\documentclass[12pt,a4paper]{amsart}
\usepackage{amssymb,amsxtra,mathtools}
\usepackage{hyperref}
\usepackage[all,cmtip]{xy}
\usepackage{tikz-cd}

\calclayout

\allowdisplaybreaks

\newcommand{\+}{\protect\nobreakdash-}
\renewcommand{\:}{\colon}

\newcommand{\rarrow}{\longrightarrow}
\newcommand{\larrow}{\longleftarrow}
\newcommand{\ot}{\otimes}
\newcommand{\ocn}{\odot}
\renewcommand{\d}{\partial}
\newcommand{\lan}{\langle}
\newcommand{\ran}{\rangle}

\newcommand{\bu}{{\text{\smaller\smaller$\scriptstyle\bullet$}}}
\newcommand{\lrarrow}{\mskip.5\thinmuskip\relbar\joinrel\relbar\joinrel
 \rightarrow\mskip.5\thinmuskip\relax}

\newcommand{\ovot}{\mathbin{\overline\otimes}}

\DeclareMathOperator{\Hom}{Hom}
\DeclareMathOperator{\Ext}{Ext}
\DeclareMathOperator{\Tor}{Tor}
\DeclareMathOperator{\Br}{Bar}
\DeclareMathOperator{\Id}{Id}
\DeclareMathOperator{\id}{id}
\DeclareMathOperator{\im}{im}
\DeclareMathOperator{\coker}{coker}
\DeclareMathOperator{\cone}{cone}
\DeclareMathOperator{\Diff}{Diff}
\newcommand{\CDiff}{\Diff^\cry}
\DeclareMathOperator{\Sym}{Sym}
\DeclareMathOperator{\Tot}{Tot}
\DeclareMathOperator{\End}{End}

\newcommand{\rings}{{\operatorname{\mathsf{--rings}}}}
\newcommand{\modl}{{\operatorname{\mathsf{--mod}}}}
\newcommand{\modr}{{\operatorname{\mathsf{mod--}}}}
\newcommand{\comodl}{{\operatorname{\mathsf{--comod}}}}
\newcommand{\comodr}{{\operatorname{\mathsf{comod--}}}}
\newcommand{\contra}{{\operatorname{\mathsf{--contra}}}}

\newcommand{\red}{{\operatorname{\mathsf{-red}}}}
\newcommand{\ind}{{\operatorname{\mathsf{-ind}}}}
\newcommand{\coind}{{\operatorname{\mathsf{-coind}}}}

\newcommand{\modrgr}{{\operatorname{\mathsf{mod_{gr}--}}}}
\newcommand{\comodrgr}{{\operatorname{\mathsf{comod_{gr}--}}}}

\newcommand{\Rings}{\mathsf{Rings}}
\newcommand{\gr}{\mathrm{gr}}
\newcommand{\q}{\mathrm{q}}
\newcommand{\rop}{{\mathrm{op}}}
\newcommand{\cry}{{\mathrm{cr}}}
\newcommand{\g}{\mathfrak g}

\newcommand{\sop}{{\mathsf{op}}}
\newcommand{\sgr}{{\mathsf{gr}}}
\newcommand{\cdg}{{\mathsf{cdg}}}
\newcommand{\fil}{{\mathsf{fil}}}
\newcommand{\wnlq}{{\mathsf{wnlq}}}
\renewcommand{\rq}{{\mathsf{rq}}}
\newcommand{\sg}{{\mathsf{sg}}}
\newcommand{\sgsm}{{\mathsf{sgsm}}}
\newcommand{\sm}{{\mathsf{sm}}}
\newcommand{\laug}{{\mathsf{laug}}}
\newcommand{\dg}{{\mathsf{dg}}}
\newcommand{\qdg}{{\mathsf{qdg}}}
\newcommand{\nlk}{{\mathsf{nlk}}}
\newcommand{\rk}{{\mathsf{rk}}}
\newcommand{\co}{{\mathsf{co}}}
\newcommand{\ctr}{{\mathsf{ctr}}}
\newcommand{\sico}{{\mathsf{sico}}}
\newcommand{\sictr}{{\mathsf{sictr}}}
\newcommand{\proj}{{\mathsf{proj}}}
\newcommand{\inj}{{\mathsf{inj}}}
\newcommand{\free}{{\mathsf{free}}}
\newcommand{\cofr}{{\mathsf{cofr}}}

\newcommand{\DG}{\mathsf{DG}}
\newcommand{\Hot}{\mathsf{Hot}}
\newcommand{\Acycl}{\mathsf{Acycl}}

\newcommand{\boQ}{\mathbb Q}
\newcommand{\boZ}{\mathbb Z}
\newcommand{\boL}{\mathbb L}
\newcommand{\boR}{\mathbb R}
\newcommand{\boM}{\mathbb M}
\newcommand{\boC}{\mathbb C}

\newcommand{\sA}{\mathsf A}
\newcommand{\sC}{\mathsf C}
\newcommand{\sD}{\mathsf D}
\newcommand{\sF}{\mathsf F}

\newcommand{\tA}{\widetilde A}
\newcommand{\tB}{\widetilde B}
\newcommand{\tC}{\widetilde C}
\newcommand{\tV}{\widetilde V}
\newcommand{\tJ}{\widetilde J}
\newcommand{\tI}{\widetilde I}

\newcommand{\hA}{\widehat A}
\newcommand{\hB}{\widehat B}
\newcommand{\hC}{\widehat C}
\newcommand{\hD}{\widehat D}
\newcommand{\hI}{\widehat I}
\newcommand{\hJ}{\widehat J}

\newcommand{\hi}{\hat\imath}
\newcommand{\hj}{\hat\jmath}
\newcommand{\hf}{\hat f}
\newcommand{\hg}{\hat g}

\newcommand{\prA}
 {\mskip.3\thinmuskip\prescript{\prime\mskip-2.5\thinmuskip}{}A}
\newcommand{\secA}
 {\mskip.3\thinmuskip\prescript{\prime\prime\mskip-2.5\thinmuskip}{}A}
\newcommand{\prI}
 {\mskip.3\thinmuskip\prescript{\prime\mskip-1.2\thinmuskip}{}I}
\newcommand{\secI}
 {\mskip.3\thinmuskip\prescript{\prime\prime\mskip-1.2\thinmuskip}{}I}
\newcommand{\prV}{\prescript{\prime}{}V}
\newcommand{\secV}{\prescript{\prime\prime}{}V}
\newcommand{\prB}{\prescript{\prime\mskip-1.4\thinmuskip}{}B}
\newcommand{\secB}{\prescript{\prime\prime\mskip-1.4\thinmuskip}{}B}
\newcommand{\prC}{\prescript{\prime\mskip-1\thinmuskip}{}C}
\newcommand{\secC}{\prescript{\prime\prime\mskip-1\thinmuskip}{}C}

\newcommand{\prtA}{\prescript{\prime\mskip-1.2\thinmuskip}
 {}{\vphantom{A}}\widetilde A}
\newcommand{\sectA}{\prescript{\prime\prime\mskip-1.2\thinmuskip}
 {}{\vphantom{A}}\widetilde A}
\newcommand{\prtV}{\prescript{\prime}{}{\vphantom{V}}\widetilde V}
\newcommand{\sectV}{\prescript{\prime\prime}{}{\vphantom{V}}\widetilde V}
\newcommand{\prhI}
 {\prescript{\prime\mskip-.4\thinmuskip}{}{\vphantom{I}}\widehat I}
\newcommand{\sechI}
 {\prescript{\prime\prime\mskip-.4\thinmuskip}{}{\vphantom{I}}\widehat I}
\newcommand{\prhB}
 {\prescript{\prime\mskip-.7\thinmuskip}{}{\vphantom{B}}\widehat B}
\newcommand{\sechB}
 {\prescript{\prime\prime\mskip-.7\thinmuskip}{}{\vphantom{B}}\widehat B}

\newcommand{\Ksp}{K\spcheck}

\newcommand{\shB}{{}^\#\!B}
\newcommand{\shb}{{}^\#\!b}
\newcommand{\shhB}{{}^\#\!\hB}
\newcommand{\shd}{{}^\#\!d}
\newcommand{\shh}{{}^\#\!h}
\newcommand{\sha}{{}^\#\!a}
\newcommand{\shdelta}{{}^\#\!\delta}

\newcommand{\Section}[1]{\bigskip\section{#1}\medskip}
\theoremstyle{plain}
\newtheorem{thm}{Theorem}[section]

\newtheorem{lem}[thm]{Lemma}
\newtheorem{prop}[thm]{Proposition}
\newtheorem{cor}[thm]{Corollary}
\theoremstyle{definition}
\newtheorem{ex}[thm]{Example}
\newtheorem{exs}[thm]{Examples}
\newtheorem{rem}[thm]{Remark}
\newtheorem{rems}[thm]{Remarks}

\begin{document}

\title{Relative nonhomogeneous Koszul duality}

\author{Leonid Positselski}

\address{Institute of Mathematics of the Czech Academy of Sciences \\
\v Zitn\'a~25, 115~67 Praha~1 (Czech Republic); and
\newline\indent Laboratory of Algebra and Number Theory \\
Institute for Information Transmission Problems \\
Moscow 127051 (Russia)}

\email{positselski@math.cas.cz}

\begin{abstract}
 This paper contains a detailed exposition of the nonhomogeneous
Koszul duality theory in the relative situation over a noncentral,
noncommutative, nonsemisimple base ring, as announced
in~\cite[Section~0.4]{Psemi}.
 We prove the Poincar\'e--Birkhoff--Witt theorem in this context and
construct the triangulated equivalences of derived Koszul duality.
 The duality between the ring of differential operators and
the de~Rham DG\+algebra, with the ring of functions as the base ring,
is the thematic example.
 The moderate generality level makes the exposition in this paper more
accessible than the very heavily technical~\cite[Chapter~11]{Psemi}.
\end{abstract}

\maketitle

\tableofcontents

\bigskip
\section*{Introduction}
\medskip

\setcounter{subsection}{-1}
\subsection{{}} \label{introd-Koszul-duality-is}
 Let $A$ be an associative ring and $R\subset A$ be a subring.
 \emph{Derived Koszul duality} is the functor $\Ext^*_A({-},R)$,
or $\Tor_*^A(R,{-})$, or $\Ext^*_A(R,{-})$, enhanced to
an equivalence of derived categories of modules.

 The above definition raises many questions.
 To begin with, $R$ is not an $A$\+module.
 So what does this $\Ext$ and $\Tor$ notation even \emph{mean}?

 Secondly, let us consider the simplest example where $R=k$ is
a field and $A=k[x]$ is the algebra of polynomials in one variable.
 Then $k$~indeed can be viewed as an $A$\+module.
 There are many such module structures, indexed by elements~$a$
of the field~$k$: given $a\in k$, one can let the generator $x\in A$
act in~$k$ by the multiplication with~$a$.
 Denote the resulting $A$\+module by~$k_a$.

 To be specific, let us choose $k=k_0$ as our preferred $A$\+module
structure on~$k$.
 Then the functors $\Ext^*_A({-},k_0)$, \ $\Tor^*_A(k_0,{-})$, and
$\Ext^*_A({-},k_0)$ are indeed well-defined on the category of
$A$\+modules.
 But these functors are far from being faithful or conservative: all
of them annihilate the $A$\+modules~$k_a$ with $a\ne0$.
 How, then, can one possibly hope to enhance such cohomological
functors to derived equivalences?

\subsection{{}}
 Koszul duality has to be distinguished from
the \emph{comodule-contramodule correspondence}, which is a different,
though related, phenomenon.

 In the simplest possible form, the comodule-contramodule correspondence
is the functor $\Ext_A^*({-},A)$ enhanced to a derived equivalence
(while Koszul duality is $\Ext_A^*({-},k)$, where $k$~is the ground
field).
 In a more realistic covariant and relative situation, comparable to
the discussion of Koszul duality in
Section~\ref{introd-Koszul-duality-is}, the comodule-contramodule
correspondence would be a derived equivalence enhancement of a functor
like $\Ext^*_A(\Hom_R(A,R),{-})$ or $\Tor_*^A({-},\Hom_R(A,R))$.

\subsection{{}}
 In the present author's research, the desire to understand Koszulity
and Koszul duality was the starting point.
 Then the separate existence and importance of comodule-contramodule
correspondence was realized, particularly in the context of
semi-infinite homological algebra~\cite{Psemi}.
 The derived nonhomogeneous Koszul duality over a field was formulated
as a ``Koszul triality'' picture, which is a triangle diagram of
derived equivalences with the comodule-contramodule correspondence
present as one side of the triangle and two versions of Koszul duality
as two other sides~\cite{Pkoszul}.

 The comodule-contramodule correspondence, its various versions,
generalizations, and philosophy, are now discussed in several books
and papers of the present author,
including~\cite{Psemi,Pkoszul,Prev,Pmgm,Pps,PS1} and others.
 On the other hand, the derived nonhomogeneous Koszul duality over
a field attracted interest of a number of authors, starting from early
works~\cite{Hin,Lef,Kel} and to very recent, such as~\cite{CLM,Maun};
there is even an operadic version of it in~\cite{HM}.

 Still, there is a void in the literature concerning \emph{relative}
nonhomogeneous Koszul duality.
 Presently, the only source of information on this topic known to
this author is his own book~\cite{Psemi}, which contains an introductory
discussion without proofs or details in~\cite[Section~0.4]{Psemi}
and a heavily technical treatment in a very general and complicated
setting in~\cite[Chapter~11]{Psemi}.
 (The memoir~\cite{Pweak} represents a very different point of view.)
 The present paper is intended to fill the void by providing
a reasonably accessible, detailed exposition on a moderate generality
level.

 Let us emphasize that relative nonhomogeneous Koszul duality is
important.
 In addition to the presence of very natural examples such as
the duality between the ring of differential operators and
the de~Rham DG\+algebra (see Section~\ref{introd-D-Omega} below),
relative nonhomogeneous Koszul duality plays a crucial role
in the semi-infinite (co)homology theory, as it was first pointed out
in~\cite{Ar}.
 This idea was subsequently developed and utilized
in~\cite[Section~11.9 and Appendix~D]{Psemi}.

 The special case of triangulated equivalences between complexes of
modules over the rings/sheaves of differential operators and DG\+modules
over the de~Rham DG\+algebra has been considered in~\cite{Kap}
and~\cite[Section~7.2]{BD2}.
 Our own treatment of it is presented in~\cite[Appendix~B]{Pkoszul}.

\subsection{{}}
 Let us start to explain the meaning of the terms involved.
 In the notation of Section~\ref{introd-Koszul-duality-is},
\emph{relative} means that $R$ is an arbitrary ring rather than
simply the ground field.
 \emph{Homogeneous} Koszul duality means that
$A=\bigoplus_{n=0}^\infty A_n$ is a nonnegatively graded ring
and $R=A_0$ is the degree-zero grading component.
 In this case, $R$ is indeed naturally both a left and a right
$R$\+module, so the meaning of the $\Ext$ and $\Tor$ notation in
Section~\ref{introd-Koszul-duality-is} is clear.
 \emph{Nonhomogeneous} Koszul duality is the situation when there is
no such grading on the ring~$A$.

 The main specific aspect of the homogeneous case is that one can
consider graded $A$\+modules with a bounding condition on the grading,
that is, only positively graded or only negatively graded modules.
 If $M$ is a positively graded left $A$\+module, then
$R\ot_AM=0$ implies $M=0$, while if $P$ is a negatively graded
left $A$\+module, then $\Hom_A(R,P)=0$ implies $P=0$.
 Hence the second problem described in
Section~\ref{introd-Koszul-duality-is} does not occur, either.

 In the nonhomogeneous situation, the solution to the second problem
from Section~\ref{introd-Koszul-duality-is} is to consider
\emph{derived categories of the second kind}.
 This means that certain complexes or DG\+modules are declared to be
nonzero objects in the derived category even though their cohomology
modules vanish.

 As to the first problem, it may well happen that $R$ has a (left or
right) $A$\+module structure even though $A$ is not graded.
 When such a module structure (extending the natural $R$\+module
structure on $R$) has been chosen, one says that the ring $A$ is
\emph{augmented}.
 In this case, the related $\Ext$ or $\Tor$ functor is well-defined.
 One wants to enhance it to a functor with values in DG\+modules
over a suitable DG\+ring in such a way that it would induce
a triangulated equivalence.

 Generally speaking, the solution to the first problem is to consider
\emph{curved DG\+modules} (\emph{CDG\+modules}), whose cohomology
modules are \emph{undefined}.
 So the $\Ext$ or $\Tor$ itself has no meaning, but the related curved
DG\+module has.
 In the augmented case, this DG\+module becomes uncurved, and indeed
computes the related $\Ext$ or $\Tor$.

\subsection{{}} \label{introd-assumptions-and-results-are}
 Let us now begin to state what our assumptions and results are.
 We assume that a ring $\tA$ is endowed with an increasing filtration
$R=F_0\tA\subset F_1\tA\subset F_2\tA\subset\dotsb$ which is exhastive
($\tA=\bigcup_nF_n\tA$) and compatible with the multiplication in~$\tA$.
 Furthermore, the successive quotients $\gr^F_n\tA=F_n\tA/F_{n-1}\tA$
are assumed to be finitely generated projective left $R$\+modules.
 Finally, the associated graded ring $A=\gr^F\tA=\bigoplus_n\gr^F_n\tA$
has to be \emph{Koszul} over its degree-zero component $A_0=R$;
this means, in particular, that the ring $A$ is generated by
its degree-one component $A_1$ over $A_0$ and defined by
relations of degree~$2$.

 In these assumptions, we assign to $(\tA,F)$ a \emph{curved DG\+ring}
(\emph{CDG\+ring}) $(B,d,h)$, which is graded by nonnegative
integers, $B=\bigoplus_{n=0}^\infty B^n$, \ $B^0=R$, has a differential
(odd derivation) $d\:B^n\rarrow B^{n+1}$ of degree~$1$, and
a \emph{curvature element} $h\in B^2$.
 The CDG\+ring $(B,d,h)$ is defined uniquely up to a unique isomorphism
of CDG\+rings, which includes the possibility of
\emph{change-of-connection} transformations.
 The grading components $B^n$ are finitely generated projective right
$R$\+modules.
 In particular, one has $B^1=\Hom_R(A_1,R)$ and
$A_1=\Hom_{R^\rop}(B^1,R)$.

 Furthermore, to any left $\tA$\+module $P$ we assign a CDG\+module
structure on the graded left $B$\+module $B\ot_RP$, and to any
right $\tA$\+module $M$ we assign a CDG\+module structure on the graded
right $B$\+module $\Hom_{R^{\rop}}(B,M)$.
 These constructions are then extended to complexes of left and
right $\tA$\+modules $P^\bu$ and $M^\bu$, assigning to them left
and right CDG\+modules $B\ot_RP^\bu$ and $\Hom_{R^\rop}(B,M^\bu)$
over $(B,d,h)$.
 A certain (somewhat counterintuitive) way to totalize bigraded modules
is presumed here.
 The resulting functors induce the derived equivalences promised in
Section~\ref{introd-Koszul-duality-is}.
 The functor $P^\bu\longmapsto B\ot_RP^\bu$ is a CDG\+enhancement of
the (possibly nonexistent) $\Ext^*_{\tA}(R,P)$, and the functor
$M^\bu\longmapsto\Hom_{R^\rop}(B,M^\bu)$ is a CDG\+enhancement of
the (possibly nonexistent) $\Tor_*^{\tA}(M,R)$.
 However, there are further caveats.

\subsection{{}} \label{introd-comodules-contramodules}
 One important feature of the nonhomogeneous Koszul duality over a field,
as developed in the memoir~\cite{Pkoszul}, is that it connects modules
with comodules or contramodules.
 In fact, the ``Koszul triality'' of~\cite{Pkoszul} connects modules with
comodules \emph{and} contramodules.
 In the context of relative nonhomogeneous Koszul duality theory in
the full generality of the present paper, the Koszul triality picture
splits into two separate dualities.
 A certain exotic derived category of right $\tA$\+modules is equivalent 
to an exotic derived category of \emph{right $B$\+comodules}, while
another exotic derived category of left $\tA$\+modules is equivalent to
an exotic derived category of \emph{left $B$\+contramodules}.
 The triality picture is then restored under some additional assumptions
(namely, two-sided locally finitely generated projectivity of
the filtration $F$ and finiteness of homological dimension of
the base ring~$R$).

 What are the ``comodules'' and ``contramodules'' in our context?
 First of all, we have complexes of $\tA$\+modules on the one side and
CDG\+modules over $B$ on the other side; so both the comodules and
the contramodules are graded $B$\+modules.
 In fact, the (graded) right $B$\+comodules are a certain \emph{full
subcategory} in the graded right $B$\+modules, and similarly
the (graded) left $B$\+contramodules are a certain \emph{full
subcategory} in the graded left $B$\+modules.

 Which full subcategory?
 A graded right $B$\+module $N$ is called a \emph{graded
right $B$\+comodule} if for every element $x\in N$ there exists
an integer $m\ge1$ such that $xB^n=0$ for all $n\ge m$.
 The definition of $B$\+contramodules is more complicated and, as
usually, involves certain infinite summation operations.
 A graded left $B$\+module $Q$ is said to be a \emph{graded left
$B$\+contramodule} if, for every integer~$j$, every sequence of
elements $q_n\in Q^{j-n}$, \ $n\ge0$, and every sequence of
elements $b_n\in B^n$, an element denoted formally by
$\sum_{n=0}^\infty b_nq_n\in Q^j$ is defined.
 One imposes natural algebraic axioms on such infinite summation
operations, and then proves that an infinite summation structure on
a given graded left $B$\+module $Q$ is unique if it exists.

 In particular, this discussion implies that (somewhat
counterintuitively), in the notation of
Section~\ref{introd-assumptions-and-results-are}, the bigraded
module $\Hom_{R^\rop}(B,M^\bu)$ has to be totalized by taking
infinite \emph{direct sums} along the diagonals (to obtain a graded
right $B$\+comodule), while the bigraded module $B\ot_RP^\bu$ needs
to be totalized by taking infinite \emph{products} along the diagonals
(to obtain a graded left $B$\+contramodule).

\subsection{{}} \label{introd-coring-is-better}
 The explanation for the counterintuitive totalization procedures
mentioned in Section~\ref{introd-comodules-contramodules}, from
our perspective, is that $B$ is a ``fake'' graded ring.
 It really ``wants'' to be a coring, but this point of view is hard
to fully develop.
 It plays a key role, however, in (at least) one of our two proofs of
the Poincar\'e--Birkhoff--Witt theorem for nonhomogeneous Koszul rings.

 The graded coring in question is $C=\Hom_{R^\rop}(B,R)$, that is,
the result of applying the dualization functor $\Hom_{R^\rop}({-},R)$
to the graded ring~$B$.
 The point is that we have already done one such dualization when
we passed from the $R$\+$R$\+bimodule $A_1$ to the $R$\+$R$\+bimodule
$B^1=\Hom_R(A_1,R)$, as mentioned in
Section~\ref{introd-assumptions-and-results-are}.
 The two dualization procedures are essentially inverse to each other,
so the passage to the coring $C$ over $R$ returns us to the undualized
world, depending covariantly-functorially on the ring~$A$.

 Experience teaches that the passage to the dual vector space is
better avoided in derived Koszul duality.
 This is the philosophy utilized in the memoir~\cite{Pkoszul}
and the book~\cite{Psemi}.
 This philosophy strongly suggests that the graded coring $C$ is
preferable to the graded ring $B$ as a Koszul dual object to
a Koszul graded ring~$A$.

 The problem arises when we pass to the nonhomogeneous setting.
 In the context of the discussion in
Section~\ref{introd-assumptions-and-results-are},
the odd derivation~$d$, which is a part of the structure of
a CDG\+ring $(B,d,h)$, is \emph{not} $R$\+linear.
 In fact, the restriction of~$d$ to the subring $R=B^0\subset B$
may well be nonzero, and in the most interesting cases it is.
 This is a distinctive feature of the relative nonhomogeneous Koszul
duality.
 So how does one apply the functor $\Hom_{R^\rop}({-},R)$ to
a non-$R$-linear map?

\subsection{{}} \label{introd-D-Omega}
 The duality between the ring of differential operators and
the de~Rham DG\+algebra of differential forms is the thematic
example of relative nonhomogeneous Koszul duality.
 Let $X$ be a smooth affine algebraic variety over a field of
characteristic~$0$ (or a smooth real manifold).
 Let $O(X)$ denote the ring of functions and $\Diff(X)$ denote
the ring of differential operators on~$X$.
 Endow the ring $\Diff(X)$ with an increasing filtation $F$ by
the order of the differential operators.
 So the associated graded ring $\Sym_{O(X)}(T(X))=\gr^F\Diff(X)$ is
the symmetric algebra of the $O(X)$\+module $T(X)$ of
vector fields on~$X$.

 In this example, $R=O(X)$ is our base ring, $\tA=\Diff(X)$ is 
our nonhomogeneous Koszul ring over $R$, and $A=\Sym_{O(X)}(T(X))$
is the related homogeneous Koszul ring.
 The graded ring Koszul to $A$ over $R$ is the graded ring of
differential forms $B=\Omega(X)$.
 There is no curvature in the CDG\+ring $(B,d,h)$ (one has $h=0$;
a nonzero curvature appears when one passes to the context of
\emph{twisted} differential operators, e.~g., differential
operators acting in the sections of a vector bundle $E$ over~$X$;
see~\cite[Section~0.4.7]{Psemi} or~\cite[Appendix~B]{Pkoszul}).
 The differential $d\:B\rarrow B$ is the de~Rham differential,
$d=d_{dR}$; so $(B,d)$ is a DG\+algebra over~$k$.

 But the de~Rham DG\+algebra is not a DG\+algebra over $O(X)$
(and neither the ring $\Diff(X)$ is an algebra over~$O(X)$).
 In fact, the restriction of the de~Rham differential to
the subring $O(X)\subset\Omega(X)$ is quite nontrivial.

\subsection{{}}
 So the example of differential operators and differential forms
is a case in point for the discussion in
Section~\ref{introd-coring-is-better}.
 In this example, $C=\Hom_{O(X)}(\Omega(X),O(X))$ is the graded
coring of polyvector fields over the ring of functions on~$X$.
 Certainly there is no de~Rham differential on polyvector fields.
 What structure on polyvector fields corresponds to the de~Rham
differential on the forms?

 Here is what we do.
 We adjoin an additional generator~$\delta$ to the de~Rham
DG\+ring $(\Omega(X),d_{dR})$, or more generally to the underlying
graded ring $B$ of a CDG\+ring $(B,d,h)$.
 The new generator~$\delta$ is subject to the relations
$[\delta,b]=d(b)$ for all $b\in B$ (where the bracket denotes
the graded commutator) and $\delta^2=h$.
 Then there is a new differential on the graded ring $\hB=B[\delta]$,
which we denote by $\d=\d/\d\delta$.

 The differential~$\d$ is $R$\+linear (and more generally,
$B$\+linear with signs), so we can dualize it, obtaining
a coring $\hC=\Hom_{R^\rop}(\hB,R)$ with the dual differential
$\Hom_{R^\rop}(\d,R)$.
 This is the structure that was called a \emph{quasi-differential
coring} in~\cite{Psemi}.
 It plays a key role in the exposition in~\cite[Chapter~11]{Psemi}.
 
 Of course, the odd derivation $\d=\d/\d\delta$ is acyclic, and
so is the dual odd coderivation on the coring~$\hC$.
 This may look strange; but in fact, this is how it should be.
 Recall that we started with a curved DG\+ring $(B,d,h)$.
 Its differential~$d$ does not square to zero, and its cohomology
is undefined.
 So there is no cohomology ring in the game, and it is not supposed
to suddenly appear from the construction.

\subsection{{}}
 Now, how does one assign a derived category of modules to the acyclic
DG\+ring $(\hB,\d)$\,?
 The related constructions are discussed in~\cite[Section~11.7]{Psemi}.
 A \emph{quasi-differential module} over $(\hB,\d)$ is simply
a graded $\hB$\+module, without any differential.
 Such modules form a DG\+category.
 In fact, a DG\+module over $(\hB,\d)$ is the same thing as
a \emph{contractible} object of the DG\+category of quasi-differential
modules.
 This point of view, adopted in~\cite[Chapter~11]{Psemi} in the context
of quasi-differential comodules and contramodules over
quasi-differential corings, is so counterintuitive that one is having
a hard time with what otherwise are very simple constructions.
 We have none of that in this paper, using the equivalent, but much
more tractable concept of a CDG\+module over the CDG\+ring $(B,d,h)$.

 Some words about the \emph{coderived} and \emph{contraderived
categories} are now in order.
 These are the most important representatives of the class of
constructions known as the ``derived categories of the second kind''.

 In the spirit of the discussion in
Section~\ref{introd-comodules-contramodules}, we consider
right CDG\+comodules and left CDG\+contramodules over $(B,d,h)$.
 These are certain full subcategories in the DG\+categories of,
respectively, right and left CDG\+modules over the CDG\+ring $(B,d,h)$.
 Following the general definitions in~\cite{Psemi,Pkoszul},
the \emph{coderived category of right CDG\+comodules} over $(B,d,h)$
is constructed as the Verdier quotient category of the homotopy category
of CDG\+comodules by its minimal triangulated subcategory containing
the total CDG\+comodules of all the short exact sequences of
CDG\+comodules and closed under infinite direct sums.
 Similarly, the \emph{contraderived category of left CDG\+contramodules}
over $(B,d,h)$ is the Verdier triangulated quotient category of
the homotopy category of CDG\+contramodules by its minimal triangulated
subcategory containing the total CDG\+contramodules of all the short
exact sequences of CDG\+contramodules over $(B,d,h)$ and closed under
infinite products.

\subsection{{}}
 When the base ring $R$ has finite right homological dimension, our
derived Koszul duality result simply states that the derived category
of right $\tA$\+modules is equivalent to the coderived category of
right CDG\+comodules over $(B,d,h)$.
 When the ring $R$ has finite left homological dimension, one similarly
has a natural equivalence between the derived category of left
$\tA$\+modules and the contraderived category of left CDG\+contramodules
over $(B,d,h)$.

 The situation gets more complicated when the homological dimension of
$R$ is infinite.
 In this case, following the book~\cite{Psemi} and the paper~\cite{Pfp},
one can consider the \emph{semiderived categories} of right and left
$\tA$\+modules, or more precisely the \emph{semicoderived category}
of right $\tA$\+modules relative to $R$ and the \emph{semicontraderived
category} of left $\tA$\+modules relative to~$R$.
 These are defined as the Verdier quotient categories of the homotopy
categories of complexes of right and left $\tA$\+modules by
the triangulated subcategories of complexes that are, respectively,
coacyclic or contraacylic \emph{as complexes of $R$\+modules}.

 Then the derived Koszul duality theorem tells that the semicoderived
category of right $\tA$\+modules is equivalent to the coderived
category of right CDG\+comodules over $(B,d,h)$; and
the semicontraderived category of left $\tA$\+modules is equivalent
to the contraderived category of left CDG\+contramodules over $(B,d,h)$.

 One can also describe the derived category of right $\tA$\+modules
as the quotient category of the coderived category of right
CDG\+comodules over $(B,d,h)$ by its minimal triangulated subcategory
closed under direct sums and containing all the CDG\+comodules $(N,d_N)$
such that $NB^i=0$ for $i>0$ and $N$ is acyclic with respect to
the differential~$d_N$ (where $d_N^2=0$ since $Nh=0$).
 Simlarly, the derived category of left $\tA$\+modules is equivalent
to the quotient category of the contraderived category of left
CDG\+contramodules over $(B,d,h)$ by its minimal triangulated
subcategory closed under products and containing all
the CDG\+contramodules $(Q,d_Q)$ such that $B^iQ=0$ for $i>0$ and
$Q$ is acyclic with respect to the differential~$d_Q$.

\subsection{{}} \label{introd-conversion}
 A basic fact of the classical theory of modules over the rings of
differential operators $\Diff(X)$ is that the abelian categories of
left and right $\Diff(X)$\+modules are naturally equivalent
to each other.
 This is a rather nontrivial equivalence, in that the ring $\Diff(X)$
is \emph{not} isomorphic to its opposite ring.

 The classical \emph{conversion functor} $\Diff(X)\modl\rarrow
\modr{\Diff(X)}$ assigns to a left $\Diff(X)$\+module $M$ a natural
right $\Diff(X)$\+module structure on the tensor product
$\Omega^m(X)\ot_{O(X)}M$, where $m=\dim X$ and $\Omega^m(X)$ is
the $O(X)$\+module of global sections of the line bundle of
differential forms of the top degree on~$X$.
 The inverse conversion $\modr{\Diff(X)}\rarrow\Diff(X)\modl$ is
performed by taking the tensor product over $O(X)$ with the (module of
global sections of) the line bundle of top polyvector fields
$\Hom_{O(X)}(\Omega^m(X),O(X))=\Lambda^m_{O(X)}(T(X))$.

 The present paper offers an interpretation of the conversion functor
in the context of relative nonhomogeneous Koszul duality.
 Let $(B,d,h)$ be a nonnegatively graded CDG\+ring with
the (possibly noncommutative) degree-zero component $B^0=R$.
 Assume that the grading components of $B$ are finitely generated
projective left and right $R$\+modules, there is an integer
$m\ge0$ such that $B^n=0$ for $n>m$, the $R$\+$R$\+bimodule
$B^m$ is invertible, and the multiplication maps $B^n\ot_RB^{m-n}
\rarrow B^m$ are perfect pairings.
 Assume further that $B$ is a Koszul graded ring over~$R$.
 Then we say that $B=(B,d,h)$ is a \emph{relatively Frobenius
Koszul CDG\+ring}.

 As the grading components of $B$ are finitely generated and
projective over $R$ on both sides, there are \emph{two} nonhomogeneous
Koszul dual filtered rings to $(B,d,h)$, one on the left side
and one on the right side; we denote them by $\tA$ and~$\tA^\#$.
 Then the claim is that, whenever $B$ is relatively Frobenius
over~$R$, the two rings $\tA$ and $\tA^\#$ are Morita equivalent.
 The tensor product with the invertible $R$\+$R$\+bimodule
$T=B^m$ transforms any left $\tA$\+module into a left $\tA^\#$\+module,
and any right $\tA^\#$\+module into a right $\tA$\+module.
 The functors $\Hom_R(T,{-})$ and $\Hom_{R^\rop}(T,{-})$ provide
the inverse transformations.
 (When the graded ring $B$ is graded commutative and $h=0$, the ring
$\tA^\#$ is simply the opposite ring to the ring~$\tA$.)

 In this context, assuming additionally that the ring $R$ has finite
left homological dimension, we even obtain a ``Koszul quadrality''
picture.
 This means a commutative diagram of triangulated equivalences between
four (conventional or exotic) derived categories: the derived category
of left $\tA^\#$\+modules, the derived category of left $\tA$\+modules,
the coderived category of left CDG\+modules over $(B,d,h)$, and
the contraderived category of left CDG\+modules over $(B,d,h)$.

\subsection{{}}
 We discuss the homogeneous quadratic duality over a base ring
in Section~\ref{quadratic-duality-secn}, flat and finitely projective
Koszul graded rings over a base ring in
Section~\ref{koszulity-secn}, relative nonhomogeneous quadratic
duality in Section~\ref{nonhomogeneous-quadratic-secn},
and the Poincar\'e--Birkhoff--Witt theorem for
nonhomogeneous Koszul rings over a base ring in
Section~\ref{pbw-secn}.
 The discussion of comodules and contramodules over graded rings
in Section~\ref{comodules-and-contramodules-secn} prepares ground for
the derived Koszul duality for module categories, which is worked out
on the comodule side in Section~\ref{comodule-side-secn} and on
the contramodule side in Section~\ref{contramodule-side-secn}.
 The comodule-contramodule correspondence, connecting the comodule
and contramodule sides of the Koszul duality, is developed in
Section~\ref{co-contra-secn}.
 The interpretation of the conversion functor in terms of Koszul
duality is discussed in Section~\ref{conversion-secn}.

 Examples of relative nonhomogeneous Koszul duality are offered
in Section~\ref{examples-secn}.
 These are various species of differential operators, to which
correspond the related species of differential forms.
 We consider algebraic differential operators over smooth affine
varieties in characteristic~$0$, crystalline differential operators
over smooth affine varieties in arbitrary characteristic, differential
operators acting in the sections of a vector bundle, and differential
operators twisted with a chosen closed $2$\+form.
 Passing from the algebraic to the analytic setting, we discuss smooth
differential operators on a smooth compact real manifold and
$\bar\d$\+differential operators in the Dolbeault theory on a compact
complex manifold.
 Returning to the algebraic context, we consider relative differential
forms and differential operators for a morphism of commutative rings,
Lie algebroids with their enveloping algebras and cohomological
Chevalley--Eilenberg complexes, and finally noncommutative differential
forms for a morphism of noncommutative rings.
 In the latter situation, the related ring of ``noncommutative
differential operators'' is simply the ring of all endomorphisms of
the bigger ring as a module over the subring (endowed with the obvious
two-step filtration).
 For the benefit of the reader, we have tried to make our exposition
of these examples from various areas of algebra and geometry reasonably
self-contained with many background details included.

\subsection*{Acknowledgment}
 Parts of the material presented in this paper go back more than
a quarter century.
 This applies to the content of
Sections~\ref{quadratic-duality-secn}\+-\ref{koszulity-secn}
and the computations in Section~\ref{nonhomogeneous-quadratic-secn}
(with the notable exception of the $2$\+category story), which
I~worked out sometime around 1992.
 The particular case of duality over a field, which is much less
complicated, was presented in the paper~\cite{Pcurv}, and
the possibility of extension to the context of a base ring
was mentioned in~\cite[beginning of Section~4]{Pcurv}.
 The main results in
Sections~\ref{comodule-side-secn}\+-\ref{contramodule-side-secn}
go back to Spring~2002.
 Subsequently, I~planned and promised several times over the years
to write up a detailed exposition.
 This paper partially fulfills that promise.
 The paper also contains some much more recent results;
this applies, first of all, to the material of
Section~\ref{comodules-and-contramodules-secn}, which is largely
based on~\cite[Section~6]{Pcoun} or~\cite[Theorem~3.1]{Pper}.
 I~would like to thank all the people, too numerous to be mentioned here
by name, whose help and encouragement contributed to my survival over
the decades.
 Speaking of more recent events, I~am grateful to Andrey Lazarev,
Julian Holstein, and Bernhard Keller for stimulating discussions and
interest to this work.
 The author was supported by research plan RVO:~67985840 and
the GA\v CR project 20-13778S when writing the paper up.

\Section{Homogeneous Quadratic Duality over a Base Ring}
\label{quadratic-duality-secn}

 All the \emph{associative rings} in this paper are unital.
 We will always presume unitality without mentioning it; so all
the left and ring \emph{modules} over associative rings are unital,
all the \emph{ring homomorphisms} take the unit to the unit, all
the \emph{subrings} contain the unit, and all the \emph{gradings} and
\emph{filtrations} are such that the unit element belongs to
the degree-zero grading/filtration component.

 Given an associative ring $R$, we denote by $R\modl$ the abelian
category of left $R$\+modules and by $\modr R$ the abelian category
of right $R$\+modules.

 Let $R$, $S$, and $T$ be three associative rings.
 For any left $R$\+modules $L$ and $M$, we denote by $\Hom_R(L,M)$
the abelian group of all left $R$\+module morphisms $L\rarrow M$.
 If $L$ is an $R$\+$S$\+bimodule and $M$ is an $R$\+$T$\+bimodule,
then the group $\Hom_R(L,M)$ acquires a natural structure of
$S$\+$T$\+bimodule.
 Similarly, for any right $R$\+modules $Q$ and $N$, the abelian group
of all right $R$\+module morphisms $Q\rarrow N$ is denoted by
$\Hom_{R^\rop}(Q,N)$ (where $R^\rop$ stands for the ring opposite
to~$R$).
 If $Q$ is an $S$\+$R$\+bimodule and $N$ is a $T$\+$R$\+bimodule,
then $\Hom_{R^\rop}(Q,N)$ is a $T$\+$S$\+bimodule.

 In particular, for any $R$\+$S$\+bimodule $U$, the abelian group
$\Hom_R(U,R)$ is naturally an $S$\+$R$\+bimodule.
 If $U$ is a finitely generated projective left $R$\+module, then
$\Hom_R(U,R)$ is a finitely generated projective right $R$\+module.
 Similarly, for any $S$\+$R$\+bimodule $M$, the abelian group
$\Hom_{R^\rop}(M,R)$ is naturally an $R$\+$S$\+bimodule.
 If $M$ is a finitely generated projective right $R$\+module, then
$\Hom_{R^\rop}(M,R)$ is a finitely generated projective left
$R$\+module.

 For any $R$\+$S$\+bimodule $U$, there is a natural morphism of
$R$\+$S$\+bimodules $U\rarrow\Hom_{R^\rop}(\Hom_R(U,R),R)$, which is
an isomorphism whenever the left $R$\+module $U$ is finitely
generated and projective.
 For any $S$\+$R$\+bimodule $M$, there is a natural morphism of
$S$\+$R$\+bimodules $M\rarrow\Hom_R(\Hom_{R^\rop}(M,R),R)$, which is
an isomorphism whenever the right $R$\+module $M$ is finitely
generated and projective.

 Let $U$ be an $R$\+$S$\+bimodule and $V$ be an $S$\+$T$\+bimodule.
 Then the left $R$\+module $U\ot_SV$ is projective whenever
the left $R$\+module $U$ and the left $S$\+module $V$ are projective.
 The left $R$\+module $U\ot_SV$ is finitely generated whenever
the left $R$\+module $U$ and the left $S$\+module $V$ are finitely
generated.
 The similar assertions apply to the projectivity and finite
generatedness on the right side.

\begin{lem} \label{tensor-dual-lemma}
\textup{(a)} Let $U$ be an $R$\+$S$\+bimodule and $V$ be
an $S$\+$T$\+bimodule.
 Then there is a natural morphism of $T$\+$R$\+bimodules
$$
 \Hom_S(V,S)\ot_S\Hom_R(U,R)\lrarrow\Hom_R(U\ot_SV,\>R),
$$
which is an isomorphism whenever the left $S$\+module $V$ is finitely
generated and projective. \par
\textup{(b)} Let $M$ be an $S$\+$R$\+bimodule and $N$ be
a $T$\+$S$\+bimodule.
 Then there is a natural morphism of $R$\+$T$\+bimodules
$$
 \Hom_{R^\rop}(M,R)\ot_S\Hom_{S^\rop}(N,S)\lrarrow
 \Hom_{R^\rop}(N\ot_SM,\>R),
$$
which is an isomorphism whenever the right $S$\+module $N$ is finitely
generated and projective.
\end{lem}

\begin{proof}
 Part~(a): the desired map takes an element $g\ot f\in \Hom_S(V,S)\ot_S
\Hom_R(U,R)$ to the map $U\ot_SV\rarrow R$ taking an element $u\ot v$ to
the element $f(ug(v))\in R$, for any $g\in\Hom_S(V,S)$, \
$f\in\Hom_R(U,R)$, \ $u\in U$, and $v\in V$.
 The second assertion does not depend on the $T$\+module structure
on $V$, so one can assume $T=\boZ$ and, passing to the finite direct
sums and direct summands in the argument $V\in S\modl$, reduce to
the obvious case $V=S$.
 Part~(b): the desired map takes an element $h\ot k\in
\Hom_{R^\rop}(M,R)\ot_S\Hom_{S^\rop}(N,S)$ to the map $N\ot_SM\rarrow R$
taking an element $n\ot m$ to the element $h(k(n)m)\in R$, for any
$h\in\Hom_{R^\rop}(M,R)$, $k\in\Hom_{S^\rop}(N,S)$, \
$n\in N$, and $m\in M$.
 The second assertion does not depend on the $T$\+module structure
on $N$, so it reduces to the obvious case $N=S$.
\end{proof}

 Let $R$ be an associative ring and $V$ be an $R$\+$R$\+bimodule.
 The \emph{tensor ring} of $V$ over $R$ (otherwise called
the ring \emph{freely generated} by an $R$\+$R$\+bimodule~$V$)
is the graded ring $T_R(V)=\bigoplus_{n=0}^\infty T_{R,n}(V)$ with
the components $T_{R,0}(V)=R$, \ $T_{R,1}(V)=V$, \ $T_{R,2}(V)=V\ot_RV$,
and $T_{R,n}(V)=V\ot_R\dotsb\ot_RV$ ($n$~factors) for $n\ge2$.
 The multiplication in $T_R(V)$ is defined by the obvious rules
$r(v_1\ot\dotsb\ot v_n)=(rv_1)\ot v_2\ot\dotsb\ot v_n$, \
$(v_1\ot\dotsb\ot v_n)s=v_1\ot\dotsb\ot v_{n-1}\ot(v_ns)$, and
$(v_1\ot\dotsb\ot v_n)(v_{n+1}\ot\dotsb\ot v_{n+m})=
v_1\ot\dotsb\ot v_{n+m}$ for all $r$, $s\in R$ and $v_i\in V$.

 Let $A=\bigoplus_{n=0}^\infty A_n$ be a nonnegatively graded ring
with the degree-zero component $A_0=R$.
 Denote by $V$ the $R$\+$R$\+bimodule $V=A_1$.
 Then there exists a unique homomorphism of graded rings $\pi_A\:T_R(V)
\rarrow A$ acting by the identity maps on the components of degrees~$0$
and~$1$.
 The ring $A$ is said to be \emph{generated by $A_1$} (over~$A_0$) if
the map~$\pi_A$ is surjective.
 Furthermore, denote by $J_A=\ker(\pi_A)$ the kernel ideal of the ring
homomorphism~$\pi_A$.
 Then $J_A$ is a graded ideal in $T_R(V)$, so we have
$J_A=\bigoplus_{n=2}^\infty J_{A,n}$, where $J_{A,n}\subset T_{R,n}(V)$.
 Set $I_A=J_{A,2}\subset V\ot_RV$.
 A graded ring $A$ generated by $A_1$ over $A_0$ is said to be
\emph{quadratic} (over $R=A_0$) if the two-sided ideal
$J_A\subset T_R(V)$ is generated by~$I_A$, that is $J_A=(I_A)$,
or explicitly
\begin{equation} \label{quadratic-ideal-component}
 J_{A,n}=\sum\nolimits_{i=1}^{n-1}
 T_{R,i-1}(V)\cdot I_A\cdot T_{R,n-i-1}(V)
 \quad\text{for all $n\ge 3$}.
\end{equation}

 Let $A=\bigoplus_{n=0}^\infty A_n$ be a quadratic graded ring
with the degree-zero component $R=A_0$.
 We will say that $A$ is \emph{$2$\+left finitely projective} if
the left $R$\+modules $A_1$ and $A_2$ are projective and finitely
generated.
 Furthermore, $A$ is \emph{$3$\+left finitely projective} if
the same applies to the left $R$\+modules $A_1$, $A_2$, and~$A_3$.
 Similarly, a quadratic graded ring $B=\bigoplus_{n=0}^\infty B_n$
with the degree-zero component $R=B_0$ is
\emph{$2$\+right finitely projective} if the right $R$\+modules
$B_1$ and $B_2$ are finitely generated projective, and $B$ is
\emph{$3$\+right finitely projective} if the same applies to
the right $R$\+modules $B_1$, $B_2$, and~$B_3$.

 Now we fix an associative ring $R$ and consider the
\emph{category of graded rings over~$R$} $R\rings_\sgr$, defined
as follows.
 The objects of $R\rings_\sgr$ are nonnegatively graded associative
rings $A=\bigoplus_{n=0}^\infty A_n$ endowed with a fixed ring
isomorphism $R\simeq A_0$.
 Morphisms $\prA\rarrow\secA$ in $R\rings_\sgr$ are graded ring
homomorphisms forming a commutative triangle diagram with
the isomorphisms $R\simeq\prA_0$ and $R\simeq\secA_0$.
 Various specific classes of graded rings defined above in this section
(and below in the next one) are viewed as full subcategories
in $R\rings_\sgr$.

\begin{prop} \label{2-fin-proj-quadratic-duality}
 There is an anti-equivalence between the categories of\/ $2$\+left
finitely projective quadratic graded rings $A$ over $R$ and\/
$2$\+right finitely projective quadratic graded rings $B$ over $R$,
called the \emph{quadratic duality} and defined by the following rules.
 Given a ring $A$, the ring $B$ is constructed as $B=T_R(B_1)/(I_B)$,
where $B_1=\Hom_R(A_1,R)$ and $I_B=\Hom_R(A_2,R)\subset B_1\ot_R B_1$.
 Then the $R$\+$R$\+bimodule $B_2$ is naturally isomorphic to
$\Hom_R(I_A,R)$.
 Conversely, given a ring $B$, the ring $A$ is constructed as
$A=T_R(A_1)/(I_A)$, where $A_1=\Hom_{R^\rop}(B_1,R)$ and
$I_A=\Hom_{R^\rop}(B_2,R)\subset A_1\ot_R A_1$.
 Then the $R$\+$R$\+bimodule $A_2$ is naturally isomorphic to
$\Hom_{R^\rop}(I_B,R)$.
\end{prop}

\begin{proof}
 The category of quadratic graded rings $A$ over $R$ is equivalent to
the category of $R$\+$R$\+bimodules $V=A_1$ endowed with
a subbimodule $I=I_A\subset V\ot_R V$.
 Here morphisms in the category of pairs $(\prV,\prI)\rarrow
(\secV,\secI)$ are defined as $R$\+$R$\+bimodule morphisms
$f\:\prV\rarrow\secV$ such that $(f\ot f)(\prI)\subset\secI$.

 A quadratic graded ring $A$ is $2$\+left finitely projective if and
only if in the related pair $(V,I)$ the left $R$\+modules $V$ and
$A_2=(V\ot_R V)/I$ are finitely generated and projective.
 Assuming the former condition, the left $R$\+module $V\ot_RV$ is then
finitely generated and projective, too, so the latter condition is
equivalent to the $R$\+$R$\+subbimodule $I\subset V\ot_RV$ being
\emph{split as a left $R$\+submodule}.

 Now we have a short exact sequence of $R$\+$R$\+bimodules
$$
 0\lrarrow I_A\lrarrow A_1\ot_R A_1\lrarrow A_2\lrarrow0,
$$
which is split as a short exact sequence of left $R$\+modules.
 Applying the functor $\Hom_R({-},R)$ and taking into account
Lemma~\ref{tensor-dual-lemma}(a), we obtain a short exact sequence
of $R$\+$R$\+bimodules
$$
 0\rarrow\Hom_R(A_2,R)\rarrow\Hom_R(A_1,R)\ot_R\Hom_R(A_1,R)
 \rarrow\Hom_R(I_A,R)\rarrow0,
$$
which is split as a short exact sequence of right $R$\+modules.
 It remains to set $B_1=\Hom_R(A_1,R)$ and $I_B=\Hom_R(A_2,R)$,
so that $B_2=\Hom_R(I_A,R)$.
 According to the discussion in the beginning of this section,
$B_1$ and $B_2$ are finitely generated projective right $R$\+modules.
 It is clear that this construction is a contravariant functor between
the categories of $2$\+left finitely projective and $2$\+right finitely
projective quadratic graded rings over $R$, and that the similar
construction with the left and right sides switched provides
the inverse functor in the opposite direction.
\end{proof}

 The $2$\+left finitely projective quadratic ring $A$ and
the $2$\+right finitely projective quadratic ring $B$ as in
Proposition~\ref{2-fin-proj-quadratic-duality} are said to be
\emph{quadratic dual} to each other.

\begin{prop} \label{3-fin-proj-quadratic-duality}
 The anti-equivalence of categories from
Proposition~\ref{2-fin-proj-quadratic-duality} takes\/ $3$\+left
finitely projective quadratic graded rings to\/ $3$\+right finitely
projective quadratic graded rings and vice versa.
\end{prop}

\begin{proof}
 For any quadratic graded ring $A$, the grading component $A_3$
is the cokernel of the map $(A_1\ot_R I_A)\oplus(I_A\ot_RA_1)
\rarrow A_1\ot_RA_1\ot_RA_1$ induced by the inclusion map
$I_A\rarrow A_1\ot_RA_1$.
 When the components $A_1$ and $A_2$ are projective as (say, left)
$R$\+modules, the maps $A_1\ot_RI_A\rarrow A_1\ot_RA_1\ot_RA_1$
and $I_A\ot_RA_1\rarrow A_1\ot_RA_1\ot_RA_1$ are injective, so
we have a four-term exact sequence of $R$\+$R$\+bimodules
\begin{multline} \label{four-term-sequence}
 0\lrarrow I_A^{(3)}\lrarrow
 (A_1\ot_R I_A)\oplus(I_A\ot_RA_1) \\ \lrarrow A_1\ot_RA_1\ot_RA_1
 \lrarrow A_3\lrarrow0,
\end{multline}
where $I_A^{(3)}=(A_1\ot_R I_A)\cap(I_A\ot_RA_1)\subset
A_1\ot_RA_1\ot_RA_1$.
 When the component $A_3$ is a projective left $R$\+module, too, we
observe that all the terms of this exact sequence, except perhaps
the leftmost one, are projective left $R$\+modules.
 It follows that the sequence~\eqref{four-term-sequence} splits as
an exact sequence of left $R$\+modules, and the leftmost term
$I_A^{(3)}$ is a projective left $R$\+module, too.

 Furthermore, when $A_1$, $A_2$, and $A_3$ are finitely generated
projective left $R$\+modules, all the terms of
the sequence~\eqref{four-term-sequence} are also finitely generated
projective left $R$\+modules.
 Applying the functor $\Hom_R({-},R)$ to~\eqref{four-term-sequence},
we obtain a four-term exact sequence of $R$\+$R$\+bimodules
\begin{multline} \label{dual-four-term-sequence}
 0\lrarrow\Hom_R(A_3,R)\lrarrow B_1\ot_R B_1\ot_R B_1 \\ \lrarrow
 (B_1\ot_RB_2)\oplus(B_2\ot_RB_1)\lrarrow\Hom_R(I_A^{(3)},R)\lrarrow0.
\end{multline}
 Now for any quadratic graded ring $B$, the cokernel of the map
$B_1\ot_RB_1\ot_RB_1\rarrow(B_1\ot_RB_2)\oplus(B_2\ot_RB_1)$ induced
by the (surjective) multiplication map $B_1\ot_R B_1\rarrow B_2$ is
the grading component~$B_3$.
 Hence we have a natural isomorphism of $R$\+$R$\+bimodules
$B_3\simeq\Hom_R(I_A^{(3)},R)$, and it follows that $B_3$ is a finitely
generated projective right $R$\+module.
\end{proof}

 The statement similar to Proposition~\ref{3-fin-proj-quadratic-duality}
does \emph{not} hold in degrees higher than~$3$ for quadratic graded
rings in general.
 It holds under the Koszulity assumption, though, as we will see in
the next section.

\begin{rem} \label{degree-zero-iso-remark}
 The above discussion of the categories of nonnegatively graded and
quadratic rings can be modified or expanded by including
\emph{non-identity isomorphisms} (in particular, automorphisms)
in the degree-zero component. 
 Denote by $\Rings_\sgr$ the category whose objects are nonnegatively
graded associative rings $A=\bigoplus_{n=0}^\infty A_n$, and morphisms
are defined as follows.
 A morphism $\prA\rarrow\secA$ in $\Rings_\sgr$ is a morphism of graded
rings $f\:\prA\rarrow\secA$ whose \emph{degree-zero component
$f_0\:\prA_0\rarrow\secA_0$ is an isomorphism}.
 Then the same classes of $2$- and $3$\+left/right finitely projective
quadratic rings as in Propositions~\ref{2-fin-proj-quadratic-duality}
and~\ref{3-fin-proj-quadratic-duality} can be viewed as full
subcategories in $\Rings_\sgr$.
 The assertions of the two propositions remain valid with this
modification.
\end{rem}

 The inclusion of the full subcategory of quadratic graded rings
over $R$ into the category of (nonnegatively) graded rings over $R$
has a right adjoint functor, which we denote by $A\longmapsto \q A$.
 For any nonnegatively graded ring $A$, the quadratic graded
ring $A'=\q A$ together with the graded ring homomorphism
$A'\rarrow A$ is characterized by the properties that the maps
$A'_0\rarrow A_0$ and $A'_1\rarrow A_1$ are isomorphisms and
the map $A'_2\rarrow A_2$ is injective.
 Explicitly, the ring $\q A$ is constructed as the ring with
degree-one generators and quadratic relations $\q A=
T_{A_0}(A_1)/(I_A)$, where $I_A\subset A_1\ot_R A_1$ is the kernel of
the multiplication map $A_1\ot_RA_1\rarrow A_2$.

\Section{Flat and Finitely Projective Koszulity}
\label{koszulity-secn}

\subsection{Graded and ungraded Ext and Tor}
\label{graded-ungraded-subsecn}
 So far in this paper we denoted a graded ring by $A=
\bigoplus_{n=0}^\infty A_n$ (and for the most part we will continue
to do so in the sequel), but this is a colloquial abuse of notation.
 A graded abelian group $U$ is properly thought of as a collection of
abelian groups $U=(U_n)_{n\in\boZ}$.
 Then there are several ways to produce an ungraded group
from a graded one.

 Two of them are important for us in this section.
 One can take the direct sum of the grading components, which we
denote by $\Sigma U=\bigoplus_{n\in\boZ}U_n$; or one can take
the product of the grading components, which we denote by
$\Pi U=\prod_{n\in\boZ}U_n$.

 In particular, let $A=(A_n)_{n\in\boZ}$ be a graded ring and
$M=(M_n)_{n\in\boZ}$ be a graded left $A$\+module.
 Then $\Sigma A=\bigoplus_{n\in\boZ} A_n$ is the underlying ungraded
ring of $A$; and there are \emph{two} underlying ungraded
$\Sigma A$\+modules associated with~$M$.
 Namely, both the abelian groups $\Sigma M=\bigoplus_{n\in\boZ}M_n$
and $\Pi M=\prod_{n\in\boZ}M_n$ have natural structures of
left $\Sigma A$\+modules.
 Denoting the category of graded left $A$\+modules by $A\modl_\sgr$,
we have two forgetful functors $\Sigma$ and $\Pi\:A\modl_\sgr
\rarrow\Sigma A\modl$.

 The tensor product of a graded right $A$\+module $N$ and a graded
left $A$\+module $M$ is naturally a graded abelian group $N\ot_AM$,
and applying the functor $\Sigma$ to $N\ot_AM$ produces the tensor
product of the ungraded $\Sigma A$\+modules $\Sigma N$ and $\Sigma M$,
\begin{equation} \label{sigma-tensor}
 \Sigma(N\ot_AM)\simeq\Sigma N\ot_{\Sigma A}\Sigma M.
\end{equation}

 Similarly, for any graded left $A$\+modules $L$ and $M$ one can
consider the graded abelian group $\Hom_A(L,M)$ with the components
$\Hom_{A,n}(L,M)$ consisting of all the homogeneous left $A$\+module
maps $L\rarrow M$ of degree~$n$.
 The purpose of introducing the functor $\Pi$ above was to formulate
the comparison between the graded and ungraded Hom, which has
the form
\begin{equation} \label{sigma-pi-hom}
 \Pi\Hom_A(L,M)\simeq\Hom_{\Sigma A}(\Sigma L,\Pi M).
\end{equation}

 Furthermore, the functor $\Sigma$ takes projective graded $A$\+modules
to projective $\Sigma A$\+modules, while the functor $\Pi$ takes
injective graded $A$\+modules to injective $\Sigma A$\+modules
(as one can see from the description of projective and injective modules
as the direct summands of the free and cofree modules, respectively).
 In addition, the functor $\Sigma$ takes flat graded $A$\+modules to
flat $\Sigma A$\+modules (as one can see from the Govorov--Lazard
description of flat modules as the filtered direct limits of finitely
generated free modules).
 We define the graded versions of $\Tor$ and $\Ext$ as the derived
functors of the graded tensor product and $\Hom$, computed in
the abelian categories of graded (right and left) modules.

 So, for any graded right $A$\+module $N$ and any graded left
$A$\+module $M$ there is a \emph{bigraded} abelian group
$$
 \Tor^A(N,M)=(\Tor^A_{i,j}(N,M))_{i,j}, \qquad i\ge0, \ j\in\boZ,
$$
where $i$ is the usual \emph{homological} grading and $j$ is
the \emph{internal} grading (induced by the grading of $A$, $N$,
and~$M$).
 In order to compute the bigraded group $\Tor^A(M,N)$, one chooses
a \emph{graded} projective (or flat) resolution of one of
the $A$\+modules $M$ and $N$ and takes its tensor product over $A$ with
the other module; then the grading~$i$ is induced by the homological
grading of the resolution and the grading~$j$ comes from the grading
of the tensor product of any two graded modules.
 In view of the above considerations concerning projective/flat graded
modules, the formula~\eqref{sigma-tensor} implies a similar formula for
the Tor groups,
$$
 \Sigma\Tor^A_i(N,M)\simeq\Tor^{\Sigma A}_i(\Sigma N,\Sigma M)
 \quad \text{for every $i\ge0$},
$$
or more explicitly,
\begin{equation} \label{graded-and-ungraded-Tor}
 \Tor^{\Sigma A}_i(\Sigma M,\Sigma N)\,\simeq\,
 \bigoplus\nolimits_{j\in\boZ}\Tor^A_{i,j}(M,N).
\end{equation}

 Similarly, for any graded left $A$\+modules $L$ and $M$ there is
a \emph{bigraded} abelian group
$$
 \Ext_A(L,M)=(\Ext_{A,n}^i(L,M))_{i,n}, \qquad i\ge 0, \ n\in\boZ,
$$
where $i$ is the usual \emph{cohomological} grading and $n$ is
the \emph{internal} grading.
 In order to compute the bigraded group $\Ext_A(L,M)$, one chooses
either a graded projective resolution of the $A$\+module $L$, or
a graded injective resolution of the $A$\+module $M$, and takes
the graded Hom; then the grading~$i$ is induced by the (co)homological
grading of the resolution and the grading~$n$ comes from the grading
of the Hom groups.

 In the context of the internal grading of the Ext, we will put $n=-j$
and use the notation $\Ext_A^{i,j}(L,M)=\Ext_{A,-j}^i(L,M)$.
 By abuse of terminology, the grading $j$~will be also called
the \emph{internal grading} of the Ext.

 In view of the above considerations concerning projective and injective
graded modules, the formula~\eqref{sigma-pi-hom} implies a similar
formula for the Ext groups,
$$
 \Pi\Ext^i_A(L,M)\simeq\Ext^i_{\Sigma A}(\Sigma L,\Pi M)
 \quad\text{for every $i\ge0$},
$$
or more explicitly,
\begin{equation} \label{graded-and-ungraded-Ext}
 \Ext_{\Sigma A}^i(\Sigma L,\Pi M)\simeq
 \prod\nolimits_{j\in\boZ}\Ext_A^{i,j}(L,M).
\end{equation}

 For any three graded left $A$\+modules $K$, $L$, and $M$, there are
natural associative, unital composition/multiplication maps
\begin{equation} \label{graded-ext-multiplication}
 \Ext_A^{i',j'}(L,M)\times\Ext_A^{i',j'}(K,L)
 \lrarrow\Ext_A^{i'+i'',j'+j''}(K,M), \quad
 i',i''\ge0, \ j',j''\in\boZ
\end{equation}
on the bigraded $\Ext$ groups.
 Whenever the graded $A$\+module $L$ only has a finite number of nonzero
grading components (so $\Sigma L=\Pi L$), the passage to the infinite
products with respect to the internal gradings~$j'$ and~$j''$ makes
the multiplications~\eqref{graded-ext-multiplication} on the graded
$\Ext$ groups agree with the similar multiplications on the ungraded
$\Ext$ (between the $\Sigma A$\+modules $\Sigma K$, \ $\Sigma L=\Pi L$,
and $\Pi M$).

\subsection{Relative bar resolution} \label{relative-bar-subsecn}
 Given an $R$\+$R$\+bimodule $V$, we will use the notation
$V^{\ot_R\,n}=T_{R,n}(V)$ for the tensor product $V\ot_R\dotsb\ot_RV$
($n$~factors).

 Let $R\rarrow A$ be an injective homomorphism of associative rings.
 Denote by $A_+$ the $R$\+$R$\+bimodule $A/R$.
 Let $L$ be a left $A$\+module.
 The \emph{reduced relative bar resolution} of $L$ is the complex of
left $A$\+modules
\begin{equation} \label{bar-resolution}
 \dotsb\lrarrow A\ot_R A_+\ot_R A_+\ot_R L\lrarrow
 A\ot_R A_+\ot_R L\lrarrow A\ot_R L\lrarrow L\lrarrow0
\end{equation}
with the differential given by the standard formula
$\d(a_0\ot\bar a_1\ot\dotsb\ot \bar a_n\ot l)=
a_0a_1\ot a_2\ot\dotsb\ot a_n\ot l-a_0\ot a_1a_2\ot a_3\ot\dotsb
\ot l+\dotsb+(-1)^na_0\ot a_1\ot\dotsb\ot a_{n-1}\ot a_nl$.
 One can easily check that the image of the right-hand side in
$A\ot_RA_+^{\ot_R\,n-1}\ot_R L$ does not depend on the arbitrary
choice of liftings $a_i\in A$ of the given elements $\bar a_i\in A_+$,
\ $1\le i\le n$; so the differential is well-defined.

 The complex~\eqref{bar-resolution} is contractible as a complex of
left $R$\+modules; the contracting homotopy is given by the formulas
$t(l)=1\ot l$, \ $l\in L$, and
$t(a_0\ot \bar a_1\ot\dotsb\ot\bar a_n\ot l)=1\ot\bar a_0\ot
\bar a_1\ot\dotsb\ot\bar a_n\ot l$, where $\bar a_i\in A_+$, \
$1\le i\le n$, and $\bar a_0\in A_+$ is the image of the element
$a_0\in A$ under the natural surjection $A\rarrow A_+$.
 Hence it follows that the complex~\eqref{bar-resolution} is acyclic.

 If both $L$ and $A_+$ are flat left $R$\+modules, then all the left
$A$\+modules $A\ot_RA_+^{\ot_R\,n}\ot_RL$ are flat, so
\eqref{bar-resolution}~is a flat resolution of the left $A$\+module~$L$.
 Similarly, \eqref{bar-resolution}~is a projective resolution of
the left $A$\+module $L$ whenever both the left $R$\+modules $L$
and $A_+$ are projective.

 Assume that $L$ and $A_+$ are flat left $R$\+modules, and let $N$
be an arbitrary right $A$\+modules.
 Then one can use the flat resolution~\eqref{bar-resolution} of the left
$A$\+module $L$ in order to compute the groups $\Tor^A_i(N,L)$.
 Thus the groups $\Tor^A_i(N,L)$ are naturally isomorphic to
the homology groups of the bar complex
\begin{equation} \label{bar-complex}
\dotsb\lrarrow N\ot_R A_+\ot_R A_+\ot_R L\lrarrow
 N\ot_R A_+\ot_R L\lrarrow N\ot_R L\lrarrow0.
\end{equation}
 Switching the roles of the left and right modules and using
the reduced relative bar resolution of $N$, we conclude that
the same bar complex~\eqref{bar-complex} computes the groups
$\Tor^A_i(N,L)$ whenever a right $A$\+module $N$ is a flat right
$R$\+module, $A_+$ is a flat right $R$\+module, and $L$ is
an arbitrary left $A$\+module.

 Let $M$ be a left $A$\+module.
 The \emph{reduced relative cobar resolution} of $M$ is the complex of
left $A$\+modules
\begin{multline} \label{cobar-resolution}
 0\lrarrow M\lrarrow\Hom_R(A,M)\lrarrow\Hom_R(A_+\ot_R A,\>M) \\
 \lrarrow\Hom_R(A_+\ot_R A_+\ot_R A,\>M)\lrarrow\dotsb
\end{multline}
with the differential given by the formula
$(\d f)(\bar a_n\ot\dotsb\ot\bar a_1\ot a_0)=a_nf(a_{n-1}\ot\dotsb
\ot a_0)-f(a_na_{n-1}\ot a_{n-2}\ot\dotsb\ot a_0)+\dotsb+
(-1)^nf(a_n\ot\dotsb\ot a_2\ot a_1a_0)$, where $f\in
\Hom_R(A_+^{\ot_R\,n-1}\ot_R A,\>M)$ and $a_i\in A$ are arbitrary
liftings of elements $\bar a_i\in A_+$, \ $1\le i\le n$.
 One easily checks that the expression in the right-hand side
vanishes on the kernel of the natural surjection
$A^{\ot_R\,n+1}\rarrow A_+^{\ot_R\,n}\ot_R A$; so the differential
is well-defined.
 The left $A$\+module structure on $\Hom_R(A_+^{\ot_R\,n}\ot_R A,\>M)$
is induced by the right $A$\+module structure on~$A$.

 The complex~\eqref{cobar-resolution} is contractible as a complex
of left $R$\+modules; the contracting homotopy is given by
the formula $t(f)(\bar a_{n-1}\ot\dotsb\ot\bar a_1\ot a_0)=
(-1)^nf(\bar a_{n-1}\ot\dotsb\ot\bar a_1\ot\bar a_0\ot 1)$.
 In particular, it follows that the complex~\eqref{cobar-resolution}
is acyclic.
 If $M$ is an injective left $R$\+module and $A_+$ is a flat right
$R$\+module, then all the left $A$\+modules
$\Hom_R(A_+^{\ot_R\,n}\ot_R A,\>M)$, \ $n\ge 0$, are injective;
so~\eqref{cobar-resolution} is an injective resolution of
the left $A$\+module~$M$.

 Let $L$ and $M$ be left $A$\+modules.
 Assuming that $L$ is a projective left $R$\+module and $A_+$ is
a projective left $R$\+module, one can use the projective
resolution~\eqref{bar-resolution} of the left $A$\+module $L$ in
order to compute the groups $\Ext_A^i(L,M)$.
 Thus the groups $\Ext_A^i(L,M)$ are naturally isomorphic to
the cohomology groups of the cobar complex
\begin{multline} \label{cobar-complex}
 0\lrarrow\Hom_R(L,M)\lrarrow\Hom_R(A_+\ot_R L,\>M) \\
 \lrarrow\Hom_R(A_+\ot_R A_+\ot_R L,\>M)\lrarrow\dotsb
\end{multline}

 Alternatively, assuming that $M$ is an injective left $R$\+module
and $A_+$ is a flat right $R$\+module, one can use the injective
resolution~\eqref{cobar-resolution} of the left $A$\+module $M$ in
order to compute the groups $\Ext_A^i(L,M)$.
 Under these assumptions, one comes to the same conclusion that
the groups $\Ext_A^i(L,M)$ are naturally isomorphic to
the cohomology groups of the cobar complex~\eqref{cobar-complex}.

 When $R$ and $A$ are graded rings, $R\rarrow A$ is a graded ring
homomorphism, and $L$ is a graded left $A$\+module, one can
interpret~\eqref{bar-resolution} as a graded resolution of
the graded $A$\+module~$L$.
 When $L$ and $A_+$ are flat graded left $R$\+modules,
\eqref{bar-resolution}~is a graded flat resolution of the graded
$A$\+module $L$; and when $L$ and $A_+$ are projective graded left
$R$\+modules, \eqref{bar-resolution} is a graded projective resolution.

 For a graded left $A$\+module $M$, one can also interpret
the $\Hom$ notation in~\eqref{cobar-resolution} as the graded $\Hom$,
obtaning a graded resolution of the graded $A$\+module~$M$.
 When $M$ is an injective graded left $R$\+module and $A_+$ is
a flat graded right $R$\+module, \eqref{cobar-resolution}~is
a graded injective resolution.

 It follows that, under the graded versions of the above flatness,
projectivity, and/or injectivity assumptions, the graded
bar complex~\eqref{bar-complex} computes the bigraded
$\Tor^A(N,L)$, and the graded version of the cobar
complex~\eqref{cobar-complex} computes the bigraded $\Ext_A(L,M)$.

 The functor $\Sigma$ transforms the graded versions
of~\eqref{bar-resolution} and~\eqref{bar-complex} (for graded rings
$R$ and $A$ and graded modules $L$ and~$N$) into the ungraded ones
(for the ungraded rings $\Sigma R$ and $\Sigma A$ and the ungraded
modules $\Sigma L$ and~$\Sigma N$).
 The functor $\Pi$ transforms the graded versions 
of~\eqref{cobar-resolution} and~\eqref{cobar-complex} (for graded rings
$R$ and $A$ and graded modules $L$ and~$M$) into the ungraded ones
(for the ungraded rings $\Sigma R$ and $\Sigma A$ and the ungraded
modules $\Sigma L$ and~$\Pi M$).

 The cobar complexes also compute the composition/multiplication on
the $\Ext$ groups.
 Let $K$, $L$, and $M$ be left $A$\+modules; assume that the left
$R$\+modules $K$, $L$, and $A_+$ are projective.
 Then the natural composition/multiplication maps on the cobar complexes
$$
 \Hom_R(A_+^{\ot_R\,i'}\ot_R L,\>M)\times
 \Hom_R(A_+^{\ot_R\,i''}\ot_R K,\>L)\lrarrow
 \Hom_R(A_+^{\ot_R\,i'+i''}\ot_R K,\>M)
$$
agree with the cobar differentials and induce the Yoneda multiplication
maps $\Ext_A^{i'}(L,M)\times\Ext_A^{i''}(K,L)\rarrow
\Ext_A^{i'+i''}(K,M)$ on the $\Ext$ groups.
 In the case of graded rings $R$ and $A$ and graded $A$\+modules $K$,
$L$, and $M$, the same assertion applies to the graded $\Ext$.

\subsection{Diagonal Tor and Ext} \label{diagonal-Tor-Ext-subsecn}
 Let $A=\bigoplus_{n=0}^\infty A_n$ be a nonnegatively graded ring
with the degree-zero component $R=A_0$.
 Then the projection onto the degree-zero component is a graded ring
homomorphism $A\rarrow R$.
 Using this homomorphism, one can consider $R$ as a left and right
graded module over~$A$.

 So we can consider the bigraded $\Tor$ groups $\Tor^A_{i,j}(R,R)$ and
the bigraded $\Ext$ ring $\Ext_A^{i,j}(R,R)$, with the (co)homological
grading~$i$ and the internal grading~$j$.
 In fact, as $R$ is an $R$\+$A$\+bimodule and an $A$\+$R$\+bimodule,
the groups $\Tor^A_{i,j}(R,R)$ have natural structures of
$R$\+$R$\+bimodules.

 First of all, we notice the connection between the graded and ungraded
$\Tor$ and $\Ext$.
 As the graded (left or right) $A$\+module $R$ is concentrated in
the internal grading~$0$, we have $\Sigma R=R=\Pi R$.
 Hence the formulas~\eqref{graded-and-ungraded-Tor}
and~\eqref{graded-and-ungraded-Ext} reduce to
\begin{align*}
 \Tor^{\Sigma A}_i(R,R)&\simeq
 \bigoplus\nolimits_{j\in\boZ}\Tor^A_{i,j}(R,R) \\
 \Ext_{\Sigma A}^i(R,R)&\simeq
 \prod\nolimits_{j\in\boZ}\Ext_A^{i,j}(R,R).
\end{align*}

\begin{prop} \label{diagonal-Tor}
 Assume that the grading components $A_n$, \,$n\ge1$, are flat left
$R$\+modules, and the multiplication map $A_1\ot_R A_1\rarrow A_2$
is surjective with the kernel~$I_A$.
 Then there are natural isomorphisms of $R$\+$R$\+bimodules \par
\textup{(a)} $\Tor_{i,j}^A(R,R)=0$ whenever $i<0$,
or $i=0$ and $j>0$, or $i>j$; \par
\textup{(b)} $\Tor_{0,0}^A(R,R)=R$, \ $\Tor_{1,1}^A(R,R)\simeq A_1$, \
$\Tor_{2,2}^A(R,R)\simeq I_A$, and
$$
 \Tor_{n,n}^A(R,R)\simeq \bigcap\nolimits_{k=1}^{n-1}
 A_1^{\ot_R\,k-1}\ot_R I_A\ot_R A_1^{\ot_R\,n-k-1}
 \subset A_1^{\ot_R\,n}, \qquad n\ge 2.
$$
\end{prop}

\begin{proof}
 The $R$\+$R$\+bimodules $\Tor^A_{i,j}(R,R)$ can be computed
as the homology bimodules of the bar complex~\eqref{bar-complex}
for $N=R=L$,
\begin{equation} \label{N=R=L}
 \dotsb\lrarrow A_+\ot_R A_+\ot_R A_+\lrarrow A_+\ot_RA_+
 \lrarrow A_+\overset0\lrarrow R\lrarrow0.
\end{equation}
 More explicitly, this means that, for every fixed $n\ge1$,
the bimodules $\Tor^A_{i,n}(R,R)$ are the homology bimodules of
the complex of $R$\+$R$\+bimodules
\begin{multline} \label{Tor-R-R-bar-complex}
 0\lrarrow A_1^{\ot_R\,n}\lrarrow\bigoplus\nolimits_{k=1}^{n-1}
 A_1^{\ot_R\,k-1}\ot_RA_2\ot_R A_1^{\ot_R\,n-k-1}\lrarrow \\ \dotsb
 \lrarrow (A_+^{\ot_R\,i})_n\lrarrow\dotsb \\ \lrarrow
 \bigoplus\nolimits_{k=1}^{n-1}A_k\ot_R A_{n-k}\lrarrow A_n\lrarrow0.
\end{multline}
 The assertion~(a) follows immediately; and to prove~(b), it remains
to compute, for every $1\le k\le n-1$, the kernel of the map
$A_1^{\ot_R\,n}\rarrow A_1^{\ot_R\,k-1}\ot_RA_2\ot_RA_1^{\ot_R\,n-k-1}$
induced by the multiplication map $A_1\ot_RA_1\rarrow A_2$.

 Specifically, we have to check that the natural short sequence
of $R$\+$R$\+bimodules
$$
 0\rarrow A_1^{\ot_R\,k-1}\ot_RI_A\ot_R A_1^{\ot_R\,n-k-1}
 \rarrow A_1^{\ot_R\,n}\rarrow
 A_1^{\ot_R\,k-1}\ot_RA_2\ot_RA_1^{\ot_R\,n-k-1}\rarrow 0
$$
is exact.
 Indeed, by assumption we have a short exact sequence of
$R$\+$R$\+bimodules
\begin{equation} \label{degree-2-sequence}
 0\lrarrow I_A\lrarrow A_1\ot_R A_1\lrarrow A_2\lrarrow0,
\end{equation}
whose terms are flat left $R$\+modules.
 Furthermore, $A_1^{\ot_R\,n-k-1}$ is a flat left $R$\+module, too.
 It follows that tensoring~\eqref{degree-2-sequence} 
(with any right $R$\+module, and in particular) with $A_1^{\ot_R\,k-1}$
on the left does not affect exactness; and neither does tensoring
(any exact sequence of right $R$\+modules, and in particular)
the resulting short sequence with $A_1^{\ot_R\,n-k-1}$ on the right.
\end{proof}

 Concerning the $\Ext$, our notation $\Ext_A(R,R)$ presumes, as above,
that the $\Ext$ is taken in the category of \emph{left} $R$\+modules.
 So, in particular, the ring $\Ext^0_A(R,R)=\Hom_R(R,R)=R^\rop$ is
the opposite ring to~$R$.

\begin{prop} \label{diagonal-Ext}
 Assume that the grading components $A_n$, \,$n\ge1$, are finitely
generated projective left $R$\+modules, and the multiplication map
$A_1\ot_R A_1\rarrow A_2$ is surjective with the kernel~$I_A$.
 Then \par
\textup{(a)} $\Ext_A^{i,j}(R,R)=0$ whenever $i<0$,
or $i=0$ and $j>0$, or $i>j$; \par
\textup{(b)} the diagonal Ext ring\/ $\bigoplus_{n=0}^\infty
\Ext_A^{n,n}(R,R)$ is naturally isomorphic, as a graded ring,
to the opposite ring $B^\rop$ to the $2$\+right finitely projective
quadratic dual ring $B=T_R(B_1)/(I_B)$, \ $B_1=\Hom_R(A_1,R)$, \
$I_B=\Hom_R(A_2,R)$ to the $2$\+left finitely projective quadratic
ring\/ $\q A=T_R(A_1)/(I_A)$.
\end{prop}

\begin{proof}
 Part~(a) does not depend on the finite generatedness assumptions
on the $R$\+modules $A_n$, and requires only the projectivity
assumptions.
 The bigraded ring $\Ext_A(R,R)$ can be computed as the cohomology
ring of the DG\+ring~\eqref{cobar-complex} for $L=R=M$.
 The latter can be obtained by applying the functor $\Hom_R({-},R)$
to the bar complex~\eqref{N=R=L},
\begin{equation} \label{L=R=M}
 0\lrarrow R\overset0\lrarrow\Hom_R(A_+,R)\lrarrow
 \Hom_R(A_+\ot_RA_+,\>R)\lrarrow\dotsb
\end{equation}
 More specifically, for every fixed $n\ge1$, the groups
$\Ext^{i,n}_A(R,R)$ are the cohomology groups of the complex
obtained by applying $\Hom_R({-},R)$ to
the complex~\eqref{Tor-R-R-bar-complex}.
 This proves part~(a).

 When all the grading components $A_j$ are finitely generated and
projective left $R$\+modules, the complex obtained by applying
the functor $\Hom_R({-},R)$ to the complex~\eqref{Tor-R-R-bar-complex}
can be computed using Lemma~\ref{tensor-dual-lemma}(a) as
\begin{multline} \label{Ext-R-R-cobar-complex}
 0\lrarrow A_n\spcheck\lrarrow
 \bigoplus\nolimits_{k=1}^{n-1}A_{n-k}\spcheck\ot_RA_k\spcheck\lrarrow
 \\ \dotsb\lrarrow(A_+\spcheck{}^{\ot_R\,i})^n\lrarrow\dotsb \\
 \lrarrow\bigoplus\nolimits_{k=1}^{n-1}A_1\spcheck{}^{\ot_R\,n-k-1}
 \ot_RA_2\spcheck\ot_RA_1\spcheck{}^{\ot_R\,k-1}\lrarrow
 A_1\spcheck{}^{\ot_R\,n}\lrarrow0,
\end{multline}
where the notation is $U\spcheck=\Hom_R(U,R)$.
 Setting $B_1=A_1\spcheck$, \ $I_B=A_2\spcheck$, and
$B=T_R(B_1)/(I_B)$, we obtain the desired isomorphism
$\Ext_A^{n,n}(R,R)\simeq B_n$.
 It follows immediately from the construction of the multiplication on
the cobar complex that this is a graded ring isomorphism between
$\bigoplus_n\Ext_A^{n,n}(R,R)$ and $B^\rop$.
\end{proof}

\begin{rem} \label{bar-complex-dg-coring}
 The bar complex $\Br_R(A)$ computing $\Tor^A(R,R)$
\,\eqref{N=R=L} has a natural structure of a (coassociative,
counital) \emph{DG\+coring} over the ring $R$, with the obvious maps
of counit $\Br_R(A)\rarrow R$ and a comultiplication $\Br_R(A)\rarrow
\Br_R(A)\ot_R\Br_R(A)$ compatible with the bar differential.
 However, as the functor of tensor product over $R$ is not left
exact, this DG\+coring structure, generally speaking, does \emph{not}
descend to a coring structure on the homology modules.
 When $A_+$ is a flat left $R$\+module \emph{and} $\Tor^A(R,R)$ is
a flat left $R$\+module, the nonexactness problem does not interfere,
and the (bi)graded $R$\+$R$\+bimodule $\Tor^A(R,R)$ is
a (coassociative, counital) coring over~$R$.
 In particular, the diagonal $\Tor$ bimodule
$C=\bigoplus_{n=0}^\infty\Tor^A_{n,n}(R)$ becomes a graded coring
over $R$ in these assumptions, with the obvious counit $C\rarrow R$
and the induced comultiplication map $C\rarrow C\ot_RC$.
\end{rem}

\subsection{Low-dimensional Tor, degree-one generators and
quadratic relations}
 As in the previous section, we consider a nonnegatively graded ring
$A=\bigoplus_{n=0}^\infty A_n$ with the degree-zero component $R=A_0$.

\begin{prop} \label{low-dimensional-Tor-prop}
 Assume that the grading components $A_n$, \,$n\ge1$, are flat left
$R$\+modules.  Then \par
\textup{(a)} the graded ring $A$ is generated by $A_1$ over $R$ if
and only if\/ $\Tor^A_{1,j}(R,R)=0$ for all $j>1$; \par
\textup{(b)} the graded ring $A$ is quadratic if and only if\/
$\Tor^A_{1,j}(R,R)=0$ for all $j>1$ and\/
$\Tor^A_{2,j}(R,R)=0$ for all $j>2$.
\end{prop}

\begin{proof}
 Part~(a): folowing the proof of Proposition~\ref{diagonal-Tor},
the group $\Tor^A_{1,n}(R,R)$ can be computed as the rightmost homology
group of the complex~\eqref{Tor-R-R-bar-complex}, that is, the cokernel
of the multiplication map $\bigoplus_{k=1}^n A_k\ot_RA_{n-k}
\rarrow A_n$.
 Now it is clear that a nonnegatively graded ring $A$ is generated by
$A_1$ over $R=A_0$ if and only if the latter map is surjective for
all $n\ge 2$.

 Part~(b): by part~(a), we can assume that $A$ is generated by $A_1$;
then we have to show that $A$ is quadratic if and only if
$\Tor^A_{2,j}(R,R)=0$ for all $j>2$.
 Once again, we compute the group $\Tor^A_{2,n}(R,R)$ using
the complex~\eqref{Tor-R-R-bar-complex}.
 So we have to show that $A$ is quadratic if and only if the short
sequence
\begin{equation} \label{low-dimensional-Tor-R-R-bar-complex}
 \bigoplus\nolimits_{k+l+m=n}^{k,l,m\ge1}A_k\ot_RA_l\ot_RA_m
 \lrarrow\bigoplus\nolimits_{k+l=n}^{k,l\ge1}A_k\ot_RA_l
 \lrarrow A_n\lrarrow0
\end{equation}
is right exact for all $n\ge3$.

 Given an $R$\+$R$\+bimodule $V$ and an $R$\+$R$\+subbimodule
$I\subset V\ot_RV$, one can construct the quadratic ring
$A'=T_R(V)/(I)$ by the following inductive procedure.
 Set $A'_0=R$, \ $A'_1=V$, and $A'_2=(V\ot_RV)/I$; then there are
the obvious multplication maps $A'_k\times A'_l\rarrow A'_{k+l}$
for $k$, $l\ge0$, \ $k+l\le2$.
 For every $n\ge3$, set $A'_n$ to be the cokernel of
the $R$\+$R$\+bimodule morphism
$$
 \bigoplus\nolimits_{k+l+m=n}^{k,l,m\ge1}A'_k\ot_RA'_l\ot_RA'_m
 \lrarrow\bigoplus\nolimits_{k+l=n}^{k,l\ge1}A'_k\ot_RA'_l.
$$
 Then we have biadditive multiplication maps $A'_k\times A'_l
\rarrow A'_n$ defined for all $k+l=n$, \,$k$, $l\ge0$ and
satisfying the associativity equations $(ab)c=a(bc)$ for all
$a\in A'_k$, \ $b=A'_l$, \ $c\in A'_m$, \ $k+l+m=n$,
\,$k$, $l$, $m\ge0$.
 So we obtain a graded ring~$A'$.
 Obviously, $A'$ is the graded associative ring freely generated
by $A'_1=V$ over $A'_0=R$ with the relations $I\subset V\ot_R V$;
so in other words, $A'\simeq T_R(V)/(I)$.

 Returning to the original graded ring $A$, put $V=A_1$ and let $I$ be
the kernel of the multiplication map $A_1\ot_RA_1\rarrow A_2$.
 Then there is a unique homomorphism of graded rings $A'\rarrow A$
acting by the identity map in degrees~$0$ and~$1$.
 The ring homomorphism $A'\rarrow A$ is surjective by assumption.
 Arguing by induction in~$n$, one can easily see that this ring
homomorphism is an isomorphism if and only if the short
sequence~\eqref{low-dimensional-Tor-R-R-bar-complex} is right
exact for every $n\ge3$.
\end{proof}

\subsection{First Koszul complex} \label{first-koszul-complex-subsecn}
 Let $A=\bigoplus_{n=0}^\infty A_n$ and $B=\bigoplus_{n=0}^\infty B_n$
be two nonnegatively graded rings with the same degree-zero component
$A_0=R=B_0$.
 Suppose that $B_1$ is a finitely generated projective right
$R$\+module, and that we are given an $R$\+$R$\+bimodule morphism
$\tau\:\Hom_{R^\rop}(B_1,R)\rarrow A_1$.

 Equivalently, instead of a map~$\tau$ one can consider an element
of the tensor product $B_1\ot_RA_1$ satisfying the equation spelled out
in the next lemma (which is to be applied to the $R$\+$R$\+bimodules
$V=\Hom_{R^\rop}(B_1,R)$ and $U=A_1$).

\begin{lem} \label{element-e-lemma}
 Let $R$ and $S$ be assocative rings. \par
\textup{(a)} Let $U$ be a left $R$\+module and $V$ be a finitely
generated projective left $R$\+module.
 Then there is a natural isomorphism of abelian groups\/
$\Hom_R(V,U)\simeq\Hom_R(V,R)\ot_R\nobreak U$. \par
\textup{(b)} Let $V$ and $U$ be two $R$\+$S$\+bimodules such that
$V$ is a finitely generated projective left $R$\+module.
 Then $R$\+$S$\+bimodule morphisms $f\:V\rarrow U$ correspond to
elements $e\in\Hom_R(V,R)\ot_RU$ satisfying the equation
$se=es$ for all $s\in S$ under the isomorphism of part~(a).
 In other words, there is a natural bijective correspondence
between $R$\+$S$\+bimodule morphisms $f\:V\rarrow U$ and
$S$\+$S$\+bimodule morphisms $e\:S\rarrow\Hom_R(V,R)\ot_RU$.
\end{lem}

\begin{proof}
 Both the assertions are easy.
 Let us only specify that the correspondence between the maps
$f\:V\rarrow U$ and the elements $e=e(1)\in\Hom_R(V,R)\ot_RU$ is
given by the rule $f(v)=\lan v,e\ran$.
 Here the pairing notation
(cf.\ Section~\ref{cdg-ring-constructed-subsecn} below)
stands for the map
$$
 \lan\ , \ \ran\:V\times\Hom_R(V,R)\lrarrow R,
 \qquad \lan v,b\ran = b(v),
$$
and, by extension, for the map
$$
 \lan\ ,\ \ran\:V\times\Hom_R(V,R)\ot_R U\lrarrow U
$$
given by the formulas $\lan v,\,b\ot u\ran = \lan v,b\ran u = b(v)u$,
where $v\in V$, \ $u\in U$, and $b\in\Hom_R(V,R)$.
\end{proof}

 Denote by $K^\tau(B,A)$ the following bigraded $A$\+$B$\+bimodule
endowed with an $A$\+$B$\+bimodule endomorphism~$\d^\tau$.
 The bigrading components of $K^\tau(B,A)$ are
$$
 K^\tau_{i,n}(B,A)=\Hom_{R^\rop}(B_i,A_{n-i}),
 \quad i\ge 0, \ n\ge i.
$$
 Here $i$~is interpreted as the homological grading and $n$~is
the internal grading.
 The homogeneous $A$\+$B$\+bimodule endomorphism $\d^\tau\:
K^\tau(B,A)\rarrow K^\tau(B,A)$ of bidegree $(i,n)=(-1,0)$ is
constructed as the composition
\begin{multline*}
 \Hom_{R^\rop}(B_i,A_j)\lrarrow\Hom_{R^\rop}(B_{i-1}\ot_R B_1,\>A_j)
 \\ \,\simeq\,\Hom_{R^\rop}(B_{i-1},\Hom_{R^\rop}(B_1,A_j))\,\simeq\,
 \Hom_{R^\rop}(B_{i-1},\>A_j\ot_R\Hom_{R^\rop}(B_1,R)) \\
 \overset\tau\lrarrow \Hom_{R^\rop}(B_{i-1},\>A_j\ot_R A_1)\lrarrow
 \Hom_{R^\rop}(B_{i-1},A_{j+1}).
\end{multline*}
 Here the map $\Hom_{R^\rop}(B_i,A_j)\rarrow
\Hom_{R^\rop}(B_{i-1}\ot_RB_1,\>A_j)$ is induced by the multiplication
map $B_{i-1}\ot_RB_1\rarrow B_i$, the map
$\Hom_{R^\rop}(B_{i-1},\>A_j\ot_R\Hom_{R^\rop}(B_1,R))\rarrow
\Hom_{R^\rop}(B_{i-1},\>A_j\ot_RA_1)$ is induced by the given map
$\tau\:\Hom_{R^\rop}(B_1,R)\rarrow A_1$, and the map
$\Hom_{R^\rop}(B_{i-1},\>A_j\ot_RA_1)\rarrow
\Hom_{R^\rop}(B_{i-1},A_{j+1})$ is induced by the multiplication map
$A_j\ot_RA_1\rarrow A_{j+1}$.

\begin{lem} \label{Koszul-differential-square-zero}
 Assume that both $B_1$ and $B_2$ are finitely generated projective
right $R$\+modules.
 Then the following conditions are equivalent: \par
\textup{(a)} $(\d^\tau)^2=0$ on the whole bigraded $A$\+$B$\+bimodule
$K^\tau(B,A)$; \par
\textup{(b)} the composition\/ $\Hom_{R^\rop}(B_2,R)\overset{\d^\tau}
\rarrow\Hom_{R^\rop}(B_1,A_1)\overset{\d^\tau}\rarrow
\Hom_{R^\rop}(R,A_2)=A_2$ vanishes; \par
\textup{(c)} the composition of maps\/ $\Hom_{R^\rop}(B_2,R)\rarrow
\Hom_{R^\rop}(B_1\ot_RB_1,\>R)\simeq\Hom_{R^\rop}(B_1,R)
\ot_R\Hom_{R^\rop}(B_1,R)\overset{\tau\ot\tau}\lrarrow A_1\ot_R A_1
\rarrow A_2$ vanishes.
\end{lem}

\begin{proof}
 (a)\,$\Longrightarrow$\,(b) is obvious.
 
 (b)\,$\Longleftrightarrow$\,(c) holds because the two maps in question
are the same.

 (c)\,$\Longrightarrow$\,(a) is straightforward, using the isomorphism
$\Hom_{R^\rop}(B_2,A_j)\simeq A_j\ot_R\Hom_{R^\rop}(B_2,R)$.
\end{proof}

 When the equivalent conditions of
Lemma~\ref{Koszul-differential-square-zero} hold, the bigraded
$A$\+$B$\+bimodule $K^\tau(B,A)$ with the differential~$\d^\tau$ can be
viewed as a complex of graded left $A$\+modules.
 We call it the \emph{first Koszul complex} and denote by
$K^\tau_\bu(B,A)$.
 Assuming that $B_i$ is a finitely generated projective right
$R$\+module for every $i\ge0$, the first Koszul complex
$K^\tau_\bu(B,A)$ has the form
\begin{equation} \label{first-koszul-complex}
 \dotsb\lrarrow A\ot_R\Hom_{R^\rop}(B_2,R)\lrarrow
 A\ot_R\Hom_{R^\rop}(B_1,R)\lrarrow A\lrarrow0.
\end{equation}

 Specifically, the following particular case is important.
 Let $A$ be a nonnegatively graded ring with the degree-zero component
$R=A_0$.
 Assume that $A_1$ and $A_2$ are finitely generated projective left
$R$\+modules and the multiplication map $A_1\ot_RA_1\rarrow A_2$
is surjective.
 Denote the kernel of the latter map by $I_A$, and consider
the $2$\+left finitely projective quadratic graded ring
$\q A=T_R(A_1)/(I_A)$.
 Let $B$ be the $2$\+right finitely projective quadratic graded ring
quadratic dual to~$\q A$.

 Then we have an isomorphism of $R$\+$R$\+bimodules
$\Hom_{R^\rop}(B_1,R)\simeq A_1$, which we use as our choice of
the map~$\tau$.
 The assumptions of Lemma~\ref{Koszul-differential-square-zero} are
satisfied, and condition~(c) holds by the construction of quadratic
duality.
 Therefore, the first Koszul complex $K^\tau_\bu(B,A)$ is well-defined.

\subsection{Dual Koszul complex} \label{dual-koszul-complex-subsecn}
 As in Section~\ref{first-koszul-complex-subsecn}, we consider two
nonnegatively graded rings $A=\bigoplus_{n=0}^\infty A_n$ and
$B=\bigoplus_{n=0}^\infty B_n$ with the same degree-zero component
$A_0=R=B_0$.
 Let $e\in B_1\ot_RA_1$ be an element satisfying the equation
$re=er$ for all $r\in R$, as in Lemma~\ref{element-e-lemma}(b).

 Denote by $\Ksp_e(B,A)$ the following bigraded $B$\+$A$\+bimodule
endowed with a $B$\+$A$\+bimodule endomorphism~$d_e$.
 The bigrading components of $\Ksp_e(B,A)$ are
$$
 \Ksp_e{}^{i,n}(B,A)=B_i\ot A_{n+i},
 \quad i\ge 0, \ n\ge -i.
$$
 Here $i$~is intepreted as the cohomological grading and $n$~is
the internal grading.
 The homogeneous $B$\+$A$\+bimodule endomorphism $d_e\:\Ksp_e(B,A)
\rarrow\Ksp_e(B,A)$ of bidegree $(i,n)=(1,0)$ consists of the components
$B_i\ot_RA_j\rarrow B_{i+1}\ot_RA_{j+1}$, which are constructed
as the compositions
$$
 B_i\ot_R A_j=B_i\ot_RR\ot_RA_j\overset e\lrarrow
 B_i\ot_R B_1\ot_R A_1\ot_R A_j\lrarrow B_{i+1}\ot_R A_{j+1}.
$$
 Here the map $B_i\ot_RR\ot_R A_j\rarrow B_i\ot_RB_1\ot_RA_1\ot_RA_j$
is induced by the $R$\+$R$\+bimodule map $e\:R\rarrow B_1\ot_R A_1$,
while the map $B_i\ot_RB_1\ot_RA_1\ot_RA_j\rarrow B_{i+1}\ot_RA_{j+1}$
is the tensor product of the multiplication maps
$B_i\ot_RB_1\rarrow B_{i+1}$ and $A_1\ot_RA_j\rarrow A_{j+1}$.

\begin{lem} \label{dual-Koszul-differential}
 Assume that both $B_1$ and $B_2$ are finitely generated projective
right $R$\+modules, and that an $R$\+$R$\+bimodule map\/
$\Hom_{R^\rop}(B_1,R)\rarrow A_1$ and an element $e\in B_1\ot_RA_1$
correspond to each other under the construction of
Lemma~\ref{element-e-lemma}.
 Then the following conditions are equivalent: \par
\textup{(a)} $d_e^2=0$ on the whole bigraded $B$\+$A$\+bimodule
$\Ksp_e(B,A)$; \par
\textup{(b)} the element $d_e^2(1\ot1)\in B_2\ot_RA_2$ vanishes; \par
\textup{(c)} the equivalent conditions of
Lemma~\ref{Koszul-differential-square-zero} hold. \par
 Furthermore, the bigraded $A$\+$B$\+bimodule $K^\tau(B,A)$ with
the endomorphism~$\d^\tau$ can be obtained by applying the functor
$\Hom_{A^\rop}({-},A)$ to the bigraded $B$\+$A$\+bimodule
$\Ksp_e(B,A)$ with the endomorphism~$d_e$.
\end{lem}

\begin{proof}
 The equivalence (a)\,$\Longleftrightarrow$\,(b) does not depend on
the assumptions of finite generatedness and projectivity of
the right $R$\+modules $B_1$ and~$B_2$.
 The equivalence (b)\,$\Longleftrightarrow$\,(c) appears to need
these assumptions.
 The last assertion of the lemma is straightforward, and only needs
the assumption of $B_1$ being a finitely generated projective
right $R$\+module.

 (a)\,$\Longrightarrow$\,(b) is obvious.
 
 (b)\,$\Longrightarrow$\,(a) holds because $d_e$~is
a $B$\+$A$\+bimodule map.

 (b)\,$\Longleftrightarrow$\,(c) holds because the map
$\Hom_{R^\rop}(B_2,R)\rarrow A_2$ from
Lemma~\ref{Koszul-differential-square-zero}(c) corresponds to
the element $d_e^2(1\ot1)\in B_2\ot_RA_2$ under the correspondence
of Lemma~\ref{element-e-lemma}.

 (a)\,$\Longrightarrow$\,(c): the condition~(a) implies
the condition of Lemma~\ref{Koszul-differential-square-zero}(a)
in view of the last assertion of the present lemma.
\end{proof}

 When the equivalent conditions of the lemma hold, the bigraded
$B$\+$A$\+bimodule $\Ksp_e(B,A)$ with the differential~$d_e$ can be
viewed as a complex of graded right $A$\+modules.
 We call it the \emph{dual Koszul complex} and denote by
$\Ksp_e{}^\bu(B,A)$.
 The dual Koszul complex has the form
\begin{equation} \label{dual-koszul-complex}
 0\lrarrow A\lrarrow B_1\ot_RA\lrarrow B_2\ot_RA\lrarrow
 B_3\ot_RA\lrarrow\dotsb.
\end{equation}
 According to the lemma, we have
$$
 K^\tau_\bu(B,A)=\Hom_{A^\rop}(\Ksp_e{}^\bu(B,A),A).
$$

 The graded Hom group $\Hom_{R^\rop}(B,R)$ with the components
$\Hom_{R^\rop}(B_i,R)$, \ $i\ge0$, has natural structures of a graded
right $B$\+module (induced by the graded left $B$\+module structure
on~$B$) and a graded left $R$\+module (induced by the left $R$\+module
structure on~$R$).
 We denote the tensor product $\Hom_{R^\rop}(B,R)\ot_B\Ksp_e(B,A)$ by
$$
 {}^\tau\!K(B,A)=\Hom_{R^\rop}(B,R)\ot_B\Ksp_e(B,A)=
 \Hom_{R^\rop}(B,R)\ot_RA.
$$
 ${}^\tau\!K(B,A)$ is a bigraded $R$\+$A$\+bimodule with the bigrading
components
$$
 {}^\tau\!K_{i,n}(B,A)=\Hom_{R^\rop}(B_i,R)\ot_RA_{n-i},
 \quad i\ge0, \ n\ge i.
$$
 Here $i$~is the homological grading and $n$~is the internal grading.
 Notice that there is \emph{no} $B$\+module structure on
${}^\tau\!K(B,A)$.

 The $B$\+$A$\+bimodule endomorphism $d_e\:\Ksp_e(B,A)\rarrow
\Ksp_e(B,A)$ induces an $R$\+$A$\+bimodule endomorphism of bidegree
$(i,n)=(-1,0)$ on ${}^\tau\!K(B,A)$, which we denote by
${}^\tau\!\d\:{}^\tau\!K(B,A)\rarrow{}^\tau\!K(B,A)$.
 When the equivalent conditions of Lemma~\ref{dual-Koszul-differential}
hold, it follows that $({}^\tau\!\d)^2=0$.
 So the bigraded $R$\+$A$\+bimodule ${}^\tau\!K(B,A)$ can be viewed
as a complex of graded right $A$\+modules.
 We call it the \emph{second Koszul complex} and denote by
${}^\tau\!K_\bu(B,A)$.
 The second Koszul complex ${}^\tau\!K_\bu(B,A)$ has the form
\begin{equation} \label{second-koszul-complex}
 \dotsb\lrarrow\Hom_{R^\rop}(B_2,R)\ot_RA\lrarrow
 \Hom_{R^\rop}(B_1,R)\ot_RA\lrarrow A\lrarrow0.
\end{equation}

 The particular case described at the end of
Section~\ref{first-koszul-complex-subsecn} is important for us.
 In this case, $\q A$ is a $2$\+left finitely projective quadratic
graded ring and $B$ is the quadratic dual $2$\+right finitely
projective quadratic graded ring.
 Choosing~$\tau$ to be the natural isomorphism $\Hom_{R^\rop}(B_1,R)
\simeq A_1$, we see that the equivalent conditions of
Lemma~\ref{dual-Koszul-differential} are satisfied (because
the condition of Lemma~\ref{Koszul-differential-square-zero}(c) holds).
 Therefore, the dual Koszul complex $\Ksp_e{}^\bu(B,A)$ and
the second Koszul complex ${}^\tau\!K_\bu(B,A)$ are well-defined.

\begin{rem} \label{B-sharp-remark}
 Assume that $B_n$ is a finitely generated projective right
$R$\+module for every $n\ge0$.
 Then the isomorphism of Lemma~\ref{tensor-dual-lemma}(b) provides
the graded $R$\+$R$\+bimodule $C$ with the components
$C_n=\Hom_{R^\rop}(B_n,R)$ with the structure of a graded coring
over~$R$.
 The comultiplication maps $C_{i+j}\rarrow C_i\ot_RC_j$ are obtained
by applying the functor $\Hom_{R^\rop}({-},R)$ to the multiplication
maps $B_j\ot B_i\rarrow B_{i+j}$.

 By construction, we have $B_n=\Hom_R(C_n,R)$.
 Put $\shB_n=\Hom_{R^\rop}(C_n,R)$ (in the spirit of the notation
in~\cite{Ar}).
 Then applying the functor $\Hom_{R^\rop}({-},R)$ to
the comultiplication maps $C_{i+j}\rarrow C_i\ot_RC_j$ and
composing the resulting maps $\Hom_{R^\rop}(C_i\ot_RC_j,\>R)\rarrow
\Hom_{R^\rop}(C_{i+j},R)$ with the natural $R$\+$R$\+bimodule
morphisms $\Hom_{R^\rop}(C_j,R)\ot_R\Hom_{R^\rop}(C_i,R)\rarrow
\Hom_{R^\rop}(C_{i+j},R)$ from Lemma~\ref{tensor-dual-lemma}(b)
produces multiplication maps $\shB_j\ot_R \shB_i\rarrow \shB_{i+j}$.
 So $\shB=\bigoplus_{n=0}^\infty \shB_n$ is a nonnegatively
graded ring with the degree-zero component $\shB_0=R$.

 One easily observes that $C=\bigoplus_{n=-\infty}^0 C_{-n}$ is
a graded $\shB$\+$B$\+bimodule.
 In fact, the left action map $\shB_i\ot_R C_{i+j}\rarrow C_j$
can be constructed as the composition
$$
 \shB_i\ot_R C_{i+j}\lrarrow \shB_i\ot_RC_i\ot_RC_j\lrarrow
 \shB_i\ot_R C_j
$$
of the map $\shB_i\ot_R C_{i+j}\rarrow \shB_i\ot_RC_i\ot_RC_j$
induced by the comultiplication map $C_{i+j}\rarrow C_i\ot_RC_j$
with the map $\shB_i\ot_RC_i\ot_RC_j\rarrow \shB_i\ot_R C_j$
induced by the evaluation map $\shB_i\ot_RC_i=
\Hom_{R^\rop}(C_i,R)\ot_RC_i\rarrow R$.

 Thus the tensor product ${}^\tau\!K(B,A)=C\ot_B\Ksp_e(B,A)=
C\ot_RA$ is naturally a graded $\shB$\+$A$\+bimodule, and
${}^\tau\!\d:{}^\tau\!K(B,A)\rarrow{}^\tau\!K(B,A)$ is its graded
$\shB$\+$A$\+bimodule endomorphism.
 So, instead of an action of the graded ring $B$, there is an action
of the graded ring $\shB$ in ${}^\tau\!K(B,A)$.
\end{rem}

 We refer to Section~\ref{revisited-subsecn} for a further
discussion of the two Koszul complexes.

\subsection{Distributive collections of subobjects}
\label{distributive-subsecn}
 Let $\sC$ be an abelian category in which all the subobjects of any
given object form a set, and let $W\in\sC$ be an object.
 A set $\Omega$ of subobjects of $W$ is said to be a \emph{lattice
of subobjects} if $0\in\Omega$, \,$W\in\Omega$, and for any
$X$, $Y\in\Omega$ one has $X\cap Y\in\Omega$ and $X+Y\in\Omega$.
 Any lattice of subobjects is \emph{modular}, i.~e., one has
$(X+Y)\cap Z=X+(Y\cap Z)$ whenever $X$, $Y$, $Z\in\Omega$ and
$X\subset Z$.
 A lattice of subobjects $\Omega$ is said to be \emph{distributive}
if the identity $(X+Y)\cap Z=X\cap Z + Y\cap Z$ holds for all
$X$, $Y$, $Z\in\Omega$, or equivalently, the identity
$(X+Y)\cap(X+Z)=X+(Y\cap Z)$ holds for all $X$, $Y$, $Z\in\Omega$.

 Let $n\ge2$ be an integer and $X_1$,~\dots, $X_{n-1}\subset W$ be
a collection of $n-1$ subobjects in~$W$.
 The lattice $\Omega$ of subobjects of $W$ \emph{generated by}
$X_1$,~\dots, $X_{n-1}$ consists of all the subobjects of $W$ that
can be obtained from $X_1$,~\dots, $X_{n-1}$ by applying iteratively
the operations of finite sum and finite intersection.

 The collection of subobjects $(X_1,\dotsc,X_{n-1})$ in $W$ is said
to be \emph{distributive} if the lattice $\Omega$ of subobjects of
$W$ generated by $X_1$,~\dots, $X_{n-1}$ is distributive.
 Any pair of subobjects $X_1$, $X_2$ ($n=3$) forms a distributive
collection, but a triple of subobjects $X_1$, $X_2$, $X_3$ ($n=4$)
does not need to be distributive.
 A collection of subobjects $(X_1,\dotsc,X_{n-1})$ in $W$ is said
to be \emph{almost distributive} if all its proper subcollections
$(X_1,\dotsc,\widehat X_k,\dotsc,X_{n-1})$, \ $1\le k\le n-1$, are
distributive.

 The following two lemmas hold in any modular lattice $\Omega$,
but we state them for lattices of subobjects only.

\begin{lem}
 A triple of subobjects $X$, $Y$, $Z\subset W$ is distributive if
and only if $(X+Y)\cap Z=X\cap Z + Y\cap Z$, and if and only if
$(X+Y)\cap(X+Z)=X+(Y\cap Z)$.
 Any permutation of $X$, $Y$, $Z$ replaces these equations by
equivalent ones.
\end{lem}

\begin{proof}
 This is an easy exercise; see~\cite[Lemma~1]{Jon} or~\cite[Lemma~6.1
in Chapter~1]{PP}.
\end{proof}

\begin{lem}
 An almost distributive collection of subobjects $(X_1,\dotsc,X_{n-1})$
in $W$ is distrubutive if and only if, for every\/ $2\le k\le n-2$,
the triple of subobjects
$$
 X_1+\dotsb+X_{k-1}, \ X_k, \ X_{k+1}\cap\dotsb\cap X_{n-1}\subset W
$$
is distributive.
\end{lem}

\begin{proof}
 This is the result of the paper~\cite{MB}, improving upon
the previous work in~\cite{Jon}.
 See also~\cite[Theorem~6.3 in Chapter~1]{PP}.
\end{proof}

 For any collection of subobjects $X_1$,~\dots, $X_{n-1}\subset W$,
we consider the following three complexes in~$\sC$.
 The \emph{bar complex} $B_\bu=B_\bu(W;X_1,\dotsc,X_{n-1})$ has
the form
\begin{multline} \label{lattice-bar-complex}
 W\lrarrow\bigoplus_t W/X_t\lrarrow\bigoplus_{t_1<t_2}
 W/(X_{t_1}+X_{t_2})\lrarrow\dotsb \\
 \lrarrow\bigoplus_{t_1<\dotsb<t_i}W/(X_{t_1}+\dotsb+X_{t_i})
 \lrarrow\dotsb\lrarrow W/(X_1+\dotsb+X_{n-1}).
\end{multline}
 Here the indices~$t_s$ range from~$1$ to~$n-1$.
 The leftmost term $W$ is placed in the homological degree~$n$,
the term with the summation over $1\le t_1<\dotsb<t_i\le n-1$ is
placed in the homological degree~$n-i$, and the rightmost term
$W/(X_1+\dotsb+X_{n-1})$ is placed in the homological degree~$1$.
 The component of the differential acting from the direct summand
$W/(X_{t_1}+\dotsb+\widehat X_{t_s}+\dotsb+X_{t_i})$ in
$B_{n-i+1}$ to the direct summand $W/(X_{t_1}+\dotsb+X_{t_i})$ in
$B_{n-i}$ is the natural epimorphism taken with the sign $(-1)^{s-1}$.

 The \emph{cobar complex} $B^\bu=B^\bu(W;X_1,\dotsc,X_{n-1})$ has
the form
\begin{multline} \label{lattice-cobar-complex}
 X_1\cap\dotsb\cap X_{n-1}\lrarrow\dotsb\lrarrow
 \bigoplus_{1\le t_1<\dotsb<t_i\le n-1}
 X_{t_1}\cap\dotsb\cap X_{t_i}\lrarrow\dotsb \\ \lrarrow
 \bigoplus_{1\le t_1<t_2\le n-1}X_{t_1}\cap X_{t_2}\lrarrow
 \bigoplus_{1\le t\le n-1} X_t\lrarrow W.
\end{multline}
 Here the leftmost term $X_1\cap\dotsb\cap X_{n-1}$ is placed in
the cohomological degree~$1$, the term with the summation over
$1\le t_1<\dotsb<t_i\le n-1$ is placed in the cohomological
degree~$n-i$, and the rightmost term $W$ is placed in
the cohomological degree~$n$.
 The component of the differential acting from the direct summand
$X_{t_1}\cap\dotsb\cap X_{t_i}$ in $B^{n-i}$ to the direct summand
$X_{t_1}\cap\dotsb\cap\widehat X_{t_s}\cap\dotsb\cap X_{t_i}$ in
$B^{n-i+1}$ is the natural monomorphism taken with the sign
$(-1)^{s-1}$.

 The \emph{Koszul complex} $K_\bu(W;X_1,\dotsc,X_{n-1})$ is
\begin{multline} \label{lattice-Koszul-complex}
 X_1\cap\dotsb\cap X_{n-1}\lrarrow X_2\cap\dotsb\cap X_{n-1}\lrarrow
 X_3\cap\dotsb\cap X_{n-1}/X_1 \\ \lrarrow\dotsb\lrarrow
 (X_{i+1}\cap\dotsb\cap X_{n-1})/(X_1+\dotsb+X_{i-1})
 \lrarrow \dotsb\lrarrow \\
 X_{n-1}/(X_1+\dotsb+X_{n-3})\lrarrow W/(X_1+\dotsb+X_{n-2})
 \lrarrow W/(X_1+\dotsb+X_{n-1}).
\end{multline}
 Here the notation is $Y/Z=Y/(Y\cap Z)=(Y+Z)/Z$.
 The term $X_1\cap\dotsb\cap X_{n-1}$ is placed in the homological
degree~$n$, and the term $W/(X_1+\dotsb+X_{n-1})$ is placed in
the homological degree~$0$.

\begin{lem} \label{almost-distributive-collections}
 Let $X_1$,~\dots, $X_{n-1}$ be an almost distributive collection of
subobjects of an object $W\in\sC$.
 Then the following conditions are equivalent:
\begin{enumerate}
\renewcommand{\theenumi}{\alph{enumi}}
\item the collection $(X_1,\dotsc,X_{n-1})$ is distributive;
\item the Koszul complex $K_\bu(W;X_1,\dotsc,X_{n-1})$ is exact; \par
\item the bar complex $B_\bu(W;X_1,\dotsc,X_{n-1})$ is exact everywhere
except for its leftmost term~$W$; \par
\renewcommand{\theenumi}{\alph{enumi}$^*$}
\setcounter{enumi}{2}
\item the cobar complex $B^\bu(W;X_1,\dotsc,X_{n-1})$ is exact
everywhere except for its rightmost term~$W$.
\end{enumerate}
\end{lem}

\begin{proof}
 This is~\cite[Proposition~7.2 in Chapter~1]{PP}.
 The assertion in~\cite{PP} is stated for collections of subspaces
in vector spaces, but the same argument applies to subobjects of
any object in an abelian category.
\end{proof}

 Let $\sF\subset\sC$ be a class of objects closed under extensions and
the passages to the kernels of epimorphisms.
 We will say that a lattice $\Omega$ of subobjects in an object $W$
is an \emph{$\sF$\+lattice} if one has $Y/Z\in\sF$ for any pair of
subobjects $Y$, $Z\in\Omega$ such that $Z\subset Y$.
 In particular, existence of an $\sF$\+lattice of subobjects in $W$
implies that $W=W/0\in\sF$.
 A collection of subobjects $X_1$,~\dots, $X_{n-1}\subset W$ is said to
be \emph{$\sF$\+distributive} if the lattice $\Omega$ of subobjects in
$W$ generated by $X_1$,~\dots, $X_{n-1}$ is a distributive $\sF$\+lattice.

\begin{lem} \label{F-distributive-lemma}
 A distributive collection of subobjects $X_1$,~\dots,
$X_{n-1}\subset W$ is\/ $\sF$\+distribu\-tive if and only if one has
$W/(X_{t_1}+\dotsb+X_{t_i})\in\sF$ for all $1\le t_1<\dotsb<t_i\le n-1$,
\ $0\le i\le n-1$.
\end{lem}

\begin{proof}
 This is a generalization of~\cite[Lemma~11.4.3.2]{Psemi}, and the same
proof applies.
 Alternatively, one can argue by induction in~$n$ in the following way.
 For any three subobjects $Z\subset Y$ and $X$ in $W$, there is a short
exact sequence $0\rarrow (X\cap Y)/(X\cap Z)\rarrow Y/Z\rarrow
(X+Y)/(X+Z)\rarrow0$.
 Taking $Y$, $Z\in\Omega$ and $X=X_{n-1}$, we observe that the object
$(X_{n-1}+Y)/(X_{n-1}+Z)$ belongs to $\sF$ by the induction assumption
applied to the collection of subobjects $(X_t+X_{n-1})/X_{n-1}\subset
W/X_{n-1}$, \ $1\le t\le n-2$.
 Since $\sF$ is closed under extensions, it suffices to show that
$(X_{n-1}\cap Y)/(X_{n-1}\cap Z)\in\sF$.
 Applying the same argument to $X=X_{n-2}$, etc., we reduce the question
to showing that the object $X_1\cap X_2\cap\dotsb\cap X_{n-1}$ belongs
to~$\sF$.
 The argument finishes similarly to the proof in~\cite{Psemi},
by invoking Lemma~\ref{almost-distributive-collections}\,%
(a)\,$\Rightarrow$\,(c) and the assumption that the class $\sF$
is closed under the kernels of epimorphisms.
\end{proof}

 We will say that a lattice $\Omega$ of subobjects in an object $W$ is
\emph{split} if for any pair of subobjects $Y$, $Z\in\Omega$ such that
$Z\subset Y$ we have that $Z$ is a split subobject of~$Y$.
 Equivalently, $\Omega$ is split if every $Z\in\Omega$ is a split
subobject of~$W$.
 A collection of subobjects $X_1$,~\dots, $X_{n-1}\subset W$ is said
to be \emph{split distributive} if the lattice $\Omega$ of subobjects
in $W$ generated by $X_1$,~\dots, $X_{n-1}$ is split and distributive.

\begin{lem} \label{split-distributive-lemma}
 A collection of subobjects $X_1$,~\dots, $X_{n-1}\subset W$ is split
distributive if and only if there exists a finite direct sum
decomposition $W=\bigoplus_\eta W_\eta$ of the object $W$ such that each
of the subobjects $X_i$ is the sum of a set of subobjects~$W_\eta$.
\end{lem}

\begin{proof}
 The proof of~\cite[Proposition~7.1\,(a)\,$\Leftrightarrow$\,(b)
in Chapter~1]{PP} is applicable.
\end{proof}

\subsection{Collections of subbimodules}
 Let $R$ and $S$ be associative rings and $W$ be an $R$\+$S$\+bimodule.
 Following the terminology of Section~\ref{distributive-subsecn}, we say
that a collection of $R$\+$S$\+subbimodules $X_1$,~\dots, $X_{n-1}
\subset W$ is \emph{distributive} if the lattice of
$R$\+$S$\+subbi\-modules $\Omega$ generated by $X_1$,~\dots, $X_{n-1}$
in $W$ is distributive.
 A collection of subbimodules $X_1$,~\dots, $X_{n-1}\subset W$
is said to be \emph{left flat distributive} if it is distributive and
for every pair of subbimodules $Y$, $Z\in\Omega$ such that $Z\subset Y$
the quotient bimodule $Y/Z$ is a flat left $R$\+module.

\begin{lem} \label{flat-tensor-products-of-collections}
 Let $R$, $S$, and $T$ be associative rings, $W$ be
an $R$\+$S$\+bimodule, and $U$ be an $S$\+$T$\+bimodule.
 Let $X_1$,~\dots, $X_{n-1}\subset W$ be a left flat distributive
collection of $R$\+$S$\+subbimodules, and let
$Y_1$,~\dots, $Y_{m-1}\subset U$ be a left flat distributive
collection of $S$\+$T$\+subbimodules.
 Then
$$
 X_1\ot_SU,\,\dotsc,\,X_{n-1}\ot_SU,\;
 W\ot_SY_1,\,\dotsc,W\ot_SY_{m-1}\,\subset\,W\ot_SU
$$
is a left flat distributive collection of $R$\+$T$\+subbimodules in
the $R$\+$T$\+bimodule $W\ot_SU$.
\end{lem}

\begin{proof}
 First of all, the map $X_i\ot_SU\rarrow W\ot_SU$ induced by
the inclusion $X_i\rarrow W$ is injective for every $1\le i\le n-1$,
since $U$ is a flat left $S$\+module.
 The map $W\ot_SY_j\rarrow W\ot_SU$ induced by the inclusion
$Y_j\rarrow U$ is also injective for every $1\le j\le m-1$, since
$U/Y_j$ is a flat left $S$\+module.
 So $X_i\ot_SU$ and $W\ot_SY_j$ are indeed subbimodules in $W\ot_SU$.

 Furthermore, arguing by induction in $n+m$, we can assume that our
collection of $n+m-2$ subbimodules in $W\ot_SU$ is almost distributive.
 Up to a homological shift by $[-1]$, the bar complex
$B_\bu=B_\bu(W\ot_SU;X_1\ot_SU,\dotsc,X_{n-1}\ot_S\nobreak U,\allowbreak
W\ot_S\nobreak Y_1,\dotsc,W\ot_SY_{m-1})$ \,\eqref{lattice-bar-complex}
is isomorphic to the tensor product of two bar complexes
$B_\bu(W;X_1,\dotsc,X_{n-1})\ot_SB_\bu(U;Y_1,\dotsc,Y_{m-1})$.
 The only nonzero homology bimodule of the bar complex
$B_\bu(U;Y_1,\dotsc,Y_{m-1})$ is $Y_1\cap\dotsc\cap Y_{m-1}\subset U$,
and it is a flat left $S$\+module by assumption.
 So are all the terms of the bar complex
$B_\bu(U;Y_1,\dotsc,Y_{m-1})$.
 Thus the bar complex $B_\bu$ is exact everywhere except for its
leftmost term $W\ot_SU$.
 Applying Lemma~\ref{almost-distributive-collections}\,%
(c)\,$\Rightarrow$\,(a), we can conclude that out collection of
$n+m-2$ subbimodules in $W\ot_SU$ is distributive.

 Finally, the quotient bimodule of $W\ot_SU$ by the sum of any subset
of $X_1\ot_SU$,~\dots, $X_{n-1}\ot_SU$, \ $W\ot_SY_1$,~\dots,
$W\ot_SY_{m-1}$ is isomorphic to the tensor product of the quotient
bimodules of $W$ and $U$ by the sums of the respective subsets of
$X_1$,~\dots, $X_{n-1}$ and $Y_1$,~\dots, $Y_{n-1}$.
 As the tensor product of an $R$\+flat $R$\+$S$\+bimodule and
an $S$\+flat $S$\+$T$\+bimodule is an $R$\+flat $R$\+$T$\+bimodule,
any such quotient bimodule of $W\ot_SU$ is a flat left $R$\+module.
 It remains to apply Lemma~\ref{F-distributive-lemma} (for the class
$\sF$ of $R$\+flat $R$\+$S$\+bimodules in the abelian category $\sC$
of $R$\+$S$\+bimodules) in order to finish the proof of the lemma.
\end{proof}

\begin{lem} \label{tensor-intersections}
 In the context of Lemma~\ref{flat-tensor-products-of-collections},
for any two subbimodules $X'$, $X''\subset W$ belonging to the lattice
of subbimodules generated by $X_1$,~\dots, $X_{n-1}$ in $W$ and
any two subbimodules $Y'$, $Y''\subset U$ belonging to the lattice of
subbimodules generated by $Y_1$,~\dots, $Y_{m-1}$ in $U$, the following
equations for subbimodules in $W\ot_SU$ hold: \par
\textup{(a)} $(X'+X'')\ot_SU=(X'\ot_SU)+(X''\ot_SU)$; \par
\textup{(b)} $(X'\cap X'')\ot_SU=(X'\ot_SU)\cap(X''\ot_SU)$; \par
\textup{(c)} $W\ot_S(Y'+Y'')=(W\ot_SY')+(W\ot_SY'')$; \par
\textup{(d)} $W\ot_S(Y'\cap Y'')=(W\ot_SY')\cap(W\ot_SY'')$; \par
\textup{(e)} $X'\ot_S Y'=(W\ot_SY')\cap(X'\ot_SU)$; \par
\textup{(f)} $(X'\cap X'')\ot_S(Y'\cap Y'')=(X'\ot_SY')\cap(X''\ot_SY'')$.
\end{lem}

\begin{proof}
 The maps $X'\ot_SU\rarrow W\ot_SU$ and $W\ot_SY'\rarrow W\ot_SU$
induced by the inclusions $X'\rarrow W$ and $Y'\rarrow U$ are injective,
as explained in the proof of
Lemma~\ref{flat-tensor-products-of-collections}.
 Furthermore, the map $X'\ot_SY'\rarrow W\ot_SU$ is injective as
the composition of injective maps $X'\ot_SY'\rarrow W\ot_SY'\rarrow
W\ot_SU$ or $X'\ot_SY'\rarrow X'\ot_SU\rarrow W\ot_SU$.
 So these are indeed subbimodules in $W\ot_SU$.
 
 Now the equations~(a) and~(c) are obvious.
 The equation~(b) holds since $U$ is flat left $S$\+module.
 To prove~(d), it suffices to consider the four-term exact sequence
of flat left $S$\+modules
$$
 0\lrarrow Y'\cap Y''\lrarrow Y'\oplus Y''\lrarrow U
 \lrarrow U/(Y'+Y'')\lrarrow0
$$
and tensor it with $W$ over $S$ on the left.

 To prove~(e), one has to check that that the short sequence
$0\rarrow X'\ot_SY'\rarrow W\ot_SU\rarrow (W/X'\ot_SU)\oplus(W\ot_SU/Y')$
is left exact.
 This follows from exactness of the sequences $0\rarrow X'\ot_SY'
\rarrow W\ot_SY'\rarrow (W/X')\ot_SY'$ and
$0\rarrow W\ot_SY'\rarrow W\ot_SU\rarrow W\ot_S(U/Y')$, and
injectivity of the map $(W/X')\ot_SY'\rarrow(W/X')\ot_SU$.
 The equation~(f) follows from~(b), (d), and~(e).
\end{proof}

 We will say that a collection of subbimodules
$X_1$,~\dots, $X_{n-1}\subset W$ in an $R$\+$S$\+bimod\-ule $W$ is
\emph{left projective distibutive} if it is distributive and for every
pair of subbimodules $Y$, $Z\in\Omega$ such that $Z\subset Y$
the quotient bimodule $Y/Z$ is a projective left $R$\+module.
 A collection of subbimodules $X_1$,~\dots, $X_{n-1}\subset W$ is said
to be \emph{left split distributive} if it is split distributive as
a collection of submodules in the left $R$\+module $W$ (in
the sense of the definition in Section~\ref{distributive-subsecn}).
 \emph{Right projective distributive} and \emph{right split
distributive} collections of subbimodules are defined similarly.

 Clearly, a collection of subbimodules in an $R$\+$S$\+bimodule $W$
is left projective distributive if and only if it is left split
distributive \emph{and} the left $R$\+module $W$ is projective.

\begin{lem}
 Let $R$, $S$, and $T$ be associative rings, $W$ be
an $R$\+$S$\+bimodule, and $U$ be an $S$\+$T$\+bimodule.
 Let $X_1$,~\dots, $X_{n-1}\subset W$ be a left projective distributive
collection of $R$\+$S$\+subbimodules, and let
$Y_1$,~\dots, $Y_{m-1}\subset U$ be a left projective distributive
collection of $S$\+$T$\+subbimodules.
 Then
$$
 X_1\ot_SU,\,\dotsc,\,X_{n-1}\ot_SU,\;
 W\ot_SY_1,\,\dotsc,W\ot_SY_{m-1}\,\subset\,W\ot_SU
$$
is a left projective distributive collection of subbimodules
in the $R$\+$T$\+bimodule $W\ot_SU$.
\end{lem}

\begin{proof}
 Similar to the proof of
Lemma~\ref{flat-tensor-products-of-collections}.
\end{proof}

\begin{lem} \label{dual-collection}
 Let $X_1$,~\dots, $X_{n-1}\subset W$ be a left split distributive
collection of subbimodules in an $R$\+$S$\+bimodule~$W$.
 Then\/ $\Hom_R(W/X_1,R)$,~\dots, $\Hom_R(W/X_{n-1},R)\allowbreak\subset
\Hom_R(W,R)$ is a right split distributive collection of subbimodules
in the $S$\+$R$\+bimod\-ule\/ $\Hom_R(W,R)$.
 If the left $R$\+module $W$ is projective and finitely generated,
then the above collection of subbimodules in\/ $\Hom_R(W,R)$ is right
projective distributive.
 In this case, the map $W\supset Y\longmapsto\Hom_R(W/Y,R)
\subset\Hom_R(W,R)$ is an anti-isomorphism between the lattice of
subbimodules in $W$ generated by $X_1$,~\dots, $X_{n-1}$ and the lattice
of subbimodules in\/ $\Hom_R(W,R)$ generated by\/
$\Hom_R(W/X_1,R)$,~\dots, $\Hom_R(W/X_{n-1},R)$ (i.~e., a bijection
of sets transforming the sums into the intersections and vice versa).
\end{lem}

\begin{proof}
 Follows easily from Lemma~\ref{split-distributive-lemma} (applied in
the category of left $R$\+modules $\sC=R\modl$).
\end{proof}

\subsection{Left flat Koszul rings} \label{flat-Koszul-subsecn}
 Let $A=\bigoplus_{n=0}^\infty A_n$ be a nonnegatively graded ring with
the degree-zero component $R=A_0$.
 The graded ring $A$ is said to be \emph{left flat Koszul} if one of
the equivalent conditions of the next theorem is satisfied.
 (\emph{Right flat Koszul} graded rings are defined similarly.)
 
\begin{thm} \label{flat-koszul-theorem}
 The following three conditions are equivalent: \par
\textup{(a)} $A_n$ is a flat left $R$\+module for every $n\ge1$ and\/
$\Tor_{i,j}^A(R,R)=0$ for all\/ $i\ne j$; \par
\textup{(b)} $A_n$ is a flat left $R$\+module for every $n\ge1$,
the graded ring $A$ is quadratic, and, setting $V=A_1$ and denoting by
$I\subset V\ot_RV$ the kernel of the multiplication map $A_1\ot_R A_1
\rarrow A_2$, for every $n\ge4$ the collection of\/ $n-1$ subbimodules
\begin{equation} \label{relation-subbimodules-collection}
 I\ot_R V\ot_R\dotsb\ot_R V, \
 V\ot_R I\ot_R V\ot_R\dotsb\ot_R V, \,\dotsc, \
 V\ot_R\dotsb\ot_R V\ot_R I
\end{equation}
in the $R$\+$R$\+bimodule $W=V^{\ot_R\,n}$ is distributive; \par
\textup{(c)} $A_1$ and $A_2$ are flat left $R$\+modules, the graded
ring $A$ is quadratic, and for every $n\ge1$ the collection of
$R$\+$R$\+subbimodules $V^{\ot_R\,k-1}\ot_RI\ot_RV^{\ot_R\,n-k-1}$,
\ $1\le k\le n-1$ \,\eqref{relation-subbimodules-collection}
in the $R$\+$R$\+bimodule $W=V^{\ot_R\,n}$ is left flat distributive.
\end{thm}

\begin{proof}
 First of all, by Proposition~\ref{low-dimensional-Tor-prop}(b),
condition~(a) implies that $A$ is quadratic.
 So we can assume that $A=T_R(V)/(I)$.
 Furthermore, any collection of less than three subobjects is
distributive; so the distributivity condition in~(b) is trivial
for $n\le3$.
 It is clear from the formula~\eqref{quadratic-ideal-component}
from Section~\ref{quadratic-duality-secn} that the left flat
distributivity condition in~(c) implies flatness of
the left $R$\+modules $A_n$ for all $n\ge1$.
 This suffices to prove the implication (c)\,$\Longrightarrow$\,(b).
 
 To deduce (b)\,$\Longrightarrow$\,(c), we observe that the quotient
bimodule of $V^{\ot_R\,n}$ by any subset of
the subbimodules~\eqref{relation-subbimodules-collection} is
isomorphic to the tensor product $A_{j_1}\ot_R\dotsb\ot_RA_{j_s}$
for some $j_1$,~\dots, $j_s\ge1$, \ $j_1+\dotsb+j_s=n$.
 Flatness of the left $R$\+modules $A_j$ implies flatness of such
tensor products (as left $R$\+modules), and it remains to invoke
Lemma~\ref{F-distributive-lemma} for the abelian category $\sC$ of
$R$\+$R$\+bimodules and the class $\sF\subset\sC$ of all
$R$\+$R$\+bimodules that are flat as left $R$\+modules.

 It remains to prove the equivalence (a)\,$\Longleftrightarrow$\,(b).
 Arguing by induction in~$n\ge1$, we can assume that the collection
of $j-1$ subbimodules~\eqref{relation-subbimodules-collection}
in the $R$\+$R$\+bimodule $V^{\ot_R\,j}$ is left flat distributive
for all $1\le j\le n-1$ and $\Tor_{i,j}^A(R,R)=0$ for all
$i\ne j\le n-1$.
 Under this induction assumption, we will prove the equivalence
of conditions~(a) and~(b) for the fixed value of $j=n$.

 By Lemma~\ref{flat-tensor-products-of-collections} (applied to
the lattice of $k-1$ subbimodules in the $R$\+$R$\+bimodule
$W=V^{\ot_R\,k}$ and the lattice of $n-k-1$ subbimodules in
the $R$\+$R$\+bimodule $U=V^{\ot_R\,n-k}$, \ $1\le k\le n-1$),
the induction assumption implies that the collection of $n-1$
subbimodules~\eqref{relation-subbimodules-collection} in
the $R$\+$R$\+bimodule $V^{\ot_R\,n}$ is almost distributive.

 Finally, the bar complex~\eqref{Tor-R-R-bar-complex} computing
the $R$\+$R$\+bimodules $\Tor^A_{i,n}(R,R)$ is isomorphic to
the lattice bar-complex $B_\bu(W;X_1,\dotsc,X_{n-1})$
\,\eqref{lattice-bar-complex} for the collection of
$n-1$ subbimodules $X_k=V^{\ot_R\,k-1}\ot_RI\ot_RV^{\ot_R\,n-k-1}$
in the $R$\+$R$\+bimodule $W=V^{\ot_R\,n}$.
 Hence Lemma~\ref{almost-distributive-collections}%
(a)\,$\Leftrightarrow$\,(c) implies the desired equivalence
(a)\,$\Longleftrightarrow$\,(b).
\end{proof}

 The construction of the first Koszul complex in
Section~\ref{first-koszul-complex-subsecn} was using double dualization:
first we passed from $A_1$ to $B_1=\Hom_R(A_1,R)$, and then set
$K_1^\tau(B,A)=\Hom_{R^\rop}(B_1,A)$.
 Therefore, the assumption that $A_1$ is a finitely generated
projective $R$\+module was needed.
 A similar double dualization was used in the construction of
the second Koszul complex in Section~\ref{dual-koszul-complex-subsecn}.
 The following alternative approach allows to produce the two Koszul
complexes for any left flat Koszul ring.

 Let $A$ be a left flat Koszul graded ring.
 Denote the kernel of the (surjective) multiplication map
$A_1\ot_RA_1\rarrow A_2$ by $I_A\subset A_1\ot_RA_1$.
 Set $I_A^{(0)}=R$, \ $I_A^{(1)}=A_1$, \ $I_A^{(2)}=I_A$, and
\begin{equation} \label{dual-coring-components}
 I_A^{(n)}=\bigcap\nolimits_{k=1}^{n-1}
 A_1^{\ot_R\,k-1}\ot_R I_A\ot_RA_1^{\ot_R\,n-k-1}
 \,\subset\,A_1^{\ot_R\,n}.
\end{equation}
 So $I_A^{(n)}$ is an $R$\+$R$\+subbimodule in $A_1^{\ot_R\,n}$.

 The intersection of all the subbimodules indexed by $k=2$,~\dots,
$n-1$ (i.~e., of all but the first one)
in~\eqref{dual-coring-components} is the subbimodule
$$
 \bigcap\nolimits_{k=2}^{n-1}
 A_1^{\ot_R\,k-1}\ot_R I_A\ot_RA_1^{\ot_R\,n-k-1}=
 A_1\ot_R I_A^{(n-1)}\,\subset\,A_1^{\ot_R\,n}
$$
by Lemma~\ref{tensor-intersections}(d).
 Hence we obtain an injective $R$\+$R$\+bimodule morphism
$I_A^{(n)}\rarrow A_1\ot_RI_A^{(n-1)}$, which is defined for
all $n\ge1$.
 For every $n\ge2$, the image of the composition $I_A^{(n)}\rarrow
A_1\ot_RI_A^{(n-1)}\rarrow A_1\ot_RA_1\ot_RI_A^{(n-2)}$ is
contained in the subbimodule $I_A\ot_RI_A^{(n-2)}\subset
A_1\ot_RA_1\ot_RI_A^{(n-2)}$.

 Put $K_{i,n}^\tau(A)=A_{n-i}\ot_RI_A^{(i)}$ for every $i\ge0$,
\,$n\ge i$ (where $\tau$~is just a placeholder or a notation for
the identity map $A_1\rarrow A_1$).
 Define the differential $\d\:K^\tau_i(A)\rarrow K^\tau_{i-1}(A)$
as the composition
$$
 A_j\ot_RI_A^{(i)}\lrarrow A_j\ot_RA_1\ot_RI_A^{(i-1)}
 \lrarrow A_{j+1}\ot_RI_A^{(i-1)}
$$
of the map $A_j\ot_RI_A^{(i)}\rarrow A_j\ot_RA_1\ot_RI_A^{(i-1)}$
induced by the map $I_A^{(i)}\rarrow A_1\ot_RI_A^{(i-1)}$ and the map
$A_j\ot_RA_1\ot_RI_A^{(i-1)}\rarrow A_{j+1}\ot_RI_A^{(i-1)}$
induced by the multiplication map $A_j\ot_RA_1\rarrow A_{j+1}$.
 Since the composition $I_A\rarrow A_1\ot_RA_1\rarrow A_2$ vanishes,
we have $\d^2=0$.
 So we obtain a complex of graded left $A$\+modules $K^\tau_\bu(A)$ with
the homological grading~$i$ and the internal grading~$n=i+j$.

\begin{prop} \label{left-flat-first-Koszul-complex}
 For any left flat Koszul graded ring $A$, the first Koszul complex
$K^\tau_\bu(A)$ is a graded flat resolution of the left $A$\+module~$R$.
\end{prop}

\begin{proof}
 The graded left $A$\+modules $A\ot_RI_A^{(i)}$ are flat, because
the left $R$\+modules $I_A^{(i)}$ are.
 Furthermore, for every $n\ge1$ the internal degree~$n$ component of
the Koszul complex $K^\tau_\bu(A)$ is isomorphic to the lattice Koszul
complex $K_\bu(W;X_1,\dotsc,X_{n-1})$ \,\eqref{lattice-Koszul-complex}
for the collection of $n-1$ subbimodules $X_k=A_1^{\ot_R\,k-1}
\ot_RI_A\ot_RA_1^{\ot_R\,n-k-1}$ in the $R$\+$R$\+bimodule
$W=A_1^{\ot_R\,n}$.
 Thus it remains to apply
Lemma~\ref{almost-distributive-collections}\,(a)\,$\Rightarrow$\,(b).
\end{proof}

 Similarly, the intersection of all the subbimodules indexed by
$k=1$,~\dots, $n-2$ (i.~e., of all but the last one)
in~\eqref{dual-coring-components} is the subbimodule
$$
 \bigcap\nolimits_{k=1}^{n-2}
 A_1^{\ot_R\,k-1}\ot_R I_A\ot_RA_1^{\ot_R\,n-k-1}=
 I_A^{(n-1)}\ot_R A_1\,\subset\,A_1^{\ot_R\,n}
$$
by Lemma~\ref{tensor-intersections}(b).
 Hence we obtain an injective $R$\+$R$\+bimodule morphism
$I_A^{(n)}\rarrow I_A^{(n-1)}\ot_RA_1$, which is defined for
all $n\ge1$.

 Put ${}^\tau\!K_{i,n}(A)=I_A^{(i)}\ot_RA_{n-i}$ for every $i\ge0$,
\,$n\ge i$.
 Define the differential $\d:{}^\tau\!K_i(A)\rarrow{}^\tau\!K_{i-1}(A)$
as the composition
$$
 I_A^{(i)}\ot_RA_j\lrarrow I_A^{(i-1)}\ot_RA_1\ot_RA_j
 \lrarrow I_A^{(i-1)}\ot_RA_{j+1}.
$$
 Similarly to the construction above, we obtain a complex of graded
right $A$\+modules ${}^\tau\!K_\bu(A)$ with the homological grading~$i$
and the internal grading $n=i+j$.

 Recall the definition of a \emph{weakly $A/R$\+flat right $A$\+module}
for a morphism of associative rings $R\rarrow A$ such that $A$ is
a flat left $R$\+module~\cite[Section~5]{Pfp}.
 A right $A$\+module $F$ is said to be weakly $A/R$\+flat (or
\emph{weakly flat relative to~$R$}) if the functor $F\ot_A{-}$ takes
short exact sequences of $R$\+flat left $A$\+modules to short exact
sequences of abelian groups.
 For any right $R$\+module $N$, the right $A$\+module $N\ot_RA$ is
weakly flat relative to~$R$.
 The main properties of the class of weakly $A/R$\+flat right
$A$\+modules are listed in~\cite[Lemma~5.3(b)]{Pfp}.

\begin{prop} \label{left-flat-second-Koszul-complex}
 For any left flat Koszul graded ring $A$, the second Koszul complex
${}^\tau\!K_\bu(A)$ is a graded weakly $A/R$\+flat resolution of
the right $A$\+module~$R$.
\end{prop}

\begin{proof}
 In view of the discussion above, the (graded) right $A$\+modules
$I_A^{(i)}\ot_RA$ are (graded) weakly $A/R$\+flat.
 This assertion does not require any assumptions about the right
$R$\+modules~$I_A^{(i)}$.
 Furthermore, for every $n\ge1$ the internal degree~$n$ component of
the Koszul complex ${}^\tau\!K_\bu(A)$ is isomorphic to the lattice
Koszul complex $K_\bu(W;X_1,\dotsc,X_{n-1})$
\,\eqref{lattice-Koszul-complex} for the collection of $n-1$
subbimodules $X_k=A_1^{\ot_R\,n-k-1}\ot_RI_A\ot_RA_1^{\ot_R\,k-1}$
in the $R$\+$R$\+bimodule $W=A_1^{\ot_R\,n}$.
 This is the same collection of subbimodules as in the proof of
Proposition~\ref{left-flat-first-Koszul-complex}, but numbered in
the opposite order.
 So it remains to apply
Lemma~\ref{almost-distributive-collections}\,(a)\,$\Rightarrow$\,(b)
in order to conclude that the internal degree~$n$ component of
the complex ${}^\tau\!K_\bu(A)$ is exact for $n>0$.
 Thus $H_i({}^\tau\!K_\bu(A))=0$ for $i>0$ and
$H_0({}^\tau\!K_\bu(A))=R$.
\end{proof}
 
\subsection{Finitely projective Koszul rings}
\label{finitely-projective-Koszul-subsecn}
 Let $A=\bigoplus_{n=0}^\infty A_n$ be a nonnegatively graded ring
with the degree-zero component $R=A_0$.
 Assume that $A_1$ and $A_2$ are finitely generated projective left
$R$\+modules.
 When the multiplication map $A_1\ot_RA_1\rarrow A_2$ is surjective,
we denote its kernel by $I\subset V\ot_RV$, where $V=A_1$, and consider
the $2$\+left finitely projective quadratic graded ring
$\q A=T_R(V)/(I)$ together with its quadratic dual $2$\+right finitely
projective quadratic graded ring~$B$.

 As above in this paper, we denote by $\Ext_A(R,R)$ the bigraded $\Ext$
ring of the graded \emph{left} $A$\+module $R$, and by
$\Ext_{B^\rop}(R,R)$ the bigraded $\Ext$ ring of the graded \emph{right}
$B$\+module~$R$.
 The graded ring $A$ is said to be \emph{left finitely projective
Koszul} if one of the equivalent conditions of the next theorem
is satisfied.

\begin{thm} \label{projective-koszul-theorem}
 For any nonnegatively graded ring $A$, the following four conditions
are equivalent: \par
\textup{(a)} $A_n$ is a finitely generated projective left $R$\+module
for every $n\ge1$ and the graded ring $A$ is left flat Koszul; \par
\textup{(b)} $A_n$ is a finitely generated projective left $R$\+module
for every $n\ge1$, the multiplication map $A_1\ot_RA_1\rarrow A_2$ is
surjective, $B_n$ is a finitely generated projective right $R$\+module
for every $n\ge1$, and\/ $\Ext_A^{i,j}(R,R)=0$ for all\/ $i\ne j$; \par
\textup{(c)} the graded ring $A$ is quadratic, $A_1$ and $A_2$ are
finitely generated projective left $R$\+modules, and for every $n\ge1$
the collection of $R$\+$R$\+subbimodules
$V^{\ot_R\,k-1}\ot_RI\ot_RV^{\ot_R\,n-k-1}$, \ $1\le k\le n-1$
in the $R$\+$R$\+bimodule $W=V^{\ot_R\,n}$ is left projective
distributive; \par
\textup{(d)} $A_1$ and $A_2$ are finitely generated projective left
$R$\+modules, the multiplication map $A_1\ot_RA_1\rarrow A_2$ is
surjective, $B_i$ is a finitely generated projective right $R$\+module
for every $i\ge1$, and the first Koszul complex $K^\tau_\bu(B,A)$
\,\eqref{first-koszul-complex} from
Section~\ref{first-koszul-complex-subsecn} is exact in all
the internal degrees $n\ge1$; \par
\textup{(e)} $A_n$ is a finitely generated projective left $R$\+module
for every $n\ge1$, the multiplication map $A_1\ot_RA_1\rarrow A_2$ is
surjective, and the second Koszul complex ${}^\tau\!K_\bu(B,A)$
\,\eqref{second-koszul-complex} from
Section~\ref{dual-koszul-complex-subsecn} is exact in all
the internal degrees $n\ge1$.
\end{thm}

\begin{proof}
 (a)\,$\Longleftrightarrow$\,(c) Take the condition of
Theorem~\ref{flat-koszul-theorem}(b) as the definition of left flat
Koszulity and argue similarly to the proof of
Theorem~\ref{flat-koszul-theorem}\,(b)\,$\Leftrightarrow$\,(c),
using Lemma~\ref{F-distributive-lemma} for the abelian category
$\sC$ of $R$\+$R$\+bimodules and the class $\sF\subset\sC$ of all
$R$\+$R$\+bimodules that are projective as left $R$\+modules.

 (a)\,$\Longleftrightarrow$\,(b) Take the condition of
Theorem~\ref{flat-koszul-theorem}(a) as the definition of left flat
Koszulity.
 The $R$\+$R$\+bimodules $\Tor^A_{i,n}(R,R)$ can be computed as
the homology of the complex~\eqref{Tor-R-R-bar-complex}, while
the $R$\+$R$\+bimodules $\Ext_A^{i,n}(R,R)$ can be computed as
the cohomology of the complex~\eqref{Ext-R-R-cobar-complex}.
 Under~(a), the graded ring $A$ is quadratic; so any one of
the conditions~(a) or~(b) implies surjectivity of the multiplication
map $A_1\ot_RA_1\rarrow A_2$.
 By Proposition~\ref{diagonal-Ext}, we have $\Ext_A^{n,n}(R,R)
\simeq B_n$.

 Now the following lemma shows that, given a nonnegatively graded
algebra $A$ with the degree-zero component $R=A_0$ such that
$A_j$ is a finitely generated projective left $R$\+module for
every $j\ge1$, and given a fixed $n\ge1$, one has
$\Tor^A_{i,n}(R,R)=0$ for all $1\le i\le n-1$ if and only if
$\Ext_A^{i,n}(R,R)=0$ for all $1\le i\le n-1$ \emph{and}
the right $R$\+module $\Ext_A^{n,n}(R,R)$ is projective.

\begin{lem}
 Let\/ $0\rarrow C_n\rarrow C_{n-1}\rarrow\dotsb\rarrow C_1\rarrow0$
be a complex of finitely generated projective left $R$\+modules,
and let\/ $0\rarrow\Hom_R(C_1,R)\rarrow\dotsb\rarrow\Hom_R(C_n,R)
\rarrow0$ be the dual complex of finitely generated projective
right $R$\+modules.
 Then the following four conditions are equivalent: \par
\textup{(a)} $H_i(C_\bu)=0$ for all\/ $1\le i\le n-1$; \par
\textup{(a$'$)} the complex of left $R$\+modules $C_\bu$ is homotopy
equivalent to the one-term complex $H_n(C_\bu)[n]$; \par
\textup{(b)} $H^i(\Hom_R(C_\bu,R))=0$ for all\/ $1\le i\le n-1$
and $H^n(\Hom_R(C_\bu,R))$ is a (finitely generated) projective
right $R$\+module; \par
\textup{(b$'$)} the complex of right $R$\+modules $\Hom_R(C_\bu,R)$
is homotopy equivalent to the one-term complex
$H^n(\Hom_R(C^\bu,R))[-n]$. \par
 If any one of these equivalent conditions holds, then $H_n(C_\bu)$ is
a finitely generated projective left $R$\+module and
$H^n(\Hom_R(C_\bu,R))\simeq\Hom_R(H_n(C_\bu),R)$. \qed
\end{lem}

 (c)\,$\Longrightarrow$\,(d),~(e) Under~(c), the results of
Lemmas~\ref{tensor-intersections} and~\ref{dual-collection} are
available for the collections of subbimodules $V^{\ot_R\,k-1}\ot_R
I\ot_RV^{\ot_R\,j-k-1}$ in the $R$\+$R$\+bimodules $V^{\ot_R\,j}$,
\ $1\le k\le j-1$, \ $j\ge 1$.
 This allows to compute the internal degree~$n$ component of the Koszul
complex $K^\tau_\bu(B,A)$ as the lattice Koszul complex
$K_\bu(W;X_1,\dotsc,X_{n-1})$ \,\eqref{lattice-Koszul-complex} for
the collection of subbimodules
$X_k=V^{\ot_R\,k-1}\ot_RI\ot_RV^{\ot_R\,n-k-1}\subset
V^{\ot_R\,n}=W$.
 Similarly, the internal degree~$n$ component of the Koszul
complex ${}^\tau\!K_\bu(B,A)$ is computed as the lattice Koszul complex
for the collection of subbimodules $X_k=V^{\ot_R\,n-k-1}\ot_RI\ot_R
V^{\ot_R\,k-1}\subset V^{\ot_R\,n}$ (the same subbimodules numbered in
the opposite order).
 Then in both cases it remains to apply
Lemma~\ref{almost-distributive-collections}\,(a)\,$\Rightarrow$\,(b).
 
 (d)\,$\Longrightarrow$\,(a),~(b) Under~(d), the Koszul complex
$K^\tau_\bu(B,A)$ is a graded projective resolution of the graded left
$A$\+module~$R$.
 Computing $\Tor^A(R,R)$ and $\Ext_A(R,R)$ in terms of this resolution
immediately yields (a) and~(b) (where (a)~is interpreted as
the condition of Theorem~\ref{flat-koszul-theorem}(a)).

 (e)\,$\Longrightarrow$\,(a) Under~(e), the Koszul complex
${}^\tau\!K_\bu(B,A)$ is a graded weakly $A/R$\+flat resolution of
the graded right $A$\+module $R$, as explained in
Section~\ref{flat-Koszul-subsecn} (see the proof of
Proposition~\ref{left-flat-second-Koszul-complex}).
 Since the graded left $A$\+module $R$ is flat over $R$, one can compute
$\Tor^A(R,R)$ in terms of this resolution (e.~g., in view
of~\cite[Lemma~5.3(b)]{Pfp}).
 This immediately implies~(a) (intepreted as the condition of
Theorem~\ref{flat-koszul-theorem}(a)).
\end{proof}

 The definition of \emph{right finitely projective Koszul} graded ring
is obtained from the above definition of a left finitely projective
Koszul ring by switching the roles of the left and right sides.

\begin{prop} \label{finitely-projective-Koszul-duality}
 Let $A$ be a $2$\+left finitely projective quadratic graded ring with
the degree-zero component $R=A_0$, and let $B$ be the $2$\+right
finitely projective quadratic graded ring quadratic dual to~$A$.
 Then the ring $A$ is left finitely projective Koszul if and only if
the ring $B$ is right finitely projective Koszul.
 If this is the case, one has
$$
 B\simeq\Ext_A(R,R)^\rop \quad\text{and}\quad
 A\simeq\Ext_{B^\rop}(R,R).
$$
\end{prop}

\begin{proof}
 The first assertion is provable using the condition of
Theorem~\ref{projective-koszul-theorem}(c) as the definition
of projective Koszulity and the result of
Lemma~\ref{dual-collection}.
 Then, by Theorem~\ref{projective-koszul-theorem}(b), we have
$\Ext_A^{i,j}(R,R)=0$ for $i\ne j$, and similarly,
$\Ext_{B^\rop}^{i,j}(R,R)=0$ for $i\ne j$.
 Finally, the graded ring isomorphism $B^\rop\simeq\bigoplus_{n\ge0}
\Ext_A^{n,n}(R,R)$ is provided by Proposition~\ref{diagonal-Ext}.
 The graded ring isomorphism $A\simeq\bigoplus_{n\ge0}
\Ext_{B^\rop}^{n,n}(R,R)$ can be obtained by applying the same
proposition to the ring $B^\rop$ and observing that the $2$\+left
finitely projective quadratic graded ring $B^\rop$ is quadratic
dual to the $2$\+right finitely projective quadratic
graded ring $A^\rop$.
\end{proof}

\begin{cor} \label{homogeneous-koszul-duality-cor}
 The anti-equivalences of categories from
Propositions~\ref{2-fin-proj-quadratic-duality}\+-%
\ref{3-fin-proj-quadratic-duality} restrict to an anti-equivalence
between the category of left finitely projective Koszul graded rings
$A$ and the category of right finitely projective Koszul graded rings 
$B$ over any fixed base ring~$R$. \qed
\end{cor}

\begin{rem} \label{dual-koszul-complex-no-good}
 Let us warn the reader that, even when the graded rings $A$ and
$B$ are Koszul (say, $A$ is left finitely projective Koszul and
$B$ is right finitely projective Koszul, as above), the dual Koszul
complex $\Ksp_e{}^\bu(B,A)$ \,\eqref{dual-koszul-complex} from
Section~\ref{dual-koszul-complex-subsecn}
has \emph{no} good exactness properties, generally speaking.
 In fact, in this case one has $\Ksp_e{}^\bu(B,A)=
\Hom_A(K^\tau_\bu(B,A),A)$, hence it follows from
Theorem~\ref{projective-koszul-theorem}(d) that the complex
$\Ksp_e{}^\bu(B,A)$ computes the graded $\Ext_A^*(R,A)$.
 This Ext can be quite complicated even for a Koszul algebra $A$ over
a field $R=k$.
\end{rem}

\Section{Relative Nonhomogeneous Quadratic Duality}
\label{nonhomogeneous-quadratic-secn}

 Nonhomogeneous quadratic rings $\tA$ can be informally described
as rings defined by nonhomogeneous quadratic relations over a fixed
base ring~$R$.
 Not every system of nonhomogeneous quadratic relations is good enough to
define a nonhomogeneous quadratic ring (see the general discussion 
in the introductions to~\cite{Pcurv,PP}, \cite[Section~5]{PP}, and
the counterexamples in~\cite[Section~5]{Pqf}
and~\cite[Section~3.4]{Pcurv}).
 For a system of nonhomogeneous quadratic relations to ``make sense'',
its coefficients must, in turn, satisfy a certain system of equations,
called the \emph{self-consistency equations}.

 Dualizing the degree-one and degree-zero parts of the nonhomogeneous
quadratic relations and imposing the equations dual to
the self-consistency equations on the relations' coefficients produces
a curved DG\+ring structure $(B,d,h)$ on the dual quadratic graded ring
$B$ to the quadratic graded ring $A$ defined by the homogeneous
quadratic parts of the nonhomogeneous quadratic relations.

\subsection{Nonhomogeneous quadratic rings}
\label{nonhomogeneous-quadratic-subsecn}
 Let $\tA$ be an associative ring with a subring $R\subset\tA$.
 Consider the $R$\+$R$\+bimodule $\tA/R$ and suppose that we have
chosen a subbimodule $V\subset\tA/R$.
 Denote by $\tV\subset\tA$ the full preimage of $V$ under
the surjective $R$\+$R$\+bimodule morphism $\tA\rarrow\tA/R$.

 Consider the tensor ring $T_R(\tV)=
\bigoplus_{n=0}^\infty T_{R,n}(\tV)$, \ $T_{R,n}(\tV)=
\tV^{\ot_R\,n}$, as in Section~\ref{quadratic-duality-secn}.
 Then the $R$\+$R$\+bimodule morphism of identity inclusion
 $\tV\rarrow\tA$ extends uniquely to a ring homomorphism
$\pi_{\tA}\:T_R(\tV)\rarrow\tA$ forming a commutative triangle diagram
with the subring inclusions $R\rarrow T_R(\tV)$ and $R\rarrow\tA$.

 Assume that the ring $\tA$ is generated by its subgroup~$\tV$.
 In other words, the ring homomorphism $\pi_{\tA}\:T_R(\tV)\rarrow\tA$
is surjective.
 Let $\tJ_{\tA}\subset T_R(\tV)$ be the kernel of~$\pi_{\tA}$.

 Define an increasing filtration $F$ on the ring $T_R(\tV)$ by
the rule $F_nT_R(\tV)=\bigoplus_{i=0}^nT_{R,i}(\tV)\subset T_R(\tV)$.
 Furthermore, put $F_n\tA=\pi_{\tA}(F_nT_R(\tV))\subset\tA$.
 So we have $F_0T_R(\tV)=R$ and $F_1T_R(\tV)=R\oplus\tV$, hence
$F_0\tA=R$ and $F_1\tA=\tV\subset\tA$.

 Clearly, $F$ is an exhaustive multiplicative filtration on the ring
$\tA$, that is $\tA=\bigcup_{n\ge0}F_n\tA$ and $F_n\tA\cdot
F_m\tA\subset F_{n+m}\tA$.
 In fact, the filtration $F$ on $\tA$ is \emph{generated by $F_1$},
which means that the equation $F_n\tA\cdot F_m\tA = F_{n+m}\tA$
holds for all $n$, $m\ge0$.
 (Sometimes we will say that the filtration $F$ on $\tA$ is
\emph{generated by $F_1\tA$ over $R$}, which means that $F$ is
generated by $F_1$ and $F_0\tA=R$.)
 The filtration $F$ on the graded ring $T_R(\tV)$ has similar
properties.

 Put $\tI_{\tA}=F_2T_R(\tV)\cap\tJ_{\tA}$; so $\tI_{\tA}$ is
an $R$\+$R$\+subbimodule in $F_2T_R(\tV)=
R\oplus\tV\oplus\tV^{\ot_R\,2}$.
 We will say that the ring $\tA$ with a fixed subring $R\subset\tA$
and a fixed $R$\+$R$\+subbimodule of generators $R\subset\tV\subset\tA$
is \emph{weak nonhomogeneous quadratic} over $R$ if the ideal
$\tJ_{\tA}\subset T_R(\tV)$ is generated by its subgroup~$\tI_{\tA}$.

 Let $\prtA$ and $\sectA$ be two weak nonhomogeneous quadratic rings
over the same base ring $R$, and let $R\subset\prtV\subset\prtA$ and
$R\subset\sectV\subset\sectA$ be their fixed subbimodules of generators.
 A \emph{morphism} of weak nonhomogeneous quadratic rings $f\:\prtA
\rarrow\sectA$ is a ring homomorphism forming a commutative triangle
diagram with the subring inclusions $R\rarrow\prtA$ and $R\rarrow\sectA$
and satisfying the condition that $f(\prtV)\subset\sectV$.

 Conversely, suppose that we are given an associative ring $\tA$
with an exhaustive multiplicative increasing filtration $0=F_{-1}\tA
\subset F_0\tA\subset F_1\tA\subset F_2\tA\subset\dotsb$.
 Then $R=F_0\tA$ is a subring in~$\tA$.
 Suppose that the filtration $F$ on $\tA$ is generated by~$F_1$,
and put $\tV=F_1\tA$.
 Then we have $R\subset\tV\subset\tA$ and the ring $\tA$ is generated
by its subgroup~$\tV$.

 Consider the associated graded ring $\gr^F\tA=\bigoplus_{n=0}^\infty
\gr^F_n\tA$, where $\gr^F_n\tA=F_n\tA/F_{n-1}\tA$.
 We will say that the filtered ring $\tA$ is \emph{nonhomogeneous
quadratic} if the ring $A=\gr^F\tA$ is quadratic
(in the sense of Section~\ref{quadratic-duality-secn}).

 Let $(\prtA,F)$ and $(\sectA,F)$ be two nonhomogeneous quadratic rings
with the same degree-zero filtration component
$F_0\sectA=R=F_0\prtA$.
 A \emph{morphism} of nonhomogeneous quadratic rings $f\:\prtA\rarrow
\sectA$ is a morphism of filtered rings (that is, a ring homomorphism
such that $f(F_n\prtA)\subset F_n\sectA$ for all $n\ge0$) forming
a commutative triangle diagram with the subring inclusions
$R\rarrow\prtA$ and $R\rarrow\sectA$.

\begin{lem}
 The above construction assigning the subbimodule of generators
$\tV=F_1\tA$ to a filtration $F$ on a ring $\tA$ defines a fully
faithful functor from the category of nonhomogeneous quadratic rings
to the category of weak nonhomogeneous quadratic rings (over any fixed
base ring~$R$).
 In particular, any nonhomogeneous quadratic ring is weak
nonhomogeneous quadratic (so our terminology is consistent).
\end{lem}

\begin{proof}
 We will only prove the second assertion.
 Let $(A,F)$ be a nonhomogeneous quadratic ring.
 We have to show that the ideal $\tJ_{\tA}\subset T_R(\tV)$ is generated
by $\tI_{\tA}=F_2T_R(\tV)\cap\tJ_{\tA}$.
 Indeed, denote the ideal generated by $\tI_{\tA}$ by
$\tJ'\subset T_R(\tV)$, and consider the ring $\tA'=T_R(\tV)/\tJ'$.
 Then $\tJ'\subset\tJ_{\tA}$, so there is a unique surjective ring
homomorphism $\tA'\rarrow\tA$ forming a commutative triangle diagram
with the surjective ring homomorphisms $T_R(\tV)\rarrow\tA'$ and
$T_R(\tV)\rarrow\tA$.

 For every $n\ge0$, denote by $F_n\tA'\subset\tA'$ the image of
the subgroup $F_nT_R(\tV)\subset T_R(\tV)$ under the ring
homomorphism $T_R(\tV)\rarrow\tA'$.
 Then $F$ is an exhastive multiplicative filtration on the ring~$\tA'$.
 Furthermore, the image of $F_n\tA'$ under the ring homomorphism
$\tA'\rarrow\tA$ coincides with $F_n\tA$.
 It is clear that the maps $F_n\tA'\rarrow F_n\tA$ are isomorphisms for
$n=0$ and~$1$.
 Moreover, we have $\tI_{\tA}\subset F_2T_R(\tV)\cap\tJ'$ by
construction, hence $F_2T_R(\tV)\cap\tJ'=F_2T_R(\tV)\cap\tJ_{\tA}$.
 Therefore, the map $F_2\tA'\rarrow F_2\tA$ is an isomorphism, too.

 We need to show that $\tJ'=\tJ_{\tA}$; equivalently, this means that
the ring homomorphism $\tA'\rarrow\tA$ is an isomorphism.
 It suffices to check that the induced homomorphism of graded rings
$\gr^F\tA'\rarrow\gr^F\tA$ is an isomorphism.
 Set $A'=\gr^F\tA'$ and $A=\gr^F\tA$.
 The graded ring $A'$ is generated by $A'_1$ over $A'_0=R$, since
the filtration $F$ on the ring $\tA'$ is generated by $F_1$ by
construction.
 The graded ring $A$ is quadratic by assumption.
 The graded rings homomorphism $A'\rarrow A$ is an isomorphism in
the degrees~$0$, $1$, and~$2$, as we have shown.
 It remains to apply the next lemma.
\end{proof}

\begin{lem} \label{quadratic-ring-iso}
 Let $f\:A'\rarrow A$ be a homomorphism of nonnegatively graded rings
such that the maps $f_n\:A'_n\rarrow A_n$ are isomorphisms for $n=0$,
$1$, and\/~$2$.
 Assume that the graded ring $A'$ is generated by $A'_1$ over $A'_0$,
while the graded ring $A$ is quadratic.
 Then the map~$f$ is an isomorphism of graded rings.
\end{lem}

\begin{proof}
 Set $A'_0=R=A_0$ and $A'_1=V=A_1$.
 Then there is a unique graded ring homomorphism $\pi_{A'}\:T_R(V)
\rarrow A'$ acting by the chosen isomorphisms $R\simeq A_0$ and
$V\simeq A'_1$ on the components of degrees~$0$ and~$1$, and
a similar unique graded ring homomorphism $\pi_A\:T_R(V)\rarrow A$.
 The triangle diagram of ring homomorphisms $T_R(V)\rarrow A'
\rarrow A$ is commutative, i.~e., $\pi_A=f\pi_{A'}$.
 Furthermore, both the graded rings $A$ and $A'$ are generated by
their degree-one components (over their degree-zero components) by
assumption, hence both the maps~$\pi_{A'}$ and $\pi_A$ are surjective.

 Denote by $J_{A'}$ and $J_A\subset T_R(V)$ the kernels of the graded
ring homomorphisms~$\pi_{A'}$ and~$\pi_A$.
 Then we have $J_{A'}\subset J_A$ and $J_{A',2}=J_{A,2}$, since
the map $f_2\:A'_2\rarrow A_2$ is an isomorphism by assumption.
 The algebra $A$ is quadratic, so the ideal $J_A$ is generated by
$I_A=J_{A,2}$.
 It follows that $J_{A'}=J_A$, hence $f_n$~is an isomorphism
for all~$n$.
\end{proof}

\begin{rem} \label{weak-Koszul-strong-remark}
 We will see below in Section~\ref{pbw-theorem-subsecn} that under
the left finitely projective Koszulity assumption the classes of weak
nonhomogeneous quadratic rings and nonhomogeneous quadratic rings
coincide.
 Specifically, if $\tA$ is a weak nonhomogeneous quadratic ring
such that the quadratic graded ring $\q\,\gr^F\tA$ is left finitely
projective Koszul (where the filtration $F$ on $\tA$ is generated
by $F_1\tA=\tV$ over $F_0\tA=R$, as above), then the graded ring
$\gr^F\tA$ is quadratic (so the filtered ring $\tA$ is
nonhomogeneous quadratic).
\end{rem}

\subsection{Curved DG-rings} \label{curved-dg-rings-subsecn}
 A \emph{CDG\+ring} (\emph{curved differential graded ring}) $B=(B,d,h)$
is a graded associative ring $B=\bigoplus_{n\in\boZ}B^n$ endowed with
a sequence of additive maps $d_n\:B^n\rarrow B^{n+1}$, \ $n\in\boZ$, and
an element $h\in B^2$ satisfying the following conditions:
\begin{enumerate}
\renewcommand{\theenumi}{\roman{enumi}}
\item $d$ is an odd derivation of $B$, that is $d(bc)=d(b)c+
(-1)^{|b|}bd(c)$ for all $b\in B^{|b|}$ and $c\in B^{|c|}$, \
$|b|$, $|c|\in\boZ$;
\item $d^2(b)=[h,b]$ for all $b\in B$ (where $[h,b]=hb-bh$ is
the commutator);
\item $d(h)=0$.
\end{enumerate}

 In the context of the present paper, all CDG\+rings will be
nonnegatively graded, that is $B=\bigoplus_{n=0}^\infty B^n$.
 We denote the grading of $B$ by upper indices, because
the differential~$d$ has degree~$1$.

 Let $\prB=(\prB,d',h')$ and $\secB=(\secB,d'',h'')$ be two
CDG\+rings.
 A \emph{morphism} of CDG\+rings $\secB\rarrow{}\prB$ is a pair
$(f,a)$ consisting of a morphism of graded rings $f\:\secB\rarrow
\prB$ and an element $a\in\prB^1$ such that
\begin{enumerate}
\renewcommand{\theenumi}{\roman{enumi}}
\setcounter{enumi}{3}
\item $f(d''(b))=d'(f(b))+[a,f(b)]$ for all $b\in\secB^{|b|}$
(where $[x,y]=xy-(-1)^{|x||y|}yx$, \ $x\in\prB^{|x|}$,
\,$y\in{}\prB^{|y|}$, is the graded commutator); and
\item $f(h'')=h'+d'(a)+a^2$.
\end{enumerate}

 The composition of morphisms is defined by the formula
$(f,a)(g,b)=(fg,a+f(b))$.
 The identity morphism is the morphism $(\id,0)$.
 These rules define the \emph{category of CDG\+rings}.
 We will denote the category of nonnegatively graded CDG\+rings
$(B,d,h)$ with the fixed degree-zero component $B^0=R$ by
$R\rings_\cdg$.
 Morphisms $\secB\rarrow\prB$ in $R\rings_\cdg$ are CDG\+ring
morphisms $(f,a)\:\secB\rarrow\prB$ such that the graded
ring homomorphism $f\:\secB\rarrow\prB$ forms a commutative
triangle diagram with the fixed isomorphisms $R\simeq\secB^0$
and $R\simeq\prB^0$.

 The element $h\in B^2$ is called the \emph{curvature element}.
 The element $a\in\prB^1$ is called the \emph{change-of-connection
element}.

 For any CDG\+ring $B=(B,d,h)$ and any element $a\in B^1$, the triple
$\prB=(B,d',h')$ with $d'=d+[a,{-}]$ and $h'=h+d(a)+a^2$ is also
a CDG\+ring.
 The CDG\+rings $B$ and $\prB$ are connected by the isomorphism
$(\id,a)\:\prB\rarrow B$.
 Such isomorphisms will be called \emph{change-of-connection
isomorphisms}, while CDG\+ring morphisms of the form $(f,0)$ will be
called \emph{strict morphisms}.
 Any morphism of CDG\+rings $(f,a)\:\secB=(\secB,d'',h'')\rarrow
\prB=(\prB,d',h')$ decomposes uniquely into a strict morphism followed
by a change-of-connection isomorphism, $(f,a)=(\id,a)(f,0)$.

 Furthermore, one can define the \emph{$2$\+category of CDG\+rings}
as follows.
 Let $(f,a)$ and $(g,b)\:\secB=(\secB,d'',h'')\rarrow
\prB=(\prB,d',h')$ be two CDG\+ring morphisms with the same
domain and codomain.
 A \emph{$2$\+morphism} $(f,a)\overset z\rarrow(g,b)$ is an invertible 
element $z\in\prB^0$ satisfying the equations
\begin{enumerate}
\renewcommand{\theenumi}{\roman{enumi}}
\setcounter{enumi}{5}
\item $g(c)=zf(c)z^{-1}$ for all $c\in\secB$; and
\item $b=zaz^{-1}-d'(z)z^{-1}$.
\end{enumerate}
 The element $z\in\prB^0$ is called the \emph{gauge transformation
element}.

 The vertical composition of two $2$\+morphisms $(f',a')\overset w
\rarrow(f'',a'') \overset z\rarrow(f''',a''')$ is the $2$\+morphism
$(f',a')\overset{zw}\rarrow(f''',a''')$.
 The identity $2$\+morphism is the $2$\+morphism $(f,a)\overset1\rarrow
(f,a)$.
 The horizontal composition of two $2$\+morphisms $(g',b')\overset w
\rarrow(g'',b'')\:(C,d_C,h_C)\rarrow(B,d_B,h_B)$ and
$(f',a')\overset z\rarrow(f'',a'')\:(B,d_B,h_B)\rarrow(A,d_A,h_A)$ is
the $2$\+morphism $(f'g',\,a'+f'(b'))\overset{z\circ w}\lrarrow
(f''g'',\,a''+f''(b''))\:(C,d_C,h_C)\allowbreak\rarrow(A,d_A,h_A)$ with
the element $z\circ w=zf'(w)=f''(w)z\in A^0$.

 All the $2$\+morphisms of CDG\+rings are invertible.
 If $(f,a)\:\secB\rarrow\prB$ is a morphism of CDG\+rings and
$z\in\prB^0$ is an invertible element, then the pair $(g,b)$
defined by the formulas (vi\+-vii) is also a morphism of CDG\+rings
$(g,b)\:\secB\rarrow\prB$.
 The morphisms $(f,a)$ and $(g,b)$ are connected by
the $2$\+isomorphism $(f,a)\overset z\rarrow(g,b)$.

 The $2$\+category $\Rings_{\cdg2}$ of nonnegatively graded CDG\+rings
is defined as the following subcategory of the $2$\+category of
CDG\+rings (cf.\ Remark~\ref{degree-zero-iso-remark}).
 The objects of $\Rings_{\cdg2}$ are nonnegatively graded CDG\+rings
$(B,d,h)$, \ $B=\bigoplus_{n=0}^\infty B^n$.
 Morphisms $(\secB,d'',h'')\rarrow(\prB,d',h')$ in $\Rings_{\cdg2}$
are morphisms of CDG\+rings $(f,a)\:(\secB,d'',h'')\rarrow(\prB,d',h')$
such that the \emph{degree-zero component} $f_0\:\secB^0\rarrow
\prB^0$ of the graded ring homomorphism~$f$ is an isomorphism
$\secB^0\simeq\prB^0$.
 \,$2$\+morphisms $(f,a)\overset z\rarrow(g,b)$ in $\Rings_{\cdg2}$
between morphisms $(f,a)$ and $(g,b)$ belonging to $\Rings_{\cdg2}$ are
arbitrary $2$\+morphisms from $(f,a)$ to $(g,b)$ in the $2$\+category
of CDG\+rings.

\subsection{Self-consistency equations} \label{self-consistency-subsecn}
 Let $\tA$ be a weak nonhomogeneous quadratic ring over its subring
$R\subset\tA$ with the $R$\+$R$\+subbimodule of generators
$\tV\subset\tA$.
 Consider the related increasing filtration $F$ on the ring $\tA$, as
constructed in Section~\ref{nonhomogeneous-quadratic-subsecn}, and
let $A=\gr^F\tA$ be the associated graded ring.
 Then $A$ is a nonnegatively graded ring generated by its degree-one
component $A_1=V=\tV/R$ over the degree-zero component $A_0=R$.

 Denote by $I\subset V\ot_RV$ the kernel of the multiplication map
$A_1\ot_RA_1\rarrow A_2$, and consider the quadratic graded ring
$\q A=T_R(V)/(I)$.
 Then we have a natural (adjunction) homomorphism of graded rings
$\q A\rarrow A$, which is an isomorphism in the degrees $n=0$, $1$,
and~$2$, and a surjective map in the degrees $n\ge3$.
 We will say that a weak nonhomogeneous quadratic ring $\tA$ is
\emph{$3$\+left finitely projective} if such is the quadratic
graded ring~$\q A=\q\,\gr^F\tA$.

 Let $\tA$ be a weak nonhomogeneous quadratic ring.
 In the rest of this section, we will assume that the $R$\+$R$\+bimodule
$V=\tV/R=A_1$ is projective as a left $R$\+module and
the $R$\+$R$\+bimodule $F_2\tA/F_1\tA=A_2$ is flat as
a left $R$\+module.

 Then, in particular, in the short exact sequence of $R$\+$R$\+bimodules
$0\rarrow R\rarrow\tV\rarrow V\rarrow0$ all the bimodules are
projective as left $R$\+modules.
 Therefore, this sequence splits as a short exact sequence of left
$R$\+modules, and we can choose a splitting $V\rarrow\tV$.
 Let $V'\subset\tV$ be the image of $V$ under such a splitting map;
so $V'$ is a left (but not a right) $R$\+submodule in
the $R$\+$R$\+bimodule $\tV\subset\tA$ such that $\tV=R\oplus V'$
as a left $R$\+module.
 We will call $V'$ a \emph{submodule of strict generators} of~$\tA$.

 In the rest of this section, we will identify $V$ with~$V'$.
 Let $v\in V$ and $r\in R$ be two elements.
 Let $rv$ and $vr\in V$ denote the elements obtained by applying to
$v\in V$ the left and right action of the element $r\in R$ in
the $R$\+$R$\+bimodule $V$, and let $r*v$ and $v*r\in\tV$ denote
the products of the elements $v\in V\simeq V' \subset\tV\subset\tA$
and $r\in R\subset\tA$ in the ring~$\tA$.
 Then we have
\begin{equation} \label{map-q}
 r*v=rv \quad\text{and}\quad v*r = vr+q(v,r)\,\in\,\tV,
\end{equation}
where $q(v,r)\in R$ is a certain uniquely defined element.
 The map $q\:V\times R\rarrow R$ is obviously biadditive, so it can be
(uniquely) extended to an abelian group homomorphism $q\:V\ot_\boZ R
\rarrow R$.
 Since $R$ is a subring in $\tA$, we also have $r*s=rs$ (where
the left-hand side denotes the product in $\tA$ and the right-hand side
is the product in~$R$) for any pair of elements $r$, $s\in R$.

 Furthermore, let $i\in I$ be an element.
 Consider the natural surjective map $V\ot_\boZ V\rarrow V\ot_RV$
from the tensor product of two copies of $V$ over the ring of integers
$\boZ$ to their tensor product over $R$, and denote by $\hI\subset
V\ot_\boZ V$ the full preimage of the subbimodule $I\subset V\ot_RV$
under the map $V\ot_\boZ V\rarrow V\ot_RV$.
 Let $\hi\in\hI$ denote some preimage of the element $i\in I$ under
the natural surjective map $\hI\rarrow I$.

 The element $i\in I\subset V\ot_RV$ can be presented as a finite sum of
decomposable tensors, $i=\sum_\alpha i_{1,\alpha}\ot_R i_{2,\alpha}$,
\ $i_{1,\alpha}$, $i_{2,\alpha}\in V$.
 We will suppress the notation for the sum over~$\alpha$, and write
simply $i=i_1\ot i_2$.
 Similarly, the element $\hi\in\hI\subset V\ot_\boZ V$ can be presented
as a finite sum $\hi=\sum_\alpha\hi_{1,\alpha}\ot_\boZ\hi_{2,\alpha}$,
where $\hi_{1,\alpha}$, $\hi_{2,\alpha}\in V$.
 Once again, we will omit the notation for the sum over~$\alpha$, and
write $\hi=\hi_1\ot\hi_2$.
 All our formulas will be biadditive in $\hi_1$ and $\hi_2$ (or, as
may be the case, appropriately $R$\+bilinear in $i_1$ and~$i_2$), so
such simplification of the notation will be harmless.

 For any two elements $u$ and $v\in V$, we identify $u$ and~$v$ with
their images under the embedding $V\simeq V'\hookrightarrow\tV
\subset\tA$ and consider their product $u*v\in\tA$ in the ring~$\tA$.
 The assignment of the element $u*v$ to a pair of elements $u$ and~$v$
can be uniquely extended to well-defined homomorphism of abelian groups
(or left $R$\+modules) $V\ot_\boZ V\rarrow\tA$, which does \emph{not},
generally speaking, factorize through $V\ot_RV$.

 In particular, for any element $i=i_1\ot i_2\in I$ and any its preimage
$\hi=\hi_1\ot\hi_2\in\hI$, we consider the element $\hi_1*\hi_2\in\tA$.
 We have $\hi_1$, $\hi_2\in V'\subset\tV=F_1\tA$, so
$\hi_1*\hi_2\in F_2\tA$.
 Moreover, the image of the element $\hi_1*\hi_2$ under the natural
surjection $F_2\tA\rarrow F_2\tA/F_1\tA$ is equal to the product
$i_1i_2$ computed in the associated graded ring $A=\gr^F\tA$.
 Now we have $i_1i_2=0$ in $F_2\tA/F_1\tA=A_2$, since the composition
$I\rarrow V\ot_RV\rarrow A_2$ vanishes by construction.
 Thus $\hi_1*\hi_2\in F_1\tA=\tV\subset\tA$.
 Therefore, there exist uniquely defined elements $p(\hi)\in V$ and
$h(\hi)\in R$ such that
\begin{equation} \label{maps-p-and-h}
 \hi_1*\hi_2=p(\hi)-h(\hi)\,\in\,\tV\subset\tA,
\end{equation}
where the element $p(\hi)\in V$ is identified with its image under
the splitting $V\simeq V'\hookrightarrow\tV$ and the element
$h(\hi)\in R$ is identified with its image under the subring inclusion
$R\hookrightarrow\tV\subset\tA$.

 Furthermore, denote by $I^{(3)}$ the intersection
$I\ot_RV\cap V\ot_RI\subset V\ot_RV\ot_RV$.
 Notice that, since the left $R$\+modules $V$ and $A_2=(V\ot_RV)/I$ are
flat by assumption, the tensor products $I\ot_RV$ and $V\ot_RI$ are
subbimodules in the triple tensor product $V\ot_RV\ot_RV$ (see the proofs
of Propositions~\ref{3-fin-proj-quadratic-duality}
and~\ref{diagonal-Tor}).

 We are also interested in ``the intersection
$\hI\ot_\boZ V\cap V\ot_\boZ\hI\subset V\ot_\boZ V\ot_\boZ V$'',
but here we cannot claim that $\hI\ot_\boZ V$ and $V\ot_\boZ\hI$
are subgroups in $V\ot_\boZ V\ot_\boZ V$.
 So the intersection as such is not well-defined.
 Instead, we have two abelian group homomorphisms $\hI\ot_\boZ V\rarrow
V\ot_\boZ V\ot_\boZ V$ and $V\ot_\boZ\hI\rarrow V\ot_\boZ V\ot_\boZ V$
induced by the inclusion $\hI\rarrow V\ot_\boZ V$.
 We denote by $\hI^{(3)}$ the fibered product of the abelian groups
$\hI\ot_\boZ V$ and $V\ot_\boZ\hI$ over the group
$V\ot_\boZ V\ot_\boZ V$. 

\begin{lem}
 The natural surjections $\hI\ot_\boZ V\rarrow I\ot_RV$, \
$V\ot_\boZ\hI\rarrow V\ot_RI$, and $V\ot_\boZ V\ot_\boZ V\rarrow
V\ot_RV\ot_RV$ induce a surjective map $\hI^{(3)}\rarrow I^{(3)}$.
\end{lem}

\begin{proof}
 Denote by $\hI\ovot_\boZ V$ and $V\ovot_\boZ\hI\subset
V\ot_\boZ V\ot_\boZ V$ the images of the maps
$\hI\ot_\boZ V\rarrow V\ot_\boZ V\ot_\boZ V$ and
$V\ot_\boZ\hI\rarrow V\ot_\boZ V\ot_\boZ V$, and let
$\overline I^{(3)}\subset V\ot_\boZ V\ot_\boZ V$ stand for
the intersection of $\hI\ovot_\boZ V$ and $V\ovot_\boZ\hI$.
 Then there is a natural surjective map $\hI^{(3)}\rarrow
\overline I^{(3)}$.
 The surjection $V\ot_\boZ V\ot_\boZ V\rarrow V\ot_R V\ot_R V$
restricts to a map $\overline I^{(3)}\rarrow I^{(3)}$, and
the map $\hI^{(3)}\rarrow I^{(3)}$ is equal to the composition
$\hI^{(3)}\rarrow\overline I^{(3)}\rarrow I^{(3)}$.
 It remains to prove that the map $\overline I^{(3)}\rarrow
I^{(3)}$ is surjective.

 We have a short exact sequence of abelian groups $0\rarrow\hI
\rarrow V\ot_\boZ V\rarrow A_2\rarrow0$.
 Hence the subgroups $X_1=\hI\ovot_\boZ V$ and $X_2=V\ovot_\boZ\hI
\subset V\ot_\boZ V\ot_\boZ V$ are the kernels of the surjective maps
$V\ot_\boZ V\ot_\boZ V\rarrow A_2\ot_\boZ V$ and
$V\ot_\boZ V\ot_\boZ V\rarrow V\ot_\boZ A_2$, respectively.
 Denote by $Y_1$ and $Y_2\subset V\ot_\boZ V\ot_\boZ V$ the kernels of
the natural surjective maps $V\ot_\boZ V\ot_\boZ V\rarrow
V\ot_RV\ot_\boZ V$ and $V\ot_\boZ V\ot_\boZ V\rarrow
V\ot_\boZ V\ot_R V$.
 Then we have $Y_1\subset X_1$ and $Y_2\subset X_2\subset
V\ot_\boZ V\ot_\boZ V$.

 Furthermore, $Y_1+Y_2$ is the kernel of the natural surjective
map $V\ot_\boZ V\ot_\boZ V\rarrow V\ot_RV\ot_RV$.
 Similarly, $X_1+Y_2$ is the kernel of the surjective map
$V\ot_\boZ V\ot_\boZ V\rarrow A_2\ot_RV$, and $Y_1+X_2$ is
the kernel of the surjective map $V\ot_\boZ V\ot_\boZ V\rarrow
V\ot_R A_2$.
 It follows that the subbimodule $I\ot_R V\subset V\ot_RV\ot_RV$
is the image of $X_1\subset V\ot_\boZ V\ot_\boZ V$ under the map
$V\ot_\boZ V\ot_\boZ V\rarrow V\ot_RV\ot_RV$, and the subbimodule
$V\ot_R I\subset V\ot_RV\ot_RV$ is the image of
$X_2\subset V\ot_\boZ V\ot_\boZ V$ under the same map.

 The assertion that the map $\overline I^{(3)}\rarrow I^{(3)}$ is
surjective is now expressed by the distributivity equation
$$
 X_1\cap X_2+(Y_1+Y_2)=(X_1+Y_2)\cap(Y_1+X_2)
$$
on subgroups in the abelian group $W=V\ot_\boZ V\ot_\boZ V$.
 What we have here is two filtrations $0\subset Y_1\subset X_1
\subset W$ and $0\subset Y_2\subset X_2\subset W$ of an abelian
group~$W$.
 It remains to observe that \emph{any two filtrations of an abelian
category object generate a distributive lattice of its subobjects},
which is a particular case of~\cite[Theorem~5]{Jon}
or~\cite[Corollary~6.4 in Chapter~1]{PP}.
\end{proof}

 Let $j\in I^{(3)}$ be an element and $\hj\in\hI^{(3)}$ be one of its
preimages.
 The element $j\in I^{(3)}\subset V^{\ot_R\,3}$ can be presented as
a finite sum of decomposable tensors,
$j=\sum_\alpha j_{1,\alpha}\ot_R j_{2,\alpha}\ot_R j_{3,\alpha}$, \
$j_{1,\alpha}$, $j_{2,\alpha}$, $j_{3,\alpha}\in V$.
 As above, we will suppress the notation for the sum over~$\alpha$,
and write simply $j=j_1\ot j_2\ot j_3$.
 Moreover, the image of the element $\hj\in\hI^{(3)}$ under the natural
map $\hI^{(3)}\rarrow V^{\ot_\boZ\,3}$ can be presented as a finite sum
$\sum_\alpha\hj_{1,\alpha}\ot_\boZ\hj_{2,\alpha}\ot_\boZ\hj_{3,\alpha}$,
where $\hj_{1,\alpha}$, $\hj_{2,\alpha}$, $\hj_{3,\alpha}\in V$.
 We will write simply $\hj=\hj_1\ot\hj_2\ot\hj_3$, omitting
the notation for the sum over~$\alpha$ and ignoring the distinction
between an element of $\hI^{(3)}$ and its image in $V^{\ot_\boZ\,3}$
in this notation.

 For any three elements $u$, $v$, and $w\in V$, we will identify~$u$,
$v$, and~$w$ with their images under the embedding $V\simeq V'
\hookrightarrow\tV\subset\tA$ and consider the triple product
$(u*v)*w=u*v*w=u*(v*w)\in\tA$ in the ring~$\tA$.
 In particular, for any element $\hj\in\hI^{(3)}$, the element
$\hj_1*\hj_2*\hj_3\in\tA$ is well-defined.

\begin{prop} \label{self-consistency-equations-prop}
 Let $\tA$ be a weak nonhomogeneous quadratic ring such that
the left $R$\+module $V=A_1$ is projective and the left $R$\+module
$A_2$ is flat.
 Suppose that a left $R$\+linear splitting $V\simeq V'\hookrightarrow
\tV$ of the surjecitve $R$\+$R$\+bimodule map $\tV\rarrow\tV/R=V$
has been chosen.
 Then the maps $q\:V\times R\rarrow R$, \ $p\:\hI\rarrow V$, and
$h\:\hI\rarrow R$ defined above in~(\ref{map-q}\+-\ref{maps-p-and-h})
satisfy the following self-consistency equations:
\begin{enumerate}
\renewcommand{\theenumi}{\alph{enumi}}
\item $q(rv,s)=rq(v,s)$ for all $r$, $s\in R$ and $v\in V$;
\item $q(v,rs)=q(vr,s)+q(v,r)s$ for all $r$, $s\in R$ and
$v\in V$;
\item $p(r\hi_1\ot\hi_2)=rp(\hi_1\ot\hi_2)$ for all $r\in R$
and $\hi\in\hI$;
\item $h(r\hi_1\ot\hi_2)=rh(\hi_1\ot\hi_2)$ for all $r\in R$
and $\hi\in\hI$;
\item $p(u\ot rv-ur\ot v)=q(u,r)v$ for all $r\in R$ and
$u$, $v\in V$;
\item $h(u\ot rv-ur\ot v)=0$ for all $r\in R$ and
$u$, $v\in V$;
\item $p(\hi_1\ot\hi_2r)=p(\hi_1\ot\hi_2)r-\hi_1q(\hi_2,r)$
for all $r\in R$ and $\hi=\hi_1\ot\hi_2\in\hI$;
\item $h(\hi_1\ot\hi_2r)=h(\hi_1\ot\hi_2)r-q(p(\hi_1\ot\hi_2),r)
+q(\hi_1,q(\hi_2,r))$ for all $r\in R$ and $\hi\in\hI$;
\item $p(\hj_1\ot\hj_2)\ot\hj_3-\hj_1\ot p(\hj_2\ot\hj_3)\,\in\,
\hI\subset V\ot_\boZ V$ for all $\hj=\hj_1\ot\hj_2\ot\hj_3\in\hI^{(3)}$;
\item $p(p(\hj_1\ot\hj_2)\ot\hj_3-\hj_1\ot p(\hj_2\ot\hj_3))
=h(\hj_1\ot\hj_2)\hj_3-\hj_1h(\hj_2\ot\hj_3)$ for all $\hj\in\hI^{(3)}$;
\item $h(p(\hj_1\ot\hj_2)\ot\hj_3-\hj_1\ot p(\hj_2\ot\hj_3))
=q(\hj_1,h(\hj_2\ot\hj_3))$ for all $\hj\in\hI^{(3)}$.
\end{enumerate}
\end{prop}

\begin{proof}
 All these equations follow, in one way or another, from
the associativity of multiplication in the ring~$\tA$.
 The specific computations proving each of the formulas are presented
below one by one.
 
 Part~(a): compare $r*v*s=(r*v)*s=(rv)*s=rvs+q(rv,s)$ with $r*v*s=
r*(v*s)=r*(vs+q(v,s))=rvs+rq(v,s)$.

 Part~(b): compare $v*r*s=v*(r*s)=v*(rs)=vrs+q(v,rs)$ with
$v*r*s=(v*r)*s=(vr+q(v,r))*s=vrs+q(vr,s)+q(v,r)s$.

 Parts~(c) and~(d): notice first of all that $\sum_\alpha
r\hi_{1,\alpha}\ot\hi_{2,\alpha}\in\hI$ whenever $r\in R$ and
$\sum_\alpha\hi_{1,\alpha}\ot\hi_{2,\alpha}\in\hI$.
 So the left-hand sides of both the equations are well-defined.
 To deduce the equations, compare
$$
 r*\hi_1*\hi_2=(r*\hi_1)*\hi_2=(r\hi_1)*\hi_2=
 p(r\hi_1\ot\hi_2)-h(r\hi_1\ot\hi_2)
$$
with
$$
 r*\hi_1*\hi_2=r*(\hi_1*\hi_2)=r*(p(\hi_1\ot\hi_2)-
 h(\hi_1\ot\hi_2))=rp(\hi_1\ot\hi_2)-rh(\hi_1\ot\hi_2)
$$
and equate the terms belonging to $V\simeq V'\subset\tV$ separately
and the terms belonging to $R\subset\tV$ separately.

 Parts~(e) and~(f): first of all, one has $u\ot rv-ur\ot v\in\hI$
for all $r\in R$ and $u$, $v\in V$.
 So the left-hand sides of the equations are well-defined.
 Furthermore,
\begin{multline*}
 0=u*(r*v)-(u*r)*v=u*(rv)-(ur)*v-q(u,r)*v \\
 =p(u\ot rv-ur\ot v)-h(u\ot rv-ur\ot v)-q(u,r)v
\end{multline*}
and it remains to equate separately the terms belonging to $V'$
and to~$R$.

 Parts~(g) and~(h): we have $\sum_\alpha
\hi_{1,\alpha}\ot\hi_{2,\alpha}r\in\hI$ whenever $r\in R$ and
$\sum_\alpha\hi_{1,\alpha}\ot\hi_{2,\alpha}\in\hI$;
so the left-hand sides of both the equations are well-defined.
 Now compare
\begin{multline*}
 \hi_1*\hi_2*r=(\hi_1*\hi_2)*r=(p(\hi_1\ot\hi_2)-h(\hi_1\ot\hi_2))*r
 \\ =p(\hi_1\ot\hi_2)r+q(p(\hi_1\ot\hi_2),r)-h(\hi_1\ot\hi_2)r
\end{multline*}
with
\begin{multline*}
 \hi_1*\hi_2*r=\hi_1*(\hi_2*r)=\hi_1*(\hi_2r+q(\hi_2,r)) \\
 = p(\hi_1\ot\hi_2r)-h(\hi_2\ot\hi_2r)+
 \hi_1q(\hi_2,r)+q(\hi_1,q(\hi_2,r))
\end{multline*}
and equate separately the terms belonging to $V'$ and to~$R$.

 Parts~(i\+-k): given an element $j\in\hI^{(3)}$, we have
$(\hj_1*\hj_2)*\hj_3=\hj_1*\hj_2*\hj_3=\hj_1*(\hj_2*\hj_3)$ in
$F_3\tA\subset\tA$.
 By construction, the value of this triple product in $\tA$ only
depends on the image of the element~$j$ in the group
$V\ot_\boZ V\ot_\boZ V$.
 We will compute the value of $(\hj_1*\hj_2)*\hj_3$ in terms
of the image of~$j$ in $\hI\ot_\boZ V$ and the value of
$\hj_1*(\hj_2*\hj_3)$ in terms of the image of~$j$ in
$V\ot_\boZ\hI$, and then equate the two expresssions.

 Specifically, we have
$$
 (\hj_1*\hj_2)*\hj_3=p(\hj_1\ot\hj_2)*\hj_3-h(\hj_1\ot\hj_2)*\hj_3
$$
and
$$
 \hj_1*(\hj_2*\hj_3)=\hj_1*p(\hj_2\ot\hj_3)-\hj_1*h(\hj_2\ot\hj_3),
$$
hence
\begin{equation} \label{main-self-consistency-eqn}
 p(\hj_1\ot\hj_2)*\hj_3-\hj_1*p(\hj_2\ot\hj_3)=
 h(\hj_1\ot\hj_2)*\hj_3-\hj_1*h(\hj_2\ot\hj_3).
\end{equation}

 Now, first of all, the right-hand side
of~\eqref{main-self-consistency-eqn} belongs to $\tV\subset\tA$,
hence so does the left-hand side.
 Both summands in the left-hand side belong to $F_2\tA$.
 So the image of the left-hand side in $A_2=F_2\tA/F_1\tA$ has to
vanish, which means that the expression $p(\hj_1\ot\hj_2)\ot_R\hj_3
-\hj_1\ot_R p(\hj_2\ot\hj_3)$ belongs to $I\subset V\ot_RV$.
 This proves part~(i).

 It remains to compute both sides of~\eqref{main-self-consistency-eqn}
as
\begin{multline*}
 p(\hj_1\ot\hj_2)*\hj_3-\hj_1*p(\hj_2\ot\hj_3) \\ =
 p(p(\hj_1\ot\hj_2)\ot\hj_3-\hj_1\ot p(\hj_2\ot\hj_3))-
 h(p(\hj_1\ot\hj_2)\ot\hj_3-\hj_1\ot p(\hj_2\ot\hj_3))
\end{multline*}
and
$$
 h(\hj_1\ot\hj_2)*\hj_3-\hj_1*h(\hj_2\ot\hj_3)=
 h(\hj_1\ot\hj_2)\hj_3-
 \hj_1h(\hj_2\ot\hj_3)-q(\hj_1,h(\hj_2\ot\hj_3)).
$$
 Comparing and equating separately the terms belonging to $V'$ and
to $R$ produces the desired formulas~(j\+-k).
\end{proof}

\subsection{The CDG-ring corresponding to a nonhomogeneous
quadratic ring} \label{cdg-ring-constructed-subsecn}
 Let $R$ and $S$ be associative rings, $U$ be an $R$\+$S$\+bimodule,
and $U\spcheck=\Hom_R(U,R)$ be the dual $S$\+$R$\+bimodule.
 We will use the notation
$$
 \lan u,f\ran = f(u) \,\in\,R \quad
 \text{for any $u\in U$ and $f\in U\spcheck$}.
$$
 Then the condition that $f\:U\rarrow R$ is a left $R$\+module
homomorphism is expressed by the identity
$$
 \lan ru,f\ran = r\lan u,f\ran \quad
 \text{for all $r\in R$, \ $u\in U$, and $f\in U\spcheck$},
$$
while the construction of the left $S$\+module structure on
$U\spcheck$ is expressed by the identity
$$
 \lan u,sf\ran = \lan us,f\ran \quad
 \text{for all $u\in U$, \ $s\in S$, and $f\in U\spcheck$},
$$
and the construction of the right $R$\+module structure on
$U\spcheck$ is expressed by
$$
 \lan u,fr\ran = \lan u,f\ran r \quad
 \text{for all $u\in U$, \ $f\in U\spcheck$, and $r\in R$}.
$$

 Furthermore, given three rings $R$, \,$S$, and $T$,
an $R$\+$S$\+bimodule $U$, and an $S$\+$T$\+bi\-module $V$,
the construction of the natural homomorphism of $T$\+$R$\+bimodules
$$
 \Hom_S(V,S)\ot_S\Hom_R(U,R)\lrarrow\Hom_R(U\ot_SV,\>R),
$$
from Lemma~\ref{tensor-dual-lemma}(a) can be expressed by
the formula
$$
 \lan u\ot v,\>g\ot f\ran =\lan u\lan v,g\ran,f\ran
 =\lan u,\lan v,g\ran f\ran
$$
for all $u\in U$, \ $v\in V$, \ $g\in\Hom_S(V,S)$,
and $f\in\Hom_R(U,R)$.

\begin{prop} \label{nonhomogeneous-dual-cdg-ring}
 Let $\tA$ be a\/ $3$\+left finitely projective weak nonhomogeneous
quadratic ring over its subring $R\subset\tA$ with
the $R$\+$R$\+bimodule of generators $R\subset\tV\subset\tA$.
 Denote by $B$ the\/ $3$\+right finitely projective quadratic graded
ring quadratic dual to the\/ $3$\+left finitely projective
quadratic graded ring\/ $\q A=\q\,\gr^F\tA$.
 Suppose that a left $R$\+linear splitting $V\simeq V'\hookrightarrow
\tV$ of the surjective $R$\+$R$\+bimodule map $\tV\rarrow\tV/R=V$
has been chosen.
 Then the formulas
\begin{equation} \label{d-zero-defined}
 \lan v,d_0(r)\ran = q(v,r)
\end{equation}
and
\begin{equation} \label{d-one-defined}
 \lan i,d_1(b)\ran =
 \lan p(\hi_1\ot\hi_2),b\ran - q(\hi_1,\lan \hi_2,b\ran)
\end{equation}
for all $r\in R$, \ $v\in V$, \ $i\in I$, and\/ $b\in B^1$, where
the maps~$q$ and~$p$ are given
by~(\ref{map-q}\+-\ref{maps-p-and-h}), specify well-defined abelian
group homomorphisms $d_0\:B^0\rarrow B^1$ and $d_1\:B^1\rarrow B^2$.
 Furthermore, the map $h\:\hI\rarrow R$ descends uniquely to
a well-defined left $R$\+linear map $I\rarrow R$, providing an element
$h\in\Hom_R(I,R)=B^2$.
 The maps $d_0$ and~$d_1$ satisfy the equations
\begin{enumerate}
\renewcommand{\theenumi}{\alph{enumi}}
\item $d_0(rs)=d_0(r)s+rd_0(s)$ for all $r$, $s\in R$;
\item $d_1(rb)=d_0(r)b+rd_1(b)$ for all $r\in R$, \ $b\in B^1$;
\item $d_1(br)=d_1(b)r-bd_0(r)$ for all $r\in R$, \ $b\in B^1$;
\item $d_1(d_0(r))=hr-rh$ for all $r\in R$;
\item $\sum_\alpha d_1(e_{1,\alpha})e_{2,\alpha}-
 \sum_\alpha e_{1,\alpha}d_1(e_{2,\alpha})=0$ in $B^3$ for all
 tensors $e=\sum_\alpha e_{1,\alpha}\ot_R e_{2,\alpha}\in B^1\ot_R B^1$
 such that the image of~$e$ vanishes in~$B^2$.
\end{enumerate}
 The formula
$$
 d_2(e_1e_2)=d_1(e_1)e_2-e_1d_1(e_2) \quad
 \text{for all\/ $e_1$, $e_2\in B^1$}
$$
specifies a well-defined abelian group homomorphism
$d_2\:B^2\rarrow B^3$, which satisfies the equations
\begin{enumerate}
\renewcommand{\theenumi}{\alph{enumi}}
\setcounter{enumi}{5}
\item $d_2(d_1(b))=hb-bh$ for all $b\in B^1$; and
\item $d_2(h)=0$.
\end{enumerate}
 In other words, the maps~$d_0$ and~$d_1$ admit a unique extension to
an odd derivation $d\:B\rarrow B$ of degree\/~$1$, and the triple
$(B,d,h)$ is a CDG\+ring.
\end{prop}

\begin{proof}
 Recall that, by the definition of quadratic duality, we have
$B^0=R$, \ $B^1=\Hom_R(V,R)$, and $B^2=\Hom_R(I,R)$ (where $I$ is
the kernel of the surjective multiplication map $A_1\ot_R A_1
\rarrow A_2\simeq(\q A)_2$; the latter isomorphism holds since
the graded ring $A$ is generated by $A_1$ over $R=A_0$).
 The grading on the ring $B$ was denoted by lower indices in
Sections~\ref{quadratic-duality-secn}\+-\ref{koszulity-secn}, but
we denote it by upper indices here; the convention is $B^n=B_n$ for
all $n\ge0$ (and $B^n=0=B_n$ for $n<0$).

 Firstly we have to check that the maps~$d_0$ and~$d_1$ are well-defined
by the formulas~(\ref{d-zero-defined}\+-\ref{d-one-defined}).
 Concerning~$d_0$, it needs to be checked that $v\longmapsto
\lan v,d_0(r)\ran$ is a left $R$\+linear map $V\rarrow R$ for every
$r\in R$.
 Indeed, we have
$$
 \lan sv,d_0(r)\ran=q(sv,r)=sq(v,r)=s\lan v,d_0(r)\ran
$$
for all $r$, $s\in R$ and $v\in V$ by
Proposition~\ref{self-consistency-equations-prop}(a).

 Concerning~$d_1$, it needs to be checked that, for every element
$b\in B^1$, the map $\hI\rarrow R$ defined by the formula
$\hi\longmapsto\lan p(\hi_1\ot\hi_2),b\ran-q(\hi_1,\lan\hi_2,b\ran)$
descends to a left $R$\+linear map $I\rarrow R$.
 Indeed, for all $u$, $v\in V$ and $r\in R$ we have
\begin{multline*}
 \lan p(u\ot rv-ur\ot v),b\ran - q(u,\lan rv,b\ran)+
 q(ur,\lan v,b\ran) \\ = 
 \lan q(u,r)v,b\ran - q(u,\lan rv,b\ran)+q(ur,\lan v,b\ran) \\
 = q(u,r)\lan v,b\ran - q(u,r\lan v,b\ran)+q(ur,\lan v,b\ran) = 0 
\end{multline*}
by Proposition~\ref{self-consistency-equations-prop}(e) and~(b),
the latter of which is being applied to the elements $u\in V$ and
$r$, $\lan v,b\ran\in R$.
 Since $\hI\rarrow I$ is a surjective map with the kernel spanned,
as an abelian group, by the elements $u\ot rv-ur\ot v$, it follows
that our map $\hI\rarrow R$ descends uniquely to a map
$d_1(b)\:I\rarrow R$.
 To prove that the latter map is left $R$\+linear, we compute
\begin{multline*}
 \lan ri,d_1(b)\ran = \lan p(r\hi_1\ot\hi_2),b\ran -
 q(r\hi_1,\lan\hi_2,b\ran) \\ =\lan rp(\hi_1\ot\hi_2),b\ran -
 rq(\hi_1,\lan\hi_2,b\ran) \\ = r\lan p(\hi_1\ot\hi_2),b\ran
 -rq(\hi_1,\lan\hi_2,b\ran) = r\lan i,d_1(b)\ran
\end{multline*}
using Proposition~\ref{self-consistency-equations-prop}(c) and~(a).

 Similarly, the map $h\:\hI\rarrow R$ descends to a well-defined
map $I\rarrow R$ by
Proposition~\ref{self-consistency-equations-prop}(f), and the latter
map is left $R$\+linear by
Proposition~\ref{self-consistency-equations-prop}(d).

 Now we have to prove the equations~(a\+-g).
 Part~(a): for every element $v\in V$, one has
\begin{multline*}
 \lan v,d_0(rs)\ran = q(v,rs) = q(vr,s)+q(v,r)s =
 \lan vr,d_0(s)\ran + \lan v,d_0(r)\ran s \\ =
 \lan v,rd_0(s)\ran + \lan v,d_0(r)s\ran =
 \lan v,\,rd_0(s)+d_0(r)s\ran
\end{multline*}
by Proposition~\ref{self-consistency-equations-prop}(b).

 Part~(b): for every element $i\in I$ and its preimage $\hi\in\hI$,
one has
\begin{multline*}
\lan i,d_1(rb)\ran = \lan p(\hi_1\ot\hi_2),rb\ran -
q(\hi_1,\lan\hi_2,rb\ran) = \lan p(\hi_1\ot\hi_2)r,b\ran -
q(\hi_1,\lan\hi_2r,b\ran) \\ = \lan \hi_1q(\hi_2,r),b\ran +
\lan p(\hi_1\ot\hi_2r),b\ran - q(\hi_1,\lan\hi_2r,b\ran) \\ =
\lan\hi_1\lan\hi_2,d_0(r)\ran,b\ran + \lan ir,d_1(b)\ran =
\lan i,d_0(r)b\ran + \lan i,rd_1(b)\ran
\end{multline*}
by Proposition~\ref{self-consistency-equations-prop}(g).

 Part~(c): for every element $i\in I$ and its preimage $\hi\in\hI$,
one has
\begin{multline*}
\lan i,d_1(br)\ran = \lan p(\hi_1\ot\hi_2),br\ran -
q(\hi_1,\lan\hi_2,br\ran) = \lan p(\hi_1\ot\hi_2),b\ran r
-q(\hi_1,\lan\hi_2,b\ran r) \\ = \lan p(\hi_1\ot\hi_2),b\ran r
-q(\hi_1,\lan\hi_2,b\ran)r - q(\hi_1\lan\hi_2,b\ran,r) \\ =
\lan i,d_1(b)\ran r - \lan\hi_1\lan\hi_2,b\ran,d_0(r)\ran =
\lan i, d_1(b)r\ran - \lan i,bd_0(r)\ran
\end{multline*}
by Proposition~\ref{self-consistency-equations-prop}(b)
applied to the elements $\hi_1\in V$ and $\lan\hi_2,b\ran$,
$r\in R$.

 Part~(d): for every element $i\in I$ and its preimage $\hi\in\hI$,
one has
\begin{multline*}
\lan i,d_1(d_0(r))\ran = \lan p(\hi_1\ot\hi_2),d_0(r)\ran -
q(\hi_1,\lan\hi_2,d_0(r)\ran) \\ = q(p(\hi_1\ot\hi_2),r) -
q(\hi_1,q(\hi_2,r)) = h(\hi_1\ot\hi_2)r - h(\hi_1\ot\hi_2r) \\
= \lan i,h\ran r - \lan ir,h\ran = \lan i,hr\ran - \lan i,rh\ran
\end{multline*}
by Proposition~\ref{self-consistency-equations-prop}(h).

 In order to prove parts~(e\+-g), we recall the natural isomorphism
$B^3\simeq\Hom_R(I^{(3)},R)$, which holds for a $3$\+left finitely
projective quadratic graded ring $\q A$ and its quadratic dual
$3$\+right finitely projective quadratic graded ring $B$ according to
the proof of Proposition~\ref{3-fin-proj-quadratic-duality}.
 In view of this isomorphism, in order to verify an equation in
the group $B^3$, it suffices to evaluate it on every element
$j\in I^{(3)}$ and check that the resulting equation in $R$ holds.
{\hbadness=1350\par}

 Furthermore, for any tensor $e=\sum_\alpha e_{1,\alpha}\ot
e_{2,\alpha}=e_1\ot e_2\in B^1\ot_RB^1$, and for every element
$j\in I^{(3)}$ and its preimage $\hj\in\hI^{(3)}$, we compute
\begin{multline} \label{d-two-computed}
 \lan j,\,d_1(e_1)e_2-e_1d_1(e_2)\ran =
 \lan j_1\lan j_2\ot j_3, d_1(e_1)\ran,e_2\ran -
 \lan j_1\ot j_2\lan j_3,e_1\ran,d_1(e_2)\ran \\
 = \lan\hj_1\lan p(\hj_2\ot\hj_3),e_1\ran,e_2\ran -
 \lan\hj_1q(\hj_2,\lan\hj_3,e_1\ran),e_2\ran \\
 -\lan p(\hj_1\ot\hj_2\lan\hj_3,e_1\ran),e_2\ran +
 q(\hj_1,\lan\hj_2\lan\hj_3,e_1\ran,e_2\ran) \\
 = \lan\hj_1\lan p(\hj_2\ot\hj_3),e_1\ran,e_2\ran
 - \lan p(\hj_1\ot\hj_2)\lan\hj_3,e_1\ran,e_2\ran +
 q(\hj_1,\lan\hj_2\lan\hj_3,e_1\ran,e_2\ran) \\ =
 \lan\hj_1\ot p(\hj_2\ot\hj_3)-p(\hj_1\ot\hj_2)\ot\hj_3,\,e_1\ot e_2\ran
 + q(\hj_1,\lan\hj_2\ot\hj_3,e_1\ot e_2\ran)
\end{multline}
by Proposition~\ref{self-consistency-equations-prop}(g) applied to
the tensor $\hi_1\ot\hi_2=\hj_1\ot\hj_2$ and the element
$r=\lan\hj_3,e_1\ran$.

 Now, tensors $e\in B^1\ot_R B^1$ whose image vanishes in $B^2$ form
the $R$\+$R$\+bimodule $I_B\subset B^1\ot_RB^1$ of quadratic relations
in the quadratic graded ring~$B$.
 By Proposition~\ref{2-fin-proj-quadratic-duality}, we have
$I_B=\Hom_R(A_2,R)=\Hom_R((V\ot_R\nobreak V)/I,\>R)$.
 To prove part~(e), it remains to observe that, in view of
Proposition~\ref{self-consistency-equations-prop}(i), both the summands
in the final expression in~\eqref{d-two-computed} involve the pairing
of an element of $I\subset V\ot_RV$ with the element~$e$.
 Therefore, both the summands vanish for $e\in I_B\subset B^1\ot_RB^1$.

 Part~(f): choose an element $e_1\ot e_2=
\sum_\alpha e_{1,\alpha}\ot e_{2,\alpha}\in B^1\ot_RB^1$ whose
image under the multiplication map $B^1\ot_RB^1\rarrow B^2$ is
equal to~$d_1(b)$.
 Then for every element $j\in I^{(3)}$ and its preimage
$\hj\in\hI^{(3)}$ we have
\begin{multline*}
 \lan j, d_2(d_1(b))\ran = \lan j,\,d_1(e_1)e_2-e_1d_1(e_2)\ran \\ =
 \lan\hj_1\ot p(\hj_2\ot\hj_3)-p(\hj_1\ot\hj_2)\ot\hj_3,\,e_1\ot e_2\ran
 + q(\hj_1,\lan\hj_2\ot\hj_3,e_1\ot e_2\ran) \\ =
 \lan\hj_1\ot p(\hj_2\ot\hj_3)-p(\hj_1\ot\hj_2)\ot\hj_3,\,d_1(b)\ran
 + q(\hj_1,\lan\hj_2\ot\hj_3,d_1(b)\ran) \\ =
 \lan p(\hj_1\ot p(\hj_2\ot\hj_3)-p(\hj_1\ot\hj_2)\ot\hj_3), b\ran \\
 - q(\hj_1,\,\lan p(\hj_2\ot\hj_3),b\ran)
 + q(p(\hj_1\ot\hj_2),\lan\hj_3,b\ran) \\
 + q(\hj_1,\,\lan p(\hj_2\ot\hj_3),b\ran)
 - q(\hj_1,q(\hj_2,\lan\hj_3,b\ran)) \\ =
 \lan\hj_1h(\hj_2\ot\hj_3),b\ran - \lan h(\hj_1\ot\hj_2)\hj_3,b\ran
 + q(p(\hj_1\ot\hj_2),\lan\hj_3,b\ran)
 - q(\hj_1,q(\hj_2,\lan\hj_3,b\ran)) \\ =
 \lan\hj_1h(\hj_2\ot\hj_3),b\ran - h(\hj_1\ot\hj_2)\lan\hj_3,b\ran
 + q(p(\hj_1\ot\hj_2),\lan\hj_3,b\ran)
 - q(\hj_1,q(\hj_2,\lan\hj_3,b\ran)) \\ =
 \lan\hj_1h(\hj_2\ot\hj_3),b\ran - h(\hj_1\ot\hj_2\lan\hj_3,b\ran) \\
 = \lan j_1\lan j_2\ot j_3,h\ran, b\ran
 - \lan j_1\ot j_2\lan j_3,b\ran, h\ran =
 \lan j,\,hb-bh\ran
\end{multline*}
by~\eqref{d-two-computed}, \eqref{d-one-defined}, and
Proposition~\ref{self-consistency-equations-prop}(j) and~(h),
the latter of which is being applied to the tensor
$\hi_1\ot\hi_2=\hj_1\ot\hj_2$ and the element $r=\lan\hj_3,b\ran$.

 Part~(g): choose an element $h_1\ot h_2=
\sum_\alpha h_{1,\alpha}\ot h_{2,\alpha}\in B^1\ot_RB^1$ whose
image under the multiplication map $B^1\ot_RB^1\rarrow B^2$ is
equal to~$h$.
 Then for every element $j\in I^{(3)}$ and its preimage
$\hj\in\hI^{(3)}$ we have
\begin{multline*}
 \lan j, d_2(h)\ran = \lan j,\,d_1(h_1)h_2-h_1d_1(h_2)\ran \\ =
 \lan\hj_1\ot p(\hj_2\ot\hj_3)-p(\hj_1\ot\hj_2)\ot\hj_3,\,h_1\ot h_2\ran
 + q(\hj_1,\lan\hj_2\ot\hj_3,h_1\ot h_2\ran) \\ =
 h(\hj_1\ot p(\hj_2\ot\hj_3)-p(\hj_1\ot\hj_2)\ot\hj_3)
 + q(\hj_1,h(\hj_2\ot\hj_3))=0
\end{multline*}
by~\eqref{d-two-computed} and
Proposition~\ref{self-consistency-equations-prop}(k).

 Finally, for any quadratic graded ring $B$, any pair of maps
$d_0\:B^0\rarrow B^1$ and $d_1\:B^1\rarrow B^2$ satisfying~(a\+-c)
and~(e) can be extended to an odd derivation $d\:B\rarrow B$
of degree~$1$ in a unique way.
 This assertion is provable, e.~g., using the next lemma (where
one takes $W=B^0\oplus B^1$).
 The equations~(d) and~(f) imply that $d(d(x))=[h,x]$ for all
$x\in B$, since $B$ is generated by $B^1$ over~$B^0$.
\end{proof}

\begin{lem} \label{odd-derivations-lemma}
\textup{(a)} Let $W=W_{\bar0}\oplus W_{\bar1}$ be
a\/ $\boZ/2\boZ$\+graded abelian group (where\/ $\boZ/2\boZ=
\{\bar0,\bar1\}$ is the group of order~$2$), and let
$T_{\boZ}(W)=\bigoplus_{n=0}^\infty W^{\ot_\boZ\,n}$ denote
the free associative ring spanned by~$W$.
 Endow $T_{\boZ}(W)$ with the\/ $\boZ/2\boZ$\+grading induced by
that of~$W$.
 Let $d_W\:W\rarrow T_{\boZ}(W)$ be an arbitrary odd homomorphism
of\/ $\boZ/2\boZ$\+graded abelian groups.
 Then there exists a unique odd derivation $d_T\:T_{\boZ}(W)\rarrow
T_{\boZ}(W)$ of the free ring $T_{\boZ}(W)$ extending the map~$d_W$
from $W\subset T_{\boZ}(W)$. \par
\textup{(b)} Let $L\subset T_{\boZ}(W)$ be
a\/ $\boZ/2\boZ$\+homogeneous subgroup and
$K\subset T_{\boZ}(W)$ be the two-sided\/ $\boZ/2\boZ$\+homogeneous
ideal generated by $L$ in $T_{\boZ}(W)$; so $K=(L)$.
 Suppose that $d_T(L)\subset K$.
 Then $d_T(K)\subset K$, and $d_T$~descends to a well-defined odd
derivation~$d$ of the quotient ring $T_{\boZ}(W)/K$.
\end{lem}

\begin{proof}
 Part~(a): put $d_T(w_1\ot\dotsb\ot w_n)=d_W(w_1)\ot w_2\ot\dotsb\ot
w_n+(-1)^{|w_1|}w_1\ot d_W(w_2)\ot w_3\ot\dotsb\ot w_n+\dotsb+
(-1)^{|w_1|+\dotsb+|w_{n-1}|}w_1\ot w_2\ot\dotsb\ot w_{n-1}
\ot d_W(w_n)$ for all $w_i\in W_{|w_i|}$, \ $1\le i\le n$, \ $n\ge0$.
 Part~(b) is obvious.
\end{proof}

\subsection{Change of strict generators}
\label{change-of-strict-gens-subsecn}
 Let $\tA$ be a weak nonhomogeneous quadratic ring over a subring
$R\subset\tA$ with the $R$\+$R$\+subbimodule of generators
$\tV\subset\tA$.
 Assume that the $R$\+$R$\+bimodule $V=\tV/R$ is projective as
a left $R$\+module and the $R$\+$R$\+bimodule $F_2\tA/F_1\tA$ is
flat as a left $R$\+module.

 Furthermore, assume that we are given two left $R$\+linear
splittings $V\simeq V''\subset\tV$ and $V\simeq V'\subset\tV$
of the surjective $R$\+$R$\+bimodule morphism $\tV\rarrow V$.
 Given an element $v\in V$, we will denote by $v''\in V''$ and
$v'\in V'$ its images under the two splittings.
 Then $v\longmapsto v''-v'$ is a left $R$\+linear map
$V\rarrow R$, which we will denote by~$a$.
 Conversely, given a left $R$\+linear splitting $V\simeq V'\subset\tV$
and a left $R$\+linear map $a\:V\rarrow R$, one can construct
a second splitting $V\simeq V''\subset\tV$ by the rule
\begin{equation} \label{change-of-strict-generators-eqn}
 V''=\{\,v'+a(v)\mid v\in V\,\}\,\subset\,\tV.
\end{equation}
 Denote the maps $q$, $p$, and~$h$ defined by
the formulas~(\ref{map-q}\+-\ref{maps-p-and-h}) using
the splitting $V'\subset\tV$ by
$$
 q'\:V\times R\rarrow R, \ \ p'\:\hI\rarrow V, \text{ and \ }
 h'\:\hI\rarrow R
$$
and the similar maps constructed using the splitting $V''\subset\tV$ by
$$
 q''\:V\times R\rarrow R, \ \ p''\:\hI\rarrow V, \text{ and \ }
 h''\:\hI\rarrow R.
$$

\begin{prop} \label{change-of-strict-generators-equations-prop}
 The maps~$q''$, $p''$, and~$h''$ can be obtained from the maps~$q'$,
$p'$, and~$h'$ and the map $a\:V\rarrow R$ by the formulas
\begin{enumerate}
\renewcommand{\theenumi}{\alph{enumi}}
\item $q''(v,r)=q'(v,r)+a(v)r-a(vr)$ for all $r\in R$ and $v\in V$;
\item $p''(\hi_1\ot\hi_2)=p'(\hi_1\ot\hi_2)+a(\hi_1)\hi_2+
\hi_1a(\hi_2)$ for all tensors $\hi=\hi_1\ot\hi_2\in\hI$;
\item $h''(\hi_1\ot\hi_2)=h'(\hi_1\ot\hi_2)+a(p'(\hi_1\ot\hi_2))
-q'(\hi_1,a(\hi_2))+a(\hi_1a(\hi_2))$ for all $\hi\in\hI$.
\end{enumerate}
\end{prop}

\begin{proof}
 In our new notation, the formulas~(\ref{map-q}\+-\ref{maps-p-and-h})
take the form
\begin{gather*}
 r*v'=(rv)' \quad\text{and}\quad r*v''=(rv)'', \\
 v'*r=(vr)'+q'(v,r) \quad\text{and}\quad
 v''*r=(vr)''+q''(v,r),
\end{gather*}
and
$$
 \hi_1'*\hi_2'=p'(\hi_1\ot\hi_2)'-h'(\hi_1\ot\hi_2)
 \quad\text{and}\quad
 \hi_1''*\hi_2''=p''(\hi_1\ot\hi_2)''-h''(\hi_1\ot\hi_2)
$$
for all $r\in R$, \ $v\in V$, and $\hi\in\hI$.
 Here the elements $rv$, $vr\in V$ correspond to the elements
$(rv)'$, $(vr)'\in V'$ and $(rv)''$, $(vr)''\in V''$.
 Similarly, $p'(\hi_1\ot\hi_2)$, $p''(\hi_1\ot\hi_2)\in V$, while
$p'(\hi_1\ot\hi_2)'\in V'$ and $p''(\hi_1\ot\hi_2)''\in V''$.

 Part~(a): one has, on the one hand,
$$
 v'*r=(vr)'+q'(v,r)=(vr)''-a(vr)+q'(v,r),
$$
and on the other hand,
$$
 v'*r=(v''-a(v))*r=(vr)''+q''(v,r)-a(v)r,
$$
since $v'=v''-a(v)$ for all $v\in V$.

 Parts~(b\+-c): we have, on the one hand
\begin{equation} \label{p-h-first-second-one-hand}
 \hi_1'*\hi_2'=p'(\hi_1\ot\hi_2)'-h'(\hi_1\ot\hi_2)=
 p'(\hi_1\ot\hi_2)''-a(p'(\hi_1\ot\hi_2))-h'(\hi_1\ot\hi_2),
\end{equation}
and on the other hand,
\begin{multline} \label{p-h-first-second-other-hand}
 \hi_1'*\hi_2'=(\hi_1''-a(\hi_1))*(\hi_2''-a(\hi_2)) \\ =
 \hi_1''*\hi_2''-a(\hi_1)*\hi_2''-\hi_1''*a(\hi_2)+a(\hi_1)*a(\hi_2) \\
 = p''(\hi_1\ot\hi_2)''-h''(\hi_1\ot\hi_2)-(a(\hi_1)\hi_2)''-
 (\hi_1a(\hi_2))''-q''(\hi_1,a(\hi_2)) + a(\hi_1)a(\hi_2).
\end{multline}
 Comparing~\eqref{p-h-first-second-one-hand}
with~\eqref{p-h-first-second-other-hand} and equating separately
the terms belonging to $V''\subset\tV$ and to $R\subset\tV$,
we obtain part~(b) as well as the equation
$$
 h''(\hi_1\ot\hi_2)=h'(\hi_1\ot\hi_2)+a(p'(\hi_1\ot\hi_2))
 -q''(\hi_1,a(\hi_2))+a(\hi_1)a(\hi_2).
$$
 In order to deduce part~(c), it remains to take into account
the equation
$$
 -q''(\hi_1,a(\hi_2))+a(\hi_1)a(\hi_2)=
 -q'(\hi_1,a(\hi_2))+a(\hi_1a(\hi_2))
$$
obtained by substituting $v=\hi_1$ and $r=a(\hi_2)$ into part~(a).
\end{proof}

\begin{prop} \label{strict-generators-connection-change}
 Let $\tA$ be a\/ $3$\+left finitely projective weak nonhomogeneous
quadratic ring over its subring $R\subset\tA$ with
the $R$\+$R$\+bimodule of generators $R\subset\tV\subset\tA$.
 Let $V\simeq V'\hookrightarrow\tV$ and $V\simeq V''\hookrightarrow\tV$
be two left $R$\+linear splittings of the surjective $R$\+$R$\+bimodule
map $\tV\rarrow\tV/R=V$.
 Denote by\/ $\prB=(B,d',h')$ and\/ $\secB=(B,d'',h'')$ the two
related CDG\+ring structures on the\/ $3$\+right finitely projective
graded ring $B$ quadratic dual to\/ $\q A$, as constructed in
Proposition~\ref{nonhomogeneous-dual-cdg-ring}.
 Let $a\in\Hom_R(V,R)=B^1$ be the element for which the two splittings
$V'\subset\tV$ and $V''\subset\tV$ are related by the rule
$V'=\{v'\mid v\in V\}$ and $V''=\{v''\mid v\in V\}$ with $v'$,
$v''\longmapsto v$ under the map $\tV\rarrow V$ and
$v''=v'+a(v)$ \,\eqref{change-of-strict-generators-eqn}.
 Then the equations
\begin{enumerate}
\renewcommand{\theenumi}{\alph{enumi}}
\item $d_0''(r)=d_0'(r)+ar-ra$ for all $r\in R$;
\item $d_1''(b)=d_1'(b)+ab+ba$ for all $b\in B^1$; and
\item $h''=h'+d'(a)+a^2$
\end{enumerate}
hold in $B$, showing that $(\id,a)\:\secB\rarrow\prB$ is
a CDG\+ring isomorphism.
\end{prop}

\begin{proof}
 In our new notation,
the formulas~(\ref{d-zero-defined}\+-\ref{d-one-defined})
take the form
\begin{gather*}
 \lan v,d_0'(r)\ran = q'(v,r)
 \quad\text{and}\quad
 \lan v,d_0''(r)\ran = q''(v,r), \\
\begin{aligned}
 \lan i,d_1'(b)\ran &=
 \lan p'(\hi_1\ot\hi_2),b\ran - q'(\hi_1,\lan \hi_2,b\ran), \\
 \lan i,d_1''(b)\ran &=
 \lan p''(\hi_1\ot\hi_2),b\ran - q''(\hi_1,\lan \hi_2,b\ran).
\end{aligned}
\end{gather*}

 Part~(a): for every element $v\in V$, one has
\begin{multline*}
 \lan v,d_0''(r)\ran = q''(v,r) = q'(v,r)+a(v)r-a(vr) \\ =
 \lan v,d_1'(r)\ran + \lan v,a\ran r - \lan vr,a\ran =
 \lan v,d_1'(r)\ran + \lan v,ar\ran - \lan v, ra\ran
\end{multline*}
by Proposition~\ref{change-of-strict-generators-equations-prop}(a).

 Part~(b): for every element $i\in I$ and its preimage $\hi\in\hI$,
one has
\begin{multline*}
 \lan i,d_1''(b)\ran = \lan p''(\hi_1\ot\hi_2),b\ran -
 q''(\hi_1,\lan\hi_2,b\ran) \\ =
 \lan p'(\hi_1\ot\hi_2),b\ran + \lan a(\hi_1)\hi_2,b\ran
 + \lan\hi_1a(\hi_2),b\ran \\ - q'(\hi_1,\lan\hi_2,b\ran)
 - a(\hi_1)\lan\hi_2,b\ran + a(\hi_1\lan\hi_2,b\ran) \\ =
 \lan i,d_1'(b)\ran + \lan \hi_1\lan \hi_2,a\ran,b\ran
 + \lan\hi_1\lan\hi_2,b\ran,a\ran =
 \lan i,d'_1(b)\ran + \lan i,ab\ran + \lan i,ba\ran
\end{multline*}
by Proposition~\ref{change-of-strict-generators-equations-prop}(b)
and~(a), the latter of which is being applied to the elements
$v=\hi_1$ and $r=\lan\hi_2,b\ran$.

 Part~(c): for every element $i\in I$ and its preimage $\hi\in\hI$,
one has
\begin{multline*}
 \lan i,h''\ran = \lan i,h'\ran + \lan p'(\hi_1\ot\hi_2),a\ran
 - q'(\hi_1,\lan\hi_2,a\ran) + \lan \hi_1\lan\hi_2,a\ran,a\ran \\ =
 \lan i,h'\ran + \lan i,d_1'(a)\ran + \lan i,a^2\ran
\end{multline*}
by Proposition~\ref{change-of-strict-generators-equations-prop}(c).

 Finally, the equations~(a) and~(b) imply that
$d''(x)=d'(x)+[a,x]$ for all $x\in B$, since $B$ is generated by
$B^1$ over $B^0=R$.
\end{proof}

\subsection{The nonhomogeneous quadratic duality functor}
\label{nonhomogeneous-duality-functor-subsecn}
 Let $R$ be an associative ring.
 We denote by $R\rings_\fil$ the category of filtered rings $(\tA,F)$
with increasing filtrations $F$ and the filtration component $F_0\tA$
identified with~$R$.

 So the objects of $R\rings_\fil$ are associative rings $\tA$ endowed
with a filtration $0=F_{-1}\tA\subset F_0\tA\subset F_1\tA\subset
F_2\tA\subset\dotsb$ such that $\tA=\bigcup_{n=0}^\infty F_n\tA$,
the filtration $F$ is compatible with the multiplication in $\tA$, and
an associative ring isomorphism $R\simeq F_0\tA$ has been chosen.
 Morphisms $(\prtA,F)\rarrow (\sectA,F)$ in $R\rings_\fil$ are ring
homomorphisms $f\:\prtA\rarrow\sectA$ such that
$f(F_n\prtA)\subset F_n\sectA$ for all $n\ge0$ and the ring homomorphism
$F_0f\:F_0\prtA\rarrow F_0\sectA$ forms a commutative triangle diagram
with the fixed isomorphisms $R\simeq F_0\prtA$ and $R\simeq F_0\sectA$.

 The \emph{category of\/ $3$\+left finitely projective weak
nonhomogeneous quadratic rings over~$R$}, denoted by $R\rings_\wnlq$,
is defined as the full subcategory in $R\rings_\fil$ whose objects
are the $3$\+left finitely projective weak nonhomogeneous quadratic
rings $R\subset\tV\subset\tA$ endowed with the filtration $F$
generated by $F_1\tA=\tV$ over $F_0\tA=R$.
 In other words, this means that a morphism $(\prtA,\prtV)\rarrow
(\sectA,\sectV)$ in $R\rings_\wnlq$ is a ring homomorphism $f\:\prtA
\rarrow\sectA$ forming a commutative triangle diagram with
the inclusions $R\simeq F_0\prtA\hookrightarrow\prtA$ and
$R\simeq F_0\sectA\hookrightarrow\sectA$ and satisfying the condition
that $f(\prtV)\subset\sectV$
(cf.\ the discussion in Section~\ref{nonhomogeneous-quadratic-subsecn}).

 Furthermore, the \emph{category of\/ $3$\+right finitely projective
quadratic CDG\+rings over~$R$}, denoted by $R\rings_{\cdg,\rq}$,
is the full subcategory in the category $R\rings_\cdg$ (as defined
in Section~\ref{curved-dg-rings-subsecn}) consisting of all
the CDG\+rings $(B,d,h)$ whose underlying nonnegatively graded ring $B$
is $3$\+right finitely projective quadratic over~$R$.

\begin{thm} \label{nonhomogeneous-duality-functor-existence-thm}
 The constructions of Propositions~\ref{nonhomogeneous-dual-cdg-ring}
and~\ref{strict-generators-connection-change} define a fully
faithful contravariant functor
\begin{equation} \label{nonhomogeneous-duality-functor}
 (R\rings_\wnlq)^\sop\lrarrow R\rings_{\cdg,\rq}
\end{equation}
from the category of\/ $3$\+left finitely projective weak nonhomogeneous
quadratic rings to the category of\/ $3$\+right finitely projective
quadratic CDG\+rings over~$R$.
\end{thm}

\begin{proof}
 Let $R\rings_\wnlq^\sg$ the category whose objects are $3$\+left
finitely projective weak nonhomogeneous quadratic rings $R\subset\tV
\subset\tA$ with a chosen submodule of strict generators
$V'\subset\tV$.
 So $V'$ is a left $R$\+submodule in $\tV$ such that $\tV=R\oplus V'$
as a left $R$\+module.
 Morphisms $(\prtA,\prtV,\prV)\rarrow(\sectA,\sectV,\secV)$ in
$R\rings_\wnlq^\sg$ are the same as morphisms $(\prtA,\prtV)\rarrow
(\sectA,\sectV)$ in $R\rings_\wnlq^\sg$; so a morphism in
$R\rings_\wnlq^\sg$ has to take $\prtV$ into $\sectV$, but it does not
need to respect the chosen submodules of strict generators
$\prV\subset\prtV$ and $\secV\subset\sectV$ in any way.
 Then the functor $R\rings_\wnlq^\sg\rarrow R\rings_\wnlq$ forgetting
the choice of the submodule of strict generators $V'\subset\tV$ is
fully faithful and surjective on objects; so it is an equivalence
of categories.

 Furthermore, let $R\rings_\wnlq^\sgsm$ be the subcategory in
$R\rings_\wnlq^\sg$ whose objects are all the objects of
$R\rings_\wnlq^\sg$ and whose morphisms are the morphisms
$f\:(\prtA,\prtV,\prV)\rarrow(\sectA,\sectV,\secV)$ in
$R\rings_\wnlq^\sg$ such that $f(\prV)\subset \secV$.
 The category $R\rings_\wnlq^\sg$ can be called the category of
$3$\+left finitely projective weak nonhomogeneous quadratic rings
\emph{with strict generators chosen}, while the category
$R\rings_\wnlq^\sgsm$ is the category of $3$\+left finitely projective
weak nonhomogeneous quadratic rings \emph{with strict generators and
strict morphisms}. {\hbadness=1300\par}

 Finally, let $R\rings_\cdg^\sm$ denote the subcategory in
$R\rings_\cdg$ whose objects are all the objects of $R\rings_\cdg$
and whose morphisms are the strict morphisms only, i.~e., all morphisms
of the form $(g,0)\:(\secB,d'',h'')\rarrow(\prB,d',h')$
in $R\rings_\cdg$.
 We denote by $R\rings_{\cdg,\rq}^\sm$ the intersection
$R\rings_\cdg^\sm\cap R\rings_{\cdg,\rq}\subset R\rings_\cdg$,
that is, the category of $3$\+right finitely projective quadratic
CDG\+rings over $R$ and strict morphisms between them.

 Then the construction of Proposition~\ref{nonhomogeneous-dual-cdg-ring}
assigns an object $(B,d,h)\in R\rings_{\cdg,\rq}$ to every object
$(\tA,\tV,V')\in R\rings_\wnlq^\sg$ in a natural way.
 Given a morphism $f\:(\prtA,\prtV,\prV)\rarrow(\sectA,\sectV,\secV)$
in $R\rings_\wnlq^\sg$ such that $f(\prV)\subset \secV$, we consider
the induced morphism of $3$\+left finitely projective quadratic
graded rings $\q\,\gr^Ff\:\q\,\gr^F\prtA\rarrow\q\,\gr^F\sectA$.
 According to Propositions~\ref{2-fin-proj-quadratic-duality}
and~\ref{3-fin-proj-quadratic-duality}, the morphism $\q\,\gr^Ff$
induces a morphism in the opposite direction between the quadratic
dual $3$\+right finitely projective quadratic graded rings,
$g\:\secB\rarrow\prB$.

 Let $(\prB,d',h')$ and $(\secB,d'',h'')$ denote the $3$\+right
finitely projective quadratic CDG\+rings assigned to the $3$\+left
finitely projective weak nonhomogeneous quadratic rings
$(\prtA,\prtV)$ and $(\sectA,\sectV)$ with the submodules of strict
generators $\prV\subset\prtV$ and $\secV\subset\sectV$ by
the construction of Proposition~\ref{nonhomogeneous-dual-cdg-ring}.
 Assigning the morphism $(g,0)\:(\secB,d'',h'')\rarrow(\prB,d',h')$
to the morphism $f\:(\prtA,\prtV,\prV)\rarrow(\sectA,\sectV,\secV)$,
we obtain a contravariant functor
\begin{equation} \label{strict-morphisms-duality-functor}
 (R\rings_\wnlq^\sgsm)^\sop\lrarrow R\rings_{\cdg,\rq}^\sm.
\end{equation}

 We still have to check that
the functor~\eqref{strict-morphisms-duality-functor} is well-defined,
i.~e., that $(g,0)\:\allowbreak (\secB,d'',h'')\rarrow
(\prB,d',h')$ is indeed a morphism of CDG\+rings.
 Simultaneously we will see that
the functor~\eqref{strict-morphisms-duality-functor} is fully faithful.

 Indeed, specifying a morphism $f\:(\prtA,\prtV,\prV)\rarrow
(\sectA,\sectV,\secV)$ in $R\rings_\wnlq^\sgsm$ means specifying
an $R$\+$R$\+bimodule map $f_1\:\prtV/R\rarrow\sectV/R$ which,
interpreted as a map $\prV\rarrow \secV$ and taken together with
the identity map $R\rarrow R$, extends (necessarily uniquely)
to a ring homomorphism $\prtA\rarrow\sectA$.
 The latter condition is equivalent to the following two:
\begin{enumerate}
\renewcommand{\theenumi}{\roman{enumi}}
\item the map $\tilde f_1=\id_R\oplus f_1\:\sectV=R\oplus \secV\rarrow R
\oplus \prV=\prtV$ agrees with the right $R$\+module structures on
$\sectV$ and $\prtV$;
\item the tensor ring homomorphism $T_R(\tilde f_1)\:T_R(\prtV)
\rarrow T_R(\sectV)$ induced by the map $\tilde f_1\:\prtV\rarrow\sectV$
takes the $R$\+$R$\+subbimodule $\tI_{\prtA}\subset R\oplus\prtV\oplus
\prtV{}^{\ot_R\,2}$ of nonhomogeneous quadratic relations in the ring
$\prtA$ into the $R$\+$R$\+subbimodule $\tI_{\sectA}\subset R\oplus
\sectV\oplus\sectV{}^{\ot_R\,2}$ of nonhomogeneous quadratic relations
in the ring~$\sectA$ (see Section~\ref{nonhomogeneous-quadratic-subsecn}
for the notation).
\end{enumerate}

 Denote the maps defined by
the formulas~(\ref{map-q}\+-\ref{maps-p-and-h}) for
$(\prtA,\prtV,\prV)$ and $(\sectA,\sectV,\secV)$ by
\begin{gather*}
 q'\:\prV\times R\rarrow R, \ \ p'\:\prhI\rarrow \prV, \ \ 
 h'\:\prhI\rarrow R, \\
 q''\:\secV\times R\rarrow R, \ \ p''\:\sechI\rarrow \secV, \ \
 h''\:\sechI\rarrow R,
\end{gather*}
where the notation $\prV=\prtV/R$ and $\secV=\sectV/R$ is presumed, while
$\prhI\subset \prV\ot_\boZ \prV$ and $\sechI\subset \secV\ot_\boZ \secV$
are the full preimages of the $R$\+$R$\+subbimodules
$\prI\subset \prV\ot_R\prV$ and $\secI\subset \secV\ot_R \secV$
of quadratic relations in the graded rings $\prA=\gr^F\prtA$
and $\secA=\gr^F\sectA$, respectively.
 Then condition~(i) is equivalent to the equation
\begin{equation} \label{q-f-commute}
 q'(v,r)=q''(f_1(v),r) \quad\text{for all $v\in \prV$, \,$r\in R$}.
\end{equation}
 Assuming~(i) or~\eqref{q-f-commute}, condition~(ii) is equivalent to
the combination of the inclusion
\begin{equation} \label{I-f-preserves}
 (f_1\ot f_1)(\prI)\,\subset\,\secI
\end{equation}
with the equations
\begin{equation} \label{p-h-f-commute}
\begin{gathered}
 f_1(p'(\hi_1\ot\hi_2))=p''(f_1(\hi_1)\ot f_1(\hi_2))
 \quad\text{for all $\hi_1\ot\hi_2\in\prhI$} \\
 h'(i_1'\ot i_2')=h''(f_1(i_1)\ot f_1(i_2))
 \quad\text{for all $i_1\ot i_2\in\prI$}.
\end{gathered}
\end{equation}

 Finally, the inclusion~\eqref{I-f-preserves} holds if and only if
the $R$\+$R$\+bimodule morphism $g_1=\Hom_R(f_1,R)\:
\secB^1=\Hom_R(\secV,R)\rarrow\Hom_R(\prV,R)=\prB^1$ together with
the identity map $\secB^0=R\rarrow R=\prB^0$ can be extended to
a graded ring homomorphism $g\:\secB\rarrow\prB$.
 Assuming~\eqref{I-f-preserves} or (equivalently) the existence of
$g=(g_n)_{n=0}^\infty$, the equations~\eqref{q-f-commute}
and~\eqref{p-h-f-commute} are equivalent to the equations
\begin{equation} \label{g-strict-CDG-morphism}
\begin{gathered}
 d'_0(r)=g_1(d''_0(r)) \quad\text{for all $r\in R$}, \\
 d'_1(g_1(b))=g_2(d''_1(b)) \quad\text{for all $b\in\secB^1$}, \\
 h'=g_2(h''),
\end{gathered}
\end{equation}
which mean that $(g,0)\:(\secB,d'',h'')\rarrow(\prB,d',h')$
is a CDG\+ring morphism.

 Conversely, given a strict morphism $(g,0)\:(\secB,d'',h'')
\rarrow(\prB,d',h')$ between $3$\+right finitely projective
quadratic CDG\+rings coming from $3$\+left finitely projective
nonhomogeneous quadratic rings $(\prtA,\prtV,\prV)$ and
$(\sectA,\sectV,\secV)$, we put
$f_1=\Hom_{R^\rop}(g_1,R)\:\prV=\Hom_{R^\rop}(\prB^1,R)
\rarrow\Hom_{R^\rop}(\secB^1,R)=\secV$ and extend the map~$f_1$
together with the identity map $\prtA\supset R\rarrow R\subset\sectA$
to a ring homomorphism $f\:\prtA\rarrow\sectA$.
 This can be done, since the combination of conditions~(i) and~(ii)
is equivalent to the existence of a (necessarily unique) graded ring
homomorphism $g\:\secB\rarrow\prB$ extending the given map
$g_1\:\secB^1\rarrow\prB^1$ together with the identity map
$\secB^0=R\rarrow R=\prB^0$ and satisfying
the equations~\eqref{g-strict-CDG-morphism}.

 Now we will construct a fully faithful functor
\begin{equation} \label{strict-generators-chosen-duality-functor}
 (R\rings_\wnlq^\sg)^\sop\lrarrow R\rings_{\cdg,\rq}
\end{equation}
extending the functor~\eqref{strict-morphisms-duality-functor}
to nonstrict morphisms.
 For a $3$\+left finitely projective weak nonhomogeneous
quadratic ring $(\tA,\tV)\in R\rings_\wnlq$ and any two choices
of a submodule of strict generators $V'$, $V''\subset\tV$,
we have two objects $(\tA,\tV,V')$ and $(\tA,\tV,V'')\in
R\rings_\wnlq^\sg$ connected by an isomorphism
$(\tA,\tV,V')\rarrow(\tA,\tV,V'')$ corresponding to the identity
map $\id\:\tA\rarrow\tA$.
 Let us call such isomorphisms in $R\rings_\wnlq^\sg$
the \emph{change-of-strict-generators isomorphisms}.
 
 Let $f\:(\prtA,\prtV,\prV)\rarrow(\sectA,\sectV,\secV)$ be
an arbitrary morphism in $R\rings_\wnlq^\sg$.
 Denote by $\prV'\subset\prtV$ the full preimage of the left
$R$\+submodule $\secV\subset\sectV$ under the $R$\+$R$\+bimodule
morphism $F_1f\:\prtV\rarrow\sectV$.
 Then one has $\prtV=R\oplus\prV'$; so $\prV'\subset\prtV$ is
another choice of a submodule of strict generators in $\prtA$,
alternative to $\prV\subset\prtV$.
 Any morphism $f\:(\prtA,\prtV,\prV)\rarrow(\sectA,\sectV,\secV)$ in
$R\rings_\wnlq^\sg$ decomposes uniquely into
a change-of-strict-generators isomorphism
$(\prtA,\prtV,\prV)\rarrow(\prtA,\prtV,\prV')$ followed by
a strict morphism $(\prtA,\prtV,\prV')\rarrow(\sectA,\sectV,\secV)$.

 The construction of
Proposition~\ref{strict-generators-connection-change} assigns
a change-of-connection isomorphism in $R\rings_{\cdg,\rq}$ to
every change-of-strict-generators isomorphism in
$R\rings_\wnlq^\sg$.
 Decomposing any morphism in $R\rings_\wnlq^\sg$ into
a change-of-strict-generators isomorphism followed by a strict morphism,
one extends the functor~\eqref{strict-morphisms-duality-functor} to
a contravariant
functor~\eqref{strict-generators-chosen-duality-functor}.
 We omit further details, which are straightforward.

 We still have to show that
the functor~\eqref{strict-generators-chosen-duality-functor} is
fully faithful.
 It is clear from the construction of
Proposition~\ref{strict-generators-connection-change} that this functor
functor restricts to a fully faithful functor from the subcategory
of change-of-strict-generators isomorphisms in $R\rings_\wnlq^\sg$ to
the subcategory of change-of-connection isomorphisms in
$R\rings_{\cdg,\rq}$.
 Since the functor~\eqref{strict-morphisms-duality-functor} is fully
faithful as well, it follows that so is
the functor~\eqref{strict-generators-chosen-duality-functor}.

 In order to construct the desired fully faithful
functor~\eqref{nonhomogeneous-duality-functor}, it remains to
choose any quasi-inverse functor to the category equivalence
$R\rings_\wnlq^\sg\rarrow R\rings_\wnlq$ and compose it
with the fully faithful
functor~\eqref{strict-generators-chosen-duality-functor}.
\end{proof}

\begin{rem}
 The counterexample in~\cite[Section~3.4]{Pcurv} shows that the fully
faithful contravariant functor in
Theorem~\ref{nonhomogeneous-duality-functor-existence-thm} is
\emph{not} an anti-equivalence of categories (even when $R=k$ is
the ground field).
 We will see below in Section~\ref{pbw-theorem-subsecn} that this
functor becomes an anti-equivalence when restricted to the full
subcategories of, respectively, left and right finitely projective
Koszul rings in $R\rings_\wnlq$ and $R\rings_{\cdg,\rq}$.
\end{rem}

\subsection{Nonhomogeneous duality 2-functor}
\label{nonhomogeneous-duality-2-functor-subsecn}
 We define the \emph{$2$\+category of filtered rings} $\Rings_{\fil2}$
as follows.
 The objects of $\Rings_{\fil2}$ are associative rings $\tA$ endowed
with an exhastive increasing filtration $0=F_{-1}\tA\subset F_0\tA
\subset F_1\tA\subset F_2\tA\subset\dotsb$ compatible with
the multiplication on~$\tA$.
 Morphisms $(\prtA,F)\rarrow(\sectA,F)$ in $\Rings_{\fil2}$ are ring
homomorphisms $f\:\prtA\rarrow\sectA$ such that $f(F_n\prtA)\subset
F_n\sectA$ for all $n\ge0$ and the map $F_0f\:F_0\prtA\rarrow F_0\sectA$
is an isomorphism.
 $2$\+morphisms $f\overset z\rarrow g$ between a pair of parallel
morphisms $f$, $g\:(\prtA,F)\rarrow(\sectA,F)$ are invertible elements
$z\in F_0\sectA$ such that $g(c)=zf(c)z^{-1}$ for all $c\in\prtA$.

 The vertical composition of two $2$\+morphisms
$f'\overset w\rarrow f''\overset z\rarrow f'''$ is the $2$\+morphism
$f'\overset{zw}\rarrow f'''$.
 The identity $2$\+morphism is the $2$\+morphism $f\overset1\rarrow f$.
 The horizontal composition of two $2$\+morphisms
$g'\overset w\rarrow g''\:(\tA,F)\rarrow(\tB,F)$
and $f'\overset z\rarrow f''\:(\tB,F)\rarrow(\tC,F)$
is the $2$\+morphism $f'g'\overset{z\circ w}\lrarrow f''g''\:
(\tA,F)\rarrow(\tC,F)$ with the element
$z\circ w=zf'(w)=f''(w)z$.

 All the $2$\+morphisms of filtered rings are invertible.
 If $f\:(\prtA,F)\rarrow(\sectA,F)$ is a morphism of filtered rings
and $z\in F_0\sectA$ is an invertible element, then $g\:c\longmapsto
zf(c)z^{-1}$ is also a morphism of filtered rings
$g\:(\prtA,F)\rarrow(\sectA,F)$.
 The morphisms~$f$ and~$g$ are connected by the $2$\+isomorphism
$f\overset z\rarrow g$.

 The condition about the map $F_0f\:F_0\prtA\rarrow F_0\sectA$ being
an isomorphism could be harmlessly dropped from the above definition
(which makes sense without this condition just as well), but we need
it for the purposes of the next definition.
 The \emph{$2$\+category of\/ $3$\+left finitely projective weak
nonhomogeneous quadratic rings}, denoted by $\Rings_{\wnlq2}$, is
defined as the following $2$\+subcategory in $\Rings_{\fil2}$.
 The objects of $\Rings_{\wnlq2}$ are the $3$\+left finitely
projective weak nonhomogeneous quadratic rings $R\subset\tV\subset\tA$
with the filtration $F$ generated by $F_1\tA=\tV$ over $F_0\tA=R$.
 All morphisms in $\Rings_{\fil2}$ between objects of $\Rings_{\wnlq2}$
are morphisms in $\Rings_{\wnlq2}$, and all $2$\+morphisms in
$\Rings_{\fil2}$ between morphisms of $\Rings_{\wnlq2}$ are
$2$\+morphisms in $\Rings_{\wnlq2}$.

 Furthermore, the \emph{$2$\+category of\/ $3$\+right finitely
projective quadratic CDG\+rings}, denoted by $\Rings_{\cdg2,\rq}$, is
the following $2$\+subcategory in the $2$\+category $\Rings_{\cdg2}$
(as defined in Section~\ref{curved-dg-rings-subsecn}).
 The objects of $\Rings_{\cdg2,\rq}$ are all the CDG\+rings
$(B,d,h)$ whose underlying nonnegatively graded ring $B$ is
$3$\+right finitely projective quadratic (over its degree-zero
component~$B^0$).
 All morphisms in $\Rings_{\cdg2}$ between objects of
$\Rings_{\cdg2,\rq}$ are morphisms in $\Rings_{\cdg2,\rq}$, and all
$2$\+morphisms in $\Rings_{\cdg2}$ between morphisms of
$\Rings_{\cdg2,\rq}$ are $2$\+morphisms in $\Rings_{\cdg2,\rq}$.

\begin{lem} \label{conjugation-gauge-transformation-lemma}
 Let $R\subset\tV\subset\tA$ be a\/ $3$\+left finitely projective
weak nonhomogeneous quadratic ring, and let $\prV\subset\tV$ be
a submodule of strict generators of~$\tA$.
 Let $(B,d,h)$ be the $3$\+right finitely projective quadratic
CDG\+ring corresponding to $(\tA,\tV,\prV)$ under the construction of
Proposition~\ref{nonhomogeneous-dual-cdg-ring}.
 Let $z\in R$ be an invertible element.
 Consider the conjugation morphism $f_{z^{-1}}\:\tA\rarrow\tA$ taking
any element $c\in\tA$ to the element $f_{z^{-1}}(c)=z^{-1}cz$.
 Then the CDG\+ring morphism $(B,d,h)\rarrow(B,d,h)$ corresponding to
the morphism~$f_{z^{-1}}$ under the duality functor of
Theorem~\ref{nonhomogeneous-duality-functor-existence-thm} is equal to
$$
 (g_z,a_z)\:(B,d,h)\lrarrow(B,d,h),
$$
where $g_z\:B\rarrow B$ is the conjugation map taking any element
$b\in B$ to the element $g_z(b)=zbz^{-1}$, and the element $a_z\in B^1$
is given by the formula $a_z=-d(z)z^{-1}$.
\end{lem}

\begin{proof}
 Strictly speaking, the assertion of the lemma does not literally make
sense as stated, and we need to make it more precise before proving it.
 The problem is that $f_{z^{-1}}\:\tA\rarrow\tA$ is \emph{not}
a morphism in the category $R\rings_\wnlq$ or $R\rings_\fil$, as defined
in Section~\ref{nonhomogeneous-duality-functor-subsecn}, because it does
not restrict to the identity map $\tA\supset R\rarrow R\subset\tA$, but
rather to the map $R\rarrow R$ of conjugation with~$z^{-1}$.
 For the same reason, $g_z\:B\rarrow B$ is not a morphism in
$R\rings_\sgr$, and consequently $(g_z,a_z)\:(B,d,h)\allowbreak
\rarrow(B,d,h)$ is not a morphism in $R\rings_\cdg$ or
$R\rings_{\cdg,\rq}$.
 So some preparatory work is needed before
the functor~\eqref{nonhomogeneous-duality-functor} could be applied
to the map~$f_{z^{-1}}$.

 Denote by $t_z\:R\rarrow R$ the conjugation map
$r\longmapsto zrz^{-1}$, \ $r\in R$.
 Let $(\tA^{(z)},\tV^{(z)})$ denote the following $3$\+left finitely
projective weak nonhomogeneous quadratic ring over~$R$.
 As an associative ring, $\tA^{(z)}$ coincides with~$\tA$.
 Denoting by $\iota\:R\rarrow \tA$ the embedding of the ring $R$ as
a subring of the ring~$\tA$, the map $\iota^{(z)}=\iota t_z\:R
\rarrow \tA=\tA^{(z)}$ makes $R$ a subring of~$\tA^{(z)}$.
 For any element $c\in \tA$, we will denote by $c^{(z)}\in \tA^{(z)}$
the corresponding element of the ring~$\tA^{(z)}$.
 So in particular, we have $\iota^{(z)}(r)=(z\iota(r)z^{-1})^{(z)}$
for all $r\in R$.
 Furthermore, as a subgroup of $\tA^{(z)}=\tA$, the group $\tV^{(z)}$
coincides with~$\tV$. 
 Then the map $f_{z^{-1}}\:\tA^{(z)}\rarrow \tA$ taking an element
$c^{(z)}\in \tA^{(z)}$ to the element $z^{-1}cz\in \tA$ is a morphism
in the category $R\rings_\wnlq$.

 Let $(B^{(z)},d^{(z)},h^{(z)})$ denote the following $3$\+right
finitely projective quadratic CDG\+ring over~$R$.
 As a graded associative ring, $B^{(z)}$ coincides with~$B$; and
both the differential $d^{(z)}\:B^{(z)}\rarrow B^{(z)}$ and
the curvature element $h\in B^{(z),2}$ coincide with
the differential~$d$ and the curvature element~$h$ in~$B$.
 However, denoting by $\iota\:R\rarrow B^0$ the identification of
the ring $R$ with the degree-zero component of the graded ring $B$,
the identification of the ring $R$ with the degree-zero component
of the graded ring $B^{(z)}$ is provided by the map $\iota^{(z)}=
\iota t_z\:R\rarrow B^{(z),0}=B^0$.

 For any element $b\in B$, we denote by $b^{(z)}\in B^{(z)}$
the corresponding element of the ring~$B^{(z)}$.
 So, in particular, we have $d^{(z)}(b^{(z)})=(d(b))^{(z)}$ for all
$b\in B$, and our notation for the element $h^{(z)}$ is consistent:
$h^{(z)}\in B^{(z),2}$ is the element corresponding to $h\in B^2$.
 Then the map $g_z\:B\rarrow B^{(z)}$ taking an element $b\in B$ to
the element $(zbz^{-1})^{(z)}\in B^{(z)}$ is a morphism in
the category $R\rings_\sgr$.
 Furthermore, let $a_z^{(z)}\in B^{(z),1}$ be the element corresponding
to the element $a_z=-d(z)z^{-1}\in B^1$ under the identity isomorphism
$B=B^{(z)}$.
 Then the pair $(g_z,a_z^{(z)})$ is a morphism of CDG\+rings
$(B,d,h)\rarrow(B^{(z)},d^{(z)},h^{(z)})$, as one can readily check.
 Therefore, it is also a morphism in the category $R\rings_\cdg$ and
in the category $R\rings_{\cdg,\rq}$.

 Now the promised precise formulation of the lemma claims that
the duality functor~\eqref{nonhomogeneous-duality-functor} takes
the morphism $f_{z^{-1}}\:\tA^{(z)}\rarrow \tA$ in the category
$R\rings_\wnlq$ to the morphism $(g_z,a_z^{(z)})\:(B,d,h)
\rarrow(B^{(z)},d^{(z)},h^{(z)})$ in the category $R\rings_{\cdg,\rq}$.

 To be even more precise, we need to establish first that
the functor~\eqref{nonhomogeneous-duality-functor} takes the object
$\tA^{(z)}\in R\rings_\wnlq$ to the object $(B^{(z)},d^{(z)},h^{(z)})
\in R\rings_{\cdg,\rq}$.
 Put $A^{(z)}=\gr^F\tA^{(z)}$; then the identity isomorphism of rings
$A^{(z)}_0=B^{(z),0}$ commuting with the identifications
$A^{(z)}_0\simeq R\simeq B^{(z),0}$ together with
the $A^{(z)}_0$\+$A^{(z)}_0$\+bimodule
isomorphism $\Hom_{A^{(z)}_0}(A^{(z)}_1,A^{(z)}_0)=B^{(z),1}$
allow to consider $B^{(z)}$ as the quadratic dual ring to~$A^{(z)}$.
 Let the subgroup $\prV^{(z)}\subset\tV^{(z)}$ coincide with
the subgroup $\prV\subset\tV$; then the construction of
Proposition~\ref{nonhomogeneous-dual-cdg-ring} takes the $3$\+left
finitely projective weak nonhomogeneous quadratic ring
$(\tA^{(z)},\tV^{(z)})$ with the submodule of strict generators
$\prV^{(z)}\subset\tV^{(z)}$ to the $3$\+right finitely projective
quadratic CDG\+ring $(B^{(z)},d^{(z)},h^{(z)})$.

 At last, we can perform the computation proving the lemma.
 We have to check that
the functor~\eqref{strict-generators-chosen-duality-functor} takes
the morphism $f_{z^{-1}}\:(\tA^{(z)},\tV^{(z)},\prV^{(z)})\rarrow
(\tA,\tV,\prV)$ to the morphism
$(g_z,a_z^{(z)})\:(B,d,h)\rarrow(B^{(z)},d^{(z)},h^{(z)})$.

 Set $\secV^{(z)}=z\prV^{(z)}z^{-1}\subset\tV^{(z)}$.
 Then our morphism $f_{z^{-1}}\:(\tA^{(z)},\tV^{(z)},\prV^{(z)})\rarrow
(\tA,\tV,\prV)$ decomposes into a change-of-strict-generators morphism
$(\tA^{(z)},\tV^{(z)},\prV^{(z)})\allowbreak\rarrow
(\tA^{(z)},\tV^{(z)},\secV^{(z)})$, which acts by the identity map on
the underlying associative ring $\tA^{(z)}$, followed by
the strict morphism $f_{z^{-1}}\:(\tA^{(z)},\tV^{(z)},\secV^{(z)})
\rarrow(\tA,\tV,\prV)$.
 Denote by $(B^{(z)},d_{(z)},h_{(z)})$ the CDG\+ring assigned to
the $3$\+left finitely projective weak nonhomogeneous quadratic ring
$(\tA^{(z)},\tV^{(z)})$ with the submodule of strict generators
$\secV^{(z)}\subset\tV^{(z)}$ by the construction of
Proposition~\ref{nonhomogeneous-dual-cdg-ring}.

 Let $q'\:V^{(z)}\times R\rarrow R$ be the map defined by
the formula~\eqref{map-q} using the splitting
$\prV^{(z)}\subset\tV^{(z)}$ of the bimodule of generators of
the $3$\+left finitely projective weak nonhomogeneous quadratic
ring~$\tA^{(z)}$.
 Denoting by $u'\in\prV^{(z)}$ and $u''\in\secV^{(z)}$ the elements
corresponding to an element $u=v^{(z)}\,\in\,V^{(z)}=\tV^{(z)}/R$,
we have
$$
 u''=z*(z^{-1}uz)'*z^{-1}=(uz)'*z^{-1}=u'+q'(uz,z^{-1}).
$$
 Furthermore, by Proposition~\ref{self-consistency-equations-prop}(b),
$$
 0=q'(u,1)=q'(u,zz^{-1})=q'(uz,z^{-1})+q'(u,z)z^{-1},
$$
hence
$$
 u''=u'-q'(u,z)z^{-1}=u'-\lan u,d_0(z)\ran z^{-1}=
 u'-\lan u,d_0(z)z^{-1}\ran = u'+a_z^{(z)}(u).
$$

 By Proposition~\ref{strict-generators-connection-change}, it follows
that $(\id,a_z^{(z)})\:(B^{(z)},d_{(z)},h_{(z)})\rarrow
(B^{(z)},d^{(z)},h^{(z)})$ is a change-of-connection morphism
of CDG\+rings over~$R$.
 Hence one can compute that $d_{(z)}(zb^{(z)}z^{-1})=
zd^{(z)}(b^{(z)})z^{-1}=z(d(b))^{(z)}z^{-1}$ for all $b\in B$,
and $h_{(z)}=zh^{(z)}z^{-1}$.
 By construction, $(\id,a_z^{(z)})\:
(B^{(z)},d_{(z)},h_{(z)})\rarrow (B^{(z)},d^{(z)},h^{(z)})$
is \emph{the} change-of-connection morphism of CDG\+rings over $R$
assigned to the change-of-strict-generators morphism
$(\tA^{(z)},\tV^{(z)},\prV^{(z)})\rarrow
(\tA^{(z)},\tV^{(z)},\secV^{(z)})$
by the functor~\eqref{strict-generators-chosen-duality-functor}.

 Finally, we need to show that
the functor~\eqref{strict-morphisms-duality-functor} takes the strict
morphism $f_{z^{-1}}\:(\tA^{(z)},\tV^{(z)},\secV^{(z)})\rarrow
(\tA,\tV,\prV)$ to the strict morphism $(g_z,0)\:(B,d,h)\rarrow
(B^{(z)},d_{(z)},h_{(z)})$.
 For this purpose, it suffices to check that the morphism
$\bar f_{z^{-1}}=\gr^Ff_{z^{-1}}\:A^{(z)}=\gr^F\tA^{(z)}\rarrow\gr^F\tA=A$
corresponds to the morphism $g_z\:B\rarrow B^{(z)}$ under
the homogeneous quadratic duality of
Propositions~\ref{2-fin-proj-quadratic-duality}
and~\ref{3-fin-proj-quadratic-duality}.
 All we need to do is to observe that
\begin{multline*}
 \lan v^{(z)},g_z(b)\ran = \lan v^{(z)},zb^{(z)}z^{-1}\ran
 =z\lan z^{-1}v^{(z)}z,b^{(z)}\ran z^{-1} \\ 
 =t_z(\lan z^{-1}v^{(z)}z,b^{(z)}\ran)
 =\lan z^{-1}vz,b\ran=
 \lan\bar f_{z^{-1}}(v^{(z)}),b\ran
\end{multline*}
for all $v^{(z)}\in V^{(z)}=A^{(z)}_1$ and $b\in B^1$.

 It remains to compute the image of the composition of our morphisms
in the category $R\rings_\wnlq^\sg$
$$
 (\tA^{(z)},\tV^{(z)},\prV^{(z)})\lrarrow
 (\tA^{(z)},\tV^{(z)},\secV^{(z)})\overset{f_{z^{-1}}}
 \lrarrow(\tA,\tV,\prV)
$$
under the functor~\eqref{strict-generators-chosen-duality-functor}.
 This is equal, by construction, to the composition of
morphisms of CDG\+rings
$$
 (\id,a_z^{(z)})\circ(g_z,0)=(g_z,a_z^{(z)}),
$$
as desired.

 In the context of
the functor~\eqref{enhanced-nonhomogeneous-duality-functor} below
instead of the functor~\eqref{nonhomogeneous-duality-functor},
the assertion of the lemma becomes literally true as stated, without
the additional discussion in the first half of the above proof.
\end{proof}

\begin{thm} \label{nonhomogeneous-duality-2-functor-theorem}
 The constructions of
Theorem~\ref{nonhomogeneous-duality-functor-existence-thm} and
Lemma~\ref{conjugation-gauge-transformation-lemma} define a fully
faithful strict contravariant\/ $2$\+functor
\begin{equation} \label{nonhomogeneous-duality-2-functor}
 (\Rings_{\wnlq2})^\sop\lrarrow\Rings_{\cdg2,\rq}
\end{equation}
from the\/ $2$\+category of\/ $3$\+left finitely projective weak
nonhomogeneous quadratic rings to the\/ $2$\+category of\/
$3$\+right finitely projective quadratic CDG\+rings.
\end{thm}

\begin{proof}
 The $2$\+functor~\eqref{nonhomogeneous-duality-2-functor} is fully
faithful in the strict sense: for any two objects $(\prtA,\prtV)$ and
$(\sectA,\sectV)\in\Rings_{\wnlq2}$ and the corresponding
CDG\+rings $(\prB,d',h')$ and $(\secB,d'',h'')\in\Rings_{\cdg2,\rq}$,
the $2$\+functor~\eqref{nonhomogeneous-duality-2-functor} induces
a bijection between morphisms $(\prtA,\prtV)\rarrow(\sectA,\sectV)$
in $\Rings_{\wnlq2}$ and morphisms $(\secB,d'',h'')\rarrow
(\prB,d',h')$ in $\Rings_{\cdg2,\rq}$.
 Furthermore, for any pair of parallel morphisms $f'$, $f''\:
(\prtA,\prtV)\rarrow(\sectA,\sectV)$ in $\Rings_{\wnlq2}$ and
the corresponding pair of parallel morphisms $g'$, $g''\:
(\secB,d'',h'')\rarrow(\prB,d',h')$ in $\Rings_{\cdg2,\rq}$,
the $2$\+functor~\eqref{nonhomogeneous-duality-2-functor} induces
a bijection between $2$\+morphisms $f'\overset z\rarrow f''$ in
$\Rings_{\wnlq2}$ and $2$\+morphisms $g''\overset w\rarrow g'$ in
$\Rings_{\cdg2,\rq}$. {\emergencystretch=1em\par}

 To construct the desired $2$\+functor, denote by
$\Rings_\wnlq\subset\Rings_{\wnlq2}$ the category whose objects are
the objects of the $2$\+category $\Rings_{\wnlq2}$ and whose morphisms
are the morphisms of the $2$\+category $\Rings_{\wnlq2}$ (but
there are no $2$\+morphisms in $\Rings_\wnlq$).
 Similarly, denote by $\Rings_{\cdg,\rq}$ the category whose objects
are the objects of the $2$\+category $\Rings_{\cdg2,\rq}$ and whose
morphisms are the morphisms of the $2$\+category $\Rings_{\cdg2,\rq}$
(but there are no $2$\+morphisms in $\Rings_{\cdg,\rq}$).
 Then essentially the same construction that was used to define
the functor~\eqref{nonhomogeneous-duality-functor} in
Theorem~\ref{nonhomogeneous-duality-functor-existence-thm} provides
a fully faithful contravariant functor
\begin{equation} \label{enhanced-nonhomogeneous-duality-functor}
 (\Rings_\wnlq)^\sop\lrarrow\Rings_{\cdg,\rq}.
\end{equation}

 In order to extend
the functor~\eqref{enhanced-nonhomogeneous-duality-functor}
to a $2$\+functor~\eqref{nonhomogeneous-duality-2-functor}, we
notice that, in the notation of the first paragraph of this proof,
for every $2$\+morphism $f'\overset z\rarrow f''$ in $\Rings_{\wnlq2}$,
both the morphisms $F_0f'\:F_0\prtA\rarrow F_0\sectA$ and
$F_0f''\:F_0\prtA\rarrow F_0\sectA$ are ring isomorphisms
and the preimages of the element $z\in F_0\sectA$ under these two
isomorphisms coincide, $(F_0f')^{-1}(z)=(F_0f'')^{-1}(z)$,
because the element~$z$ is preserved by the conjugation with~$z$.
 Similarly, for any $2$\+morphism $g''\overset w\rarrow g'$ in
$\Rings_{\cdg2,\rq}$, both the morphisms $g'_0\:\secB^0\rarrow
\prB^0$ and $g''_0\:\secB^0\rarrow\prB^0$ are ring isomorphisms and
the preimages of the element $w\in\prB^0$ under these two isomorphisms
coincide, $g'_0{}^{-1}(w)=g''_0{}^{-1}(w)$.
 By construction, we have $\prB^0=F_0\prtA$ and $\secB^0=F_0\sectA$.
 The maps $F_0f'\:F_0\prtA\rarrow F_0\sectA$ and $g'_0\:\secB^0\rarrow
\prB^0$ are mutually inverse under this identification, and so are
the two maps $F_0f''\:F_0\prtA\rarrow F_0\sectA$ and
$g''_0\:\secB^0\rarrow\prB^0$, that is $g'_0=(F_0f')^{-1}$
and $g''_0=(F_0f'')^{-1}$.

 We assign a $2$\+morphism $g''\overset w\rarrow g'$ in
$\Rings_{\cdg2,\rq}$ to a $2$\+morphism $f'\overset z\rarrow f''$
in $\Rings_{\wnlq2}$ if $g'_0(z)=w=g''_0(z)$, or equivalently,
if $(F_0f')(w)=z=(F_0f'')(w)$.
 It only needs to be checked that $g''\overset w\rarrow g'$ is
a $2$\+morphism in $\Rings_{\cdg2,\rq}$ if and only if
$f'\overset z\rarrow f''$ is a $2$\+morphism in $\Rings_{\wnlq2}$.
 Then the compatibility with the vertical and horizontal
compositions of $2$\+morphisms will be clear from the construction of
such compositions in the beginning of this section and in
Section~\ref{curved-dg-rings-subsecn}.

 Let us define \emph{basic\/ $2$\+morphisms} in the $2$\+category
$\Rings_{\wnlq2}$ as $2$\+morphisms of the form $f_{z^{-1}}
\overset z\rarrow \id_{(\tA,\tV)}$, where $(\tA,\tV)$ is an object
of $\Rings_{\wnlq}$, \ $z\in F_0\tA$ is an invertible element,
$f_{z^{-1}}\:(\tA,\tV)\rarrow(\tA,\tV)$ is the morphism taking
an element $c\in\tA$ to the element $f_{z^{-1}}(c)=z^{-1}cz\in\tA$,
and $\id_{(\tA,\tV)}\:(\tA,\tV)\rarrow(\tA,\tV)$ is the identity
morphism.
 Then any $2$\+morphism in $\Rings_{\wnlq2}$ decomposes uniquely as
a morphism followed by a basic $2$\+morphism, and also as a basic
$2$\+morphism followed by a morphism.
 Specifically, a $2$\+morphism $f'\overset z\rarrow f''$ as above
is the composition of the morphism $f''$ followed by the basic
$2$\+morphism $f_{z^{-1}}\overset z\rarrow\id_{(\sectA,\sectV)}$;
and it is also the composition of the basic $2$\+morphism
$f_{w^{-1}}\overset w\rarrow\id_{(\prtA,\prtV)}$ followed by
the $2$\+morphism~$f''$.

 Similarly, we define \emph{basic\/ $2$\+morphisms} in the $2$\+category
$\Rings_{\cdg2,\rq}$ as $2$\+morphisms of the form
$(\id_B,0)\overset z\rarrow(g_z,a_z)\:(B,d,h)\rarrow(B,d,h)$,
where $(B,d,h)$ is an object of $\Rings_{\cdg,\rq}$, \ $z\in B^0$ is
an invertible element, $g_z\:B\rarrow B$ is the graded ring
homomorphism taking an element $b\in B$ to the element $g_z(b)=
zbz^{-1}\in B$, and $a_z=-d(z)z^{-1}\in B^1$.
 Then any $2$\+morphism in $\Rings_{\wnlq2}$ decomposes uniquely as
a morphism followed by a basic $2$\+morphism, and also as a basic
$2$\+morphism followed by a morphism.
 Specifically, a $2$\+morphism $g''\overset w\rarrow g'$ as above is
the composition of the morphism~$g''$ followed by the basic
$2$\+morphism $(\id_{\prB},0)\overset w\rarrow (g_w,a_w)\:
(\prB,d',h')\rarrow(\prB,d',h')$, and it is also the composition
of the basic $2$\+morphism $(\id_{\secB},0)\overset z\rarrow (g_z,a_z)
\:(\secB,d'',h'')\rarrow(\secB,d'',h'')$ followed by
the morphism~$g''$.

 In other words, $f'\overset z\rarrow f''$ is a $2$\+morphism in
$\Rings_{\wnlq2}$ if and only if $f'=f_{z^{-1}}f''$, or equivalently,
$f'=f''f_{w^{-1}}$.
 Similarly, $g''\overset w\rarrow g''$ is a $2$\+morphism in
$\Rings_{\cdg2,\rq}$ if and only if $g'=(g_w,a_w)\circ g''$, or
equivalently, $g'=g''\circ(g_z,a_z)$.

 It remains to refer to
Lemma~\ref{conjugation-gauge-transformation-lemma} for the assertion
that, for any object $(\tA,\tV)\in\Rings_{\wnlq2}$ and
the corresponding object $(B,d,h)\in\Rings_{\cdg2,\rq}$, basic
$2$\+morphisms $f_{z^{-1}}\overset z\rarrow\id_{(\tA,\tV)}$ in
$\Rings_{\wnlq2}$ correspond to basic $2$\+morphisms $(\id_B,0)
\overset z\rarrow (g_z,a_z)\:(B,d,h)\rarrow(B,d,h)$ in
$\Rings_{\cdg2,\rq}$.
\end{proof}

\subsection{Augmented nonhomogeneous quadratic rings}
\label{augmented-subsecn}
 Let $\tA$ be an associative ring and $R\subset\tA$ be a subring.
 A \emph{left augmentation of\/ $\tA$ over} $R$ is a left ideal
$\tA^+\subset\tA$ such that $\tA=R\oplus\tA^+$.
 Equivalently, a left augmentation is a left action of $\tA$ in $R$
extending the regular left action of $R$ in itself.

 Given a left augmentation ideal $\tA^+\subset\tA$, such a left
action of $\tA$ in $R$ is obtained by identifying $R$ with
the quotient left $\tA$\+module $\tA/\tA^+=R$.
 Conversely, given a left augmentation action of $\tA$ in $R$,
the left augmentation ideal $\tA^+\subset\tA$ is recovered as
the annihilator of the element $1\in R$.

 The \emph{category of left augmented rings over~$R$}, denoted by
$R\rings^\laug$, is defined as follows.
 The objects of $R\rings^\laug$ are associative rings $\tA$ endowed
with a subring identified with $R$ and a left augmentation ideal
$\tA^+\subset\tA$.
 Morphisms $(\prtA,\prtA^+)\rarrow(\sectA,\sectA^+)$ in
$R\rings^\laug$ are ring homomorphisms $f\:\prtA\rarrow\sectA$
forming a commutative triangle diagram with the embeddings
$R\rarrow\prtA$ and $R\rarrow\sectA$ and satisfying the condition of
compatibility with the augmentations, namely, that
$f(\prtA^+)\subset\sectA^+$.
 Equivalently, both the conditions on~$f$ can be expressed by
saying that the left action of $\prtA$ in $R$ coincides with
the action obtained from the left action of $\sectA$ in $R$ by
the restriction of scalars via~$f$.

 Moreover, one can define the \emph{$2$\+category of left augmented
rings}, denoted by $\Rings^{\laug2}$, in the following way.
 The objects of $\Rings^{\laug2}$ are associative rings $\tA$
endowed with a subring $F_0\tA\subset\tA$ and a left ideal
$\tA^+\subset\tA$ such that $\tA=F_0\tA\oplus\tA^+$.
 Morphisms $(\prtA,F_0\prtA,\prtA^+)\rarrow(\sectA,F_0\sectA,\sectA^+)$
in $\Rings^{\laug2}$ are ring homomorphisms $f\:\prtA\rarrow\sectA$
such that $f$~restricts to an isomorphism $F_0f\:F_0\prtA\rarrow
F_0\sectA$ and $f(\prtA^+)\subset\sectA^+$.
 \,$2$\+morphisms $f\overset z\rarrow g$ between a pair of parallel
morphisms $f$, $g\:(\prtA,F_0\prtA,\prtA^+)\rarrow
(\sectA,F_0\sectA,\sectA^+)$ are invertible elements $z\in F_0\prtA$
such that $g(c)=f(zcz^{-1})$ for all $c\in\prtA$.

 Notice that it follows from the latter condition that
$z\prtA^+z^{-1}=\prtA^+$, or equivalently, $\prtA^+z=\prtA^+$.
 Indeed, $zcz^{-1}=r+a$ for $c$, $a\in\prtA^+$ and $r\in F_0\prtA$
implies $g(c)=f(zcz^{-1})=f(r)+f(a)\in\sectA$ with $f(r)\in
F_0\sectA$ and $f(a)\in\sectA^+$, hence $f(r)=0$ and $r=0$.
 Similarly, $z^{-1}cz=r+a$ for $c$, $a\in\prtA^+$ and $r\in F_0\prtA$
implies $f(c)=g(z^{-1}cz)=g(r)+g(a)\in\sectA$, hence $r=0$.

 The vertical composition of two $2$\+morphisms $f'\overset w\rarrow f''
\overset z\rarrow f'''$ is the $2$\+morphism
$f'\overset{wz}\rarrow f'''$.
 The identity $2$\+morphism is the $2$\+morphism $f\overset1\rarrow f$.
 The horizontal composition of two $2$\+morphisms $g'\overset w\rarrow
g''\:(\tA,F_0\tA,\tA^+)\rarrow(\tB,F_0\tB,\tB^+)$ and $f'\overset z
\rarrow f''\:(\tB,F_0\tB,\tB^+)\rarrow(\tC,F_0\tC,\tC^+)$ is
the $2$\+morphism $f'g'\overset{z\circ w}\lrarrow f''g''\:
(\tA,F_0\tA,\tA^+)\rarrow(\tC,F_0\tC,\tC^+)$ with the element
$z\circ w=(F_0g')^{-1}(z)w=w(F_0g'')^{-1}(z)$, where $(F_0g')^{-1}$,
$(F_0g'')^{-1}\:F_0\tB\rarrow F_0\tA$ are the inverse maps to
the ring isomorphisms $F_0g'$, $F_0g''\:F_0\tA\rarrow F_0\tB$.

 All the $2$\+morphisms of left augmented rings are invertible.
 If $f\:(\prtA,F_0\prtA,\prtA^+)\rarrow(\sectA,F_0\sectA,\sectA^+)$
is a morphism of left augmented rings and $z\in F_0\prtA$ is
an invertible element such that $\prtA^+z=\prtA^+$, then
$g\:c\longmapsto f(zcz^{-1})$ is also a morphism of left augmented rings
$g\:(\prtA,F_0\prtA,\prtA^+)\rarrow(\sectA,F_0\sectA,\sectA^+)$.
 The morphisms $f$ and~$g$ are connected by the $2$\+isomorphism
$f\overset z\rarrow g$.

 Let $(\tA,F)$ be a filtered ring with an increasing filtration
$0=F_{-1}\tA\subset F_0\tA\subset F_1\tA\subset F_2\tA\subset\dotsb$
(which, as above, is presumed to be exhastive and compatible with
the multiplication in~$\tA$).
 The filtered ring $(\tA,F)$ is said to be \emph{left augmented}
if the ring $\tA$ is left augmented over its subring $F_0\tA$.
 In other words, this means that a left ideal $\tA^+\subset\tA$ is
chosen such that $\tA=F_0\tA\oplus\tA^+$.

 We denote by $R\rings_\fil^\laug$ the category of left augmented
filtered rings with the filtration component $F_0\tA$ identified
with~$R$.
 So the objects of $R\rings_\fil^\laug$ are left augmented filtered
rings $(\tA,F,\tA^+)$ for which a ring isomorphism $R\simeq
F_0\tA$ has been chosen.
 Morphisms $(\prtA,F,\prtA^+)\rarrow(\sectA,F,\sectA^+)$ in
$R\rings_\fil^\laug$ are ring homomorphisms $f\:\prtA\rarrow
\sectA$ such that $f(F_n\prtA)\subset F_n\sectA$ for all $n\ge0$,
\ $f(\prtA^+)\subset\sectA^+$, and the ring homomorphism
$F_0f\:F_0\prtA\rarrow F_0\sectA$ forms a commutative triangle
diagram with the fixed isomorphisms $R\simeq F_0\prtA$ and
$R\simeq F_0\sectA$.

 The definition of the \emph{$2$\+category of left augmented filtered
rings}, denoted by $\Rings_\fil^{\laug2}$, is similar to the above.
 The objects of $\Rings_\fil^{\laug2}$ are left augmented filtered
rings $(\tA,F,\tA^+)$.
 Morphisms $(\prtA,F,\prtA^+)\rarrow(\sectA,F,\sectA^+)$ are ring
homomorphisms $f\:\prtA\rarrow\sectA$ such that $f(F_n\prtA)\subset
F_n\sectA$ for all $n\ge0$, the map $F_0f\:F_0\prtA\rarrow F_0\sectA$
is an isomorphism, and $f(\prtA^+)\subset\sectA^+$.
 \,$2$\+morphisms $f\overset z\rarrow g$ between a pair of parallel
morphisms $f$, $g\:(\prtA,F,\prtA^+)\rarrow (\sectA,F,\sectA^+)$ are
invertible elements $z\in F_0\prtA$ such that $g(c)=f(zcz^{-1})$
for all $c\in\prtA$.
 The composition of $2$\+morphisms in $\Rings_\fil^{\laug2}$ is
defined in the same way as in the category $\Rings^{\laug2}$; so
there is an obvious forgetful strict $2$\+functor
$\Rings_\fil^{\laug2}\rarrow\Rings^{\laug2}$.

 A weak nonhomogeneous quadratic ring $R\subset\tV\subset\tA$ is said
to be \emph{left augmented} if the ring $\tA$ is endowed with a left
augmentation over its subring~$R$.
 The \emph{category of\/ $3$\+left finitely projective left
augmented weak nonhomogeneous quadratic rings over~$R$}, denoted by
$R\rings_\wnlq^\laug$, is defined as the full subcategory in
$R\rings_\fil^\laug$ whose objects are the left augmented
weak nonhomogeneous quadratic rings over $R$ that are $3$\+left
finitely projective as weak nonhomogeneous quadratic rings.

 The \emph{$2$\+category of\/ $3$\+left finitely projective left
augmented weak nonhomogeneous quadratic rings}, denoted by
$\Rings_\wnlq^{\laug2}$, is defined as the following $2$\+subcategory
in $\Rings_\fil^{\laug2}$.
 The objects of $\Rings_\fil^{\laug2}$ are the $3$\+left finitely
projective left augmented weak nonhomogeneous quadratic rings
$R\subset\tV\subset\tA\supset\tA^+$ with the filtration $F$ generated
by $F_1\tA$ over $F_0\tA$. 
 All morphisms in $\Rings_\fil^{\laug2}$ between objects of
$\Rings_\wnlq^{\laug2}$ are morphisms in $\Rings_\wnlq^{\laug2}$,
and all $2$\+morphisms in $\Rings_\fil^{\laug2}$ between morphisms of
$\Rings_\wnlq^{\laug2}$ are $2$\+morphisms in $\Rings_\wnlq^{\laug2}$.

 A \emph{DG\+ring} $(B,d)$ is a graded associative ring $B=
\bigoplus_{n\in\boZ}B^n$ endowed with an odd derivation $d\:B\rarrow B$
of degree~$1$ such that $d^2=0$.
 In other words, one can say that a DG\+ring is a CDG\+ring $(B,d,h)$
with $h=0$.
 In this section, we consider nonnegatively graded DG\+rings, that is
$B=\bigoplus_{n=0}^\infty B^n$.

 A \emph{morphism of DG\+rings} $f\:(\secB,d'')\rarrow(\prB,d')$ is
a morphism of graded rings $f\:\secB\rarrow\prB$ such that
$fd''=d'f$.
 In other words, one can say that a morphism of DG\+rings is
a morphism of CDG\+rings $(f,a)\:(\secB,d'',0)\rarrow(\prB,d',0)$
with $a=0$.
 Notice that there exist CDG\+ring morphisms $(f,a)$ with $a\ne0$
both the domain and codomain of which are DG\+rings.
 In other words, DG\+rings form a subcategory in CDG\+rings, but it
is \emph{not} a full subcategory.

 We will denote the category of nonnegatively graded DG\+rings
$(B,d)$ with the fixed degree-zero component $B^0=R$ by
$R\rings_\dg$.
 Morphisms $f\:(\secB,d'')\rarrow(\prB,d')$ in $R\rings_\dg$ are
DG\+ring morphisms such that the graded ring homomorphism
$f\:\secB\rarrow\prB$ forms a commutative triangle diagram with
the fixed isomorphisms $R\simeq\secB^0$ and $R\simeq\prB^0$.

 One can define the \emph{$2$\+category of DG\+rings} as follows.
 Let $f$, $g\:(\secB,d'')\rarrow(\prB,d')$ be a pair of parallel
morphisms of DG\+rings.
 A $2$\+morphism $f\overset z\rarrow g$ is an invertible element
$z\in\prB^0$ such that $d'(z)=0$ and $g(c)=zf(c)z^{-1}$ for
all $c\in\secB$.

 The vertical composition of two $2$\+morphisms $f'\overset w\rarrow f''
\overset z\rarrow f'''$ is the $2$\+morphism $f'\overset{zw}
\rarrow f'''$.
 The identity $2$\+morphism is the $2$\+morphism $f\overset 1\rarrow f$.
 The horizontal composition of two $2$\+morphisms $g'\overset w
\rarrow g''\:(C,d_C)\rarrow (B,d_B)$ and $f'\overset z\rarrow f''\:
(B,d_B)\rarrow(A,d_A)$ is the $2$\+morphism $f'g'\overset{z\circ w}
\lrarrow f''g''\:(C,d_C)\rarrow(A,d_A)$ with the element
$z\circ w=zf'(w)=f''(w)z\in A^0$.

 All the $2$\+morphisms of DG\+rings are invertible.
 If $f\:(\secB,d'')\rarrow(\prB,d')$ is a morphism of DG\+rings and
$z\in\prB^0$ is an invertible element such that $d'(z)=0$, then
$g\:c\longmapsto zf(c)z^{-1}$ is also a morphism of DG\+rings
$g\:(\secB,d'')\rarrow(\prB,d')$.
 The morphisms $f$ and~$g$ are connected by the $2$\+isomorphism
$f\overset z\rarrow g$.

 It is clear from these definitions that the $2$\+category of
DG\+rings is a $2$\+subcategory of the $2$\+category of CDG\+rings.
 Notice the difference, however: the $2$\+morphisms of CDG\+rings
correspond to arbitrary invertible elements $z\in\prB^0$.
 The $2$\+morphisms of DG\+rings correspond to invertible
\emph{cocycles} $z\in\prB^0$, \ $d'(z)=0$.

 The $2$\+category $\Rings_{\dg2}$ of nonnegatively graded
DG\+rings is defined as the following subcategory of the $2$\+category
of DG\+rings.
 The objects of $\Rings_{\dg2}$ are nonnegatively graded DG\+rings
$(B,d)$, \ $B=\bigoplus_{n=0}^\infty B^n$.
 Morphisms $f\:(\secB,d'')\rarrow(\prB,d')$ in $\Rings_{\dg2}$ are
morphisms of DG\+rings such that the map $f_0\:\secB^0\rarrow
\prB^0$ is an isomorphism.
 \,$2$\+morphisms $f\overset z\rarrow g$ in $\Rings_{\dg2}$ between
morphisms $f$ and $g$ belonging to $\Rings_{\dg2}$ are arbitrary
$2$\+morphisms from $f$ to~$g$ in the $2$\+category of DG\+rings.

 The \emph{category of\/ $3$\+right finitely projective quadratic
DG\+rings over~$R$}, denoted by $R\rings_{\dg,\rq}$, is the full
subcategory in the category $R\rings_\dg$ consisting of all
the DG\+rings $(B,d)$ whose underlying nonnegatively graded ring $B$
is $3$\+right finitely projective quadratic over~$R$.

 The \emph{$2$\+category of\/ $3$\+right finitely projective quadratic
DG\+rings}, denoted by $\Rings_{\dg2,\rq}$, is the similar
$2$\+subcategory in the $2$\+category $\Rings_{\dg2}$.
 The objects of $\Rings_{\dg2,\rq}$ are all the DG\+rings $(B,d)$
whose underlying nonnegatively graded ring $B$ is $3$\+right finitely
projective quadratic over~$B^0$.
 All morphisms in $\Rings_{\dg2}$ between objects of $\Rings_{\dg2,\rq}$
are morphisms in $\Rings_{\dg2,\rq}$, and all $2$\+morphisms in
$\Rings_{\dg2}$ between morphisms of $\Rings_{\dg2,\rq}$ are
$2$\+morphisms in $\Rings_{\dg2,\rq}$.

\begin{thm} \label{augmented-duality-functor-existence-thm}
 The nonhomogeneous quadratic duality functor of
Theorem~\ref{nonhomogeneous-duality-functor-existence-thm}
restricts to a fully faithful contravariant functor
\begin{equation} \label{augmented-duality-functor}
 (R\rings_\wnlq^\laug)^\sop\lrarrow R\rings_{\dg,\rq}
\end{equation}
from the category of\/ $3$\+left finitely projective left augmented
weak nonhomogeneous quadratic rings to the category of\/
$3$\+right finitely projective quadratic DG\+rings over~$R$.
\end{thm}

\begin{proof}
 It is clear from the above discussion that $R\rings_{\dg,\rq}$ is
a subcategory in $R\rings_{\cdg,\rq}$.
 Moreover, the category of $3$\+right finitely projective quadratic
DG\+rings $R\rings_{\dg,\rq}$ is a full subcategory in the category
$R\rings_{\cdg,\rq}^\sm$ of $3$\+right finitely projective quadratic
CDG\+rings over $R$ and strict morphisms between them (which was
introduced in the proof of
Theorem~\ref{nonhomogeneous-duality-functor-existence-thm}).

 Similarly, we observe that the category of $3$\+left finitely
projective left augmented weak nonhomogeneous quadratic rings
$R\rings_\wnlq^\laug$ is a full subcategory in the category
$R\rings_\wnlq^\sgsm$ of $3$\+left finitely projective
weak nonhomogeneous quadratic rings $(\tA,\tV)$ with a fixed
submodule of strict generators $V'\subset\tV$ and morphisms
$f\:(\prtA,\prtV,\prV)\rarrow(\sectA,\sectV,\secV)$ preserving
the submodule of strict generators.
 Indeed, given a $3$\+left finitely projective left augmented weak
nonhomogeneous quadratic ring $(\tA,\tV,\tA^+)$ over $R$, we
choose the left $R$\+submodule $V'=\tA^+\cap\tV\subset\tV$
as the submodule of strict generators of~$\tA$.
 The left augmentation ideal $\tA^+\subset\tA$ can be then recovered
as the left ideal (equivalently, the subring without unit)
generated by $V'$ in~$\tA$.

 Moreover, the essential image of the fully faithful functor
$R\rings_\wnlq^\laug\rarrow R\rings_\wnlq^\sgsm$ can be explicitly
described as follows.
 Given an object $(\tA,\tV,V')\in R\rings_\wnlq^\sgsm$, consider
the related maps $q\:V\times R\rarrow R$, \ $p\:\hI\rarrow R$,
and $h\:\hI\rarrow R$ defined by
the formulas~(\ref{map-q}\+-\ref{maps-p-and-h}).
 Then the object $(\tA,\tV,V')$ corresponds to a ($3$\+left
finitely projective) left augmented weak nonhomogeneous quadratic
ring if and only if one has $\hi_1*\hi_2\in V'\subset\tV$ for all
$\hi\in\hI$, that is, $h=0$.

 Indeed, the ``only if'' assertion is obvious.
 To prove the ``if'', one observes that, for any weak nonhomogeneous
quadratic ring $(\tA,\tV)$ satisfying the assumptions of
Section~\ref{self-consistency-subsecn} and any chosen submodule of
strict generators $V'\subset\tV$, the ring $\tA$ is generated
by the ring $R$ and the abelian group $V'$ with the defining
relations~(\ref{map-q}\+-\ref{maps-p-and-h}).
 It is clear from the form of these relations that the left ideal
generated by $V'$ in $\tA$ does not intersect $R$ whenever $h=0$.

 The latter condition means exactly that the quadratic CDG\+ring
$(B,d,h)$ assigned to $(\tA,\tV,V')$ by
the functor~\eqref{strict-morphisms-duality-functor} is
a DG\+ring.
 The desired fully faithful contravariant
functor~\eqref{augmented-duality-functor} can be now obtained
as a restriction of the fully faithful contravariant
functor~\eqref{strict-morphisms-duality-functor} to the full
subcategory $R\rings_\wnlq^\laug\subset R\rings_\wnlq^\sgsm$.
\end{proof}

\begin{thm} \label{augmented-duality-2-functor-existence-thm}
 The nonhomogeneous quadratic duality\/ $2$\+functor of
Theorem~\ref{nonhomogeneous-duality-2-functor-theorem} restricts
to a fully faithful strict contravariant\/ $2$\+functor
\begin{equation} \label{augmented-duality-2-functor}
 (\Rings_\wnlq^{\laug2})^\sop\lrarrow
 \Rings_{\dg2,\rq}
\end{equation}
from the\/ $2$\+category of\/ $3$\+left finitely projective left
augmented weak nonhomogeneous quadratic rings to the\/ $2$\+category
of\/ $3$\+right finitely projective quadratic DG\+rings.
\end{thm}

\begin{proof}
 Similarly to the functor~\eqref{nonhomogeneous-duality-2-functor},
the functor~\eqref{augmented-duality-2-functor} is fully faithful
in the strict sense.
 For any two objects $(\prtA,\prtV,\prtA^+)$ and
$(\sectA,\sectV,\sectA^+)\in\Rings_\wnlq^{\laug2}$ and the corresponding
DG\+rings $(\prB,d')$ and $(\secB,d'')\in\Rings_{\dg2,\rq}$,
the $2$\+functor~\eqref{augmented-duality-2-functor} induces
a bijection between morphisms $(\prtA,\prtV,\prtA^+)\rarrow
(\sectA,\sectV,\sectA^+)$ in $\Rings_\wnlq^{\laug2}$ and
morphisms $(\secB,d'')\rarrow(\prB,d')$ in $\Rings_{\dg2,\rq}$.
 Furthermore, for any pair of parallel morphisms
$f'$, $f''\:(\prtA,\prtV,\prtA^+)\rarrow(\sectA,\sectV,\sectA^+)$ in
$\Rings_\wnlq^{\laug2}$ and the corresponding pair of parallel morphisms
$g'$, $g''\:(\secB,d'')\rarrow(\prB,d')$ in $\Rings_{\dg2,\rq}$,
the $2$\+functor~\eqref{augmented-duality-2-functor} induces a bijection
between $2$\+morphisms $f'\overset z\rarrow f''$ in
$\Rings_\wnlq^{\laug2}$ and $2$\+morphisms $g''\overset z\rarrow g'$
in $\Rings_{\dg2,\rq}$.

 It is clear from the discussion above in this section that
$\Rings_{\dg2,\rq}$ is a $2$\+subcategory in $\Rings_{\cdg2,\rq}$.
 Similarly, the $2$\+category $\Rings_\wnlq^{\laug2}$ can be viewed
as a $2$\+subcategory of the $2$\+category $\Rings_{\wnlq2}$ in
the following way.
 To any object $(\tA,\tV,\tA^+)\in\Rings_\wnlq^{\laug2}$ one assigns
the object $(\tA,\tV)\in\Rings_{\wnlq2}$, and to any
morphism $f\:(\prtA,\prtV,\prtA^+)\rarrow(\sectA,\sectV,\sectA^+)$ in
$\Rings_\wnlq^{\laug2}$ one assigns the morphism
$f\:(\prtA,\prtV)\rarrow(\sectA,\sectV)$ in $\Rings_{\wnlq2}$.
 Finally, to any $2$\+morphism $f\overset z\rarrow g\:
(\prtA,\prtV,\prtA^+)\rarrow(\sectA,\sectV,\sectA^+)$ in
$\Rings_\wnlq^{\laug2}$ one assigns the $2$\+morphism
$f\overset w\rarrow g\:(\prtA,\prtV)\rarrow(\sectA,\sectV)$ in
$\Rings_{\wnlq2}$ with the element $w=(F_0f)(z)=(F_0g)(z)$.
 Here $z\in F_0\prtA$ is an invertible element such that
$\prtA^+z=\prtA^+$ and $w\in F_0\sectA$ is an invertible element,
while $F_0f$ and $F_0g\:F_0\prtA\rarrow F_0\sectA$ are two
ring isomorphisms whose values coincide on the element~$z$.

 Denote by $\Rings_\wnlq^\laug\subset\Rings_\wnlq^{\laug2}$ the category
whose objects are the objects of the $2$\+category
$\Rings_\wnlq^{\laug2}$ and whose morphisms are the morphisms of
the $2$\+category $\Rings_\wnlq^{\laug2}$ (but
there are no $2$\+morphisms in $\Rings_\wnlq^\laug$).
 Similarly, denote by $\Rings_{\dg,\rq}$ the category whose objects
are the objects of the $2$\+category $\Rings_{\dg2,\rq}$ and whose
morphisms are the morphisms of the $2$\+category $\Rings_{\dg2,\rq}$
(but there are no $2$\+morphisms in $\Rings_{\dg,\rq}$).
 Then essentially the same argument that was used to restrict
the functor~\eqref{nonhomogeneous-duality-functor}
to the functor~\eqref{augmented-duality-functor} in
Theorem~\ref{augmented-duality-functor-existence-thm} shows that
the fully faithful
functor~\eqref{enhanced-nonhomogeneous-duality-functor} restricts
to a fully faithful functor
\begin{equation} \label{enhanced-augmented-duality-functor}
 (\Rings_\wnlq^\laug)^\sop\lrarrow\Rings_{\dg,\rq}.
\end{equation}

 To deduce the existence of a fully faithful
functor~\eqref{augmented-duality-2-functor} from the existence
of the fully faithful functors~\eqref{nonhomogeneous-duality-2-functor}
and~\eqref{enhanced-augmented-duality-functor}, one can observe
that both the embeddings of $2$\+categories
$\Rings_\wnlq^{\laug2}\rarrow\Rings_{\wnlq2}$ and
$\Rings_{\dg2,\rq}\rarrow\Rings_{\cdg2,\rq}$ are \emph{fully
faithful on the level of\/ $2$\+morphisms}.
 In other words, this means that any $2$\+morphism in
$\Rings_{\wnlq2}$ between a pair of parallel morphisms in
$\Rings_\wnlq^{\laug2}$ belongs to $\Rings_\wnlq^{\laug2}$, and
similarly, any $2$\+morphism in $\Rings_{\cdg2,\rq}$ between
a pair of parallel morphisms in $\Rings_{\dg2,\rq}$ belongs
to $\Rings_{\dg2,\rq}$.
 (More generally, $\Rings_\fil^{\laug2}$ is a $2$\+subcategory in
$\Rings_{\fil2}$ such that any $2$\+morphism in $\Rings_{\fil2}$
between a pair of parallel morphisms in $\Rings_\fil^{\laug2}$
belongs to $\Rings_\fil^{\laug2}$; and $\Rings_{\dg2}$ is
a $2$\+subcategory in $\Rings_{\cdg2}$ such that any $2$\+morphism
in $\Rings_{\cdg2}$ between a pair of parallel morphisms in
$\Rings_{\dg2}$ belongs to $\Rings_{\dg2}$.)
 This suffices to prove the theorem.

 Alternatively, one can construct the fully faithful strict
$2$\+functor~\eqref{augmented-duality-2-functor} in the way similar
to the construction of the fully faithful strict
$2$\+functor~\eqref{nonhomogeneous-duality-2-functor} in the proof
of Theorem~\ref{nonhomogeneous-duality-2-functor-theorem}.
 For this purpose, one defines the ``basic $2$\+morphisms'' in
the $2$\+categories $\Rings_\wnlq^{\laug2}$ and $\Rings_{\dg2,\rq}$,
and observes that any $2$\+morphism in $\Rings_\wnlq^{\laug2}$
decomposes uniquely as a basic $2$\+morphism followed by
a morphism (but \emph{not} in the other order), while any
$2$\+morphism in $\Rings_{\dg2,\rq}$ decomposes uniquely as
a morphism followed by a basic $2$\+morphism (but \emph{not}
in the other order).

 The reason is, essentially, that the cocycle equation
$d'(z)=0$ in the definition of a $2$\+morphism of DG\+rings
$f\overset z\rarrow g\:(\secB,d'')\rarrow(\prB,d')$ does \emph{not}
imply the equation $d''(w)=0$ for the element
$w=f_0^{-1}(z)=g_0^{-1}(z)$.
 The latter equation is stronger than the former one, and does not
need to hold.
 Similarly, the condition $\prtA^+z=\prtA^+$ related to the definition
a $2$\+morphism of left augmented rings $f\overset z\rarrow g\:
(\prtA,F_0\prtA,\prtA^+)\rarrow(\sectA,F_0\sectA,\sectA^+)$ 
does \emph{not} imply the condition $\sectA^+w=\sectA^+$ for
the element $w=(F_0f)(z)=(F_0g)(z)$.
 The latter condition is stronger than the former one, and does not
need to hold.
\end{proof}

\Section{The Poincar\'e--Birkhoff--Witt Theorem} \label{pbw-secn}

\subsection{Central element theorem}
 Let $\hA=\bigoplus_{n=0}^\infty\hA_n$ be a nonnegatively graded ring
with the degree-zero component $R=\hA_0$, and let $t\in\hA_1$ be
a central element.
 Let $A=\hA/\hA t$ denote the quotient ring of $\hA$ by the homogeneous
ideal generated by~$t$.
 So $A=\bigoplus_{n=0}^\infty A_n$ is also a nonnegatively graded ring
with the degree-zero component $A_0=R$, the degree-one component
$A_1=\hA_1/Rt$, and the degree~$n$ component $A_n=\hA_n/\hA_{n-1}t$
for all $n\ge1$.
 We will say that $t$~is a \emph{nonzero-divisor} in $\hA$ if $at=0$
implies $a=0$ for any $a\in\hA_n$, \,$n\ge0$.

\begin{prop} \label{central-quotient-reflects-quadraticity}
 Let $\hA$ be a nonnegatively graded ring and $t\in\hA_1$ be a central
element.  Then \par
\textup{(a)} the graded ring\/ $\hA$ is generated by $\hA_1$ over
$\hA_0$ if and only if the graded ring $A=\hA/\hA t$ is generated by
$A_1$ over $A_0$; \par
\textup{(b)} assuming that $t$ is a nonzero-divisor in $\hA$,
the graded ring\/ $\hA$ is quadratic if and only if the graded ring
$A=\hA/\hA t$ is quadratic.
\end{prop}

\begin{proof}
 Part~(a): for any nonnegatively graded ring
$C=\bigoplus_{n=0}^\infty C_n$ generated by $C_1$ over $C_0$, and for
any homogeneous ideal $H\subset C$, the quotient ring $A=C/H$ is
generated by $A_1$ over~$A_0$.
 This proves the implication ``only if''.
 
 To prove the ``if'', suppose that we are given a nonnegatively graded
ring $C$ and a homogeneous ideal $H\subset C$ which is generated,
as a two-sided ideal, by its degree-one component~$H_1$.
 Suppose further that the quotient ring $A=C/H$ is generated by $A_1$
over~$A_0$.
 Let $C'\subset C$ denote the subring in $C$ generated by $C'_1=C_1$
over $C'_0=C_0$.
 Then the composition $C'\rarrow C\rarrow A$ is surjective, so we have
$C=C'+H$.
 Arguing by induction, we will prove that $C'_n=C_n$ for every $n\ge2$.
 Indeed, assume that $C'_k=C_k$ for all $k\le n-1$.
 Then $H_n=\sum_{k=1}^nC_{k-1}H_1C_{n-k}=
\sum_{k=1}^nC'_{k-1}H_1C'_{n-k}\subset C'_n$, hence $C_n=C'_n+H_n=C'_n$.
 Thus $C'=C$, so $C$ is generated by $C_1$ over~$C_0$.

 Part~(b): for any quadratic graded ring $C=\bigoplus_{n=0}^\infty C_n$
and any homogeneous ideal $H\subset C$ that is generated, as
a two-sided ideal, by its components $H_1$ and $H_2$, the quotient
ring $A=C/H$ is quadratic.
 This proves the implication ``only if'' (which does not depend on
the assumption that $t$~is a nonzero-divisor).

 To prove the ``if'', suppose that $t\in\hA_1$ is a central
nonzero-divisor and the graded ring $A=\hA/\hA t$ is quadratic.
 Then, by part~(a), the graded ring $\hA$ is generated by $\hA_1$
over $R=\hA_0$.

 Let $\hI\subset\hA_1\ot_R\hA_1$ be the kernel of the multiplication
map $\hA_1\ot_R\hA_1\rarrow\hA_2$ and $\q\hA=T_R(\hA_1)/(\hI)$
be the quadratic graded ring generated by $\hA_1$ with
the relations $\hI$ over~$R$.
 Then we have a unique surjective homomorphism of graded rings
$\q\hA\rarrow\hA$ acting by the identity maps on the components
of degree~$0$ and~$1$.
 By construction, the graded ring map $\q\hA\rarrow\hA$ is also
an isomorphism in degree~$2$.

 The isomorphism $\q\hA_1\simeq\hA_1$ allows to consider~$t$ as
an element of the ring $\q\hA$.
 Moreover, $t\in\q\hA_1$ is a central element, since $\q\hA$ is
generated by $\q\hA_1$ over $R$ and the relations of commutativity
of~$t$ with the elements of $R$ and $\q\hA_1$ have degree~$\le2$,
so they hold in $\q\hA$ whenever they hold in~$\hA$.

 Furthermore, by the ``only if'' assertion (which we have already
explained) the quotient ring $A'=\q\hA/(\q\hA)t$ is quadratic.
 We have the induced homomorphism of graded rings
$A'=\q\hA/(\q\hA)t\rarrow\hA/\hA t=A$.
 Since the map $\q\hA\rarrow\hA$ is an isomorphism in degree~$\le 2$,
so is the map $A'\rarrow A$.
 Since the graded ring $A$ is quadratic by assumption, it follows
by virtue of Lemma~\ref{quadratic-ring-iso} that the map
$A'\rarrow A$ is an isomorphism of graded rings.

 It follows that the kernel $H\subset\q\hA$ of the graded ring
homomorphism $\q\hA\rarrow\hA$ is contained in $(\q\hA)t\subset\q\hA$.
 Now we will prove by induction in~$n\ge3$ that $H_n=0$.
 Indeed, let $h\in H_n$ be an element.
 Then $h=h't$ for some $h'\in\q\hA_{n-1}$.
 The image of~$h$ under the ring homomorphism $\q\hA\rarrow\hA$
vanishes, and since $t$~is a nonzero-divisor in $\hA$, it follows
that the image of~$h'$ under the same homomorphism vanishes as well.
 Hence $h'\in H_{n-1}=0$ by the induction assumption and
$h=h't=0$.

 We have shown that $\q\hA\rarrow\hA$ is an isomorphism of graded
rings, and it follows that the graded ring $\hA$ is quadratic.
\end{proof}

\begin{lem} \label{hA-left-projective}
 Let $\hA$ be a nonnegatively graded ring, $t\in\hA_1$ be a central
nonzero-divisor, and $A=\hA/\hA t$ be the quotient ring.
 Let $n\ge0$ be an integer.
 Assume that $A_j$ is a finitely generated projective
(projective, or flat) left module over the ring $R=\hA_0=A_0$
for all\/ $0\le j\le n$.
 Then $\hA_j$ is a finitely generated projective (resp., projective or
flat) left $R$\+module for all\/ $0\le j\le n$.
\end{lem}

\begin{proof}
 Provable by induction in~$n$ using the short exact sequences of
$R$\+$R$\+bimodules $0\rarrow\hA_{n-1}\overset t\rarrow\hA_n
\rarrow A_n\rarrow0$.
\end{proof}

 The following theorem extends to the relative context a very
specific particular case of the result of~\cite[second assertion
of Theorem~6.1]{Pbogom}.

\begin{thm} \label{central-nonzerodivisor-theorem}
 Let $\hA$ be a nonnegatively graded ring and $t\in\hA_1$ be a central
nonzero-divisor.
 Assume that $A_n=\hA_n/\hA_{n-1}t$ is a flat left $R$\+module for
every $n\ge1$.
 Then the graded ring $\hA$ is left flat Koszul if and only if
the graded ring $A$ is left flat Koszul.
\end{thm}

\begin{proof}
 For any homomorphism of (graded) rings $C\rarrow A$, any (graded)
right $C$\+module $N$, and any (graded) left $A$\+module $M$,
the isomorphism of left derived functors of tensor product
$$
 (N\ot_C^\boL A)\ot_A^\boL M\simeq N\ot_C^\boL M
$$
on the derived categories of modules leads to a spectral sequence
of (internally graded) abelian groups
$$
 E^2_{p,q}=\Tor^A_p(\Tor^C_q(N,A),M)\,\Longrightarrow\,
 E^\infty_{p,q}=\gr^F_p\Tor^C_{p+q}(N,M)
$$
with the differentials $d^r_{p,q}\:E^r_{p,q}\rarrow E^r_{p-r,q+r-1}$.

 In particular, for any homomorphism of nonnegatively graded rings
$C\rarrow A$ acting by the identity map on their degree-zero
components $C_0=R=A_0$, we have a spectral sequence of
internally graded $R$\+$R$\+bimodules
\begin{equation} \label{change-of-rings-spectral-sequence}
 E^2_{p,q}=\Tor^A_p(\Tor^C_q(R,A),R)\,\Longrightarrow\,
 E^\infty_{p,q}=\gr^F_p\Tor^C_{p+q}(R,R).
\end{equation}

 In the situation at hand with $C=\hA$ and $A=\hA/t$, where $t\in\hA_1$
is a central nonzero-divisor, we have
$$
 \Tor^{\hA}_q(R,A)=
 \begin{cases}
 R & \text{for $q=0$}, \\
 Rt & \text{for $q=1$}, \\
 0 & \text{for $q\ge2$},
 \end{cases}
$$
where the $R$\+$R$\+bimodule $\Tor^{\hA}_0(R,A)=R$ is situated in
the internal degree~$0$ and the $R$\+$R$\+bimodule
$\Tor^{\hA}_1(R,A)=Rt$ is situated in the internal degree~$1$.

 Furthermore, the assumption that $A_n$ is a flat left $R$\+module
for every $n\ge1$ implies that $\hA_n$ is a flat left $R$\+module as
well, as one can show arguing by induction in~$n$ and using
the short exact sequences of $R$\+$R$\+bimodules $0\rarrow \hA_{n-1}
\overset t\rarrow\hA_n\rarrow A_n\rarrow0$.
 Now if $\Tor^A_{i,j}(R,R)=0$ for all $i\ne j$, then
every term $E^2_{p,q}$ of the spectral
sequence~\eqref{change-of-rings-spectral-sequence} is concentrated
in the internal degree $p+q$; hence so is the term $E^\infty_{p,q}$.
 It follows immediately that
$\Tor^{\hA}_{i,j}(R,R)=0$ for all $i\ne j$.
 So the conditions of Theorem~\ref{flat-koszul-theorem}(a) hold
for $\hA$ whenever they hold for~$A$.
 This proves the implication ``if''.

 To prove the ``only if'', one can proceed by induction in~$i$.
 Assume that the graded $R$\+$R$\+bimodule $\Tor^A_p(R,R)$ is
concentrated in the internal degree~$j=p$ for all $p\le i-1$.
 Then the terms $E^2_{p,q}$ are concentrated in the internal
degree~$p+q$ for all $p\le i-1$.
 Furthermore, if the graded $R$\+$R$\+bimodule $\Tor^{\hA}_i(R,R)$ is
concentrated in the internal degree~$i$, then so are
the $R$\+$R$\+bimodules $E^\infty_{p,q}$ for all $p+q=i$.
 In particular, the term $E^\infty_{i,0}$ is concentrated in
the internal degree~$i$.

 The only possibly nontrivial differentials passing through
$E^r_{i,0}$ with $r\ge2$ are $d^2_{i,0}\:E^2_{i,0}\rarrow E^2_{i-2,1}$.
 As the term $E^2_{i-2,1}$ is concentrated in the internal degree $i-1$
by the induction assumption and the above discussion, and the term
$E^\infty_{i,0}$ is concentrated in the internal degree~$i$,
it follows that the term $E^2_{i,0}=\Tor^A_i(R,R)$ can only have nonzero
components in the internal degrees $j=i-1$ and~$i$.
 It remains to recall that $\Tor^A_{i,j}(R,R)=0$ for $j<i$
by Proposition~\ref{diagonal-Tor}(a).
 Thus the graded $R$\+$R$\+bimodule $\Tor^A_i(R,R)$ is concentrated
in the internal degree $j=i$.
\end{proof}

 The next result is a kind of Poincar\'e--Birkhoff--Witt theorem
(see first proof of Theorem~\ref{pbw-theorem-thm}
in Section~\ref{pbw-theorem-subsecn}).

\begin{thm} \label{central-element-pbw}
 Let $\hA$ be a quadratic graded ring and $t\in\hA_1$ be a central
element.
 Assume that the left $R$\+modules $\hA_n$ are flat for all $n\ge1$
and the graded ring $A=\hA/\hA t$ is left flat Koszul.
 Assume further that the three maps $R\overset t\rarrow \hA_1
\overset t\rarrow\hA_2\overset t\rarrow\hA_3$ are injectve.
 Then the central element~$t$ is a nonzero-divisor in~$\hA$.
\end{thm}

\begin{proof}
 For any graded module $M$ over a graded ring $C$, let us denote by
$M(1)$ the same module with the shifted grading, $M(1)_n=M_{n-1}$.
 Denote by $H\subset\hA$ the kernel of the multiplication map
$\hA\overset t\rarrow\hA$.
 Then we have a four-term exact sequence of graded
$\hA$\+$\hA$\+bimodules
$$
 0\lrarrow H(1)\lrarrow\hA(1)\overset t\lrarrow\hA\lrarrow A\lrarrow0.
$$
 It follows that $\Tor_2^{\hA}(R,A)\simeq R\ot_{\hA}H(1)$ and
$\Tor_0^{\hA}(\Tor_2^{\hA}(R,A),R)\simeq R\ot_{\hA}H(1)\ot_{\hA}R$
as an internally graded $R$\+$R$\+bimodule.
 We also have $\Tor_1^{\hA}(R,A)=Rt\simeq R(1)$ and
$\Tor_0^{\hA}(R,A)=R$.

 By assumption, we have $H_n=0$ for $n\le 2$.
 Arguing by induction in $n\ge3$, we will prove that $H_n=0$ for
all~$n$.
 Assuming that $H_{n-1}=0$, the $R$\+$R$\+bimodule $H_n$ is
isomorphic to the degree~$n$ component of the $R$\+$R$\+bimodule
$R\ot_{\hA}H\ot_{\hA}R$, or which is the same, the degree $n+1$
component of the $R$\+$R$\+bimodule $R\ot_{\hA}H(1)\ot_{\hA}R$.

 We will make use of the spectral
sequence~\eqref{change-of-rings-spectral-sequence} for the graded ring
homomorphism $C=\hA\rarrow A$.
 In the induction assumption as above, we need to check that
the degree $n+1$ component of the term $E^2_{0,2}$ vanishes.
 Indeed, we have $\Tor_{2,j}^{\hA}(R,R)=0$ for $j\ge3$ by
Proposition~\ref{low-dimensional-Tor-prop}(b) (since the graded ring
$\hA$ is quadratic and its grading components are flat left
$R$\+modules by assumption), hence the term $E^\infty_{p,q}$
has no grading components in the internal degrees $j\ge3$ when $p+q=2$.
 In particular, this applies to the term $E^\infty_{0,2}$.

 The only possibly nontrivial differentials passing through
$E^r_{0,2}$ with $r\ge2$ are $d^2_{2,1}\:E^2_{2,1}\rarrow E^2_{0,2}$
and $d^3_{3,0}\:E^3_{3,0}\rarrow E^3_{0,2}$.
 Since the graded ring $A$ is left flat Koszul by assumption, we have
$\Tor_{2,j}^A(R,R)=0$ for $j\ge3$ and $\Tor_{3,j}^A(R,R)=0$ for $j\ge4$
by Theorem~\ref{flat-koszul-theorem}(a).
 Therefore, both the terms $E^2_{2,1}$ and $E^2_{3,0}$ have no
grading components in the internal degrees $j\ge4$.
 It follows that the term $E^2_{0,2}$ cannot have a nonzero grading
component in the degree $n+1\ge4$, and we are done.
\end{proof}

\subsection{Quasi-differential graded rings and CDG-rings}
\label{qdg-cdg-subsecn}
 A \emph{quasi-differential graded ring} $(\hB,\d)$ is a graded
associative ring $\hB=\bigoplus_{n\in\boZ}\hB^n$ endowed with an odd
derivation $\d\:\hB\rarrow\hB$ of degree~$-1$ such that $\d^2=0$
and the homology ring $H_\d(\hB)=\ker\d/\im\d$ vanishes.
 The latter condition holds (for an odd derivation~$\d$ of degree~$-1$
satisfying $\d^2=0$) if and only if the homology class of the unit
element $1\in\hB^0$ vanishes, that is, $1\in\d(\hB^1)$.
{\hbadness=1450\par}

 The \emph{underlying graded ring} $B$ of a quasi-differential graded
ring $(\hB,\d)$ is defined as the kernel of the differential
$\d\:\hB\rarrow\hB$, that is $B=\ker\d\subset\hB$.
 A \emph{quasi-differential structure} on a graded ring $B$ is
the datum of a quasi-differential graded ring $(\hB,\d)$ together
with a graded ring isomorphism $B\simeq\ker\d\subset\hB$.

 A \emph{morphism of quasi-differential graded rings}
$\hf\:(\sechB,\d'')\rarrow(\prhB,\d')$ is a morphism of graded rings
$\hf\:\sechB\rarrow\prhB$ such that $\hf\d''=\d'\hf$.
 The composition of morphisms of quasi-differential graded rings is
defined in the obvious way.
 These rules define the \emph{category of quasi-differential graded
rings}.

 A quasi-differential graded ring $(\hB,\d)$ is said to be
\emph{nonnegatively graded} if its underlying graded ring $B$ is
nonnegatively graded, $B=\bigoplus_{n=0}^\infty B^n$, or equivalently,
the graded ring $\hB$ is nonnegatively graded,
$\hB=\bigoplus_{n=0}^\infty\hB^n$.
 For a nonnegatively graded quasi-differential ring $(\hB,\d)$, one
has $B^0=\hB^0$.

 We will denote the category of nonnegatively graded quasi-differential
rings $(\hB,\d)$ with the fixed degree-zero component $\hB^0=B^0=R$
by $R\rings_\qdg$.
 Morphisms $\hf\:(\sechB,\d'')\rarrow(\prhB,\d')$ in $R\rings_\qdg$ are
morphisms of quasi-differential rings such that the graded ring
homomorphism $\hf\:\sechB\rarrow\prhB$ forms a commutative triangle
diagram with the fixed isomorphisms $R\simeq\sechB^0$ and
$R\simeq\prhB^0$, or equivalently, the induced homomorphism
$f\:\secB\rarrow\prB$ between the graded rings $\secB=\ker\d''
\subset\sechB$ and $\prB=\ker\d'\subset\prhB$ forms a commutative
diagram with the fixed isomorphisms $R\simeq\secB^0$ and
$R\simeq\prB^0$.

\begin{thm} \label{cdg-qdg-equivalence}
 The category of quasi-differential graded rings is equivalent to
the category of CDG\+rings.
 In particular, for any fixed ring $R$, the category $R\rings_\qdg$
is equivalent to the category $R\rings_\cdg$.
\end{thm}

\begin{proof}
 To a CDG\+ring $(B,d,h)$, one assigns the following quasi-differential
graded ring $(\hB,\d)$.
 The graded ring $\hB=B[\delta]$ is obtained by adjoining to the graded
ring $B$ a new element $\delta\in\hB^1$ satisfying the relations
\begin{equation} \label{commutator-delta-b}
 [\delta,b]=\delta b-(-1)^{|b|}b\delta=d(b)
 \quad\text{for all $b\in B$}
\end{equation}
and
\begin{equation} \label{delta-squared}
 \delta^2=h.
\end{equation}
 In other words, the elements of the grading component $\hB^n$ are all
the formal expressions $b+c\delta$ with $b\in B^n$ and $c\in B^{n-1}$.
 The unit element in $\hB$ is $1_{\hB}=1_B+0\delta$.
 The multiplication in $\hB$ is given by the formula
$$
 (b'+c'\delta)(b''+c''\delta)=(b'b''+c'd(b'')+(-1)^{|c''|}c'c''h)+
 (b'c''+(-1)^{|b''|}c'b''+c'd(c''))\delta.
$$
 The odd derivation~$\d$ on $\hB$ can be described informally as
the partial derivative $\d=\d/\d\delta$.
 Explicitly, we put $\d(b+c\delta)=(-1)^{|c|}c+0\delta$.
 So the kernel $\ker\d\subset\hB$ clearly coincides with the subring
$B\subset\hB$ embedded by the obvious rule $b\longmapsto b+0\delta$.

 Let $(\id,a)\:(B,d'',h'')\rarrow(B,d,h)$ be a change-of-connection
isomorphism of CDG\+rings; so $d''(b)=d'(b)+[a,b]$ for all $b\in B$
and $h''=h'+d'(a)+a^2$.
 Let $\sechB=B[\delta'']$ and $\prhB=B[\delta']$ denote
the quasi-differential graded rings corresponding to
$(B,d'',h'')$ and $(B,d',h')$, respectively.
 Then the rules $b\longmapsto b$ for all $b\in B$ and $\delta''
\longmapsto\delta'+a$ define an isomorphism of graded rings $\sechB
\simeq\prhB$ forming a commutative square diagram with the odd
derivations $\d''\:\sechB\rarrow\sechB$ and $\d'\:\prhB\rarrow\prhB$.

 Generally, to a morphism of CDG\+rings
$(f,a)\:(\secB,d'',h'')\rarrow(\prB,d',h')$ one assigns the morphism of
quasi-differential graded rings $\hf\:(\sechB,\d'')\rarrow(\prhB,\d')$,
where the graded ring homomorphism $\hf\:\sechB=\secB[\delta'']\rarrow
\prB[\delta']=\prhB$ takes any element $b\in\secB\subset\sechB$ to
the element $\hf(b)=f(b)\in\prB \subset\prhB$ and the element $\delta''
\in\sechB$ to the element $\hf(\delta'')=\delta'+a\in\prhB$.
 This construction defines a functor from the category of CDG\+rings
to the category of quasi-differential graded rings, and one can
easily see that this functor is fully faithful.

 To construct the inverse functor, one needs to choose, for each
quasi-differential graded ring $(\hB,\d)$, an element $\delta\in\hB^1$
such that $\d(\delta)=1$.
 Then the CDG\+ring $(B,d,h)$ is recovered by the rules $B=\ker\d
\subset\hB$, \ $d(b)=[\delta,b]$ for all $b\in B\subset\hB$,
and $h=\delta^2$.
 One has $\d([\delta,b])=[\d(\delta),b]-[\delta,\d(b)]=
[1,b]-[\delta,0]=0$ for all $b\in B$ and
$\d(\delta^2)=[\d(\delta),\delta]=[1,\delta]=0$, hence
$[\delta,b]\in B$ and $\delta^2\in B$, as desired.
 The construction of the morphism of CDG\+rings assigned to a morphism
of quasi-differential graded rings is obvious from the above.
\end{proof}

 One can define the \emph{$2$\+category of quasi-differential graded
rings} in the way similar to the definition of the $2$\+category
of DG\+rings in Section~\ref{augmented-subsecn}.
 Let $\hf$, $\hg\:(\sechB,\d'')\rarrow(\prhB,\d')$ be a pair of parallel
morphisms of quasi-differential graded rings.
 A $2$\+morphism $\hf\overset z\rarrow\hg$ is an invertible element
$z\in\prB^0=\ker(\d_0\:\prhB^0\to\prhB^{-1})$ such that
$\hg(c)=z\hf(c)z^{-1}$ for all $c\in\sechB$.
 The composition of $2$\+morphisms is defined in the same way as
in Section~\ref{augmented-subsecn}.

 The $2$\+category $\Rings_{\qdg2}$ of nonnegatively graded
quasi-differential rings is defined as the following subcategory of
the $2$\+category of quasi-differential graded rings.
 The objects of $\Rings_{\qdg2}$ are nonnegatively graded
quasi-differential rings $(\hB,\d)$.
 Morphisms $\hf\:(\sechB,\d'')\rarrow(\prhB,\d')$ in $\Rings_{\qdg2}$
are morphisms of quasi-differential graded rings such that the map
$\hf_0\:\sechB^0\rarrow\prhB^0$ is an isomorphism.
 $2$\+morphisms $\hf\overset z\rarrow\hg$ in $\Rings_{\qdg2}$ between
morphisms $\hf$ and~$\hg$ belonging to $\Rings_{\qdg2}$ are arbitrary
$2$\+morphisms from $\hf$ to~$\hg$ in the $2$\+category of
quasi-differential graded rings.

\begin{thm} \label{cdg-qdg-2cat-equivalence}
 The equivalence of categories from Theorem~\ref{cdg-qdg-equivalence}
can be extended to a strict equivalence between the\/ $2$\+category
of CDG\+rings and the\/ $2$\+category of quasi-differential graded
rings.
 In particular, the\/ $2$\+category of nonnegatively graded CDG\+rings\/
$\Rings_{\cdg2}$ is strictly equivalent to the\/ $2$\+category of
quasi-differential graded rings\/ $\Rings_{\qdg2}$.
\end{thm}

\begin{proof}
 Let $(f,a)$ and $(g,b)\:(\secB,d'',h'')\rarrow(\prB,d',h')$ be
a pair of parallel morphisms of CDG\+rings and
$\hf$, $\hg\:\sechB=\secB[\delta'']\rarrow\prB[\delta']=\prhB$ be
the corresponding pair of parallel morphisms of quasi-differential
graded rings.
 Let $z\in\prB^0$ be an invertible element.
 We leave it to the reader to check that $(f,a)\overset z\rarrow(g,b)$
is a $2$\+morphism in the $2$\+category of CDG\+rings if and only if
$\hf\overset z\rarrow\hg$ is a $2$\+morphism in the $2$\+category of
quasi-differential graded rings.
\end{proof}

\subsection{Quadratic quasi-differential graded rings}
\label{quadratic-qdg-subsecn}
 We will say that a nonnegatively graded quasi-differential ring
$(\hB,\d)$ is \emph{quadratic} if its underlying graded ring
$B=\ker\d\subset\hB$ is quadratic.

\begin{lem} \label{adding-delta-preserves-quadraticity}
 Let $(\hB,\d)$ be a nonnegatively graded quasi-differential ring
and $B=\ker\d\subset\hB$ be its underlying graded ring.
 Then \par
\textup{(a)} the graded ring $\hB$ is generated by $\hB^1$ over $\hB^0$
whenever the graded ring $B$ is generated by $B^1$ over~$B^0$; \par
\textup{(b)} the graded ring $\hB$ is quadratic whenever
the graded ring $B$ is quadratic.
\end{lem}

\begin{proof}
 Part~(a): following the proof of Theorem~\ref{cdg-qdg-equivalence},
we have $\hB=B[\delta]$, where $\delta\in\hB^1$.
 So the ring $\hB$ is generated by $B$ and~$\delta$; and if $B$ is
generated by $B^1$ over $B^0$, then it follows that $\hB$ is
generated by $\hB^1$ over $\hB^0=B^0$.

 Part~(b): the ring $\hB=B[\delta]$ is generated by the ring $B$
and the adjoined element~$\delta$ with the defining
relations~(\ref{commutator-delta-b}\+-\ref{delta-squared}).
 The relation~\eqref{delta-squared} has degree~$2$; and it remains to
observe that, whenever the ring $B$ is generated by $B^1$ over $B^0$,
it suffices to impose the relation~\eqref{commutator-delta-b} for
elements $b\in B^0$ and $B^1$ only.
 This makes the ring $\hB$ defined by relations of degree~$\le 2$
between generators of degree~$\le 1$ provided that the ring $B$
can be so defined.
\end{proof}

\begin{lem} \label{removing-delta-often-preserves-quadraticity}
 Let $(\hB,\d)$ be a nonnegatively graded quasi-differential ring
and $B=\ker\d\subset\hB$ be its underlying graded ring. \par
\textup{(a)} Assume that the graded ring $\hB$ is generated by
$\hB^1$ over~$\hB^0$.
 Then the graded ring $B$ is generated by $B^1$ over $B^0$ if and
only if the component $B^2$ is generated by $B^1$, that is,
the multiplication map $B^1\ot_{B^0}B^1\rarrow B^2$ is surjective. \par
\textup{(b)} Assume that the graded ring $B$ is generated by $B^1$
over $B^0$ and the graded ring $\hB$ is quadratic.
 Then the graded ring $B$ is quadratic if and only if it has no
relations of degree\/~$3$, that is, in other words, the natural
homomorphism of graded rings\/ $\q B\rarrow B$ is injective
in degree\/~$3$.
\end{lem}

\begin{proof}
 The ``only if'' assertion in both (a) and~(b) is obvious.
 To prove the ``if'', choose an element $\delta\in\hB^1$ such that
$\d(\delta)=1$ and consider the CDG\+ring $(B,d,h)$ corresponding
to $(\hB,\d)$ under the equivalence of categories from
Theorem~\ref{cdg-qdg-equivalence}.

 Part~(a): let $\prB\subset B$ denote the graded subring in $B$
generated by $\prB^1=B^1$ over $\prB^0=B^0$.
 By assumption, we have $\prB^2=B^2$.
 Hence the curvature element $h\in B^2$ belongs to~$\prB$.
 Moreover, we have $d(\prB^1)\subset\prB^2$, and it follows that
$d(\prB)\subset\prB\subset B$.
 Thus $(\prB,d|_{\prB},h)$ is a CDG\+ring.
 Let $\prhB=\prB[\delta]$ be the quasi-differential graded ring
corresponding to the CDG\+ring $(\prB,d|_{\prB},h)$ under
the construction of Theorem~\ref{cdg-qdg-equivalence}.
 Then $\prhB$ is naturally a graded subring in~$\hB$.
 Furthermore, we have $\prhB^0=\hB^0$ and $\prhB^1=\hB^1$, since
$\prB^0=B^0$ and $\prB^1=B^1$.
 Since the graded ring $\hB$ is generated by $\hB^1$ over $\hB^0$
by assumption, it follows that $\prhB=\hB$.
 In view of the equivalence of categories from
Theorem~\ref{cdg-qdg-equivalence}, we can conclude that $\prB=B$.

 Part~(b): the natural homomorphism of graded rings $\q B\rarrow B$,
which we will denote by~$f$, is surjective in our assumptions.
 If it is injective in degree~$3$, this means that it is an isomorphism
in degrees~$2$ and~$3$.
 So we can consider the curvature element $h\in B^2$ as an element
of $\q B^2$.

 Furthermore, the components $d_0\:B^0\rarrow B^1$ and $d_1\:B^1\rarrow
B^2$ of the odd derivation $d\:B\rarrow B$ can be considered as
maps $d'_0\:\q B^0\rarrow\q B^1$ and $d'_1\:\q B^1\rarrow\q B^2$.
 The possibility of extending these maps to an odd derivation
$d'\:\q B\rarrow\q B$ presents itself if and only if they respect
the defining relations of the graded ring $\q B$, which all live
in the degrees~$\le 2$.
 The same relations hold in the graded ring~$B$.
 Applying the maps~$d'_0$ and~$d'_1$ to the defining relations of $\q B$,
one obtains equations for elements of the components $\q B^1$, $\q B^2$,
and~$\q B^3$.
 These are isomorphic to~$B^1$, $B^2$, and $B^3$, respectively, via
the graded ring homomorphism $f\:\q B\rarrow B$.
 Since the odd derivation $d\:B\rarrow B$ exists, it follows that
one can indeed extend the maps~$d'_0$ and~$d'_1$ to an odd
derivation $d'\:\q B\rarrow \q B$.

 Finally, the equation $d'(h)=0$ holds in $\q B^3$ since the equation
$d(h)=0$ holds in~$B^3$.
 Concerning the equation $d'(d'(b))=[h,b]$ for all $b\in\q B$, it
suffices to check it for $b\in\q B^0$ and $\q B^1$, since the graded
ring $\q B$ is generated by $\q B^1$ over $\q B^0$.
 Once again, the equations need to hold in the components $\q B^2$
and $\q B^3$, which are isomorphic to $B^2$ and~$B^3$ via
the graded ring homomorphism~$f$.
 Since $(B,d,h)$ is a CDG\+ring, it follows that $(\q B,d',h)$
is a CDG\+ring as well.
 We have a strict morphism of CDG\+rings $(f,0)\:(\q B,d',h)\rarrow
(B,d,h)$.

 Let $\prhB=(\q B)[\delta]$ denote the quasi-differential graded ring
corresponding to the CDG\+ring $(\q B,d',h)$ under the construction
of Theorem~\ref{cdg-qdg-equivalence}; and let
$\hf\:\prhB\rarrow\hB$ denote the morphism of quasi-differential
graded rings corresponding to the CDG\+ring morphism $(f,0)$.
 Then the graded ring homomorphism~$\hf$ is an isomorphism in
degrees~$0$, $1$, and~$2$, since so is the graded ring
homomorphism~$f$.
 Since the graded ring $\hB$ is quadratic by assumption and
the graded ring $\prhB$ is quadratic by
Lemma~\ref{adding-delta-preserves-quadraticity}, we can conclude
that the graded ring homomorphism $\hf\:\prhB\rarrow\hB$ is
an isomorphism (in all the degrees) by Lemma~\ref{quadratic-ring-iso}.
 In view of the equivalence of categories from
Theorem~\ref{cdg-qdg-equivalence}, it follows that
the graded ring homomorphism $f\:\q B\rarrow B$ is an isomorphism, too.
 Thus the graded ring $B$ is quadratic.
\end{proof}

\begin{lem} \label{B-right-projective}
 Let $(\hB,\d)$ be a nonnegatively graded quasi-differential ring and
$B=\ker\d\subset\hB$ be its underlying graded ring.
 Let $n\ge0$ be an integer.
 Then $\hB^j$ is a finitely generated projective (projective, or flat)
right module over the ring $R=\hB^0=B^0$ for all\/ $0\le j\le n$ if and
only if $B^j$ is a finitely generated projective (resp., projective or
flat) right $R$\+module for all\/ $0\le j\le n$.
\end{lem}

\begin{proof}
 First we observe that, for any the odd derivation~$\d$ of degree~$-1$
on a nonnegatively graded ring $\hB$, the component
$\d_0\:\hB^0\rarrow\hB^{-1}$ vanishes, since $\hB^{-1}=0$.
 In other words, we have $\d(R)=0$; hence $\d$~is an $R$\+$R$\+bimodule
map.

 Returning to the situation at hand, the case $n=0$ is obvious.
 Arguing by induction in~$n$, it suffices to assume that $B^{n-1}$ is
a finitely generated projective (resp., projective or flat) right
$R$\+module, and prove that $\hB^n$ is a finitely generated projective
(resp., projective or flat) right $R$\+module if and only if $B^n$ is.
 This is clear from the short exact sequence of $R$\+$R$\+bimodules
$0\rarrow B^n\rarrow\hB^n\overset\d\rarrow B^{n-1}\rarrow0$.
\end{proof}

 The following theorem is the version of a very specific particular
case of~\cite[Corollary~6.2(b)]{Pbogom}.

\begin{thm} \label{adding-delta-preserves-Koszulity}
 Let $(\hB,\d)$ be a nonnegatively graded quasi-differential ring and
$B=\ker\d\subset\hB$ be its underlying graded ring.  Then \par
\textup{(a)} the graded ring $\hB$ is right flat Koszul whenever
the graded ring $B$ is; \par
\textup{(b)} the graded ring $\hB$ is right finitely projective Koszul
whenever the graded ring $B$ is.
\end{thm}

\begin{proof}
 In view of Lemma~\ref{B-right-projective}, it suffices to prove
part~(a).
 Then, by the same lemma, we know that $\hB^n$ is a flat right
$R$\+module for every $n\ge0$.

 In the notation of the proof of
Theorem~\ref{central-element-pbw}, we have a short exact sequence
of $R$\+$B$\+bimodules $0\rarrow B\rarrow\hB\overset\d\rarrow B(1)
\rarrow0$.
 So the graded right $B$\+module $\hB$ is projective (in fact, free).
 Hence the spectral
sequence~\eqref{change-of-rings-spectral-sequence} for the morphism
of graded rings $B\rarrow\hB$ (with the roles of the left and right
sides switched) degenerates to an isomorphism
$$
 \Tor_i^B(R,R)\simeq\Tor_i^{\hB}(R,\Tor^B_0(\hB,R))
 \quad\text{for all $i\ge0$}.
$$

 Furthermore, the internally graded $R$\+$R$\+bimodule
$\Tor^B_0(\hB,R)$ can be computed as
$$
 \hB/\hB B^+=\Tor^B_0(\hB,R)=R\oplus R(1),
$$
where $B^+=\bigoplus_{n=1}^\infty B^n$.
 The $R$\+$R$\+bimodule direct summand $R(1)\subset \Tor^B_0(\hB,R)$
is a left $\hB$\+subbimodule, since the ring $\hB$ is nonnegatively
graded.
 So we have a short exact sequence of graded $\hB$\+$R$\+bimodules
$$
 0\lrarrow R(1)\lrarrow\Tor^B_0(\hB,R)\lrarrow R\lrarrow0,
$$
where $\hB$ acts in $R(1)$ and in $R$ via the augmentation map
$B\rarrow B/B^+=R$.
 We are arriving to a long exact sequence of $R$\+$R$\+bimodules
\begin{multline*}
 \dotsb\lrarrow\Tor_{i,j-1}^{\hB}(R,R)\lrarrow\Tor^B_{i,j}(R,R)
 \\ \lrarrow\Tor^{\hB}_{i,j}(R,R)\lrarrow\Tor^{\hB}_{i-1,j-1}(R,R)
 \lrarrow\dotsb,
\end{multline*}
where the index $i\ge0$ denotes the homological grading, while
the index $j\ge0$ denotes the internal grading of
the $\Tor$ induced by the grading of $B$ and~$\hB$.

 Now if $\Tor_{i,j}^B(R,R)=0$ for $i\ne j$ and
$\Tor_{i-1,j-1}^{\hB}(R,R)=0$, then it follows that
$\Tor_{i,j}^{\hB}(R,R)=0$.
 Arguing by induction in~$i\ge0$ (or in $j\ge0$), one proves that
$\Tor_{i,j}^{\hB}(R,R)=0$ for $i\ne j$.
\end{proof}

\subsection{Central elements and acyclic derivations}
 Let $\hA=\bigoplus_{n=0}^\infty\hA_n$ be a $2$\+left finitely
projective quadratic graded ring with the degree-zero component
$R=\hA_0$, and let $\hB=\bigoplus_{n=0}^\infty\hB^n$, \ $\hB^0=R$,
be the $2$\+right finitely projective quadratic graded ring
quadratic dual to~$\hA$.
 The aim of this section is to establish and study a correspondence
between central elements $t\in\hA_1$ and odd derivations
$\d\:\hB\rarrow\hB$ of degree~$-1$.

 First of all, to every element $t\in A_1\simeq\Hom_{R^\rop}(\hB^1,R)$
we assign a right $R$\+module morphism $\d_1\:\hB^1\rarrow R$, and
vice versa.
 This is done in the obvious way except for a sign rule.
 In the pairing notation of
Section~\ref{cdg-ring-constructed-subsecn}, we put
$$
 \d_1(b) = -\lan t,b\ran \,\in\,R
 \qquad\text{for all $b\in\hB^1$}.
$$

\begin{lem} \label{central-iff-derivation}
 For any element $t\in\hA_1$, consider the related right $R$\+module
morphism $\d_1\:\hB^1\rarrow R$.
 Then \par
\textup{(a)} the element~$t$ commutes with all the elements of
$\hA_0=R$ in $\hA$ if and only if the map~$\d_1$ is
an $R$\+$R$\+bimodule morphism; \par
\textup{(b)} the element~$t$ is central in $\hA$ if and only if
the map~$\d_1$ can be (uniquely) extended to an odd derivation
$\d\:\hB\rarrow\hB$ of degree~$-1$.
\end{lem}

\begin{proof}
 Part~(a): we have $\d_1(b)=-\lan t,b\ran$ for all $b\in\hB^1$.
 Furthermore,
$$
 \lan rt,b\ran=r\lan t,b\ran
 \quad\text{and}\quad
 \lan tr,b\ran=\lan t,rb\ran
 \qquad\text{for all $r\in R$ and $b\in\hB^1$}.
$$
 Hence one has $rt=tr$ for all $r\in R$ if and only if
$r\d_1(b)=\d_1(rb)$ for all $r\in R$ and $b\in\hB^1$, i.~e.,
if and only if $\d_1$~is a left $R$\+module map.

 Part~(b): any odd derivation $\d\:\hB\rarrow\hB$ of degree~$-1$ is
an $R$\+$R$\+bimodule map, as it was explained in the proof of
Lemma~\ref{B-right-projective}.
 Moreover, the odd derivation~$\d$ is uniquely determined by
its component~$\d_1$, since the graded ring $\hB$ is generated by
$\hB^1$ over $\hB^0$.
 Similarly, the element~$t$ is central in $\hA$ if and only if it
commutes with all the elements of $\hA_0$ and $\hA_1$, since
the graded ring $\hA$ is generated by $\hA_1$ over $\hA_0$.

 So we can assume that both the equivalent conditions of part~(a)
hold, and it remains to prove that the element~$t$ commutes
with all the elements of $\hA_1$ if and only if the $R$\+$R$\+bimodule
map $\d_1\:\hB^1\rarrow\hB^0$ can be extended to an odd derivation
of degree~$-1$.
 The latter condition can be equivalently restated as follows.
 Consider the $R$\+$R$\+bimodule morphism
$$
 \tilde\d_2\:\hB^1\ot_R\hB^1\lrarrow\hB^1
$$
defined by the formula
$$
 \tilde\d_2(b_1\ot b_2)=\d_1(b_1)b_2-b_1\d_1(b_2)
 \quad\text{for all $b_1$, $b_2\in\hB^1$}.
$$
 For any $R$\+$R$\+bimodule map $\d_1\:\hB^1\rarrow R$,
the corresponding map~$\tilde\d_2$ is well-defined.
 Since the graded ring $\hB$ is quadratic, the map~$\d_1$ extends
to an odd derivation of $\hB$ if and only if the map~$\tilde\d_2$
annihilates the kernel $\hJ\subset\hB^1\ot_R\hB^1$ of
the surjective multiplication map $\hB^1\ot_R\hB^1\rarrow\hB^2$.
 This is provable using Lemma~\ref{odd-derivations-lemma}.

 On the other hand, the condition that the element~$t$ commutes
with all the elements of $\hA_1$ can be restated as follows.
 Consider the $R$\+$R$\+bimodule morphism
$$
 \tilde t\:\hA_1\lrarrow\hA_1\ot_R\hA_1
$$
defined by the formula
$$
 \tilde t(a)=a\ot t-t\ot a.
$$
 Then the element $t\in\hA_1$ commutes with all the elements of
$\hA_1$ if and only if the composition of the map~$\tilde t$ with
the surjective multiplication map $\hA_1\ot_R\hA_1\rarrow\hA_2$
vanishes.
 In order to prove the desired equivalence, it remains to recall
and observe that the contravariant functor $\Hom_R({-},R)$ transforms
the $R$\+$R$\+bimodule $\hA_1$ into the $R$\+$R$\+bimodule
$\hB^1$, the $R$\+$R$\+bimodule $\hA_1\ot_R\hA_1$ into
the $R$\+$R$\+bimodule $\hB^1\ot_R\hB^1$, the map~$\tilde t$ into
the map~$-\tilde\d_2$, the $R$\+$R$\+bimodule $\hA_2$ into
the $R$\+$R$\+bimodule $\hJ$, and the surjection $\hA_1\ot_R\hA_1\rarrow
\hA_2$ into the inclusion $\hJ\rarrow\hB^1\ot_R\hB^1$.
 Besides, the $R$\+$R$\+bimodules $\hA_1$, \ $\hA_1\ot_R\hA_1$, and
$\hA_2$ are finitely generated and projective as left $R$\+modules,
while the $R$\+$R$\+bimodules $\hB^1$, \ $\hB^1\ot_R\hB^1$, and
$\hJ$ are finitely generated and projective as right $R$\+modules;
so the contravariant functor $\Hom_{R^\rop}({-},R)$ performs
the inverse transformation.
\end{proof}

\begin{lem} \label{nonzerodivisor-in-degree-one}
 Let $\d\:\hB\rarrow\hB$ be an odd derivation of degree~$-1$.  Then \par
\textup{(a)} $\d^2=0$; \par
\textup{(b)} the homology ring $H_\d(\hB)$ vanishes if and only if,
for the central element $t\in\hA_1$ corresponding to~$\d$,
the multiplication map $\hA_0\overset t\rarrow\hA_1$ is injective
and its cokernel $A_1=\hA_1/\hA_0 t$ is a projective left $R$\+module.
\end{lem}

\begin{proof}
 Part~(a): the compositions $\d_{-1}\d_0\:\hB^0\rarrow\hB^{-1}\rarrow
\hB^{-2}$ and $\d_0\d_1\:\hB^1\rarrow\hB^0\rarrow\hB^{-1}$ vanish,
since $\hB^{-1}=0$.
 Since the square of an odd derivation is a derivation and the graded
ring $\hB$ is generated by $\hB^1$ over $\hB^0$, it follows that
$\d^2=0$ on the whole graded ring~$\hB$.
 Part~(b) is a particular case of the assertion that a morphism of
finitely generated projective left $R$\+modules $t\:\hA_0\rarrow\hA_1$ is
injective with a projective cokernel if and only if the dual morphism
of finitely generated projective right $R$\+modules $-\d_1\:\hB^1\rarrow
\hB^0$ is surjective.
\end{proof}

\begin{prop} \label{nonzerodivisor-in-degree-two}
 Let $\hA$ be a\/ $2$\+left finitely projective quadratic graded ring
and $\hB$ be the quadratic dual\/ $2$\+right finitely projective
quadratic graded ring.
 Let $t\in\hA_1$ be a central element and $\d\:\hB\rarrow\hB$ be
the corresponding odd derivation of degree~$-1$
(see Lemma~\ref{central-iff-derivation}(b)).
 Assume that the equivalent conditions of
Lemma~\ref{nonzerodivisor-in-degree-one}(b) are satisfied.
 Let $B=\ker\d\subset\hB$ be the underlying graded ring of
the quasi-differential graded ring $(\hB,\d)$, and let $A=\hA/\hA t$
be the quotient ring.
 Then \par
\textup{(a)} the graded ring $B$ is generated by $B^1$ over $B^0$ if
and only if the multiplication map $\hA_1\overset t\rarrow\hA_2$ is
injective and its cokernel $A_2=\hA_2/\hA_1 t$ is a projective left
$R$\+module; \par
\textup{(b)} if the equivalent conditions of part~(a) hold, then
the\/ $2$\+right finitely projective quadratic graded ring\/ $\q B$ is
quadratic dual to the\/ $2$\+left finitely projective quadratic
graded ring~$A$.
 The composition of graded ring homomorphisms\/
$\q B\rarrow B\rarrow\hB$ is the morphism of\/ $2$\+right finitely
projective quadratic graded rings corresponding to the surjective
morphism of\/ $2$\+left finitely projective quadratic graded rings
$\hA\rarrow A$ under the equivalence of categories from
Proposition~\ref{2-fin-proj-quadratic-duality}.
\end{prop}

\begin{proof}
 Part~(a): according to
Lemma~\ref{removing-delta-often-preserves-quadraticity}(a),
the graded ring $B$ is generated by $B^1$ over $B^0$ if and only if
the multiplication map $B^1\ot_RB^1\rarrow B^2$ is surjective.
 In the notation of the proof of Lemma~\ref{central-iff-derivation}(b),
we have
$$
 B^1\ot_RB^1\,\subset\,\ker\tilde\d_2\,\subset\,\hB^1\ot_R\hB^1.
$$
 Here the map $B^1\ot_RB^1\rarrow\hB^1\ot_R\hB^1$ induced by
the inclusion of the $R$\+$R$\+sub\-bi\-mod\-ule $B^1$ into
the $R$\+$R$\+bimodule $\hB^1$ is injective, because the right
$R$\+module $\hB^1$ is projective and the quotient bimodule
$\hB^1/B^1$, being naturally isomorphic to $R$ via the surjective
$R$\+$R$\+bimodule map $\d_1\:\hB^1\rarrow\hB^0$, is also projective
as a right $R$\+module.
 The inclusion $B^1\ot_R B^1\subset\ker\tilde\d_2$ is clear from
the construction of the map~$\tilde\d_2$.

 Furthermore, following the proof of
Lemma~\ref{central-iff-derivation}(b), we have
$\hJ\subset\ker\tilde\d_2$, since $\d$~is an odd derivation of
$\hB$ by assumption.
 The two $R$\+$R$\+bimodule maps
$$
 \hB^1\ot_R\hB^1\rightrightarrows\hB^1, \qquad
 b_1\ot_R b_2\,\longmapsto\, \d(b_1)b_2\, \text{ and } \,b_1\d(b_2)
$$
restrict to one and the same map $u\:\ker\tilde\d_2\rarrow\hB_1$.
 The subbimodule $B^1\ot_RB^1\subset\ker\tilde\d_2$ is the kernel
of the map~$u$.

 The map $\tilde\d_2$ is equal to the composition of the surjective
multiplication map $\hB^1\ot_R\hB^1\rarrow\hB^2$ with the map
$\d_2\:\hB^2\rarrow\hB^1$.
 By the definition, $B^2\subset\hB^2$ is the kernel of the latter map.
 It follows that the multiplication map $B^1\ot_RB^1\rarrow B^2$ is
surjective if and only if one has $\ker\tilde\d_2=\hJ+B^1\ot_RB^1$.
 The latter condition holds if and only if the map
$$
 \bar u=u|_{\hJ}\:\hJ\rarrow\hB^1
$$
is surjective.

 It remains to observe that the functor $\Hom_R({-},R)$ transforms
the map $\hA_1\overset t\rarrow\hA_2$ into the map~$-\bar u$.
 Thus the map~$\bar u$ is surjective if and only if the map~$t$
is injective and its cokernel $A_2$ is a projective left $R$\+module.

 In part~(b), the graded ring $A$ is quadratic by (the proof of)
Proposition~\ref{central-quotient-reflects-quadraticity}(b) and
$2$\+left finitely projective by assumptions; while the quadratic
graded ring $\q B$ is $2$\+right finitely projective by
Lemma~\ref{B-right-projective}.
 Furthermore, whenever the equivalent conditions of
Lemma~\ref{nonzerodivisor-in-degree-one}(b) hold, the functor
$\Hom_R({-},R)$ transforms the $R$\+$R$\+bimodule $A_1$ into
the $R$\+$R$\+bimodule~$B^1$.
 This is clear from the proof of 
Lemma~\ref{nonzerodivisor-in-degree-one}(b).
 Finally, whenever the equivalent conditions of part~(a) of
the present proposition hold, the functor $\Hom_R({-},R)$ transforms
the $R$\+$R$\+bimodule $A_2$ into the $R$\+$R$\+bimodule
$$
 \ker\bar u =\hJ\cap(B^1\ot_RB^1)\subset\hB^1\ot_R\hB^1.
$$
 The latter $R$\+$R$\+bimodule is the kernel of the surjective
multiplication map $B^1\ot_RB^1\rarrow B^2$.
 Then one has to check that the same functor transforms the surjective
multiplication map $A_1\ot_RA_1\rarrow A_2$ into the inclusion map
$\ker\bar u\rarrow B^1\ot_RB^1$.
\end{proof}

\begin{prop} \label{nonzerodivisor-in-degree-three}
 Let $\hA$ be a\/ $3$\+left finitely projective quadratic graded ring
and $\hB$ be the quadratic dual\/ $3$\+right finitely projective
quadratic graded ring.
 Let $t\in\hA_1$ be a central element and $\d\:\hB\rarrow\hB$ be
the corresponding odd derivation of degree~$-1$.
 Assume that the equivalent conditions of
Lemma~\ref{nonzerodivisor-in-degree-one}(b) are satisfied and
the equivalent conditions of
Proposition~\ref{nonzerodivisor-in-degree-two}(a) are satisfied
as well.
 Then the graded ring $B$ is quadratic if and only if
the multiplication map $\hA_2\overset t\rarrow\hA_3$ is injective
and its cokernel $A_3=\hA_3/\hA_2 t$ is a projective left $R$\+module.
\end{prop}

\begin{proof}
 By Lemma~\ref{B-right-projective}, the right $R$\+module $B^3$ is
projective in our assumptions.
 So, if the graded ring $B$ is quadratic, then it is $3$\+right
finitely projective quadratic.
 Since $\q B$ is quadratic dual to $A$ by
Proposition~\ref{nonzerodivisor-in-degree-two}(b), it then follows
by virtue of Proposition~\ref{3-fin-proj-quadratic-duality} that
the quadratic graded ring $A$ is $3$\+left finitely projective.
 Thus we can assume that $A_3$ is a projective left $R$\+module in
all cases.

 Then, again by Propositions~\ref{nonzerodivisor-in-degree-two}(b)
and~\ref{3-fin-proj-quadratic-duality}, it follows that the quadratic
graded ring $\q B$ is $3$\+right finitely projective.
 According to
Lemma~\ref{removing-delta-often-preserves-quadraticity}(b),
the graded ring $B$ is quadratic if and only if the surjective
homomorphism of graded rings $\q B\rarrow B$ is an isomorphism
in degree~$3$.
 Equivalently, this means that the composition $\q B^3\rarrow
B^3\rarrow\hB^3$ is an injective map.
 Furthermore, since both $\q B^3$ and $B^3$ are projective right
$R$\+modules, and since the right $R$\+module $\hB^3/B^3\simeq B^2$
is projective as well, we can conclude that the map $\q B^3\rarrow B^3$
is an isomorphism if and only if applying the functor
$\Hom_{R^\rop}({-},R)$ to the composition $\q B^3\rarrow B^3
\rarrow\hB^3$ produces a surjective map
$f\:\Hom_{R^\rop}(\hB^3,R)\rarrow\Hom_{R^\rop}(\q B^3,R)$.

 Denote by $I\subset A_1\ot_RA_1$ and $\hI\subset\hA_1\ot_R\hA_1$
the kernels of the multiplication maps $A_1\ot_RA_1\rarrow A_1$
and $\hA_1\ot_R\hA_1\rarrow\hA_1$.
 Both the graded rings $A$ and $\hA$ are $3$\+left finitely projective
quadratic in our assumptions; $\q B$ and $\hB$, respectively, are
their quadratic dual $3$\+right finitely projective quadratic rings.
 Put
$$
 I^{(3)}=(I\ot_R A_1)\cap(A_1\ot_R I)\subset A_1\ot_RA_1\ot_RA_1
$$
and
$$
 \hI^{(3)}=(\hI\ot_R\hA_1)\cap(\hA_1\ot_R\hI)\subset
 \hA_1\ot_R\hA_1\ot_R\hA_1.
$$
 Following the proof of Proposition~\ref{3-fin-proj-quadratic-duality},
we have natural isomorphisms of $R$\+$R$\+bimod\-ules
$\Hom_{R^\rop}(\q B^3,R)\simeq I^{(3)}$ and $\Hom_{R^\rop}(\hB^3,R)
\simeq\hI^{(3)}$.
 The map
$$
 f\:\hI^{(3)}=\Hom_{R^\rop}(\hB^3,R)\lrarrow
 \Hom_{R^\rop}(\q B^3,R)=I^{(3)}
$$
that we are interested in is induced by the surjective morphism of
quadratic graded rings $\hA\rarrow A$ (which is quadratic dual to
the morphism of quadratic graded rings $\q B\rarrow\hB$ according to
Proposition~\ref{nonzerodivisor-in-degree-two}(b)).

 We have shown that the graded ring $B$ is quadratic if and only if
the natural map $f\:\hI^{(3)}\rarrow I^{(3)}$ is surjective.
 Let us show that the latter condition holds if and only if
the multiplication map $\hA_2\overset t\rarrow\hA_3$ is injective.

 Consider the three-term complex
$$
 \hA_1\ot_R\hA_1\ot_R\hA_1\lrarrow
 (\hA_1\ot_R\hA_2)\oplus(\hA_2\ot_R\hA_1)\lrarrow\hA_3\lrarrow0
$$
and denote it by~$\hC_\bu$.
 For the sake of certainty of notation, let us place the complex
$\hC_\bu$ in the homological degrees~$1$, $2$, and~$3$, so that
$\hC_3=\hA_1^{\ot_R\,3}$ and $\hC_1=\hA_3$.
 Endow the graded ring~$\hA$ with a multiplicative decreasing
filtration by homogeneous ideals $\hA=G^0\hA\supset G^1\hA\supset
G^2\hA\supset\dotsb$ defined by the rule $G^p\hA=t^p\hA$, and endow
the complex $\hC_\bu$ with the induced decreasing filtration.
 So, in particular, one has $\hA/G^1\hA=A$, and the complex
$\hC_\bu/G^1\hC_\bu$ is isomorphic to the complex
$$
 A_1\ot_R A_1\ot_R A_1\lrarrow (A_1\ot_RA_2)\oplus(A_2\ot_RA_1)
 \lrarrow A_3\lrarrow0,
$$
which we denote by~$C_\bu=C_\bu(A)$.

 Since the graded rings $\hA$ and $A$ are quadratic, we have
$H_i(\hC_\bu)=0=H_i(C_\bu)$ for all $i\ne 3$, while
$H_3(\hC_\bu)=\hI^{(3)}$ and $H_3(C_\bu)=I^{(3)}$.
 It is clear from the homological long exact sequence related
to the short exact sequence of complexes
\begin{equation} \label{G-1-sequence}
 0\lrarrow G^1\hC_\bu\lrarrow\hC_\bu\lrarrow C_\bu\lrarrow0
\end{equation}
that the natural map $f\:\hI^{(3)}\rarrow I^{(3)}$ is surjective
if and only if $H_2(G^1\hC_\bu)=0$.

 Consider the nonnegatively graded ring $\boZ[\bar t]$ of polynomials
in one variable~$\bar t$ of degree~$1$ with integer coefficients.
 Let $A[\bar t]=A\ot_\boZ\boZ[\bar t]$ be the tensor product of
the graded rings $A$ and $\boZ[\bar t]$, taken over the ring $\boZ$
and endowed with the induced grading.
 Obviously, the graded ring $A[\bar t]$ is quadratic.
 Furthermore, there is a natural morphism of graded rings
\begin{equation} \label{bar-t-gr-G-map}
 A[\bar t]\lrarrow\gr_G\hA =
 \bigoplus\nolimits_{p=0}^\infty G^p\hA/G^{p+1}\hA
\end{equation}
whose restriction to $A\subset A[\bar t]$ is equal to the inclusion
$A\simeq \hA/G^1\hA\hookrightarrow\gr_G\hA$ and which takes
the element $\bar t\in A[\bar t]$ to the coset $t+G^2\hA\in
G^1\hA/G^2\hA$.

 One easily observes that the map~\eqref{bar-t-gr-G-map} is surjective.
 Furthermore, the map $A[\bar t]_n\rarrow(\gr_G\hA)_n$ is injective
(equivalently, an isomorphism) if and only if the map $\hA_{j-1}\overset
t\rarrow\hA_j$ is injective for all $j\le n$.
 In our assumptions, we know that this injectivity holds in
the internal degrees $j=1$ and~$2$, and we are interested in knowing
whether it holds for $j=3$.
 So the map $\hA_2\overset t\rarrow\hA_3$ is injective if and only if
the map $A[\bar t]_3\rarrow \gr_G\hA_3$ is.

 Put $\tA=A[\bar t]$ and consider the three-term complex
$$
 \tA_1\ot_R\tA_1\ot_R\tA_1\lrarrow
 (\tA_1\ot_R\tA_2)\oplus(\tA_2\ot_R\tA_1)\lrarrow\tA_3\lrarrow0,
$$
denoted by $\tC_\bu=C_\bu(\tA)$.
 Since the graded ring $A$ is $2$\+left finitely projective,
the complex $\gr_G\hC_\bu$ is isomorphic to the complex
$C_\bu(\gr_G\hA)$.
 The surjective morphism of graded rings $\tA\rarrow\gr_G\hA$ induces
a surjective morphism of complexes $\tC_\bu\rarrow\gr_G\hC_\bu$.

 Moreover, the maps $\tC_i\rarrow\gr_G\hC_i$ are isomorphisms for
$i\ne 1$, since the map $\tA_j\rarrow\gr_G\hA_j$ is an isomorphism
for $j<3$.
 Furthermore, the complex $\tC_i$ is acyclic in the homological
degrees $i\ne 3$, since the graded ring $\tA=A[\bar t]$ is quadratic. 
 Since the map $\tC_3\rarrow\gr_G\hC_3$ is surjective, it follows
that the map $A[\bar t]_3\rarrow\gr_G\hA_3$ is an isomorphism if and
only if $H_2(\gr_G\hC_\bu)=0$.

 It remains to show that $H_2(G^1\hC_\bu)=0$ if and only if
$H_2(\gr_G\hC_\bu)=0$.
 The implication ``if'' is obvious.
 To prove the ``only if'', we assume that $H_2(G^1\hC_\bu)=0$, the map
$f\:\hI^{(3)}\rarrow I^{(3)}$ is surjective,
and the graded ring $B$ is quadratic.

 Consider the additional grading~$p$ on the ring $\gr_G\hA$ and
the complex $\gr_G\hC_\bu$ induced by the indexing of
the filtration~$G$.
 The related grading~$p$ on the ring $\tA=A[\bar t]$ and the complex
$\tC_\bu=C_\bu(\tA)$ is induced by the grading of
the ring $\boZ[\bar t]$.
 The additional grading~$p$ on the rings $\gr_G\hA$ and $\tA$
takes values in the monoid of nonnegative integers $p\ge0$.
 On the complexes $\gr_G\hC_\bu$ and $\tC_\bu$, the additional
grading takes values $0\le p\le 3$.

 In the additional grading $p=0$, the morphism of graded rings
$A[\bar t]\rarrow\gr_G\hA$ is an isomorphism, since $A=\hA/\hA t$.
 Therefore, the morphism of complexes $\tC_\bu\rarrow\gr_G\hC_\bu$ is
an isomorphism in the additional grading $p=0$.
 Any possible homology classes in $H_2(\gr_G\hC_\bu)$ would occur in
the additional grading $p=1$, $2$, or~$3$.

 We have $H_1(\gr_G\hC_\bu)=0$, as the graded ring $\gr_G\hA$ is
generated by its component $A\oplus R[\bar t]$ of degree $n=1$ over
its component $R$ of degree $n=0$.
 Let us compute the $p$\+graded abelian group ($R$\+$R$\+bimodule)
$\tI^{(3)}=H_3(\tC_\bu)=H_3(\gr_G\hC_\bu)$.

 The quadratic graded ring $\tA=A[\bar t]$ (in the grading~$n$) is
$3$\+left finitely projective, and its quadratic dual $3$\+right
finitely projective graded ring $\tB$ can be computed as
the graded ring $\tB=B[\bar\delta]=B\ot_\boZ^{-1}\boZ[\bar\delta]$
obtained by adjoining to $B$ a generator~$\bar\delta$ of degree $n=1$
(and $p=1$) with the relations $\bar\delta b+b\bar\delta=0$
for all $b\in B^1$ and $\bar\delta^2=0$.
 Here the notation $B\ot_\boZ^q\boZ[\bar\delta]$ stands for
the relations $\bar\delta b=qb\bar\delta$, \,$b\in B^1$,
with $q=-1$ \cite[Section~1 of Chapter~3]{PP}.
 Consequently, $\tB^3=B^3\oplus B^2\bar\delta$ and
$\tI^{(3)}\simeq\Hom_{R^\rop}(\tB^3,R)=I^{(3)}\oplus\bar\epsilon I$,
where $\bar\epsilon$ is a variable dual to~$\bar\delta$.
 The direct summand $I^{(3)}$ of the $p$\+graded group $\tI^{(3)}$
sits in the additional grading $p=0$ and the direct summand
$\bar\epsilon I\simeq I$ sits in the additional grading $p=1$.

 Thus we have
$$
 H_3(G^p\hC_\bu/G^{p+1}\hC_\bu)=
 \begin{cases}
  I^{(3)} & \text{for $p=0$} \\ 
  \bar\epsilon I &\text{for $p=1$} \\
  0 & \text{for $p\ge2$}.
 \end{cases}
$$
 In view of the spectral sequence connecting the homology groups
of the complexes $G^1\hC_\bu$ and $\gr_G\hC^\bu$, it remains to
check that the map
\begin{equation} \label{map-from-kernel}
 H_3(G^1\hC_\bu)\lrarrow
 H_3(G^1\hC_\bu/G^2\hC_\bu)=\bar\epsilon I
\end{equation}
is surjective (or equivalently, an isomorphism).
 Here the group $H_3(G^1\hC_\bu)$ is the kernel of the natural map
$f\:\hI^{(3)}\rarrow I^{(3)}$ (see short exact
sequence~\eqref{G-1-sequence}).

 Following the discussion in the beginning of this proof, the map
$f\:\hI^{(3)}\rarrow I^{(3)}$ can be obtained by applying the functor
$\Hom_{R^\rop}({-},R)$ to the composition of maps $\q B^3\rarrow
B^3\rarrow\hB^3$.
 In our present assumptions, the map $f$~is surjective and the graded
ring $B$ is quadratic.
 The kernel of the map~$f$ is
$\ker(f)\simeq\Hom_{R^{\rop}}(\hB^3/B^3,R)=
\Hom_{R^{\rop}}(B^2\bar\delta,R)\simeq\bar\epsilon I$, as desired.
 Notice that existence of a well-defined odd derivation $\d=\d/\d\delta$
on $\hB$ ensures injectivity (equivalently, bijectivity) of
the multiplication map $B^2\overset{\bar\delta}\rarrow\hB^3/B^3$.
\end{proof}

\subsection{Nonhomogeneous quadratic duality via quasi-differential
graded rings}
 The construction of the nonhomogeneous quadratic duality functors
$$
 (R\rings_\wnlq)^\sop\lrarrow R\rings_{\cdg,\rq}
 \quad\text{and}\quad
 (\Rings_{\wnlq2})^\sop\lrarrow\Rings_{\cdg2,\rq}
$$
in Theorems~\ref{nonhomogeneous-duality-functor-existence-thm}
and~\ref{nonhomogeneous-duality-2-functor-theorem} are based on
the computations in Sections~\ref{self-consistency-subsecn}\+-%
\ref{change-of-strict-gens-subsecn}, which are beautiful, but quite
involved.
 The definitions and results above in Section~\ref{pbw-secn} allow to
produce the nonhomogeneous duality functors in a more conceptual fashion.

 A nonnegatively graded quasi-differential ring $(\hB,\d)$ is said to be
\emph{$3$\+right finitely projective quadratic} if its underlying graded
ring $B=\ker\d\subset\hB$ is $3$\+right finitely projective quadratic.
 The results of Section~\ref{quadratic-qdg-subsecn} explain how such
properties of the graded ring $B$ are related to the similar properties
of the graded ring~$\hB$.

 Let $R$ be an associative ring.
 The \emph{category of\/ $3$\+right finitely projective quadratic
quasi-differential graded rings over~$R$}, denoted by
$R\rings_{\qdg,\rq}$, is the full subcategory in the category
$R\rings_\qdg$ (as defined in Section~\ref{qdg-cdg-subsecn}) consisting
of all the $3$\+right finitely projective quadratic
quasi-differential graded rings.

 Furthermore, the \emph{$2$\+category of\/ $3$\+right finitely
projective quadratic quasi-differential graded rings}, denoted by
$\Rings_{\qdg2,\rq}$, is the following $2$\+subcategory in
the $2$\+category $\Rings_{\qdg2}$ (which was also defined in
Section~\ref{qdg-cdg-subsecn}).
 The objects of $\Rings_{\qdg2,\rq}$ are all the $3$\+right finitely
projective quadratic quasi-differential graded rings.
 All morphisms in $\Rings_{\qdg2}$ between objects of
$\Rings_{\qdg2,\rq}$ are morphisms in $\Rings_{\qdg2,\rq}$, and
all $2$\+morphisms in $\Rings_{\qdg2}$ between morphisms of
$\Rings_{\qdg2,\rq}$ are $2$\+morphisms in $\Rings_{\qdg2,\rq}$.
{\hbadness=2700\par}

 Restricting the equivalence of categories $R\rings_\qdg\simeq
R\rings_\cdg$ provided by Theorem~\ref{cdg-qdg-equivalence} to
the full subcategories of $3$\+right finitely projective quadratic
CDG\+rings and quasi-differential graded rings over $R$, we obtain
an equivalence of categories
\begin{equation} \label{quadratic-qdg-cdg-equivalence}
 R\rings_{\qdg,\rq}\,\simeq\,R\rings_{\cdg,\rq}.
\end{equation}
 Similarly, a restriction of the strict equivalence of $2$\+categories
$\Rings_{\qdg2}\,\simeq\Rings_{\cdg2}$ provided by
Theorem~\ref{cdg-qdg-2cat-equivalence} produces a strict equvalence
between the $2$\+categories of $3$\+right finitely projective
quadratic CDG\+rings and quasi-differential graded rings,
\begin{equation} \label{quadratic-qdg-cdg-2cat-equivalence}
 \Rings_{\qdg2,\rq}\,\simeq\,\Rings_{\cdg2,\rq}.
\end{equation}

 Let $(\tA,F)$ be a filtered ring with an increasing filtration
$0=F_{-1}\tA\subset F_0\tA\subset F_1\tA\subset F_2\tA\subset\dotsb$.
 Consider two graded rings related to such a filtration: the associated
graded ring $A=\gr^F\tA$ and the Rees ring
$\hA=\bigoplus_{n=0}^\infty F_n\tA$.

 The unit element $1\in F_0\tA$, viewed as an element of $F_1\tA$,
represents a central nonzero-divisor $t\in\hA_1$.
 The quotient ring $\hA/\hA t$ is naturally isomorphic to the associated
graded ring~$A$ of the filtration $F$ on the ring~$\tA$.
 By Proposition~\ref{central-quotient-reflects-quadraticity}, the graded
ring $\hA$ is generated by $\hA_1$ over $\hA_0$ if and only if
the graded ring $A$ is generated by $A_1$ over $A_0$; and moreover,
the graded ring $\hA$ is quadratic if and only if the graded ring $A$ is.

 Let $\tA$ be a weak nonhomogeneous quadratic ring over a subring
$R\subset\tA$ with the $R$\+$R$\+bimodule of generators $\tV\subset\tA$.
 Let $F$ be the related increasing filtration on the ring $\tA$,
as defined in Section~\ref{nonhomogeneous-quadratic-subsecn}.
 Then the graded ring $\hA$ is generated by $\hA_1$ over $\hA_0$, and
the graded ring $A$ is generated by $A_1$ over~$A_0$.
 The graded rings $\hA$ and $A$ do \emph{not} need to be quadratic (as
we only assume $\tA$ be a \emph{weak} nonhomogeneous quadratic ring);
so let us consider the quadratic graded rings $\q\hA$ and~$\q A$.

\begin{lem} \label{quadratic-part-central-element}
 Let $\hA$ be a nonnegatively graded ring, and let $t\in\hA$ be a central
element.
 Let $A=\hA/\hA t$ be the quotient ring.
 Then \par
\textup{(a)} $t\in\q\hA_1=\hA_1$ is a central element in the quadratic
graded ring\/ $\q\hA$; \par
\textup{(b)} for each $n=1$, $2$, or\/~$3$, the multiplication map\/
$\q\hA_{n-1}\overset t\rarrow\q\hA_n$ is injective whenever
the multipication map $\hA_{n-1}\overset t\rarrow\hA_n$ is; \par
\textup{(c)} there is a natural isomorphism of quadratic graded rings\/
$\q\hA/(\q\hA)t\simeq\q A$.
\end{lem}

\begin{proof}
 In part~(a), the element $t\in\q\hA_1$ commutes with all the elements
of $\q\hA_0$ and $\q\hA_1$, since the graded ring homomorphism
$\q\hA\rarrow\hA$ is an isomorphism in degree $n\le 1$ and a monomorphism
in degree $n=2$, and the element $t\in\hA_1$ is central in~$\hA$.
 Since the graded ring $\q\hA$ is generated by $\q\hA_1$ over $\q\hA_0$,
it follows that the element $t\in\q\hA_1$ is central in~$\q\hA$.

 Part~(b) follows from injectivity of the map $\q\hA_{n-1}\rarrow
\hA_{n-1}$ for $n\le3$ and commutativity of the square diagram
$\q\hA_{n-1}\overset t\rarrow\q\hA_n\rarrow\hA_n$, \
$\q\hA_{n-1}\rarrow\hA_{n-1}\overset t\rarrow\hA_n$.

 In part~(c), the graded ring homomorphism $\q\hA\rarrow\hA$ induces
a graded ring homomorphism between the quotient rings
$\q\hA/(\q\hA)t\rarrow\hA/\hA t=A$.
 Since the graded ring $\q\hA/(\q\hA)t$ is quadratic by
Proposition~\ref{central-quotient-reflects-quadraticity}, the latter
morphism, in turn, induces the desired graded ring homomorphism
$f\:\q\hA/(\q\hA)t\rarrow\q A$.
 The map~$f$ is an isomorphism in the degrees $n\le2$ by
construction, hence it is an isomorphism of quadratic graded rings by
Lemma~\ref{quadratic-ring-iso}.
\end{proof}

\begin{prop} \label{qdg-nonhomogeneous-duality-construction}
 Let $R\subset\tV\subset\tA$ be a\/ $3$\+left finitely projective
weak nonhomogeneous quadratic ring.
 Consider the corresponding graded rings $\hA=\bigoplus_{n=0}^\infty
F_n\tA$ and $A=\gr^F\tA$.
 Then the quadratic graded rings\/ $\q\hA$ and\/ $\q A$ are\/ $3$\+left
finitely projective.

 Let $\hB$ be the\/ $3$\+right finitely projective quadratic graded
ring quadratic dual to\/~$\q\hA$, and let $\d\:\hB\rarrow\hB$ be
the odd derivation of degree~$-1$ corresponding to the central
element $t\in\hA_1$, as in Lemma~\ref{central-iff-derivation}(b).
 Then $(\hB,\d)$ is a quasi-differential graded ring.
 The underlying graded ring $B=\ker\d\subset\hB$ is\/ $3$\+right
finitely projective quadratic and quadratic dual to\/~$\q A$.
\end{prop}

\begin{proof}
 The quadratic graded ring $\q A$ is $3$\+left finitely projective by
the definition of what it means for a weak nonhomogeneous quadratic ring
$\tA$ to be $3$\+left finitely projective
(see Section~\ref{self-consistency-subsecn}).
 To prove that the left $R$\+module $\q\hA_n$ is finitely generated
projective for $n\le3$, one argues by induction in $0\le n\le 3$, using
the short exact sequences of $R$\+$R$\+bimodules
$0\rarrow\q\hA_{n-1}\overset t\rarrow\q\hA_n\rarrow\q A_n\rarrow0$
provided by Lemma~\ref{quadratic-part-central-element}.
 Here the multiplication maps $\q\hA_{n-1}\overset t\rarrow\q\hA_n$ are
injective for $n\le 3$ by Lemma~\ref{quadratic-part-central-element}(b),
since $t$~is a nonzero-divisor in~$\hA$.

 The pair $(\hB,\d)$ is a quasi-differential graded ring by
Lemma~\ref{nonzerodivisor-in-degree-one} (applied to the quadratic
graded ring $\q\hA$ with its central element~$t$).
 The graded ring $B=\ker\d\subset\hB$ is generated by $B^1$ over $B^0$
by Proposition~\ref{nonzerodivisor-in-degree-two}(a) and quadratic
by Proposition~\ref{nonzerodivisor-in-degree-three}.
 The quadratic graded ring $B$ is $3$\+right finitely projective by
Lemma~\ref{B-right-projective}.
 It is quadratic dual to $\q A$ by
Proposition~\ref{nonzerodivisor-in-degree-two}(b).
\end{proof}

\begin{thm} \label{qdg-nonhomogeneous-duality-theorem}
 The construction of
Proposition~\ref{qdg-nonhomogeneous-duality-construction} defines
a fully faithful contravariant functor
\begin{equation} \label{qdg-nonhomogeneous-duality-functor}
 (R\rings_\wnlq)^\sop\lrarrow R\rings_{\qdg,\rq}
\end{equation}
from the category of\/ $3$\+left finitely projective weak nonhomogeneous
quadratic rings to the category of\/ $3$\+right finitely projective
quasi-differential quadratic graded rings over~$R$.
 The functor~\eqref{qdg-nonhomogeneous-duality-functor} forms
a commutative triangle diagram with the fully faithful contravariant
functor~\eqref{nonhomogeneous-duality-functor} and
the equivalence of categories~\eqref{quadratic-qdg-cdg-equivalence}.

 The same construction also defines a fully faithful strict
contravariant\/ $2$\+functor
\begin{equation} \label{qdg-nonhomogeneous-duality-2-functor}
 (\Rings_{\wnlq2})^\sop\lrarrow \Rings_{\qdg2,\rq}
\end{equation}
from the\/ $2$\+category of\/ $3$\+left finitely projective weak
nonhomogeneous quadratic rings to the\/ $2$\+category of\/ $3$\+right
finitely projective quasi-differential quadratic graded rings.
 The strict\/ $2$\+functor~\eqref{qdg-nonhomogeneous-duality-2-functor}
forms a commutative triangle diagram with the fully faithful strict
contravariant\/ $2$\+functor~\eqref{nonhomogeneous-duality-2-functor}
and the strict equivalence of\/
$2$\+categories~\eqref{quadratic-qdg-cdg-2cat-equivalence}.
\end{thm}

\begin{proof}
 The proof is straightforward.
\end{proof}

\subsection{PBW theorem} \label{pbw-theorem-subsecn}
 The Poincar\'e--Birkhoff--Witt theorem in nonhomogeneous quadratic duality
tells that, when restricted to the left/right finitely projective Koszul
rings on both sides, the fully faithful contravariant nonhomogeneous
quadratic duality functors
from Theorems~\ref{nonhomogeneous-duality-functor-existence-thm},
\ref{nonhomogeneous-duality-2-functor-theorem},
\ref{augmented-duality-functor-existence-thm},
\ref{augmented-duality-2-functor-existence-thm},
and~\ref{qdg-nonhomogeneous-duality-theorem} become anti-equivalences
of categories (or strict anti-equivalences of $2$\+categories).
 In other words, \emph{every right finitely projective Koszul CDG\+ring
(or quasi-differential graded ring) arises from a left finitely
projective nonhomogeneous Koszul ring}.

 A (weak) nonhomogeneous quadratic ring $R\subset \tV\subset \tA$ is
said to be \emph{left finitely projective Koszul} if the quadratic
graded ring $\q A=\q\gr^F\tA$ is left finitely projective Koszul
(in the sense of the definition in
Section~\ref{finitely-projective-Koszul-subsecn}).
 A nonnegatively graded CDG\+ring $(B,d,h)$ is said to be \emph{right
finitely projective Koszul} if the nonnegatively graded ring $B$ is
right finitely projective Koszul.
 A nonnegatively graded quasi-differential ring $(\hB,\d)$ is said to
be \emph{right finitely projective Koszul} if its underlying graded
ring $B=\ker\d\subset\hB$ is right finitely projective Koszul.

 The following result was mentioned in
Remark~\ref{weak-Koszul-strong-remark} in
Section~\ref{nonhomogeneous-quadratic-subsecn}.

\begin{thm} \label{weak-nonhomogeneous-Koszul-is-strong}
 If a weak nonhomogeneous quadratic ring is left finitely projective
Koszul, then it is nonhomogeneous quadratic.
\end{thm}

 So we will call such filtered rings as in
Theorem~\ref{weak-nonhomogeneous-Koszul-is-strong}
\emph{left finitely projective nonhomogeneous Koszul rings}.
 Notice that a filtered ring $(\tA,F)$ with an increasing filtration
$F$ such that $A_n=\gr^F_n\tA$ is a projective left $A_0$\+module
for all $n\ge0$ is left finitely projective Koszul (i.~e., the graded
ring $A=\gr^F\tA$ is left finitely projective Koszul) if and only
if the Rees ring $\hA=\bigoplus_{n=0}^\infty F_n\tA$ is left
finitely projective Koszul.
 This is clear from Lemma~\ref{hA-left-projective}
and Theorem~\ref{central-nonzerodivisor-theorem}.

 The proof of Theorem~\ref{weak-nonhomogeneous-Koszul-is-strong}
will be given at the end of Section~\ref{pbw-theorem-subsecn}.

\begin{thm} \label{pbw-theorem-thm}
 Every right finitely projective Koszul CDG\+ring arises from a left
finitely projective nonhomogeneous Koszul ring via the construction
of Proposition~\ref{nonhomogeneous-dual-cdg-ring}.
 Equivalently, every right finitely projective Koszul quasi-differential 
graded ring arises from a left finitely projective nonhomogeneous Koszul
ring via the construction of
Proposition~\ref{qdg-nonhomogeneous-duality-construction}.
\end{thm}

 The two assertions in Theorem~\ref{pbw-theorem-thm} are equivalent
in view of Theorems~\ref{cdg-qdg-equivalence}
and~\ref{qdg-nonhomogeneous-duality-theorem}.
 We will prove the second assertion.
 
\begin{proof}[First proof of Theorem~\ref{pbw-theorem-thm}]
 Let $(\hB,\d)$ be a right finitely projective Koszul quasi-differential
graded ring.
 Then the graded ring $B$ is right finitely projective Koszul by
definition and the graded ring $\hB$ is right finitely projective
Koszul by Theorem~\ref{adding-delta-preserves-Koszulity}(b).
 Let $\hA$ be the left finitely projective Koszul ring quadratic dual
to~$\hB$ (see Proposition~\ref{finitely-projective-Koszul-duality}),
and let $t\in\hA_1$ be the central element corresponding to
the odd derivation $\d\:\hB\rarrow\hB$
(see Lemma~\ref{central-iff-derivation}).
 Let $A=\hA/\hA t$ be the quotient ring.

 By Lemma~\ref{nonzerodivisor-in-degree-one}, the multiplication map
$\hA_0\overset t\rarrow\hA_1$ is injective and $A_1$ is a projective
left $R$\+module.
 By Proposition~\ref{nonzerodivisor-in-degree-two}(a),
the multiplication map $\hA_1\overset t\rarrow\hA_2$ is injective
and $A_2$ is a projective left $R$\+module.
 By Proposition~\ref{nonzerodivisor-in-degree-two}(b),
the right finitely projective Koszul ring $B$ is quadratic dual to
the $2$\+left finitely projective quadratic ring~$A$.
 Therefore, the quadratic graded ring $A$ is left finitely
projective Koszul.
 By Proposition~\ref{nonzerodivisor-in-degree-three},
the multiplication map $\hA_2\overset t\rarrow\hA_3$ is injective.

 Applying Theorem~\ref{central-element-pbw}, we conclude that $t$~is
a nonzero-divisor in~$\hA$.
 It follows that $\hA$ is the Rees ring of the filtered ring
$\tA=\hA/\hA(t-1)$ with the filtration $F_n\tA=\hA_n+\hA(t-1)$,
and $A$ is the associated graded ring, $A\simeq\gr^F\tA$.
 We have obtained the desired left finitely projective nonhomogeneous
Koszul ring~$\tA$.
\end{proof}

\begin{proof}[Second proof of Theorem~\ref{pbw-theorem-thm}]
 This is a particular case of the argument in~\cite[proof of
Theorem~11.6]{Psemi}.
 Let $(\hB,\d)$ be a right finitely projective Koszul quasi-differential
graded ring and $B=\ker\d\subset\hB$ be its underlying right finitely
projective Koszul graded ring.
 Then the right $R$\+modules $B^n$ are finitely generated projective by
definition and the right $R$\+modules $\hB^n$ are finitely generated
projective by Lemma~\ref{B-right-projective}.

 The (essentially) spectral sequence argument below is to be compared
with, and distinguished from, a quite different (and simpler)
spectral sequence argument proving the Poincar\'e--Birkhoff--Witt
theorem for nonhomogeneous quadratic algebras \emph{over the ground
field} in~\cite[Section~3.3]{Pcurv} and~\cite[Proposition~7.2 in
Chapter~5]{PP}.

 Put $D_n=\Hom_{R^\rop}(B^n,R)$.
 Then $D=\bigoplus_{n=0}^\infty D_n$ is a graded coring over the ring
$R$ with the counit map $\varepsilon\:D\rarrow D_0\simeq R$ dual to
the inclusion map $R\simeq B^0\rarrow B$ and the comultiplication maps
$\mu_{i,j}\:D_{i+j}\rarrow D_i\ot_R D_j$ obtained by dualizing
the multiplication maps $B^j\ot_R B^i\rarrow B^{i+j}$ in the graded
ring $B$,
\begin{multline*}
 D_{i+j}=\Hom_{R^\rop}(B^{i+j},R)\lrarrow
 \Hom_{R^\rop}(B^j\ot_RB^i,\>R) \\ \simeq
 \Hom_{R^\rop}(B^i,R)\ot_R\Hom_{R^\rop}(B^j,R)
 =D_i\ot_R D_j.
\end{multline*}
 In the pairing notation of Section~\ref{cdg-ring-constructed-subsecn},
we have $\lan \mu(f),b_1\ot b_2\ran =\lan f,b_1b_2\ran$
for all $f\in D$ and $b_1$, $b_2\in B$. 

 Similarly, we set $\hD_n=\Hom_{R^\rop}(\hB^n,R)$; so
$\hD=\bigoplus_{n=0}^\infty\hD_n$ is also a graded coring over~$R$.
 The odd derivation $\d\:\hB\rarrow\hB$ is an $R$\+$R$\+bilinear
map, so it dualizes to an $R$\+$R$\+bilinear map
$\Hom_{R^\rop}(\d,R)\:\hD\rarrow\hD$, which we denote for brevity
also by~$\d$.
 There is a sign rule involved: in the pairing notation, we put
$\lan \d(f),b\ran = (-1)^{|f|+1}\lan f,\d(b)\ran$ for $f\in\hD$ and
$b\in\hB$.
 The map $\d\:\hD\rarrow\hD$ is an odd coderivation of degree~$1$
on the coring $\hD$, in the sense that its components act as
$\d_n\:\hD_n\rarrow\hD_{n+1}$ and, in the symbolic notation
$\mu(f)=\mu_1(f)\ot\mu_2(f)$ for the comultiplication, one has
$$
 \d(\mu(f))=\d\mu_1(f)\ot\mu_2(f)+
 (-1)^{|\mu_1(f)|}\mu_1(f)\ot\d\mu_2(f)
 \quad\text{for all $f\in\hD$}.
$$
 The odd coderivation~$\d$ on the graded coring $\hD$ has zero
square, $\d^2=0$ (so one can say that $\hD$ is a \emph{DG\+coring
over\/~$R$}).
 Furthermore, the odd coderivation $\d$ on $\hD$ has vanishing
cohomology, $H_\d(\hD)=0$.
 The cokernel of~$\d$ is the quotient coring $D$ of the coring $\hD$,
that is $\hD\twoheadrightarrow\coker(\d)\simeq D$.

 Consider the bigraded $R$\+$R$\+bimodule $K$ with the components
$K^{p,q}=\hD_{q-p}$ for $p\le0$, $q\le0$, and $K^{p,q}=0$
otherwise.
 The $R$\+$R$\+bimodule $K$ can be viewed as a bigraded coring
over~$R$ with the comultiplication inherited from the comultiplication
of~$\hD$.
 Considered as a graded coring in the total grading $p+q$,
the coring $K$ has an odd coderivation $\d_K$ of degree~$1$
with the components $\d_K^{p,q}\:K^{p,q}\rarrow K^{p,q+1}$
given by the rule $\d_K^{p,q}=\d_{q-p}\:\hD_{q-p}\rarrow\hD_{q-p+1}$.
 So for every fixed $p=-n$, \,$n\ge0$, the components $K^{p,q}$
with varying~$q$ form a complex of $R$\+$R$\+bimodules
$$
 0\lrarrow K^{p,p}\overset{\d_K^{p,p}}\lrarrow K^{p,p+1}
 \xrightarrow{\d_K^{p,p+1}}\dotsb\xrightarrow{\d_K^{p,-2}}K^{p,-1}
 \xrightarrow{\d_K^{p,-1}}K^{p,0}\lrarrow0
$$
isomorphic to
$$
 0\lrarrow\hD_0\overset{\d_0}\lrarrow\hD_1\overset{\d_1}\lrarrow
 \dotsb\overset{\d_{n-2}}\lrarrow
 \hD_{n-1}\overset{\d_{n-1}}\lrarrow\hD_n\lrarrow0.
$$
 The only nontrivial cohomology bimodule of this complex occurs at
its rightmost term ($q=0$) and is isomorphic to~$D_n$.
 So there is a surjective morphism of bigraded corings $K\rarrow D$
inducing an isomorphism of the bimodules/corings of cohomology, where
the coring $D$ is placed in the bigrading $D^{p,0}=D_n=D_{-p}$ and
endowed with the zero differential.

 Denote by $K_+$ the cokernel of the inclusion $R\simeq K^{0,0}
\rarrow K$.
 Consider the tensor ring $T_R(K_+)$ of the bigraded
$R$\+$R$\+bimodule~$K_+$.
 By the definition, $T_R(K_+)$ is a trigraded ring with the gradings
$p\le0$ and $q\le0$ inherited from the bigrading of $K_+$ and
the additional grading $r\ge0$ by the number of tensor factors.
 We will consider $T_R(K_+)$ as a graded ring in the total grading
$p+q+r$.
 The graded ring $T_R(K_+)$ is endowed with three odd derivations
of total degree~$1$, which we will now introduce.

 For any (graded) $R$\+$R$\+bimodule $V$, derivations of the tensor
ring $T(V)$ (say, odd derivations with respect to the total parity
on $T(V)$) annihilating the subring $R\subset T(V)$ are uniquely
determined by their restriction to the subbimodule $V\subset T(V)$,
which can be an arbitrary $R$\+$R$\+bilinear map $V\rarrow T(V)$
(that is homogeneous of the prescribed degrees).
 In the situation at hand, we let $\d_T$~be the only odd derivation of
$T_R(K_+)$ which preserves the subbimodule $K_+\subset T_R(K_+)$
and whose restriction to $K_+$ is equal to~$-\d_K$.
 Let $d_T$ be the only odd derivation of $T_R(K_+)$ which maps
$K_+$ into $K_+\ot_RK_+$ by the comultiplication map~$\mu$ with
the sign rule
$$
 d_T(k)=(-1)^{p_1+q_1}\mu_1(k)\ot\mu_2(k), \quad
 k\in K_+, \ \mu_1(k)\in K^{p_1,q_1}.
$$
 Finally, let $\delta_T$ be the only odd derivation of $T_R(K_+)$
whose restiction to $K_+$ is the identity map of the component
$K^{-1,-1}=R$ to the unit component $T_R(K_+)^{0,0,0}=R$ and
zero on all the remaining components of~$K_+$.

 All the three derivations are constructed to be odd derivations
in the parity $p+q+r$.
 The derivations $\d_T$, $d_T$, and~$\delta_T$ have tridegrees
$(0,1,0)$, \ $(0,0,1)$, and $(1,1,-1)$, respectively, in
the trigrading $(p,q,r)$ of the tensor ring $T_R(K_+)$.
 All the three derivations have zero squares, and they pairwise
anti-commute.

 There is an increasing filtration $F$ on the graded ring $T_R(K_+)$
whose component $F_nT_R(K_+)$ is the direct sum of all
the trigrading components $T_R(K_+)^{p,q,r}$ with $-p\le n$.
 This filtration is compatible with the differentials $\d_T$,
$d_T$, and~$\delta_T$; the graded ring $\gr^FT_R(K_+)=
\bigoplus_{n=0}^\infty F_nT_R(K_+)/F_{n-1}T_R(K_+)$ with
the differential induced by $\d_T+d_T+\delta_T$ is naturally
isomorphic to $T_R(K_+)$ with the differential $\d_T+d_T$.

 Let us also consider the tensor ring $T_R(D_+)$ of
the $R$\+$R$\+bimodule $D_+=D/R$.
 Similarly to the above, we endow $T_R(D_+)$ with the grading~$p$
coming from the grading $D^p=D_n=D_{-p}$ of $D_+$ and
the grading~$r$ by the number of tensor factors.
 As to the grading~$q$, we set it to be identically zero on $T_R(D_+)$.
 We will consider $T_R(D_+)$ as a graded ring in the total grading
$p+r$.
 Let~$d'_T$ be the only odd derivation of $T_R(D_+)$ which maps
the subbimodule $D_+\subset T_R(D_+)$ to $D_+\ot_RD_+$ by
the comultiplication map~$\mu$ with the sign rule similar to
the above.

 The DG\+ring $T_R(D_+)$ is naturally isomorphic to the cobar
complex~\eqref{L=R=M} (with the roles of the left and right sides
switched) computing the bigraded $\Ext_{B^{\rop}}(R,R)$.
 There is some sign rule involved in this isomorphism, which is
discussed in~\cite[proof of Theorem~11.6]{Psemi}.
 By Proposition~\ref{finitely-projective-Koszul-duality}, it follows
that the bigraded ring of cohomology of the DG\+ring $T_R(D_+)$ is
isomorphic to the left finitely projective Koszul graded ring
$A=\bigoplus_{n=0}^\infty A_n$ quadratic dual to~$B$, placed in
the diagonal bigrading $-p=r=n$.
 In the total grading $p+r$, the cohomology ring of the DG\+ring
$T_R(D_+)$ is the whole ring $A$ placed in the degree $p+r=0$.

 Consider the morphism of DG\+rings $(T_R(K_+),\>\d_T+d_T)\rarrow
(T_R(D_+),d'_T)$ induced by the surjective morphism of graded
corings $K\rarrow D$.
 The morphism of DG\+corings $K\rarrow D$ is a quasi-isomorphism,
so the induced morphism of DG\+rings is a quasi-isomorphism, too,
due to the right flatness/projectivity conditions imposed on
the $R$\+$R$\+bimodule $D$ and the presence of a nonpositive
(essentially, negative) internal grading~$p$.
 The fact that the components of fixed grading~$p$ in
the DG\+coring $K$ are finite complexes (of projective right
$R$\+modules, with projective right $R$\+modules of cohomology)
is relevant here.
 Thus we have $H_{\d_T+d_T}^0(T_R(K_+))\simeq A$ and
$H_{\d_T+d_T}^i(T_R(K_+))=0$ for $i\ne0$.

 Finally, we put $\tA=H_{\d_T+d_T+\delta_T}^0(T_R(K_+))$.
 Then the ring $\tA$ is endowed with an increasing filtration $F$
induced by the filtration $F$ of the DG\+ring
$(T_R(K_+),\>\d_T+d_T+\delta_T)$.
 Since $H_{\d_T+d_T}^i(T_R(K_+))=0$ for $i\ne0$, we can conclude
that the associated graded ring of the ring $\tA$ is naturally
isomorphic to the graded ring $A$, that is $\gr^F\tA\simeq A$,
while $H_{\d_T+d_T+\delta_T}^i(T_R(K_+))=0$ for $i\ne0$.

 Since the graded ring $A=\gr^F\tA$ is left finitely projective
Koszul, the graded ring $\hA=\bigoplus_{n=0}^\infty F_n\tA$ is
left finitely projective Koszul as well (by
Lemma~\ref{hA-left-projective} and
Theorem~\ref{central-nonzerodivisor-theorem}).
 Let $\hB'$ be the right finitely projective Koszul graded ring
quadratic dual to~$\hA$.
 By Proposition~\ref{qdg-nonhomogeneous-duality-construction},
the graded ring $\hB'$ is endowed with an odd derivation
$\d'\:\hB'\rarrow\hB'$ of degree~$-1$, making it a right finitely
projective Koszul quasi-differential graded ring.
 The underlying graded ring $B'=\ker\d'\subset\hB'$ of
the quasi-differential graded ring $(\hB',\d')$ is quadratic dual
to $A$, so we have $B'\simeq B$.
 It remains to construct a natural isomorphism of quasi-differential
graded rings $(\hB',\d')\simeq(\hB,\d)$.

 For this purpose, let us consider the dual graded corings
$D'=\Hom_{R^\rop}(B',R)$ and $\hD'=\Hom_{R^\rop}(\hB',R)$.
 As in the beginning of this proof, we have an odd coderivation
$\d'\:\hD'\rarrow\hD'$ of degree~$1$ dual to the odd derivation
$\d'\:\hB'\rarrow\hB'$.
 It suffices to construct a natural isomorphism $(\hD,\d)\rarrow
(\hD',\d')$ of DG\+corings over~$R$.

 The embedding of the component $\hD_1=T_R(K_+)^{-1,0,1}\rarrow
T_R(K_+)$ induces an isomorphism of $R$\+$R$\+bimodules
$\hD_1\simeq F_1\tA$.
 The composition $\hD_2\rarrow\hD_1\ot_R\hD_1\simeq F_1\tA\ot_R
F_1\tA\rarrow F_2\tA$ of the comultiplication and multiplication
maps vanishes, being killed by the differential
$(\d_T+d_T+\delta_T)^{-2,0,1}=d_T^{-2,0,1}\:T_R(K_+)^{-2,0,1}
\rarrow T_R(K_+)^{-2,0,2}$.
 So there is a natural morphism of graded corings $\hD\rarrow\hD'$.
 Since the embedding $R=F_0\tA\rarrow F_1\tA$ corresponds to
the map $\d_0\:R=\hD_0\rarrow\hD_1$ under the isomorphisms
$F_0\tA=R=\hD_0$ and $F_1\tA\simeq\hD_1$, the graded coring
morphism $\hD\rarrow\hD'$ forms a commutative square diagram with
the differentials $\d$ on $\hD$ and $\d'$ on~$\hD'$.

 The induced morphism $\coker(\d)\rarrow\coker(\d')$ coincides with
the natural isomorphism $D\rarrow D'$ on the components of
degree~$1$, and hence on the other components as well.
 Therefore, the morphism of corings $\hD\rarrow\hD'$ is also
an isomorphism.
\end{proof}

\begin{proof}[Proof of
Theorem~\ref{weak-nonhomogeneous-Koszul-is-strong}]
 Let $R\subset\tV\subset\tA$ be a weak nonhomogeneous quadratic ring
and $A=\gr^F\tA$ be its associted graded ring with respect to
the filtration $F$ generated by $F_1\tA=\tV$ over $F_0\tA=R$.
 Assume that the quadratic graded ring $\q A=\q\gr^F\tA$ is left
finitely projective Koszul.
 Let $(B,d,h)$ be the CDG\+ring corresponding to $\tA$ under
the construction of Proposition~\ref{nonhomogeneous-dual-cdg-ring}.
 Equivalently, one can consider the quasi-differential graded ring
$(\hB,\d)$ with the underlying graded ring $B=\ker\d\subset\hB$
corresponding to $\tA$ under the construction of
Proposition~\ref{qdg-nonhomogeneous-duality-construction}.

 Whichever one of these two points of view one takes, the graded ring
$B$ is quadratic and quadratic dual to the quadratic graded ring $\q A$
by construction.
 Hence the graded ring $B$ is right finitely projective Koszul by
Proposition~\ref{finitely-projective-Koszul-duality}.
 Applying Theorem~\ref{pbw-theorem-thm}, we see that the CDG\+ring
$(B,d,h)$ or the quasi-differential ring $(\hB,\d)$ comes from
a left finitely projective nonhomogeneous Koszul ring $(\tA',F)$.
 Following either one of the two proofs of
Theorem~\ref{pbw-theorem-thm}, the graded ring $\gr^F\tA'$ is quadratic
and left finitely projective Koszul.

 It remains to observe that the nonhomogeneous quadratic duality
functor assigns the same CDG\+ring $(B,d,h)$ (or the same
quasi-differential graded ring $(\hB,\d)$) to the two $3$\+left
finitely projective (weak) nonhomogeneous quadratic rings $\tA$
and~$\tA'$.
 Since the nonhomogeneous quadratic duality functor is fully faithful
by Theorem~\ref{nonhomogeneous-duality-functor-existence-thm}
or~\ref{qdg-nonhomogeneous-duality-theorem}, it follows that
the two (weak) nonhomogeneous quadratic rings $(\tA,F)$ and $(\tA',F)$
are isomorphic.
 Hence the graded ring $A=\gr^F\tA\simeq\gr^F\tA'$ is quadratic
(and left finitely projective Koszul).
 In other words, the weak nonhomogeneous quadratic ring $R\subset\tV
\subset\tA$ is actually nonhomogeneous quadratic.
\end{proof}

\subsection{Anti-equivalences of Koszul ring categories}
\label{anti-equivalences-koszul-rings-subsecn}
 Let $R$ be an associative ring.
 The \emph{category of left finitely projective nonhomogeneous Koszul
rings over~$R$}, denoted by $R\rings_\nlk$, is defined as
the full subcategory of the category of filtered rings $R\rings_\fil$
(see Section~\ref{nonhomogeneous-duality-functor-subsecn}) whose objects
are the left finitely projective nonhomogeneous Koszul rings $(\tA,F)$.
 The \emph{category of right finitely projective Koszul CDG\+rings
over~$R$}, denoted by $R\rings_{\cdg,\rk}$, is the full subcategory in
the category of nonnegatively graded CDG\+rings $R\rings_\cdg$ whose
objects are the right finitely projective Koszul CDG\+rings.
 The \emph{category of right finitely projective Koszul
quasi-differential rings over~$R$}, denoted by $R\rings_{\qdg,\rk}$, is
the full subcategory in the category of nonnegatively graded
quasi-differential rings $R\rings_\qdg$ whose objects are
the right finitely projective Koszul quasi-differential graded rings.

\begin{cor} \label{nonhomogeneous-koszul-duality-anti-equivalence}
 The constructions of
Theorems~\ref{nonhomogeneous-duality-functor-existence-thm},
\ref{cdg-qdg-equivalence},
and~\ref{qdg-nonhomogeneous-duality-theorem} define natural
(anti)-equivalences
$$
 (R\rings_\nlk)^\sop\,\simeq\,R\rings_{\qdg,\rk}\,\simeq\,
 R\rings_{\cdg,\rk}
$$
between the categories of left finitely projective nonhomogeneous
Koszul rings, right finitely projective Koszul quasi-differential rings,
and right finitely projective Koszul CDG\+rings over~$R$.
\end{cor}

\begin{proof}
 Follows from the mentioned theorems and
Theorem~\ref{pbw-theorem-thm}.
\end{proof}

 The \emph{$2$\+category of left finitely projective nonhomogeneous
Koszul rings}, denoted by $\Rings_{\nlk2}$, is defined as
the following $2$\+subcategory in the $2$\+category of filtered
rings $\Rings_{\fil2}$ (see
Section~\ref{nonhomogeneous-duality-2-functor-subsecn}).
 The objects of $\Rings_{\nlk2}$ are the left finitely projective
nonhomogeneous Koszul rings $(\tA,F)$.
 All morphisms in $\Rings_{\fil2}$ between objects of $\Rings_{\nlk2}$
are morphisms in $\Rings_{\nlk2}$, and all $2$\+morphisms in
$\Rings_{\fil2}$ between morphisms of $\Rings_{\nlk2}$ are
$2$\+morphisms in $\Rings_{\nlk2}$.

 The \emph{$2$\+category of right finitely projective Koszul CDG\+rings},
denoted by $\Rings_{\cdg2,\rk}$, is the following $2$\+subcategory in
the $2$\+category of nonnegatively graded CDG\+rings $\Rings_{\cdg2}$.
 The objects of $\Rings_{\cdg2,\rk}$ are all the right finitely
projective Koszul CDG\+rings $(B,d,h)$.
 All morphisms in $\Rings_{\cdg2}$ between objects of
$\Rings_{\cdg2,\rk}$ are morphisms in $\Rings_{\cdg2,\rk}$, and all
$2$\+morphisms in $\Rings_{\cdg2}$ between morphisms of
$\Rings_{\cdg2,\rk}$ are $2$\+morphisms in $\Rings_{\cdg2,\rk}$.

 The \emph{$2$\+category of right finitely projective Koszul
quasi-differential graded rings}, denoted by $\Rings_{\qdg2,\rk}$, is
the following $2$\+subcategory in the $2$\+category of nonnegatively
graded quasi-differential rings $\Rings_{\qdg2}$ (see
Section~\ref{qdg-cdg-subsecn}).
 The objects of $\Rings_{\qdg2,\rk}$ are all the right finitely
projective Koszul quasi-differential graded rings $(\hB,\d)$.
 All morphisms in $\Rings_{\qdg2}$ between objects of
$\Rings_{\qdg2,\rk}$ are morphisms in $\Rings_{\qdg2,\rk}$, and
all $2$\+morphisms in $\Rings_{\qdg2}$ between morphisms of
$\Rings_{\qdg2,\rk}$ are $2$\+morphisms in $\Rings_{\qdg2,\rk}$.

\begin{cor}
 The constructions of
Theorems~\ref{nonhomogeneous-duality-2-functor-theorem},
\ref{cdg-qdg-2cat-equivalence},
and~\ref{qdg-nonhomogeneous-duality-theorem} define natural
strict (anti)-equivalences
$$
 (\Rings_{\nlk2})^\sop\,\simeq\,\Rings_{\qdg2,\rk}\,\simeq\,
 \Rings_{\cdg2,\rk}
$$
between the\/ $2$\+categories of left finitely projective nonhomogeneous
Koszul rings, right finitely projective Koszul quasi-differential rings,
and right finitely projective Koszul CDG\+rings.
\end{cor}

\begin{proof}
 Follows from the mentioned theorems and
Theorem~\ref{pbw-theorem-thm}.
\end{proof}

 A left finitely projective nonhomogeneous Koszul ring $(\tA,F)$ is
said to be \emph{left augmented} if the ring $\tA$ is left augmented
over its subring~$F_0\tA$.
 In other words, this means that a left ideal $\tA^+\subset\tA$ is
chosen such that $\tA=F_0\tA\oplus\tA^+$
(see Section~\ref{augmented-subsecn}).
 The \emph{category of left augmented left finitely projective
nonhomogeneous Koszul rings over~$R$}, denoted by $R\rings_\nlk^\laug$,
is defined as the full subcategory in the category of left augmented
filtered rings $R\rings_\fil^\laug$ whose objects are the left
augmented left finitely projective nonhomogeneous Koszul rings.

 A nonnegatively graded DG\+ring $(B,d)$ (in the sense of
Section~\ref{augmented-subsecn}) is said to be \emph{right finitely
projective Koszul} if the graded ring $B$ is right finitely
projective Koszul.
 The \emph{category of right finitely projective Koszul DG\+rings
over~$R$}, denoted by $R\rings_{\dg,\rk}$ is defined as the full
subcategory in the category of nonnegatively graded DG\+rings
$R\rings_\dg$ whose objects are the right finitely projective
Koszul DG\+rings.

\begin{cor} \label{left-augmented-koszul-duality-anti-equivalence}
 The construction of
Theorem~\ref{augmented-duality-functor-existence-thm}
defines a natural anti-equivalence
$$
 (R\rings_\nlk^\laug)^\sop\,\simeq\,R\rings_{\dg,\rk}
$$
between the category of left augmented left finitely projective
nonhomogeneous Koszul rings and the category of right finitely
projective Koszul DG\+rings.
\end{cor}

\begin{proof}
 Follows from Theorem~\ref{pbw-theorem-thm} and the discussion
in the proof of Theorem~\ref{augmented-duality-functor-existence-thm}.
\end{proof}

 The \emph{$2$\+category of left augmented left finitely projective
nonhomogeneous Koszul rings}, denoted by $\Rings_\nlk^{\laug2}$, is
defined as the following $2$\+subcategory in the $2$\+category of
left augmented filtered rings $\Rings_\fil^{\laug2}$
(see Section~\ref{augmented-subsecn}).
 The objects of $\Rings_\nlk^{\laug2}$ are the left augmented left
finitely projective nonhomogeneous Koszul rings $(\tA,F,\tA^+)$.
 All morphisms in $\Rings_\fil^{\laug2}$ between objects of
$\Rings_\nlk^{\laug2}$ are morphisms in $\Rings_\nlk^{\laug2}$,
and all $2$\+morphisms in $\Rings_\fil^{\laug2}$ between morphisms of
$\Rings_\nlk^{\laug2}$ are $2$\+morphisms in $\Rings_\nlk^{\laug2}$.

 The \emph{$2$\+category of right finitely projective Koszul DG\+rings},
denoted by $\Rings_{\dg2,\rk}$, is defined as the following
$2$\+subcategory in the $2$\+category of nonnegatively graded
DG\+rings $\Rings_{\dg2}$.
 The objects of $\Rings_{\dg2,\rk}$ are the right finitely projective
Koszul DG\+rings $(B,d)$.
 All morphisms in $\Rings_{\dg2}$ between objects of $\Rings_{\dg2,\rk}$
are morphisms in $\Rings_{\dg2,\rk}$, and all $2$\+morphisms in
$\Rings_{\dg2}$ between morphisms of $\Rings_{\dg2,\rk}$ are
$2$\+morphisms in $\Rings_{\dg2,\rk}$.

\begin{cor}
 The construction of
Theorem~\ref{augmented-duality-2-functor-existence-thm}
defines a natural strict anti-equivalence
$$
 (\Rings_\nlk^{\laug2})^\sop\,\simeq\,\Rings_{\dg2,\rk}
$$
between the\/ $2$\+category of left augmented left finitely projective
nonhomogeneous Koszul rings and the\/ $2$\+category of right finitely
projective Koszul DG\+rings.  \qed
\end{cor}

\Section{Comodules and Contramodules over Graded Rings}
\label{comodules-and-contramodules-secn}

\subsection{Ungraded comodules over nonnegatively graded rings}
\label{ungraded-comodules-subsecn}
 Let $B=\bigoplus_{n=0}^\infty B_n$ be a nonnegatively graded ring.
 We denote the underlying ungraded ring of $B$ by the same letter
$B=\Sigma B$ (see Section~\ref{graded-ungraded-subsecn}), and consider
ungraded right $B$\+modules.
 An ungraded right $B$\+module $M$ is said to be a \emph{$B$\+comodule}
(or an \emph{ungraded right $B$\+comodule}) if for every element
$x\in M$ there exists an integer $m\ge0$ such that $xb=0$ in $M$
for all $b\in B_n$, \,$n>m$.
 Notice that any right $B$\+comodule has a natural structure of right
module over the ring $\Pi B=\prod_{n=0}^\infty B_n$.
{\hbadness=1275\par}

 We will denote the full subcategory of right $B$\+comodules by
$\comodr B\subset\modr B$.
 The following examples explain the terminology.

\begin{exs} \label{comodules-examples}
 (1)~Let $B$ be a nonnegatively graded associative algebra over a field
$k=B_0$ with finite-dimensional grading components, $\dim_kB_n<\infty$
for all $n\ge0$.
 Then the graded dual vector space $C=\bigoplus_{n=0}^\infty B_n^*$ to
$B$ has a natural structure of graded coassociative coalgebra over~$k$.
 An ungraded right $B$\+comodule in the sense of the above definition is
the same thing as an ungraded right $C$\+comodule.

 (2)~More generally, let $B$ be a nonnegatively graded ring with
the degree-zero grading component $R=B_0$.
 Assume that $B_n$ is a finitely generated projective right $R$\+module
for every $n\ge0$, and put $C=\bigoplus_{n=0}^\infty C_n$, where
$C_n=\Hom_{R^\rop}(B_n,R)$.
 Then, in view of Lemma~\ref{tensor-dual-lemma}(b), applying
the functor $\Hom_{R^\rop}({-},R)$ to the multiplication maps
$B_i\ot_R B_j\rarrow B_{i+j}$ produces comultiplication maps
$C_{i+j}\rarrow C_j\ot_R C_i$ endowing the graded $R$\+$R$\+bimodule $C$
with a natural structure of graded coring over~$R$.
 An ungraded right $B$\+comodule $M$ in the sense of the above
definition is the same thing as an ungraded right $C$\+comodule,
i.~e., a right $R$\+module $M$ endowed with a coassociative, counital
comultiplication map $M\rarrow M\ot_RC$ (see
Example~\ref{co-monad-examples}(3) below for a discussion of
corings and comodules).

 (3)~In the context of~(1), one can also say that ungraded left
$B$\+comodules are the same thing as ungraded left $C$\+comodules.
 But let us \emph{warn} the reader that, in the context of~(2),
left $B$\+comodules are, generally speaking, entirely \emph{unrelated}
to left $C$\+comodules.
 Rather, in order to describe left $C$\+comodules, one needs to consider
the graded ring $\shB$ defined in Remark~\ref{B-sharp-remark}.
 Then, assuming that $C_n$ is a finitely generated projective right
$R$\+module for every $n\ge0$, an ungraded left $\shB$\+comodule is
the same thing as an ungraded left $C$\+comodule (and the same holds
for graded comodules as defined in
Section~\ref{graded-comodules-subsecn} below;
cf.\ Example~\ref{graded-comodules-example}).
\end{exs}

 Obviously, the full subcategory $\comodr B\subset\modr B$ is closed
under subobjects, quotients, and infinite direct sums.
 In fact, $\comodr B$ is a hereditary pretorsion class in $\modr B$
corresponding to the filter/topology of right ideals in $B$ in which
the two-sided ideals $B_{\ge m}=\bigoplus_{n\ge m}B_n\subset B$,
\,$m\ge1$, form a base~\cite[Section~VI.4]{St},
\cite[Sections~2.3\+-2.4]{Pcoun}.
 So $\comodr B$ is a Grothendieck abelian category and its inclusion
$\comodr B\rarrow\modr B$ is an exact functor preserving infinite
direct sums.

 Let us show that, under certain assumptions, the full subcategory
$\comodr B$ is also closed under extensions in $\modr B$ (in other
words, it is a \emph{hereditary torsion class}).

\begin{lem} \label{generated-by-several-components}
 Let $B=\bigoplus_{n=0}^\infty B_n$ be a nonnegatively graded ring
and $m\ge1$ be an integer.
 Then the following five conditions are equivalent:
\begin{enumerate}
\renewcommand{\theenumi}{\alph{enumi}}
\item there exists an integer $k\ge1$ such that the graded right
$B$\+module $B_{\ge1}=\bigoplus_{n=1}^\infty B_n$ is generated by
elements of degree~$\le k$;
\item there exists an integer $l\ge m$ such that the graded right
$B$\+module $B_{\ge m}=\bigoplus_{n=m}^\infty B_n$ is generated by
elements of degree~$\le l$;
\item there exists an integer $k\ge1$ such that for
every $n>k$ the multiplication map\/ $\bigoplus_{i+j=n}^{i,j\ge 1}
B_i\ot_R B_j\rarrow B_n$ is surjective;
\item there exists an integer $l\ge 2m-1$ such that for
every $n>l$ the multiplication map\/ $\bigoplus_{i+j=n}^{i,j\ge m}
B_i\ot_R B_j\rarrow B_n$ is surjective.
\item there exists an integer $k\ge 1$ such that the ring $B$ is
generated by its subgroup\/ $\bigoplus_{n=0}^k B_k\subset B$.
\end{enumerate}
\end{lem}

\begin{proof}
 The equivalences (a)\,$\Longleftrightarrow$\,(c)%
\,$\Longleftrightarrow$\,(e) (for the same~$k$) are easy.
 So are the implications (d)\,$\Longrightarrow$\,(b)%
\,$\Longrightarrow$\,(c) (for $k=l$).
 The implication (e)\,$\Longrightarrow$\,(d) holds for $l=2m-2+k$.
\end{proof}

 Notice that it follows from Lemma~\ref{generated-by-several-components}
that the conditions~(a\+-b) are left-right symmetric (because
the conditions~(c\+-e) are).
 It also follows that the conditions~(b) and~(d) do not depend on~$m$
(because the conditions (a), (c), and~(e) don't).

\begin{lem} \label{finitely-generated-graded-ideal}
 Let $B=\bigoplus_{n=0}^\infty B_n$ be a nonnegatively graded ring
and $m\ge1$ be an integer.
 Then, for $m=1$, the following two conditions are equivalent:
\begin{enumerate}
\renewcommand{\theenumi}{\roman{enumi}}
\item $B_{\ge m}=\bigoplus_{n=m}^\infty B_n$ is a finitely generated
(graded) right $B$\+module;
\item $B_n$ is a finitely generated right $B_0$\+module for every
$n\ge m$, and any one of the equivalent conditions of
Lemma~\ref{generated-by-several-components} holds.
\end{enumerate}
 For any $m\ge1$, the implication (ii)\,$\Longrightarrow$\,(i) is true.
\end{lem}

\begin{proof}
 (i)\,$\Longrightarrow$\,(ii)  Clearly, (i)\,$\Longrightarrow$\,(b)
(for the same~$m$); so it remains to show that (i)~implies the first
part of~(ii).
 Indeed, let $\{b_{s,n}\in B_n\mid n\ge m, \,1\le s\le t_n\}$ be
a finite set of homogeneous generators of the right
$B$\+module~$B_{\ge m}$.
 Put $R=B_0$.
 Then, for every $n\ge m$, the cokernel of the multiplication map
$\bigoplus_{i+j=n}^{i\ge m,j\ge1}B_i\ot_R B_j\rarrow B_n$ is
generated by the cosets of the elements $b_{s,n}$, $1\le s\le t_n$
as a right $R$\+module.
 For any $R$\+$R$\+bimodules $U$ and $V$, if the right $R$\+modules
$U$ and $V$ are finitely generated, then so is the right $R$\+module
$U\ot_RV$.
 This allows to prove by induction in~$n$ that the right $R$\+module
$B_n$ is finitely generated for $n\ge m$.

 (ii)\,$\Longrightarrow$\,(i)  Choose~$l$ as in~(b), and for every
$m\le n\le l$ choose a finite set of generators $b_{s,n}$,
$1\le s\le t_n$ of the right $R$\+module~$B_n$.
 Then the right $B$\+module $B_{\ge m}$ is generated by the elements
$\{b_{s,n}\in B_n\mid m\le n\le l, \,1\le s\le t_n\}$.
\end{proof}

\begin{lem} \label{augmentation-implies-component-ideals}
 Let $B=\bigoplus_{n=0}^\infty B_n$ be a nonnegatively graded ring
such that the augmentation ideal $B_{\ge1}$ is finitely generated
as a right ideal in~$B$.
 Then, for every $m\ge1$, the right ideal $B_{\ge m}\subset B$ is also
finitely generated.
\end{lem}

\begin{proof}
 Follows from Lemmas~\ref{generated-by-several-components}
and~\ref{finitely-generated-graded-ideal}.
\end{proof}

 The following particular case is important for our purposes.
 Let $B$ be a nonnegatively graded ring generated by $B_1$ over~$B_0$.
 Assume the right $B_0$\+module $B_1$ is finitely generated.
 Then it is clear from Lemma~\ref{generated-by-several-components}
(take $k=1$) and Lemma~\ref{finitely-generated-graded-ideal}
that the augmentation ideal $B_{\ge1}\subset B$ is finitely
generated as a right ideal in~$B$.
 In particular, this holds for any $2$\+right finitely projective
quadratic graded ring~$B$.

\begin{prop} \label{comodules-closed-under-extensions-prop}
 Assume that the augmentation ideal $B_{\ge1}=
\bigoplus_{n=1}^\infty B_n$ of a nonnegatively graded ring $B$ is
finitely generated as a right ideal in~$B$.
 Then the full subcategory of ungraded right $B$\+comodules\/
$\comodr B$ is closed under extensions in the category of
right $B$\+modules\/ $\modr B$.
\end{prop}

\begin{proof}
 Let $0\rarrow K\rarrow L\rarrow M\rarrow 0$ be a short exact sequence
of (ungraded) right $B$\+modules.
 Assume that $K$ and $M$ are $B$\+comodules.
 Let $x\in L$ be an element, and let $y\in M$ be the image of~$x$
under the surjective $B$\+module morphism $L\rarrow M$.
 By assumption, there exists $m\ge1$ such that $yB_{\ge m}=0$ in~$M$.

 By Lemma~\ref{augmentation-implies-component-ideals}, the right
ideal $B_{\ge m}\subset B$ is finitely generated.
 Let $\{b_s\mid 1\le s\le t\}$ be a finite set of its homogeneous
generators; so $B_{\ge m}=\sum_{s=1}^t b_sB$.
 Let $k_s$~be the homogeneous degree of the element $b_s\in B$, and
let $k\ge m$ be an integer such that $k\ge k_s$ for all $1\le s\le t$.
 Then we have $B_n=\sum_{s=1}^t b_sB_{n-k_s}$ for all $n\ge k$.

 For every $1\le s\le t$, we have $z_s=xb_s\in K\subset L$.
 By assumption, there exists $j_s\ge1$ such that $z_s B_{\ge j_s}=0$
in~$K$.
 Let $j\ge1$ integer such that $j\ge j_s$ for all $1\le s\le t$.
 Then for every $n\ge k+j$ we have $xB_n=\sum_{s=1}^t xb_sB_{n-k_s}
=\sum_{s=1}^t z_sB_{n-k_s}=0$ in $L$, since $n-k_s\ge j_s$.
 Thus $xB_{\ge k+j}=0$, as desired.
\end{proof}

\subsection{Right exact monads on abelian categories}
\label{monads-subsecn}
 Let $\sA$ be a category.
 A \emph{monad} on $\sA$ is a monoid object in the monoidal category
of endofunctors of $\sA$ (with respect to the composition of functors).
 In other words, this means that a monad $\boM\:\sA\rarrow\sA$ is
a covariant functor endowed with two natural transformations of
\emph{monad unit} $e\:\Id_\sA\rarrow\boM$ and
\emph{monad multiplication} $m\:\boM\circ\boM\rarrow\boM$ satisfying
the associativity and unitality axioms.
 We refer to~\cite[Chapter~VI]{McL} for the definitions of monads and
algebras over monads.

 Since we are interested in abelian categories arising from monads,
we prefer to use the term ``modules over a monad'' for what are usually
called algebras over a monad.
 Given a monad $\boM\:\sA\rarrow\sA$, we denote the category of
algebras/modules over $\boM$ (known also as the ``Eilenberg--Moore
category of~$\boM$'') by $\boM\modl$.
 The category $\boM\modl$ comes together with a faithful forgetful
functor $\boM\modl\rarrow\sA$, which has a left adjoint functor
$\sA\rarrow\boM\modl$.
 The latter functor takes any object $A\in\sA$ to the object
$\boM(A)\in\sA$, which is endowed with a natural structure of
an $\boM$\+module, making it an object of $\boM\modl$.
 The object $\boM(A)\in\boM\modl$ is called the \emph{free\/
$\boM$\+module} spanned by~$A$.

\begin{lem} \label{monad-additivity}
 Let\/ $\sA$ be an additive category and\/ $\boM\:\sA\rarrow\sA$ be
a monad.
 Then the following three conditions are equivalent:
\begin{enumerate}
\renewcommand{\theenumi}{\alph{enumi}}
\item the category\/ $\boM\modl$ is additive;
\item the category\/ $\boM\modl$ is additive, and both the forgetful
functor\/ $\boM\modl\rarrow\sA$ and the free module functor\/
$\sA\rarrow\boM\modl$ are additive;
\item the functor\/ $\boM\:\sA\rarrow\sA$ is additive.
\end{enumerate}
\end{lem}

\begin{proof}
 (a)\,$\Longleftrightarrow$\,(b)\,$\Longrightarrow$\,(c)
 The functor $\boM\:\sA\rarrow\sA$ is the composition of two adjoints
$\sA\rarrow\boM\modl\rarrow\sA$.
 Any functor between additive categories having an adjoint functor
is additive; and the composition of two additive functors is additive.

 (c)\,$\Longrightarrow$\,(b)  If the functor $\boM\:\sA\rarrow\sA$
is additive, then for any two $\boM$\+modules $P$ and $Q$
the set $\Hom_\boM(P,Q)$ of morphisms $P\rarrow Q$ in $\boM\modl$
is a subgroup in the group $\Hom_\sA(P,Q)$ of morphisms $P\rarrow Q$
in the category~$\sA$.
 This shows that the category $\boM\modl$ is preadditive and
the faithful forgetful functor $\boM\modl\rarrow\sA$ is additive.
 Finally, for any finite collection of objects $P_i\in\boM\modl$,
\ $i=1$,~\dots,~$n$, the direct sum ${}^\sA\bigoplus_{i=1}^n P_i$
taken in the category $\sA$ is naturally endowed with an $\boM$\+module
structure, which makes it the direct sum of the objects $P_i$ in
the category $\boM\modl$.
 Thus finite (co)products exist in $\boM\modl$.
\end{proof}

\begin{lem} \label{monad-exactness}
 Let\/ $\sA$ be an abelian category and\/ $\boM\:\sA\rarrow\sA$ be
a monad.
 Then the following two conditions are equivalent:
\begin{enumerate}
\renewcommand{\theenumi}{\alph{enumi}}
\item the category\/ $\boM\modl$ is abelian \emph{and} the forgetful
functor\/ $\boM\modl\rarrow\sA$ is exact;
\item the functor\/ $\boM\:\sA\rarrow\sA$ is right exact.
\end{enumerate}
\end{lem}

\begin{proof}
 Notice first of all that any exact or one-sided exact functor is
additive by definition; so Lemma~\ref{monad-additivity} is applicable.
 This allows us to presume that all the categories and functors involved
are additive.

 (a)\,$\Longrightarrow$\,(b)  Any left adjoint functor preserves
colimits; hence any left adjoint functor between abelian categories
is right exact.
 So the free module functor $\sA\rarrow\boM\modl$ is right exact.
 If the forgetful functor $\boM\modl\rarrow\sA$ is exact, it follows
that the composition $\boM\:\sA\rarrow\sA$ of two adjoints
$\sA\rarrow\boM\modl\rarrow\sA$ is right exact.

 (b)\,$\Longrightarrow$\,(a) For any monad $\boM\:\sA\rarrow\sA$ and
any morphism $f\:P\rarrow Q$ in $\boM\modl$, the kernel $\ker^\sA(f)$
of the morphism~$f$ taken in hte category $\sA$ is naturally endowed
with an $\boM$\+module structure, which makes it the kernel of~$f$
in the category $\boM\modl$.
 When the functor $\boM\:\sA\rarrow\sA$ is right exact, the same
holds true for the cokernel of~$f$.
 Both the properties in~(a) follow easily.
\end{proof}

 The following examples are illuminating.

\begin{exs} \label{co-monad-examples}
 (1)~Let $g\:R\rarrow S$ be an associative ring homomorphism.
 Put $\sA=R\modl$, and consider the functor $\boM=S\ot_R{-}\:
R\modl\rarrow R\modl$.
 The map~$g$ is an $R$\+$R$\+bimodule morphism, so it induces
a natural transformation $\Id_\sA\rarrow\boM$.
 The multiplication map $S\ot_RS\rarrow S$ is an $R$\+$R$\+bimodule
morphism, too, and it induces a natural transformation
$\boM\circ\boM\rarrow\boM$.
 These two morphisms of functors make $\boM$ a monad on $\sA=R\modl$.
 The category of $\boM$\+modules $\boM\modl$ is equivalent to $S\modl$.
 The forgetful functor $\boM\modl\rarrow\sA$ is the functor of
restriction of scalars $S\modl\rarrow R\modl$, while the free module
functor $\sA\rarrow\boM\modl$ is the functor of extension of
scalars $S\ot_R{-}\:R\modl\rarrow S\modl$.

 (2)~The notion of a \emph{comonad} on a category $\sA$ is dual to
that of a monad.
 Specifically, a comonad $\boC$ on $\sA$ is a covariant endofunctor
$\boC\:\sA\rarrow\sA$ endowed with two natural transformations of
\emph{counit} $\boC\rarrow\Id_\sA$ and \emph{comultiplication}
$\boC\rarrow\boC\circ\boC$ satisfying the coassociativity and
counitality axioms.
 One can say that a comonad on $\sA$ is the same thing as a monad on
$\sA^\sop$.
 Denote by $\boC\comodl$ the category of comodules (usually called
``coalgebras'') over a comonad $\boC$ on~$\sA$.
 
 Let $\sA$ be an abelian category and $\boC$ be a comonad on~$\sA$.
 Then the dual version of Lemma~\ref{monad-exactness} tells that
the category $\boC\comodl$ is abelian with exact forgetful functor
$\boC\comodl\rarrow\sA$ if and only if the functor $\boC\:\sA\rarrow\sA$
is left exact.

 (3)~Let $R$ be a ring and $C$ be a (coaccosiative, counital) coring
over $R$, that is, a comonoid object in the monoidal category of
$R$\+$R$\+bimodules.
 Then the functor $\boC={-}\ot_RC\:\modr R\rarrow \modr R$ is
a comonad on the category of right $R$\+modules $\modr R$.
 The counit map $C\rarrow R$ of the coring $C$ induces the counit
morphism $\boC\rarrow\Id_{\modr R}$, and the comultiplication map
$C\rarrow C\ot_RC$ induces the comultiplication morphism
$\boC\rarrow\boC\circ\boC$.
 Comodules over the comonad $\boC$ are known as \emph{right comodules}
over the coring~$C$.
 In other words, a right $C$\+comodule $N$ is a right $R$\+module
endowed with a \emph{right coaction} map $N\rarrow N\ot_RC$, which
must be a right $R$\+module morphism satisfying the coassociativity and
counitality equations involving the comultiplication and counit maps
of the coring~$C$.
 Similarly one defines \emph{left $C$\+comodules} (using the comonad
induced by $C$ on the category of left $R$\+modules).

 Let $\comodr C=\boC\comodl$ denote the category of right
$C$\+comodules.
 Then, as a particular case of~(2), one obtains the assertion that
the category $\comodr C$ is abelian with exact forgetful functor
$\comodr C\rarrow\modr R$ if and only if $C$ is a flat left
$R$\+module~\cite[Section~1.1.2]{Psemi}, \cite[Section~2.5]{Prev}.

 (4)~Let $R$ be a ring and $C$ be a coring over $R$, as in~(3).
 Then the functor $\boM=\Hom_R(C,{-})\:R\modl\rarrow R\modl$
is a monad on the category of left $R$\+modules.
 The counit map of the coring $C$ induces the unit morphism
$\Id_{R\modl}\rarrow\boM$, and the comultiplication map of
the coring $C$ induces the multiplication morphism
$\boM\circ\boM\rarrow\boM$.
 Modules over the monad $\boM$ are known as \emph{left contramodules}
over the coring~$C$.
 In other words, a left $C$\+contramodule $P$ is a left $R$\+module
endowed with a \emph{left contraaction} map $\Hom_R(C,P)\rarrow P$,
which must be a left $R$\+module morphism satisfying
(contra)associativity and (contra)unitality equations involving
the comultiplication and counit of the coring~$C$.
 We refer to~\cite[Section~3.1]{Psemi} or~\cite[Section~2.5]{Prev}
for details on contramodules over corings.

 Let $C\contra=\boM\modl$ denote the category of left
$C$\+contramodules.
 Then it follows from Lemma~\ref{monad-exactness} that the category
$C\contra$ is abelian with exact forgetful functor $C\contra\rarrow
R\modl$ if and only if $C$ is a projective left
$R$\+module~\cite[Section~3.1.2]{Psemi}, \cite[Proposition~2.5]{Prev}.
\end{exs}

\subsection{Ungraded contramodules over nonnegatively graded rings}
\label{ungraded-contramodules-subsecn}
 Let $B=\bigoplus_{n=0}^\infty B_n$ be a nonnegatively graded ring,
and let $K$ be an associative ring endowed with a ring homomorphism
$K\rarrow B_0$.
 Consider the following monad $\boM=\boM_K\:K\modl\rarrow K\modl$ on
the category of left $K$\+modules.

 To any left $K$\+module $L$, the monad $\boM$ assigns the left
$K$\+module
$$
 \boM_K(L)=\prod\nolimits_{n=0}^\infty \left(B_n\ot_K L\right).
$$
 The monad unit map
$$
 e_L\:L\lrarrow\prod\nolimits_{n=0}^\infty B_n\ot_K L
$$
is the composition
$$
 L\lrarrow B_0\ot_K L\lrarrow \prod\nolimits_{n=0}^\infty B_n\ot_K L
$$
of the map $L\rarrow B_0\ot_KL$ induced by the ring homomorphism
$K\rarrow B_0$ with the inclusion of the $(n=0)$\+indexed component
$B_0\ot_KL\rarrow\prod_{n=0}^\infty B_n\ot_KL$.
 The monad multiplication map
$$
 m_L\:\prod\nolimits_{p=0}^\infty B_p\ot_K
 \left(\prod\nolimits_{q=0}^\infty B_q\ot_K L\right)
 \lrarrow\prod\nolimits_{n=0}^\infty B_n\ot_K L
$$
is the composition
$$
 \prod\nolimits_{p=0}^\infty B_p\ot_K
 \left(\prod\nolimits_{q=0}^\infty B_q\ot_K L\right)\lrarrow
 \prod\nolimits_{p,q=0}^\infty B_p\ot_K B_q\ot_K L\lrarrow
 \prod\nolimits_{n=0}^\infty B_n\ot_K L.
$$
 Here the leftmost arrow is the product over $p\ge0$ of the maps
$B_p\ot_K\left(\prod_{q=0}^\infty B_q\ot_K L\right)\allowbreak\rarrow
\prod_{q=0}^\infty B_p\ot_K B_q\ot_K L$ whose components are
the direct summand projections $B_p\ot_K\left(\prod_{k=0}^\infty
B_k\ot_KL\right)\rarrow B_p\ot_K B_q\ot_KL$.
 The rightmost arrow is induced by the multiplication maps
$\prod_{p+q=n}^{p,q\ge0} B_p\ot_K B_q\rarrow B_n$.

 By the definition, a \emph{left $B$\+contramodule} (or
an \emph{ungraded left $B$\+contramodule}) is a module over
the monad $\boM_K\:K\modl\rarrow K\modl$.
 In other words, a left $B$\+contramodule $P$ is a left $K$\+module
endowed with a \emph{left $B$\+contraaction map}
$$
 \pi_P\:\boM_K(P)=\prod\nolimits_{n=0}^\infty B_n\ot_K P\lrarrow P,
$$
which must be a morphism of left $K$\+modules satisfying
the (contra)associativity and (contra)unitality equations involving
the multiplication and unit maps of the monad~$\boM_K$.
 This means that the compositions
$$
 \prod\nolimits_{i=0}^\infty B_i\ot_K\left(\prod\nolimits_{j=0}^\infty
 B_j\ot_K P\right)\,\rightrightarrows\,
 \prod\nolimits_{n=0}^\infty B_n\ot_K P\rarrow P
$$
of the monad multiplication map~$m_P$ and the map $\boM_K(\pi_P)$ with
the contraaction map~$\pi_P$ are equal to each other,
while the composition
$$
 P\lrarrow\prod\nolimits_{n=0}^\infty B_n\ot_KP\lrarrow P
$$
of the monad unit map~$e_P$ and the contraaction map~$\pi_P$ is equal
to the identity map~$\id_P$.

 We denote the category of ungraded left $B$\+contramodules by
$B\contra=\boM_K\modl$.
 Given two ungraded left $B$\+contramodules $P$ and $Q$, the group
of morphisms $P\rarrow Q$ in $B\contra$ is denoted by
$\Hom^B(P,Q)$.

 The next proposition shows that the notion of an ungraded left
$B$\+contramodule does not depend on the choice of a ring~$K$.
 One can consider two polar cases.
 On the one hand, one can take $K=\boZ$ and the unique ring
homomorphism $\boZ\rarrow B_0$.
 On the other hand, one can take $K=B_0$ and the identity map
$K\rarrow B_0$.

\begin{prop} \label{ungraded-B-contramodules-prop}
 There are natural equivalences of categories
$$
 \boM_\boZ\modl\simeq\boM_K\modl\simeq\boM_{B_0}\modl,
$$
making the notation $B\contra$ unambiguous.
 The category $B\contra$ of ungraded left $B$\+contramodules is abelian
with enough projective objects.
 The forgetful functor $B\contra\rarrow K\modl$ can be naturally lifted
to a forgetful functor $B\contra\rarrow B\modl$ taking values in
the category of ungraded left $B$\+modules.
 The forgetful functor $B\contra\rarrow B\modl$ is exact and preserves
infinite products.
\end{prop}

\begin{proof}
 To lift the forgetful functor $\boM_K\modl\rarrow K\modl$ to
a functor $\boM_K\modl\rarrow B\modl$, consider an $\boM_K$\+module
$P$ and restrict the contraction map~$\pi_P$ to the left
$K$\+submodule
$$
 B\ot_KP\,=\,\bigoplus\nolimits_{n=0}^\infty B_n\ot_KP\,\subset\,
 \prod\nolimits_{n=0}^\infty B_n\ot_K P.
$$
 The resulting map $B\ot_K P\rarrow P$ extends the left $K$\+module
structure on $P$ to a left $B$\+module structure, as explained in
Example~\ref{co-monad-examples}(1).

 The functor $\boM_K\:K\modl\rarrow K\modl$ is right exact, since
the infinite product functor is exact in $K\modl$ and the tensor
product functors $B_n\ot_K{-}$ are right exact.
 By Lemma~\ref{monad-exactness}, it follows that the category
$\boM_K\modl$ is abelian and the forgetful functor $\boM_K\modl\rarrow
K\modl$ is exact.
 The latter functor, being a right adjoint, also preserves all limits.
 Since the forgetful functor $B\modl\rarrow K\modl$ is faithul, exact,
and preserves all limits, it follows that the forgetful functor
$\boM_K\modl\rarrow B\modl$ is also exact and preserves all limits.

 The free object functor $K\modl\rarrow\boM_K\modl$, being left adjoint
to an exact functor, takes projective objects of $K\modl$ to projective
objects of $\boM_K\modl$.
 One can easily check that any object $P\in\boM_K\modl$ is a quotient
object of an $\boM_K$\+module $\boM(K^{(X)})$ obtained by applying
the free $\boM_K$\+module functor to a free left $K$\+module with $X$
generators, for some set~$X$.
 Indeed, if the underlying left $K$\+module of $P$ is a quotient of
a free left $K$\+module $K^{(X)}$, then the $\boM_K$\+module $P$ is
a quotient of the $\boM_K$\+module $\boM_K(K^{(X)})$.

 It remains to prove the first assertion.
 For this purpose, we will show that, for any $\boM_K$\+module $P$,
the contraaction map~$\pi_P$ factorizes through the natural
surjective map $\prod_{n=0}^\infty B_n\ot_KP\rarrow
\prod_{n=0}^\infty B_n\ot_{B_0}P$.
 Indeed, passing to the product over $n\ge0$ of the right exact
sequences (bar-complex fragments) $B_n\ot_K B_0\ot_K P\rarrow B_n\ot_K P
\rarrow B_n\ot_{B_0}P\rarrow0$, one obtains a right exact sequence
$$
 \prod\nolimits_{n=0}^\infty B_n\ot_K B_0\ot_KP \lrarrow
 \prod\nolimits_{n=0}^\infty B_n\ot_KP \lrarrow
 \prod\nolimits_{n=0}^\infty B_n\ot_{B_0}P \lrarrow0.
$$
 The map $\prod_{n=0}^\infty B_n\ot_K B_0\ot_K P
\rarrow\prod_{n=0}^\infty B_n\ot_K P$ is constructed as
the difference of two maps, one of which is induced by
the multiplication maps $B_n\ot_K B_0\rarrow B_n$ and
the other one by the action maps $B_0\ot_K P\rarrow P$.

 This pair of maps can be obtained by restricting the pair of maps
$m_P$, $\boM_K(\pi_P)\:\prod_{i=0}^\infty \left(B_i\ot_K
\left(\prod_{j=0}^\infty B_j\ot_KP\right)\right)\rightrightarrows
\prod_{n=0}^\infty B_n\ot_KP$ to the direct summand
$\prod_{i=0}^\infty B_i\ot_K B_0\ot_K P\subset\prod_{i=0}^\infty
\left(B_i\ot_K \left(\prod_{j=0}^\infty B_j\ot_KP\right)\right)$
corresponding to the value of the index $j=0$.
 Thus the contraassociativity equation for the $\boM_K$\+module $P$
implies the desired factorization.
 Now it is straightforward to check that the natural functors
$\boM_{B_0}\modl\rarrow\boM_K\modl\rarrow\boM_\boZ\modl$ are
category equivalences.  \hbadness=1175
\end{proof}

 Explicitly, the $\boM_K$\+module $\boM_K(K^{(X)})$ mentioned in
the above proof has the form
$$
 \boM_K(K^{(X)})=\prod\nolimits_{n=0}^\infty B_n\ot_K K^{(X)}
 =\prod\nolimits_{n=0}^\infty B_n^{(X)}=\Pi(B^{(X)}),
$$
where the notation $A^{(X)}$ stands for the direct sum of $X$ copies
of a (graded) abelian group~$A$.
 The left $B$\+contramodule $\boM_K(K^{(X)})$ does not depend on
the choice of a ring $K$, but only on the nonnegatively graded ring $B$
and a set~$X$.
 The $\boM_K$\+modules $\boM_K(K^{(X)})$ are called the \emph{free
ungraded left $B$\+contramodules}.
 Since there are enough free ungraded left $B$\+contramodules,
it follows that every projective ungraded left $B$\+contramodule is
a direct summand of a free one.

 The forgetful functor lifting assertion in
Proposition~\ref{ungraded-B-contramodules-prop} can be strengthened
as follows.
 In fact, the underlying left $B$\+module structure of any ungraded
left $B$\+contramodule can be extended naturally to a structure of
left module over the ring $\Pi B=\prod_{n=0}^\infty B_n$.
 To define the underlying left $\Pi B$\+module structure on a left
$B$\+contramodule $P\in\boM_K\modl$, it suffices to compose
the contraaction map~$\pi_P$ with the natural map
$$
 \Pi B\ot_K P\,=\,\left(\prod\nolimits_{n=0}^\infty B_n\right)\ot_K P
 \lrarrow\prod\nolimits_{n=0}^\infty\left(B_n\ot_K P\right) =\boM_K(P)
$$
and use Example~\ref{co-monad-examples}(1).

\begin{exs} \label{contramodules-examples}
 (1)~Let $B$ be a nonnegatively graded algebra over a field $k=B_0$
with finite-dimensional components $B_n$, as in
Example~\ref{comodules-examples}(1), and let $C=\bigoplus_{n=0}^\infty
B_n^*$ be the graded dual coalgebra over~$k$.
 Then an ungraded left $B$\+contramodule in the sense of the above
definition is the same thing as an ungraded left $C$\+contramodule
in the sense of~\cite[Section~1.1]{Prev} and~\cite[Appendix~A]{Psemi}.

 (2)~More generally, let $B$ be a nonnegatively graded ring with
the degree-zero component $R=B_0$ such that $B_n$ is a finitely
generated projective right $R$\+module for every~$n$, as
in Example~\ref{comodules-examples}(2).
 Let $C=\bigoplus_{n=0}^\infty C_n$, where $C_n=\Hom_{R^\rop}(B_n,R)$,
be the graded (right) dual coring.
 Then an ungraded left $B$\+contramodule in the sense of the above
definition is the same thing as an ungraded left $C$\+contramodule
in the sense of~\cite[Section~3.1.1]{Psemi}
and~\cite[Section~2.5]{Prev} (see Example~\ref{co-monad-examples}(4)).
\end{exs}

\begin{ex} \label{ungraded-co-contra-dualization-example}
 Let $E$ be an associative ring, $N$ be an $E$\+$B$\+bimodule, and
$U$ be a left $E$\+module.
 Assume that the right $B$\+module $N$ is a right $B$\+comodule.
 Then the left $B$\+module structure on the Hom group
$\Hom_E(N,U)$ can be extended to a left $B$\+contramodule structure
in a natural way.

 To construct the contraaction map
$$
 \pi\:\prod\nolimits_{n=0}^\infty B_n\ot_K \Hom_E(N,U)
 \lrarrow \Hom_E(N,U),
$$
consider an element $y\in N$, and choose an integer $m\ge0$ such that
$yB_{\ge m+1}=0$.
 Let $w\in\prod\nolimits_{n=0}^\infty B_n\ot_K\Hom_E(N,U)$ be
an arbitrary element.
 In order to evaluate the element $\pi(w)\in\Hom_E(N,U)$ on
the element $y\in N$, consider the composition
$$
 \prod\nolimits_{n=0}^\infty B_n\ot_K\Hom_E(N,U)\lrarrow
 \bigoplus\nolimits_{n=0}^m B_n\ot_K\Hom_E(N,U)\lrarrow
 \Hom_E(N,U)
$$
of the direct summand projection $\prod_{n=0}^\infty B_n\ot_K\Hom_E(N,U)
\rarrow\bigoplus_{n=0}^m B_n\ot_K\Hom_E(N,U)$ and the map
$\bigoplus_{n=0}^m B_n\ot_K\Hom_E(N,U)\rarrow\Hom_E(N,U)$ provided by
the left action of $B$ in $\Hom_E(N,U)$.
 This composition of maps needs to be applied to the element~$w$,
and the resuting element in $\Hom_E(N,U)$ needs to be evaluated on
the element $y\in N$, producing the desired element $\pi(w)(y)\in U$.
\end{ex}

\begin{rem} \label{two-definitions-of-contramodules-remark}
 Following the discussion in Section~\ref{ungraded-comodules-subsecn},
the hereditary pretorsion class $\comodr B\subset\modr B$ corresponds
to a right linear (in fact, two-sided linear) topology on the ring
$B=\bigoplus_{n=0}^\infty B_n$ with a base of neighborhoods of zero
formed by the ideals $B_{\ge m}=\bigoplus_{n\ge m} B_n$.
 The completion of $B$ with respect to this topology is the topological
ring $\Pi B=\prod_{n=0}^\infty B_n$ with the product topology.
 The abelian category $B\contra$ of ungraded left $B$\+contramodules,
as defined above, is in fact equivalent to the abelian category of left
contramodules over the topological ring $\Pi B$, as defined
in~\cite[Section~1.2]{Pweak}, \cite[Section~2.1]{Prev},
\cite[Sections~1.2 and~5]{PR}, \cite[Example~1.3(2)]{Pper},
\cite[Section~6.2]{PS1}, \cite[Section~2.7]{Pcoun}, etc.

 Let us explain the connection between the two definitions.
 The difference is that in the above references, unlike in the present
section, we were working with \emph{monads on the category of sets}.
 This allowed for an extra flexibility, in that one could easily define
contramodules over topological rings of arithmetic flavor, like
the $p$\+adic integers~\cite[Section~1.4]{Prev}.
 In the context of the present paper, such additional flexibility or
extra generality is not needed or relevant, so we are mostly restricting
ourselves to the more straightforward definition above.

 Still it is interesting and important for our purposes to establish
a comparison between the two definitions, i.~e., show that they are
indeed equivalent for the topological ring~$\Pi B$.
 This is provable along the lines of the argument
in~\cite[Section~1.10]{Pweak} and~\cite[Section~2.3]{Prev}.
 Let us spell out some details.
 
 The category of left contramodules over the topological ring~$\Pi B$
is, by the definition, the category of modules over the monad
$X\longmapsto\Pi B[[X]]$ on the category of sets, assigning to a set
$X$ the set of all infinite formal linear combinations of elements
of $X$ with the coefficients in $\Pi B$, where the family of
coefficients converges to zero in the topology of $\Pi B$.
 For any set $X$, we have a natural isomorphism of left $K$\+modules
$$
 \Pi B[[X]]\simeq\boM_K(K[X]).
$$
 For any left $K$\+module $A$, there is an obvious surjective map of
sets (or left $K$\+modules)
$$
 \gamma_A\:\Pi B[[A]]\lrarrow\boM_K(A).
$$
 Composing the map~$\gamma_P$ with the contraaction map $\pi_P\:
\boM_K(P)\rarrow P$, one defines a left $\Pi B$\+contramodule structure
on every ungraded left $B$\+contramodule~$P$.
 This construction defines a functor $B\contra\rarrow\Pi B\contra$,
where $\Pi B\contra$ denotes the category of left contramodules over
the topological ring~$\Pi B$.
 This functor is fully faithful, because the map $\gamma_P$ is
surjective.
 In order to show that this is an equivalence of categories,
it remains to check that, for any left $\Pi B$\+contramodule $C$,
the contraaction map $\Pi B[[C]]\rarrow C$ factorizes through
the surjection~$\gamma_C$.

 First of all, we recall that the underlying set of any contramodule
over the topological ring $\Pi B$ has a natural left $\Pi B$\+module
structure, hence also a natural left $K$\+module structure.
 The two forgetful functors $B\contra\rarrow\Pi B\modl$ and
$\Pi B\contra\rarrow\Pi B\modl$ form a commutative triangle diagram
with the above comparison functor $B\contra\rarrow\Pi B\contra$.

 For any set $X$, we denote by $K[X]=K^{(X)}$ the free left $K$\+module
with $X$~generators.
 Then for any left $K$\+module $A$ one can construct the monadic
bar-resolution related to the forgetful functor from $K\modl$
to the category of sets.
 We are interested in the fragment
\begin{equation} \label{ab-sets-resolution}
 K[K[A]]\rightrightarrows K[A]\rarrow A.
\end{equation}
 Here for any left $K$\+module $A$ there is a natural surjective map
$p_A\:K[A]\rarrow A$; this is the rightmost map
in~\eqref{ab-sets-resolution}.
 The leftmost pair of maps is formed by the maps $p_{K[A]}$
and $K[p_A]$.
 Passing to the difference of the leftmost pair of maps, one obtains
a right exact sequence of left $K$\+modules
\begin{equation} \label{ab-sets-resol-additive}
 K[K[A]]\lrarrow K[A]\lrarrow A\lrarrow0.
\end{equation}

 Applying the right exact functor $\boM_K$
to~\eqref{ab-sets-resolution} and~\eqref{ab-sets-resol-additive},
we obtain a right exact sequence of left $K$\+modules
\begin{equation} \label{kernel-of-gamma}
 \Pi B[[K[A]]]\lrarrow \Pi B[[A]]\lrarrow\boM_K(A)\lrarrow0,
\end{equation}
where the rightmost map in~\eqref{kernel-of-gamma} is~$\gamma_A$.
 We have obtained a description of the kernel of~$\gamma_A$ as
the image of the leftmost map in~\eqref{kernel-of-gamma},
which is the difference of two maps $\boM_K(p_{K[A]})$ and
$\boM_K(K[p_A])=\Pi B[[p_A]]$.

 Returning to a left $\Pi B$\+contramodule $C$, we recall that, by
the definition, its contraaction map $\Pi B[[C]]\rarrow C$ has to
satisfy the (contra)associativity equation of a module over
the monad $X\longmapsto\Pi B[[X]]$ on the category of sets.
 This equation, claiming that two compositions of maps of sets
$$
 \Pi B[[\Pi B[[C]]]]\,\rightrightarrows\,\Pi B[[C]]\rarrow C
$$
are equal to each other, can be also expressed by the vanishing of
the composition of two maps in the sequence of left $K$\+module
morphisms
\begin{equation} \label{topological-contraassociativity}
 \Pi B[[\Pi B[[C]]]]\lrarrow \Pi B[[C]]\lrarrow C\lrarrow 0.
\end{equation}
 It remains to observe that there is a natural commutative square
of a morphism from the leftmost map in~\eqref{kernel-of-gamma} for
$A=C$ to the leftmost map in~\eqref{topological-contraassociativity}.
 Hence the induced map on the rightmost terms $\boM_K(C)\rarrow C$,
as desired.
\end{rem}

 Let $P$ be an ungraded left $B$\+contramodule.
 One can define a natural decreasing filtration $P=G^0P\supset G^1P
\supset G^2P\supset\dotsb$ on $P$ by the rule
$$
 G^mP=\im\left(\left(\prod\nolimits_{n\ge m}B_n\ot_KP\right)
 \overset{\pi_P}\lrarrow P\right).
$$
 Here the infinite product $\prod_{n\ge m}B_n\ot_KP$ is viewed as
a $K$\+submodule in $\boM_K(P)=\prod_{n\ge0}B_n\ot_KP$,
the contraaction map $\pi_P\:\boM_K(P)\rarrow P$ is restricted onto
this submodule, and the image of the restriction is taken.

 It follows from the contraassociativity axiom for $\pi_P$ that
$G^mP$ is a $B$\+subcontra\-module in $P$ for every $m\ge0$ and
the successive quotient contramodules $G^mP/G^{m+1}P$ have
\emph{trivial} ungraded left $B$\+contramodule structures.
 Here an ungraded left $B$\+con\-tramodule $Q$ is said to be trivial
if the contraaction map $\pi_Q$ vanishes on the submodule
$\prod_{n=1}^\infty B_n\ot_KQ\subset\prod_{n=0}^\infty B_n\ot_KQ$
(in other words, this means that $G^1Q=0$).
 The full subcategory of trivial ungraded left $B$\+contramodules
in $B\contra$ is equivalent to the category of left $B_0$\+modules.

 An ungraded left $B$\+contramodule $P$ is said to be \emph{separated}
if\/ $\bigcap_{m\ge0}G^mP=0$, that is, in other words, the natural
$B$\+contramodule map $\lambda_{B,P}\:P\rarrow\varprojlim_{m\ge0}
P/G^mP$ is injective.
 A $B$\+contramodule $P$ is said to be \emph{complete} if the map
$\lambda_{B,P}$ is surjective.

 The following result is essentially well-known for contramodules over
a topological ring with a countable base of neighborhoods of zero.

\begin{prop} \label{separated-ungraded-contramodules}
 All ungraded left $B$\+contramodules are complete (but not necessarily
separated).
 All projective ungraded left $B$\+contramodules are separated.
 Every ungraded left $B$\+contramodule is the cokernel of an injective
morphism of separated ungraded left $B$\+contramodules.
\end{prop}

\begin{proof}
 By Remark~\ref{two-definitions-of-contramodules-remark}, we can
consider left contramodules over the topological ring $\Pi B$ in lieu of
the ungraded left $B$\+contramodules.
 Then the first assertion is~\cite[Lemma~6.3(b)]{PR}
(see also~\cite[Theorem~5.3]{Pcoun} or~\cite[Lemma~A.2.3 and
Remark~A.3]{Psemi}).
 The second assertion is clear from the explicit description of
projective ungraded left $B$\+contramodules above
(cf.~\cite[Lemma~6.9]{PR}).
 The third assertion follows from the second one together with existence
of sufficiently many projective contramodules, as all subcontramodules
of separated contramodules are separated (in particular, any
subcontramodule of a projective ungraded left $B$\+contramodule
is separated).
\end{proof}

\begin{thm} \label{ungraded-contramodules-fully-faithful}
 Assume that the augmentation ideal $B_{\ge1}=\bigoplus_{n=1}^\infty
B_n$ of a nonnegatively graded ring $B$ is finitely generated as
a right ideal in~$B$.
 Then the forgetful functor $B\contra\rarrow B\modl$ from the category
of ungraded left $B$\+contramodules to the category of ungraded left
$B$\+modules is fully faithful.
\end{thm}

\begin{proof}
 By Proposition~\ref{comodules-closed-under-extensions-prop},
the full subcategory $\comodr B\subset\modr B$ is closed under
extensions under our assumptions.
 In other words, $\comodr B$ is a hereditary torsion class in
$\modr B$.
 In the terminology of~\cite[Section~VI.5]{St}
and~\cite[Section~2.4]{Pcoun}, this means that the topology on $B$
with a base formed by the ideals $B_{\ge m}$, \,$m\ge1$, is
a \emph{right Gabriel topology}.
 By Lemma~\ref{augmentation-implies-component-ideals}, the right
ideals $B_{\ge m}\subset B$ are finitely generated.

 Now we can apply~\cite[Corollary~6.7]{Pcoun} in order to conclude
that the forgetful functor $\Pi B\contra\rarrow B\modl$ is fully
faithful.
 Alternatively, it is straightforward to see that the right ideal
$B_{\ge m}\subset B$ is strongly (finitely) generated in
the sense of~\cite[Section~6]{Pcoun};
so~\cite[Theorem~6.2\,(ii)\,$\Rightarrow$\,(iii)]{Pcoun}
(or even~\cite[Theorem~3.1]{Pper}) is directly applicable.
 Finally, Remark~\ref{two-definitions-of-contramodules-remark}
identifies the category $\Pi B\contra$ of left contramodules over
the topological ring $\Pi B$ with the category $B\contra$ which
we are interested in.
\end{proof}

\subsection{Weak Koszulity and the Ext comparison}
\label{weak-koszulity-subsecn}
 Let $B=\bigoplus_{n=0}^\infty B_n$ be a nonnegatively graded ring
with the degree-zero component $R=B_0$.
 We will say that $B$ is \emph{weakly right finitely projective Koszul}
if the graded right $B$\+module $R$ has a graded $B$\+module
resolution of the form
\begin{equation} \label{weakly-Koszul-resolution}
 \dotsb\lrarrow V_2\ot_R B\lrarrow V_1\ot_R B\lrarrow B
 \lrarrow R\lrarrow 0,
\end{equation}
where $V_i$, \,$i\ge1$, are finitely generated projective graded
right $R$\+modules.
 This means that $V_i=\bigoplus_{j\in\boZ} V_{i,j}$ is a graded right
$R$\+module concentrated in a finite set of internal degrees~$j$,
and for every $j\in\boZ$ the right $R$\+module $V_{i,j}$ is finitely
generated and projective.

 By the opposite version of Theorem~\ref{projective-koszul-theorem}(d)
(with the roles of the rings $A$ and $B$ switched), any right
finitely projective Koszul graded ring is weakly right finitely
projective Koszul.
 For any weakly right finitely projective Koszul graded ring $B$,
the augmentation ideal $B_{\ge1}=\bigoplus_{n=1}^\infty B_n$ is
finitely generated as a right ideal in $B$.
 In fact, any finite set of homogeneous generators of the graded
right $R$\+module $V_1$ gives rise to a finite set of homogeneous
generators of the right ideal $B_{\ge1}\subset B$.

\begin{lem} \label{Sigma-Pi-hom-tensor-lemma}
\textup{(a)} For any graded right $R$\+module $V$ and graded left
$R$\+module $H$, there is a natural map of ungraded abelian groups
$$
 \Sigma V\ot_R\Pi H\lrarrow \Pi(V\ot_R H).
$$
 This map is an isomorphism whenever $V$ is a finitely generated
projective graded right $R$\+module. \par
\textup{(b)} For any graded right $R$\+modules $V$ and $J$, there
is a natural map of ungraded abelian groups
$$
 \Sigma\Hom_{R^\rop}(V,J)\lrarrow \Hom_{R^\rop}(\Sigma V, \Sigma J).
$$
 This map is an isomorphism whenever $V$ is a finitely generated
projective graded right $R$\+module. \qed
\end{lem}

\begin{thm} \label{ungraded-contramodule-Ext-comparison-thm}
 For any weakly right finitely projective Koszul graded ring
$B=\bigoplus_{n=0}^\infty B_n$, the exact, fully faithful forgetful
functor of abelian categories $B\contra\rarrow B\modl$ induces
isomorphisms on all the Ext groups.
\end{thm}

\begin{proof}
 According to the discussion above, the augmentation ideal $B_{\ge1}
\subset B$ is finitely generated as a right ideal in~$B$.
 By Theorem~\ref{ungraded-contramodules-fully-faithful},
it follows that the forgetful functor $B\contra\rarrow B\modl$ is
indeed fully faithful.
 We have to show that it also induces isomorphisms of the groups
$\Ext^i$ for $i>0$.

 There are enough projective objects in the abelian category $B\contra$.
 Hence it suffices to prove that $\Ext_B^i(P,Q)=0$ for all projective
objects $P\in B\contra$, all objects $Q\in B\contra$, and all $i>0$
(where, as usually, $\Ext_B^*$ denotes the Ext functor in the category
$B\modl$).
 By Proposition~\ref{separated-ungraded-contramodules}, any object of
$B\contra$ is the cokernel of an injective morphism of separated
contramodules.
 Hence we can assume that $Q$ is a separated ungraded left
$B$\+contramodule.

 Any separated ungraded left $B$\+contramodule $Q$ is an infinitely
iterated extension, in the sense of projective limit, of trivial
ungraded left $B$\+contramodules $G^mQ/G^{m+1}Q$ (because, by the same
proposition, all ungraded left $B$\+contramodules are complete).
 In view of the dual Eklof lemma~\cite[Proposition~18]{ET}, we can
assume that $Q$ is a trivial ungraded left $B$\+contramodule
(that is, just a left $R$\+module in which the augmentation ideal
$B_{\ge1}\subset B$ acts by zero).

 Now the left $R$\+module $Q$ has a resolution by injective left
$R$\+modules of the form $F^+=\Hom_\boZ(F,\boQ/\boZ)$, where $F$
ranges over free right $R$\+modules.
 This reduces the question to the case of a left $R$\+module of
the form $Q=F^+$.
 Furthermore, we have $\Ext_B^i(P,F^+)\simeq\Tor^B_i(F,P)^+$
(where $B_{\ge1}$ acts by zero in~$F$).
 Thus it suffices to show that $\Tor^B_i(R,P)=0$ for all projective
objects $P\in B\contra$ and all integers $i>0$ (where $R$ is viewed
as a right $B$\+module with $B_{\ge1}$ acting by zero).

 We have done a chain of reductions for an object $Q\in B\contra$;
now we need to deal with a projective object $P\in B\contra$.
 According to the discussion after
Proposition~\ref{ungraded-B-contramodules-prop} in
Section~\ref{ungraded-contramodules-subsecn}, the underlying left
$B$\+module of $P$ is a direct summand of an ungraded left $B$\+module
$\Pi H$ for some free graded left $B$\+module $H=B^{(X)}$ with
generators in degree zero.

 Finally, we compute the groups $\Tor^B_i(R,\Pi H)$ using the projective
resolution~\eqref{weakly-Koszul-resolution} of the right
$B$\+module~$R$.
 More precisely, we use the ungraded projective resolution obtained by
applying $\Sigma$ to~\eqref{weakly-Koszul-resolution}.
 We obtain the complex of ungraded abelian groups
\begin{equation} \label{Sigma-Pi-Tor-complex}
 \dotsb\lrarrow \Sigma V_2\ot_R\Pi H\lrarrow \Sigma V_1\ot_R\Pi H
 \lrarrow \Pi H\lrarrow0
\end{equation}
computing $\Tor^B_*(R,\Pi H)$.
 According to Lemma~\ref{Sigma-Pi-hom-tensor-lemma}(a),
the complex~\eqref{Sigma-Pi-Tor-complex} is isomorphic to
\begin{equation} \label{Pi-Tor-complex-computed}
 \dotsb\lrarrow \Pi(V_2\ot_R H) \lrarrow \Pi (V_1\ot_R H)
 \lrarrow \Pi(H)\lrarrow0.
\end{equation}
 The complex~\eqref{Pi-Tor-complex-computed} can be obtained by
applying $\Pi$ to the complex
\begin{equation}
 \dotsb\lrarrow V_2\ot_R H \lrarrow V_1\ot_R H \lrarrow H\lrarrow0
\end{equation}
computing the graded abelian groups $\Tor^B_*(R,H)$.
 It remains to observe that $\Tor^B_i(R,H)=0$ for $i>0$ and
the functor $\Pi$ is exact.
 Thus $\Tor^B_i(R,\Pi H)\simeq\Pi\Tor^B_i(R,H)=0$ for $i>0$ under
our assumptions.
\end{proof}

 The assertion of the next theorem can be equivalently restated by
saying that, under its assumptions, the full subcategory of right
$B$\+comodules $\comodr B$ is a \emph{weakly stable torsion class} in
$\modr B$ in the sense of the paper~\cite{VY}.

\begin{thm} \label{ungraded-comodule-Ext-comparison-thm}
 For any weakly right finitely projective Koszul graded ring
$B=\bigoplus_{n=0}^\infty B_n$, the exact, fully faithful inclusion
functor of abelian categories\/ $\comodr B\rarrow\modr B$ induces
isomorphisms on all the Ext groups. \hbadness=1750
\end{thm}

\begin{proof}
 There are enough injective objects in a Grothendieck abelian category
of ungraded right $B$\+comodules $\comodr B$.
 Hence it suffices to prove that $\Ext_{B^\rop}^i(M,L)=0$ for all
objects $M\in\comodr B$, all injective objects $L\in\comodr B$, and
all $i>0$ (where $\Ext_{B^\rop}^*$ denotes the Ext functor in
the category $\modr B$).

 Any object $M\in\comodr B$ has a natural increasing filtration
$0=G_{-1}M\subset G_0 M\subset G_1 M\subset G_2M\subset \dotsb$
by right $B$\+submodules $G_mM\subset M$, where $G_mM$ is the subset
of all elements $x\in M$ such that $xB_{\ge m+1}=0$.
 By the definition, we have $M=\bigcup_{m\ge0}G_mM$.
 The successive quotient modules $G_mM/G_{m-1}M$ are \emph{trivial}
ungraded right $B$\+comodules, i.~e., right $B$\+modules in which
the augmentation ideal $B_{\ge1}$ acts by zero.
 The category of trivial ungraded right $B$\+comodules is equivalent
to the category of right $R$\+modules $\modr R$.

 In view of the Eklof lemma~\cite[Lemma~1]{ET}, we can assume that
$M$ is a trivial ungraded right $B$\+comodule.
 Now the right $R$\+module $M$ has a resolution by free right
$R$\+modules, which reduces the question to the case of a free right
$R$\+module $M=R^{(X)}$ (with $B_{\ge1}$ acting by zero in~$M$).
 Thus it suffices to show that $\Ext_{B^\rop}^i(R,L)=0$ for all
injective objects $L\in\comodr B$ and all integers $i>0$ (where $R$
is viewed as a right $B$\+module with $B_{\ge1}$ acting by zero).

 Injective objects of the category $\comodr B$ can be described as
follows.
 The inclusion functor $\comodr B\rarrow\modr B$ has a right adjoint
functor $\Gamma_B\:\modr B\rarrow\comodr B$ assigning to an ungraded
right $B$\+module $N$ its submodule $\Gamma_B(N)\subset N$ consisting
of all the elements $x\in N$ for which there exists $m\ge1$ such that
$xB_{\ge m}=0$.
 The injective objects of $\comodr B$ are precisely the direct summands
of the objects $\Gamma_B(N)$, where $N$ ranges over the injective
objects of $\modr B$.

 Injective ungraded right $B$\+modules can be further described as
follows.
 Let $H=B^{(X)}$ be a free graded left $B$\+module with generators in
degree zero, and let $H^+=\Hom_\boZ(B^{(X)},\boQ/\boZ)$ be the graded
character module of~$H$.
 So $H^+$ is a graded right $B$\+module concentrated in the nonpositive
degrees.
 Then $\Pi(H^+)$ is an injective ungraded right $B$\+module, and all
injective right $B$\+modules are direct summands of $B$\+modules of
this form.
 One can easily see that $\Gamma_B(\Pi(H^+))=\Sigma(H^+)$.
 
 Finally, we compute the groups $\Ext_{B^\rop}^i(R,\Sigma(H^+))$ using
the projective resolution~\eqref{weakly-Koszul-resolution} of
the right $B$\+module~$R$.
 More precisely, we use the projective resolution of ungraded right
$B$\+module $R$ obtained by applying $\Sigma$
to~\eqref{weakly-Koszul-resolution}.
 We obtain a complex of ungraded abelian groups
\begin{equation} \label{Sigma-Sigma-Ext-complex}
 0\lrarrow\Sigma(H^+)\lrarrow\Hom_{R^\rop}(\Sigma V_1,\Sigma(H^+))
 \lrarrow\Hom_{R^\rop}(\Sigma V_2,\Sigma(H^+))\lrarrow\dotsb
\end{equation}
computing $\Ext_{B^\rop}^*(R,\Sigma(H^+))$.
 According to Lemma~\ref{Sigma-Pi-hom-tensor-lemma}(b),
the complex~\eqref{Sigma-Sigma-Ext-complex} is isomorphic to
\begin{equation} \label{Sigma-Ext-complex-computed}
 0\lrarrow\Sigma(H^+)\lrarrow\Sigma\Hom_{R^\rop}(V_1,H^+)
 \lrarrow\Sigma\Hom_{R^\rop}(V_2,H^+)\lrarrow\dotsb
\end{equation} 

 The complex~\eqref{Sigma-Ext-complex-computed} can be obtained by
applying $\Sigma$ to the complex of graded abelian groups
\begin{equation}
 0\lrarrow H^+\lrarrow\Hom_{R^\rop}(V_1,H^+)\lrarrow
 \Hom_{R^\rop}(V_2,H^+)\lrarrow\dotsb
\end{equation}
computing the graded $\Ext_{B^\rop}^*(R,H^+)$ (that is, the derived
functor of graded $\Hom$ of graded right $B$\+modules).
 It remains to observe that $H^+$ is an injective object of
the category of graded right $B$\+modules and the functor $\Sigma$
is exact.
 Thus $\Ext_{B^\rop}^i(R,\Sigma(H^+))\simeq
\Sigma\Ext_{B^\rop}^i(R,H^+)=0$ for $i>0$.
\end{proof}

\subsection{Graded comodules over nonnegatively graded rings}
\label{graded-comodules-subsecn}
 Let $B=\bigoplus_{n=0}^\infty B_n$ be a nonnegatively graded ring.
 We consider graded right $B$\+modules $M=\bigoplus_{n\in\boZ}M_n$.
 A graded right $B$\+module $M$ is said to be a \emph{$B$\+comodule}
(or a \emph{graded right $B$\+comodule}) if its underlying ungraded
right $B$\+module $\Sigma M$ is an ungraded right $B$\+comodule.
 Equivalently, this means that for every (homogeneous) element
$x\in M$ there exists an integer $m\ge0$ such that $xb=0$ in $M$
for all $b\in B_n$, \,$n>m$.

 Clearly, any \emph{nonpositively} graded right $B$\+module
$M=\bigoplus_{n=-\infty}^0 M_n$ is a $B$\+comod\-ule.
 But the free graded right $B$\+module $B$ is usually \emph{not}
a $B$\+comodule (in fact, $B$ is a $B$\+comodule if and only if
$B_n=0$ for $n\gg0$).

 We denote the category of graded right $B$\+modules by
$\modrgr B$ and the full subcategory of graded right $B$\+comodules
by $\comodrgr B\subset\modrgr B$.
 Obviously, the full subcategory $\comodrgr B\subset\modrgr B$ is
closed under subobjects, quotients, and infinite direct sums (in other
words, $\comodrgr B$ is a hereditary pretorsion class in $\modrgr B$).
 Hence $\comodrgr B$ is a Grothendieck abelian category and its
inclusion $\comodrgr B\rarrow\modrgr B$ is an exact functor preserving
infinite direct sums.

\begin{ex} \label{graded-comodules-example}
 Let $B$ be a nonnegatively graded ring with the degree-zero grading
component $R=B_0$.
 Assume that $B_n$ is a finitely generated projective right
$R$\+module for every $n\ge0$, and consider the graded right dual
coring $C$ as in Example~\ref{comodules-examples}(2).
 Then a graded right $B$\+comodule in the sense of the above definition
is the same thing as a graded right $C$\+comodule.
\end{ex}

\begin{prop} \label{graded-comodules-closed-under-extensions}
 Assume that the augmentation ideal $B_{\ge1}=
\bigoplus_{n=1}^\infty B_n$ of a nonnegatively graded ring $B$
is finitely generated as a right ideal in~$B$.
 Then the full subcategory of graded right $B$\+comodules\/
$\comodrgr B$ is closed under extensions in the category of
graded right $B$\+modules\/ $\modrgr B$.
\end{prop}

\begin{proof}
 This is a graded version of
Proposition~\ref{comodules-closed-under-extensions-prop}.
 It can be either proved by the same computation, or deduced
formally from the ungraded version.
\end{proof}

\subsection{Graded contramodules over nonnegatively graded rings}
\label{graded-contramodules-subsecn}
 Let $B=\bigoplus_{n=0}^\infty B_n$ be a nonnegatively graded ring,
and let $K$ be an associative ring endowed with a ring homomorphism
$K\rarrow B_0$.
 Consider the following monad $\boM=\boM_K^\gr\:K\modl_\sgr\rarrow
K\modl_\sgr$ on the category of graded left $K$\+modules $K\modl_\sgr$
(where $K$ is viewed as a graded ring concentrated in the grading~$0$).

 To any graded left $K$\+module $L=(L_j)_{j\in\boZ}$, the monad $\boM$
assigns the graded left $K$\+module with the components
$$
 \boM_K^\gr(L)_j=\prod\nolimits_{n=0}^\infty
 \left(B_n\ot_K L_{j-n}\right).
$$
 The monad unit map $e_L\:L\rarrow\boM_K^\gr(L)$ is the graded
$K$\+module morphism whose degree~$j$ component is the composition
$$
 L_j\lrarrow B_0\ot_K L_j\lrarrow
 \prod\nolimits_{n=0}^\infty B_n\ot_K L_{j-n}
$$
of the map $L_j\rarrow B_0\ot_KL_j$ induced by the ring homomorphism
$K\rarrow B_0$ with the inclusion of the $(n=0)$\+indexed summand
$B_0\ot_K L_j\rarrow\prod_{n=0}^\infty B_n\ot_KL_{j-n}$.
 The monad multiplication map $m_L\:\boM_K^\gr(\boM_K^\gr(L))\rarrow
\boM_K^\gr(L)$ is the graded $K$\+module morphism whose degree~$j$
component is the composition
$$
 \prod_{p=0}^\infty B_p\ot_K
 \left(\prod_{q=0}^\infty B_q\ot_K L_{j-p-q}\right) \lrarrow
 \prod_{p,q=0}^\infty B_p\ot_K B_q\ot_K L_{j-p-q} \lrarrow
 \prod_{n=0}^\infty B_n\ot_K L_{j-n}.
$$
 Here the leftmost arrow is the product over $p\ge0$ of the maps
$$
 B_p\ot_K\left(\prod\nolimits_{q=0}^\infty B_q\ot_K L_{j-p-q}\right)
 \lrarrow\prod\nolimits_{q=0}^\infty B_p\ot_K B_q\ot_K L_{j-p-q}
$$
whose $q$\+components are the direct summand projections
$B_p\ot_K\left(\prod_{k=0}^\infty B_k\ot_K L_{j-p-k}\right)
\allowbreak\rarrow B_p\ot_K B_q\ot_K L_{j-p-q}$.
 The rightmost arrow is induced by the multiplication maps
$\prod_{p+q=n}^{p,q\ge0} B_p\ot_K B_q\rarrow B_n$.

 By the definition, a \emph{graded left $B$\+contramodule} is
a module over the monad $\boM_K^\gr\:K\modl_\sgr\rarrow K\modl_\sgr$.
 In other words, a graded left $B$\+contramodule $P$ is a graded
left $K$\+module endowed with a \emph{left $B$\+contraaction map}
$$
 \pi_P\:\boM_K^\gr(P)\lrarrow P,
$$
which must be a morphism of graded left $K$\+modules satisfying
the (contra)associativity and (contra)unitality equations involving
the multiplication and unit maps of the monad~$\boM_K^\gr$.
 This means that, in every degree $j\in\boZ$, the two {compositions
\hbadness=1850
$$
 \prod\nolimits_{p=0}^\infty B_p\ot_K
 \left(\prod\nolimits_{q=0}^\infty B_q\ot_K P_{j-p-q}\right)
 \,\rightrightarrows\,\prod\nolimits_{n=0}^\infty B_n\ot_K P_{j-n}
 \rarrow P_j
$$
of} the degree~$j$ components of the monad multiplication map~$m_P$
and the map $\boM_K^\gr(\pi_P)$ with the degree~$j$ component
of the contraaction map~$\pi_P$ are equal to each other, while
the composition
$$
 P_j\lrarrow\prod\nolimits_{n=0}^\infty B_n\ot_K P_{j-n}
 \lrarrow P_j
$$
of the degree~$j$ components of the monad unit map~$e_P$ and
the contraaction map~$\pi_P$ is equal to the identity map~$\id_{P_j}$.

 We denote the category of graded left $B$\+contramodules by
$B\contra_\sgr=\boM_K^\gr\modl$.
 The category of graded left $B$\+modules $M=\bigoplus_{j\in\boZ}M_j$
is denoted by $B\modl_\sgr$.

 For any graded left $B$\+contramodule $P$, the ungraded
left $K$\+module $\Pi P=\prod_{j\in\boZ}P_j$ has a natural structure
of ungraded left $B$\+contramodule.
 Indeed, the functors $\boM_K^\gr\:K\modl_\sgr\rarrow K\modl_\sgr$
and $\boM_K\:K\modl\rarrow K\modl$ form a commutative square
diagram with the functor $\Pi\:K\modl_\sgr\rarrow K\modl$; and
the monad multiplication and unit maps of the two monads agree.
 So we obtain a faithful functor of forgetting the grading
$$
 \Pi\:B\contra_\sgr\lrarrow B\contra.
$$
 It is clear from the next proposition that the functor $\Pi\:
B\contra_\sgr\rarrow B\contra$ is exact and preserves infinite
products (because the functor $\Pi\:K\modl_\sgr\rarrow K\modl$
has such properties).

\begin{prop} \label{graded-B-contramodules-proposition}
 There are natural equivalences of categories
$$
 \boM_\boZ^\gr\modl\simeq\boM_K^\gr\modl\simeq\boM_{B_0}^\gr\modl,
$$
making the notation $B\contra_\sgr$ unambiguous.
 The category $B\contra_\sgr$ of graded left $B$\+contramodules is
abelian with enough projective objects.
 The forgetful functor $B\contra_\sgr\rarrow K\modl_\sgr$ can be
naturally lifted to a forgetful functor $B\contra_\sgr\rarrow
B\modl_\sgr$ taking values in the category of graded left $B$\+modules.
 The forgetful functor $B\contra_\sgr\rarrow B\modl_\sgr$ is exact
and preserves infinite products.
\end{prop}

\begin{proof}
 This is a graded version of
Proposition~\ref{ungraded-B-contramodules-prop}.
 To prove the first assertion, notice that for every graded left
$B_0$\+module $L$ there is a natural surjective map of graded
left $K$\+modules $\boM_K^\gr(L)\rarrow\boM_{B_0}^\gr(L)$, and
for every graded left $K$\+module $L$ there is a natural surjective
map of graded abelian groups $\boM_\boZ^\gr(L)\rarrow\boM_K^\gr(L)$.
 Using these maps, one constructs a natural $\boM_K^\gr$\+module
structure on the underlying graded left $K$\+module of every
$\boM_{B_0}^\gr$\+module, and a natural $\boM_\boZ^\gr$\+module
structure on the underlying graded abelian group of every
$\boM_K^\gr$\+module.
 The resulting functors $\boM_{B_0}^\gr\modl\rarrow\boM_K^\gr\modl
\rarrow\boM_\boZ^\gr\modl$ are fully faithful, since the above
maps of graded modules are surjective.

 To prove these functors are essentially surjective, it remains to
show that, for every $\boM_K^\gr$\+module $P$, the contraaction map
$\pi_P\:\boM_K^\gr(P)\rarrow P$ factorizes through the natural
surjection $\boM_K^\gr(P)\rarrow\boM_{B_0}^\gr(P)$.
 This can be either checked directly in the way similar to the proof
of the ungraded version in
Proposition~\ref{ungraded-B-contramodules-prop}, or deduced formally
from the ungraded version.
 The latter approach works as follows: the map
$\Pi\pi_P\:\Pi\boM_K^\gr(P)\rarrow\Pi P$ is naturally identified
with the contraaction map $\pi_{\Pi P}\:\boM_K(\Pi P)\rarrow\Pi P$.
 The latter map factorizes though $\boM_{B_0}(\Pi P)$, as explained
in the proof of Proposition~\ref{ungraded-B-contramodules-prop};
this means that the map $\Pi\pi_P$ factorizes though
$\Pi\boM_{B_0}^\gr(P)$.
 Hence the map~$\pi_P$ factorizes through $\boM_{B_0}^\gr(P)$.

 The remaining assertions are provable in the way similar to
Proposition~\ref{ungraded-B-contramodules-prop}.
 It is important that the functor $\boM_K^\gr\:K\modl_\sgr\rarrow
K\modl_\sgr$ is right exact (as a combination of tensor products
and direct products).
 We skip the obvious details, and will only explain what
the projective objects of the category $B\contra_\sgr$ are.
 Enough of these can be obtained by applying the free
$\boM_K^\sgr$\+module functor $K\modl_\sgr\rarrow\boM_K^\gr\modl$
to free graded left $K$\+modules (with generators possibly sitting
in all the degrees).

 Explicitly, let $X$ be a \emph{graded set}, which means a sorted set
with the sorts indexed by the integers, or equivalently, a set presented
as a disjoint union $X=\coprod_{i\in\boZ} X_i$.
 Let $K[X]$ denote the graded left $K$\+module with the components
$K[X]_i=K^{(X_i)}$.
 Then, for every $j\in\boZ$, the degree~$j$ component of the graded
left $B$\+contramodule $\boM_K^\sgr(K[X])$ has the form
$$
 \boM_K^\sgr(K[X])_j=
 \prod\nolimits_{n=0}^\infty B_n\ot_K K[X]_{j-n}=
 \prod\nolimits_{n=0}^\infty B_n^{(X_{j-n})}.
$$
 The graded left $B$\+contramodule $\boM_K(K[X])$ does not depend
on the choice of a ring $K$, but only on the nonnegatively graded
ring $B$ and a graded set~$X$.
 The $\boM_K^\gr$\+module $\boM_K(K[X])$ is called the \emph{free
graded left $B$\+contramodule} spanned by a graded set~$X$.
 Since there are enough free graded left $B$\+contamodules, it follows
that every projective graded left $B$\+contramodule is a direct
summand of a free one.
\end{proof}

\begin{exs} \label{graded-contramodules-examples}
 (1)~Let $B$ be a nonnegatively graded algebra over a field $k=B_0$
with finite-dimensional components~$B_n$, as in
Examples~\ref{comodules-examples}(1)
and~\ref{contramodules-examples}(1), and let
$C=\bigoplus_{n=0}^\infty B_n^*$ be the graded dual coalgebra over~$k$.
 Then a graded left $B$\+con\-tramodule in the sense of the above
definition is the same thing as a graded left $C$\+contramodule
in the sense of~\cite[Section~2.2]{Pkoszul}.
 (We ignore the issue of the sign rule, which only becomes important
when one considers differentials on coalgebras, comodules, and
contramodules.)

 (2)~More generally, let $B$ be a nonnegatively graded ring with
the degree-zero component $R=B_0$ such that $B_n$ is a finitely
generated projective right $R$\+module for every~$n$, and let
$C=\bigoplus_{n=0}^\infty\Hom_{R^\rop}(B_n,R)$ be the graded (right)
dual coring, as in Examples~\ref{comodules-examples}(2)
and~\ref{contramodules-examples}(2).
 Then a graded left $B$\+contramodule in the sense of the above
definition is the same thing as a graded left $C$\+contramodule
in the sense of~\cite[Section~11.1.1]{Psemi}.
\end{exs}

\begin{exs} \label{graded-contramodules-bounded-grading-examples}
 (1)~Any \emph{nonnegatively} graded left $B$\+module
$P=\bigoplus_{j=0}^\infty P_j$ has a unique graded left
$B$\+contramodule structure compatible with its $B$\+module structure.
 The point is that the relevant direct products are essentially finite
in this case, and therefore coincide with the similar direct sums,
$$
 \prod\nolimits_{n=0}^\infty B_n\ot_K P_{j-n}=
 \bigoplus\nolimits_{n=0}^\infty B_n\ot_K P_{j-n}
 \qquad\text{for all $j\in\boZ$}.
$$

 (2)~For the same reason, there is no difference between graded left
$B$\+modules and graded left $B$\+contramodules when $B_n=0$
for $n\gg0$.
 (There is also no difference between ungraded left $B$\+modules and
ungraded left $B$\+contramodules in this case; and similarly for
graded or ungraded right $B$\+comodules.)

 (3)~However, outside of the degenerate case~(2), the cofree graded
left $B$\+module $B^+=\Hom_\boZ(B,\boQ/\boZ)$ usually does not admit
an extension of its graded $B$\+module structure to a graded
$B$\+contramodule structure.
 In order to prove this for a concrete graded ring $B$, it suffices
to show that the cofree ungraded left $B$\+module $\Pi(B^+)$ does not
admit an exension of its $B$\+module structure to an ungraded
$B$\+contramodule structure.
 Let us give a specific example.

 Choose a field~$k$, and let $B$ be the polynomial algebra $B=k[x]$
in one variable~$x$ of the degree $\deg x=1$.
 Then the (graded or ungraded) $B$\+contramodules are the same thing as
the (respectively, graded or ungraded) contramodules over the coalgebra
$C$ over~$k$ graded dual to~$B$ (as in
Examples~\ref{contramodules-examples}(1)
and~\ref{graded-contramodules-examples}(1)).
 Ungraded $C$\+contramodules are described
in~\cite[Remark~A.1.1]{Psemi}, \cite[Sections~1.3 and~1.6]{Prev},
\cite[Theorem~B.1.1]{Pweak}, and~\cite[Theorem~3.3]{Pcta} as
$k[x]$\+modules $P$ such that
$\Hom_{k[x]}(k[x,x^{-1}],P)=0=\Ext^1_{k[x]}(k[x,x^{-1}],P)$.
 In particular, it follows that $P$ has no nonzero $x$\+divisible
$k[x]$\+submodules~\cite[Lemma~3.2]{Pcta}.
 Therefore, the $k[x]$\+module $\Pi(B^+)$ is not a $C$\+contramodule,
because it is nonzero and $x$\+divisible (i.~e., the operator
$x\:\Pi(B^+)\rarrow\Pi(B^+)$ is surjective).
\end{exs}

\begin{ex} \label{graded-co-contra-dualization-example}
 Let $E=\bigoplus_{n\in\boZ}E_n$ be a graded associative ring,
$N=\bigoplus_{n\in\boZ}N_n$ be a graded $E$\+$B$\+bimodule, and
$U=\bigoplus_{n\in\boZ}U_n$ be a graded left $E$\+module.
 Assume that the graded right $B$\+module $N$ is a graded right
$B$\+comodule.
 Then the graded left $B$\+module structure on the graded Hom group
$\Hom_E(N,U)$ can be extended to a graded left $B$\+contramodule
structure in a natural way.
 The construction of the contraaction is similar to that in
Example~\ref{ungraded-co-contra-dualization-example}.

 Specifically, for every $j\in\boZ$ we need to construct
the degree~$j$ component of the contraaction map
$$
 \pi_j\:\prod\nolimits_{n=0}^\infty B_n\ot_K\Hom_E(N,U)_{j-n}
 \lrarrow\Hom_E(N,U)_j.
$$
 Let $w\in\prod_{n=0}^\infty B_n\ot_K\Hom_E(N,U)_{j-n}$ be
an arbitrary element.
 Consider an element $y\in N_i$, \ $i\in\boZ$, and choose an integer
$m\ge0$ such that $yB_{\ge m+1}=0$.
 In order to evaluate element $\pi_j(w)\in\Hom_E(N,U)_j$ on
the element $y\in N_i$, consider the composition
$$
 \prod_{n=0}^\infty B_n\ot_K\Hom_E(N,U)_{j-n}
 \lrarrow\bigoplus_{n=0}^m B_n\ot_K\Hom_E(N,U)_{j-n}
 \lrarrow\Hom_E(N,U)_j
$$
of the direct summand projection $\prod_{n=0}^\infty B_n\ot_K
\Hom_E(N,U)_{j-n}\rarrow\bigoplus_{n=0}^m B_n\ot_K\Hom_E(N,U)_{j-n}$
and the map $\bigoplus_{n=0}^m B_n\ot_K\Hom_E(N,U)_{j-n}
\rarrow\Hom_E(N,U)_j$ provided by the left action of $B$
in $\Hom_E(N,U)$.
 This composition of maps has to be applied to the element~$w$,
and the resulting element in $\Hom_E(N,U)_j$ has to be evaluated
on the element $y\in N_i$, producing the desired element
$\pi_j(w)(y)\in U_{i+j}$.
\end{ex}

 Let $P$ be a graded left $B$\+contramodule.
 The natural decreasing filtration $P=G^0P\supset G^1P\supset G^2P
\supset\dotsb$ on $P$ is defined in the following way.
 For every $m\ge0$, the graded left $K$\+submodule $G^mP\subset P$
has the grading components
$$
 G^mP_j=\im\left(\left(\prod\nolimits_{n\ge m}B_n\ot_KP_{j-n}\right)
 \overset{\pi_P}\lrarrow P_j\right).
$$
 Here the degree~$j$ component of the contraaction map
$\pi_{P,j}\:\boM_K^\gr(P)_j\rarrow P_j$ is restricted onto
the $K$\+submodule (in fact, a direct summand) $\prod_{n\ge m}
B_n\ot_KP_{j-n}\subset\prod_{n=0}^\infty B_n\ot_KP_{j-n}$,
and the image of the restriction is taken.

 The key observation is that $\bigoplus_{n\ge m}B_m$ is a homogeneous
two-sided ideal in~$B$.
 Hence it follows from the contraassociativity axiom for~$\pi_P$
that $G^mP$ is a graded $B$\+subcontramodule in $P$ for every $m\ge0$
and the successive quotient graded contramodules $G^mP/G^{m+1}P$
have \emph{trivial} graded left $B$\+contramodule structures.
 Here a graded left $B$\+contramodule $Q$ is said to be trivial if
$G^1Q=0$.
 The full subcategory of trivial graded left $B$\+contramodules in
$B\contra$ is equivalent to the category of graded left
$B_0$\+modules $B_0\modl_\sgr$.

 A graded left $B$\+contramodule $P$ is said to be \emph{separated}
if\/ $\bigcap_{m\ge0} G^mP=0$, that is, in other words, the natural
graded $B$\+contramodule map $\lambda_{B,P}\:P\rarrow
\varprojlim_{m\ge0}P/G^mP$ is injective.
 Here the projective limit is taken in the category of graded left
$K$\+modules, which agrees with the projective limit in the category
of graded left $B$\+contramodules.
 A graded left $B$\+contramodule $P$ is said to be \emph{complete}
if the map~$\lambda_{B,P}$ is surjective.

 The functor $\Pi\:B\contra_\sgr\rarrow B\contra$ takes the canonical
decreasing filtration $G$ on a graded left $B$\+contramodule $P$ to
the canonical decreasing filtration $G$ on the ungraded left
$B$\+contramodule $\Pi P$ (which was constructed in
Section~\ref{ungraded-contramodules-subsecn}).
 Consequently, the functor $\Pi\:B\contra_\sgr\rarrow B\contra$ takes
the morphism~$\lambda_{B,P}$ to the morphism $\lambda_{B,\Pi P}$
constructed in Section~\ref{ungraded-contramodules-subsecn}.
 Hence a graded left $B$\+contramodule $P$ is separated (resp.,
complete) if and only if the ungraded left $B$\+contramodule $\Pi P$
is separated (resp., complete).

\begin{prop} \label{separated-graded-contramodules-prop}
 All graded left $B$\+contramodules are complete (but not necessarily
separated).
 All projective graded left $B$\+contramodules are separated.
 Every graded left $B$\+contramodule is the cokernel of an injective
morphism of separated graded left $B$\+contramodules.
\end{prop}

\begin{proof}
 This is a graded version of
Proposition~\ref{separated-ungraded-contramodules}.
 The first assertion follows from the ungraded version in view of
the discussion above.
 The second assertion is clear from the explicit description of
projective graded contramodules as direct summands of the free ones
(similarly to the ungraded version).
 The third assertion is provable similarly to the ungraded version
explained in Proposition~\ref{separated-ungraded-contramodules}.
\end{proof}

\begin{ex} \label{nonseparated-graded-contramodule-counterex}
 Here is an example of a nonseparated graded contramodule.
 It is a modification of a now-classical ungraded example, which
appeared in various guises in~\cite[Example~2.5]{Sim},
\cite[Example~3.20]{Yek1}, and~\cite[Section~A.1.1]{Psemi}
(see~\cite[Section~2.5]{Prev} and~\cite[Example~2.7(1)]{Pcta}
for a further discussion).

 Let $B=k[x]$ be the polynomial algebra in one variable~$x$ of
the degree $\deg x=1$ over a field~$k$.
 Let $E$ be the free graded $B$\+contramodule with a countable
graded set of generators $e_1$, $e_2$, $e_3$,~\dots, all of them
of the degrees $\deg e_i=0$, and let $F$ be the free graded
$B$\+contramodule with a countable set of generators $f_1$, $f_2$,
$f_3$,~\dots\ of degrees $\deg f_i=-i$.
 Let $g\:E\rarrow F$ be the unique graded $B$\+contramodule
morphism taking $e_i$ to $x^if_i$ for every $i\ge1$.
 Consider the quotient contramodule $P=F/g(E)$.
 Then the formal expression $f=\sum_{i=1}^\infty x^if_i$ defines
an element of degree~$0$ in $F$ whose image~$p$ in $P$ belongs to
$\bigcap_{m\ge0}G^mP$.
 Still $f\notin g(E)$, hence $p\ne0$ in~$P$.
\end{ex}

\begin{thm} \label{graded-contramodules-fully-faithful-thm}
 Assume that the augmentation ideal $B_{\ge1}=\bigoplus_{n=1}^\infty
B_n$ of a nonnegatively graded ring $B$ is finitely generated as
a right ideal in~$B$.
 Then the forgetful functor $B\contra_\sgr\rarrow B\modl_\sgr$ from
the category of graded left $B$\+contramodules to the category of
graded left $B$\+modules is fully faithful.
\end{thm}

\begin{proof}
 Let $P$ and $Q$ be graded left $B$\+contramodules, and let
$f\:P\rarrow Q$ be a graded left $B$\+module morphism.
 Consider the ungraded left $B$\+contramodules $\Pi P$ and $\Pi Q$, and
the morphism of ungraded left $B$\+modules $\Pi f\:\Pi P\rarrow\Pi Q$.
 By Theorem~\ref{ungraded-contramodules-fully-faithful}, \,$\Pi f$ is
a morphism of ungraded left $B$\+contramodules.
 It follows easily that $f$~is a morphism of graded left
$B$\+contramodules.
\end{proof}

\subsection{The graded Ext comparison}
\label{graded-ext-comparison-subsecn}
 Let $K$ be a ring.
 Given a graded right $K$\+module $M$ and a graded left $K$\+module $L$,
we will denote by $M\ot_K^\Pi L$ the graded abelian group with
the grading components
$$
 (M\ot_K^\Pi L)_n=\prod\nolimits_{p+q=n} M_p\ot_K L_q.
$$
 Similarly, given two graded right $K$\+modules $L$ and $M$, we denote
by $\Hom_{K^\rop}^\Sigma(L,M)$ the graded abelian group with the grading
components
$$
 \Hom_{K^\rop}^\Sigma(L,M)_n=\bigoplus\nolimits_{p-q=n}
 \Hom_{K^\rop}(L_q,M_p).
$$
 When one of the graded modules involved is concentrated in a finite
set of degrees (so there is no difference between the direct sum and
product), we will drop the superindices $\Pi$ and $\Sigma$ from
the above notation.
 The superindices $\Pi$ and $\Sigma$ will be also absent when the usual
totalization rule is applied (i.~e., the bigraded group of tensor
products is totalized using the direct sums, and the bigraded Hom group
is totalized using the direct products).

\begin{lem} \label{unusual-totalization-Hom-tensor-lemma}
\textup{(a)} Let $K$ and $R$ be associative rings, $V$ be a graded
right $R$\+module, $D$ be a graded $K$\+$R$\+bimodule, and
$H$ be a graded left $K$\+module.
 Then there is a natural map of graded abelian groups
$$
 V\ot_R(D\ot_K^\Pi H)\lrarrow (V\ot_RD)\ot_K^\Pi H.
$$
 This map is an isomorphism whenever $V$ is a finitely generated
projective graded right $R$\+module. \par
\textup{(a)} Let $K$ and $R$ be associative rings, $V$ be a graded
right $R$\+module, $D$ be a graded $R$\+$K$\+bimodule, and
$J$ be a graded right $K$\+module.
 Then there is a natural map of graded abelian groups
$$
 \Hom_{K^\rop}^\Sigma(V\ot_R D,\>J)\lrarrow
 \Hom_{R^\rop}(V,\,\Hom_{K^\rop}^\Sigma(D,J)).
$$
 This map is an isomorphism whenever $V$ is a finitely generated
projective graded right $R$\+module. \qed
\end{lem}

 Let $B=\bigoplus_{n=0}^\infty B_n$ be a nonnegatively graded ring
with the degree-zero component $R=B_0$.
 In the proofs below, we denote by $\Tor^B_*$ and $\Ext_B^*$
the (bi)graded Tor and Ext of graded $B$\+modules (see
Section~\ref{graded-ungraded-subsecn} for a discussion).
 So for any graded left $B$\+module $L$ and graded right $B$\+module
$M$ we have $\Tor^B_i(M,L)=(\Tor^B_i(M,L)_n)_{n\in\boZ}$, and
similarly, for any graded left $B$\+modules $L$ and $M$ we have
$\Ext_B^i(L,M)=(\Ext_B^i(L,M)_n)_{n\in\boZ}$.
 In particular, the degree-zero components
$$
 \Ext^i_B(L,M)_0=\Ext^i_{B\modl_\sgr}(L,M).
$$
are the Ext groups in the abelian category of graded left $B$\+modules
$B\modl_\sgr$.
 More generally, for every $n\in\boZ$,
$$
 \Ext^i_B(L,M)_n=\Ext^i_{B\modl_\sgr}(L,M(-n)),
$$
where $M(-n)$ denotes the left $B$\+module $M$ with the shifted grading,
$M(-n)_j=M_{j+n}$.

 Similarly, as above in this
Section~\ref{comodules-and-contramodules-secn},
given a graded module/abelian group $F$, the notation
$F^+=\Hom_\boZ(F,\boQ/\boZ)$ stands for the graded module/abelian group
with the components $F^+_n=\Hom_\boZ(F_{-n},\boQ/\boZ)$.
 In particular, the natural isomorphism of graded abelian groups
$$
 \Ext^i_B(L,M^+)\simeq\Tor_i^B(M,L)^+
$$
holds for any graded right $B$\+module $M$ and graded left
$B$\+module~$L$.

\begin{thm} \label{graded-contramodule-Ext-comparison}
 Assume that $B$ is a weakly right finitely projective Koszul
graded ring.
 Then the exact, fully faithful forgetful functor of abelian categories
$B\contra_\sgr \rarrow B\modl_\sgr$ induces isomorphisms on all
the Ext groups. \hbadness=1450
\end{thm}

\begin{proof}
 This is a graded version of
Theorem~\ref{ungraded-contramodule-Ext-comparison-thm}.
 The forgetful functor $B\contra_\sgr\rarrow B\modl_\sgr$ is fully
faithful by Theorem~\ref{graded-contramodules-fully-faithful-thm}.
 We have to show that it induces isomorphisms of the groups $\Ext^i$
for $i>0$.

 There are enough projective objects in the abelian category
$B\contra_\sgr$.
 Hence it suffices to prove the graded Ext group vanishing
$\Ext_B^i(P,Q)=0$ for all projective objects $P\in B\contra_\sgr$,
all objects $Q\in B\contra_\sgr$, and all $i>0$.

 By Proposition~\ref{separated-graded-contramodules-prop}, any
graded left $B$\+contramodule is the cokernel of an injective
morphism of separated graded left $B$\+contramodules.
 Furthermore, any separated graded left $B$\+contramodule $Q$ is
an infinitely iterated extension, in the sense of the projective limit,
of trivial graded left $B$\+contramodules $G^mQ/G^{m+1}Q$.
 Using a graded version of the dual Eklof
lemma~\cite[Proposition~18]{ET} (see~\cite[Lemma~4.5]{PR} for
a far-reaching generalization), we can reduce the question to
the case when $Q$ is a trivial graded left $B$\+contramodule.
 In other words, $Q$ is just a graded left $R$\+module in which
the augmentation ideal $B_{\ge1}\subset B$ acts by zero.

 The graded left $R$\+module $Q$ has a resolution by injective
graded left $R$\+modules of the form $F^+=\Hom_\boZ(F,\boQ/\boZ)$,
where $F$ ranges over free graded right $R$\+modules.
 Arguing further as in the proof of
Theorem~\ref{ungraded-contramodule-Ext-comparison-thm},
we see that it suffices to prove the graded Tor group vanishing
$\Tor^B_i(R,P)=0$ for all projective objects $P\in B\contra_\sgr$
and all integers $i>0$ (where $R$ is viewed as a graded right
$B$\+module with $B_{\ge1}$ acting by zero).

 A description of projective graded left $B$\+contramodules was given
in the proof of Proposition~\ref{graded-B-contramodules-proposition}.
 In the notation above, we can rephrase it as follows.
 The projective graded left $B$\+contramodules are precisely
the direct summands of graded left $B$\+contramodules of the form
$\boM_K^\gr(H)=B\ot_K^\Pi K[X]$, where $K[X]$ is the free graded
left $K$\+module spanned by a graded set~$X$.
 We put $H=K[X]$ for brevity.
 
 Now we compute the graded groups $\Tor^B_i(R,\>B\ot_K^\Pi H)$ using
the graded projective resolution~\eqref{weakly-Koszul-resolution}
of the graded right $B$\+module~$R$.
 Taking the tensor product of~\eqref{weakly-Koszul-resolution}
with $B\ot_K^\Pi H$ over~$B$, we obtain the complex of graded
abelian groups
\begin{equation} \label{Pi-tensor-complex}
 \dotsb\lrarrow V_2\ot_R(B\ot_K^\Pi H)\lrarrow
 V_1\ot_R(B\ot_K^\Pi H)\lrarrow B\ot_K^\Pi H\lrarrow0.
\end{equation}
 According to Lemma~\ref{unusual-totalization-Hom-tensor-lemma}(a),
the complex~\eqref{Pi-tensor-complex} is isomorphic to
\begin{equation} \label{Pi-tensor-complex-rewritten}
 \dotsb\lrarrow (V_2\ot_RB)\ot_K^\Pi H\lrarrow
 (V_1\ot_RB)\ot_K^\Pi H\lrarrow B\ot_K^\Pi H\lrarrow0.
\end{equation}
 The complex~\eqref{Pi-tensor-complex-rewritten} can be obtained by
applying the exact functor ${-}\ot_K^\Pi H$ to
the resolution~\eqref{weakly-Koszul-resolution}.
 Thus $\Tor^B_i(R,\>B\ot_K^\Pi H)=0$ for $i>0$, as desired.
\end{proof}

\begin{thm} \label{graded-comodule-Ext-comparison}
 Assume that $B$ is a weakly right finitely projective Koszul
graded ring.
 Then the exact, fully faithful inclusion functor of abelian
categories\/ $\comodrgr B \rarrow\modrgr B$ induces isomorphisms
on all the Ext groups. \hbadness=1450
\end{thm}

\begin{proof}
 This is a graded version of
Theorem~\ref{ungraded-comodule-Ext-comparison-thm}.
 There are enough injective objects in a Grothendieck abelian category
$\comodrgr B$.
 Hence it suffices to prove the graded Ext group vanishing
$\Ext_{B^\rop}^i(M,L)=0$ for all objects $M\in\comodrgr B$, all
injective objects $L\in\comodrgr B$, and all $i>0$.

 Every graded right $B$\+comodule $M$ has a natural increasing
filtration $0=G_{-1}M\subset G_0M\subset G_1M\subset G_2M\subset\dotsb$
by graded right $B$\+submodules $G_mM\subset M$ with the grading
component $G_mM_j\subset M_j$ consisting of all the elements
$x\in M_j$ such that $xB_{\ge m+1}=0$ in~$M$.
 By the definition, we have $M=\bigcup_{m\ge0}G_mM$.
 The successive quotient modules $G_mM/G_{m-1}M$ are \emph{trivial}
graded right $B$\+comodules, that is graded right $B$\+modules in
which $B_{\ge1}$ acts by zero.
 The category of trivial graded right $B$\+comodules is equivalent to
the category of graded right $R$\+modules $\modrgr R$.

 Using a graded version of the Eklof lemma~\cite[Lemma~1]{ET}, we can
assume that $M$ is a trivial graded right $B$\+comodule.
 Then the graded right $R$\+module $M$ has a resolutions by free graded
right $R$\+modules.
 Similarly to the proof of
Theorem~\ref{ungraded-comodule-Ext-comparison-thm}, we see that it
suffices to prove the graded Ext group vanishing
$\Ext^i_{B^\rop}(R,L)=0$ for all injective objects $L\in\comodrgr B$
and all integers $i>0$ (where $R$ is viewed as a graded right
$R$\+module with $B_{\ge1}$ acting by zero).

 Injective objects of the category $\comodrgr B$ can be described as
follows.
 The inclusion functor $\comodrgr B\rarrow\modrgr B$ has a right
adjoint functor $\Gamma_B\:\modrgr B\rarrow\comodrgr B$ assigning
to a graded right $B$\+module $N$ its graded submodule $\Gamma_B(N)
\subset N$ whose component $\Gamma_B(N)_j$ consists of all elements
$x\in N_j$ for which there exists $m\ge1$ such that $xB_{\ge m}=0$.
 The injective objects of $\comodrgr B$ are precisely the direct
summands of the objects $\Gamma_B(N)$, where $N$ ranges over
the injective objects of $\modrgr B$.

 The following description of injective graded right $B$\+modules is
convenient for our purposes.
 Let $H=K[X]$ be the free graded left $K$\+module spanned by a graded
set $X$, and let $J=K[X]^+$ be the graded character module of~$H$.
 Then $J$ is an injective graded right $K$\+module, and all injective
graded right $K$\+modules are direct summands of modules of this form.
 The graded Hom module $\Hom_{K^\rop}(B,J)$ is an injective graded
right $B$\+module, and all injective graded right $B$\+modules are
direct summands of modules of this form.
 One can easily check that $\Gamma_B(\Hom_{K^\rop}(B,J))=
\Hom_{K^\rop}^\Sigma(B,J)$ (see Lemma~\ref{gamma-hom-hom-sigma} below).

 Now we compute the graded groups
$\Ext_{B^\rop}^i(R,\Hom_{K^\rop}^\Sigma(B,J))$ using the graded
projective resolution~\eqref{weakly-Koszul-resolution} of
the graded right $B$\+module~$R$.
 Taking the graded right $B$\+module Hom
from~\eqref{weakly-Koszul-resolution} to $\Hom_{K^\rop}^\Sigma(B,J)$,
we obtain the complex of graded abelian groups
\begin{multline} \label{Sigma-Hom-complex}
 0\lrarrow\Hom_{K^\rop}^\Sigma(B,J)\lrarrow
 \Hom_{R^\rop}(V_1,\,\Hom_{K^\rop}^\Sigma(B,J)) \\
 \lrarrow\Hom_{R^\rop}(V_2,\,\Hom_{K^\rop}^\Sigma(B,J))\lrarrow\dotsb
\end{multline}
 According to Lemma~\ref{unusual-totalization-Hom-tensor-lemma}(b),
the complex~\eqref{Sigma-Hom-complex} is isomorphic to
\begin{multline} \label{Sigma-Hom-complex-rewritten}
 0\lrarrow\Hom_{K^\rop}^\Sigma(B,J)\lrarrow
 \Hom_{K^\rop}^\Sigma(V_1\ot_R B,\>J) \\ \lrarrow
 \Hom_{K^\rop}^\Sigma(V_2\ot_R B,\>J)\lrarrow\dotsb
\end{multline}
 The complex~\eqref{Sigma-Hom-complex-rewritten} can be obtained by
applying the exact functor $\Hom_{K^\rop}^\Sigma({-},J)$ to
the resolution~\eqref{weakly-Koszul-resolution}.
 Thus $\Ext_{B^\rop}^i(R,\Hom_{K^\rop}^\Sigma(B,J))=0$ for
$i>0$, as desired.
\end{proof}

\Section{Relative Nonhomogeneous Derived Koszul Duality: \\
the Comodule Side}
\label{comodule-side-secn}

\subsection{CDG-modules} \label{cdg-modules-subsecn}
 Let $(B,d,h)$ be a curved DG\+ring, as defined in
Section~\ref{curved-dg-rings-subsecn}, and let $(\hB,\d)$ be
the related quasi-differential graded ring constructed in
Theorem~\ref{cdg-qdg-equivalence}.
 A \emph{left CDG\+module over $(B,d,h)$} can be simply defined as
a graded left $\hB$\+module.

 Equivalently, a left CDG\+module $(M,d_M)$ over $(B,d,h)$ is a graded
left $B$\+module $M=\bigoplus_{n\in\boZ} M^n$ endowed with a sequence
of additive maps $d_{M,n}\:M^n\rarrow M^{n+1}$, \ $n\in\boZ$, satisfying
following equations:
\begin{enumerate}
\renewcommand{\theenumi}{\roman{enumi}}
\item $d_M$~is an odd derivation of the graded left module $M$
compatible with the odd derivation~$d$ on the graded ring $B$,
that is $d_M(bx)=d(b)x+(-1)^{|b|}bd_M(x)$ for all $b\in B^{|b|}$ and
$x\in M^{|x|}$, \,$|b|$, $|x|\in\boZ$;
\item $d_M^2(x)=hx$ for all $x\in M$.
\end{enumerate}

 A \emph{right CDG\+module over $(B,d,h)$} can be simply defined as
a graded right $\hB$\+module.
 Equivalently, a right CDG\+module $(N,d_N)$ over $(B,d,h)$ is
a graded right $B$\+module $N=\bigoplus_{n\in\boZ}N^n$ endowed with
a sequence of additive maps $d_{N,n}\:N^n\rarrow N^{n+1}$, \ $n\in\boZ$,
satisfying the following equations:
\begin{enumerate}
\renewcommand{\theenumi}{\roman{enumi}}
\item $d_N$~is an odd derivation of the graded right module $N$
compatible with the odd derivation~$d$ on the graded ring $B$,
that is $d_N(yb)=d_N(y)b+(-1)^{|y|}yd(b)$ for all $b\in B^{|b|}$ and
$y\in N^{|y|}$, \,$|b|$, $|y|\in\boZ$;
\item $d_N^2(y)=-yh$ for all $y\in N$.
\end{enumerate}

 Notice that the graded $B$\+$B$\+bimodule $B$ with its odd
derivation~$d$ is \emph{neither} a left \emph{nor} a right CDG\+module
over $(B,d,h)$, because the formula~(ii) for $d^2(b)$, \ $b\in B$ in
Section~\ref{curved-dg-rings-subsecn} does not agree with
the formulas~(ii) in the two definitions above.
 However, the diagonal graded $B$\+$B$\+bimodule $B$ is naturally
a \emph{CDG\+bimodule} over the CDG\+ring $(B,d,h)$, in the sense
of the following definition.

 Let $B=(B,d_B,h_B)$ and $C=(C,d_C,h_C)$ be two CDG\+rings.
 A \emph{CDG\+bimodule $K$ over $B$ and~$C$} is a graded
$B$\+$C$\+bimodule $K=\bigoplus_{n\in\boZ}K^n$ endowed with a sequence
of additive maps $d_{K,n}\:K^n\rarrow K^{n+1}$, \ $n\in\boZ$,
satisfying the following conditions: {\hbadness=1100
\begin{enumerate}
\renewcommand{\theenumi}{\roman{enumi}}
\item $d_K$ is \emph{both} an odd derivation of the graded left
$B$\+module $K$ compatible with the odd derivation~$d_B$ on
the graded ring~$B$ \emph{and} an odd derivation of the graded right
$C$\+module $K$ compatible with the odd derivation~$d_C$ on
the graded ring~$C$;
\item $d_K^2(z)=h_Bz-zh_C$ for all $z\in K$.
\end{enumerate}

 Notice} that a CDG\+bimodule over $B$ and $C$ is \emph{not} the same
thing as a graded $\hB$\+$\hC$\+bimodule.
 We leave it to the reader to construct a CDG\+ring structure on
the graded ring $B\ot_\boZ C^\rop$ for which a CDG\+bimodule over $B$
and $C$ would be the same thing as a left CDG\+module over
$B\ot_\boZ C^\rop$.

 Given a graded left $B$\+module $M$ and an integer $n\in\boZ$,
we denote by $M[n]$ the graded left $B$\+module with the grading
components $M[n]^i=M^{n+i}$, \ $i\in\boZ$, and the left action of $B$
transformed with the sign rule $b(x[n])=(-1)^{n|b|}(bx)[n]$ for
any element $x\in M^{n+i}$ and the corresponding element
$x[n]\in M[n]^i$.
 Given a graded right $B$\+module $N$ and an integer $n\in\boZ$,
we let $N[n]$ denote the graded right $B$\+module $N$ with
the components $N[n]^i=N^{n+i}$ and the right action of $B$
unchanged, $(y[n])b=(yb)[n]$ for $y\in N^{n+i}$ and $y[n]\in N[n]^i$.
 If a differential $d_M$ or $d_N$ is given, a grading shift transforms
it with the usual sign rule $d_{M[n]}(x[n])=(-1)^n(d_M(x))[n]$
or $d_{N[n]}(y[n])=(-1)^n(d_N(y))[n]$.
 These sign rules are best understood with the hint that we write
the shift-of-degree symbol $[n]$ to the right of our modules and their
elements, as it is traditionally done, but tacitly presume it to be
a \emph{left} operator (as all our operators, such as the differentials,
homomorphisms, etc., are generally left operators).
 So one may want to write $[n]x\in [n]M$ instead of $x[n]\in M[n]$
for $x\in M$ and similarly $[n]y\in [n]N$ for $y\in N$, which would
explain the signs.

 Given two graded left $B$\+modules $L$ and $M$, we denote by
$\Hom_B^n(L,M)$ the abelian group of all graded left $B$\+module
homomorphisms $L\rarrow M[n]$.
 For any two left CDG\+modules $L$ and $M$ over a CDG\+ring $(B,d,h)$,
there is a natural differential~$d$ on the graded abelian group
$\Hom_B(L,M)$ defined by the usual rule $d(f)(l)=d_M(f(l))-
(-1)^{|f|}f(d_L(l))$ for all $f\in\Hom_B^{|f|}(L,M)$, \ $|f|\in\boZ$,
and $l\in L$.
 Similarly, given two graded right $B$\+modules $K$ and $N$, we denote
by $\Hom_{B^\rop}^n(K,N)$ the abelian group of all graded right
$B$\+module homomorphisms $K\rarrow N[n]$.
 For any two right CDG\+modules $K$ and $N$ over a CDG\+ring $(B,d,h)$,
there is a natural differential~$d$ on the graded abelian group
$\Hom_{B^\rop}(K,N)$ defined by the formula $d(g)(k)=d_N(g(k))-
(-1)^{|g|}g(d_K(k))$ for all $g\in\Hom_{B^\rop}^{|g|}(K,N)$, \
$|g|\in\boZ$, and $k\in K$.

 One easily checks that $d^2(f)=0=d^2(g)$ for all $f$ and~$g$ in
the notation above.
 Hence for any two left CDG\+modules $L$ and $M$ over $B$,
the \emph{complex of morphisms} $\Hom_B(L,M)$ is defined; and similarly,
for any two right CDG\+modules $K$ and $N$ over $B$, there is
the complex of morphisms $\Hom_{B^\rop}(K,N)$.
 These constructions produce the \emph{DG\+category of left
CDG\+modules} $\DG(B\modl)$ and the similar \emph{DG\+category of
right CDG\+modules} $\DG(\modr B)$ \emph{over~$B$}.
 We refer to~\cite[Sections~1.2 and~3.1]{Pkoszul} for a discussion
of these DG\+categories and their triangulated \emph{homotopy
categories} $\Hot(B\modl)=H^0\DG(B\modl)$ and
$\Hot(\modr B)=H^0\DG(\modr B)$.

 Given a right graded $B$\+module $N$ and a left graded $B$\+module
$M$, the graded tensor product $N\ot_BM$ is a graded abelian group
constructed in the usual way.
 So we consider the bigraded abelian group $N\ot_\boZ M$ with
the components $N^i\ot_\boZ M^j$, totalize it by taking infinite
direct sums along the diagonals $i+j=n$, and consider the graded
quotient group by the graded subgroup spanned by the elements
$yb\ot x-y\ot bx$, where $y\in N$, \ $b\in B$, and $x\in M$.
 For any right CDG\+module $N$ and any left CDG\+module $M$ over
$(B,d,h)$, the tensor product $N\ot_BM$ is naturally a complex of
abelian groups with the differential $d(y\ot x)=d_N(y)\ot x+
(-1)^{|y|}y\ot d_M(x)$.
 One can easily check that $d^2(y\ot x)=0$ in $N\ot_BM$.

 More generally, let $(A,d_A,h_A)$, \ $(B,d_B,h_B)$, and $(C,d_C,h_C)$
be three CDG\+rings.
 Given a graded $C$\+$B$\+bimodule $N$ and a graded $B$\+$A$\+bimodule
$K$, the graded $C$\+$A$\+bimodule structure on the tensor product
$N\ot_BK$ is defined by the usual rule $c(y\ot z)a=(cy)\ot(za)$
for all $c\in C$, \ $y\in N$, \ $z\in K$, and $a\in A$.
 Now let $N$ be a CDG\+bimodule over $C$ and $B$, and let $K$ be
a CDG\+bimodule over $B$ and~$A$.
 Then the graded $C$\+$A$\+bimodule $N\ot_BK$ with the differential
$d(y\ot z)=d_N(y)\ot z+(-1)^{|y|}y\ot d_K(z)$, where $y\in N^{|y|}$
and $z\in K^{|z|}$, is a CDG\+bimodule over $C$ and~$A$.

 Given a graded $B$\+$A$\+bimodule $K$ and a graded $B$\+$C$\+bimodule
$L$, the graded $A$\+$C$\+bimodule structure on the graded abelian
group $\Hom_B(K,L)$ is defined by the formula $(afc)(z)=
(-1)^{|a||f|+|a||z|+|c||z|}f(za)c$ for all $a\in A^{|a|}$, \
$c\in C^{|c|}$, \ $z\in K^{|z|}$, and $f\in\Hom_B^{|f|}(K,L)$.
 Now let $K$ be a CDG\+bimodule over $B$ and $A$, and let $L$ be
a CDG\+bimodule over $B$ and~$C$.
 Then the graded $A$\+$C$\+bimodule $\Hom_B(K,L)$ with the differential
$d(f)(z)=d_L(f(z))-(-1)^{|f|}f(d_K(z))$ is a CDG\+bimodule
over $A$ and~$C$.

 Given a graded $B$\+$A$\+bimodule $K$ and a graded $C$\+$A$\+bimodule
$M$, the graded $C$\+$B$\+bimodule structure on the graded abelian
group $\Hom_{A^\rop}(K,M)$ is defined by the formula
$(cgb)(z)=cg(bz)$ for all $c\in C$, \ $b\in B$, \ $z\in K$, and
$g\in\Hom_{A^\rop}(K,M)$.
 Now let $K$ be a CDG\+bimodule over $B$ and $A$, and let $M$ be
a CDG\+bimodule over $C$ and~$A$.
 Then the graded $C$\+$B$\+bimodule $\Hom_{A^\rop}(K,M)$ with
the differential $d(g)(z)=d_M(g(z))-(-1)^{|g|}g(d_K(z))$, where
$z\in K^{|z|}$ and $g\in\Hom_{A^\rop}^{|g|}(K,M)$, is
a CDG\+bimodule over $C$ and~$B$.

 Let $A$, $B$, $C$, and $D$ be four CDG\+rings.
 Then for any CDG\+bimodule $L$ over $A$ and $B$, any CDG\+bimodule $M$
over $B$ and $C$, and any CDG\+bimodule $N$ over $C$ and $D$,
the natural isomorphism of graded $A$\+$D$\+bimodules
$$
 (L\ot_BM)\ot_CN\simeq L\ot_B(M\ot_CN)
$$
is a closed isomorphism of CDG\+bimodules over $A$ and~$D$.

 Similarly, for any CDG\+bimodule $L$ over $B$ and $A$, any
CDG\+bimodule $M$ over $C$ and $B$, and any CDG\+bimodule $N$ over $C$
and $D$, the natural isomorphism of graded $A$\+$D$\+bimodules
$$
 f'\in\Hom_C(M\ot_BL,\>N)\,\simeq\,\Hom_B(L,\Hom_C(M,N))\ni f''
$$
given by the rule $f''(z)(x)=(-1)^{|x||z|}f'(x\ot z)$ for all
$z\in L^{|z|}$ and $x\in M^{|x|}$ is a closed isomorphism of
CDG\+bimodules over $A$ and~$D$.

 For any CDG\+bimodule $L$ over $A$ and $B$, any CDG\+bimodule $M$
over $B$ and $C$, and any CDG\+bimodule $N$ over $D$ and $C$,
the natural isomorphism of graded $D$\+$A$\+bimodules
$$
 g'\in\Hom_{C^\rop}(L\ot_BM,\>N)\,\simeq\,
 \Hom_{B^\rop}(L,\Hom_{C^\rop}(M,N))\ni g''
$$
given by the rule $g''(z)(x)=g'(z\ot x)$ for all $z\in L$ and $x\in M$
is a closed isomorphism of CDG\+bimodules over $D$ and~$A$.

 In particular, let $K$ be a CDG\+bimodule over $B$ and $A$.
 Then the DG\+functor
$$
 K\ot_A{-}\:\DG(A\modl)\lrarrow\DG(B\modl)
$$
is left adjoint to the DG\+functor
$$
 \Hom_B(K,{-})\:\DG(B\modl)\lrarrow\DG(A\modl).
$$
 Here the DG\+functor $K\ot_A{-}$ acts on morphisms by the rule
$(K\ot_Af)(k\ot l)=(-1)^{|f||k|}k\ot f(l)$ for all
$L$, $M\in\DG(A\modl)$, \ $f\in\Hom_A^{|f|}(L,M)$, \ $k\in K^{|k|}$,
and $l\in L^{|l|}$.
 The DG\+functor $\Hom_B(K,{-})$ acts on morphisms by the rule
$\Hom_B(K,f)(g)(k)=f(g(k))$ for all $L$, $M\in\DG(B\modl)$, \
$f\in\Hom_B(L,M)$, \ $g\in\Hom_B(K,L)$, and $k\in K$.

 Hence the induced triangulated functors $K\ot_A{-}$ and $\Hom_B(K,{-})$
between the homotopy categories $\Hot(A\modl)$ and $\Hot(B\modl)$
are also adjoint.

 Similarly, the DG\+functor
$$
 {-}\ot_BK\:\DG(\modr B)\lrarrow\DG(\modr A)
$$
is left adjoint to the DG\+functor
$$
 \Hom_{A^\rop}(K,{-})\:\DG(\modr A)\lrarrow\DG(\modr B).
$$
 Here the DG\+functor ${-}\ot_BK$ acts on morphisms by the rule
$(f\ot_BK)(l\ot k)=f(l)\ot k$ for all $L$, $M\in\DG(\modr B)$, \
$f\in\Hom_{B^\rop}(L,M)$, \ $l\in L$, and $k\in K$.
 The DG\+functor $\Hom_{A^\rop}(K,{-})$ acts on morphisms by
the fule $\Hom_{A^\rop}(K,f)(g)(k)=f(g(k))$ for all
$L$, $M\in\DG(\modr A)$, \ $f\in\Hom_{A^\rop}(L,M)$, \
$g\in\Hom_{A^\rop}(K,L)$, and $k\in K$.

 Hence the induced triangulated functors ${-}\ot_BK$ and
$\Hom_{A^\rop}(K,{-})$ between the homotopy categories $\Hot(\modr B)$
and $\Hot(\modr A)$ are also adjoint.

 Now let $(f,a)\:\secB=(\secB,d'',h'')\rarrow\prB=(\prB,d',h')$ be
a morphism of CDG\+rings (as defined in
Section~\ref{curved-dg-rings-subsecn}).
 Let $(M,d'_M)$ be a left CDG\+module over $(\prB,d',h')$.
 The restriction of scalars via the graded ring homomorphism
$f\:\secB\rarrow\prB$ defines a structure of graded left
$\secB$\+module on~$M$.
 Put $d''_M(x)=d'_M(x)+ax$ for all $x\in M$.
 Then $(M,d''_M)$ is a left CDG\+module over $(\secB,d'',h'')$.
 This construction defines the DG\+functor $\DG(\prB\modl)\rarrow
\DG(\secB\modl)$ of restriction of scalars via~$(f,a)$.

 Let $(N,d'_N)$ be a right CDG\+module over $(\prB,d',h')$.
 The restriction of scalars via the graded ring homomorphism
$f\:\secB\rarrow\prB$ defines a structure of graded right
$\secB$\+module on~$N$.
 Put $d''_N(y)=d'_N(y)-(-1)^{|y|}ya$ for all $y\in N^{|y|}$,
\ $|y|\in\boZ$.
 Then $(N,d''_N)$ is a right CDG\+module over $(\secB,d'',h'')$.
 This construction defines the DG\+functor $\DG(\prB\modl)\rarrow
\DG(\secB\modl)$ of restriction of scalars via~$(f,a)$.

 Similarly, let $(f,a_B)\:\secB=(\secB,d''_B,h''_B)\rarrow
\prB=(\prB,d'_B,h'_B)$ and $(g,a_C)\:\secC=(\secC,d''_C,h''_C)
\rarrow\prC=(\prC,d''_C,h''_C)$ be two morphisms of CDG\+rings.
 Let $(K,d'_K)$ be a CDG\+bimodule over $(\prB,d'_B,h'_B)$ and
$(\prC,d'_C,h'_C)$.
 The restriction of scalars via the graded ring homomorphisms
$f\:\secB\rarrow\prB$ and $g\:\secC\rarrow\prC$ defines a structure
of graded $\secB$\+$\secC$\+bimodule on~$K$.
 Put $d''_K(z)=d'_K(z)+a_Bz-(-1)^{|z|}za_C$ for all $z\in K^{|z|}$,
\ $|z|\in\boZ$.
 Then $(K,d''_K)$ is a CDG\+bimodule over $(\secB,d''_B,h''_B)$ and
$(\secC,d''_C,h''_C)$.
 All the constructions of the tensor product complex and the tensor
product and Hom CDG\+bimodules above in this section agree with
the transformations of CDG\+(bi)modules induced by isomorphisms of
CDG\+rings.

\subsection{Nonhomogeneous Koszul complex/CDG-module}
\label{nonhomogeneous-koszul-cdg-module-subsecn}
 Let $R\subset\tV\subset\tA$ be a $3$\+left finitely projective weak
nonhomogeneous quadratic ring, as defined in
Sections~\ref{nonhomogeneous-quadratic-subsecn}
and~\ref{self-consistency-subsecn}, and let $A=\gr^F\tA$ be its
associated graded ring with respect to the filtration~$F$.
 As in Sections~\ref{koszulity-secn}\+-%
\ref{nonhomogeneous-quadratic-secn}, we put $V=A_1$ and denote
by $B$ the $3$\+right finitely projective quadratic graded ring
quadratic dual to~$\q A$.
 Choose a left $R$\+submodule of strict generators $V'\subset\tV=
F_1\tA$, and let $(B,d,h)=(B,d',h')$ denote the related CDG\+ring
structure on the graded ring $B$, as constructed in
Proposition~\ref{nonhomogeneous-dual-cdg-ring}.

 Denote by $e'\in\Hom_R(V,R)\ot_R\tV=B^1\ot_RF_1\tA$ the element
corresponding to the injective left $R$\+module map
$V\simeq V'\hookrightarrow \tV$ under the construction of
Lemma~\ref{element-e-lemma}.
 Our aim is to construct the \emph{dual nonhomogeneous Koszul
CDG\+module} $\Ksp(B,\tA)=\Ksp_{e'}(B,\tA)$, which has the form
\begin{equation} \label{dual-nonhomogeneous-koszul-cdg-module}
 0\lrarrow\tA\lrarrow B^1\ot_R\tA\lrarrow B^2\ot_R\tA\lrarrow
 B^3\ot_R\tA\lrarrow\dotsb.
\end{equation}
 The differential on $\Ksp_{e'}(B,\tA)$ does not square to zero when
$h\ne0$.
 Rather, $\Ksp_{e'}(B,\tA)$ is a left CDG\+module over the CDG\+ring
$(B,d,h)$, with the obvious underlying graded left $B$\+module
structure of the tensor product $B\ot_R\tA$.
 However, the differential on $\Ksp_{e'}(B,\tA)$ preserves the right
$\tA$\+module structure of the tensor product $B\ot_R\tA$.
 Summarizing, one can say that $\Ksp_{e'}(B,\tA)$ is going to be
a CDG\+bimodule over the CDG\+rings $B=(B,d,h)$ and $\tA=(\tA,0,0)$
(where $\tA$ is viewed as a graded ring concentrated entirely
in degree zero).

 When $\tA=A=\q A$ is a $3$\+left finitely projective quadratic
graded ring and $V'=A_1\hookrightarrow F_1\tA=R\oplus A_1$ is
the obvious splitting of the direct summand projection $R\oplus A_1
\rarrow A_1$ (so the related CDG\+ring structure on $B$ is trivial,
$d=0=h$), the construction of differential~$d_{e'}$ in
the next proposition reduces to that of the differential~$d_e$
on the dual Koszul complex $\Ksp_e{}^\bu=B\ot_RA$ from
Section~\ref{dual-koszul-complex-subsecn}.

\begin{prop}
 The formula
\begin{equation} \label{dual-nonhomogeneous-koszul-differential}
 d_{e'}(b\ot c)=d(b)\ot c+(-1)^{|b|}be'c
\end{equation}
for all $b\in B^{|b|}$ and $c\in\tA$, where $e'\in B^1\ot_R F_1\tA$
is the element mentioned above, endows the graded left $B$\+module
$\Ksp(B,\tA)=B\ot_R\tA$ with a well-defined structure of
left CDG\+module over the CDG\+ring $B=(B,d,h)$.
\end{prop}

\begin{proof}
 Neither one of the two summands in
the formula~\eqref{dual-nonhomogeneous-koszul-differential}
produces a well-defined map $B\ot_R\tA\rarrow B\ot_R\tA$.
 It is only the sum that is well-defined as an endomorphism of
the tensor product over~$R$.
 However, the separate summands are well-defined as maps
$B\ot_\boZ\tA\rarrow B\ot_R\tA$.
 Indeed, the first summand is obviously well-defined as a map
$B\ot_\boZ\tA\rarrow B\ot_\boZ\tA$, and one can consider
the composition with the natural surjection $B\ot_\boZ\tA
\rarrow B\ot_R\tA$.
 The grading/filtration components of the second summand (without
the $\pm$~sign) are the compositions
$$
 B^n\ot_\boZ F_j\tA\overset{e'}\lrarrow B^n\ot_\boZ B^1\ot_R F_1\tA
 \ot_\boZ F_j\tA\lrarrow B^{n+1}\ot_R F_{j+1}\tA.
$$
 Here the map $B^n\ot_\boZ F_j\tA\rarrow B^n\ot_\boZ B^1\ot_R F_1\tA
\ot_\boZ F_j\tA$ is induced by the abelian group element or map
$e'\:\boZ\rarrow B^1\ot_RF_1\tA$, while the map $B^n\ot_\boZ B^1
\ot_R F_1\tA\ot_\boZ F_j\tA\lrarrow B^{n+1}\ot_R F_{j+1}\tA$ is
the tensor product of two multiplication maps $B^n\ot_\boZ B^1\rarrow
B^{n+1}$ and $F_1\tA\ot_\boZ F_j\tA\rarrow F_{j+1}\tA$.

 Let us check that the map~$d_{e'}$ is well-defined on
$B\ot_R\tA$, that is, for any $b\in B$, \ $r\in R$, and $c\in\tA$
one has
$$
 d(br)\ot c+(-1)^{|b|}bre'c=d(b)\ot rc+(-1)^{|b|}be'rc
$$
in $B\ot_R\tA$.
 In the pairing notation from the proof of Lemma~\ref{element-e-lemma},
by the definition of the element $e'\in\Hom_R(V,R)\ot_R\tV$ we have
$\lan v,e'\ran = v'$ for every $v\in V$, where $v'\in V'\subset\tV$
denotes the element corresponding to~$v$ under the left $R$\+module
identification $V\simeq V'$.
 Hence, in particular,
$$
 \lan v,re'\ran=\lan vr,e'\ran=(vr)'
$$
while
$$
 \lan v,e'r\ran=\lan v,e'\ran*r=v'*r=(vr)'+q(v,r)=
 (vr)'+\lan v,d(r)\ran.
$$
 Here, as in Section~\ref{self-consistency-subsecn}, the symbol~$*$
denotes the multiplication in $\tA$, and in particular, the right
action of $R$ in~$\tV$.
 The notation $q(v,r)$ was defined in the equation~\eqref{map-q}, and
the equation $q(v,r)=\lan v,d(r)\rangle$ comes from the definition
of the differential~$d$ on the graded ring~$B$ in
Proposition~\ref{nonhomogeneous-dual-cdg-ring}.
 Therefore,
$$
 \lan v, e'r \ran = \lan v, re' \ran +
 \lan v, d(r) \ran
 \qquad\text{for all $v\in V$}.
$$
 We have deduced the equation
\begin{equation} \label{e-prime-left-right-R-actions}
 e'r=re'+d(r)\ot 1
\end{equation}
comparing the left and right actions of the ring $R$ on the element
$e'\in\Hom_R(V,R)\ot_R\tV$ (where $1\in R\subset\tV$ denotes the unit
element).
 Now we can compute that
\begin{multline*}
 d(br)\ot c+(-1)^{|b|}bre'c=d(b)r\ot c+(-1)^{|b|}bd(r)\ot c+
 (-1)^{|b|}bre'c \\ =d(b)\ot rc + (-1)^{|b|}be'rc \,\in\, B\ot_R\tA,
\end{multline*}
as desired.
 So the map $d_{e'}\:B\ot_R\tA\rarrow B\ot_R\tA$ is well-defined.

 It is obvious that $d_{e'}$~is a right $\tA$\+module map.
 It is also easy to check that $d_{e'}$~is an odd derivation of
the graded left $B$\+module $B\ot_R\tA$ compatible with
the odd derivation~$d$ on the graded ring~$B$.
 Indeed, for all $b$, $b'\in B$ and $c\in\tA$ we have
\begin{multline*}
 d_{e'}(b'b\ot c) = d(b'b)\ot c + (-1)^{|b'|+|b|}b'be'c \\ =
 d(b')b\ot c+(-1)^{|b'|}b'd(b)\ot c +(-1)^{|b'|+|b|}b'be'c \\ =
 d(b')b\ot c+(-1)^{|b'|}b'\left(d(b)\ot c+(-1)^{|b|}be'c\right)=
 d(b')b\ot c+(-1)^{|b'|}b'd_{e'}(b\ot c),
\end{multline*}
as desired.

 It remains to check the equation
$$
 d_{e'}^2(b\ot c)=hb\ot c
 \qquad\text{for all $b\in B$ and $c\in\tA$}.
$$
 For this purpose, we choose two finite collections of elements
$v_\alpha\in V$ and $u^\alpha\in B^1=\Hom_R(V,R)$, indexed by the same
indices~$\alpha$, such that $e=\sum_\alpha u^\alpha\ot v_\alpha
\in \Hom_R(V,R)\ot_R V$.
 Then the assertion of Lemma~\ref{element-e-lemma}(b) can be written
as the equation
\begin{equation} \label{left-right-actions-on-e-equation}
 \sum\nolimits_\alpha su^\alpha\ot v_\alpha=
 \sum\nolimits_\alpha u^\alpha\ot v_\alpha s
 \qquad\text{for all $s\in R$}.
\end{equation}
 For any element $v\in V$, we have
$\sum_\alpha\lan v,u^\alpha\ran v_\alpha=v$.

 By the definition,
$$
 d_{e'}(b\ot c)=d(b)\ot c+
 (-1)^{|b|}\sum\nolimits_\alpha bu^\alpha\ot v_\alpha'*c,
$$
hence
\begin{multline*}
 d_{e'}^2(b\ot c)= d^2(b)\ot c +
 (-1)^{|b|}\sum\nolimits_\alpha d(bu^\alpha)\ot v_\alpha'*c \\
 +(-1)^{|b|+1}\sum\nolimits_\alpha d(b)u^\alpha\ot v_\alpha'*c 
 -\sum\nolimits_\alpha\sum\nolimits_\beta
 bu^\alpha u^\beta\ot v_\beta'*v_\alpha'*c \\
 =(hb-bh)\ot c + \sum\nolimits_\alpha bd(u^\alpha)\ot v_\alpha'*c
 -\sum\nolimits_\alpha\sum\nolimits_\beta
 bu^\alpha u^\beta\ot v_\beta'*v_\alpha'*c.
\end{multline*}
 So we only need to check that the equation
\begin{equation} \label{nonhomogeneous-duality-equation-with-indices}
 -h\ot 1+\sum\nolimits_\alpha d(u^\alpha)\ot v_\alpha'
 -\sum\nolimits_\alpha\sum\nolimits_\beta
 u^\alpha u^\beta\ot v_\beta'*v_\alpha' = 0
\end{equation}
holds in $B^2\ot_R F_2\tA$.

 We have $B^2=\Hom_R(I,R)$.
 Hence, in order to prove
the equation~\eqref{nonhomogeneous-duality-equation-with-indices},
it suffices to check that the pairing of the left-hand side with
any element $i\in I$ vanishes as an element of $F_2\tA$.
 As in Section~\ref{self-consistency-subsecn}, we denote by
$\hI\subset V\ot_\boZ V$ the full preimage of
the $R$\+$R$\+subbimodule $I\subset V\ot_RV$ under the natural
surjective map $V\ot_\boZ V\rarrow V\ot_RV$.
 Consider an element $i\in I$ and choose its preimage $\hi\in\hI$.
 As in Section~\ref{self-consistency-subsecn}, we write symbolically
$i=i_1\ot i_2$ and $\hi=\hi_1\ot\hi_2$.

 Now we deal we the three summands
in~\eqref{nonhomogeneous-duality-equation-with-indices} one by one.
 Firstly, $\lan i,h\ran=h(i)$.
 Secondly,
\begin{multline*}
 \lan i,\sum\nolimits_\alpha d(u^\alpha)\ot v_\alpha'\ran
 =\sum\nolimits_\alpha\lan i,d(u^\alpha)\ran v_\alpha' \\
 \,\overset{\eqref{d-one-defined}}=\,
 \sum\nolimits_\alpha\lan p(\hi_1\ot\hi_2),u^\alpha\ran v_\alpha'
 -\sum\nolimits_\alpha q(\hi_1,\lan \hi_2,u^\alpha\ran)v_\alpha' \\
 = p(\hi_1\ot\hi_2)'
 -\sum\nolimits_\alpha q(\hi_1,\lan \hi_2,u^\alpha\ran)v_\alpha'
\end{multline*}
according to the formula~\eqref{d-one-defined} in
Proposition~\ref{nonhomogeneous-dual-cdg-ring}.
 Thirdly,
\begin{multline*}
 \lan i,\,\sum\nolimits_\alpha\sum\nolimits_\beta
 u^\alpha u^\beta\ot v_\beta'*v_\alpha'\ran =
 \sum\nolimits_\alpha\sum\nolimits_\beta
 \lan i_1\ot i_2,u^\alpha u^\beta\ran
 v_\beta'*v_\alpha' \\
 = \sum\nolimits_\alpha\sum\nolimits_\beta
 \lan i_1,\lan i_2,u^\alpha\ran u^\beta\ran v_\beta'*v_\alpha'
 \,\overset{\eqref{left-right-actions-on-e-equation}}=\,
 \sum\nolimits_\alpha\sum\nolimits_\beta
 \lan i_1,u^\beta\ran (v_\beta\lan i_2,u^\alpha\ran)'*v_\alpha' \\
 \,\overset{\eqref{map-q}}=\,
 \sum\nolimits_\alpha\sum\nolimits_\beta
 \lan\hi_1,u^\beta\ran v_\beta'*\lan\hi_2,u^\alpha\ran*v_\alpha'
 - \sum\nolimits_\alpha\sum\nolimits_\beta
 \lan\hi_1,u^\beta\ran q(v_\beta,\lan\hi_2,u^\alpha\ran)*v_\alpha' \\
 \,\overset{\eqref{map-q}}=\,
 \sum\nolimits_\alpha\sum\nolimits_\beta
 \lan\hi_1,u^\beta\ran v_\beta'*\lan\hi_2,u^\alpha\ran v_\alpha'
 - \sum\nolimits_\alpha\sum\nolimits_\beta
 \lan\hi_1,u^\beta\ran q(v_\beta,\lan\hi_2,u^\alpha\ran)v_\alpha' \\
 =\hi_1'*\hi_2'
 -\sum\nolimits_\alpha q(\hi_1,\lan\hi_2,u^\alpha\ran)v_\alpha',
\end{multline*}
where the equation~\eqref{left-right-actions-on-e-equation} is
being applied to the expression $e=\sum_\beta u^\beta\ot v_\beta$
and the elements $s=\lan i_2,u^\alpha\ran\in R$.

 Finally, we recall that
$$
 -h(i)+p(\hi_1\ot\hi_2)'-\hi_1'*\hi_2'=0
$$
in $F_2\tA$ by the equation~\eqref{maps-p-and-h}.
 This finishes the proof of the desired
equation~\eqref{nonhomogeneous-duality-equation-with-indices} and
of the whole proposition.
\end{proof}

 Let us show that the left CDG\+module $\Ksp(B,\tA)$ which we have
constructed does not depend on any arbitrary choices.
 Let $V''\subset\tV$ be another left $R$\+submodule of strict
generators, related to our original choice of the left $R$\+submodule
$V'\subset\tV$ by the rule~\eqref{change-of-strict-generators-eqn}.
 Let $\prB=(B,d',h')$ and $\secB=(B,d'',h'')$ denote the CDG\+ring
 structures on the graded ring $B$ corresponding to
$V'\subset\tV$ and $V''\subset\tV$, as in
Proposition~\ref{strict-generators-connection-change}.
 Let $d_{e'}$ and~$d_{e''}$ denote the related two differentials on
the tensor product $B\ot_R\tA$, endowing it with the structures of
a left CDG\+module over $\prB$ and~$\secB$.

\begin{lem}
 The restriction of scalars (as defined in
Section~\ref{cdg-modules-subsecn}) with respect to the CDG\+ring 
isomorphism $(\id,a)\:\secB\rarrow\prB$ from
Proposition~\ref{strict-generators-connection-change} transforms
the left CDG\+module $\Ksp_{e'}(B,\tA)=(B\ot_R\tA,\>d_{e'})$ over
the CDG\+ring\/ $\prB$ into the left CDG\+module
$\Ksp_{e''}(B,\tA)=(B\ot_R\tA,\>d_{e''})$ over the CDG\+ring\/~$\secB$.
\end{lem}

\begin{proof}
 We use the notation of Section~\ref{change-of-strict-gens-subsecn}:
given an element $v\in V$, the corresponding elements in $V'$ and
$V''$ are denoted by $v'\in V'\subset\tV$ and $v''\in V''\subset\tV$.
 Let $e'$ and $e''\in\Hom_R(V,R)\ot_R\tV$ denote the elements
corresponding to the injective $R$\+module maps
$V\simeq V'\hookrightarrow\tV$ and $V\simeq V''\hookrightarrow\tV$.
 Then for every $v\in V$ we have
$$
 \lan v,e'\ran =v'\quad
 \text{and}\quad
 \lan v,e''\ran = v''=v'+a(v)
$$
by formula~\eqref{change-of-strict-generators-eqn}, hence
$$
 e''=e'+a\ot 1\,\in\,B^1\ot_RF_1\tA.
$$
 Now for all $b\in B$ and $c\in\tA$ we can compute that
\begin{multline*}
 d_{e''}(b\ot c)=d''(b)\ot c+(-1)^{|b|}be''c \\
 = d'(b)\ot c+(ab-(-1)^{|b|}ba)\ot c+
 (-1)^{|b|}be'c+(-1)^{|b|}ba\ot c \\
 =d'(b)\ot c+ab\ot c+(-1)^{|b|}be'c=d_{e'}(b\ot c)+ab\ot c,
\end{multline*}
as desired (where the middle equation takes into account
the equation~(iv) from the definition of a morphism of CDG\+rings in
Section~\ref{curved-dg-rings-subsecn}).
\end{proof}

\begin{rem} \label{nonhomogeneous-Koszul-complexes-remark}
 The dual nonhomogeneous Koszul CDG\+module $\Ksp_{e'}(B,\tA)$ which
we have constructed is a nonhomogeneous generalization of the dual
Koszul complex $\Ksp_e(B,A)$ from
Section~\ref{dual-koszul-complex-subsecn}.
 The CDG\+(bi)module $\Ksp_{e'}(B,\tA)$ will be very convenient to use 
below in Section~\ref{koszul-duality-functors-comodule-side-subsecn}
for the purpose of constructing the nonhomogeneous Koszul duality
DG\+functors which we are really interested in.
 But the CDG\+module $\Ksp_{e'}(B,\tA)$ itself, even when it happens
to be a complex, has no good exactness properties
(see Remark~\ref{dual-koszul-complex-no-good}).
 However, applying the functor $\Hom_{\tA^\rop}({-},\tA)$ to
$\Ksp_{e'}(B,\tA)$ produces a right CDG\+module
$\Hom_{\tA^\rop}(\Ksp_{e'}(B,\tA),\tA)$ over $(B,d,h)$ (in fact,
a CDG\+bimodule over $(\tA,0,0)$ and $(B,d,h)$), whose underlying
graded $\tA$\+$B$\+bimodule is $\tA\ot_R\Hom_{R^\rop}(B,R)$.
 This CDG\+module, which is a nonhomogeneous generalization of
the first Koszul complex $K^\tau_\bu(B,A)$ from
Section~\ref{first-koszul-complex-subsecn}, will appear below in
the proof of Theorem~\ref{comodule-side-koszul-duality-theorem}
and the subsequent discussion in
Section~\ref{comodule-side-triangulated-equivalence-subsecn}. 
\end{rem}

\subsection{Semicoderived category of modules}
\label{semicoderived-category-subsecn}
 Let $R$ be an associative ring.
 The following definition of the \emph{coderived category of right
$R$\+modules} goes back to~\cite[Section~2.1]{Psemi}
and~\cite[Section~3.3]{Pkoszul} (see~\cite[Section~2]{Kra}
and~\cite[Proposition~1.3.6(2)]{Bec} for an alternative approach).

 As a particular case of the notation of
Section~\ref{cdg-modules-subsecn}, we denote by $\Hot(\modr R)$
the homotopy category of (unbounded complexes of infinitely generated)
right $R$\+modules.
 Let us consider short exact sequences of complexes of right
$R$\+modules $0\rarrow L^\bu\rarrow M^\bu\rarrow N^\bu\rarrow0$.
 Any such short exact sequence can be viewed as a bicomplex with
three rows, and the totalization (the total complex) $\Tot(L^\bu
\to M^\bu\to N^\bu)$ of such a bicomplex can be constructed
in the usual way.
 A complex of right $R$\+modules is said to be \emph{coacyclic} if it
belongs to the minimal full triangulated subcategory of $\Hot(\modr R)$
containing all the totalizations of short exact sequences of complexes
of right $R$\+modules and closed under infinite direct sums.

 The full triangulated subcategory of coacyclic complexes of right
$R$\+modules is denoted by $\Acycl^\co(\modr R)\subset\Hot(\modr R)$.
 The \emph{coderived category of right $R$\+modules} is the triangulated
Verdier quotient category of the homotopy category $\Hot(\modr R)$ by
the triangulated subcategory of coacyclic complexes,
$$
 \sD^\co(\modr R)=\Hot(\modr R)/\Acycl^\co(\modr R).
$$
 Obviously, any coacyclic complex of $R$\+modules is acyclic,
$\Acycl^\co(\modr R)\subset\Acycl(\modr R)$; but the converse
does not hold in general~\cite[Examples~3.3]{Pkoszul}.
 So the conventional derived category $\sD(\modr R)$ is a quotient
category of $\sD^\co(\modr R)$.
 All acyclic complexes of right $R$\+modules are coacyclic when
the right homological dimension of the ring $R$ is
finite~\cite[Remark~2.1]{Psemi}.

 The following lemma is useful.

\begin{lem} \label{complex-filtration-coacyclic-lemma}
 Let\/ $0=F_{-1}M^\bu\subset F_0M^\bu\subset F_1M^\bu\subset\dotsb
\subset M^\bu$ be a complex of $R$\+modules endowed with an exhaustive
increasing filtration by subcomplexes of $R$\+modules; so
$M^\bu=\bigcup_{n\ge0}F_nM^\bu$.
 Assume that, for every $n\ge0$, the complex of $R$\+modules
$F_nM^\bu/F_{n-1}M^\bu$ is coacyclic.
 Then the complex of $R$\+modules $M^\bu$ is coacyclic, too.
\end{lem}

\begin{proof}
 First one proves by induction in~$n$ that the complex of
$R$\+modules $F_nM^\bu$ is coacyclic for every $n\ge0$, using
the fact that the totalizations of short exact sequences $0\rarrow
F_{n-1}M^\bu\rarrow F_nM^\bu\rarrow F_nM^\bu/F_{n-1}M^\bu\rarrow0$
are coacyclic.
 Then one concludes that the complex of $R$\+modules
$\bigoplus_{n=0}^\infty F_nM^\bu$ is coacyclic.
 Finally, one deduces coacyclicity of the complex $M^\bu$ from
the fact that the totalization of the telescope short exact
sequence $0\rarrow\bigoplus_{n=0}^\infty F_nM^\bu\rarrow
\bigoplus_{n=0}^\infty F_nM^\bu\rarrow M^\bu\rarrow0$ is
a coacyclic complex of $R$\+modules.
\end{proof}

 Now let $R\rarrow A$ be a homomorphism of associative rings.
 The following definition of the \emph{$A/R$\+semicoderived category
of right $A$\+modules} can be found in~\cite[Section~5]{Pfp}.
 The semicoderived category
$$
 \sD^\sico_R(\modr A)=\Hot(\modr A)/\Acycl^\sico_R(\modr A)
$$
is the triangulated Verdier quotient category of the homotopy
category of right $A$\+modules by its full triangulated subcategory
$\Acycl^\sico_R(\modr A)$ consisting of all the complexes of right
$A$\+modules that are \emph{coacyclic as complexes of
right $R$\+modules}.

 Obviously, any coacyclic complex of $A$\+modules is coacyclic
as a complex of $R$\+modules; but the converse does not need to be true.
 So the semicoderived category is intermediate between the derived and
the coderived category of $A$\+modules: there are natural triangulated
Verdier quotient functors $\sD^\co(\modr A)\rarrow\sD^\sico_R(\modr A)
\rarrow\sD(\modr A)$.
 (We refer to the paper~\cite{Pps} for a more general discussion of
such intermediate triangulated quotient categories, called
\emph{pseudo-coderived categories} in~\cite{Pps}.)
 The semicoderived category $\sD^\sico_R(\modr A)$ can be thought of
as a mixture of the coderived category ``along the variables from~$R$''
and the conventional derived category ``in the direction of $A$
relative to~$R$''.

 For example, taking $R=\boZ$, one obtains the conventional derived
category of $A$\+modules
$$
 \sD^\sico_\boZ(\modr A)=\sD(\modr A)
$$
(because all the acyclic complexes of abelian groups are coacyclic),
while taking $R=A$ one obtains the coderived category of $A$\+modules
$$
 \sD^\sico_A(\modr A)=\sD^\co(\modr A).
$$

\subsection{Coderived category of CDG-comodules}
\label{coderived-cdg-comodules-subsecn}
 Let $B=(B,d,h)$ be a nonnegatively graded CDG\+ring (as defined in
Section~\ref{curved-dg-rings-subsecn}).
 So $B=\bigoplus_{n=0}^\infty B^n$ is a nonnegatively graded
associative ring, $d\:B\rarrow B$ is an odd derivation of degree~$1$,
and $h\in B^2$ is a curvature element for~$d$.

 Let $(N,d_N)$ be a right CDG\+module over $(B,d,h)$ (as defined in
Section~\ref{cdg-modules-subsecn}).
 So $N=\bigoplus_{n\in\boZ} N^n$ is a graded right $B$\+module,
$d_N\:N\rarrow N$ is an odd derivation of degree~$1$ compatible with
the derivation~$d$ on~$B$, and the equation $d_N^2(y)=-yh$ is
satisfied for all $y\in N$.
 We will say that the CDG\+module $(N,d_N)$ is a \emph{CDG\+comodule}
(or a \emph{right CDG\+comodule}) over $(B,d,h)$ if the graded right
$B$\+module $N$ is a graded right $B$\+comodule (in the sense of
Section~\ref{graded-comodules-subsecn}).

 In other words, a right CDG\+comodule $N$ over $(B,d,h)$ is
a graded right $B$\+comod\-ule endowed with an odd coderivation
$d_N\:N\rarrow N$ of degree~$1$ compatible with the derivation~$d$
on $B$ and satisfying the above equation for~$d_N^2$.
 Here by an \emph{odd coderivation} of a graded $B$\+comodule $N$
(compatible with a given odd derivation~$d$ on~$B$) we mean an odd
derivation of $N$ as a graded $B$\+module.
 Equivalently, a right CDG\+comodule over $(B,d,h)$ can be simply
defined as a graded right $\hB$\+comodule (where $(\hB,\d)$ is
the quasi-differential graded ring corresponding to $(B,d,h)$,
as in Theorem~\ref{cdg-qdg-equivalence} and
Section~\ref{cdg-modules-subsecn}).

 The DG\+category $\DG(\modr B)$ of right CDG\+modules over $B=(B,d,h)$
was defined in Section~\ref{cdg-modules-subsecn}.
 The \emph{DG\+category of right CDG\+comodules} $\DG(\comodr B)$
\emph{over $(B,d,h)$} is defined as the full DG\+subcategory of
the DG\+category $\DG(\modr B)$ whose objects are the right
CDG\+comodules over $(B,d,h)$.
 Obviously, the full DG\+subcategory $\DG(\comodr B)\subset
\DG(\modr B)$ is closed under shifts, twists, and infinite direct
sums; in particular, it is closed under the cones of closed morphisms
(see~\cite[Section~1.2]{Pkoszul} for the terminology).

 Let $\Hot(\comodr B)=H^0\DG(\modr B)$ denote the homotopy category
(of the DG\+category) of right CDG\+comodules over~$B$.
 So $\Hot(\comodr B)$ is a full triangulated subcategory in
$\Hot(\modr B)$.
 The following definition of the \emph{coderived category of
right CDG\+comodules over~$B$} is a variation on the definitions
in~\cite[Sections~11.7.1\+-2]{Psemi}
and~\cite[Sections~3.3 and~4.2]{Pkoszul}.

 Let us consider short exact sequences $0\rarrow L\overset i\rarrow M
\overset p\rarrow N\rarrow0$ of right CDG\+comodules over $(B,d,h)$.
 Here $i$ and~$p$ are closed morphisms of degree zero in
the DG\+category $\DG(\comodr B)$ and the short sequence $0\rarrow L
\rarrow M\rarrow N\rarrow 0$ is exact in the abelian category
$\comodrgr B$.
 Using the construction of the cone of a closed morphism (or
the shift, finite direct sum, and twist) in the DG\+category
$\DG(\comodr B)$, one can produce the totalization (or
the \emph{total CDG\+comodule}) $\Tot(L\to M\to N)$ of the finite
complex of CDG\+comodules $L\rarrow M\rarrow N$.

 A right CDG\+comodule over $B$ is said to be \emph{coacyclic} if it
belongs to the minimal full triangulated subcategory of
$\Hot(\comodr B)$ containing all the totalizations of short exact
sequences of right CDG\+comodules over $B$ and closed under infinite
direct sums.
 The full triangulated subcategory of coacyclic CDG\+comodules is
denoted by $\Acycl^\co(\comodr B)\subset\Hot(\comodr B)$.
 The \emph{coderived category of right CDG\+comodules over~$B$}
is the triangulated Verdier quotient category of the homotopy
category $\Hot(\comodr B)$ by its triangulated subcategory of
coacyclic CDG\+comodules,
$$
 \sD^\co(\comodr B)=\Hot(\comodr B)/\Acycl^\co(\comodr B).
$$

 Notice that, generally speaking, the differential $d_N$ on
a CDG\+(co)module $(N,d_N)$ has nonzero square (if $B$ is not
a DG\+ring, that is $h\ne0$ in~$B$).
 So CDG\+(co)modules are \emph{not complexes} and their cohomology
(co)modules are \emph{undefined}.
 For this reason, one \emph{cannot} speak about ``acyclic
CDG\+(co)modules'', \emph{nor} about ``quasi-isomorphisms of
CDG\+(co)modules'' in the usual sense of the word ``quasi-isomorphism''.
 Therefore, the conventional construction of the derived category
of DG\+modules is \emph{not} applicable to curved DG\+modules or
curved DG\+comodules.
 There is only the coderived category and its variations (known as
\emph{derived categories of the second kind}).

\begin{lem} \label{cdg-comodule-filtration-coacyclic-lemma}
 Let $B=(B,d,h)$ be a nonnegatively graded CDG\+ring and $N=(N,d_N)$
be a right CDG\+comodule over~$B$.
 Let\/ $0=F_{-1}N\subset F_0N\subset F_1N\subset\dotsb\subset N$
be an exhaustive increasing filtration of $N$ by CDG\+subcomodules,
that is $(F_nN)B\subset F_nN$, \ $d_N(F_nN)\subset F_nN$, and
$N=\bigcup_{n=0}^\infty F_n N$.
 Assume that, for every $n\ge0$, the CDG\+comodule $F_nN/F_{n-1}N$
is coacyclic over~$B$.
 Then the CDG\+comodule $N$ is coacyclic, too.
\end{lem}

\begin{proof}
 This is a generalization of
Lemma~\ref{complex-filtration-coacyclic-lemma}, provable by
the same argument.
\end{proof}

 There is an important particular case when the nonnegatively graded
ring $B$ has only finitely many grading components, that is $B^n=0$
for $n\gg0$.
 In this case, all the (ungraded or graded) $B$\+modules are
(respectively, ungraded or graded) $B$\+comodules, so
$\comodrgr B=\modrgr B$.
 Accordingly, the full DG\+subcategory of (right) CDG\+comodules over
$B$ coincides with the whole ambient DG\+category of
(right) CDG\+modules, $\DG(\comodr B)=\DG(\modr B)$.
 In this case, we will write simply $\sD^\co(\modr B)$ instead of
$\sD^\co(\comodr B)$.

\subsection{Koszul duality functors for modules and comodules}
\label{koszul-duality-functors-comodule-side-subsecn}
 Let $R\subset\tV\subset\tA$ be a $3$\+left finitely projective
weak nonhomogeneous quadratic ring, and let $B=(B,d,h)$ be
the nonhomogeneous quadratic qual CDG\+ring.
 The dual nonhomogeneous Koszul CDG\+module $\Ksp(B,\tA)=
\Ksp_{e'}(B,\tA)=(B\ot_R\tA,\>d_{e'})$ constructed in
Section~\ref{nonhomogeneous-koszul-cdg-module-subsecn} plays a key
role below.
 We recall from the discussion in
Section~\ref{nonhomogeneous-koszul-cdg-module-subsecn} that
$\Ksp(B,\tA)$ is a CDG\+bimodule over the CDG\+rings $B=(B,d,h)$
and $\tA=(\tA,0,0)$.

 According to Section~\ref{cdg-modules-subsecn}, we have the induced
pair of adjoint DG\+functors between the DG\+category $\DG(\modr\tA)$
of complexes of right $\tA$\+modules and the DG\+category $\DG(\modr B)$
of right CDG\+modules over $(B,d,h)$.
 Here the left adjoint DG\+functor
$$
 {-}\ot_B\Ksp(B,\tA)\:\DG(\modr B)\lrarrow\DG(\modr\tA)
$$
takes a right CDG\+module $N=(N,d_N)$ over $(B,d,h)$ to the complex
of right $\tA$\+modules whose underlying graded $\tA$\+module is
$$
 N\ot_B\Ksp(B,\tA)=N\ot_B(B\ot_R\tA)=N\ot_R\tA.
$$
 The right adjoint DG\+functor
$$
 \Hom_{\tA^\rop}(\Ksp(B,\tA),{-})\:\DG(\modr\tA)\lrarrow\DG(\modr B)
$$
takes a complex of right $\tA$\+modules $M=M^\bu$ to the right
CDG\+module over $(B,d,h)$ whose underlying graded $B$\+module is
$$
 \Hom_{\tA^\rop}(\Ksp(B,\tA),M)=\Hom_{\tA^\rop}(B\ot_R\tA,\>M)
 =\Hom_{R^\rop}(B,M).
$$

 Our aim is to restrict this adjoint pair to the DG\+subcategory
$\DG(\comodr B)\subset\DG(\modr B)$ of CDG\+comodules over $(B,d,h)$,
as defined in Sections~\ref{graded-comodules-subsecn}
and~\ref{coderived-cdg-comodules-subsecn}.

 We recall the notation $\Hom_{R^\rop}^\Sigma(L,M)\subset
\Hom_{R^\rop}(L,M)$ for the direct sum totalization of the bigraded
Hom group of two graded right $R$\+modules $L$ and $M$, which was
introduced in Section~\ref{graded-ext-comparison-subsecn}.
 Let $\Gamma_B\:\modrgr B\rarrow\comodrgr B$ denote the functor
assigning to every graded right $B$\+module its maximal (graded)
submodule which is a (graded) right $B$\+comodule (see the proof of
Theorem~\ref{graded-comodule-Ext-comparison}).

\begin{lem} \label{gamma-hom-hom-sigma}
 Let $B=\bigoplus_{n=0}^\infty B_n$ be a nonnegatively graded ring
with the degree-zero component $R=B_0$.
 Then for any graded right $R$\+module $M$ one has
$$
 \Gamma_B(\Hom_{R^\rop}(B,M))=\Hom_{R^\rop}^\Sigma(B,M).
$$ 
\end{lem}

\begin{proof}
 We have $\Hom_{R^\rop}^\Sigma(B,M)=\bigoplus_{n\in\boZ}
\Hom_{R^\rop}(B,M_n(n))$, where $M\longmapsto M(n)$ denotes
the functor of grading shift.
 Since the nonpositively graded right $B$\+module $\Hom_{R^\rop}(B,L)$
is a graded right $B$\+comodule for any ungraded right $R$\+module~$L$,
and the full subcategory of graded right $B$\+comodules $\comodrgr B$
is closed under direct sums in $\modrgr B$, it follows that
$\Hom_{R^\rop}^\Sigma(B,M)$ is a graded right $B$\+comodule.
 Hence $\Hom_{R^\rop}^\Sigma(B,M)\subset\Gamma_B(\Hom_{R^\rop}(B,M))$.
 Conversely, let $f\in\Hom_{R^\rop}(B,M)_n$ be a homogeneous
element of degree~$n$ such that $fB_{\ge m}=0$ for some $m\ge1$.
 Then for every $b\in B_{\ge m}$ we have $f(b)=f(b\cdot 1)=(fb)(1)=0$,
where $1\in R\subset B$ is the unit element.
 It follows immediately that $f\in\Hom_{R^\rop}^\Sigma(B,M)_n
\subset\Hom_{R^\rop}(B,M)_n$.
\end{proof}

 For any complex of right $A$\+modules $M^\bu$, the graded right
$B$\+submodule
$$
 \Hom_{\tA^\rop}^\Sigma(\Ksp(B,\tA),M^\bu)=
 \Hom_{R^\rop}^\Sigma(B,M^\bu)\subset\Hom_{R^\rop}(B,M^\bu)
$$
is preserved by the differential of the CDG\+module
$\Hom_{\tA^\rop}(\Ksp(B,\tA),M^\bu)=\Hom_{R^\rop}(B,M^\bu)$
over $(B,d,h)$ (so it is a CDG\+submodule).
 In view of the above discussion and Lemma~\ref{gamma-hom-hom-sigma},
we obtain a pair of adjoint DG\+functors between the DG\+category
$\DG(\modr\tA)$ of complexes of right $\tA$\+modules and
the DG\+category $\DG(\comodr B)$ of right
CDG\+comodules over~$(B,d,h)$.

 Here the left adjoint DG\+functor
$$
 {-}\ot_B\Ksp(B,\tA)\:\DG(\comodr B)\lrarrow\DG(\modr\tA)
$$
is simply the restriction of the above DG\+functor to the full
DG\+subcategory $\DG(\comodr B)\subset\DG(\modr B)$.
 We denote it by
$$
 N\longmapsto N\ot_B\Ksp(B,\tA)=N\ot_R^{\tau'}\tA
$$
where the placeholder~$\tau'$ is understood as a notation for
the injective left $R$\+module map $\Hom_{R^\rop}(B^1,R)=
F_1\tA/F_0\tA=V\simeq V'\hookrightarrow \tV=F_1\tA$.

 The right adjoint DG\+functor
$$
 \Hom_{\tA^\rop}^\Sigma(\Ksp(B,\tA),{-})\:\DG(\modr\tA)
 \lrarrow\DG(\comodr B)
$$
takes a complex of right $\tA$\+modules $M=M^\bu$ to the right
CDG\+comodule over $(B,d,h)$ whose underlying graded right
$B$\+comodule is
$$
 \Hom_{\tA^\rop}^\Sigma(B\ot_R\tA,\>M)=\Hom_{R^\rop}^\Sigma(B,M)
 =M\ot_R\Hom_{R^\rop}(B,R)=M\ot_RC,
$$
where $C$ is a notation for the graded $R$\+$B$\+bimodule
(and a graded right $B$\+comodule) $C=\Hom_{R^\rop}(B,R)$.
 We denote this right CDG\+comodule over $(B,d,h)$ by
$$
 \Hom_{\tA^\rop}^\Sigma(\Ksp(B,\tA),M^\bu)=
 M^\bu\ot_R^{\sigma'}\Hom_{R^\rop}(B,R)=M^\bu\ot_R^{\sigma'}C,
$$
where $\sigma'$~is, for the time being, just a placeholder.
 The purpose of this placeholder is to remind us of the nontrivial
differential~$d_{e'}$ on $\Ksp(B,\tA)$, which is incorporated
into the construction of the differential on the CDG\+comodule
$M^\bu\ot_R^{\sigma'}C$ (cf.\ the twisting cochain notation
in the simpler context of~\cite[Section~6]{Pkoszul}).
 See~\eqref{sigma-prime-explained-eqn} and
Section~\ref{revisited-subsecn} below for our suggested interpretation
of~$\sigma'$.

 The pair of adjoint DG\+functors that we have constructed induces
a pair of adjoint triangulated functors
$$
 {-}\ot_B\Ksp(B,\tA)={-}\ot_R^{\tau'}\tA\:
 \Hot(\comodr B)\lrarrow\Hot(\modr\tA)
$$
and
$$
 \Hom_{\tA^\rop}^\Sigma(\Ksp(B,\tA),{-})={-}\ot_R^{\sigma'}C\:
 \Hot(\modr\tA)\lrarrow\Hot(\comodr B)
$$
between the homotopy category $\Hot(\modr\tA)$ of complexes of
right $\tA$\+modules and the homotopy category $\Hot(\comodr B)$
of right CDG\+comodules over~$(B,d,h)$.

\subsection{Triangulated equivalence}
\label{comodule-side-triangulated-equivalence-subsecn}
 Let $R\subset\tV\subset\tA$ be a left finitely projective
nonhomogeneous Koszul ring (as defined in
Section~\ref{pbw-theorem-subsecn}), and let $(B,d,h)$ be
the corresponding right finitely projective Koszul CDG\+ring
(as constructed in Proposition~\ref{nonhomogeneous-dual-cdg-ring}; see
also Corollary~\ref{nonhomogeneous-koszul-duality-anti-equivalence}).
 The following theorem is the main result of
Section~\ref{comodule-side-secn}.

\begin{thm} \label{comodule-side-koszul-duality-theorem}
 The pair of adjoint triangulated functors\/ ${-}\ot_R^{\sigma'}C\:
\Hot(\modr\tA)\rarrow\Hot(\comodr B)$ and\/ ${-}\ot_R^{\tau'}\tA\:
\Hot(\comodr B)\rarrow\Hot(\modr\tA)$ defined in
Section~\ref{koszul-duality-functors-comodule-side-subsecn}
induces a pair of adjoint triangulated functors between
the $\tA/R$\+semicoderived category of right $\tA$\+modules\/
$\sD^\sico_R(\modr\tA)$, as defined
in Section~\ref{semicoderived-category-subsecn}, and
the coderived category\/ $\sD^\co(\comodr B)$ of right CDG\+comodules
over $(B,d,h)$, as defined in
Section~\ref{coderived-cdg-comodules-subsecn}.
 Under the above Koszulity assumption, the latter two functors are
mutually inverse triangulated equivalences,
$$
 \sD_R^\sico(\modr\tA)\simeq\sD^\co(\comodr B).
$$
\end{thm}

\begin{proof}
 This is essentially a generalization of~\cite[Theorem~B.2(a)]{Pkoszul}
and a particular case of~\cite[Theorem~11.8(b)]{Psemi}.
 The proof is similar.
 Let us spell out some details.
 
 Recall from the proof of Theorem~\ref{graded-comodule-Ext-comparison}
that a graded right $B$\+comodule $N$ is said to be \emph{trivial} if
$NB^{\ge1}=0$ (where the notation is
$B^{\ge m}=\bigoplus_{n\ge m} B^n\subset B$ for any integer $m\ge0$).
 A right CDG\+comodule $(N,d_N)$ over $(B,d,h)$ is said to be
\emph{trivial} if its underlying graded $B$\+comodule is trivial.
 The differential~$d_N$ on a trivial CDG\+comodule $N$ squares to zero,
as the curvature element~$h$ acts by zero in~$N$.
 The DG\+category of trivial right CDG\+comodules over $(B,d,h)$ is
equivalent to the DG\+category of complexes of right $R$\+modules.

 The ring $\tA$ is endowed with an increasing filtration~$F$.
 Let us also introduce an increasing filtration $F$ on
the graded $R$\+$B$\+bimodule $C=\Hom_{R^\rop}(B,R)$.
 We put $F_nC=\bigoplus_{i\le n}\Hom_{R^\rop}(B^i,R)\subset C$
for every $n\in\boZ$.
 So $F_{-1}C=0$ and $C=\bigcup_{n\ge0}F_nC$.
 Clearly, $F_nC$ is a graded $R$\+$B$\+subbimodule in~$C$.
 Notice that both the left $R$\+modules $F_n\tA/F_{n-1}\tA$ and
the graded left $R$\+modules $F_nC/F_{n-1}C$ are finitely generated
and projective (it is important for the argument below that
they are flat).

 Let $M^\bu$ be a complex of right $\tA$\+modules that is coacyclic as
a complex of right $R$\+modules.
 We need to show that the right CDG\+comodule $M^\bu\ot_R^{\sigma'}C$
is coacyclic over~$(B,d,h)$.
 Let $F$ be the increasing filtration on the tensor product
$M\ot_RC$ induced by the increasing filtration $F$ on~$C$, that is
$F_n(M\ot_RC)=M\ot_RF_nC$.
 Then $F$ is a filtration of the CDG\+comodule $M^\bu\ot_R^{\sigma'}C$
by CDG\+subcomodules over $(B,d,h)$.
 The successive quotient CDG\+comodules $(M^\bu\ot_R^{\sigma'}F_nC)/
(M^\bu\ot_R^{\sigma'}F_{n-1}C)$ are trivial.
 Viewed as complexes of right $R$\+comodules, these are the tensor
products $M^\bu\ot_R(F_nC/F_{n-1}C)$, where $F_nC/F_{n-1}C$ is
a graded $R$\+$R$\+bimodule concentrated in the cohomological
degree~$-n$ and endowed with the zero differential.

 Since the graded left $R$\+module $F_nC/F_{n-1}C$ is flat, tensoring
with it preserves short exact sequences of graded right $R$\+modules.
 The tensor product functor also preserves infinite direct sums.
 It follows that tensoring with the graded $R$\+$R$\+bimodule
$F_nC/F_{n-1}C$ takes coacyclic complexes of right $R$\+modules to
coacyclic complexes of right $R$\+modules.
 Clearly, any coacyclic complex of right $R$\+modules is also coacyclic
as a trivial CDG\+comodule over $(B,d,h)$.
 It remains to make use of
Lemma~\ref{cdg-comodule-filtration-coacyclic-lemma} in order to
conclude that the right CDG\+comodule $M^\bu\ot_R^{\sigma'}C$ over
$(B,d,h)$ is coacyclic.

 For any coacyclic right CDG\+comodule $N$ over $(B,d,h)$, the complex
of right $\tA$\+modules $N\ot_R^{\tau'}\tA$ is coacyclic not only as
a complex of right $R$\+modules, but even as a complex of right
$\tA$\+modules.
 This follows immediately from the fact that $\tA$ is a flat left
$R$\+module; so the functor ${-}\ot_R^{\tau'}\tA$ takes short exact
sequences of right CDG\+comodules over $(B,d,h)$ to short exact
sequences of complexes of right $\tA$\+modules (and this functor also
preserves infinite direct sums).

 Let $M^\bu$ be arbitrary complex of right $\tA$\+modules.
 We need to show the cone $Y^\bu$ of the natural closed morphism of
complexes of right $\tA$\+modules (the adjunction morphism)
$M^\bu\ot_R^{\sigma'}C\ot_R^{\tau'}\tA\rarrow M^\bu$ is coacyclic
as a complex of right $R$\+modules.
 Endow the graded right $R$\+module $M\ot_RC\ot_R\tA$ with
the increasing filtration $F$ induced by the increasing filtrations $F$
on $C$ and $\tA$, that is $F_n(M\ot_RC\ot_R\tA)=\sum_{i+j=n}
M\ot_RF_iC\ot_RF_j\tA\subset M\ot_RC\ot_R\tA$.
 Then $F$ is an exhaustive filtration of
$M^\bu\ot_R^{\sigma'}C\ot_R^{\tau'}\tA$ by subcomplexes of right
$R$\+modules (but not of right $\tA$\+modules!).
 The complex $M^\bu$ is endowed with the trivial filtration
$F_{-1}M^\bu=0$, \ $F_0M^\bu=M^\bu$; and the complex $Y^\bu$ is
endowed with the induced increasing filtration~$F$.

 Then the associated graded complex of right $R$\+modules
$$
 \gr^F(M^\bu\ot_R^{\sigma'}C\ot_R^{\tau'}\tA)=
 \bigoplus\nolimits_{n=0}^\infty
 F_n(M^\bu\ot_R^{\sigma'}C\ot_R^{\tau'}\tA)\big/
 F_{n-1}(M^\bu\ot_R^{\sigma'}C\ot_R^{\tau'}\tA)
$$
is naturally isomorphic to the tensor product
$M^\bu\ot_R\,{}^\tau\!K_\bu(B,A)$ of the complex of right
$R$\+modules $M^\bu$ and the complex of (graded) $R$\+$R$\+bimodules
${}^\tau\!K_\bu(B,A)=C\ot^\tau_RA$ constructed in
Section~\ref{dual-koszul-complex-subsecn} (where it is called
``the second Koszul complex'').
 Here, as usually, $A=\gr^F\tA$ is the left finitely projective
Koszul graded ring associated with~$\tA$.
 By Proposition~\ref{left-flat-second-Koszul-complex} or
Theorem~\ref{projective-koszul-theorem}(e), every internal
degree~$n\ge1$ component of the complex ${}^\tau\!K_\bu(B,A)$ is
a finite acyclic complex of $R$\+$R$\+bimodules whose terms
are finitely generated and projective left $R$\+modules.
 The internal degree $n=0$ component of the complex
${}^\tau\!K_\bu(B,A)$ is isomorphic to the $R$\+$R$\+bimodule~$R$.

 It follows that the associated graded complex $\gr^F\,Y^\bu=
\bigoplus_{n\ge0}F_nY^\bu/F_{n-1}Y^\bu$ is isomorphic to
the tensor product $M^\bu\ot_R\cone({}^\tau\!K_\bu(B,A)\to R)$
of the complex of right $R$\+modules $M^\bu$ and the cone of
the natural morphism of complexes of (graded) $R$\+$R$\+bimodules
${}^\tau\!K_\bu(B,A)\rarrow R$.
 For every $n\in\boZ$, the internal degree~$n$ component of
the complex $\cone({}^\tau\!K_\bu(B,A)\to R)$ is a finite acyclic
complex of $R$\+$R$\+bimodules which are finitely generated and
projective as left $R$\+modules.
 Consequently, the complex $\cone({}^\tau\!K_\bu(B,A)\to R)$
is ``coacyclic with respect to the exact category of left
$R$\+flat (or even left $R$\+projective) $R$\+$R$\+bimodules''
in the sense of~\cite[Section~2.1]{Psemi}.
 Since the functor $M^\bu\ot_R{-}$ preserves exactness of short exact
sequences of left $R$\+flat $R$\+$R$\+bimodules, we can conclude that
$\gr^F\,Y^\bu\simeq M^\bu\ot_R\cone({}^\tau\!K_\bu(B,A)\to R)$ is
a coacyclic complex of right $R$\+modules.
 By Lemma~\ref{complex-filtration-coacyclic-lemma}, it follows that
$Y^\bu$ is a coacyclic complex of right $R$\+modules, as desired.

 Let $N=(N,d_N)$ be an arbitrary right CDG\+comodule over $(B,d,h)$.
 We need to show that the cone $Z$ of the natural closed morphism
of right CDG\+comodules (the adjunction morphism)
$N\rarrow N\ot_R^{\tau'}\tA\ot_R^{\sigma'}C$ is coacyclic as
a CDG\+comodule over $(B,d,h)$.
 Endow the graded right $B$\+comodule $N$ with the canonical
increasing filtration by graded subcomodules $0=G_{-1}N\subset G_0N
\subset G_1N\subset G_2N\subset\dotsb\subset N$ as in the proof of
Theorem~\ref{graded-comodule-Ext-comparison}: the grading
component $G_mN^i\subset N^i$ consists of all the elements $y\in N^i$
such that $yB^{\ge m+1}=0$ in~$N$.
 We have $d_N(G_mN)\subset G_mN$, as
$d_N(y)b=d_N(yb)-(-1)^{|y|}yd(b)=0$ for all $y\in G_mN^{|y|}$ and
$b\in B^{\ge m+1}$.
 So $G$ is an exhaustive increasing filtration of $N$ by
CDG\+subcomodules over $(B,d,h)$.

 Consider the induced increasing filtrations $G$ of the right
CDG\+comodule $N\ot_R^{\tau'}\tA\ot_R^{\sigma'}C$ and the right
CDG\+comodule $Z$ over $(B,d,h)$; so
$G_m(N\ot_R^{\tau'}\tA\ot_R^{\sigma'}C)=
(G_mN)\ot_R^{\tau'}\tA\ot_R^{\sigma'}C$ and
$G_mZ=\cone(G_mN\to G_mN\ot_R^{\tau'}\tA\ot_R^{\sigma'}C)$.
 By Lemma~\ref{cdg-comodule-filtration-coacyclic-lemma}, it suffices
to check that the CDG\+comodules $G_mZ/G_{m-1}Z$ are coacyclic over
$(B,d,h)$ for all $m\ge0$.
 Then, since $\tA$ and $C$ are flat (graded) left $R$\+modules,
the CDG\+comodule $G_mZ/G_{m-1}Z$ is naturally isomorphic to
the cone of the adjunction morphism $G_mN/G_{m-1}N\rarrow
(G_mN/G_{m-1}N)\ot_R^{\tau'}\tA\ot_R^{\sigma'}C$. {\hbadness=1600\par}

 The successive quotient CDG\+comodules $G_mN/G_{m-1}N$ are annihilated
by $B^{\ge1}$; so they are trivial right CDG\+comodules over $(B,d,h)$
or, which is the same, just complexes of right $R$\+modules.
 Thus we have reduced our problem to the case of a trivial right
CDG\+comodule $N$, and we can assume that $N=N^\bu$ is simply a complex
of right $R$\+modules with $B^{\ge1}$ acting by zero.
 Then the right CDG\+comodule $N\ot_R^{\tau'}\tA\ot_R^{\sigma'}C=
N^\bu\ot_R\tA\ot_R^{\sigma'}C$ over $(B,d,h)$ is simply the tensor
product of the complex of right $R$\+modules $N^\bu$ and
the right CDG\+comodule $\tA\ot_R^{\sigma'}C=
\Hom_{\tA^\rop}(\Ksp_{e'}(B,\tA),\tA)$ (which is, in fact,
a CDG\+bimodule over $\tA=(\tA,0,0)$ and $B=(B,d,h)$, hence
in particular a CDG\+bimodule over $R=(R,0,0)$ and~$B$).

 Endow the graded right $R$\+module $N\ot_R\tA\ot_RC$ with
the increasing filtration $F$ induced by the increasing filtrations $F$
on $\tA$ and $C$, that is $F_n(N\ot_R\tA\ot_RC)=\sum_{i+j=n}
N\ot_RF_j\tA\ot_RF_iC\subset N\ot_R\tA\ot_RC$.
 Then $F$ is an exhaustive filtration of
$N^\bu\ot_R\tA\ot_R^{\sigma'}C$ by CDG\+subcomodules over $(B,d,h)$.
 Let the complex $N^\bu$ be endowed with the trivial filtration
$F_{-1}N^\bu=0$, \ $F_0N^\bu=N^\bu$, and let the CDG\+comodule
$Z=\cone(N^\bu\to N^\bu\ot_R\tA\ot_R^{\sigma'}C)$ be endowed with
the induced filtration.

 Then the associated graded CDG\+comodule
$$
 \gr^F(N^\bu\ot_R\tA\ot_R^{\sigma'}C)
 =\bigoplus\nolimits_{n=0}^\infty F_n(N^\bu\ot_R\tA\ot_R^{\sigma'}C)
 \big/F_{n-1}(N^\bu\ot_R\tA\ot_R^{\sigma'}C)
$$
is annihilated by the action of $B^{\ge1}$; so it is a trivial
right CDG\+comodule over $(B,d,h)$.
 As a complex of right $R$\+modules, it is naturally isomorphic to
the tensor product $N^\bu\ot_RK^\tau_\bu(B,A)$ of the complex of
right $R$\+modules $N^\bu$ and the complex of (graded)
$R$\+$R$\+bimodules $K^\tau_\bu(B,A)=A\ot_R^\sigma C$ constructed
in Section~\ref{first-koszul-complex-subsecn} (where it is called
``the first Koszul complex'').
 Here, once again, the superindex~$\sigma$ is just a placeholder for
the time being (see Section~\ref{revisited-subsecn} below for
an explanation).
 By Proposition~\ref{left-flat-first-Koszul-complex} or
Theorem~\ref{projective-koszul-theorem}(d), every internal degree
$n\ge1$ component of the complex $K^\tau_\bu(B,A)$ is a finite
acyclic complex of $R$\+$R$\+bimodules with finitely generated
projective underlying left $R$\+modules (while the internal degree
$n=0$ component is the $R$\+$R$\+bimodule~$R$).

 It follows that the associated graded CDG\+comodule $\gr^FZ=
\bigoplus_{n\ge0}F_nZ/F_{n-1}Z$ is a trivial right CDG\+comodule
over $(B,d,h)$ which, as a complex of right $R$\+modules, is
naturally isomorphic to the tensor product
$N^\bu\ot_R\cone(R\to K^\tau_\bu(B,A))$.
 For every $n\in\boZ$, the internal degree~$n$ component of
the complex $\cone(R\to K^\tau_\bu(B,A))$ is a finite acyclic
complex of $R$\+$R$\+bimodules which are finitely generated and
projective (hence flat) as left $R$\+modules.
 As above, we can conclude that $\gr^FZ\simeq
N^\bu\ot_R\cone(R\to K^\tau_\bu(B,A))$ is a coacyclic complex of
right $R$\+modules.
 Hence $\gr^FZ$ is also coacyclic as a right CDG\+comodule
over $(B,d,h)$.
 By Lemma~\ref{cdg-comodule-filtration-coacyclic-lemma}, it follows
that $Z$ is a coacyclic CDG\+comodule over $(B,d,h)$, too.
\end{proof}

\begin{exs} \label{bimodule-resolution-etc-examples}
 (1)~The free right $R$\+module $R$ can be considered as a one-term
complex of right $R$\+modules, concentrated in the cohomological
degree~$0$ and endowed with the zero differential.
 This complex, in turn, can be considered as a trivial right
CDG\+comodule over $(B,d,h)$.
 The functor ${-}\ot_R^{\tau'}\tA={-}\ot_B\Ksp_{e'}(B,\tA)$ takes
this trivial right CDG\+comodule $R$ over $(B,d,h)$ to the free
right $\tA$\+module $\tA$ (viewed as a one-term complex concentrated
in the cohomological degree~$0$).

 (2)~Applying the functor ${-}\ot_R^{\sigma'}C=
\Hom_{\tA^\rop}(\Ksp_{e'}(B,\tA),{-})$ to the free right $\tA$\+module
$\tA$ produces the right CDG\+comodule $\tA\ot_R^{\sigma'}C=
\Hom_{\tA^\rop}(\Ksp_{e'}(B,\tA),\tA)$,
which was mentioned in the above proof and in
Remark~\ref{nonhomogeneous-Koszul-complexes-remark}.
 The left $\tA$\+module structure on $\tA$ induces a left action of
$\tA$ in $\tA\ot_R^{\sigma'}C$, making it a CDG\+bimodule over
$(\tA,0,0)$ and $(B,d,h)$.
 The adjunction map $R\rarrow\tA\ot_R^{\sigma'}C$ is a closed morphism
of right CDG\+comodules over $(B,d,h)$ (or, if one wishes, of
CDG\+bimodules over $(R,0,0)$ and $(B,d,h)$).
 According to the above argument, its cone is a coacyclic right
CDG\+comodule over $(B,d,h)$.

 (3)~Applying the functor ${-}\ot_R^{\tau'}\tA$ to the right
CDG\+comodule $\tA\ot_R^{\sigma'}C$ over $(B,d,h)$ produces a complex
of right $\tA$\+modules $\tA\ot_R^{\sigma'}C\ot_R^{\tau'}\tA$.
 The left action of $\tA$ in $\tA\ot_R^{\sigma'}C$ induces a left
action of $\tA$ in the complex $\tA\ot_R^{\sigma'}C\ot_R^{\tau'}\tA$,
making it a complex of $\tA$\+$\tA$\+bimodules.
 The adjunction map $\tA\ot_R^{\sigma'}C\ot_R^{\tau'}\tA\rarrow\tA$
is a morphism of complexes of $\tA$\+$\tA$\+bimodules.
 According to the above proof, its cone is coacyclic as a complex
of right $R$\+modules (following the details of the argument,
one can see that this cone is, in fact, coacyclic as a complex of
$\tA$\+$R$\+bimodules).

 In particular, the cone is an acyclic complex; so
$\tA\ot_R^{\sigma'}C\ot_R^{\tau'}\tA$ is a bimodule resolution of
the diagonal $\tA$\+$\tA$\+bimodule~$\tA$.
 The terms of this resolution are \emph{not} projective
$\tA\ot_\boZ\tA$\+modules, of course; and $\tA\ot_R\tA$ is,
generally speaking, not even a ring!
 Still one can say that the terms of this resolution
are the $\tA$\+$\tA$\+bimodules $\tA\ot_RC_i\ot_R\tA$ induced from
the $R$\+$R$\+bimodules $C_i=\Hom_{R^\rop}(B^i,R)$ (which are finitely
generated and projective as left $R$\+modules in our assumptions).

 (4)~Notice that the bimodule resolution
$\tA\ot_R^{\sigma'}C\ot_R^{\tau'}\tA\rarrow\tA$ from~(3) is
a quasi-isomorphism of bounded above complexes of projective
left $\tA$\+modules; hence it is a homotopy equivalence of complexes
of left $\tA$\+modules.
 It is also a quasi-isomorphism of bounded above complexes of
weakly $\tA/R$\+flat right $\tA$\+modules in the sense
of~\cite[Section~5]{Pfp} (cf.\ Section~\ref{flat-Koszul-subsecn}).
 In view of~\cite[Lemma~5.3(b)]{Pfp}, it follows that the tensor
product of this resolution with any $R$\+flat left $\tA$\+module is
still an exact complex.

 We refer to Section~\ref{revisited-subsecn} below for a further
discussion of this bimodule resolution.

 (5)~Now assume that $\tA$ is a left augmented left finitely
projective nonhomogeneous Koszul ring over $R$, in the sense of
the definitions in Sections~\ref{augmented-subsecn}
and~\ref{anti-equivalences-koszul-rings-subsecn}.
 Let $\tA^+\subset\tA$ be the augmentation ideal; so $\tA^+$ is
a left ideal in $\tA$ such that $\tA=R\oplus\tA^+$.
 Choose the left $R$\+submodule of strict generators $V'\subset\tV$
as $V'=\tA^+\cap\tV$.
 Then, by Theorem~\ref{augmented-duality-functor-existence-thm} (see
also Corollary~\ref{left-augmented-koszul-duality-anti-equivalence}),
we have $h=0$; so $(B,d,h)=(B,d,0)$ is a DG\+ring.
 Accordingly, the CDG\+bimodule $\tA\ot_R^{\sigma'}C$ from~(2) is
a DG\+bimodule, i.~e., a complex.
 This complex can be called the \emph{nonhomogeneous Koszul complex}
of a left augmented left finitely projective nonhomogeneous Koszul
ring~$\tA$.
 
 The left augmentation of~$\tA$ endows the base ring $R$ with
a structure of left $\tA$\+module provided by the identification
$R\simeq\tA/\tA^+$.
 Taking the tensor product of the bimodule resolution from~(3) with
this left $\tA$\+module, we obtain a morphism of complexes of
left $\tA$\+modules $\tA\ot_R^{\sigma'}C\rarrow R$, which is
a quasi-isomorphism by the argument in~(4).
 Here $\tA\ot_R^{\sigma'}C$ is the above DG\+bimodule over
$(\tA,0)$ and $(B,d)$.

 According to~(2), we also have a natural map in the opposite
direction $R\rarrow\tA\ot_R^{\sigma'}C$, whose cone is coacyclic
(hence acyclic) as a right DG\+bimodule over $(B,d)$.
 The composition $R\rarrow\tA\ot_R^{\sigma'}C\rarrow R$ is
the identity map.
 This provides another way to prove that the map
$\tA\ot_R^{\sigma'}C\rarrow R$ is a quasi-isomorphism.

 Clearly, the terms of the complex $\tA\ot_R^{\sigma'}C$ are finitely
generated projective left $\tA$\+modules in our assumptions.
 Thus the nonhomogeneous Koszul complex $\tA\ot_R^{\sigma'}C$ is
a resolution of the left $\tA$\+module $R$ by finitely generated
projective left $\tA$\+modules.
\end{exs}

\begin{rems} \label{comodule-side-koszul-duality-functors-interpreted}
 (1)~Consider the trivial right CDG\+comodule $R$ over $(B,d,h)$,
as in Example~\ref{bimodule-resolution-etc-examples}(1).
 Since the triangulated equivalence ${-}\ot_R^{\tau'}\tA\:
\sD^\co(\comodr B)\rarrow\sD^\sico_R(\modr\tA)$ takes $R$ to
the free right $\tA$\+module $\tA$, it follows that the graded ring
of endomorphisms $R\rarrow R[i]$ in $\sD^\co(\comodr B)$, \ $i\in\boZ$,
is naturally isomorphic to ring $\tA$ (concentrated in the cohomological
grading $i=0$).

 Let $N$ be a right CDG\+comodule over $(B,d,h)$.
 Then the complex of right $\tA$\+modules $N\ot_R^{\tau'}\tA$
computes the $\tA$\+modules $\Hom_{\sD^\co(\comodr B)}(R,N[i])$.
 Indeed, we have
$$
 \Hom_{\sD^\co(\comodr B)}(R,N[i])\simeq
 \Hom_{\sD^\sico_R(\modr\tA)}(\tA,\>N\ot_R^{\tau'}\tA[i])=
 H^i(N\ot_R^{\tau'}\tA),
$$
because
$$
 \Hom_{\sD^\sico_R(\modr\tA)}(\tA,M^\bu[i])=
 \Hom_{\Hot(\modr\tA)}(\tA,M^\bu[i])=H^i(M^\bu)
$$
for any complex of right $\tA$\+modules~$M^\bu$.
 The latter isomorphism holds since all complexes in
$\Acycl^\sico_R(\modr\tA)$ are acyclic.

 (2)~Now assume that $\tA$ is left augmented over $R$, as in
Example~\ref{bimodule-resolution-etc-examples}(5).
 Then, for any complex of right $\tA$\+modules $M^\bu$,
the CDG\+(co)module $M^\bu\ot_R^{\sigma'}C$ is a DG\+module
over $(B,d,h)=(B,d,0)$, i.~e., a complex.
 We have a natural isomorphism
$M^\bu\ot_R^{\sigma'}C\simeq M^\bu\ot_{\tA}(\tA\ot_R^{\sigma'}C)$.
 Since $\tA\ot_R^{\sigma'}C$ is a projective resolution of the left
$\tA$\+module $R$, it follows that the complex of abelian groups
$M^\bu\ot_R^{\sigma'}C$ computes the derived tensor product
$M^\bu\ot_{\tA}^\boL R$.
 In fact, the DG\+ring $(B,d)$ computes (the opposite ring to) the
graded ring of endomorphisms $R\rarrow R[i]$ in $\sD(\tA\modl)$, \
$i\in\boZ$, and the DG\+module $M^\bu\ot_R^{\sigma'}C$ computes
the homology $\Tor^{\tA}_*(M^\bu,R)=H_*(M^\bu\ot_{\tA}^\boL R)$ as
a graded module over this graded ring of endomorphisms.
\end{rems}

\subsection{Reduced coderived category}
\label{reduced-coderived-category-subsecn}
 We keep the assumptions of
Section~\ref{comodule-side-triangulated-equivalence-subsecn}.
 So $R\subset\tV\subset\tA$ is a left finitely projective
nonhomogeneous Koszul ring and $(B,d,h)$ is the corresponding
right finitely projective Koszul CDG\+ring.
 The following result is a special case of
Theorem~\ref{comodule-side-koszul-duality-theorem}.

\begin{cor} \label{fin-dim-base-koszul-duality-comodule-side}
 Assume additionally that the right homological dimension of
the ring $R$ (that is, the homological dimension of the abelian
category\/ $\modr R$) is finite.
 Then the pair of adjoint triangulated functors\/
${-}\ot_R^{\sigma'}C\:\Hot(\modr\tA)\rarrow\Hot(\comodr B)$ and\/
${-}\ot_R^{\tau'}\tA\:\Hot(\comodr B)\rarrow\Hot(\modr\tA)$ induces
mutually inverse triangulated equivalences
\begin{equation} \label{fin-dim-base-comodule-side}
 \sD(\modr\tA)\simeq\sD^\co(\comodr B)
\end{equation}
between the derived category of right $\tA$\+modules and
the coderived category of right CDG\+comodules over $(B,d,h)$.
\end{cor}

\begin{proof}
 When the right homological dimension of the ring $R$ is finite,
all acyclic complexes of right $R$\+modules are coacyclic
by~\cite[Remark~2.1]{Psemi}.
 So the $\tA/R$\+semi\-coderived category of right $\tA$\+modules
coincides with their derived category,
$\sD^\sico_R(\modr\tA)=\sD(\modr\tA)$.
 This means, specifically, that the right CDG\+comodule
$M^\bu\ot_R^{\sigma'}C$ over $(B,d,h)$ is coacyclic for any
acyclic complex of right $\tA$\+modules~$M^\bu$.
 The rest is explained in
Theorem~\ref{comodule-side-koszul-duality-theorem} and its proof.
\hbadness=1250
\end{proof}

 We would like to obtain a description of the derived category
$\sD(\modr\tA)$ in terms of some kind of exotic derived category
of right CDG\+comodules over $(B,d,h)$ \emph{without} assuming
finiteness of the right homological dimension of~$R$.

 For this purpose, let us consider the full subcategory of
acyclic complexes of right $R$\+modules $\Acycl(\modr R)\subset
\Hot(\modr R)$.
 All complexes of right $R$\+modules, and in particular acyclic
complexes of right $R$\+modules, can be viewed as trivial right
CDG\+comodules over $(B,d,h)$.

 By a kind of abuse of notation, let us denote simply by
$\lan\Acycl(\modr R)\ran_\oplus$ the minimal full triangulated
subcategory in $\sD^\co(\comodr B)$ containing all the acyclic
complexes of right $R$\+modules (viewed as trivial
CDG\+comodules) and closed under infinite direct sums.
 Consider the triangulated Verdier quotient category
$$
 \sD^\co_{R\red}(\comodr B)=
 \sD^\co(\comodr B)/\lan\Acycl(\modr R)\ran_\oplus;
$$
let us call it the \emph{reduced coderived category of
right CDG\+comodules over $(B,d,h)$ relative to~$R$}.

\begin{thm} \label{reduced-koszul-duality-comodule-side-thm}
 The pair of adjoint triangulated functors\/
${-}\ot_R^{\sigma'}C\:\Hot(\modr\tA)\rarrow\Hot(\comodr B)$ and\/
${-}\ot_R^{\tau'}\tA\:\Hot(\comodr B)\rarrow\Hot(\modr\tA)$ induces
mutually inverse triangulated equivalences
$$
 \sD(\modr\tA)\simeq\sD^\co_{R\red}(\comodr B)
$$
between the derived category of right $\tA$\+modules and
the reduced coderived category of right CDG\+comodules
over $(B,d,h)$ relative to~$R$.
\end{thm}

\begin{proof}
 In view of Theorem~\ref{comodule-side-koszul-duality-theorem},
only two things still need to be checked.
 Firstly, we have to show that, for any acyclic complex of right
$\tA$\+modules $M^\bu$, the CDG\+comodule $M^\bu\ot_R^{\sigma'}C$
over $(B,d,h)$ belongs to the triangulated subcategory
$\lan\Acycl(\modr R)\ran_\oplus\subset\sD^\co(\comodr B)$.
 Secondly, it needs to be demonstrated that, for any
CDG\+comodule $N\in\lan\Acycl(\modr R)\ran_\oplus$, the complex
of right $\tA$\+modules $N\ot_R^{\tau'}\tA$ is acyclic.

 Firstly, let $M=M^\bu$ be an acyclic complex of right $\tA$\+modules.
 Consider the increasing filtration $F$ on the graded right $R$\+module
$M\ot_RC$ induced by the increasing filtration $F$ on
the $R$\+$R$\+bimodule $C$, as in the proof of
Theorem~\ref{comodule-side-koszul-duality-theorem}.
 Then $F$ is a filtration of the CDG\+comodule $M^\bu\ot_R^{\sigma'}C$
over $(B,d,h)$ by CDG\+subcomodules.
 The successive quotient CDG\+comodules are trivial CDG\+comodules
which, viewed as complexes of $R$\+modules, can be computed as
the tensor products $M^\bu\ot_R(F_nC/F_{n-1}C)$.
 Here $F_nC/F_{n-1}C$ is a graded $R$\+$R$\+bimodule, finitely generated
and projective as a left $R$\+module, and concentrated in
the single cohomological degree~$-n$.

 It follows that the complexes of right $R$\+modules
$M^\bu\ot_R(F_nC/F_{n-1}C)$ are acyclic.
 Now the iterated extension and telescope sequence argument from
the proofs of Lemmas~\ref{complex-filtration-coacyclic-lemma}
and~\ref{cdg-comodule-filtration-coacyclic-lemma} shows that
the right CDG\+comodule $M^\bu\ot_R^{\sigma'}C$ over $(B,d,h)$ belongs
to the triangulated subcategory $\lan\Acycl(\modr R)\ran_\oplus$
in $\sD^\co(\comodr B)$.

 Secondly, since the the functor ${-}\ot_R^{\tau'}\tA$ is triangulated
and preserves infinite direct sums, it suffices to check that
for any acyclic complex of right $R$\+modules $N^\bu$, viewed as
a trivial right CDG\+comodule over $(B,d,h)$, the complex of right
$\tA$\+modules $N^\bu\ot_R^{\tau'}\tA$ is acyclic.
 As $N^\bu$ is a trivial CDG\+comodule, the complex of right
$\tA$\+modules $N^\bu\ot_R^{\tau'}\tA=N^\bu\ot_R\tA$ is simply obtained
from the complex of right $R$\+modules $N^\bu$ by applying the tensor
product functor ${-}\ot_R\tA$ termwise.
 Since $\tA$ is a projective (hence flat) left $R$\+module, it follows
that the complex $N^\bu\ot_R\tA$ is acyclic.
\end{proof}

\Section{Relative Nonhomogeneous Derived Koszul Duality: \\
the Contramodule Side}
\label{contramodule-side-secn}

\subsection{Semicontraderived category of modules}
\label{semicontraderived-category-subsecn}
 Let $R$ be an associative ring.
 The following definition of the \emph{contraderived category of left
$R$\+modules} is dual to the definition of the coderived category in
Section~\ref{semicoderived-category-subsecn}.
 It goes back to~\cite[Section~4.1]{Psemi}
and~\cite[Section~3.3]{Pkoszul} (see~\cite{Jor}
and~\cite[Proposition~1.3.6(1)]{Bec} for an alternative approach).

 Similarly to Section~\ref{semicoderived-category-subsecn}, we consider
short exact sequences of left $R$\+modules $0\rarrow L^\bu\rarrow
M^\bu\rarrow N^\bu\rarrow0$ and their totalizations
$\Tot(L^\bu\to M^\bu\to N^\bu)$.
 A complex of left $R$\+modules is said to be \emph{contraacyclic} if it
belongs to the minimal full triangulated subcategory of the homotopy
category of complexes of left $R$\+modules $\Hot(R\modl)$ containing
all the totalizations of short exact sequences of left $R$\+modules
and closed under infinite products.

 The full subcategory of contraacyclic complexes of left $R$\+modules
is denoted by $\Acycl^\ctr(R\modl)\subset\Hot(R\modl)$.
 The \emph{contraderived category of left $R$\+modules} is
the triangulated Verdier quotient category of the homotopy category
$\Hot(R\modl)$ by the triangulated subcategory of contraacyclic
complexes,
$$
 \sD^\ctr(R\modl)=\Hot(R\modl)/\Acycl^\ctr(R\modl).
$$
 Clearly, any contraacyclic complex of $R$\+modules is acyclic,
$\Acycl^\ctr(R\modl)\subset\Acycl(R\modl)$; but the converse does not
hold in general~\cite[Examples~3.3]{Pkoszul}.
 So the derived category $\sD(R\modl)$ is a quotient category of
$\sD^\ctr(R\modl)$.
 All acyclic complexes of left $R$\+modules are contraacyclic when
the left homological dimension of the ring $R$ is
finite~\cite[Remark~2.1]{Psemi}.

 The following lemma is dual-analogous to and very slightly more
complicated than Lemma~\ref{complex-filtration-coacyclic-lemma}.

\begin{lem} \label{complex-filtration-contraacyclic-lemma}
 Let\/ $P^\bu=F^0P^\bu\supset F^1 P^\bu\supset F^2P^\bu\supset\dotsb$
be a complex of $R$\+modules that is separated and complete with
respect to a decreasing filtration by subcomplexes of $R$\+modules,
that is $P^\bu=\varprojlim_{n\ge0}P^\bu/F^{n+1}P^\bu$.
 Assume that, for every $n\ge0$, the complex of $R$\+modules
$F^nP^\bu/F^{n+1}P^\bu$ is contraacyclic.
 Then the complex of $R$\+modules $P^\bu$ is contraacyclic.
\end{lem}

\begin{proof}
 First one proves by induction in~$n$ that the complex of $R$\+modules
$P^\bu/F^{n+1}P^\bu$ is contraacyclic for every $n\ge0$, using
the fact that the totalizations of short exact sequences
$0\rarrow F^nP^\bu/F^{n+1}P^\bu\rarrow P^\bu/F^{n+1}P^\bu\rarrow
P^\bu/F^nP^\bu\rarrow0$ are contraacyclic.
 Then one concludes that the complex of $R$\+modules
$\prod_{n=0}^\infty P^\bu/F^{n+1}P^\bu$ is contraacyclic.
 Finally, one observes that the telescope sequence $0\rarrow
P^\bu\rarrow\prod_{n=0}^\infty P^\bu/F^{n+1}P^\bu\rarrow
\prod_{n=0}^\infty P^\bu/F^{n+1}P^\bu\rarrow0$ of the projective
system $(P^\bu/F^{n+1}P^\bu)_{n\ge0}$ is exact, since the transition
maps $P^\bu/F^{n+1}P^\bu\rarrow P^\bu/F^nP^\bu$ are surjective
(so the derived projective limit of this projective system vanishes,
$\varprojlim_{n\ge0}^1 P^\bu/F^{n+1}P^\bu=0$).
 Hence the totalization of the telescope sequence is a contraacyclic
complex of $R$\+modules, and contraacyclicity of the complex $P^\bu$
follows.
\end{proof}

 Now let $R\rarrow A$ be a homomorphism of associative rings.
 The following definition of the \emph{$A/R$\+semicontraderived category
of left $A$\+modules} can be found in~\cite[Section~5]{Pfp}.
 The semicontraderived category
$$
 \sD^\sictr_R(A\modl)=\Hot(A\modl)/\Acycl^\sictr_R(A\modl)
$$
is the triangulated Verdier quotient category of the homotopy category
of left $A$\+modules by its full triangulated subcategory
$\Acycl^\sictr_R(A\modl)$ consisting of all the complexes of left
$A$\+modules that are \emph{contraacyclic as complexes of
left $R$\+modules}.

 Clearly, any contraacyclic complex of $A$\+modules is contraacyclic
as a complex of $R$\+modules; but the converse does not need to be true.
 So the semicontraderived category is intermediate between the derived
and the contraderived category of $A$\+modules: there are natural
triangulated Verdier quotient functors $\sD^\ctr(A\modl)\rarrow
\sD^\sictr_R(A\modl)\rarrow\sD(A\modl)$.
 (See the paper~\cite{Pps} for a more general discussion of such
intermediate triangulated quotient categories, called
\emph{pseudo-contraderived categories} in~\cite{Pps}.)
 The semicontraderived category $\sD^\sictr_R(A\modl)$ can be thought
of as a mixture of the contraderived category ``along the variables
from~$R$'' and the conventional derived category ``in the direction
of $A$ relative to~$R$''.

 For example, taking $R=\boZ$, one obtains the conventional derived
category of $A$\+modules
$$
 \sD^\sictr_\boZ(A\modl)=\sD(A\modl)
$$
(because all the acyclic complexes of abelian groups are contraacyclic),
while taking $R=A$ one obtains the contraderived category of
$A$\+modules
$$
 \sD^\sictr_A(A\modl)=\sD^\ctr(A\modl).
$$

\subsection{Contraderivations} \label{contraderivations-subsecn}
 Let $B=\bigoplus_{n=0}^\infty B^n$ be a nonnegatively graded ring
and $K$ be a ring endowed with a ring homomorphism $K\rarrow B^0$.

 Let $k\in\boZ$ be a fixed integer and $d\:B\rarrow B$ be a fixed
derivation of degree~$k$ (even or odd, depending on the parity of~$k$).
 This means that we have additive maps $d_n\:B^n\rarrow B^{n+k}$
given for all $n\ge0$ and satisfying the equation
$d(bc)=d(b)c+(-1)^{k|b|}bd(c)$ for all
$b\in B^{|b|}$ and $c\in B^{|c|}$.
 Assume that composition $K\rarrow B^0\rarrow B^k$ of the maps
$K\rarrow B^0$ and $d_0\:B^0\rarrow B^k$ vanishes (e.~g., one
can always take $K=\boZ$).
 So $d\:B\rarrow B$ is a left and right $K$\+linear map.

 Let $M=\bigoplus_{n\in\boZ}M^n$ be a graded left $B$\+module.
 A homogeneous additive map $d_M\:M\rarrow M$ of degree~$k$ (that is,
$d_M(M^n)\subset M^{n+k}$ for all $n\in\boZ$) is said to be
an (\emph{even} or \emph{odd}) \emph{derivation of $M$ compatible with
the derivation~$d$ on~$B$} if the equation
$d_M(bx)=d(b)x+(-1)^{k|b|}bd_M(x)$ holds
for all $b\in B^{|b|}$ and $x\in M^{|x|}$.

 Let $C$ and $D$ be graded right $K$\+modules, $P$ and $Q$ be graded
left $K$\+modules, $g\:C\rarrow D$ be a homogenenous $K$\+linear map
of degree~$l$, and $f\:P\rarrow Q$ be a homogeneous $K$\+linear map
of degree~$m$.
 Consider the graded abelian groups $C\ot_K^\Pi P$ and $D\ot_K^\Pi Q$
(in the notation of Section~\ref{graded-ext-comparison-subsecn}),
and construct the homogeneous map
$$
 g\ot^\Pi f\:C\ot_K^\Pi P\lrarrow D\ot_K^\Pi Q
$$
of degree $l+m$ by taking the infinite products along the diagonals of
the bigrading for the bihomogeneous map
$g\ot f\:C\ot_K P\rarrow D\ot_K Q$ with the components
$g\ot f\:C^i\ot_K P^j\rarrow D^{i+l}\ot_K Q^{j+m}$ given by the rule
$(g\ot f)(c\ot p)=(-1)^{mi}g(c)\ot f(p)$ for all
$c\in C^i$ and $p\in P^j$.

 Let $P=\bigoplus_{n\in\boZ}P^n$ be a graded left $B$\+contramodule
(in the sense of the definition in
Section~\ref{graded-contramodules-subsecn}), and let
$d_P\:P\rarrow P$ be a homogeneous $K$\+linear map of degree~$k$.
 We will say that $d_P$~is an (\emph{even} or \emph{odd})
\emph{contraderivation of $P$ compatible with the derivation~$d$
on~$B$} if the contraaction map $\pi_P\:\boM_K^\gr(P)=
B\ot_K^\Pi P\rarrow P$ forms a commutative square diagram with
the maps
$$
 d\ot^\Pi{\id_P}+{\id_B}\ot^\Pi d_P\:B\ot_K^\Pi P\lrarrow B\ot_K^\Pi P
$$
and $d_P\:P\rarrow P$, that is
$$
 d_P\circ\pi_P=\pi_P\circ(d\ot^\Pi{\id_P}+{\id_B}\ot^\Pi d_P).
$$
 Obviously, any contraderivation of a graded left $B$\+contramodule
$P$ compatible with the given derivation~$d$ on $B$ is a derivation
of the underlying graded left $B$\+module of $P$ compatible with
the same derivation~$d$ on~$B$.

\begin{prop} \label{contra-derivations-prop}
 Assume that the forgetful functor $B\contra_\sgr\rarrow B\modl_\sgr$
from the category of graded left $B$\+contramodules to the category
of graded left $B$\+modules is fully faithful.
 Let $Q$ be a graded left $B$\+contramodule.
 Then any derivation of the underlying graded left $B$\+module of $Q$
(compatible with the given derivation~$d$ on~$B$) is a contraderivation
of $Q$ (compatible with the same derivation~$d$ on~$B$).
\end{prop}

\begin{proof}
 The argument is based on a buildup of auxiliary definitions.
 Given a graded left $B$\+module $M$ and an integer~$n$, the graded
left $B$\+module $M[n]$ has the grading components $M[n]^i=M^{n+i}$;
the action of $B$ in $M[n]$ is defined in terms of the left action
of $B$ in $M$ with the sign rule written down in
Section~\ref{cdg-modules-subsecn}.
 Given a graded left $B$\+contramodule $P$, the graded left
$B$\+contramodule $P[n]$ has the grading components $P[n]^i=P^{n+i}$.
 We leave it to the reader to spell out the construction of the left
contraaction of $B$ in $P[n]$ involving the same sign rule.
 
 Let $f\:M\rarrow N$ be a morphism of graded left $B$\+modules,
and let $d_{N,M}\:M\rarrow N$ be a homogeneous additive map of
degree~$k$.
 We will say that $d_{N,M}$ is a \emph{relative derivation compatible
with the derivation~$d$ on $B$ and the morphism~$f$} if the equation
$d_{N,M}(bx)=d(b)f(x)+(-1)^{k|b|}bd_{N,M}(x)$ holds in $N^{|b|+|x|+k}$
for all $b\in B^{|b|}$ and $x\in M^{|x|}$.
 A derivation $d_M\:M\rarrow M$ compatible with the given derivation~$d$
on $B$ is the same thing as a relative derivation $d_{M,M}\:M\rarrow M$
compatible with~$d$ and the identity morphism $\id_M\:M\rarrow M$.

 For any two relative derivations $d'_{N,M}$ and $d''_{N,M}\:M\rarrow N$
compatible with the same derivation~$d$ on $B$ and the same
morphism $f\:M\rarrow N$, the difference $d''_{N,M}-d'_{N,M}$ is
a morphism of graded left $B$\+modules $M\rarrow N[k]$.
 Conversely, the sum $d'_{N,M}=d_{N,M}+g$ of a relative derivation
$d_{N,M}\:M\rarrow N$ compatible with $d$ and~$f$ and a morphism of
graded left $B$\+modules $g\:M\rarrow N[k]$ is a relative derivation
compatible with $d$ and~$f$.

 Let $f\:M\rarrow N$ and $g\:L\rarrow M$ be two morphisms of graded
left $B$\+modules, and let $d_{N,M}\:M\rarrow N$ be a relative
derivation compatible with the derivation~$d$ on $B$ and
the morphism~$f$.
 Then $d_{N,M}\circ g\:L\rarrow N$ is a relative derivation compatible
with the derivation~$d$ on $B$ and the morphism $f\circ g\:L\rarrow N$.

 Let $f\:P\rarrow Q$ be a morphism of graded left $B$\+contramodules,
and let $d_{Q,P}\:P\rarrow Q$ be a homogeneous $K$\+linear map of
degree~$k$.
 Consider the square diagram formed by the map
$$
 (d\ot^\Pi f+{\id_B}\ot^\Pi d_{Q,P})\:B\ot^\Pi P\lrarrow B\ot^\Pi Q
$$
together with the map $d_{Q,P}\:P\rarrow Q$ and the contraaction maps
$\pi_P\:B\ot^\Pi P\rarrow P$ and $\pi_Q\:B\ot^\Pi Q\rarrow Q$.
 We will say that $d_{Q,P}$ is a \emph{relative contraderivation
compatible with the derivation~$d$ on $B$ and the morphism~$f$} if
this square diagram is commutative, that is
$$
 d_{Q,P}\circ\pi_P=\pi_Q\circ(d\ot^\Pi f+{\id_B}\ot^\Pi d_{Q,P}).
$$
 A contraderivation $d_P\:P\rarrow P$ compatible with the given
derivation~$d$ on $B$ is the same thing as a relative contraderivation
$d_{P,P}\:P\rarrow P$ compatible with~$d$ and the identity morphism
$\id_P\:P\rarrow P$.

 For any two relative contraderivations $d'_{Q,P}$ and
$d''_{Q,P}\:P\rarrow Q$ compatible with the same derivation~$d$ on $B$
and the same morphism $f\:P\rarrow Q$, the difference
$d''_{Q,P}-d'_{Q,P}$ is a morphism of graded left $B$\+contramodules
$P\rarrow Q[k]$.
 Conversely, the sum $d'_{Q,P}=d_{Q,P}+g$ of a relative contraderivation
$d_{Q,P}\:P\rarrow Q$ compatible with $d$ and~$f$ and a morphism of
graded left $B$\+contramodules $g\:P\rarrow Q[k]$ is a relative
contraderivation compatible with $d$ and~$f$.

 Let $f\:P\rarrow Q$ and $g\:S\rarrow P$ be two morphisms of graded
left $B$\+contramodules, and let $d_{Q,P}\:P\rarrow Q$ be a relative
contraderivation compatible with the derivation~$d$ on $B$
and the morphism~$f$.
 Then $d_{Q,P}\circ g\:S\rarrow Q$ is a relative contraderivation
compatible with the derivation~$d$ on $B$ and the morphism
$f\circ g\:S\rarrow Q$.

 Let $L$ be a graded left $K$\+module and $M=B\ot_KL$ be the induced
graded left $B$\+module (so there is a natural graded left
$K$\+module morphism $L\rarrow M$).
 Let $f\:M\rarrow N$ a morphism of graded left $B$\+modules.
 Then a relative derivation $d_{N,M}\:M\rarrow N$ compatible with
the given derivation~$d$ on $B$ and the morphism~$f$ is uniquely
determined by its composition with the map $L\rarrow M$.

 The latter composition $L\rarrow M\overset{d_{N,M}}\lrarrow N$ can be
an arbitrary homogeneous $K$\+linear map $t\:L\rarrow N$ of degree~$k$.
 Given such a map~$t$, and denoting by $f'\:L\rarrow N$
the composition $L\rarrow M\overset f\rarrow N$ (which is
a homogeneous $K$\+linear map of degree~$0$), the related relative
derivation $d_{N,M}\:M\rarrow N$ can be recovered by the formula
$d_{N,M}(b\ot l)=d(b)f'(l)+(-1)^{k|b|}bt(l)$ for all
$b\in B^{|b|}$ and $l\in K^{|l|}$.

 In particular, let $X=\coprod_{i\in\boZ}X_i$ be a graded set and
$K[X]$ be the free graded left $K$\+module with the components
$K[X]^i=K^{(X_i)}$ as in the the proof of
Proposition~\ref{graded-B-contramodules-proposition}.
 Let $M=B[X]=B\ot_KK[X]$ be the free graded left $B$\+module spanned
by~$X$, and let $f\:M\rarrow N$ be a morphism of graded left
$B$\+modules.
 Then a relative derivation $d_{N,M}\:M\rarrow N$ compatible with
the given derivation~$d$ on $B$ and the morphism~$f$ is uniquely
determined by its restriction to the subset of free generators
$X\subset M$, and this restriction can be an arbitrary homogeneous
map $X\rarrow N$ of degree~$k$ (i.~e., a collection of maps
$X_i\rarrow N^{i+k}$, \ $i\in\boZ$).

 Consider the graded left $B$\+contramodule $P=\boM_K^\gr(L)=
B\ot_K^\Pi L$ freely generated by~$L$.
 Then, once again, there is a natural graded left $K$\+module
morphism $L\rarrow P$.
 Let $f\:P\rarrow Q$ be a morphism of graded left $B$\+contramodules.
 Then a relative contraderivation $d_{Q,P}\:P\rarrow Q$ compatible
with the given derivation~$d$ on $B$ and the morphism~$f$ is
uniquely determined by its composition with the map $L\rarrow P$.

 The latter composition $L\rarrow P\overset{d_{Q,P}}\lrarrow Q$ can be
an arbitrary homogeneous $K$\+linear map $t\:L\rarrow Q$ of degree~$k$.
 Given such a map~$t$, and denoting by $f'\:L\rarrow Q$
the composition $L\rarrow P\overset f\rarrow Q$ (which is a homogeneous
$K$\+linear map of degree~$0$), the related relative contraderivation
$d_{Q,P}\:P\rarrow Q$ can be recovered by the formula
$d_{Q,P}=\pi_Q\circ(d\ot^\Pi f'+{\id_B}\ot^\Pi t)$,
$$
 P\lrarrow B\ot_K^\Pi Q\lrarrow Q.
$$

 In particular, let $P=\boM_K^\gr(K[X])=B\ot_K^\Pi K[X]$ be
the free graded left $B$\+contramodule spanned by~$X$, and let
$f\:P\rarrow Q$ be a morphism of graded left $B$\+contramodules.
 Then a relative contraderivation $d_{Q,P}\:P\rarrow Q$ compatible
with the given derivation~$d$ on $B$ and the morphism~$f$ is
uniquely determined by its restriction to the subset of free generators
$X\subset P$, and this restriction can be an arbitrary homogeneous
map $X\rarrow Q$ of degree~$k$.

 Now that we are done with the auxiliary material, we can prove
the proposition.
 Let $d_Q\:Q\rarrow Q$ be a derivation of the underlying graded left
$B$\+module of $Q$ compatible with the given derivation~$d$ on~$B$.
 Choose a free graded left $B$\+contramodule $P$ (spanned by some
graded set $X$) together with a surjective graded left
$B$\+contramodule morphism $f\:P\rarrow Q$.
 Then the composition $d_Q\circ f\:P\rarrow Q$ is a relative
derivation of the underlying graded $B$\+modules of $P$ and $Q$
compatible with the derivation~$d$ on $B$ and
the morphism $f\:P\rarrow Q$.

 Consider the restriction of the map $d_Q\circ f\:P\rarrow Q$ to
the subset of free generators $X\subset P$, and extend the resulting
map $X\rarrow Q$ to a relative contraderivation $d_{Q,P}\:P\rarrow Q$
of the graded $B$\+contramodules $P$ and $Q$ compatible with
the derivation~$d$ on $B$ and the contramodule morphism
$f\:P\rarrow Q$.
 Then the difference $d_Q\circ f-d_{Q,P}$ is a morphism of underlying
graded $B$\+modules $P\rarrow Q[k]$ of the graded left
$B$\+contramodules $P$ and~$Q[k]$.

 By the assumption of the proposition, any morphism between
the underlying graded left $B$\+modules of two graded left
$B$\+contramodules is a graded left $B$\+contramodule morphism.
 Hence $d_Q\circ f-d_{Q,P}$ must be, in fact, a graded left
$B$\+contramodule morphism.
 Since $d_{Q,P}$ is a relative contraderivation compatible with
the derivation~$d$ on $B$ and the morphism~$f$, it follows that
$d_Q\circ f$ is a relative contraderivation compatible with $d$
and~$f$, too.

 Finally, since the map~$f$ is surjective, we can conclude that $d_Q$~is
a contraderivation of $Q$ compatible with the derivation~$d$ on~$B$.
\end{proof}

\subsection{Contraderived category of CDG-contramodules}
\label{contraderived-cdg-contramodules-subsecn}
 Let $(B,d,h)$ be a nonnegatively graded CDG\+ring, as defined
in Section~\ref{curved-dg-rings-subsecn}, and let $(\hB,\d)$ be
the related quasi-differential graded ring constructed in
Theorem~\ref{cdg-qdg-equivalence}.
 So $h\in B^2$ is a curvature element for an odd derivation
$d\:B\rarrow B$ of degree~$1$.

 A \emph{left CDG\+contramodule over $(B,d,h)$} can be simply defined
as a graded left $\hB$\+contramodule.
 Equivalently, a left CDG\+contramodule $(P,d_P)$ over $(B,d,h)$ is
a graded left $B$\+contramodule $P=\bigoplus_{n\in\boZ}P^n$ endowed
with a sequence of additive maps $d_{P,n}\:P^n\rarrow P^{n+1}$
satisfying the following conditions:
\begin{enumerate}
\renewcommand{\theenumi}{\roman{enumi}}
\item $d_P$~is an odd contraderivation of the graded left
$B$\+contramodule $P$ compatible with the odd derivation~$d$ on $B$
(in the sense of Section~\ref{contraderivations-subsecn});
\item $d_P^2(x)=hx$ for all $x\in P$.
\end{enumerate}

 According to Proposition~\ref{contra-derivations-prop}, when
the forgetful functor $B\contra_\sgr\rarrow B\modl_\sgr$ is fully
faithful, condition~(i) is equivalent to the seemingly weaker
condition that $d_P$~is an odd derivation of the underlying graded
left $B$\+module of $P$ compatible with the odd derivation~$d$ on~$B$.
 In particular, by
Theorem~\ref{graded-contramodules-fully-faithful-thm}, this holds
whenever the augmentation ideal $B^{\ge1}=\bigoplus_{n=1}^\infty B^n$
is finitely generated as a right ideal in~$B$.
 This includes all $2$\+right finitely projective quadratic graded
rings~$B$ (see Section~\ref{ungraded-comodules-subsecn}).

 Given two graded left $B$\+contramodules $P$ and $Q$, we denote by
$\Hom^{B,n}(P,Q)$ the abelian group of all graded left
$B$\+contramodule homomorphisms $P\rarrow Q[n]$ (where a sign rule
is needed in the definition of the graded left $B$\+contramodule
structure on the grading shift $Q[n]$ of a graded left $B$\+contramodule
$Q$; see the proof of Proposition~\ref{contra-derivations-prop} for
a brief discussion).
 For any two left CDG\+contramodules $P$ and $Q$ over $(B,d,h)$,
there is a natural differential~$d$ on the graded abelian group
$\Hom^B(P,Q)$.
 In fact, $\Hom^B(P,Q)$ is a subcomplex in the complex $\Hom_B(P,Q)$
of morphisms between the underlying graded left $B$\+modules of $P$
and $Q$ (which was constructed in Section~\ref{cdg-modules-subsecn}).

 Of course, one has $\Hom^B(P,Q)=\Hom_B(P,Q)$ for all left
CDG\+contramodules $P$ and $Q$ over $(B,d,h)$ whenever the forgetful
functor $B\contra_\sgr\rarrow B\modl_\sgr$ is fully faithful.
 But generally speaking, the conditions that $d_P$ and $d_Q$ are
contraderivations (rather than just derivations of the underlying
graded left $B$\+modules of $P$ and $Q$) is needed in order to check
that $d(f)\in\Hom^{B,n+1}(P,Q)\subset\Hom_B^{n+1}(P,Q)$ for every
$f\in\Hom^{B,n}(P,Q)\subset\Hom_B^n(P,Q)$.

 The above construction produces the \emph{DG\+category of left
CDG\+contramodules} $\DG(B\contra)$ \emph{over~$B=(B,d,h)$}.
 The DG\+category of left CDG\+contramodules comes endowed with
the forgetful DG\+functor $\DG(B\contra)\rarrow\DG(B\modl)$.
 All the shifts, twists, and infinite products (in the sense of
the discussion in~\cite[Section~1.2]{Pkoszul}) exist in
the DG\+category $\DG(B\contra)$; in particular, all the cones
of closed morphisms exist.
 Passing to the degree-zero cohomology groups of the complexes of
morphisms, we obtain the triangulated \emph{homotopy category of
CDG\+contramodules} $\Hot(B\contra)=H^0\DG(B\contra)$
\emph{over $(B,d,h)$} and the triangulated forgetful functor
$\Hot(B\contra)\rarrow\Hot(B\modl)$.

 Next, as in Section~\ref{coderived-cdg-comodules-subsecn}, we consider
short exact sequences $0\rarrow P\overset i\rarrow Q\overset p\rarrow S
\rarrow 0$ of left CDG\+contramodules over $(B,d,h)$.
 Here $i$ and~$p$ are closed morphisms of degree zero in
the DG\+category $\DG(B\contra)$ and the short sequence $0\rarrow P
\rarrow Q\rarrow S\rarrow 0$ is exact in the abelian category
$B\contra_\sgr$.
 Using the construction of the cone of a closed morphism (or the shift,
finite direct sum, and twist) in the DG\+category $\DG(B\contra)$, one
can produce the totalization (or the \emph{total CDG\+contramodule})
$\Tot(P\to Q\to S)$ of the finite complex of CDG\+contramodules
$P\rarrow Q\rarrow S$.

 A left CDG\+contramodule over $B$ is said to be \emph{contraacyclic}
if it belongs to the minimal full triangulated subcategory of
$\Hot(B\contra)$ containing all the totalizations of short exact
sequences of left CDG\+contramodules over $B$ and closed under
infinite products.
 The full triangulated subcategory of contraacyclic CDG\+contramodules
is denoted by $\Acycl^\ctr(B\contra)\subset\Hot(B\contra)$.
 The \emph{contraderived category of left CDG\+contramodules over~$B$}
is the triangulated Verdier quotient category of the homotopy category
$\Hot(B\contra)$ by its triangulated subcategory of contraacyclic
CDG\+contramodules,
$$
 \sD^\ctr(B\contra)=\Hot(B\contra)/\Acycl^\ctr(B\contra).
$$

 Similarly to the discussion in
Section~\ref{coderived-cdg-comodules-subsecn}, CDG\+contramodules,
generally speaking, are \emph{not complexes} (if $h\ne0$ in~$B$)
and their cohomology (contra)modules are \emph{undefined}.
 So the conventional construction of the derived category is \emph{not}
applicable to curved DG\+contramodules.
 There is only the contraderived category and its variations (known
as derived categories of the second kind).

\begin{lem} \label{cdg-contramodule-filtration-contraacyclic-lemma}
 Let $B=(B,d,h)$ be a nonnegatively graded CDG\+ring and $Q=(Q,d_Q)$
be a left CDG\+contramodule over~$B$.
 Let $Q=F^0Q\supset F^1Q\supset F^2Q\supset\dotsb$ be a decreasing
filtration of $Q$ by CDG\+subcontramodules $F^nQ\subset Q$ such that
$Q$ is separated and complete with respect to this filtration, that is
$Q=\varprojlim_{n\ge0} Q/F^{n+1}Q$.
 Assume that, for every $n\ge0$, the CDG\+contramodule $F^nQ/F^{n+1}Q$
is contraacyclic over~$B$.
 Then the CDG\+contramodule $Q$ is contraacyclic, too.
\end{lem}

\begin{proof}
 This is a generalization of
Lemma~\ref{complex-filtration-contraacyclic-lemma}, provable by
the same argument.
 Let us only point out that the forgetful functor $\DG(B\contra)
\rarrow\DG(B\modl)$ preserves infinite products
(cf.\ Proposition~\ref{graded-B-contramodules-proposition}); so
the projective limits of CDG\+contramodules over $B$ agree with
those of CDG\+modules over $B$ (or just of graded abelian groups).
 Hence the derived countable projective limit vanishes on sequences
of surjective morphisms of CDG\+contramodules, which makes
our argument work.
\end{proof}

 There is an important particular case when the nonnegatively graded
ring $B$ has only finitely many grading components, that is $B^n=0$
for $n\gg0$.
 In this case, the forgetful functors from the categories of (ungraded
or graded) $B$\+contramodules to the categories of (ungraded or
graded) $B$\+modules are equivalences of categories, so in particular
$B\contra_\sgr\simeq B\modl_\sgr$.
 Consequently, the forgetful DG\+functor $\DG(B\contra)\rarrow
\DG(B\modl)$ is an equivalence of DG\+categories.
 In this case, we will write simply $\sD^\ctr(B\modl)$ instead of
$\sD^\ctr(B\contra)$.

\subsection{Koszul duality functors for modules and contramodules}
\label{koszul-duality-functors-contramodule-side-subsecn}
 Let $R\subset\tV\subset\tA$ be a $3$\+left finitely projective
weak nonhomogeneous quadratic ring, and let $B=(B,d,h)$ be
the nonhomogeneous quadratic dual CDG\+ring, as constructed in
Proposition~\ref{nonhomogeneous-dual-cdg-ring}.
 As in Section~\ref{koszul-duality-functors-comodule-side-subsecn},
the dual nonhomogeneous Koszul CDG\+module $\Ksp(B,\tA)=
\Ksp_{e'}(B,\tA)=(B\ot_R\tA,d_{e'})$ constructed in
Section~\ref{nonhomogeneous-koszul-cdg-module-subsecn} plays a key role.

 According to Section~\ref{cdg-modules-subsecn}, the CDG\+bimodule
$\Ksp(B,\tA)$ over the DG\+rings $B=(B,d,h)$ and $\tA=(\tA,0,0)$
induces a pair of adjoint DG\+functors between the DG\+cat\-egory
$\DG(\tA\modl)$ of complexes of left $\tA$\+modules and
the DG\+category $\DG(B\modl)$ of left CDG\+modules over $(B,d,h)$.
 Here the right adjoint DG\+functor
$$
 \Hom_B(\Ksp(B,\tA),{-})\:\DG(B\modl)\lrarrow\DG(\tA\modl)
$$
takes a left CDG\+module $Q=(Q,d_Q)$ over $(B,d,h)$ to the complex
of left $\tA$\+modules whose underlying graded $\tA$\+module is
$$
 \Hom_B(\Ksp(B,\tA),Q)=\Hom_B(B\ot_R\tA,\>Q)=\Hom_R(\tA,Q).
$$
 The left adjoint DG\+functor
$$
 \Ksp(B,\tA)\ot_{\tA}{-}\,\:\DG(\tA\modl)\lrarrow\DG(B\modl)
$$
takes a complex of left $\tA$\+modules $P=P^\bu$ to the left
CDG\+module over $(B,d,h)$ whose underlying graded $B$\+module is
$$
 \Ksp(B,\tA)\ot_{\tA}P=(B\ot_R\tA)\ot_{\tA}P=B\ot_RP.
$$

 Our aim is to replace the DG\+category of CDG\+modules $\DG(B\modl)$
with the DG\+category of CDG\+contramodules $\DG(B\contra)$ in this
adjoint pair.
 In fact, the augmentation ideal $B^{\ge1}$ is finitely generated as
a right ideal in $B$ for any $2$\+right finitely projective quadratic
graded ring~$B$ (according to the discussion in
Section~\ref{ungraded-comodules-subsecn}).
 Therefore, the forgetful DG\+functor $\DG(B\contra)\rarrow\DG(B\modl)$
is fully faithful in our assumptions (see the discussion in
Section~\ref{contraderived-cdg-contramodules-subsecn}).
 So the situation is formally somewhat similar to that in
Section~\ref{koszul-duality-functors-comodule-side-subsecn}.

\begin{lem} \label{delta-tensor-tensor-pi}
 Let $B=\bigoplus_{n=0}^\infty B_n$ be a nonnegatively graded ring
with the degree-zero component $R=B_0$.
 Then the forgetful functor (between the abelian categories)
$B\contra_\sgr\rarrow B\modl_\sgr$ has a left adjoint functor
$\Delta_B\:B\modl_\sgr\rarrow B\contra_\sgr$.
 Furthermore, for any graded left $R$\+module $L$ one has
$$
 \Delta_B(B\ot_RL)=B\ot_R^\Pi L.
$$
\end{lem}

\begin{proof}
 For any graded left $B$\+contramodule $Q$, the group of morphisms
$B\ot_RL\rarrow Q$ in $B\modl_\sgr$ is isomorphic to the group of
morphisms $L\rarrow Q$ in $R\modl_\sgr$.
 The latter group is isomorphic to the group of morphisms
$B\ot_R^\Pi L=\boM_R^\gr(L)\rarrow Q$ in $B\contra_\sgr$
(see the discussion in Section~\ref{monads-subsecn}
and Proposition~\ref{graded-B-contramodules-proposition}).
 This shows that the functor $\Delta_B$ is defined on the full
subcategory in $B\modl_\sgr$ formed by the graded left $B$\+modules
$B\ot_RL$ induced from graded left $R$\+modules~$L$,
and proves the second assertion of the lemma.
 To prove the first assertion, present an arbitrary graded left
$B$\+module $N$ as the cokernel of a morphism
$f\:B\ot_RL\rarrow B\ot_RM$ of graded left $B$\+modules induced from
graded left $R$\+modules $L$ and $M$ (e.~g., one can always take $L$
and $M$ to be free graded left $R$\+modules).
 Then the graded left $B$\+contramodule $\Delta_B(N)$ can be
computed as the cokernel of the related morphism of graded
left $B$\+contramodules $\Delta_B(f)\:B\ot_R^\Pi L
\rarrow B\ot_R^\Pi M$ (it is helpful to keep in mind that the functor
$\Delta_B$, being a left adjoint, has to preserve cokernels).
\end{proof}

\begin{lem} \label{delta-cdg-lemma}
 Let $(B,d,h)$ be a nonnegatively graded CDG\+ring.
 Then, for any left CDG\+module $(M,d_M)$ over $(B,d,h)$,
the graded left $B$\+contramodule $\Delta_B(M)$ admits a unique
odd contraderivation $d_{\Delta_B(M)}$ of degree~$1$ compatible
with the derivation~$d$ on $B$ and forming a commutative square
diagram with the odd derivation~$d_M$ and the adjunction morphism
$M\rarrow\Delta_B(M)$.
 The pair $(\Delta_B(M),d_{\Delta_B(M)})$ is a CDG\+contramodule
over $(B,d,h)$.
 The assignment $(M,d_M)\longmapsto(\Delta_B(M),d_{\Delta_B(M)})$ is
a DG\+functor $\DG(B\modl)\rarrow\DG(B\contra)$ left adjoint to
the forgetful DG\+functor $\DG(B\contra)\rarrow\DG(B\modl)$.
\end{lem}

\begin{proof}
 To check uniqueness, suppose that $d'_{\Delta_B(M)}$ and
$d''_{\Delta_B(M)}$ are two odd contraderivations of degree~$1$
on $\Delta_B(M)$ compatible with the derivation~$d$ on~$B$.
 Then the difference $d''_{\Delta_B(M)}-d'_{\Delta_B(M)}$ is
a morphism of graded left $B$\+contramodules $\Delta_B(M)\rarrow
\Delta_B(M)[1]$.
 If both $d'_{\Delta_B(M)}$ and $d''_{\Delta_B(M)}$ form commutative
square diagrams with the odd derivation~$d_M$ and the adjunction
morphism $M\rarrow\Delta_B(M)$, then the difference
$d''_{\Delta_B(M)}-d'_{\Delta_B(M)}$ is annihilated by the composition
with the adjunction morphism.
 In view of the adjunction, it follows that
$d''_{\Delta_B(M)}-d'_{\Delta_B(M)}=0$.

 Now let $(Q,d_Q)$ be a left CDG\+contramodule over $(B,d,h)$.
 Then the group of all graded left $B$\+contramodule morphisms
$\Delta_B(M)\rarrow Q[n]$ is naturally isomorphic to the group
of all graded left $B$\+module morphisms $M\rarrow Q[n]$ for every
$n\in\boZ$.
 Hence we obtain an isomorphism between the underlying graded abelian
groups of the complexes $\Hom^B(\Delta_B(M),Q)$ and $\Hom_B(M,Q)$.
 Since the adjunction map $M\rarrow\Delta_B(M)$ commutes with
the differentials, it follows that the $\Hom^B(\Delta_B(M),Q)
\simeq\Hom_B(M,Q)$ is an isomorphism of complexes of abelian groups.
 This establishes the DG\+adjunction, which implies
the DG\+functoriality.

 It remains to prove existence of the desired CDG\+contramodule
structure on $\Delta_B(M)$.
 For this purpose, consider the quasi-differential ring $(\hB,\d)$
assigned to the CDG\+ring $(B,d,h)$ by the construction of
Theorem~\ref{cdg-qdg-equivalence}.
 The key observation is that the functors $\Delta_{\hB}$ and $\Delta_B$
agree with each other, that is, they form a commutative square diagram
with the functors of restriction of scalars
$\hB\modl_\sgr\rarrow B\modl_\sgr$ and
$\hB\contra_\sgr\rarrow B\contra_\sgr$.
 The latter two functors, essentially, assign to a CDG\+(contra)module
over $(B,d,h)$ its underlying graded $B$\+(contra)module.

 Since both the functors $\Delta_{\hB}$ and $\Delta_B$ are right exact
(while the functors of restriction of scalars are exact), it suffices
to check that the two functors agree on free graded $\hB$\+modules.
 It is helpful to observe that the restriction of scalars
$\hB\modl_\sgr\rarrow B\modl_\sgr$ takes free graded $\hB$\+modules
to free graded $B$\+modules.
 The rest is a straightforward computation based on
Lemma~\ref{delta-tensor-tensor-pi}.

 Finally, given a CDG\+module $(M,d_M)$ over $(B,d,h)$, in order
to produce the induced differential~$d_{\Delta_B(M)}$ on
the graded $B$\+contramodule $\Delta_B(M)$, one simply applies
the functor $\Delta_{\hB}$ to the graded $\hB$\+module~$M$.
\end{proof}

 Now we consider the composition of the pair of adjoint DG\+functors
$\Hom_B(\Ksp(B,\allowbreak\tA),{-})\: \DG(B\modl)\rarrow\DG(\tA\modl)$
and $\Ksp(B,\tA)\ot_{\tA}{-}\: \DG(\tA\modl)\rarrow\DG(B\modl)$ from
the above discussion with the pair of adjoint DG\+functors
$\DG(B\contra)\rarrow\DG(B\modl)$ and $\DG(B\modl)\rarrow\DG(B\contra)$
provided by Lemma~\ref{delta-cdg-lemma}.
 This produces a pair of adjoint DG\+functors between the DG\+category
$\DG(\tA\modl)$ of complexes of left $\tA$\+modules and the DG\+category
$\DG(B\contra)$ of left CDG\+contramodules over $(B,d,h)$.
{\hbadness=1125\par}

 Here the right adjoint DG\+functor
$$
 \Hom_B(\Ksp(B,\tA),{-})\:\DG(B\contra)\lrarrow\DG(\tA\modl)
$$
is the composition $\DG(B\contra)\rarrow\DG(B\modl)\rarrow
\DG(\tA\modl)$ of the forgetful DG\+functor $\DG(B\contra)\rarrow
\DG(B\modl)$ with the functor $\Hom_B(\Ksp(B,\tA),{-})\:\allowbreak
\DG(B\modl)\rarrow\DG(\tA\modl)$.
 We denote it by
$$
 Q\longmapsto\Hom_B(\Ksp(B,\tA),Q)=\Hom_R^{\tau'}(\tA,Q),
$$
where $\tau'$~is the same placeholder as in
Section~\ref{koszul-duality-functors-comodule-side-subsecn}.

 The left adjoint DG\+functor
$$
 \Ksp(B,\tA)\ot_{\tA}^\Pi{-}\,\:\DG(\tA\modl)\lrarrow\DG(B\contra)
$$
is the composition $\DG(\tA\modl)\rarrow\DG(B\modl)\rarrow
\DG(B\contra)$ of the DG\+functor $\Ksp(B,\tA)\ot_{\tA}{-}\:
\DG(\tA\modl)\rarrow\DG(B\modl)$ with the DG\+functor
$\Delta_B\:\DG(B\modl)\rarrow\DG(B\contra)$.
 This DG\+functor takes a complex of left $\tA$\+modules $P=P^\bu$
to the left CDG\+contramodule over $(B,d,h)$ whose underlying
graded left $B$\+contramodule is {\hbadness=1175
$$
 \Ksp(B,\tA)\ot_{\tA}^\Pi P=B\ot_R^\Pi P=
 \Hom_R(\Hom_{R^\rop}(B,R),P)=\Hom_R(C,P),
$$
where $C$ is} a notation for the graded $R$\+$B$\+bimodule
$C=\Hom_{R^\rop}(B,R)$, as in
Section~\ref{koszul-duality-functors-comodule-side-subsecn}.
 We denote this CDG\+contramodule over $(B,d,h)$ by
$$
 \Ksp(B,\tA)\ot_{\tA}^\Pi P^\bu=
 \Hom_R^{\sigma'}(\Hom_{R^\rop}(B,R),P^\bu)=
 \Hom_R^{\sigma'}(C,P^\bu),
$$
where $\sigma'$~is also the same placeholder as in
Section~\ref{koszul-duality-functors-comodule-side-subsecn}
(see Section~\ref{revisited-subsecn} below for a further discussion of
these placeholders).

 Here the computation of $\Delta_B(B\ot_R P)$ as $B\ot_R^\Pi P$
is provided by Lemma~\ref{delta-tensor-tensor-pi}.
 It is only the differential on $\Ksp(B,\tA)\ot_{\tA}^\Pi P^\bu=
\Delta_B(\Ksp(B,\tA)\ot_{\tA}P^\bu)$ that remains to be explained.
 Lemma~\ref{delta-cdg-lemma} establishes existence of a unique
odd contraderivation of degree~$1$ on $B\ot_R^\Pi P=
\Ksp(B,\tA)\ot_{\tA}^\Pi P$ compatible with the derivation~$d$ on $B$
and agreeing with the differential on $\Ksp(B,\tA)\ot_{\tA}P^\bu$.
 Now it is straightforward to check that the differential
$d_{e'}\ot^\Pi{\id_P}+{\id_{\Ksp(B,\tA)}}\ot^\Pi d_{P^\bu}$
on $\Ksp(B,\tA)\ot_{\tA}^\Pi P^\bu$ satisfies these conditions.
 Hence we have
$$
 d_{\Delta_B(\Ksp(B,\tA)\ot_{\tA}P^\bu)}=
 d_{e'}\ot^\Pi{\id_P}+{\id_{\Ksp(B,\tA)}}\ot^\Pi d_{P^\bu},
$$
as one would expect (where $d_{e'}$~is the differential on
the dual nonhomogeneous Koszul CDG\+module $\Ksp(B,\tA)$ and
$d_{P^\bu}$ denotes the differential on the complex~$P^\bu$).
 The placeholder~$\sigma'$ is supposed to remind us of the first
summand in this formula.

 The pair of adjoint DG\+functors that we have constructed induces
a pair of adjoint triangulated functors
$$
 \Hom_B(\Ksp(B,\tA),{-})=\Hom_R^{\tau'}(\tA,{-})\:
 \Hot(B\contra)\lrarrow\Hot(\tA\modl)
$$
and
$$
 \Ksp(B,\tA)\ot^\Pi_{\tA}{-}\,=\,\Hom_R^{\sigma'}(C,{-})\:
 \Hot(\tA\modl)\lrarrow\Hot(B\contra)
$$
between the homotopy category $\Hot(\tA\modl)$ of complexes of
left $\tA$\+modules and the homotopy category $\Hot(B\contra)$
of left CDG\+contramodules over $(B,d,h)$.

\subsection{Triangulated equivalence}
\label{contramodule-side-triangulated-equivalence-subsecn}
 Let $R\subset\tV\subset\tA$ be a left finitely projective
nonhomogeneous Koszul ring, and let $(B,d,h)$ be the corresponding
right finitely projective Koszul CDG\+ring (see
Sections~\ref{cdg-ring-constructed-subsecn}
and~\ref{pbw-theorem-subsecn}\+-%
\ref{anti-equivalences-koszul-rings-subsecn}).
 The following theorem is the main result of
Section~\ref{contramodule-side-secn}.

\begin{thm} \label{contramodule-side-koszul-duality-theorem}
 The pair of adjoint triangulated functors\/ $\Hom_R^{\sigma'}(C,{-})\:
\Hot(\tA\modl)\allowbreak\rarrow\Hot(B\contra)$ and\/
$\Hom_R^{\tau'}(\tA,{-})\:\Hot(B\contra)\rarrow\Hot(\tA\modl)$ defined
in Section~\ref{koszul-duality-functors-contramodule-side-subsecn}
induces a pair of adjoint triangulated functors between
the $\tA/R$\+semicon\-traderived category of left $\tA$\+modules\/
$\sD^\sictr_R(\tA\modl)$, as defined in
Section~\ref{semicontraderived-category-subsecn}, and the contraderived
category\/ $\sD^\ctr(B\contra)$ of left CDG\+contramodules over
$(B,d,h)$, as defined in
Section~\ref{contraderived-cdg-contramodules-subsecn}.
 Under the above Koszulity assumption, the latter two funtors are
mutually inverse triangulated equivalences,
$$
 \sD^\sictr_R(\tA\modl)\simeq\sD^\ctr(B\contra).
$$
\end{thm}

\begin{proof}
 This is a generalization of~\cite[Theorem~B.2(b)]{Pkoszul} and
a particular case of~\cite[Theorem~11.8(c)]{Psemi}.
 The proof is similar to the proofs of these results, and dual-analogous
to the proof of Theorem~\ref{comodule-side-koszul-duality-theorem}.

 Recall from the discussion in
Section~\ref{graded-contramodules-subsecn} that any graded left
$B$\+contramodule $Q$ is endowed with a canonical decreasing filtration
by graded $B$\+subcontramodules $Q=G^0Q\supset G^1Q\supset G^2Q\supset
\dotsb$, where $G^mQ\subset Q$ is the image of the restriction
$$
 B^{\ge m}\ot_R^\Pi Q\overset{\pi_Q}\lrarrow Q
$$
of the contraaction map $\pi_Q\:B\ot_R^\Pi Q\rarrow Q$ to the graded
submodule $B^{\ge m}\ot_R^\Pi\nobreak Q\subset B\ot_R^\Pi Q$.
 Here, as usually, the notation $B^{\ge m}$ stands for the homogeneous
two-sided ideal $\bigoplus_{n\ge m}B^n\subset B$.
 The canonical decreasing filtration $G$ on a graded $B$\+contramodule
$Q$ is always complete, but it does not have to be separated; see
Proposition~\ref{separated-graded-contramodules-prop}
and Example~\ref{nonseparated-graded-contramodule-counterex}.
 The successive quotient contramodules $G^mQ/G^{m+1}Q$ have trivial
graded $B$\+contramodule structures, in the sense explained in
Section~\ref{graded-contramodules-subsecn}.
 
 A left CDG\+contramodule $(Q,d_Q)$ over $(B,d,h)$ is said to be
\emph{trivial} if its underlying graded $B$\+contramodule  is trivial.
 The square of the differential~$d_Q$ on a trivial CDG\+contramodule $Q$
is zero, as the curvature element $h\in B^2$ acts by zero in~$Q$.
 The DG\+category of trivial left CDG\+contramodules over $(B,d,h)$ is
equivalent to the DG\+category of complexes of left $R$\+modules.
 For any left CDG\+contramodule $(Q,d_Q)$ over $(B,d,h)$, the canonical
decreasing filtration $G$ of the underlying graded $B$\+contramodule
of $Q$ is a filtration by CDG\+subcontramodules with trivial
CDG\+contramodule successive quotients.

 The ring $\tA$ is endowed with an increasing filtration $F$, which
was defined in Section~\ref{nonhomogeneous-quadratic-subsecn}.
 We also endow the graded $R$\+$B$\+bimodule $C=\Hom_{R^\rop}(B,R)$
with the increasing filtration $F$ induced by the grading, as
in the proof of Theorem~\ref{comodule-side-koszul-duality-theorem}.
 It is important for the argument below that both $F_n\tA/F_{n-1}\tA$
and $F_nC/F_{n-1}C$ are projective (graded) left $R$\+modules for all
$n\ge0$ (in fact, they are finitely generated projective left
$R$\+modules in our assumptions).
 In particular, $\tA$ is a projective left $R$\+module and $C$ is
a projective graded left $R$\+module.
 Hence both the functors $\Hom_R^{\tau'}(\tA,{-})$ and
$\Hom_R^{\sigma'}(C,{-})$ preserve short exact sequences (of
left CDG\+contramodules over $(B,d,h)$ and of complexes of
left $\tA$\+modules, respectively).
 Both the functors also obviously preserve infinite products.

 Let $P^\bu$ be a complex of left $\tA$\+modules that is contraacyclic
as a complex of left $R$\+modules.
 We need to show that the left CDG\+contramodule
$\Hom_R^{\sigma'}(C,P^\bu)$ over $(B,d,h)$ is contraacyclic.
 Let $F$ be the decreasing filtration on the Hom module $\Hom_R(C,P)$
induced by the increasing filtration $F$ on $C$, that is
$F^n\Hom_R(C,P)=\Hom_R(C/F_{n-1}C,P)\subset\Hom_R(C,P)$.
 Then $F$ is a complete, separated filtration of the CDG\+contramodule
$\Hom_R^{\sigma'}(C,P^\bu)$ by CDG\+subcontramodules over $(B,d,h)$.
 The successive quotient CDG\+contramodules
$F^n\Hom_R^{\sigma'}(C,P^\bu)/F^{n+1}\Hom_R^{\sigma'}(C,P^\bu)$
are trivial.
 Viewed as complexes of left $R$\+modules, these are the Hom complexes
$\Hom_R(F_nC/F_{n-1}C,P^\bu)$, where $F_nC/F_{n-1}C$ is a graded
$R$\+$R$\+bimodule concentrated in the cohomological degree~$-n$ and
endowed with the zero differential.

 Since the graded left $R$\+module $F_nC/F_{n-1}C$ is projective,
the functor $\Hom_R(F_nC/\allowbreak F_{n-1}C,{-})$ takes
contraacyclic complexes of left $R$\+modules to contraacyclic
complexes of left $R$\+modules.
 Clearly, any contraacyclic complex of left $R$\+modules is also
contraacyclic as a CDG\+contramodule over $(B,d,h)$ with the trivial
contramodule structure.
 It remains to use
Lemma~\ref{cdg-contramodule-filtration-contraacyclic-lemma} in order
to conclude that $\Hom_R^{\sigma'}(C,P^\bu)$ is a contraacyclic left
CDG\+contramodule over $(B,d,h)$.

 For any contraacyclic left CDG\+contramodule $Q$ over $(B,d,h)$,
the complex of left $\tA$\+modules $\Hom_R^{\tau'}(\tA,Q)$ is
contraacyclic not only as a complex of left $R$\+modules, but even
as a complex of left $\tA$\+modules.
 This follows immediately from the fact that $\tA$ is a projective
left $R$\+module; so the functor $\Hom_R^{\tau'}(\tA,{-})$ takes
short exact sequences of left CDG\+contramodules over $(B,d,h)$ to
short exact sequences of complexes of left $\tA$\+modules (and this
functor also preserves infinite products).

 Let $P^\bu$ be an arbitrary complex of left $\tA$\+modules.
 We need to show that the cone $Y^\bu$ of the adjunction morphism of
complexes of left $\tA$\+modules
$P^\bu\rarrow\Hom_R^{\tau'}(\tA,\Hom_R^{\sigma'}(C,P^\bu))$
is contraacyclic as a complex of left $R$\+modules.
 Endow the graded left $R$\+module $\Hom_R(\tA,\Hom_R(C,P))$ with
the decreasing filtration $F$ induced by the increasing filtrations
$F$ on $C$ and $\tA$, that is $F^n\Hom_R(\tA,\Hom_R(C,P))=
\sum_{i+j=n-1}\Hom_R(\tA/F_j\tA,\Hom_R(C/F_iC,P))
\subset\Hom_R(\tA,\Hom_R(C,P))$.
 Then $F$ is a complete, separated filtration of the complex
$\Hom_R^{\tau'}(\tA,\Hom_R^{\sigma'}(C,P^\bu))$ by subcomplexes of
left $R$\+modules (but not of left $\tA$\+modules).
 The complex $P^\bu$ is endowed with the trivial filtration
$F^0P^\bu=P^\bu$, \ $F^1P^\bu=0$; and the complex $Y^\bu$ is endowed
with the induced decreasing filtration~$F$.

 Then the associated graded complex of left $R$\+modules
\begin{multline*}
 \gr_F\Hom_R^{\tau'}(\tA,\Hom_R^{\sigma'}(C,P^\bu)) \\ =
 \prod\nolimits_{n=0}^\infty
 F^n\Hom_R^{\tau'}(\tA,\Hom_R^{\sigma'}(C,P^\bu))\big/
 F^{n+1}\Hom_R^{\tau'}(\tA,\Hom_R^{\sigma'}(C,P^\bu))
\end{multline*}
is naturally isomorphic to the Hom complex
$\Hom_R({}^\tau\!K_\bu(B,A),P^\bu)$ from the complex of
(graded) $R$\+$R$\+bimodules ${}^\tau\!K_\bu(B,A)=C\ot_R^\tau A$
constructed in Section~\ref{dual-koszul-complex-subsecn} into
the complex of left $R$\+modules~$P^\bu$.
 Here $A=\gr^F\tA$ is the left finitely projective Koszul graded ring
associated with~$\tA$.
 Arguing further as in the proof of
Theorem~\ref{comodule-side-koszul-duality-theorem}, we conclude that
$\gr_FY^\bu=\Hom_R(\cone({}^\tau\!K_\bu(B,A)\to R)[-1],\,P^\bu)$ is
a contraacyclic complex of left $R$\+modules.
 By Lemma~\ref{complex-filtration-contraacyclic-lemma}, it follows that
$Y^\bu$ is a contraacyclic complex of left $R$\+modules, too.

 Let $Q=(Q,d_Q)$ be an arbitrary left CDG\+contramodule over $(B,d,h)$.
 We need to show that the cone $Z$ of the adjunction morphism of
left CDG\+contramodules $\Hom_R^{\sigma'}(C,\Hom_R^{\tau'}(\tA,Q))
\rarrow Q$ is contraacyclic as a CDG\+contramodule over $(B,d,h)$.
{\emergencystretch=0em\par}

 First of all we observe that the class of all left CDG\+contramodules
$Q$ having the desired property is preserved by the passages to
the cokernels of injective closed morphisms of CDG\+contramodules
(as well as to the kernels of surjective closed morphisms; but we need
the cokernels).
 This holds because the class of all contraacyclic left
CDG\+contramodules over $(B,d,h)$ is closed under the cokernels of
injective closed morphisms (as well as the kernels of surjective
closed morphisms).

 Consider our CDG\+contramodule $Q$ as a graded left contramodule over
the graded ring $\hB$ from Theorem~\ref{cdg-qdg-equivalence}, and
present it as a quotient of a free graded $\hB$\+contramodule.
 All free graded contramodules are separated, all subcontramodules of
separated graded contramodules are separated, and the restriction of
scalars takes separated graded $\hB$\+contramodules to separated
graded $B$\+contramodules (besides, it takes free graded
$\hB$\+contramodules to free graded $B$\+contramodules).
 So any subcontramodule of a free graded $\hB$\+contramodule is
separated as a graded $B$\+contramodule.
 This allows to present an arbitrary CDG\+contramodule $Q$ over
$(B,d,h)$ as the cokernel of an injective morphism of
CDG\+contramodules whose underlying graded $B$\+contramodules are
separated.
 Thus we can assume that $Q$ is a separated graded $B$\+contramodule.

 Now we endow $Q$ with the canonical decreasing filtration
$Q=G^0Q\supset G^1Q\supset G^2Q\supset\dotsb$; so, as mentioned above,
$G^mQ$ are CDG\+subcontramodules in $Q$ and $G^mQ/G^{m+1}Q$ are
trivial CDG\+contramodules over $(B,d,h)$.
 The filtration $G$ is separated and complete.
 Consider the induced decreasing filtrations $G$ on the left
CDG\+contramodule $\Hom_R^{\sigma'}(C,\Hom_R^{\tau'}(\tA,Q))$ and
the left CDG\+contramodule $Z$ over $(B,d,h)$.
 In particular, $G^m\Hom_R^{\sigma'}(C,\Hom_R^{\tau'}(\tA,Q))=
\Hom_R^{\sigma'}(C,\Hom_R^{\tau'}(\tA,G^mQ))$.
 These decreasing filtrations are also separated and complete.
 The successive quotient CDG\+contramodules 
$\Hom_R^{\sigma'}(C,\Hom_R^{\tau'}(\tA,G^mQ))/
\Hom_R^{\sigma'}(C,\Hom_R^{\tau'}(\tA,G^{m+1}Q))$ are computable
as the CDG\+contramodules
$\Hom_R^{\sigma'}(C,\Hom_R^{\tau'}(\tA,G^mQ/G^{m+1}Q))$, and
similarly for the CDG\+contramodules $G^mZ/G^{m+1}Z$.
 Arguing as in the proof of
Theorem~\ref{comodule-side-koszul-duality-theorem} and
using Lemma~\ref{cdg-contramodule-filtration-contraacyclic-lemma},
we reduce the question to the case of a trivial left CDG\+contramodule
$Q$ over $(B,d,h)$.
 So we can assume that $Q=Q^\bu$ is simply a complex of left
$R$\+modules endowed with the trivial graded $B$\+contramodule
structure.

 Next we endow the graded left $R$\+module $\Hom_R(C,\Hom_R(\tA,Q))$
with the decreasing filtration $F$ induced by the increasing filtrations
$F$ on $\tA$ and $C$, that is $F^n\Hom_R(C,\Hom_R(\tA,Q))=
\sum_{i+j=n-1}\Hom_R(C/F_iC,\Hom_R(\tA/F_j\tA,Q))
\subset\Hom_R(C,\Hom_R(\tA,Q))$.
 Then $F$ is a complete, separated filtration of the CDG\+contramodule
$\Hom_R^{\sigma'}(C,\Hom_R^{\tau'}(\tA,Q^\bu))$ by
CDG\+subcontramodules with trivial CDG\+contramodule successive
quotients.
 As a complex of left $R$\+modules, the trivial CDG\+contramodule
$\gr_F\Hom_R^{\sigma'}(C,\Hom_R^{\tau'}(\tA,Q^\bu))$ is naturally
isomorphic to the Hom complex $\Hom_R(K^\tau_\bu(B,A),Q^\bu)$ from
the complex of (graded) $R$\+$R$\+bimodules $K^\tau_\bu(B,A)=
A\ot_R^\sigma C$ constructed in
Section~\ref{first-koszul-complex-subsecn} into
the complex of left $R$\+modules~$Q^\bu$. {\hbadness=2650\par}

 Continuing to argue as in the proof of
Theorem~\ref{comodule-side-koszul-duality-theorem}, we conclude
that $\gr_FZ=\prod_{n=0}^\infty F^nZ/F^{n+1}Z\simeq
\Hom_R(\cone(R\to K^\tau_\bu(B,A))[-1],\,Q^\bu)$ is
a contraacyclic complex of left $R$\+modules.
 Here it is important that the complex $\cone(R\to K^\tau_\bu(B,A))$ is
``coacyclic with respect to the exact category of left $R$\+projective
$R$\+$R$\+bimodules'' in the sense of~\cite[Section~2.1]{Psemi}.
 Hence $\gr_FZ$ is also contraacyclic as a left CDG\+contramodule
over $(B,d,h)$.
 By Lemma~\ref{cdg-contramodule-filtration-contraacyclic-lemma},
it follows that $Z$ is a contraacyclic CDG\+contramodule over
$(B,d,h)$ as well.
\end{proof}

\begin{exs} \label{contramodule-side-koszul-duality-examples}
 (1)~The free left $R$\+module $R$ can be considered as a one-term
complex of left $R$\+modules, concentrated in the cohomological
degree~$0$ and endowed with the zero differential.
 This one-term complex of left $R$\+modules can be then considered
as a trivial left CDG\+contramodule over $(B,d,h)$.
 The functor $\Hom_R^{\tau'}(\tA,{-})=\Hom_B(\Ksp_{e'}(B,\tA),{-})$
takes this trivial left CDG\+contramodule $R$ over $(B,d,h)$ to
the left $\tA$\+module $\Hom_R(\tA,R)$ (viewed as a one-term
complex concentrated in the cohomological degree~$0$).

 (2)~Assume that $\tA$ is a left augmented left finitely projective
nonhomogeneous Koszul ring, in the sense of the definitions in
Sections~\ref{augmented-subsecn}
and~\ref{anti-equivalences-koszul-rings-subsecn}.
 Then the left augmentation of $\tA$ endows the base ring $R$ with
a structure of left $\tA$\+module.

 Furthermore, by Theorem~\ref{augmented-duality-functor-existence-thm}
we have $h=0$, so $(B,d,h)=(B,d)$ is a DG\+ring.
 As such, $(B,d)$ is naturally a left DG\+module over itself.
 Being nonnegatively cohomologically graded, this DG\+module has
a unique structure of DG\+contramodule over $(B,d)$ compatible with
its DG\+module structure (see
Example~\ref{graded-contramodules-bounded-grading-examples}(1)).
 (Here, by definition, a \emph{DG\+contramodule} over a nonnegatively
graded DG\+ring $(B,d)$ is the same thing as a CDG\+contramodule
over $(B,d,0)$.)

 In this setting, the functor $\Hom_R^{\sigma'}(C,{-})=
\Ksp_{e'}(B,\tA)\ot_{\tA}^\Pi{-}$ takes the left $\tA$\+module $R$
to the left DG\+contramodule $(B,d)$ over the DG\+ring $(B,d)$.
\end{exs}

\begin{rem} \label{contramodule-side-koszul-dual-complex-computes}
 Assume that $\tA$ is left augmented over $R$, as in
Remark~\ref{comodule-side-koszul-duality-functors-interpreted}(2)
and Example~\ref{contramodule-side-koszul-duality-examples}(2).
 Then, for any complex of left $\tA$\+modules $P^\bu$,
the CDG\+(contra)module $\Hom_R^{\sigma'}(C,P^\bu)$ is
a DG\+module over $(B,d,h)=(B,h)$, i.~e., a complex.
 We have a natural isomorphism of complexes of abelian groups
$\Hom_R^{\sigma'}(C,P^\bu)\simeq
\Hom_{\tA}(\tA\ot_R^{\sigma'}C,\>P^\bu)$.
 Since $\tA\ot_R^{\sigma'}C$ is a projective resolution of
the left $\tA$\+module $R$ by
Example~\ref{bimodule-resolution-etc-examples}(5), it follows
that the complex of abelian groups $\Hom_R^{\sigma'}(C,P^\bu)$
computes $\boR\Hom_{\tA}(R,P^\bu)$.
 In fact, the DG\+module $\Hom_R^{\sigma'}(C,P^\bu)$ computes
the cohomology $\Ext_{\tA}^i(R,P^\bu)=\Hom_{\sD(\tA\modl)}(R,P^\bu[i])$,
\ $i\in\boZ$, as a graded module over the graded ring of endomorphisms
of the derived category object $R\in\sD(\tA\modl)$.
\end{rem}

\subsection{Reduced contraderived category}
\label{reduced-contraderived-category-subsecn}
 We keep the assumptions of
Section~\ref{contramodule-side-triangulated-equivalence-subsecn}
(which coincide with the assumptions of
Sections~\ref{comodule-side-triangulated-equivalence-subsecn}\+-%
\ref{reduced-coderived-category-subsecn}).
 The following corollary is a special case of
Theorem~\ref{contramodule-side-koszul-duality-theorem}.

\begin{cor} \label{fin-dim-base-koszul-duality-contramodule-side}
 Assume additionally that the left homological dimension of the ring $R$
(that is, the homological dimension of the abelian category $R\modl$)
is finite.
 Then the pair of adjoint triangulated functors\/
$\Hom_R^{\sigma'}(C,{-})\:\Hot(\tA\modl)\rarrow\Hot(B\contra)$ and\/
$\Hom_R^{\tau'}(\tA,{-})\:\Hot(B\contra)\rarrow\Hot(\tA\modl)$
induces mutually inverse triangulated equivalences
\begin{equation} \label{fin-dim-base-contramodule-side}
 \sD(\tA\modl)\simeq\sD^\ctr(B\contra)
\end{equation}
between the derived category of left $\tA$\+modules and
the contraderived category of left CDG\+contramodules over $(B,d,h)$.
\end{cor}

\begin{proof}
 When the left homological dimension of the ring $R$ is finite, all
acyclic complexes of left $R$\+modules are contraacyclic
by~\cite[Remark~2.1]{Psemi}.
 So the $\tA/R$\+semicontraderived category of left $\tA$\+modules
coincides with their derived category, $\sD^\sictr_R(\tA\modl)=
\sD(\tA\modl)$, and the assertion follows from
Theorem~\ref{contramodule-side-koszul-duality-theorem}.
\end{proof}

 We would like to obtain a description of the derived category
$\sD(\tA\modl)$ in terms of some kind of exotic derived category of
left CDG\+contramodules over $(B,d,h)$ \emph{without} assuming
finiteness of the left homological dimension of~$R$.

 For this purpose, let us consider the full subcategory of acyclic
complexes of left $R$\+modules $\Acycl(R\modl)\subset\Hot(R\modl)$.
 All complexes of left $R$\+modules, and in particular acyclic
complexes of left $R$\+modules, can be viewed as trivial left
CDG\+contramodules over $(B,d,h)$.

 By an abuse of notation, let us denote simply by
$\lan\Acycl(R\modl)\ran_\sqcap$ the minimal full triangulated
subcategory in $\sD^\ctr(B\contra)$ containing all the acyclic
complexes of left $R$\+modules (viewed as trivial CDG\+contramodules)
and closed under infinite products.
 Consider the triangulated Verdier quotient category
$$
 \sD^\ctr_{R\red}(B\contra)=\sD^\ctr(B\contra)/
 \lan\Acycl(R\modl)\ran_\sqcap;
$$
let us call it the \emph{reduced contraderived category of left
CDG\+contramodules over $(B,d,h)$ relative to~$R$}.

\begin{thm} \label{reduced-koszul-duality-contramodule-side-thm}
 The pair of adjoint triangulated functors\/
$\Hom_R^{\sigma'}(C,{-})\:\Hot(\tA\modl)\allowbreak\rarrow
\Hot(B\contra)$ and\/ $\Hom_R^{\tau'}(\tA,{-})\:\Hot(B\contra)\rarrow
\Hot(\tA\modl)$ induces mutually inverse triangulated equivalences
$$
 \sD(\tA\modl)\simeq\sD^\ctr_{R\red}(B\contra)
$$
between the derived category of left $\tA$\+modules and the reduced
contraderived category of left CDG\+contramodules over $(B,d,h)$
relative to~$R$.
\end{thm}

\begin{proof}
 This is the dual-analogous assertion to
Theorem~\ref{reduced-koszul-duality-comodule-side-thm}.
 In view of Theorem~\ref{contramodule-side-koszul-duality-theorem},
only two things still need to be checked.
 Firstly, one has to show that, for any acyclic complex of left
$\tA$\+modules $P^\bu$, the CDG\+contramodule
$\Hom_R^{\sigma'}(C,P^\bu)$ over $(B,d,h)$ belongs to the triangulated
subcategory $\lan\Acycl(R\modl)\ran_\sqcap\subset\sD^\ctr(B\contra)$.
 Secondly, it needs to be established that, for any CDG\+contramodule
$Q\in\lan\Acycl(R\modl)\ran_\sqcap$, the complex of left $\tA$\+modules
$\Hom_R^{\tau'}(\tA,Q)$ is acyclic. {\hbadness=1075\par}

 Firstly, let $P=P^\bu$ be an acyclic complex of left $\tA$\+modules.
 Consider the decreasing filtration $F$ on the graded left
$R$\+module $\Hom_R(C,P)$ induced by the increasing filtration $F$
on the $R$\+$R$\+bimodule $C$, that is
$F^n\Hom_R(C,P)=\Hom_R(C/F_{n-1}C,P)$ for all $n\ge0$.
 Then $F$ is a complete, separated filtration of the CDG\+contramodule
$\Hom_R^{\sigma'}(C,P^\bu)$ over $(B,d,h)$ by its
CDG\+subcontramodules.
 The successive quotient CDG\+contramodules
$F^n\Hom_R^{\sigma'}(C,P)/F^{n+1}\Hom_R^{\sigma'}(C,P)$ are trivial
CDG\+contramodules which, viewed as complexes of left $R$\+modules,
can be computed as the Hom complexes $\Hom_R(F_nC/F_{n-1}C,P^\bu)$.

 As in the proof of
Theorem~\ref{reduced-koszul-duality-comodule-side-thm}, one can see
that the complexes of left $R$\+modules $\Hom_R(F_nC/F_{n-1}C,P^\bu)$
are acyclic.
 Now the iterated extension and telescope sequence argument from
the proofs of Lemmas~\ref{complex-filtration-contraacyclic-lemma}
and~\ref{cdg-contramodule-filtration-contraacyclic-lemma} shows
that the left CDG\+contramodule $\Hom_R^{\sigma'}(C,P^\bu)$ over
$(B,d,h)$ belongs to the triangulated subcategory
$\lan\Acycl(R\modl)\ran_\sqcap\subset\sD^\ctr(R\modl)$.

 Secondly, since the functor $\Hom_R^{\tau'}(\tA,{-})$ is triangulated
and preserves infinite products, it suffices to check that for any
acyclic complex of left $R$\+modules $Q^\bu$, viewed as a trivial
left CDG\+contramodule over $(B,d,h)$, the complex of left
$\tA$\+modules $\Hom_R^{\tau'}(\tA,Q^\bu)$ is acyclic.
 As $Q^\bu$ is a trivial CDG\+contramodule, the complex of left
$\tA$\+modules $\Hom_R^{\tau'}(\tA,Q^\bu)=\Hom_R(\tA,Q^\bu)$ is simply
obtained from the complex of left $R$\+modules $Q^\bu$ by applying
the Hom functor $\Hom_R(\tA,{-})$ termwise.
 Since $\tA$ is a projective left $R$\+module, it follows that
the complex $\Hom_R(\tA,Q^\bu)$ is acyclic.
\end{proof}

\Section{The Co-Contra Correspondence}
\label{co-contra-secn}

\subsection{Contratensor product} \label{contratensor-subsecn}
 Let $B=\bigoplus_{n=0}^\infty B_n$ be a nonnegatively graded ring,
and let $K$ be a ring endowed with a ring homomorphism $K\rarrow B_0$.
 We start from the case of ungraded comodules and contramodules before
passing to the graded ones.

 Let $N$ be an ungraded right $B$\+comodule (in the sense of
Section~\ref{ungraded-comodules-subsecn}) and $P$ be an ungraded left
$B$\+contramodule (in the sense of
Section~\ref{ungraded-contramodules-subsecn}).
 The \emph{contratensor product} $N\ocn_B P$ is an abelian group
constructed as the cokernel of (the difference of) a natural pair
of abelian group homomorphisms
$$
 N\ot_K\prod\nolimits_{n=0}^\infty\left(B_n\ot_KP\right)
 \,\rightrightarrows\,N\ot_KP.
$$
 Here the first map is simply $N\ot_K\pi_P\:N\ot_K
\prod_{n=0}^\infty\left(B_n\ot_KP\right)\rarrow N\ot_KP$, where
$\pi_P\:\prod_{n=0}^\infty B_n\ot_KP \rarrow P$ is the contraaction map.
 The second map is constructed using the right (co)action of $B$ in~$N$.

 Specifically, consider a pair of elements $y\in N$ and
$w\in\prod_{n=0}^\infty B_n\ot_KP$.
 Choose an integer $m\ge0$ such that $yB_{\ge m+1}=0$ in~$N$.
 Then the image of $y\ot w$ under the second map
$N\ot_K\left(\prod_{n=0}^\infty B_n\ot_KP\right)\rarrow N\ot_KP$
is, by definition, equal to the image of $y\ot w$ under the composition
\begin{multline*}
 N\ot_K\prod\nolimits_{n=0}^\infty\left(B_n\ot_KP\right)
 \lrarrow N\ot_K\bigoplus\nolimits_{n=0}^m\left(B_n\ot_KP\right) \\
 \,=\, \bigoplus\nolimits_{n=0}^m N\ot_K B_n\ot_K P\lrarrow
 N\ot_K P.
\end{multline*}
 Here the map $N\ot_K\prod\nolimits_{n=0}^\infty\left(B_n\ot_KP\right)
\rarrow N\ot_K\bigoplus\nolimits_{n=0}^m\left(B_n\ot_KP\right)$ is
the direct summand projection, while the map
$\bigoplus\nolimits_{n=0}^m N\ot_K B_n\ot_K P\lrarrow N\ot_K P$ is
induced by the right action maps $N\ot_K B_n\rarrow N$.

 Notice that, by construction, the contratensor product $N\ocn_BP$ is
a quotient group of the tensor product $N\ot_BP$.
 In other words, there is a natural surjective map of abelian groups
\begin{equation} \label{ungraded-tensor-contratensor-surjection}
 N\ot_BP\lrarrow N\ocn_BP.
\end{equation}
 One can also easily see that the contratensor product $N\ocn_BP$ does
not depend on the choice of a ring $K$ (so our notation is unambiguous,
and one can take, at one's convenience, $K=\boZ$, or $K=B_0$, etc.)

\begin{lem} \label{ungraded-contratensor-hom-adjunction}
 Let $E$ be an associative ring and $N$ be an $E$\+$B$\+bimodule such
that the underlying right $B$\+module of $N$ is an (ungraded) right
$B$\+comodule.
 Then the functor
$$
 \Hom_E(N,{-})\:E\modl\lrarrow B\contra
$$
constructed in Example~\ref{ungraded-co-contra-dualization-example}
is right adjoint to the contratensor product functor
$$
 N\ocn_B{-}\:B\contra\lrarrow E\modl.
$$
 In other words, for any ungraded left $B$\+contramodule $P$ and any
left $E$\+module $U$ there is a natural isomorphism of abelian groups
$$
 \Hom_E(N\ocn_BP,\>U)\simeq\Hom^B(P,\Hom_E(N,U)).
$$
\end{lem}

\begin{proof}
 The claim is that, under the adjunction isomorphism
$\Hom_E(N\ot_BP,\>U)\simeq\Hom_B(P,\Hom_E(N,U))$, the subgroup
$\Hom_E(N\ocn_BP,\>U)\subset\Hom_E(N\ot_BP,\>U)$ provided by
the surjective left $E$\+module
map~\eqref{ungraded-tensor-contratensor-surjection} corresponds to
the subgroup $\Hom^B(P,\Hom_E(N,U))\subset\Hom_B(P,\Hom_E(N,U))$
provided by the faithful forgetful functor $B\contra\rarrow B\modl$.
 One can see this by comparing the above construction of
the contratensor product $N\ocn_BP$ with the construction of
the $B$\+contramodule structure on the Hom group $\Hom_E(N,U)$ in
Example~\ref{ungraded-co-contra-dualization-example}.
\end{proof}

 Now let $N$ be a graded right $B$\+comodule (in the sense of
Section~\ref{graded-comodules-subsecn}) and $P$ be a graded left
$B$\+contramodule (in the sense of
Section~\ref{graded-contramodules-subsecn}).
 Then the \emph{contratensor product} $N\ocn_BP$ is a graded abelian
group constructed as the cokernel of (the difference of)
a natural pair of homomorphisms of graded abelian groups
$$
 N\ot_K(B\ot_K^\Pi P)\,\rightrightarrows\,N\ot_KP.
$$
 Here the first map is $N\ot_K\pi_P\:N\ot_K(B\ot_K^\Pi P)
\rarrow N\ot_KP$, where $\pi_P\:B\ot_K^\Pi P=\boM_K^\gr(P)\rarrow P$
is the contraaction map.
 The second map is constructed using the right (co)action of $B$ in~$N$.

 Specifically, consider a pair of elements $y\in N_i$ and
$w\in(B\ot_K^\Pi P)_j=\prod_{n\in\boZ}\left(B_n\ot_K P_{j-n}\right)$, \
$i$, $j\in\boZ$. 
 Choose an integer $m\ge0$ such that $yB_{\ge m+1}=0$.
 Then the image of $y\ot w$ under the second map $N\ot_K(B\ot_K^\Pi P)
\rarrow N\ot_KP$ is, by definition, equal to the image of $y\ot w$
under the composition
\begin{multline*}
 N_i\ot_K\prod\nolimits_{n=0}^\infty\left(B_n\ot_K P_{j-n}\right)
 \lrarrow N_i\ot_K\bigoplus\nolimits_{n=0}^m
 \left(B_n\ot_K P_{j-n}\right) \\
 \,=\,\bigoplus\nolimits_{n=0}^m N_i\ot_K B_n\ot_K P_{j-n}\lrarrow
 (N\ot_KP)_{i+j}.
\end{multline*}
 Here the map $N_i\ot_K\prod_{n=0}^\infty\left(B_n\ot_K P_{j-n}\right)
\rarrow N_i\ot_K\bigoplus_{n=0}^m\left(B_n\ot_K P_{j-n}\right)$
is the direct summand projection, while the map $\bigoplus_{n=0}^m
N_i\ot_K B_n\ot_K P_{j-n}\rarrow(N\ot_KP)_{i+j}$ is induced by
the right action maps $N_i\ot_K B_n\rarrow N_{i+n}$.

 By construction, the contratensor product $N\ocn_BP$ is a homogeneous
quotient group of the tensor product $N\ot_BP$.
 In other words, there is a natural surjective map of graded abelian
groups $N\ot_BP\rarrow N\ocn_BP$.
 One can also easily see that the contratensor product $N\ocn_BP$ does
not depend on the choice of a ring $K$ (so our notation is unambigous).

 For any graded left $B$\+contramodules $P$ and $Q$, we denote by
$\Hom^B(P,Q)$ the graded $B$\+contramodule Hom group from $P$ to $Q$,
that is, the graded abelian group with the components
$$
 \Hom^B(P,Q)_n=\Hom_{B\contra_\sgr}(P,Q(-n)),
$$
where $Q(-n)$ denotes the left $B$\+contramodule $Q$ with the shifted
grading, $Q(-n)_j=Q_{j+n}$ (cf.\ the notation in
Section~\ref{graded-ext-comparison-subsecn}).

\begin{lem} \label{graded-contratensor-hom-adjunction}
 Let $E=\bigoplus_{n\in\boZ}E_n$ be a graded associative ring and
$N=\bigoplus_{n\in\boZ}N_n$ be a graded $E$\+$B$\+bimodule such that
the underlying graded right $B$\+module of $N$ is a (graded) right
$B$\+comodule.
 Then the functor
$$
 \Hom_E(N,{-})\:E\modl_\sgr\lrarrow B\contra_\sgr
$$
constructed in Example~\ref{graded-co-contra-dualization-example}
is right adjoint to the contratensor product functor
$$
 N\ocn_B{-}\:B\contra_\sgr\lrarrow E\modl_\sgr.
$$
 In other words, for any graded left $B$\+contramodule $P$ and any
graded left $E$\+module $U$ there is a natural isomorphism of graded
abelian groups
$$
 \Hom_E(N\ocn_BP,\>U)\simeq\Hom^B(P,\Hom_E(N,U)).
$$
\end{lem}

\begin{proof}
 Similar to the proof of
Lemma~\ref{ungraded-contratensor-hom-adjunction}.
\end{proof}

\begin{prop} \label{contratensor-tensor-comparison}
\textup{(a)} Assume that the forgetful functor $B\contra\rarrow B\modl$
is fully faithful.
 Then the natural map of abelian groups $N\ot_BP\rarrow N\ocn_BP$
\,\eqref{ungraded-tensor-contratensor-surjection} is an isomorphism
for any ungraded right $B$\+comodule $N$ and any ungraded left
$B$\+contramodule~$P$. \par
\textup{(b)} Assume that the forgetful functor $B\contra_\sgr\rarrow
B\modl_\sgr$ is fully faithful.
 Then the natural map of graded abelian groups $N\ot_BP\rarrow N\ocn_BP$
is an isomorphism for any graded right $B$\+comodule $N$ and any
graded left $B$\+contramodule~$P$.
\end{prop}

\begin{proof}
 We will prove part~(b); part~(a) is similar.
 Let $U$ be a graded abelian group.
 Then we have
\begin{align*}
 \Hom_\boZ(N\ot_BP,\>U)&\simeq\Hom_B(P,\Hom_\boZ(N,U)), \\
 \Hom_\boZ(N\ocn_BP,\>U)&\simeq\Hom^B(P,\Hom_\boZ(N,U)).
\end{align*}
 Applying the functor $\Hom_\boZ({-},U)$ to the natural map
$N\ot_BP\rarrow N\ocn_BP$ produces the map
$\Hom^B(P,\Hom_\boZ(N,U))\rarrow\Hom_B(P,\Hom_\boZ(N,U))$
induced by the forgetful functor $B\contra_\sgr\rarrow B\modl_\sgr$.
 Since the latter map is an isomorphism for all graded abelian
groups $U$ by assumption, it follows that the former map is
an isomorphism, too.
\end{proof}

\subsection{CDG\+contramodules of the induced type}
\label{cdg-contramodules-induced-type-subsecn}
 Let $B=\bigoplus_{n=0}^\infty B^n$ be a nonnegatively graded ring,
and let $K$ be a fixed associative ring endowed with a ring homomorphism
$K\rarrow B^0$.
 Let $L=\bigoplus_{n\in\boZ} L^n$ be a graded left $K$\+module.
 By the graded left $B$\+contramodule \emph{induced from} (or
\emph{freely generated by}) a graded left $K$\+module $L$ we mean
the graded left $B$\+contramodule $\boM_K^\gr(L)=B\ot_K^\Pi L$.
 These were called the ``free modules over the monad $\boM=\boM_K^\gr$''
in Section~\ref{monads-subsecn}.
 Subsequently we described such graded $B$\+contramodules in
Lemma~\ref{delta-tensor-tensor-pi}.

 Assume that $B$ is a flat graded right $K$\+module.
 Then, given a short exact sequence of graded left $K$\+modules
$0\rarrow L'\rarrow L\rarrow L''\rarrow0$, one can apply the functor
$\boM_K^\gr$, producing a short exact sequence of graded left
$B$\+contramodules $0\rarrow\boM_K^\gr(L')\rarrow\boM_K^\gr(L)\rarrow
\boM_K^\gr(L'')\rarrow0$.
 The resulting short exact sequence of graded left $B$\+contramodules
is said to be \emph{induced from} the original short exact sequence of
graded left $K$\+modules.

 Now let $B=(B,d,h)$ be a nonnegatively graded CDG\+ring.
 We will say that a left CDG\+contramodule $(Q,d_Q)$ over $(B,d,h)$ is
\emph{of the induced type} if its underlying graded left
$B$\+contramodule $Q$ is (isomorphic to a graded left $B$\+contramodule) induced from a graded left $K$\+module.
 (Of course, the meaning of this definition depends on the choice of
a ring $K$ and a ring homomorphism $K\rarrow B^0$, which we presume
to be fixed.)

 Similarly, assume that $B$ is a flat graded right $K$\+module.
 Then we will say that a short exact sequence $0\rarrow (Q',d_{Q'})
\rarrow (Q,d_Q)\rarrow (Q'',d_{Q''})\rarrow0$ of left CDG\+contramodules
over $(B,d,h)$ (with closed morphisms between them) is a \emph{left
CDG\+contramodule short exact sequence of the induced type} if
the underlying short exact sequence of graded left $B$\+contramodules
$0\rarrow Q'\rarrow Q\rarrow Q''\rarrow0$ is isomorphic to a short exact
sequence of graded left $B$\+contramodules induced from a short exact
sequence of graded left $K$\+modules.

 As in the proof of
Proposition~\ref{graded-B-contramodules-proposition},
by the \emph{free} graded left $B$\+contramodule spanned by a graded
set $X$ we mean the graded left $B$\+contramodule $\boM_K^\gr(K[X])$
induced from the free graded left $K$\+module $K[X]$.
 The \emph{projective} graded left $B$\+contramodules are
the projective objects of the abelian category $B\contra_\sgr$;
these are the direct summands of the free graded left
$B$\+contramodules.
 A CDG\+contramodule $(P,d_P)$ over $(B,d,h)$ is said to be \emph{of
the projective} (resp., \emph{free}) \emph{type} if its underlying
graded left $B$\+contramodule is projective (resp., free).

 Denote the full DG\+subcategories in the DG\+category $\DG(B\contra)$
formed by CDG\+contramodules of the induced, free, and projective type
by $\DG(B\contra_{K\ind})$, \ $\DG(B\contra_\free)$, and
$\DG(B\contra_\proj)$, respectively.
 Clearly, all the three DG\+sub\-categories are closed under shifts,
twists, and finite direct sums; in particular, they are closed under
the cones of closed morphisms.
 Therefore, their homotopy categories $\Hot(B\contra_{K\ind})=
H^0\DG(B\contra_{K\ind})$, \ $\Hot(B\contra_\free)=
H^0\DG(B\contra_\free)$, and $\Hot(B\contra_\proj)=
H^0\DG(B\contra_\proj)$ are full triangulated subcategories in
$\Hot(B\contra)$.

 The main results of this
Section~\ref{cdg-contramodules-induced-type-subsecn}, following
below, constitute our version of~\cite[Theorem~5.5(b)]{Psemi}.

\begin{lem} \label{projective=free-cancellation-trick}
 The two full triangulated subcategories\/ $\Hot(B\contra_\free)$ and\/
$\Hot(B\contra_\proj)\subset\Hot(B\contra)$ coincide. \hbadness=3500
\end{lem}

\begin{proof}
 The argument is based on a version of Eilenberg's cancellation trick.
 Let $(Q,d_Q)$ be a left CDG\+contramodule over $(B,d,h)$ whose
underlying graded left $B$\+contramodule $Q$ is a direct summand of
the free graded left $B$\+contramodule $B\ot_K^\Pi K[X]$ spanned
by a graded set~$X$.
 Let $Y=X\times\boZ$ be the graded set with the components $Y_n=
X_n\times\boZ$.
 Consider the quasi-differential graded ring $(\hB,\d)$ corresponding
to the CDG\+ring $(B,d,h)$, and let $P=\hB\ot_K^\Pi K[Y]$ be the free
graded contramodule over $\hB$ spanned by the graded set~$Y$.
 Then $P$ can be viewed as a CDG\+contamodule $(P,d_P)$ over $(B,d,h)$.
 Moreover, $(P,d_P)$ is a contractible CDG\+contramodule (i.~e., it
represents a zero object in the homotopy category $\Hot(B\contra)$);
essentially, the homogeneous map $\d\:P\rarrow P$ of degree~$-1$ induced
by the differential~$\d$ on $\hB$ provides a contracting homotopy.
 As a graded $B$\+contramodule, $P$ is isomorphic to
$B\ot_K^\Pi K[Y]\oplus (B\ot_K^\Pi K[Y])[-1]$.
 Hence $P\oplus Q$ is a free graded $B$\+contramodule, and
the object $(Q,d_Q)\in\Hot(B\contra_\proj)$ is isomorphic to
the object $(Q,d_Q)\oplus(P,d_P)\in\Hot(B\contra_\free)$.
\end{proof}

\begin{prop} \label{homotopy-free-type-absderived-induced-prop}
 Assume that $B$ is a flat graded right $K$\+module and the left
homological dimension of the ring $K$ is finite.
 Then the inclusion\/ $\Hot(B\contra_\free)\rarrow
\Hot(B\contra_{K\ind})$ induces a triangulated equivalence between\/
$\Hot(B\contra_\free)$ and the triangulated quotient category of\/
$\Hot(B\contra_{K\ind})$ by its minimal full triangulated subcategory
containing all the totalizations of left CDG\+contramodule short exact
sequences of the induced type over $(B,d,h)$.
\end{prop}

\begin{proof}
 Essentially, the claim is that the full subcategory of
CDG\+contramodules of the free type and the minimal full triangulated
subcategory containing the totalizations of CDG\+contramodule short
exact sequences of the induced type form a semiorthogonal decomposition
of $\Hot(B\contra_{K\ind})$.
 Firstly, this means that the complex of abelian groups $\Hom^B(P,Z)$ is
acyclic for any $P\in\Hot(B\contra_\free)$ and the totalization $Z$ of
any CDG\+contramodule short exact sequence $0\rarrow Q'\rarrow Q\rarrow
Q''\rarrow0$ of the induced type.
 This is a particular case of a standard observation
(cf.~\cite[Theorem~3.5(b)]{Pkoszul}), and the assumption that the short
exact sequence is of the induced type is not needed for it;
see the proof of Theorem~\ref{contraderived-category-theorem} below.

 Secondly, it is claimed that for any object
$Q\in\Hot(B\contra_{K\ind})$ there exists a distinguished triangle
$P\rarrow Q\rarrow Z\rarrow P[1]$ in $\Hot(B\contra_{K\ind})$ with
$P\in\Hot(B\contra_\free)$ and $Z$ belonging to the minimal triangulated
subcategory of $\Hot(B\contra_{K\ind})$ containing the totalizations
of CDG\+contramodule short exact sequences of the induced type.
 This is provable by a variation of the standard argument
from~\cite[proof of Theorem~3.6]{Pkoszul}.

 Let $(Q,d_Q)$ be a CDG\+contramodule structure on the graded left
$B$\+contramodule $Q=B\ot_K^\Pi L$ induced from a graded left
$K$\+module~$L$.
 Choose a free graded left $K$\+module $G=K[X]$ together with
a surjective morphism of graded left $K$\+modules $g\:G\rarrow L$.
 Consider the free graded left $\hB$\+contramodule $F=\hB\ot_K^\Pi G$
induced from~$G$.
 The CDG\+contramodule $(Q,d_Q)$ over $(B,d,h)$ can be also viewed
as a graded $\hB$\+contramodule.
 Let $f\:F\rarrow Q$ be the graded left $\hB$\+contramodule morphism
corresponding to the composition of graded left $K$\+module maps
$G\rarrow L\rarrow Q$.
 Then the restriction of~$f$ to the graded $B$\+subcontramodule
$B\ot_K^\Pi G\subset\hB\ot_K^\Pi G$ is the graded $B$\+contramodule
map $B\ot_K^\Pi G\rarrow B\ot_K^\Pi L=Q$ obtained by applying
the functor $B\ot_K^\Pi{-}$ to the graded left $K$\+module map
$g\:G\rarrow L$.

 The key observation is that the whole map $f\:F\rarrow Q$, viewed
as a morphism of graded left $B$\+contramodules, is induced from
a morphism of graded left $K$\+modules.
 Put $F'=B\ot_K^\Pi G\subset\hB\ot_K^\Pi G=F$.
 Then the quotient graded left $B$\+contramodule $F''=F/F'$ is
isomorphic to $F'[-1]$; so it is also a free graded left
$B$\+contramodule.
 Hence the short exact sequence of graded left $B$\+contramodules
$0\rarrow F'\rarrow F\rarrow F''\rarrow0$ splits; we can choose
a splitting and consider $F''$ as a graded left $B$\+subcontramodule
in~$F$.
 Put $f'=f|_{F'}\:F'\rarrow Q$ and $f''=f|_{F''}\:F''\rarrow Q$;
so $f=(f',f'')\:F'\oplus F''\rarrow Q$.
 Then $f'=B\ot_K^\Pi g$ is a surjective morphism of graded left
$B$\+contramodules (as $g\:G\rarrow L$ is a surjective morphism
of graded left $K$\+modules).
 Since $F''$ is a projective graded left $B$\+contramodule,
there exists a graded left $B$\+contramodule morphism
$t\:F''\rarrow F'$ such that $f''=f'\circ t$.
 Using the map~$t$, one can construct an automorphism of the graded
left $B$\+contramodule $F=F'\oplus F''$ such that the composition
$F'\oplus F''\rarrow F'\oplus F''\overset f\rarrow Q$ is equal to
the map $(f',0)\:F'\oplus F''\rarrow Q$.
 
 Put $F_0=F$ and $G_0=G\oplus G[-1]$; and denote by $Q_1$ the kernel
of the graded left $\hB$\+contramodule morphism $f\:F_0\rarrow Q$.
 Then $0\rarrow Q_1\rarrow F_0\rarrow Q\rarrow0$ is a short exact
sequence of graded left $\hB$\+contramodules.
 We have shown that, as a short exact sequence of graded left
$B$\+contramodules, it is induced from a short exact sequence of
graded left $K$\+modules $0\rarrow L_1\rarrow G_0\rarrow L\rarrow0$
(where $L_1$ is the kernel of the morphism $(g,0)\:G_0=G\oplus G[-1]
\rarrow L$).
 Now we can consider $0\rarrow Q_1\rarrow F_0\rarrow Q\rarrow0$ as
a CDG\+contramodule short exact sequence of the induced type over
$(B,d,h)$.

 It remains to iterate this construction, applying it to
the CDG\+contramodule $Q_1$ in place of $Q$, etc.
 Proceeding in this way, we produce an exact sequence of
CDG\+contramodules $\dotsb\rarrow F_2\rarrow F_1\rarrow F_0\rarrow Q
\rarrow0$ such that its underlying exact sequence of graded left
$B$\+contramodules is induced from an exact sequence of graded left
$K$\+modules $\dotsb\rarrow G_2\rarrow G_1\rarrow G_0\rarrow L\rarrow0$.
 Here $G_i$ are free graded left $K$\+modules by construction.
 Since the left homological dimension of the ring $K$ is finite by
assumption, there exists $k\ge0$ such that the image $L_k$ of
the morphism $G_k\rarrow G_{k-1}$ is a projective graded left
$K$\+module.
 
 Consider the finite exact sequence of CDG\+contramodules
$0\rarrow Q_k\rarrow F_{k-1}\rarrow\dotsb\rarrow F_0\rarrow Q\rarrow0$
over $(B,d,h)$, where $Q_k$ is the image of the CDG\+contramodule
morphism $F_k\rarrow F_{k-1}$.
 Denote by $(Z,d_Z)$ the total CDG\+contramodule of this finite exact
sequence of CDG\+contramodules; and denote by $(P,d_P)$ the total
CDG\+contramodule of the complex of CDG\+contramodules $Q_k\rarrow
F_{k-1}\rarrow\dotsb\rarrow F_0$.
 Then we have a distinguished triangle $P\rarrow Q\rarrow Z\rarrow P[1]$
in $\Hot(B\contra_{K\ind})$.

 The underlying exact sequence of graded left $B$\+contramodules
$0\rarrow Q_k \rarrow F_{k-1}\rarrow\dotsb\rarrow F_0\rarrow Q\rarrow0$
is induced from the exact sequence of graded left $K$\+modules
$0\rarrow L_k\rarrow G_{k-1}\rarrow\dotsb\rarrow G_0\rarrow L\rarrow0$.
 Hence it follows that the CDG\+contramodule $(Z,d_Z)$ over $(B,d,h)$
belongs to the minimal triangulated subcategory in
$\Hot(B\contra_{K\ind})$ containing the totalizations of
CDG\+contramodule short exact sequences of the induced type.
 Finally, the underlying graded left $B$\+contramodule $P$ of
the CDG\+contramodule $(P,d_P)$ is projective, since the graded
left $B$\+contramodule $Q_k$ is induced from the projective graded
left $K$\+module~$L_k$.
 By Lemma~\ref{projective=free-cancellation-trick},
the CDG\+contramodule $(P,d_P)$ is homotopy equivalent to
a CDG\+contramodule of the free type over $(B,d,h)$.
\end{proof}

\begin{prop} \label{contraderived-induced-prop}
 Assume that the right $K$\+module $B^n$ is finitely generated and
projective for every $n\ge0$.
 Then the inclusion\/ $\Hot(B\contra_{K\ind})\rarrow\Hot(B\contra)$
induces a triangulated equivalence between the contraderived category\/
$\sD^\ctr(B\contra)$ (as defined in
Section~\ref{contraderived-cdg-contramodules-subsecn}) and
the triangulated quotient category of the homotopy category\/
$\Hot(B\contra_{K\ind})$ by its intersection with the full triangulated
subcategory\/ $\Acycl^\ctr(B\contra)$ of contraacyclic
CDG\+contramodules in\/ $\Hot(B\contra)$,
$$
 \Hot(B\contra_{K\ind})/
 (\Acycl^\ctr(B\contra)\cap\Hot(B\contra_{K\ind}))
 \,\simeq\,\sD^\ctr(B\contra).
$$
\end{prop}

\begin{proof}
 In view of a well-known lemma (\cite[Lemma~2.6]{Psemi}
or~\cite[Lemma~1.6(a)]{Pkoszul}), it suffices to show that for any
left CDG\+contramodule $(S,d_S)$ over $(B,d,h)$ there exists a left
CDG\+contramodule $(Q,d_Q)$ of the induced type together with a closed
morphism of CDG\+contramodules $(Q,d_Q)\rarrow (S,d_S)$ whose cone is
a contraacyclic CDG\+contramodule over $(B,d,h)$.
 A standard argument (cf.~\cite[Theorem~3.8]{Pkoszul}) reduces
the question to showing that countable products of free graded left
$B$\+contramodules are induced graded left $B$\+contramodules.

 Let us spell out some details.
 Consider the CDG\+contramodule $(S,d_S)$ as a graded left
$\hB$\+contramodule.
 Let $\dotsb\rarrow P_2\rarrow P_1\rarrow P_0\rarrow S\rarrow0$ be
a resolution of the object $S\in\hB\contra_\sgr$ by free graded
left $\hB$\+contramodules~$P_i$.
 Consider this resolution as a complex of left CDG\+contramodules
$\dotsb\rarrow P_2\rarrow P_1\rarrow P_0\rarrow0$ over $(B,d,h)$ and
totalize it by taking infinite products along the diagonals.
 We refer to~\cite[Section~1.2]{Pkoszul} for a discussion of
totalizations of complexes in DG\+categories.
 Notice that the infinite products in $B\contra_\sgr$ agree with
those in $K\modl_\sgr$, so there is no ambiguity in how
the totalization construction needs to be applied.
 Let $(Q,d_Q)$ be the resulting total CDG\+contramodule over $(B,d,h)$.
 Then the underlying graded left $B$\+contramodule $Q$ of $(Q,d_Q)$
is simply the countable product $Q=\prod_{i=0}^\infty P_i[i]$.
 Once again we recall that the forgetful functor $\hB\contra_\sgr
\rarrow B\contra_\sgr$ takes free graded left $\hB$\+contramodules
to free graded left $B$\+contramodules.

 The cone $(Y,d_Y)$ of the natural closed morphism of CDG\+contramodules
$(Q,d_Q)\allowbreak\rarrow (S,d_S)$ is the infinite product totalization
of the exact complex $\dotsb\rarrow P_2\rarrow P_1\rarrow P_0\rarrow S
\rarrow0$.
 As such, it admits a complete, separated decreasing filtration $F$ by
the totalizations of the subcomplexes of canonical filtration of
the complex $\dotsb\rarrow P_2\rarrow P_1\rarrow P_0\rarrow S\rarrow0$.
 The successive quotient CDG\+contramodules $F^nY/F^{n+1}Y$ are
contractible.
 In fact, the CDG\+contramodules $F^nY/F^{n+1}Y$ are the cones of
identity endomorphisms of the CDG\+contramodules of cocycles of
our acyclic complex of CDG\+contramodules (up to some shifts).
 By Lemma~\ref{cdg-contramodule-filtration-contraacyclic-lemma}, it
follows that the CDG\+contramodule $(Y,d_Y)$ is contraacyclic.

 It remains to check that the class of all induced graded left
$B$\+contramodules is closed under infinite products in $B\contra_\sgr$.
 This is where we use the assumption that $B^n$ is a finitely
generated projective right $K$\+module for every $n\ge0$.
 Under this assumption, one can easily see that the functor
$B\ot_K^\Pi{-}$ (viewed as a functor $K\modl_\sgr\rarrow K\modl_\sgr$
or $K\modl_\sgr\rarrow B\contra_\sgr$) preserves infinite products.
 In fact, it suffices that $B^n$ be a finitely presentable right
$K$\+module for every $n\ge0$.
\end{proof}

\begin{thm} \label{contraderived-category-theorem}
 Assume that $B^n$ is a finitely generated projective right $K$\+module
for every $n\ge0$ and the left homological dimension of the ring $K$ is
finite.
 Then the composition of the fully faithful inclusion\/
$\Hot(B\contra_\proj)\rarrow\Hot(B\contra)$ with the triangulated
quotient functor\/ $\Hot(B\contra)\rarrow \sD^\ctr(B\contra)$ is
a triangulated equivalence
$$
 \Hot(B\contra_\proj)\simeq\sD^\ctr(B\contra)
$$
between the homotopy category of left CDG\+contramodules of
the projective type over $(B,d,h)$ and the contraderived category
of left CDG\+contramodules over $(B,d,h)$.
 Furthermore, the intersection of two full subcategories
$$
 \Acycl^\ctr(B\contra)\cap\Hot(B\contra_{K\ind})\subset\Hot(B\contra)
$$
coincides with the minimal full triangulated subcategory in\/
$\Hot(B\contra_{K\ind})$ containing all the totalizations of
CDG\+contramodule short exact sequences of the induced type
over $(B,d,h)$.
\end{thm}

\begin{proof}
 In the first assertion of the theorem, the claim is that the full
subcategory of CDG\+contramodules of the projective type
$\Hot(B\contra_\proj)$ and the full subcategory of contraacyclic
CDG\+contramodules $\Acycl^\ctr(B\contra)$ form a semiorthogonal
decomposition of the homotopy category $\Hot(B\contra)$.
 First of all, this means that the complex of abelian groups
$\Hom^B(P,Z)$ is acyclic for any $P\in\Hot(B\contra_\proj)$ and
$Z\in\Acycl^\ctr(B\contra)$.
 A particular case of this observation was already mentioned in
the proof of
Proposition~\ref{homotopy-free-type-absderived-induced-prop}.

 Indeed, let us fix a CDG\+contramodule $(P,d_P)\in
\Hot(B\contra_\proj)$.
 Then the class of all CDG\+contramodules $(Z,d_Z)\in\Hot(B\contra)$
for which the complex $\Hom^B(P,Z)$ is acyclic is closed under
shifts, cones, and infinite products in $\Hot(B\contra)$.
 Thus it remains to check that the totalization
$Z=\Tot(Q'\to Q\to Q'')$ of any short exact sequence of
CDG\+contramodules $0\rarrow (Q',d_{Q'})\rarrow (Q,d_Q)\rarrow
(Q'',d_{Q''})\rarrow0$ over $(B,d,h)$ has the desired property.
 Here it suffices to notice that the complex of abelian groups
$\Hom^B(P,Z)$ is the totalization of the bicomplex with three rows
$\Hom^B(P,Q')\rarrow\Hom^B(P,Q)\rarrow\Hom^B(P,Q'')$.
 The short sequence of complexes of abelian groups
$0\rarrow\Hom^B(P,Q')\rarrow\Hom^B(P,Q)\rarrow\Hom^B(P,Q'')\rarrow0$
is exact, since the graded left $B$\+contramodule $P$ is projective.
 The totalization of any short exact sequence of complexes of
abelian groups is an exact complex.

 Now we can prove the second assertion of the theorem.
 Let $(Z,d_Z)$ be a left CDG\+contramodule of the induced type over
$(B,d,h)$ which is contraacyclic as a left CDG\+contramodule
over $(B,d,h)$.
 Then, according to the above argument, the complex $\Hom^B(P,Z)$ is
acyclic for any left CDG\+contramodule $(P,d_P)$ of the free type
over $(B,d,h)$.
 In view of the semiorthogonal decomposition of the category
$\Hot(B\contra_{K\ind})$ constructed in the proof of
Proposition~\ref{homotopy-free-type-absderived-induced-prop},
it follows that $(Z,d_Z)$ belongs to the minimal triangulated
subcategory of $\Hot(B\contra_{K\ind})$ containing the totalizations
of CDG\+contramodule short exact sequences of the induced type.

 Finally, the assertion that the functor $\Hot(B\contra_\proj)
\rarrow\sD^\ctr(B\contra)$ is a triangulated equivalence can be now
obtained by comparing the results of
Lemma~\ref{projective=free-cancellation-trick},
Proposition~\ref{homotopy-free-type-absderived-induced-prop},
and Proposition~\ref{contraderived-induced-prop}.
 This also implies the semiorthogonal decomposition promised in
the beginning of this proof (given that the semiorthogonality of
the two subcategories in question has been established already).
\end{proof}

\subsection{CDG\+comodules of the coinduced type}
\label{cdg-comodules-coinduced-type-subsecn}
 Let $B=\bigoplus_{n=0}^\infty B^n$ be a nonnegatively graded ring,
and let $K$ be a fixed ring endowed with a ring homomorphism
$K\rarrow B^0$.
 Similarly to the definition of graded right $B$\+comodules in
Section~\ref{graded-comodules-subsecn}, one can define graded left
$B$\+comodules.
 So a graded left $B$\+module $M=\bigoplus_{n\in\boZ}M^n$ is said
to be a \emph{$B$\+comodule} (or a \emph{graded left $B$\+comodule})
if for every element $x\in M^j$, \ $j\in\boZ$, there exists
an integer $m\ge0$ such that $B^{\ge m+1}x=0$.
 We denote the full subcategory of graded left $B$\+comodules
by $B\comodl_\sgr\subset B\modl_\sgr$.

 Let $L=\bigoplus_{n\in\boZ}L^n$ be a graded left $K$\+module.
 By the graded left $B$\+comodule \emph{coinduced from} (or
\emph{cofreely cogenerated by}) a graded left $K$\+module $L$ we mean
the graded left $B$\+comodule $\Hom_K^\Sigma(B,L)$, in the notation
of Section~\ref{graded-ext-comparison-subsecn}
(cf.\ Lemma~\ref{gamma-hom-hom-sigma}).

 Assume that $B$ is a projective graded left $K$\+module.
 Then, given a short exact sequence of graded left $K$\+modules
$0\rarrow L'\rarrow L\rarrow L''\rarrow 0$, one can apply the functor
$\Hom_K^\Sigma(B,{-})$, producing a short exact sequence of graded
left $B$\+comodules $0\rarrow\Hom_K^\Sigma(B,L')\rarrow
\Hom_K^\Sigma(B,L)\rarrow\Hom_K^\Sigma(B,L'')\rarrow0$.
 The resulting short exact sequence of graded left $B$\+comodules is
said to be \emph{coinduced from} the original short exact sequence
of graded left $K$\+modules.

 Now let $B=(B,d,h)$ be a nonnegatively graded CDG\+ring.
 Similarly to the definition of right CDG\+comodules over $(B,d,h)$
in Section~\ref{coderived-cdg-comodules-subsecn}, one can define
left CDG\+comodules over $(B,d,h)$.
 So a left CDG\+module $(M,d_M)$ over $(B,d,h)$ is said to be
a \emph{CDG\+comodule} (or a \emph{left CDG\+comodule}) over $(B,d,h)$
if the graded left $B$\+module $M$ is a graded left $B$\+comodule.
 The full DG\+subcategory in $\DG(B\modl)$ whose objects are
the left CDG\+comodules over $(B,d,h)$ is called the \emph{DG\+category
of left CDG\+comodules over $(B,d,h)$} and denoted by $\DG(B\comodl)$.
 Its triangulated homotopy category is denoted by
$\Hot(B\comodl)=H^0\DG(B\comodl)$.

 We will say that a left CDG\+comodule $(M,d_M)$ over $(B,d,h)$ is
\emph{of the coinduced type} if its underlying graded left
$B$\+comodule $M$ is (isomorphic to a graded left $B$\+comodule)
coinduced from a graded left $K$\+module.
 (The meaning of this definition depends on the choice of a ring $K$
and a ring homomorphism $K\rarrow B^0$, which we presume to be fixed.)

 Assume that $B$ is a projective graded left $K$\+module.
 Then we will say that a short exact sequence $0\rarrow(M',d_{M'})
\rarrow(M,d_M)\rarrow(M'',d_{M''})\rarrow0$ of left CDG\+comodules
over $(B,d,h)$ (with closed morphisms between them) is
a \emph{left CDG\+comodule short exact sequence of the coinduced
type} if the underlying short exact sequence of graded left
$B$\+comodules $0\rarrow M'\rarrow M\rarrow M''\rarrow0$ is
isomorphic to a short exact sequence of graded left $B$\+comodules
coinduced from a short exact sequence of graded left $K$\+modules.

 By the \emph{cofree} graded left $B$\+comodule spanned by a graded
set $X$ we mean the graded left $B$\+comodule
$\Hom_K^\Sigma(B,K[X]^+)$ coinduced from the cofree graded left
$K$\+module $K[X]^+$.
 Here $F^+$ is the notation for the graded character group
$\Hom_\boZ(F,\boQ/\boZ)$ of an abelian group $F$ (see
Section~\ref{graded-ext-comparison-subsecn}).
 The \emph{injective} graded left $B$\+comodules are the injective
objects of the abelian category $B\comodl_\sgr$; these are
the direct summands of the cofree graded left $B$\+comodules
(see the discussion in the proof of
Theorem~\ref{graded-comodule-Ext-comparison}).

 A CDG\+comodule $(J,d_J)$ over $(B,d,h)$ is said to be
\emph{of the injective} (resp., \emph{cofree}) \emph{type} if its
underlying graded left $B$\+comodule is injective (resp., cofree).
 Denote the full DG\+subcategories in the DG\+category $\DG(B\comodl)$
formed by CDG\+comodules of the coinduced, cofree, and injective type
by $\DG(B\comodl_{K\coind})$, \ $DG(B\comodl_\cofr)$, and
$\DG(B\comodl_\inj)$, respectively.
 All the three DG\+subcategories are closed under shifts, twists, and
finite direct sums; in particular, they are closed under the cones of
closed morphisms.
 Therefore, their homotopy categories $\Hot(B\comodl_{K\coind})=
H^0\DG(B\comodl_{K\coind})$, \ $\Hot(B\comodl_\cofr)=
H^0\DG(B\comodl_\cofr)$, and $\Hot(B\comodl_\inj)=
H^0\DG(B\comodl_\inj)$ are full triangulated subcategories in
$\Hot(B\comodl)$. {\hbadness=1425\par}

 The main results of this
Section~\ref{cdg-comodules-coinduced-type-subsecn}, following below,
form our version of~\cite[Theorem~5.5(a)]{Psemi}.

\begin{lem} \label{injective=cofree-cancellation-trick}
 The two full triangulated subcategories\/ $\Hot(B\comodl_\cofr)$ and\/
$\Hot(B\comodl_\inj)\subset\Hot(B\comodl)$ coincide. \hbadness=3000
\end{lem}

\begin{proof}
 This is a dual-analogous version of
Lemma~\ref{projective=free-cancellation-trick}, provable by
a similar argument using the cancellation trick.
\end{proof}

\begin{prop} \label{homotopy-cofree-type-absderived-coinduced-prop}
 Assume that $B$ is a projective graded left $K$\+module and the left
homological dimension of the ring $K$ is finite.
 Then the inclusion\/ $\Hot(B\comodl_\cofr)\rarrow
\Hot(B\comodl_{K\coind})$ induces a triangulated equivalence between\/
$\Hot(B\comodl_\cofr)$ and the triangulated quotient category of\/
$\Hot(B\comodl_{K\coind})$ by its minimal full triangulated subcategory
containing all the totalizations of left CDG\+comodule short exact
sequences of the coinduced type over $(B,d,h)$.
\end{prop}

\begin{proof}
 Essentially, the claim is that the full subcategory of CDG\+comodules
of the cofree type and the minimal full triangulated subcategory
containing the totalizations of CDG\+comodule short exact sequences of
the coinduced type form a semiorthogonal decomposition of
$\Hot(B\comodl_{K\coind})$.
 Firstly, this means that the complex of abelian groups
$\Hom_B(Z,J)$ is acyclic for any $J\in\Hot(B\comodl_\cofr)$ and
the totalization $Z$ of any CDG\+comodule short exact sequence
$0\rarrow M'\rarrow M\rarrow M''\rarrow0$ of the coinduced type.
 This is a particular case of a standard observation
(cf.~\cite[Theorem~3.6(a)]{Pkoszul}), and the assumption that
the short exact sequence is of the coinduced type is not needed for it;
see the proof of Theorem~\ref{coderived-category-theorem} below.

 Secondly, it is claimed that for any object
$N\in\Hot(B\comodl_{K\coind})$ there exists a distinguished triangle
$Z\rarrow N\rarrow J\rarrow Z[1]$ in $\Hot(B\comodl_{K\coind})$ with
$J\in\Hot(B\comodl_\cofr)$ and $Z$ belonging to the minimal
triangulated subcategory of $\Hot(B\comodl_{K\coind})$ containing
the totalizations of CDG\+comodule short exact sequences of
the coinduced type.
 This is provable by a variation of the standard argument
form~\cite[proof of Theorem~3.6]{Pkoszul}.
 We skip the further details, which are very similar to the proof of
the dual-analogous result in
Proposition~\ref{homotopy-free-type-absderived-induced-prop}.
\end{proof}

\begin{prop} \label{coderived-coinduced-prop}
 Assume that the left $K$\+module $B^n$ is finitely generated
for every $n\ge0$.
 Then the inclusion\/ $\Hot(B\comodl_{K\coind})\rarrow\Hot(B\comodl)$
induces a triangulated equivalence between the coderived category\/
$\sD^\co(B\comodl)$ (as defined in
Section~\ref{coderived-cdg-comodules-subsecn}) and the triangulated
quotient category of the homotopy category\/ $\Hot(B\comodl_{K\coind})$
by its intersection with the full triangulated subcategory\/
$\Acycl^\co(B\comodl)$ of coacyclic CDG\+comodules
in\/ $\Hot(B\comodl)$,
$$
 \Hot(B\comodl_{K\coind})/
 (\Acycl^\co(B\comodl)\cap\Hot(B\comodl_{K\coind})
 \,\simeq\,\sD^\co(B\comodl).
$$
\end{prop}

\begin{proof}
 In view of~\cite[Lemma~2.6]{Psemi} or~\cite[Lemma~1.6(b)]{Pkoszul},
it suffices to show that for any left CDG\+comodule $(M,d_M)$
over $(B,d,h)$ there exists a left CDG\+comodule $(N,d_N)$ of
the coinduced type together with a closed morphism of CDG\+comodules
$(M,d_M)\rarrow(N,d_N)$ whose cone is a coacyclic CDG\+comodule
over~$(B,d,h)$.
 A standard argument using
Lemma~\ref{cdg-comodule-filtration-coacyclic-lemma}
(cf.~\cite[proof of Theorem~3.7]{Pkoszul}) reduces the question
to showing that countable direct sums of cofree graded
left $B$\+comodules are coinduced graded left $B$\+comodules.
 We skip the details, which are very similar to the ones in
the proof of the dual-analogous result in
Proposition~\ref{contraderived-induced-prop}.

 Hence it remains to check that the class of all coinduced graded left
$B$\+comodules is closed under infinite direct sums in $B\comodl_\sgr$.
 This is where we use the assumption that $B^n$ is a finitely generated
left $K$\+module for every $n\ge0$.
 Under this assumption, one can easily see that the functor
$\Hom_K^\Sigma(B,{-})$ (viewed as a functor $K\modl_\sgr\rarrow
K\modl_\sgr$ or $K\modl_\sgr\rarrow B\comodl_\sgr$) preserves infinite
direct sums.
\end{proof}

\begin{thm} \label{coderived-category-theorem}
 Assume that $B^n$ is a finitely generated projective left $K$\+module
for every $n\ge0$ and the left homological dimension of the ring $K$
is finite.
 Then the composition of the fully faithful inclusion\/
$\Hot(B\comodl_\inj)\rarrow\Hot(B\comodl)$ with the triangulated
quotient functor\/ $\Hot(B\comodl)\rarrow\sD^\co(B\comodl)$ is
a triangulated equivalence
$$
 \Hot(B\comodl_\inj)\simeq\sD^\co(B\comodl)
$$
between the homotopy category of left CDG\+comodules of the injective
type over $(B,d,h)$ and the coderived category of left CDG\+comodules
over $(B,d,h)$.
 Furthermore, the intersection of two full subcategories
$$
 \Acycl^\co(B\comodl)\cap\Hot(B\comodl_{K\coind})\subset
 \Hot(B\comodl)
$$
coincides with the minimal full triangulated subcategory in\/
$\Hot(B\comodl_{K\coind})$ containing all the totalizations of
CDG\+comodule short exact sequences of the coinduced type
over $(B,d,h)$.
\end{thm}

\begin{proof}
 In the first assertion of the theorem, the claim is that the full
subcategory of CDG\+comodules of the injective type
$\Hot(B\comodl_\inj)$ and the full subcategory of coacyclic
CDG\+comodules $\Acycl^\co(B\comodl)$ form a semiorthogonal
decomposition of the homotopy category $\Hot(B\comodl)$.
 First of all, this means that the complex of abelian groups
$\Hom_B(Z,J)$ is acyclic for any $Z\in\Acycl^\co(B\comodl)$ and
$J\in\Hot(B\comodl_\inj)$.
 A particular case of this observation was already mentioned in
the proof of
Proposition~\ref{homotopy-cofree-type-absderived-coinduced-prop}.

 The proof of the theorem is based on
Lemma~\ref{injective=cofree-cancellation-trick},
Proposition~\ref{homotopy-cofree-type-absderived-coinduced-prop},
and Proposition~\ref{coderived-coinduced-prop}.
 We skip the details, which are very similar to the proof of
the dual-analogous Theorem~\ref{contraderived-category-theorem}.
\end{proof}

\subsection{The diagonal CDG\+bicomodule}
\label{diagonal-cdg-bicomodule-subsecn}
 Let $B=\bigoplus_{n=0}^\infty B^n$ be a nonnegatively graded ring,
and let $K$ be a ring endowed with a ring homomorphism $K\rarrow B^0$.
 Put $C_n=C^{-n}=\Hom_{K^\rop}(B^n,K)$.
 We will use the pairing notation
\begin{gather*}
 \lan c,b\ran = c(b) \in K \quad
 \text{for all $c\in C_n$, \,$b\in B^n$, \,$n\ge0$}, \\
 \lan c,b\ran = 0 \quad \text{for all $c\in C_j$, \,$b\in B^i$,
 \,$i\ne j$}.
\end{gather*}
 The graded abelian group $C=\bigoplus_{n\le0} C^n$ has a natural
structure of graded $K$\+$B$\+bimodule, with the left action of $K$
given by the rule
$$
 \lan kc,b\ran = k\lan c,b\ran \quad
 \text{for all $k\in K$, \,$c\in C^{-n}$, \,$b\in B^n$}
$$
and the right action of $B$ given by the rule
$$
 \lan cb,b'\ran = \lan c,bb'\ran \quad
 \text{for all $c\in C$, \ $b$, $b'\in B$}.
$$

 Furthermore, assume that $B^n$ is a finitely generated projective right
$K$\+module for every $n\ge0$.
 Then a straightforward generalization of the construction from
Remark~\ref{B-sharp-remark} endows $C$ with the structure of a graded
coring over~$K$.
 (The discussion in Remark~\ref{B-sharp-remark} presumes
the base ring $K=R=B^0$.)
 The grading components $\mu_{i,j}\:C_{i+j}\rarrow C_i\ot_K C_j$,
\ $i$, $j\ge0$, of the comultiplication map $\mu\:C\rarrow C\ot_KC$ are
obtained by applying the functor $\Hom_{K^\rop}({-},K)$ to
the multiplication maps $B^j\ot_K B^i\rarrow B^{i+j}$ and taking into
account the natural isomorphisms of Lemma~\ref{tensor-dual-lemma}(b).
 In the pairing notation,
$$
 \lan \mu_{i,j}(c),\,b''\ot b'\ran = \lan c,b''b'\ran
 \quad\text{for all $c\in C_{i+j}$, \,$b'\in B^i$, \,$b''\in B^j$},
$$
where the pairing $\lan\ ,\ \ran\:C_i\ot_K C_j\ot_K B^j\ot_K B^i
\rarrow K$ is defined by the rule
$$
 \lan c'\ot c'',\,b''\ot b'\ran = \lan c'\lan c'',b''\ran,b'\ran =
 \lan c',\lan c'',b''\ran b'\ran. 
$$
 The counit map $\varepsilon\:C_0\rarrow K$ is obtained by applying
the functor $\Hom_{K^\rop}({-},K)$ to the ring homomorphism
$K\rarrow B^0$,
$$
 \varepsilon(c)=\lan c,1\ran \quad\text{for all $c\in C$},
$$
where $1\in B^0$ is the unit element.

 Generalizing Remark~\ref{B-sharp-remark} again, we put
$\shB^n=\Hom_{K^\rop}(C_n,K)$ for every $n\ge0$.
 The related pairing notation is
\begin{gather*}
 \lan\shb,c\ran=\shb(c) \quad
 \text{for all $\shb\in\shB^n$, \,$c\in C_n$, \,$n\ge0$}, \\
 \lan\shb,c\ran=0 \quad
 \text{for all $\shb\in\shB^i$, \,$c\in C_j$, \,$i\ne j$}.
\end{gather*}
 Then $\shB^n$ is a $K$\+$K$\+bimodule with the left and right
actions of $K$ given by the usual rules
$$
 \lan k\shb,c\ran=k\lan\shb,c\ran, \quad
 \lan\shb k,c\ran=\lan\shb,kc\ran \quad
 \text{for all $k\in K$, \,$\shb\in\shB^n$, \,$c\in C_n$}.
$$
 Moreover, $\shB=\bigoplus_{n=0}^\infty\shB^n$ is a nonnegatively
graded ring with the multiplication given by the rule
$$
 \lan\shb'\shb'',c\ran=\lan\shb'\ot\shb'',\,\mu(c)\ran \quad
 \text{for all $\shb'$, $\shb''\in\shB$, \ $c\in C$},
$$
where the pairing
$\lan\ ,\ \ran:\shB^j\ot_K\shB^i\ot_KC_i\ot_KC_j\rarrow K$ is
defined by
$$
 \lan\shb''\ot\shb',\,c'\ot c''\ran=
 \lan\shb''\lan\shb',c'\ran,c''\ran=
 \lan\shb'',\lan\shb',c'\ran c''\ran.
$$
 Applying the functor $\Hom_{K^{\rop}}({-},K)$ to the counit map
$C_0\rarrow K$, we obtain a ring homomorphism $K\rarrow\shB^0$
(which, together with the graded ring structure on $\shB$, induces
the $K$\+$K$\+bimodule structures on $\shB^n$ mentioned above).
 In the pairing notation, the formula
$$
 \lan 1,c\ran =\varepsilon(c) \quad
 \text{for all $c\in C$}
$$
defines the unit element $1\in\shB^0$.

 Furthermore, the graded left $K$\+module structure on
$C=\bigoplus_{n\ge0}C^n$ extends naturally to a graded left
$\shB$\+module structure, making $C$ a graded $\shB$\+$B$\+bimodule
(as it was mentioned already in Remark~\ref{B-sharp-remark}).
 In the pairing notation, the right action of $B$ in $C$ (defined
above) can be expressed in terms of the comultiplication in $C$ by
the formula
$$
 cb=\lan \mu_{i,j}(c),b\ran \quad
 \text{for all $c\in C_{i+j}$ and $b\in B^j$},
$$
where $\lan\ ,\ \ran$ denotes the map $C_i\ot_K C_j\ot_K B^j\rarrow C_i$
induced by the pairing $C_j\ot_K B^j\rarrow K$.
 Similarly, the left action of $\shB$ in $C$ is defined by the formula
$$
 \shb c=\lan\shb,\mu_{j,i}(c)\ran \quad
 \text{for all $c\in C_{i+j}$ and $\shb\in\shB^j$},
$$
where $\lan\ ,\ \ran$ denotes the map $\shB^j\ot_KC_j\ot_KC_i
\rarrow C_i$ induced by the pairing $\shB^j\ot_K C_j\rarrow K$.

 By construction, we know that $C_n$ is a finitely generated projective
left $K$\+module for every $n\ge0$.
 Now we make an additional assumption.
 
\begin{prop} \label{sharp-cdg-structures-correspondence}
 Assume additionally that $C_n$ is a finitely generated projective
right $K$\+module for every $n\ge0$.
 Then there is a natural bijective correspondence between
CDG\+ring structures $(B,d,h)$ on the graded ring $B$, viewed up to
change-of-connection isomorphisms $(id,a)\:(B,d',h')\rarrow(B,d,h)$,
and CDG\+ring structures $(\shB,\shd,\shh)$ on the graded ring\/
$\shB$, viewed up to change-of-connection isomorphisms $(id,\sha)\:
(\shB,\shd',\shh')\rarrow(\shB,\shd,\shh)$.
\end{prop}

\begin{proof}
 We use the result of Theorem~\ref{cdg-qdg-equivalence} and
identify CDG\+ring structures $(B,d,h)$ on $B$ with quasi-differential
structures $(\hB,\d)$ (that is, quasi-differential graded rings
$(\hB,\d)$ endowed with an injective graded ring homomorphism
$B\rarrow\hB$ with the image equal to $\ker\d\subset\hB$).
 Similarly, we identify CDG\+ring structures $(\shB,\shd,\shh)$
on the graded ring $\shB$ with quasi-differential structures
$(\shhB,\d)$ on~$\shB$, and construct a bijection between
the isomorphism classes of $(\hB,\d)$ and~$(\shhB,\d)$.

 Let $(\hB,\d)$ be a quasi-differential structure on~$B$.
 Recall that the differential $\d\:\hB^n\rarrow\hB^{n-1}$ is a morphism
of $B^0$\+$B^0$\+bimodules (since $\d(B^0)\subset B^{-1}=0$), hence also
a morphism of $K$\+$K$\+bimodules.
 Therefore, in view of the short exact sequences $0\rarrow B^n\rarrow
\hB^n\overset\d\rarrow B^{n-1}\rarrow0$, the graded right $K$\+module
$\hB^n$ is also finitely generated and projective for every $n\ge0$.

 Applying the above construction to the graded ring $\hB$, we produce
a graded coring $\hC$ with the components
$\hC_n=\hC^{-n}=\Hom_{K^\rop}(\hB^n,K)$.
 Moreover, the differential $\d\:\hB^n\rarrow\hB^{n-1}$ induces
a differential $\d\:\hC_{n-1}\rarrow\hC_n$, \ $n\ge0$, which we define
with the sign rule so that
$$
 \lan\d(c),b\ran+(-1)^{|c|}\lan c,\d(b)\ran=0 \quad
 \text{for all $c\in \hC^{|c|}$ and $b\in\hB^{|b|}$}.
$$
 By construction, $\d\:\hC\rarrow\hC$ is a homogeneous
$K$\+$K$\+bimodule map.
 As explained in the second proof of Theorem~\ref{pbw-theorem-thm}
in Section~\ref{pbw-theorem-subsecn}, it is also an odd coderivation
of the coring $\hC$ (of degree~$-1$ in the upper indices).
 Since $\d^2=0$ on $\hB$, we also have $\d^2=0$ on~$\hC$.
 The inclusion of graded rings $B\rarrow\hB$ induces a graded coring
homomorphism $\hC\rarrow C$, whose composition with the odd
coderivation~$\d$ on $\hC$ obviously vanishes.

 Furthermore, since the short exact sequences of $K$\+$K$\+bimodules
$0\rarrow B^n\rarrow\hB^n\overset\d\rarrow B^{n-1}\rarrow0$ are split
exact as short sequences of right $K$\+modules, we have short
exact sequences of $K$\+$K$\+bimodules
$0\rarrow C_{n-1}\rarrow\hC_n\rarrow C_n\rarrow0$.
 So the graded coring homomorphism $\hC\rarrow C$ is surjective
and its kernel is equal to the image of the odd coderivation~$\d$
on~$C$, as well as to the kernel of the odd coderivation~$\d$.

 Continuing to apply the above construction, we produce a graded ring
$\shhB$ with the components $\shhB^n=\Hom_{K^\rop}(\hC_n,K)$.
 The differential $\d\:\hC_{n-1}\rarrow\hC_n$ induces a differential
$\d\:\shhB^n\rarrow\shhB^{n-1}$, which we define using the sign rule
$$
 \lan\d(\shb),c\ran+(-1)^{|\shb|}\lan\shb,\d(c)\ran=0 \quad
 \text{for all $\shb\in\shhB^{|\shb|}$ and $c\in\hC^{|c|}$}.
$$
 The map $\d\:\shhB\rarrow\shhB$ is an odd derivation of degree~$-1$
on the graded ring~$\shhB$.
 Since $\d^2=0$ on $\hC$, we also have $\d^2=0$ on~$\shhB$.
 The surjective homomorphism of graded corings $\hC\rarrow C$ induces
an injective homomorphism of graded rings $\shB\rarrow\shhB$, whose
composition with the coderivation~$\d$ on $\shhB$ vanishes.

 By the additional assumption, all the $K$\+$K$\+bimodules in the short
exact sequence $0\rarrow C_{n-1}\rarrow\hC_n\rarrow C_n\rarrow0$ are
finitely generated and projective as right $K$\+modules, so the sequence
splits as a short exact sequence of right $K$\+modules.
 Hence we have short exact sequnces of $K$\+$K$\+bimodules
$0\rarrow\shB^n\rarrow\shhB^n\overset\d\rarrow\shB^{n-1}\rarrow0$.
 Thus the kernel of the odd derivation~$\d$ on~$\shhB$ is equal to
the image of the injective ring homomorphism $\shB\rarrow\shhB$, as
well as the to the image of the odd derivation $\d\:\shhB\rarrow
\shhB$.

 We have constructed the quasi-differential graded ring $(\shhB,\d)$
corresponding to a quasi-differential graded ring $(\hB,\d)$.
 To check that this is a bijective correspondence, it suffices to
reverse the construction, recovering the graded coring $\hC$ with
its odd coderivation~$\d$ as the graded Hom bimodule
$\hC=\Hom_K(\shhB,K)$, and the graded ring $\hB$ with its odd
derivation~$\d$ as the graded Hom group/bimodule $\hB=\Hom_K(\hC,K)$.
\end{proof} 
 
\begin{lem} \label{diagonal-cdg-bicomodule-lemma}
 In the context of
Proposition~\ref{sharp-cdg-structures-correspondence}, for any pair of
CDG\+ring structures $B=(B,d,h)$ and\/ $\shB=(\shB,\shd,\shh)$
corresponding to each other under the construction of
the proposition, there is a natural structure of CDG\+bimodule over\/
$\shB$ and $B$ on the graded\/ $\shB$\+$B$\+bimodule~$C$.
\end{lem}

\begin{proof}
 Saying that $B=(B,d,h)$ and $\shB=(\shB,\shd,\shh)$ correspond
to each other under the construction of
Proposition~\ref{sharp-cdg-structures-correspondence} means
the following.
 Let $(\hB,\d)$ and $(\shhB,\d)$ be a pair of quasi-differential graded
rings corresponding to each other under the construction from
the proof of the proposition.
 (This presumes that $\hB$ and $\shhB$ are nonnegatively graded, and
that we are given ring homomorphisms $K\rarrow\hB^0$ and
$K\rarrow\shhB^0$ such that $\hB^n$ is a finitely generated projective
right $K$\+module and $\shhB^n$ is a finitely generated projective
left $K$\+module for every $n\ge0$.)

 Choose arbitrary (unrelated) elements $\delta\in\hB^1$ such that
$\d(\delta)=1$ in $\hB^0$ and $\shdelta\in\shhB^1$ such that
$\d(\shdelta)=1$ in $\shhB^0$.
 Define the differential~$d$ on $B=\ker\d\subset\hB$ and the curvature
element $h\in B^2$ as in the proof of Theorem~\ref{cdg-qdg-equivalence},
$$
 d(b)=[\delta,b] \quad\text{for all $b\in B$},
 \quad\text{and}\quad h=\delta^2,
$$
where the bracket $[\ ,\ ]$ denotes the graded commutator.
 Similarly, consider the graded ring $\shB=\ker\d\subset\shhB$, and put
$$
 \shd(\shb)=[\shdelta,\shb] \quad
 \text{for all $\shb\in\shB$},
 \quad\text{and}\quad \shh=\shdelta^2.
$$

 In this context, we define the differential $d_{\hC}$ on the graded
$\shhB$\+$\hB$\+bimodule $\hC$ by the rule
$$
 d_{\hC}(c)=\shdelta c-(-1)^{|c|}c\delta \quad
 \text{for all $c\in\hC$}.
$$
 The claim is that the map~$d_{\hC}\:\hC\rarrow\hC$ induces
a well-defined map $d_C\:C\rarrow C$.
 In other words, there exists an (obviously unique) homogeneous map
$d_C\:C\rarrow C$ (of degree~$1$ in the upper indices) forming
a commutative square diagram with the map $d_{\hC}$ and the surjective
coring homomorphism $\hC\rarrow C$.

 In order to prove this assertion, we will check that $d_{\hC}$
preserves the kernel $\d(C)$ of the surjective map $\hC\rarrow C$,
that is $d_{\hC}(\d(C))\subset\d(C)$.
 Moreover, we will show that the two differentials $\d$ and $d_{\hC}$
on $\hC$ anti-commute, that is
\begin{equation} \label{d-and-delta-anti-commute}
 \d\circ d_{\hC}+d_{\hC}\circ\d=0 \quad \text{on $\hC$}.
\end{equation}

 First of all, let us check that the differential~$\d$ on the graded
right $\hB$\+module $\hC$ is an odd derivation compatible with
the odd derivation~$\d$ on $\hB$, that is
$$
 \d(cb)=\d(c)b+(-1)^{|c|}c\d(b)\in\hC \quad
 \text{for all $c\in\hC^{|c|}$ and $b\in\hB^{|b|}$}.
$$
 Indeed, for any $b'\in\hB$ we have
\begin{multline*}
 \lan \d(cb),b'\ran=-(-1)^{|c|+|b|}\lan cb,\d(b')\ran=
 -(-1)^{|c|+|b|}\lan c,b\d(b')\ran \\
 =-(-1)^{|c|}\lan c,\d(bb')\ran+
 (-1)^{|c|}\lan c,\d(b)b'\ran \\
 =\lan \d(c),bb'\ran+
 (-1)^{|c|}\lan c,\d(b)b'\ran
 =\lan \d(c)b,b'\ran+
 (-1)^{|c|}\lan c\d(b),b'\ran.
\end{multline*}
 Similarly one can check that the differential~$\d$ on the graded
left $\shhB$\+module $\hC$ is also an odd derivation compatible with
the odd derivation~$\d$ on $\shhB$, that is
$$
 \d(\shb c)=\d(\shb)c+(-1)^{|\shb|}\shb\d(c)\in\hC \quad
 \text{for all $\shb\in\shhB^{|\shb|}$ and $c\in\hC^{|c|}$}.
$$

 Now for any $c\in\hC$ we have
$$
 \d(\shdelta c)=\d(\shdelta)c-\shdelta\d(c)=c-\shdelta\d(c)
$$
and
$$
 \d(c\delta)=\d(c)\delta+(-1)^{|c|}c\d(\delta)=
 (-1)^{|c|}c+\d(c)\delta,
$$
hence
$$
 \d(d_{\hC}(c))=\d(\shdelta c-(-1)^{|c|}c\delta)
 =-\shdelta\d(c)-(-1)^{|c|}\d(c)\delta=
 -d_{\hC}(\d(c)),
$$
and the equation~\eqref{d-and-delta-anti-commute} is deduced.

 Thus the desired differential $d_C\:C\rarrow C$ is well-defined.
 Checking that $d_C$~is an odd derivation on the graded left
$\shB$\+module $C$ and the graded right $B$\+module $C$ compatible
with the odd derivations $\shd=[\shdelta,{-}]$ on $\shB$
and $d=[\delta,{-}]$ on $B$ is easy.
 The same applies to checking the equation
$d_C^2(c)=\shh c-ch$ for all $c\in C$.
\end{proof}

\subsection{Comodule-contramodule correspondence}
 We formulate two lemmas of rather general character before specializing
to the situation we are interested in.

\begin{lem} \label{coinduced-hom-induced-contratensor}
 Let $B=\bigoplus_{n=0}^\infty B^n$ be a nonnegatively graded ring,
and let $K$ be an associative ring endowed with a ring homomorphism
$K\rarrow B^0$.
 Then \par
\textup{(a)} For any graded left $K$\+module $L$ and any graded left
$B$\+comodule $M$, the graded abelian group of graded left $B$\+comodule
(or graded left $B$\+module) Hom from $M$ to the coinduced graded left
$B$\+comodule\/ $\Hom_K^\Sigma(B,L)$ can be computed as
$$
 \Hom_B(M,\Hom_K^\Sigma(B,L))\simeq\Hom_K(L,M).
$$ \par
\textup{(b)} For any graded left $K$\+module $L$ and any graded right
$B$\+comodule $N$, the contratensor product of $N$ with the induced
graded left $B$\+contramodule $\boM_K^\gr(L)=B\ot_K^\Pi L$ can be
computed as
$$
 N\ocn_B(B\ot_K^\Pi L)\simeq N\ot_KL.
$$
\end{lem}

\begin{proof}
 The natural isomorphism of graded abelian groups in part~(a) follows
from the similar natural isomorphism for the graded Hom group of
graded $B$\+modules, $\Hom_B(M,\Hom_K(B,L))\simeq\Hom_K(L,M)$,
together with Lemma~\ref{gamma-hom-hom-sigma}.

 To construct the natural isomorphism of graded abelian groups in
part~(b), let $U$ be a graded abelian group.
 Then we have natural isomorphisms
\begin{multline*}
 \Hom_\boZ(N\ocn_B(B\ot_K^\Pi L),\>U)
 \,\overset{\ref{graded-contratensor-hom-adjunction}}\simeq\,
 \Hom^B(B\ot_K^\Pi L,\>\Hom(N,U)) \\
 \,\overset{\ref{delta-tensor-tensor-pi}}\simeq\,
 \Hom_K(L,\Hom(N,U))\simeq\Hom_K(N\ot_KL,\>U),
\end{multline*}
where the first isomorphism is provided by
Lemma~\ref{graded-contratensor-hom-adjunction} (applied to the graded
ring $E=\boZ$ concentrated in degree~$0$), while the second one
holds in view of the discussion in Section~\ref{monads-subsecn} or
by Lemma~\ref{delta-tensor-tensor-pi}.
 It remains to say that any isomorphism of representable functors
comes from a unique isomorphism of the representing objects.
 (In fact, the isomorphism we have constructed is provided by the map
$N\ot_KL\rarrow N\ocn_B(B\ot_K^\Pi L)$ induced by the unit element
of~$B$.)
\end{proof}

\begin{lem} \label{cdg-hom-contratensor-functors}
 Let $E=(E,d_E,h_E)$ be a CDG\+ring and $B=(B,d,h)$ be a nonnegatively
graded CDG\+ring. 
 Let $N=(N,d_N)$ be a CDG\+bimodule over $E$ and $B$ whose underlying
graded right $B$\+module is a graded right $B$\+comodule.
 Then \par
\textup{(a)} for any left CDG\+module $M=(M,d_M)$ over~$E$, the graded
left $B$\+contramodule $\Hom_E(N,M)$ constructed in
Example~\ref{graded-co-contra-dualization-example} endowed with
the differential of the left CDG\+module $\Hom_E(N,M)$ over $B$
constructed in Section~\ref{cdg-modules-subsecn} is a left
CDG\+contra\-module over $B$ in the sense of the definition
in Section~\ref{contraderived-cdg-contramodules-subsecn}. \par
\textup{(b)} for any left CDG\+contramodule $(Q,d_Q)$ over~$B$,
the differential of the left CDG\+module $N\ot_BQ$ over $E$ constructed
in Section~\ref{cdg-modules-subsecn} induces a well-defined differential
on the graded quotient $E$\+module $N\ocn_BQ$ of $N\ot_BQ$ constructed
in Section~\ref{contratensor-subsecn}, making $N\ocn_BQ$
a left CDG\+module over~$E$.
\end{lem}

\begin{proof}
 The proofs of both~(a) and~(b) are straightforward.
\end{proof}

 Now we assume that $B=(B,d,h)$ and $\shB=(\shB,\shd,\shh)$
are two CDG\+rings corresponding to each other under the construction
of Proposition~\ref{sharp-cdg-structures-correspondence}.
 This means, in particular, that $B=\bigoplus_{n=0}^\infty B^n$
and $\shB=\bigoplus_{n=0}^\infty \shB^n$ are nonnegatively graded,
that we are given a ring $K$ together with ring homomorphisms
$K\rarrow B^0$ and $K\rarrow\shB^0$, that the right $K$\+module $B^n$
is finitely generated and projective for every $n\ge0$, and that
the left $K$\+module $\shB^n$ is finitely generated and projective
for every $n\ge0$.

 The following theorem is the first main result of
Section~\ref{co-contra-secn}.
 It is a generalization of~\cite[Theorem~B.3]{Pkoszul}.
 It also solves a particular case of the problem posed
in~\cite[Question at the end of Section~11]{Psemi}.

\begin{thm} \label{co-contra-correspondence-theorem}
 Assume additionally that the left homological dimension of the ring
$K$ is finite.
 Then there is a natural triangulated equivalence between the coderived
category of left CDG\+comodules over\/ $\shB=(\shB,\shd,\shh)$
and the contraderived category of left CDG\+contramodules over
$(B,d,h)$,
\begin{equation} \label{co-contra-correspondence-eqn}
 \boR\Hom_{\shB}(C,{-})\:\sD^\co(\shB\comodl)
 \simeq\sD^\ctr(B\contra):\!C\ocn_B^\boL{-},
\end{equation}
provided by derived functors of comodule Hom and contratensor product
with the CDG\+bi(co)module $(C,d_C)$ over\/ $\shB$ and $B$ constructed
in Lemma~\ref{diagonal-cdg-bicomodule-lemma}.
\end{thm}

\begin{proof}
 First of all, the CDG\+bimodule $(C,d_C)$ is nonpositively
cohomologically graded.
 According to the discussion in Section~\ref{graded-comodules-subsecn},
it follows that $C$ is a graded left $\shB$\+comodule and a graded right
$B$\+comodule.

 The DG\+functor
\begin{equation} \label{contra-hom-dg-functor}
 \Hom_{\shB}(C,{-})\:\DG(\shB\comodl)\lrarrow\DG(B\contra)
\end{equation}
assigns to a left CDG\+comodule $M$ over $(\shB,\shd,\shh)$
the graded left $B$\+contramodule $\Hom_{\shB}(C,M)$, as constructed
in Example~\ref{graded-co-contra-dualization-example}, endowed with
the structure of a left CDG\+contramodule over $(B,d,h)$ described
in Lemma~\ref{cdg-hom-contratensor-functors}(a).
 Since $\Hom_{\shB}(C,{-})$ is a DG\+functor $\DG(\shB\modl)\rarrow
\DG(B\modl)$ according to the discussion in
Section~\ref{cdg-modules-subsecn}, it follows that it is also
a DG\+functor $\DG(\shB\comodl)\rarrow\DG(B\contra)$.

 The DG\+functor
\begin{equation} \label{contratensor-dg-functor}
 C\ocn_B{-}\:\DG(B\contra)\lrarrow\DG(\shB\comodl)
\end{equation}
assigns to a left CDG\+contramodule $Q$ over $(B,d,h)$ the graded
abelian group of the contratensor product $C\ocn_BQ$, as constructed in
Section~\ref{contratensor-subsecn}, endowed with the obvious graded
left $\shB$\+comodule structure induced by the graded left
$\shB$\+comodule structure on $C$ and with the structure of
a left CDG\+comodule over $(\shB,\shd,\shh)$ described in
Lemma~\ref{cdg-hom-contratensor-functors}(b).
 Since $C\ot_B{-}$ is a DG\+functor $\DG(B\modl)\rarrow\DG(\shB\modl)$
according to the discussion in Section~\ref{cdg-modules-subsecn},
it follows that $C\ocn_B{-}$ is a DG\+functor
$\DG(B\contra)\rarrow\DG(\shB\comodl)$.

 The adjunction of Lemma~\ref{graded-contratensor-hom-adjunction}
makes the DG\+functors \eqref{contra-hom-dg-functor}
and~\eqref{contratensor-dg-functor} adjoint to each other.
 Hence the induced triangulated functors between the homotopy categories
\begin{equation}\label{contra-hom-homotopy-functor}
 \Hom_{\shB}(C,{-})\:\Hot(\shB\comodl)\lrarrow\Hot(B\contra)
\end{equation}
and
\begin{equation} \label{contratensor-homotopy-functor}
 C\ocn_B{-}\:\Hot(B\contra)\lrarrow\Hot(\shB\comodl)
\end{equation}
are also adjoint.

 In order to construct the derived functor $\boR\Hom_{\shB}(C,{-})\:
\sD^\co(\shB\comodl)\rarrow\sD^\ctr(B\contra)$, we restrict
the triangulated functor $\Hom_{\shB}(C,{-})$ 
\,\eqref{contra-hom-homotopy-functor} to the full triangulated
subcategory of CDG\+comodules of the coinduced type
$$
 \Hot(\shB\comodl_{K\coind})\subset\Hot(\shB\comodl),
$$
which was defined in Section~\ref{cdg-comodules-coinduced-type-subsecn}.
 Similarly, in order to construct the derived functor
$C\ocn_B^\boL{-}\:\sD^\ctr(B\contra)\rarrow\sD^\co(B\comodl)$, we
restrict the triangulated functor $C\ocn_B{-}$
\,\eqref{contratensor-homotopy-functor} to the full triangulated
subategory of CDG\+contramodules of the induced type
$$
 \Hot(B\contra_{K\ind})\subset\Hot(B\contra),
$$
which was defined in
Section~\ref{cdg-contramodules-induced-type-subsecn}.

 Now it is claimed that the functors $\Hom_{\shB}(C,{-})$ and
$C\ocn_B{-}$ take the full subcategories
$\Hot(\shB\comodl_{K\coind})\subset\Hot(\shB\comodl)$ and
$\Hot(B\contra_{K\ind})\subset\Hot(B\contra)$ into each other
and induce mutually inverse triangulated equivalences between them
\begin{equation} \label{induced-coinduced-homotopy-equivalence}
 \Hom_{\shB}(C,{-})\:\Hot(\shB\comodl_{K\coind})
 \simeq\Hot(B\contra_{K\ind}):\!C\ocn_B{-}.
\end{equation}
 Moreover, we have mutually inverse DG\+equivalences
$$
 \Hom_{\shB}(C,{-})\:\DG(\shB\comodl_{K\coind})
 \simeq\DG(B\contra_{K\ind}):\!C\ocn_B{-}.
$$

 It suffices to check the latter assertion on the level of the full
subcategory of coinduced graded comodules $\shB\comodl_{K\coind}
\subset\shB\comodl$ in the abelian category of graded left
$\shB$\+comodules $\shB\comodl$ and the full subcategory of induced
graded contramodules $B\contra_{K\ind}\subset B\contra$ in the abelian
category of graded left contramodules $B\contra$.
 So we need to check that the adjunction of
Lemma~\ref{graded-contratensor-hom-adjunction} restricts to a pair
of mutually inverse equivalences of additive categories
$$
 \Hom_{\shB}(C,{-})\:\shB\comodl_{K\coind}
 \simeq B\contra_{K\ind}:\!C\ocn_B{-}.
$$
 This amounts to a simple computation based on
Lemma~\ref{coinduced-hom-induced-contratensor}.
 To show that two adjoint functors are mutually inverse equivalences,
it suffices to check that the adjunction morphisms (i.~e.,
the adjunction unit and counit) are isomorphisms.

 Let $L$ be a graded left $K$\+module, $\Hom_K^\Sigma(\shB,L)$ be
the graded left $\shB$\+comodule coinduced from $L$, and
$B\ot_K^\Pi L$ be the graded left $B$\+contramodule induced from~$L$.
 Then we have
$$
 \Hom_{\shB}(C,\Hom_K^\Sigma(\shB,L))
 \,\,\overset{\ref{coinduced-hom-induced-contratensor}\text{(a)}}
 \simeq\,\,\Hom_K(C,L)\,\simeq\,\Hom_K(C,K)\ot_K^\Pi L\,=\,B\ot_K^\Pi L,
$$
since $C^{-n}=\Hom_{K^\rop}(B^n,K)$ is a finitely generated projective
left $K$\+module for every $n\ge0$.
 Similarly,
$$
 C\ocn_B(B\ot_K^\Pi L)
 \,\,\overset{\ref{coinduced-hom-induced-contratensor}\text{(b)}}
 \simeq\,\,C\ot_KL\,\simeq\,\Hom_K^\Sigma(\Hom_{K^\rop}(C,K),L)\,=\,
 \Hom_K^\Sigma(\shB,L),
$$
since $C^{-n}=\Hom_K(\shB^n,K)$ is a finitely generated projective
right $K$\+module for every $n\ge0$.
 This finishes the proof of the triangulated
equivalence~\eqref{induced-coinduced-homotopy-equivalence}.

 Finally, in order to show that the triangulated
equivalence~\eqref{induced-coinduced-homotopy-equivalence} descends
to the desired triangulated
equivalence~\eqref{co-contra-correspondence-eqn}, we use
the descriptions of the coderived category $\sD^\co(\shB\comodl)$
and the contraderived category $\sD^\ctr(B\contra)$ provided by
Theorems~\ref{coderived-category-theorem}
and~\ref{contraderived-category-theorem} together with
Propositions~\ref{coderived-coinduced-prop}
and~\ref{contraderived-induced-prop}.

 Specifically, by Proposition~\ref{coderived-coinduced-prop}
and the second assertion of Theorem~\ref{coderived-category-theorem},
the coderived category $\sD^\co(\shB\comodl)$ is equivalent to
the triangulated quotient category of $\Hot(\shB\comodl_{K\coind})$
by its minimal triangulated subcategory containing all the totalizations
of left CDG\+comodule short exact sequences of the coinduced type
over $(\shB,\shd,\shh)$.
 Similarly, by Proposition~\ref{contraderived-induced-prop} and
the second assertion of Theorem~\ref{contraderived-category-theorem},
the contraderived category $\sD^\ctr(B\contra)$ is equivalent to
the triangulated quotient category of $\Hot(B\contra_{K\ind})$
by its minimal triangulated subcategory containing all the totalizations
of left CDG\+contramodule short exact sequences of the induced type
over $(B,d,h)$.

 In view of these results, it remains to observe that the DG\+functors
$\Hom_{\shB}(C,{-})$ and $C\ocn_B{-}$ take CDG\+comodule short exact
sequences of the coinduced type over $(\shB,\shd,\shh)$ to
CDG\+contramodule short exact sequences of the induced bype over
$(B,d,h)$ and back.
 This is also an assertion on the level of the additive categories
$\shB\comodl_{K\coind}$ and $B\contra_{K\ind}$, and it follows from
the above computation based on
Lemma~\ref{coinduced-hom-induced-contratensor}.
\end{proof}

 Let us recall that, whenever the forgetful functor
$B\contra_\sgr\rarrow B\modl_\sgr$ is fully faitfhul,
the contratensor product functor $C\ocn_B{-}$ appearing in
the above theorem is isomorphic to the tensor product functor
$C\ot_B{-}$ by Proposition~\ref{contratensor-tensor-comparison}(b).
 In particular, by
Theorem~\ref{graded-contramodules-fully-faithful-thm}, this holds
whenever the augmentation ideal $B^{\ge1}=\bigoplus_{n=1}^\infty B^n$
is finitely generated as a right ideal in~$B$.

\subsection{Two-sided finitely projective Koszul rings}
\label{two-sided-finitely-projective-subsecn}
 Let $A=\bigoplus_{n=0}^\infty A_n$ be a nonnegatively graded ring
with the degree-zero component $R=A_0$.
 We will say that $A$ is \emph{two-sided finitely projective Koszul}
if it is \emph{both} left finitely projective Koszul and right finitely
projective Koszul (in the sense of
Section~\ref{finitely-projective-Koszul-subsecn}).

 Let $A$ be a two-sided finitely projective Koszul graded ring.
 Denote by $B$ the right finitely projective Koszul graded ring
quadratic dual to the left finitely projective Koszul ring $A$,
and denote by $\shB$ the left finitely projective Koszul gradedring
quadratic dual to the right finitely projective Koszul ring $A$
(as in Section~\ref{quadratic-duality-secn} and
Proposition~\ref{finitely-projective-Koszul-duality}).
 So we have
$$
 B=\Ext_A(R,R)^\rop \quad\text{and}\quad \shB=\Ext_{A^\rop}(R,R),
$$
where $\Ext_A$ is taken in the category of graded left $A$\+modules
and $\Ext_{A^\rop}$ in the category of graded right $A$\+modules,
as usually.

\begin{lem} \label{koszul-dual-rings-two-sides}
 For any two-sided finitely projective Koszul graded ring $A$,
the nonnegatively graded rings $B$ and\/ $\shB$ are related to each
other by the construction of Remark~\ref{B-sharp-remark} and
Section~\ref{diagonal-cdg-bicomodule-subsecn} (with $K=R$).
 In other words, the grading components of $B$ and\/ $\shB$ can be
obtained from each other as the $R$\+$R$\+bimodules
$$
 \shB^n=\Hom_{R^\rop}(\Hom_{R^\rop}(B^n,R),R) \quad
 \text{and} \quad B^n=\Hom_R(\Hom_R(\shB^n,R),R)
$$
for every $n\ge0$, and the multiplicative structures are related
accordingly.
\end{lem}

\begin{proof}
 Put $V=A_1$, and denote by $I\subset V\ot_RV$ the kernel of
the multiplication map $A_1\ot_RA_1\rarrow A_2$.
 The graded $R$\+$R$\+bimodule $C$ with the components
$$
 C_n=\Tor^A_n(R,R)=\Tor^A_{n,n}(R,R)
$$
plays a key role.
 By Proposition~\ref{diagonal-Tor}(b), we have $C_0=R$, \ $C_1=V$, \
$C_2=I$, and
$$
 C_n=\bigcap\nolimits_{k=1}^{n-1}
 V^{\ot_R\,k-1}\ot_R I\ot_R V^{\ot_R\,n-k-1}
 \subset V^{\ot_R\,n}, \qquad n\ge2.
$$

 Since $A$ is a left finitely projective Koszul graded ring,
Theorem~\ref{projective-koszul-theorem}(c) tells, in particular,
that $C_n$ is a finitely generated projective left $R$\+module
for every $n\ge0$.
 As $A$ is also a right finitely projective Koszul graded ring,
$C_n$ is also a finitely generated projective right $R$\+module.

 According to Remark~\ref{bar-complex-dg-coring}, $C$ has a natural
graded coring structure.
 The grading components of the comultiplication map
$\mu_{i,j}\:C_{i+j}\rarrow C_i\ot_RC_j$, \ $i$, $j\ge0$, can be
described explicitly as the $R$\+$R$\+subbimodule inclusions
\begin{multline*}
 C_{i+j}=\bigcap\nolimits_{k=1}^{i+j-1}
 V^{\ot_R\,k-1}\ot_R I\ot_R V^{\ot_R\,i+j-k-1} \\ \lrarrow
 \bigcap\nolimits_{1\le k\le i+j-1}^{k\ne i}
 V^{\ot_R\,k-1}\ot_R I\ot_R V^{\ot_R\,i+j-k-1}
 \,\simeq\,C_i\ot_RC_j,
\end{multline*}
where the latter isomorphism holds by
Lemma~\ref{tensor-intersections}(b,d,e).

 Computing the Ext in terms of the relative cobar complex as in
Proposition~\ref{diagonal-Ext}, one obtains the isomorphisms
$$
 B^n=\Ext^n_A(R,R)=\Ext^{n,n}_A(R,R)\simeq\Hom_R(C_n,R)
$$
and
$$
 \shB^n=\Ext^n_{A^\rop}(R,R)=\Ext^{n,n}_{A^\rop}(R,R)
 \simeq\Hom_{R^\rop}(C_n,R),
$$
and these isomorphisms are compatible with the multiplicative
structures.
\end{proof}

 Now let $(\tA,F)$ be an associative ring with an increasing filtration
$0=F_{-1}\tA\subset R=F_0\tA\subset\tV=F_1\tA\subset F_2\tA\subset
\dotsb$.
 We recall that the filtered ring $(\tA,F)$ is said to be left finitely
projective nonhomogeneous Koszul if the graded ring
$A=\gr^F\tA=\bigoplus_{n=0}^\infty F_n\tA/F_{n-1}\tA$ is left
finitely projective Koszul (see Section~\ref{pbw-theorem-subsecn}).

 Dually, let us say that $(\tA,F)$ is \emph{right finitely projective
nonhomogeneous Koszul} if $(\tA^\rop,F^\rop)$ is left finitely
projective nonhomogeneous Koszul, that is, in other words, the graded
ring $A$ is right finitely projective Koszul.
 We will say that $(\tA,F)$ is \emph{two-sided finitely projective
nonhomogeneous Koszul} if it is both left and right finitely projective
nonhomogeneous Koszul; in other words, this means that the graded
ring $A$ is two-sided finitely projective Koszul.

 Let $B$ and $\shB$ be the Koszul dual rings to $A$ on the two sides,
as in Lemma~\ref{koszul-dual-rings-two-sides}.
 Recall the notation~$\tau'$ introduced in
Section~\ref{koszul-duality-functors-comodule-side-subsecn} for
a left $R$\+linear splitting
$$
 \Hom_{R^\rop}(B^1,R)\simeq F_1\tA/F_0\tA=V\,\hookrightarrow\,
 \tV=F_1\tA
$$
of the surjective $R$\+$R$\+bimodule map
$\tV=F_1\tA\rarrow F_1\tA/F_0\tA=V$.
 Such a splitting was first considered (under the assumption of
projectivity of the left $R$\+module $V$, guaranteeing its
existence) in Section~\ref{self-consistency-subsecn}.

 For a two-sided finitely projective nonhomogeneous Koszul ring $\tA$,
one can also choose a right $R$\+linear splitting $V\rarrow\tV$ of
the surjective $R$\+$R$\+bimodule map $\tV\rarrow\tV/R=V$.
 Denote by~$\sigma'$ the resulting right $R$\+module map
\begin{equation} \label{sigma-prime-explained-eqn}
 \Hom_R(\shB^1,R)\simeq F_1\tA/F_0\tA=V\hookrightarrow\tV=F_1\tA.
\end{equation}

 Having chosen a one-sided splitting~$\tau'$, one can apply
the construction of Sections~\ref{self-consistency-subsecn}\+-%
\ref{cdg-ring-constructed-subsecn} to the left finitely projective
nonhomogeneous Koszul ring $\tA$, and produce a CDG\+ring structure
$(B,d,h)$ on the right finitely projective Koszul graded ring~$B$.
 Having chosen a splitting~$\sigma'$ on the other side, one can
apply the opposite version of the same construction (with the left
and right sides switched) to the right finitely projective
nonhomogeneous Koszul ring $\tA$, and produce a CDG\+ring structure
$(\shB,\shd,\shh)$ on the left finitely projective Koszul
graded ring~$\shB$.

\begin{lem} \label{nonhomogeneous-koszul-dual-two-sides}
 For any two-sided finitely projective nonhomogeneous Koszul ring $\tA$,
the nonnegatively graded CDG\+rings $(B,d,h)$ and
$(\shB,\shd,\shh)$ are related by the construction of
Proposition~\ref{sharp-cdg-structures-correspondence} (with $K=R$).
\end{lem}

\begin{proof}
 We use the interpretation of nonhomogeneous quadratic duality in
terms of quasi-differential corings, as formulated in
Theorem~\ref{qdg-nonhomogeneous-duality-theorem}.
 The graded ring $\hA=\bigoplus_{n=0}^\infty\tA_n$ is left finitely
projective Koszul by Lemma~\ref{hA-left-projective} and
Theorem~\ref{central-nonzerodivisor-theorem}.
 By the opposite assertion, the graded ring $\hA$ is also right
finitely projective Koszul.

 Let $\hB$ be the right finitely projective Koszul graded ring
quadratic dual to the left finitely projective Koszul graded ring
$\hA$, and let $\shhB$ be the left finitely projective Koszul graded
ring quadratic dual to the right finitely projective Koszul graded
ring~$\hA$.
 By Lemma~\ref{koszul-dual-rings-two-sides} applied to the two-sided
finitely projective Koszul graded ring~$\hA$, the nonnegatively
graded rings $\hB$ and $\shhB$ are related to each other by
the construction from the beginning of
Section~\ref{diagonal-cdg-bicomodule-subsecn}.

 It remains to observe that the acyclic odd derivations~$\d$ of
degree~$-1$ on the graded rings $\hB$ and $\shhB$, both originating
from the canonical central element $t\in\hA_1$ via the construction
of Lemma~\ref{central-iff-derivation}(b) and its opposite version,
correspond to each other under the construction from the proof of
Proposition~\ref{sharp-cdg-structures-correspondence}.
\end{proof}

\subsection{Bimodule resolution revisited} \label{revisited-subsecn}
 Now, under the more restrictive assumption of two-sided (rather
than just left) finite projective Koszulity, we can offer an explanation
of the ``twisting cochain~$\sigma'$ notation'' of
Sections~\ref{koszul-duality-functors-comodule-side-subsecn}
and~\ref{koszul-duality-functors-contramodule-side-subsecn}.
 Let us start with the ``twisting cochain~$\sigma$'', which is
the homogeneous Koszul version (it was briefly mentioned in
the proofs of Theorems~\ref{comodule-side-koszul-duality-theorem}
and~\ref{contramodule-side-koszul-duality-theorem}).

 Let $A$ be a left finitely projective Koszul graded ring and $B$ be
the quadratic dual right finitely projective Koszul graded ring.
 In this context, in the discussion of homogeneous Koszul complexes
in Sections~\ref{first-koszul-complex-subsecn}\+-%
\ref{dual-koszul-complex-subsecn}, the symbol~$\tau$ was used as
a notation for the structure isomorphism of quadratic duality
$$
 \tau\:\Hom_{R^\rop}(B^1,R)\overset\simeq\lrarrow A_1.
$$
 Assume that $A$ is also a right finitely projective Koszul graded
ring, and let $\shB$ be the quadratic dual left finitely projective
Koszul graded ring.
 Denote by~$\sigma$ the structure isomorphism
$$
 \sigma\:\Hom_R(\shB^1,R)\overset\simeq\lrarrow A_1.
$$

 Moreover, in the context of the discussion of the dual Koszul
complex $\Ksp_e{}^\bu(B,A)=(B\ot_RA,\>d_e)$ in
Section~\ref{dual-koszul-complex-subsecn}, we denoted by
$e\in B^1\ot_RA_1$ the element corresponding to the map~$\tau$
under the construction of Lemma~\ref{element-e-lemma}.
 Let us now denote by $u\in A_1\ot_R\shB^1$ the element corresponding
to the map~$\sigma$ under the opposite version of the same
construction.
 Denote the related dual Koszul complex by $\Ksp_u{}^\bu(A,\shB)=
(A\ot_R\shB,\>d_u)$.

 For any left $A$\+module $M$ and any graded right $B$\+module $N$
(which we would prefer to assume additionally to be a $B$\+comodule)
we define the complex $N\ot_R^\tau M$ as the tensor product
$$
 N\ot_R^\tau M=N\ot_B\Ksp_e{}^\bu(B,A)\ot_AM
$$
with the differential induced by the differential~$d_e$ on
$\Ksp_e{}^\bu$.
 So $N\ot_R^\tau M$ is a complex of abelian groups whose underlying
graded abelian group is $N\ot_R M$.
 Additional structures on a left $A$\+module $M$ and a graded right
$B$\+comodule $N$ may induce additional structures on
the complex $N\ot_R^\tau M$, as usually.

 Similarly, for any right $A$\+module $M$ and any graded left
$\shB$\+module $N$ (which we would rather assume to be a graded left
$\shB$\+comodule) we define the complex $M\ot_R^\sigma N$ as
the tensor product
$$
 M\ot_R^\sigma N=M\ot_A\Ksp_u{}^\bu(A,\shB)\ot_{\shB}N
$$
with the differential induced by the differential~$d_u$ on
$\Ksp_u{}^\bu$.
 So $M\ot_R^\sigma N$ is a complex of abelian groups whose underlying
graded abelian group is $M\ot_RN$.

 Now we recall the notation $C$ for the graded $R$\+$R$\+bimodule
(in fact, a graded coring) $\Tor^A(R,R)$ which was used and
discussed in the proof of Lemma~\ref{koszul-dual-rings-two-sides}
and the references therein.
 In particular, according to Remark~\ref{B-sharp-remark}
and Section~\ref{diagonal-cdg-bicomodule-subsecn}, \,$C$ is a graded
$\shB$\+$B$\+bimodule (in fact, even a comodule on both sides, as
explained in the beginning of the proof of
Theorem~\ref{co-contra-correspondence-theorem}).

 With these preparations, we can interpret the second Koszul complex
${}^\tau\!K_\bu(B,A)=(\Hom_{R^\rop}(B,R)\ot_RA,\>{}^\tau\!\d)$
\,\eqref{second-koszul-complex} as
$$
 {}^\tau\!K_\bu(B,A)=C\ot_R^\tau A.
$$
 The first Koszul complex $K^\tau_\bu(B,A)=(\Hom_{R^\rop}(B,A),\d^\tau)
=(A\ot_R\Hom_{R^\rop}(B,R),\>\d^\tau)$ \,\eqref{first-koszul-complex}
gets interpreted as
$$
 K^\tau_\bu(B,A)=A\ot_R^\sigma C,
$$
as mentioned in the proofs of
Theorems~\ref{comodule-side-koszul-duality-theorem}
and~\ref{contramodule-side-koszul-duality-theorem}.
 Looking on $A$ as primarily a right finitely projective Koszul graded
ring rather than a left one, the first and second Koszul complexes
switch their roles, of course.

 Let us turn to the nonhomogeneous situation.
 Let $(\tA,F)$ be a two-sided finitely projective nonhomogeneous
Koszul ring, and let $(B,d,h)$ and $(\shB,\shd,\shh)$ be
its two nonhomogeneous quadratic dual CDG\+rings on the right and
left sides, as in Section~\ref{two-sided-finitely-projective-subsecn}.
 The notation~$\sigma'$ for the right $R$\+linear map between
$R$\+$R$\+bimodules $\Hom_R(\shB^1,R)\rarrow F_1\tA$ was already
introduced in~\eqref{sigma-prime-explained-eqn}.

 Denote by $u'\in F_1\tA\ot_R\shB^1$ the element
corresponding to~$\sigma'$ under the opposite version of
the construction of Lemma~\ref{element-e-lemma}.
 The opposite version of the construction of
Section~\ref{nonhomogeneous-koszul-cdg-module-subsecn}
produces the dual nonhomogeneous Koszul CDG\+module
$\Ksp(\tA,\shB)=\Ksp_{u'}(\tA,\shB)$, which has the form
$$
 0\lrarrow\tA\lrarrow\tA\ot_R\shB^1\lrarrow\tA\ot_R\shB^2\lrarrow\dotsb
$$
 The differential on $\Ksp_{u'}(\tA,\shB)$ does not square to
zero when $\shh\ne0$.
 Rather, $\Ksp_{u'}(\tA,\shB)$ is a right CDG\+module over
$(\shB,\shd,\shh)$, and in fact a CDG\+bimodule over
$(\tA,0,0)$ and $(\shB,\shd,\shh)$.
 The underlying graded $\tA$\+$\shB$\+bimodule of
$\Ksp_{u'}(\tA,\shB)$ is the bimodule tensor product $\tA\ot_R\shB$.

 Now let $D=(D,d_D,h_D)$ and $E=(E,d_E,h_E)$ be two CDG\+rings.
 Let $M$ be a CDG\+bimodule over $(\tA,0,0)$ and $(E,d_E,h_E)$, and
let $N$ be a CDG\+bimodule over $(D,d_D,h_D)$ and $(B,d,h)$.
 We would prefer to assume $N$ to be a graded right $B$\+comodule
rather than an arbitrary graded right $B$\+module.
 In this setting, we define the CDG\+bimodule $N\ot_R^{\tau'}M$
over $D$ and $E$ as the tensor product of CDG\+bimodules
$$
 N\ot_R^{\tau'}M=N\ot_B\Ksp_{e'}(B,\tA)\ot_{\tA}M.
$$
 The underlying graded $D$\+$E$\+bimodule of $N\ot_R^{\tau'}M$ is
the bimodule tensor product $N\ot_RM$.
 If $D$ is nonnegatively graded and $N$ is a $D$\+comodule, or if
$E$ is nonnegatively graded and $M$ is an $E$\+comodule, then
the tensor product $N\ot_RM$ is obviously a comodule on
the respective side.

 Similarly, let $M$ be a CDG\+bimodule over $(E,d_E,h_E)$ and
$(\tA,0,0)$, and let $N$ be a CDG\+bimodule over
$(\shB,\shd,\shh)$ and $(D,d_D,h_D)$.
 We would rather assume the graded left $\shB$\+module $N$ to be
a $\shB$\+comodule.
 Then we define the CDG\+bimodule $M\ot_R^{\sigma'}N$ over $E$ and $D$
as the tensor product of CDG\+bimodules
$$
 M\ot_R^{\sigma'}N=M\ot_{\tA}\Ksp_{u'}(\tA,\shB)\ot_{\shB}N.
$$
 The underlying graded $E$\+$D$\+bimodule of $M\ot_R^{\sigma'}N$
is the bimodule tensor product $M\ot_RN$.
 If $D$ is nonnegatively graded and $N$ is a $D$\+comodule, or if
$E$ is nonnegatively graded and $M$ is an $E$\+comodule, then
the tensor product $M\ot_RN$ is obviously a comodule on
the respective side.

 Now we recall that, by Lemma~\ref{diagonal-cdg-bicomodule-lemma}
(which is applicable to the situation at hand in view of
Lemma~\ref{nonhomogeneous-koszul-dual-two-sides}), there is
a natural differential~$d_C$ on the graded $\shB$\+$B$\+bimodule
$C$, making $(C,d_C)$ a CDG\+bimodule over $(\shB,\shd,\shh)$
and $(B,d,h)$.
 With these preparations, we can interpret the nonhomogeneous Koszul
duality functors
$$
 M^\bu\longmapsto M^\bu\ot_R^{\sigma'}C
 \quad\text{and}\quad
 N\longmapsto N\ot_R^{\tau'}\tA
$$ 
of Section~\ref{koszul-duality-functors-comodule-side-subsecn} as
particular cases of the above contsructions of twisted tensor products.
 This explains the meaning of the ``placeholder~$\sigma'$\,'' in
Section~\ref{koszul-duality-functors-comodule-side-subsecn}.

 Let $M$ be a CDG\+bimodule over $(\tA,0,0)$ and $(D,d_D,h_D)$, and
let $Q$ be a left CDG\+module over $(B,d,h)$.
 We would prefer to assume $Q$ to be a left CDG\+con\-tramodule
over $(B,d,h)$.
 In this setting, we define the left CDG\+module $\Hom_R^{\tau'}(M,Q)$
over $(D,d_D,h_D)$ as the Hom CDG\+module
$$
 \Hom_R^{\tau'}(M,Q)=\Hom_B(\Ksp_{e'}(B,\tA)\ot_{\tA}M,\>Q)
 =\Hom_{\tA}(M,\Hom_B(\Ksp_{e'}(B,\tA),Q)).
$$
 The underlying graded left $D$\+module of $\Hom_R^{\tau'}(M,Q)$ is
the Hom module $\Hom_R(M,Q)$ with the left $D$\+module structure
induced by the right $D$\+module structure on~$M$.
 If $D$ is nonnegatively graded and $M$ is a $D$\+comodule, then
$\Hom_R^{\tau'}(M,Q)$ is a left CDG\+contramodule over
$(D,d_D,h_D)$ by Lemma~\ref{cdg-hom-contratensor-functors}(a). {\hbadness=1950\par}

 Let $N$ be a CDG\+bimodule over $(\shB,\shd,\shh)$ and
$(D,d_D,h_D)$, and let $P=P^\bu$ be a complex of left $\tA$\+modules.
 We would rather assume the graded left $\shB$\+module $N$ to be
a $\shB$\+comodule.
 Then we define the left CDG\+module $\Hom_R^{\sigma'}(N,P^\bu)$
over $(D,d_D,h_D)$ as the Hom CDG\+module
$$
 \Hom_R^{\sigma'}(N,P^\bu)=\Hom_{\tA}(\Ksp_{u'}(\tA,\shB)\ot_{\shB}N,
 \>P^\bu)=\Hom_{\shB}(N,\Hom_{\tA}(\Ksp_{u'}(\tA,\shB),P^\bu)).
$$
 The underlying graded left $D$\+module of $\Hom_R^{\sigma'}(N,P^\bu)$
is the Hom module $\Hom_R(N,P)$ with the left $D$\+module
structure induced by the right $D$\+module structure on~$N$.
 If $D$ is nonnegatively graded and $N$ is a $D$\+comodule, then
$\Hom_R^{\sigma'}(N,P^\bu)$ is a left CDG\+contramodule over
$(D,d_D,h_D)$ by Lemma~\ref{cdg-hom-contratensor-functors}(a).
{\hbadness=1600\par}

 With these preparations, we can interpret the nonhomogeneous
Koszul duality functors
$$
 P^\bu\longmapsto\Hom_R^{\sigma'}(C,P^\bu)
 \quad\text{and}\quad
 Q\longmapsto\Hom_R^{\tau'}(\tA,Q)
$$
of Section~\ref{koszul-duality-functors-contramodule-side-subsecn}
as particular cases of the above constructions of twisted Hom
CDG\+modules.
 This explains the meaning of the ``placeholder~$\sigma'$\,'' in
Section~\ref{koszul-duality-functors-contramodule-side-subsecn}.

 To conclude, let us say a few words about the bimodule resolution
$\tA\ot_R^{\sigma'}C\ot_R^{\tau'}\tA$ of the diagonal
$\tA$\+$\tA$\+bimodule $\tA$, which was constructed in
Examples~\ref{bimodule-resolution-etc-examples}(3\+-4).
 It was explained in Examples~\ref{bimodule-resolution-etc-examples}
why this is indeed a bimodule resolution, but the construction itself
was decidedly not left-right symmetric, and the meaning of
the placeholder~$\sigma'$ remained mysterious.

 Now this is explained, and the symmetry is restored (under
the additional assumption that $\tA$ is two-sided finitely projective
nonhomogeneous Koszul).
 As a corollary of the above discussion, we have
$$
 \tA\ot_R^{\sigma'}C\ot_R^{\tau'}\tA \,=\,
 \Ksp_{u'}(\tA,\shB)\ot_{\shB}C\ot_B\Ksp_{e'}(B,\tA),
$$
where $C$ is a CDG\+bimodule over $(\shB,\shd,\shh)$
and $(B,d,h)$, as explained above, and the tensor products in
the right-hand side are the CDG\+bimodule tensor products of
Section~\ref{cdg-modules-subsecn}.

\subsection{Koszul triality} \label{koszul-triality-subsecn}
 A ``Koszul triality'' is a commutative triangle diagram of
triangulated equivalences between a derived category of modules,
a coderived category of comodules, and a contraderived category
of contramodules.
 This phenomenon was first observed in~\cite{Pkoszul}
(see~\cite[Sections~6.3\+-6.5]{Pkoszul};
cf.~\cite[Sections~6.7\+-6.8]{Pkoszul}).
 It is \emph{not} supposed to manifest itself in the more general
setting of~\cite[Chapter~11]{Psemi}, which has coalgebra variables in
the base (it is expected to be a ``quadrality'' there; cf.\
Section~\ref{frobenius-quadrality-subsecn}).
 But it still exists in the context of the present paper, under
the assumptions listed in Theorem~\ref{koszul-triality-theorem}
in this section.

 Let $R\subset\tV\subset\tA$ be a right finitely projective
nonhomogeneous Koszul ring, and let $\shB=(\shB,\shd,\shh)$
be the corresponding left finitely projective Koszul CDG\+ring
under the opposite version of the construction of
Proposition~\ref{nonhomogeneous-dual-cdg-ring} and
Corollary~\ref{nonhomogeneous-koszul-duality-anti-equivalence}.
 Then, as a particular case of the discussion in
Section~\ref{revisited-subsecn}, there are DG\+functors
$$
 \tA\ot_R^{\sigma'}{-}\,\:\DG(\shB\comodl)\lrarrow\DG(\tA\modl)
$$
and
$$
 C\ot_R^{\tau'}{-}\,\:\DG(\tA\modl)\lrarrow\DG(\shB\comodl).
$$
 By the opposite version of the discussion in
Section~\ref{koszul-duality-functors-comodule-side-subsecn},
the DG\+functor $\tA\ot_R^{\sigma'}{-}$ is left adjoint to
the DG\+functor $C\ot_R^{\tau'}{-}$.
 Hence the induced triangulated functors between the homotopy
categories,
$$
 \tA\ot_R^{\sigma'}{-}\,\:\Hot(\shB\comodl)\lrarrow\Hot(\tA\modl)
$$
and
$$
 C\ot_R^{\tau'}{-}\,\:\Hot(\tA\modl)\lrarrow\Hot(\shB\comodl),
$$
are also adjoint.

 The following assertions are the opposite versions of
Theorem~\ref{comodule-side-koszul-duality-theorem}
and Corollary~\ref{fin-dim-base-koszul-duality-comodule-side}.

\begin{thm}
 The above pair of adjoint triangulated functors $C\ot_R^{\tau'}{-}$
and $\tA\ot_R^{\sigma'}{-}$ induces a pair of adjoint triangulated
functors between the $\tA/R$\+semicoderived category of left
$\tA$\+modules\/ $\sD^\sico_R(\tA\modl)$ and the coderived category\/
$\sD^\co(\shB\comodl)$ of left CDG\+comodules over
$(\shB,\shd,\shh)$, which are mutually inverse triangulated
equivalences
$$
 \sD^\sico_R(\tA\modl)\simeq\sD^\co(\shB\comodl).
$$ \qed
\end{thm}

\begin{cor} \label{fin-dim-base-koszul-duality-left-comodule-side}
 Assume additionally that the left homological dimension of the ring
$R$ is finite.
 Then the pair of adjoint triangulated functors
$\tA\ot_R^{\sigma'}{-}\,\:\allowbreak\Hot(\shB\comodl)\rarrow
\Hot(\tA\modl)$ and $C\ot_R^{\tau'}{-}\,\:\Hot(\tA\modl)\rarrow
\Hot(\shB\comodl)$ induces mutually inverse triangulated equivalences
\hbadness=1500
\begin{equation} \label{fin-dim-base-left-comodule-side}
 \sD(\tA\modl)\simeq\sD^\co(\shB\comodl)
\end{equation}
between the derived category of left $\tA$\+modules and
the coderived category of left CDG\+comodules over
$(\shB,\shd,\shh)$.  \qed
\end{cor}

 The next theorem is the second main result of
Section~\ref{co-contra-secn}.

\begin{thm} \label{koszul-triality-theorem}
 Let $R\subset\tV\subset\tA$ be a two-sided finitely projective
nonhomogeneous Koszul ring, let $B=(B,d,h)$ be the right finitely
projective Koszul CDG\+ring nonhomogeneous quadratic dual to $\tA$,
and let\/ $\shB=(\shB,\shd,\shh)$ be the left finitely
projective Koszul CDG\+ring nonhomogeneous quadratic dual to $\tA$,
as in Lemma~\ref{nonhomogeneous-koszul-dual-two-sides}.
 Assume that the left homological dimension of the ring $R$ is finite.
 Then the triangulated equivalence of derived nonhomogeneous Koszul
duality on the comodule side~\eqref{fin-dim-base-left-comodule-side}
of Corollary~\ref{fin-dim-base-koszul-duality-left-comodule-side}
(cf.\ the right (co)module side version
in~\eqref{fin-dim-base-comodule-side} of
Corollary~\ref{fin-dim-base-koszul-duality-comodule-side}),
$$
 C\ot_R^{\tau'}{-}\,\:\sD(\tA\modl)\simeq\sD^\co(\shB\comodl)
 :\!\tA\ot_R^{\sigma'}{-}\,,
$$
the triangulated equivalence of derived nonhomogeneous Koszul duality
on the contramodule side~\eqref{fin-dim-base-contramodule-side}
of Corollary~\ref{fin-dim-base-koszul-duality-contramodule-side},
$$
 \Hom_R^{\sigma'}(C,{-})\:\sD(\tA\modl)\simeq\sD^\ctr(B\contra)
 :\!\Hom_R^{\tau'}(\tA,{-}),
$$
and the triangulated equivalence of derived comodule-contramodule
correspondence~\eqref{co-contra-correspondence-eqn} of
Theorem~\ref{co-contra-correspondence-theorem}, \hbadness=2250
$$
 \boR\Hom_{\shB}(C,{-})\:\sD^\co(\shB\comodl)
 \simeq\sD^\ctr(B\contra):\!C\ocn_B^\boL{-}
$$
form a commutative triangle diagram of triangulated category
equivalences,
\begin{equation}
\begin{tikzcd}
&&& \sD^\co(\shB\comodl) \arrow[dd, Leftrightarrow, no head, no tail] \\
\sD(\tA\modl) \arrow[rrru, Leftrightarrow, no head, no tail]
\arrow[rrrd, Leftrightarrow, no head, no tail] \\
&&& \sD^\ctr(B\contra)
\end{tikzcd}
\end{equation}
\end{thm}

\begin{proof}
 As a preliminary observation, it is worth pointing out that
the contratensor product functor $C\ocn_B{-}$ appearing in the third
displayed formula is isomorphic to the tensor product functor
$C\ot_B{-}$ in our present assumptions.
 This is the result of
Proposition~\ref{contratensor-tensor-comparison}(b),
which is applicable in view of
Theorem~\ref{graded-contramodules-fully-faithful-thm}.
 
 To prove the theorem, we notice first of all that, for any complex
of left $\tA$\+modules $M=M^\bu$, one has
$$
 \boR\Hom_{\shB}(C,\>C\ot_R^{\tau'}M^\bu)=
 \Hom_{\shB}(C,\>C\ot_R^{\tau'}M^\bu)
$$
and similarly,
$$
 C\ocn_B^\boL\Hom_R^{\sigma'}(C,M^\bu)=
 C\ocn_B\Hom_R^{\sigma'}(C,M^\bu),
$$
that is, the derived functors $\boR\Hom_{\shB}(C,{-})$ and
$C\ocn_B^\boL{-}$ coincide with the respective underived functors
on these objects.
 Indeed, the graded left $\shB$\+module $C\ot_RM=\Hom_R^\Sigma(\shB,M)$
is coinduced (from the graded left $R$\+module $M$), and the graded
left $B$\+contramodule $\Hom_R(C,M)=B\ot_R^\Pi M$ is induced
(from the graded left $R$\+module~$M$).
 So the objects in question are adjusted to the respective functors
in view of the construction of the derived functors in
the proof of Theorem~\ref{co-contra-correspondence-theorem}.
 
 Finally, the main claim is that, for any complex of left $\tA$\+modules
$M=M^\bu$, there is a natural isomorphism of left CDG\+contramodules
over $(B,d,h)$
$$
 \Hom_{\shB}(C,\>C\ot_R^{\tau'}M^\bu)\simeq\Hom_R^{\sigma'}(C,M^\bu),
$$
as well as a natural isomorphism of left CDG\+comodules over
$(\shB,\shd,\shh)$
$$
 C\ocn_B\Hom_R^{\sigma'}(C,M^\bu)\simeq C\ot_R^{\tau'}M^\bu.
$$
 Indeed, on the level of the underlying graded co/contramodules we have
$$
 \Hom_{\shB}(C,\>C\ot_RM)\simeq\Hom_{\shB}(C,\Hom_R^\Sigma(\shB,M))
 \simeq\Hom_R(C,M)
$$
and
$$
 C\ocn_B\Hom_R(C,M)\simeq C\ocn_B(B\ot_R^\Pi M)\simeq C\ot_R M
$$
(cf.\ Lemma~\ref{coinduced-hom-induced-contratensor} and
the computation in the proof of
Theorem~\ref{co-contra-correspondence-theorem}).
 We leave it to the reader to check that the differentials agree.
\end{proof}

\Section{Koszul Duality and Conversion Functor} \label{conversion-secn}

\subsection{Relatively Frobenius rings} \label{rel-frobenius-subsecn}
 Let $K$ be an associative ring and $B$ be a $K$\+$K$\+bimodule.
 Assume that $B$ is a finitely generated projective left $K$\+module.
 Then for any left $K$\+module $L$ there is a natural isomorphism of
left $K$\+modules $\Hom_K(B,L)\simeq\Hom_K(B,K)\ot_KL$.

 Let $T$ be another $K$\+$K$\+bimodule which is finitely generated
and projective as a left $K$\+module.
 Then there is a natural morphism of $K$\+$K$\+bimodules
$K\rarrow\Hom_K(T,T)\simeq\Hom_K(T,K)\ot_KT$ induced by the right
action of $K$ in~$T$.
 There is also a natural evaluation morphism of $K$\+$K$\+bimodules
$T\ot_K\Hom_K(T,K)\rarrow K$ given by the rule $t\ot f\longmapsto f(t)$
for all $t\in T$ and $f\in\Hom_K(T,K)$.
 The pairing notation $\lan t,f\ran=f(t)$ is convenient.

 We will say that a $K$\+$K$\+bimodule $T$ is \emph{invertible} if
$T$ is projective and finitely generated as a left $K$\+module, and
both the above morphisms $K\rarrow\Hom_K(T,K)\ot_KT$ and
$T\ot_K\Hom_K(T,K)\rarrow K$ are isomorphisms.
 In this case, the functors $T\ot_K\nobreak{-}\,\:\allowbreak
K\modl\rarrow K\modl$ and $\Hom_K(T,{-})=\Hom_K(T,K)\ot_K\nobreak{-}\,\:
\allowbreak K\modl\rarrow K\modl$ are mutually inverse auto-equivalences
of the category of left $K$\+modules.
 Similarly, the functors ${-}\ot_KT\:\modr K\rarrow\modr K$ and
${-}\ot_K\Hom_K(T,K)\:\modr K\rarrow\modr K$ are mutually inverse
auto-equivalences of the category of right $K$\+modules.
 It follows that $T$ is also a finitely generated projective
right $K$\+module, and the notion of an invertible $K$\+$K$\+bimodule
is left-right symmetric.
 Moreover, we obtain a natural isomorphism of $K$\+$K$\+bimodules
$\Hom_K(T,K)\simeq\Hom_{K^\rop}(T,K)$.

 Furthermore, for any invertible $K$\+$K$\+bimodule $T$ and any
finitely generated projective left $K$\+module $L$ the natural
evaluation map $L\rarrow\Hom_{K^\rop}(\Hom_K(L,T),T)$ is an isomorphism
of left $K$\+modules (as one can see, e.~g., by reducing the question
to the case of a free module with one generator $L=K$).

 Let $B'$ and $B''$ be $K$\+$K$\+bimodules, $T$ be an invertible
$K$\+$K$\+bimodule, and $\phi\:B'\ot_KB''\rarrow T$ be
a $K$\+$K$\+bimodule morphism.
 We will say that $\phi$~is a \emph{perfect pairing} if $B'$ is
a finitely generated projective left $K$\+module and
the $K$\+$K$\+bimodule map $\check\phi\:B''\rarrow\Hom_K(B',T)
\simeq\Hom_K(B',K)\ot_KT$ induced by $\phi$~is an isomorphism.
 Then $\Hom_K(B',K)$ is a finitely generated projective right
$K$\+module, and it follows that $B''$ is a finitely generated
projective right $K$\+module.
 Furthermore, the $K$\+$K$\+bimodule map $\hat\phi\:B'\rarrow
\Hom_{K^\rop}(B'',T)\simeq T\ot_K\Hom_{K^\rop}(B'',K)$
induced by~$\phi$ is an isomorphism, because $\hat\phi\simeq
\Hom_{K^\rop}(\check\phi,T)$ (while $\check\phi\simeq
\Hom_K(\hat\phi,T)$).
 Therefore, the notion of a perfect pairing~$\phi$ is left-right
symmetric as well.

 Let $K\rarrow B$ be a ring homomorphism, $T$ be an invertible
$K$\+$K$\+bimodule, and $t\:B\rarrow T$ be a morphism of
$K$\+$K$\+bimodules.
 Take $\phi\:B\ot_KB\rarrow K$ to be equal to the composition of
the multiplication map $B\ot_KB\rarrow B$ with the map~$t$.
 Then the map $\check\phi\:B\rarrow\Hom_K(B,T)$ is a morphism of
left $B$\+modules, while the map $\hat\phi\:B\rarrow
\Hom_{K^\rop}(B,T)$ is a morphism of right $B$\+modules.
 We will say that the ring $B$ is \emph{relatively Frobenius over~$K$}
(with respect to the map~$t$) if the pairing $\phi\:B\ot_KB\rarrow T$
induced by the map $t\:B\rarrow T$ is perfect.
 Following the above discussion, any relatively Frobenius ring $B$
over $K$ is a finitely generated and projective left and right
$K$\+module, and both the maps $\check\phi$ and~$\hat\phi$ are
isomorphisms in this case.

 Let $B$ be a relatively Frobenius ring over $K$ with respect to
a $K$\+$K$\+bimodule morphism $t\:B\rarrow T$.
 Then for any left $K$\+module $L$ there is a natural isomorphism of
left $B$\+modules $B\ot_KL\simeq\Hom_K(B,T)\ot_KL\simeq
\Hom_K(B,K)\ot_KT\ot_KL\simeq\Hom_K(B,\>T\ot_KL)$.
 Consequently, for any left $K$\+module $M$ there is a natural
isomorphism of left $B$\+modules $\Hom_K(B,M)\simeq
B\ot_K\Hom_K(T,M)$.
 So the class of left $B$\+modules induced from left $K$\+modules
coincides with that of left $B$\+modules coinduced from left
$K$\+modules.

\subsection{Relatively Frobenius graded rings}
\label{frobenius-graded-subsecn}
 Let $B=\bigoplus_{n=0}^\infty B^n$ be a nonnegatively graded ring,
and let $K\rarrow B^0$ be a ring homomorphism.
 Let $T$ be an invertible $K$\+$K$\+bimodule.
 Assume that there is an integer $m\ge0$ such that $B^n=0$ for all
$n>m$, and a morphism of $K$\+$K$\+bimodules $t\:B^m\rarrow T$ is given.

 Then we will say that $B$ is a \emph{relatively Frobenius graded ring
over~$K$} (with respect to the map~$t$) if the underlying ungraded ring
of $B$ is relatively Frobenius, in the sense of the definition in
Section~\ref{rel-frobenius-subsecn}, with respect to the composition
of the direct summand projection $B=\bigoplus_{n=0}^mB^n\rarrow B^m$
with the map~$t$.
 Equivalently, this means that the composition $\phi_n\:B^n\ot_KB^{m-n}
\rarrow T$ of the multiplication map $B^n\ot_KB^{m-n}\rarrow B^m$
with the map $t\:B^m\rarrow T$ is a perfect pairing (in the sense of
the definition in Section~\ref{rel-frobenius-subsecn}) for every
$0\le n\le m$.
 In this case, $B^n$ is a finitely generated projective left and right
$K$\+module for every $0\le n\le m$.
 
 Since $B^n=0$ for all $n>m$, all (graded or ungraded, right or left)
$B$\+modules are $B$\+comodules, as mentioned at the end of
Section~\ref{coderived-cdg-comodules-subsecn}.
 So, in particular, we have $B\comodl_\sgr=B\modl_\sgr$,
in the notation of Section~\ref{cdg-comodules-coinduced-type-subsecn},
and therefore we will use $\sD^\co(B\modl)$ instead of
$\sD^\co(B\comodl)$ as the notation for the coderived category of
left CDG\+(co)modules over $B$ when $B$ is endowed with a CDG\+ring
structure $B=(B,d,h)$.
 Similarly, the forgetful functors from the categories of (graded or
ungraded) $B$\+contamodules to the similar categories of $B$\+modules
are category equivalences, as mentioned at the end of
Section~\ref{contraderived-cdg-contramodules-subsecn}.
 So we will use $\sD^\ctr(B\modl)$ instead of $\sD^\ctr(B\contra)$ as
the notation for the contraderived category of
left CDG\+(contra)modules over $B$ when $B$ is endowed with
a CDG\+ring structure.

 The definitions of the induced graded contramodules and coinduced
graded comodules given in
Sections~\ref{cdg-contramodules-induced-type-subsecn}
and~\ref{cdg-comodules-coinduced-type-subsecn} simplify similarly:
for any graded left $K$\+module $L$ one has $B\ot_K^\Pi L=B\ot_KL$
and $\Hom_K^\Sigma(B,L)=\Hom_K(B,L)$.
 So we will speak of \emph{induced} and \emph{coinduced graded
left $B$\+modules} instead of induced graded left $B$\+contramodules
and coinduced graded left $B$\+comodules.

 Finally, since $B$ is a relatively Frobenius graded ring over $K$,
the class of left graded $B$\+modules induced from graded left
$K$\+modules coincides with that of graded left $B$\+modules coinduced
from graded left $K$\+modules, essentially as explained in
Section~\ref{rel-frobenius-subsecn}.
 More precisely, for any graded left $K$\+modules $L$ and $M$ there
are natural isomorphisms of graded left $B$\+modules
$B\ot_KL\simeq\Hom_K(B,\>T\ot_KL[-m])$ and
$\Hom_K(B,M)\simeq B\ot_K\Hom_K(T,M[m])$, where $[i]$, \,$i\in\boZ$, 
denotes the usual (cohomological) grading shift.

\subsection{Relatively Frobenius co-contra correspondence}
 Let $B=(B,d,h)$ be a nonnegatively graded CDG\+ring, so
$B=\bigoplus_{n=0}^\infty B^n$, and let $K\rarrow B^0$ be a ring
homomorphism.
 We will say that the CDG\+ring $(B,d,h)$ is \emph{relatively
Frobenius over~$K$} if the graded ring $B$ is relatively Frobenius
(with respect to some morphism of $K$\+$K$\+bimodules
$t\:B^m\rarrow T$).

 Let $(B,d,h)$ be a relatively Frobenius CDG\+ring over~$K$.
 Similarly to Sections~\ref{cdg-contramodules-induced-type-subsecn}
and~\ref{cdg-comodules-coinduced-type-subsecn}, one can speak of
\emph{left CDG\+modules of the induced type} and \emph{left
CDG\+modules of the coinduced type} over $(B,d,h)$.
 Moreover, following the discussion in
Section~\ref{frobenius-graded-subsecn}, these two classes of
left CDG\+modules coincide.

 Furthermore, similarly to
Sections~\ref{cdg-contramodules-induced-type-subsecn}
and~\ref{cdg-comodules-coinduced-type-subsecn}, one can speak of
\emph{left CDG\+module short exact sequences of the induced type}
and \emph{of the coinduced type} over~$(B,d,h)$.
 Once again, following the discussion in
Section~\ref{frobenius-graded-subsecn}, these two clases of short
exact sequences of left CDG\+modules over $(B,d,h)$ coincide.

 Hence we arrive to the following theorem, which is the first
main result of Section~\ref{conversion-secn}.
 It is to be compared with~\cite[Theorems~3.9\+-3.10]{Pkoszul},
and it is a generalization of~\cite[Theorem~B.3]{Pkoszul}.
 It is also a simplified version of
Theorem~\ref{co-contra-correspondence-theorem} above.

\begin{thm} \label{frobenius-co-contra-thm}
 Let $B=(B,d,h)$ be a relatively Frobenius CDG\+ring over a ring~$K$.
 Assume that the left homological dimension of $K$ is finite.
 Then there is a natural triangulated equivalence between
the coderived and contraderived categories of left CDG\+modules
over $(B,d,h)$,
\begin{equation} \label{frobenius-co-contra-eqn}
 \sD^\co(B\modl)\simeq\sD^\ctr(B\modl).
\end{equation}
\end{thm}

\begin{proof}
 By Proposition~\ref{coderived-coinduced-prop}
and the second assertion of Theorem~\ref{coderived-category-theorem},
the coderived category $\sD^\co(B\modl)$ is equivalent to
the triangulated quotient category of the full triangulated
subcategory in $\Hot(B\modl)$ formed by left CDG\+modules of
the coinduced type over $(B,d,h)$ by its minimal triangulated
subcategory containing all the totalizations of left CDG\+module
short exact sequences of the coinduced type.
 Similarly, by Proposition~\ref{contraderived-induced-prop} and
the second assertion of Theorem~\ref{contraderived-category-theorem},
the contraderived category $\sD^\ctr(B\modl)$ is equivalent to
the triangulated quotient category of the full subcategory in
$\Hot(B\modl)$ formed by left CDG\+modules of the induced type
over $(B,d,h)$ by its minimal triangulated subcategory containing
all the totalizations of left CDG\+module short exact sequences of
the induced type.
 Following the above discussion, it is one and the same triangulated
quotient category.
\end{proof}

\subsection{Relatively Frobenius Koszul graded rings}
\label{relatively-frobenius-koszul-graded-subsecn}
 Let $A=\bigoplus_{n=0}^\infty A_n$ be a left finitely projective
Koszul graded ring with the degree-zero component $R=A_0$, and
let $B=\bigoplus_{n=0}^\infty B^n$ be the quadratic dual
right finitely projective Koszul graded ring, as per the construction
of Section~\ref{quadratic-duality-secn} and
Proposition~\ref{finitely-projective-Koszul-duality}.

 Assume that there is an integer $m\ge0$ such that $B^n=0$ for $n>m$,
the $R$\+$R$\+bimodule $T=B^m$ is invertible, and the graded ring $B$
is relatively Frobenius over the ring $K=R$ with respect to
the identity morphism of $R$\+$R$\+bimodules $B^m\rarrow T$.
 In other words, this means that the multiplication map $B^n\ot_R
B^{m-n}\rarrow B^m$ is a perfect pairing for every $0\le n\le m$.

 Consider the dual Koszul complex $\Ksp_e{}^\bu(B,A)$
\,\eqref{dual-koszul-complex} from
Section~\ref{dual-koszul-complex-subsecn}:
\begin{equation} \label{frobenius-dual-koszul-complex}
 0\lrarrow A\lrarrow B^1\ot_RA\lrarrow\dotsb\lrarrow
 B^{m-1}\ot_RA\lrarrow B^m\ot_RA\lrarrow0,
\end{equation}
and consider also the second Koszul complex
${}^\tau\!K_\bu(B,A)=C\ot_B\Ksp_e{}^\bu(B,A)=C\ot^\tau_RA$
from formula~\eqref{second-koszul-complex}
and Remark~\ref{B-sharp-remark} (where $C=\Hom_{R^\rop}(B,R)$).

 According to Theorem~\ref{projective-koszul-theorem}(e),
the complex ${}^\tau\!K_\bu(B,A)$ is a graded right $A$\+module
resolution of the graded right $A$\+module~$R$.
 In view of the isomorphism of graded right $B$\+modules $B\simeq
\Hom_{R^\rop}(B,T)\simeq T\ot_R\Hom_{R^\rop}(B,R)[-m]=T\ot_RC[-m]$ (see
Sections~\ref{rel-frobenius-subsecn}\+-\ref{frobenius-graded-subsecn}),
the dual Koszul complex~\eqref{frobenius-dual-koszul-complex}
is a graded right $A$\+module resolution of the graded right
$A$\+module $T=B^m$ (placed in the cohomological degree~$m$).
 In other words, the complex~\eqref{frobenius-dual-koszul-complex}
is acyclic at every term except the rightmost one, and at
the rightmost term $B^m\ot_RA$ its cohomology module is
the internal grading component $B^m\ot_RA_0=T$.

\subsection{Two-sided Koszul CDG-rings}
\label{two-sided-koszul-cdg-rings-subsecn}
 Let $B=\bigoplus_{n=0}^\infty B^n$ be a two-sided finitely projective
Koszul graded ring (in the sense of
Section~\ref{two-sided-finitely-projective-subsecn})
with the degree-zero component $R=B^0$.
 Denote by $A$ the left finitely projective Koszul graded ring quadratic
dual to the right finitely projective Koszul ring $B$, and denote by
$A^\#$ the right finitely projective Koszul graded ring quadratic
dual to the left finitely projective Koszul ring $B$.
 So, according to Proposition~\ref{finitely-projective-Koszul-duality},
we have
$$
 A=\Ext_{B^\rop}(R,R) \quad\text{and}\quad
 A^\#=\Ext_B(R,R)^\rop,
$$
where $\Ext_{B^\rop}$ is taken in the category of graded right
$B$\+modules and $\Ext_B$ is taken in the category of graded left
$B$\+modules.
 Following Lemma~\ref{koszul-dual-rings-two-sides}, there are
natural isomorphisms of $R$\+$R$\+bimodules
$$
 A^\#_n=\Hom_R(\Hom_R(A_n,R),R) \quad\text{and}\quad
 A_n=\Hom_{R^\rop}(\Hom_{R^\rop}(A^\#_n,R),R)
$$
for all $n\ge0$, and the multiplicative structures on $A$ and
$A^\#$ are related by the construction of
Remark~\ref{B-sharp-remark} and
Section~\ref{diagonal-cdg-bicomodule-subsecn}.

 Now let $(B,d,h)$ be a nonnegatively graded CDG\+ring.
 We recall that $(B,d,h)$ is said to be right finitely projective
Koszul if the graded ring $B=\bigoplus_{n=0}^\infty B^n$ is right
finitely projective Koszul.
 Dually, let us say that $(B,d,h)$ is \emph{left finitely projective
Koszul} if the graded ring $B$ is left finitely projective Koszul.
 We will say that a CDG\+ring $(B,d,h)$ is \emph{two-sided finitely
projective Koszul} if the graded ring $B$ is two-sided finitely
projective Koszul.

 Let $(B,d,h)$ be a two-sided finitely projective Koszul CDG\+ring.
 Applying the Poincar\'e--Birkhoff--Witt Theorem~\ref{pbw-theorem-thm}
(see also
Corollary~\ref{nonhomogeneous-koszul-duality-anti-equivalence})
to the right finitely projective Koszul CDG\+ring $(B,d,h)$, we
obtain a left finitely projective nonhomogeneous Koszul ring $\tA$
with an increasing filtration $0=F_{-1}\tA\subset R=F_0\tA\subset
\tV=F_1\tA\subset F_2\tA\subset\dotsb$ and a natural isomorphism
of graded rings $A\simeq\gr^F\tA=\bigoplus_{n=0}^\infty
F_n\tA/F_{n-1}\tA$.
 The filtered ring $\tA$ comes together with a left $R$\+linear
splitting
$$
 \Hom_{R^\rop}(B^1,R)\simeq F_1\tA/F_0\tA\simeq A_1=V\hookrightarrow
 \tV=F_1\tA,
$$
which we denote by~$\tau'$, as in
Sections~\ref{koszul-duality-functors-comodule-side-subsecn}
and~\ref{koszul-duality-functors-contramodule-side-subsecn}.

 Applying the opposite version of Theorem~\ref{pbw-theorem-thm}
and Corollary~\ref{nonhomogeneous-koszul-duality-anti-equivalence}
to the left finitely projective Koszul CDG\+ring $(B,d,h)$, we obtain
a right finitely projective nonhomogeneous Koszul ring $\tA^\#$ with
an increasing filtration $0=F_{-1}\tA^\#\subset R=F_0\tA^\#\subset
\tV^\#=F_1\tA^\#\subset F_2\tA^\#\subset\dotsb$ and a natural
isomorphism of graded rings $A^\#\simeq\gr^F\tA^\#=
\bigoplus_{n=0}^\infty F_n\tA^\#/F_{n-1}\tA^\#$.
 The filtered ring $\tA^\#$ comes together with a right $R$\+linear
splitting
$$
 \Hom_R(B^1,R)\simeq F_1\tA^\#/F_0\tA^\#\simeq A^\#_1=V^\#
 \hookrightarrow\tV^\#=F_1\tA^\#,
$$
which we will denote by~$\rho'$.

 The dual nonhomogeneous Koszul CDG\+module $\Ksp(B,\tA)$ was
constructed in
Section~\ref{nonhomogeneous-koszul-cdg-module-subsecn}.
 We recall that $\Ksp(B,\tA)$ is a CDG\+bimodule over the CDG\+rings
$B=(B,d,h)$ and $\tA=(\tA,0,0)$ whose underlying graded
$B$\+$\tA$\+bimodule is $B\ot_R\tA$.
 The opposite construction produces a CDG\+bimodule $\Ksp(\tA^\#,B)$
over the CDG\+rings $\tA^\#=(\tA^\#,0,0)$ and $B=(B,d,h)$.
 So the underlying graded $\tA^\#$\+$B$\+bimodule of $\Ksp(\tA^\#,B)$
is $\tA^\#\ot_RB$.

 From now on we will assume for simplicity that $B^n=0$ for $n\gg0$;
so there is no difference between $B$\+modules and $B$\+comodules,
and $B$ is a finitely generated projective graded left and right
$R$\+module.
 Let $D=(D,d_D,h_D)$ and $E=(E,d_E,h_E)$ be two CDG\+rings.
 Let $M$ be a CDG\+bimodule over $(\tA,0,0)$ and $(E,d_E,h_E)$, and
let $N$ be a CDG\+bimodule over $(D,d_D,h_D)$ and $(B,d,h)$.
 Following the constructions of
Sections~\ref{koszul-duality-functors-comodule-side-subsecn}
and~\ref{revisited-subsecn}, the CDG\+bimodule $N\ot_R^{\tau'}M$
over $D$ and $E$ is defined as the tensor product of CDG\+bimodules
$$
 N\ot_R^{\tau'}M = N\ot_B\Ksp(B,\tA)\ot_{\tA}M.
$$

 Furthermore, let $M$ be a CDG\+bimodule over $(E,d_E,h_E)$ and
$(\tA,0,0)$.
 Then the CDG\+bimodule $M\ot_R^{\sigma'}C$ over $(E,d_E,h_E)$
and $(B,d,h)$ is defined as the Hom CDG\+bimodule
$$
 M\ot_R^{\sigma'}C=\Hom_{\tA^\rop}(\Ksp(B,\tA),M),
$$
where $\sigma'$~is the placeholder introduced in
Section~\ref{koszul-duality-functors-comodule-side-subsecn}
and explained in Section~\ref{revisited-subsecn},
while $C$ is the graded coring $C=\Hom_{R^\rop}(B,R)$.

 Let us introduce notation for the opposite constructions.
 Let $M$ be a CDG\+bimod\-ule over $(E,d_E,h_E)$ and $(\tA^\#,0,0)$,
and let $N$ be a CDG\+bimodule over $(B,d,h)$ and $(D,d_D,h_D)$.
 We define the CDG\+bimodule $M\ot_R^{\rho'}N$ over $E$ and $D$
as the tensor product of CDG\+bimodules
$$
 M\ot_R^{\rho'}N=M\ot_{\tA^\#}\Ksp(\tA^\#,B)\ot_BN.
$$
 So the underlying graded $E$\+$D$\+bimodule of $M\ot_R^{\rho'}N$
is $M\ot_RN$.

 Furthermore, let $M$ be a CDG\+bimodule over $(\tA^\#,0,0)$ and
$(E,d_E,h_E)$.
 Then the CDG\+bimodule $C^\#\ot_R^{\pi'}M$ over $(B,d,h)$
and $(E,d_E,h_E)$ is defined as the Hom CDG\+bimodule
$$
 C^\#\ot_R^{\pi'}M=\Hom_{\tA^\#}(\Ksp(\tA^\#,B),M),
$$
where $C^\#$ is the graded coring $C^\#=\Hom_R(B,R)$ over $R$,
while $\pi'$~is another placeholder.
 The underlying graded $B$\+$E$\+bimodule of $C^\#\ot_R^{\pi'}M$
is $\Hom_{\tA^\#}(\tA^\#\ot_RB,\>M)\simeq\Hom_R(B,M)\simeq
C^\#\ot_RM$.

\subsection{Conversion bimodule} \label{conversion-bimodule-subsecn}
 Let $B=(B,d,h)$ be a two-sided finitely projective Koszul
nonnegatively graded CDG\+ring (in the sense of
Section~\ref{two-sided-koszul-cdg-rings-subsecn}) with the degree-zero
component $R=B^0$.
 We will say that $(B,d,h)$ is a \emph{relatively Frobenius Koszul
CDG\+ring} if the graded ring $B$ is relatively Frobenius over $R$ in
the sense of Section~\ref{relatively-frobenius-koszul-graded-subsecn}.
 This means that there is an integer $m\ge0$ such that $B^n=0$ for
$n>m$, the $R$\+$R$\+bimodule $T=B^m$ is invertible, and
the multiplication map $B^n\ot_RB^{m-n}\rarrow B^m$ is a perfect
pairing for every $0\le n\le m$ (in the sense of
Section~\ref{rel-frobenius-subsecn}).

 Let $(B,d,h)$ be a relatively Frobenius Koszul CDG\+ring.
 We will use the notation of
Section~\ref{two-sided-koszul-cdg-rings-subsecn}.
 In particular, $\tA$ is the left finitely projective nonhomogeneous
Koszul ring corresponding to the right finitely projective Koszul
CDG\+ring $(B,d,h)$, while $\tA^\#$ is the right finitely projective
nonhomogeneous Koszul ring correponding to the left finitely projective
Koszul CDG\+ring $(B,d,h)$.

\begin{lem} \label{frobenius-nonhomogeneous-dual-koszul}
\textup{(a)} There is a natural morphism $T[-m]\rarrow B\ot_R^{\tau'}
\tA=\Ksp(B,\tA)$ of left CDG\+modules over $(B,d,h)$, where
the one-term complex of left $R$\+modules $T[-m]$ is endowed with
the trivial structure of left CDG\+(co)module over $(B,d,h)$.
 The cone of this morphism is a coacyclic left CDG\+(co)module over
$(B,d,h)$. \par
\textup{(b)} There is a natural morphism $T[-m]\rarrow \tA^\#\ot_R^
{\rho'}B=\Ksp(\tA^\#,B)$ of right CDG\+modules over $(B,d,h)$, where
the one-term complex of right $R$\+modules $T[-m]$ is endowed with
the trivial structure of right CDG\+(co)module over $(B,d,h)$.
 The cone of this morphism is a coacyclic right CDG\+(co)module over
$(B,d,h)$.
\end{lem}

\begin{proof}
 Let us prove part~(a); part~(b) is the opposite version.
 The left finitely projective nonhomogeneous Koszul ring $\tA$ is
endowed with an increasing filtration~$F$.
 Define a decreasing filtration $F$ on the graded ring $B$ by
the rule $F^nB=\bigoplus_{i\ge n} B^i[-i]\subset B$; so $F$ is
the decreasing filtration on $B$ induced by the grading.
 Define an increasing filtration $F$ on $B$ by the rule
$F_{-n}B=F^nB$; so $0=F_{-m-1}B\subset F_{-m}B\subset\dotsb
\subset F_0B=B$.
 Consider the increasing filtration $F$ on $\Ksp(B,\tA)=B\ot_R\tA$
induced by the increasing filtrations $F$ on $B$ and $\tA$,
i.~e., $F_n(B\ot_R\tA)=\sum_{i+j=n}F_iB\ot_R F_j\tA$.

 It is clear from the definition of the differential $d_{e'}$ on
$\Ksp(B,\tA)$ that $F_n\Ksp(B,\tA)$ is a CDG\+submodule of
the left CDG\+module $\Ksp(B,\tA)$ over $(B,d,h)$ for every
$-m\le n<\infty$.
 Furthermore, the CDG\+submodule $F_{-m}\Ksp(B,\tA)=B^m[-m]\ot_RF_0A
=B^m[-m]\subset\Ksp(B,\tA)$ is isomorphic to $T[-m]$; hence
the desired (injective) morphism of CDG\+modules $T[-m]\rarrow
\Ksp(B,\tA)$.
 In view of Lemma~\ref{cdg-comodule-filtration-coacyclic-lemma},
in order to show that the cone (equivalently, the cokernel) of
this morphism is coacyclic, it suffices to check that
the CDG\+(co)module $F_n\Ksp(B,\tA)/F_{n-1}\Ksp(B,\tA)$ over
$(B,d,h)$ is coacyclic for every $n\ge-m+1$.

 In fact, $F_n\Ksp(B,\tA)/F_{n-1}\Ksp(B,\tA)$ is a trivial
CDG\+comodule over $(B,d,h)$, that is, the ideal $B^{\ge1}\subset B$
acts by zero in this quotient CDG\+module.
 The whole associated graded CDG\+module $\gr^F\Ksp(B,\tA)=
\bigoplus_{n=-m}^\infty F_n\Ksp(B,\tA)/F_{n-1}\Ksp(B,\tA)$ is
the complex of left $R$\+modules $\Ksp_e{}^\bu(B,A)=B\ot_RA$
\,\eqref{dual-koszul-complex} from
Section~\ref{dual-koszul-complex-subsecn}, viewed as a left CDG\+module
over $(B,d,h)$ with the trivial CDG\+(co)module structure.
 In view of the discussion in
Section~\ref{relatively-frobenius-koszul-graded-subsecn}, this means
that, for every $n\ge-m+1$, the complex of left $R$\+modules
$F_n\Ksp(B,\tA)/F_{n-1}\Ksp(B,\tA)$ is acyclic.
 This is a finite complex of finitely generated projective left
$R$\+modules; so it is contractible as a complex of left
$R$\+modules, hence it is also contractible (and consequently,
coacyclic) as a trivial left CDG\+(co)module over $(B,d,h)$.
\end{proof}

\begin{lem} \label{frobenius-quasi-isomorphism-lemma}
\textup{(a)} For any complex of right $\tA^\#$\+modules $M^\bu$,
there is a natural quasi-isomorphism of complexes of right $R$\+modules
$$
 M^\bu\ot_RT[-m] \lrarrow M^\bu\ot_R^{\rho'}B\ot_R^{\tau'}\tA,
$$
where $M^\bu\ot_R^{\rho'}B\ot_R^{\tau'}\tA=
M^\bu\ot_{\tA^\#}\Ksp(\tA^\#,B)\ot_B\Ksp(B,\tA)$ is the complex
of right $\tA$\+modules provided by the constructions
of Section~\ref{two-sided-koszul-cdg-rings-subsecn}. \par
\textup{(b)} For any complex of left $\tA$\+modules $L^\bu$,
there is a natural quasi-isomorphism of complexes of left $R$\+modules
$$
 T\ot_RL^\bu[-m] \lrarrow \tA^\#\ot_R^{\rho'}B\ot_R^{\tau'}L^\bu,
$$
where $\tA^\#\ot_R^{\rho'}B\ot_R^{\tau'}L^\bu=
\Ksp(\tA^\#,B)\ot_B\Ksp(B,\tA)\ot_{\tA}L^\bu$ is the complex
of left $\tA^\#$\+modules provided by the constructions
of Section~\ref{two-sided-koszul-cdg-rings-subsecn}.
\end{lem}

\begin{proof}
 Let us explain part~(a); part~(b) is opposite.
 The desired morphism is obtaned by applying the functor
$M^\bu\ot_R^{\rho'}{-}$ to the morphism of CDG\+modules in
Lemma~\ref{frobenius-nonhomogeneous-dual-koszul}(a).
 It is important here that the latter morphism is right $R$\+linear;
so it is, in fact, a morphism of CDG\+bimodules over $(B,d,h)$
and $(R,0,0)$.
 Furthermore, following the proof of
Lemma~\ref{frobenius-nonhomogeneous-dual-koszul}(a), the cone of
the morphism $T[-m]\rarrow B\ot_R^{\tau'}\tA$ is coacyclic as
a CDG\+bimodule over $(B,d,h)$ and $(R,0,0)$, and in fact even as
a CDG\+bimodule over $(B,d,h)$ and $(R,0,0)$ with a projective
underlying graded left $R$\+module (because any finite acyclic
complex of $R$\+$R$\+bimodules is coacyclic).
 Hence the tensor product functor $M^\bu\ot_R^{\rho'}{-}$ preserves
coacyclicity of the cone, so the cone of the resulting morphism
$M^\bu\ot_RT[-m]\rarrow M^\bu\ot_R^{\rho'}B\ot_R^{\tau'}\tA$ is not
only acyclic but even coacyclic as a complex of right $R$\+modules.
\end{proof}

 We are interested in the finite complex of
$\tA^\#$\+$\tA$\+bimodules $\tA^\#\ot_R^{\rho'}B\ot_R^{\tau'}\tA$,
or more specifically, in its top cohomology bimodule
$E=H^m(\tA^\#\ot_R^{\rho'}B\ot_R^{\tau'}\tA)$.
 According to Lemma~\ref{frobenius-quasi-isomorphism-lemma}
applied to the one-term complexes $M^\bu=\tA^\#$ and $L^\bu=\tA$,
we have $H^n(\tA^\#\ot_R^{\rho'}B\ot_R^{\tau'}\tA)=0$ for all $n\ne m$,
while the $\tA^\#$\+$\tA$\+bimodule $E$ is naturally isomorphic to
$\tA^\#\ot_RT$ as an $\tA^\#$\+$R$\+bimodule and naturally isomorphic
to $T\ot_R\tA$ as an $R$\+$\tA$\+bimodule.
 The following theorem tells that $E$ is a Morita equivalence bimodule
for the rings $A$ and~$A^\#$.

\begin{thm} \label{conversion-equivalence-theorem}
 Let $(B,d,h)$ be a relatively Frobenius Koszul CDG\+ring, and let
$\tA$ and $\tA^\#$ be the two nonhomogeneous Koszul rings quadratic
dual to $B$ on the two sides, as above.
 Then \par
\textup{(a)} the tensor product and Hom functors
$$
 E\ot_{\tA}{-}\,\:\tA\modl\,\simeq\,
 \tA^\#\modl\,:\!\Hom_{\tA^\#}(E,{-})
$$
are mutually inverse equivalences of abelian categories; \par
\textup{(b)} the tensor product and Hom functors
$$
 {-}\ot_{\tA^\#}E\:\modr\tA^\#\,\simeq\,
 \modr\tA\,:\!\Hom_{\tA^\rop}(E,{-})
$$
are mutually inverse equivalences of abelian categories.
\end{thm}

\begin{proof}
 The left $\tA^\#$\+module $E\simeq\tA^\#\ot_RT$ is a finitely
generated projective generator of $\tA^\#\modl$, since the left
$R$\+module $T$ is a finitely generated projective generator of
$R\modl$.
 Similarly, the right $\tA$\+module $E\simeq T\ot_R\tA$ is
a finitely generated projective generator of $\modr\tA$, since
the right $R$\+module $T$ is a finitely generated projective
generator of $\modr R$.
 Finally, one computes that $\Hom_{\tA^\#}(E,E)\simeq
\Hom_{\tA^\#}(\tA^\#\ot_R\nobreak T,\allowbreak\>E)\simeq
\Hom_R(T,\>E)\simeq\Hom_R(T,\>T\ot_R\tA)\simeq\Hom_R(T,T)\ot_R\tA
\simeq\tA^\rop$, and similarly $\Hom_{\tA^\rop}(E,E)\simeq\tA^\#$.
 These observations imply the assertions of the theorem.
\end{proof}

 The result of Theorem~\ref{conversion-equivalence-theorem} can be
rephrased as follows.
 For any left $\tA$\+module $L$, the tensor product $T\ot_RL$ has
a natural left $\tA^\#$\+module structure.
 For any right $\tA^\#$\+module $M$, the tensor product
$M\ot_RT$ has a natural right $\tA$\+module structure.
 For any left $\tA^\#$\+module $N$, the left $R$\+module
$\Hom_R(T,N)$ has a natural left $\tA$\+module structure.
 For any right $\tA$\+module $N$, the right $R$\+module
$\Hom_{R^\rop}(T,N)$ has a natural right $\tA^\#$\+module structure.

\begin{rem} \label{commutative-cdg-ring-remark}
 When $(B,d,h)$ is a graded commutative DG\+ring (so the graded ring
$B$ is graded commutative and $h=0$), the two filtered rings
$\tA$ and $\tA^\#$ are naturally opposite to each other, that is
$\tA^\#\simeq\tA^\rop$.
 This isomorphism is induced by the natural identity isomorphism
between the CDG\+ring $B=(B,d,h)$ and its opposite CDG\+ring
$B^\rop=(B^\rop,d^\rop,-h^\rop)$.
 So the Morita equivalence bimodule $E$ can be viewed as having
two right $\tA$\+module structures which commute with each other,
and Theorem~\ref{conversion-equivalence-theorem} provides
an equivalence between the categories of left and right $\tA$\+modules,
$\tA\modl\simeq\modr\tA$, in this case.

 Here, for an arbitrary CDG\+ring $(B,d,h)$, the CDG\+ring
$(B^\rop,d^\rop,-h^\rop)$ is defined as follows.
 As a graded abelian group, $B^\rop$ is identified with $B$ by
the map denoted by $b\longmapsto b^\rop\:B^n\rarrow B^{\rop,n}$,
\ $n\in\boZ$.
 The multiplication in $B^\rop$ is given by the formula
$b^\rop c^\rop=(-1)^{|b||c|}(cb)^\rop$ for all $b$, $c\in B$,
the differential is $d^\rop(b^\rop)=d(b)^\rop$, and
the curvature element is $-h^\rop\in B^{\rop,2}$.
 The construction of the opposite CDG\+ring $(B^\rop,d^\rop,-h^\rop)$ 
has the expected property that the DG\+category of right CDG\+modules
over $(B,d,h)$ is equivalent to the DG\+category of left CDG\+modules
over $(B^\rop,d^\rop,-h^\rop)$ and vice versa.
 Notice that the formulas above imply that the map 
$b\longmapsto b^\rop\:B\rarrow B^\rop$ is \emph{not} an isomorphism
of CDG\+rings when the ring $B$ is graded commutative, but
$2h\ne0$ in~$B^2$.
 Therefore, there is \emph{no} natural isomorphism between the rings
$\tA^\#$ and $\tA^\rop$ in the situation with a nonzero curvature
in characteristic different from~$2$ (generally speaking).
\end{rem}

\subsection{Relatively Frobenius Koszul quadrality}
\label{frobenius-quadrality-subsecn}
 We keep the assumptions and notation of
Sections~\ref{two-sided-koszul-cdg-rings-subsecn}\+-%
\ref{conversion-bimodule-subsecn}.
 So $B=(B,d,h)$ is a relatively Frobenius Koszul CDG\+ring, $\tA$ is
the left finitely projective nonhomogeneous Koszul ring nonhomogeneous
quadratic dual to the right finitely projective Koszul CDG\+ring
$(B,d,h)$, and $\tA^\#$ is the right finitely projective nonhomogeneous
Koszul ring nonhomogeneous quadratic dual to the left finitely
projective Koszul CDG\+ring $(B,d,h)$.
 The $\tA^\#$\+$\tA$\+bimodule $E$ was constructed in
the paragraph before Theorem~\ref{conversion-equivalence-theorem}.

 Following the discussion in
Section~\ref{two-sided-koszul-cdg-rings-subsecn}, there are DG\+functors
$$
 \tA^\#\ot_R^{\rho'}{-}\,\:\DG(B\modl)\lrarrow\DG(\tA^\#\modl)
$$
and
$$
 C^\#\ot_R^{\pi'}{-}\,\:\DG(\tA^\#\modl)\lrarrow\DG(B\modl)
$$
between the DG\+category of complexes of left $\tA^\#$\+modules
and the DG\+category of left CDG\+modules over $(B,d,h)$.
 Following the discussion in the beginning of
Section~\ref{koszul-triality-subsecn}, the DG\+functor
$\tA^\#\ot_R^{\rho'}{-}$ is left adjoint to the DG\+functor
$C^\#\ot_R^{\pi'}{-}$.
 Hence the induced triangulated functors between the homotopy
categories,
$$
 \tA^\#\ot_R^{\rho'}{-}\,\:\Hot(B\modl)\lrarrow\Hot(\tA^\#\modl)
$$
and
$$
 C^\#\ot_R^{\pi'}{-}\,\:\Hot(\tA^\#\modl)\lrarrow\Hot(B\modl),
$$
are also adjoint on the respective sides.

\begin{lem} \label{quadrality-commutativity-lemma}
\textup{(a)} For any complex of right $\tA^\#$\+modules $M^\bu$, there
is a natural closed isomorphism of right CDG\+modules over $(B,d,h)$
$$
 M^\bu\ot_{\tA^\#}E[-m] \ot_R^{\sigma'}C \simeq
 M^\bu\ot_R^{\rho'}B.
$$ \par
\textup{(b)} For any complex of left $\tA$\+modules $L^\bu$, there is
a natural closed isomorphism of left CDG\+modules over $(B,d,h)$
$$
 C^\#\ot_R^{\pi'} E[-m]\ot_{\tA}L^\bu \simeq
 B\ot_R^{\tau'}L^\bu.
$$
\end{lem}

\begin{proof}
 Let us prove part~(b); part~(a) is opposite.
 We will construct a closed isomorphism $B\ot_R^{\tau'}\tA\rarrow
C^\#\ot_R^{\pi'}E[-m]$ of CDG\+bimodules over $(B,d,h)$ and
$(\tA,0,0)$.
 Then the desired closed isomorphism of left CDG\+modules in~(b)
will be produced by applying the DG\+functor ${-}\ot_{\tA}L^\bu$.

 Recall that, by the definition, the $\tA^\#$\+$\tA$\+bimodule $E$
is constructed as $E=H^m(\tA^\#\ot_R^{\rho'}B\ot_R^{\tau'}\tA)$.
 Moreover, there is a natural quasi-isomorphism of finite complexes
of $\tA^\#$\+$\tA$\+bimodules $\tA^\#\ot_R^{\rho'}B\ot_R^{\tau'}\tA
\rarrow E[-m]$.
 Applying the DG\+functor $C^\#\ot_R^{\pi'}{-}$, we obtain a closed
morphism
\begin{equation} \label{induced-by-taking-top-cohomology}
 C^\#\ot_R^{\pi'}\tA^\#\ot_R^{\rho'}B\ot_R^{\tau'}\tA
 \lrarrow C^\#\ot_R^{\pi'}E[-m]
\end{equation}
of CDG\+bimodules over $(B,d,h)$ and $(\tA,0,0)$.
 In fact, the cone of the closed
morphism~\eqref{induced-by-taking-top-cohomology} is a coacyclic
CDG\+bimodule (since $C^\#$ is a projective graded right
$R$\+module), but we will not need to use this observation.

 For any left CDG\+module $N$ over $(B,d,h)$, there is a natural
adjunction morphism $N\rarrow C^\#\ot_R^{\pi'}\tA^\#\ot_R^{\rho'}N$,
which is a closed morphism of left CDG\+modules over $(B,d,h)$.
 In particular, we are interested in the adjunction morphism
\begin{equation} \label{particular-case-of-dg-adjunction}
 B\ot_R^{\tau'}\tA \lrarrow
 C^\#\ot_R^{\pi'}\tA^\#\ot_R^{\rho'}B\ot_R^{\tau'}\tA,
\end{equation}
which is a closed morphism of CDG\+bimodules over
$(B,d,h)$ and $(\tA,0,0)$.
 The composition of~\eqref{particular-case-of-dg-adjunction}
and~\eqref{induced-by-taking-top-cohomology} is the desired closed
morphism of CDG\+bimodules
\begin{equation} \label{closed-isomorphism-of-CDG-bimodules}
 B\ot_R^{\tau'}\tA \lrarrow C^\#\ot_R^{\pi'}E[-m].
\end{equation}

 In order to show that~\eqref{closed-isomorphism-of-CDG-bimodules}
is an isomorphism, one can define filtrations on the left-hand and
the right-hand side and check that the associated graded map is
an isomorphism.
 The nonhomogeneous Koszul rings $\tA$ and $\tA^\#$ are endowed with
increasing filtrations~$F$.
 A finite increasing filtration $F$ on the graded ring $B$ was defined
in the proof of Lemma~\ref{frobenius-nonhomogeneous-dual-koszul}.
 Hence one obtains an induced increasing filtration on the complex
of $\tA^\#$\+$\tA$\+bimodules $\tA^\#\ot_R^{\rho'}B\ot_R^{\tau'}\tA$,
and consequently on its cohomology $\tA^\#$\+$\tA$\+bimodule~$E$.
 Finally, a finite increasing filtration $F$ on the graded coring $C$
was defined in the proof of
Theorem~\ref{comodule-side-koszul-duality-theorem} (as the filtration
induced by the grading); a filtration $F$ on the graded coring $C^\#$
is constructed similarly.

 Both the left-hand and the right-hand sides of
the maps~\eqref{induced-by-taking-top-cohomology}
and~\eqref{particular-case-of-dg-adjunction} acquire induced
filtrations, which are preserved by both the maps.
 It is straightforward to check that
the composition~\eqref{closed-isomorphism-of-CDG-bimodules} becomes
an isomorphism after the passage to the associated graded bimodules
with respect to~$F$.
\end{proof}

 The following assertion is a restatement of
Corollary~\ref{fin-dim-base-koszul-duality-left-comodule-side}
with $\tA$ replaced by $\tA^\#$ and $\shB$ replaced by~$B$,
and with the notation taking into account the assumption that
$B^n=0$ for $n>m$.
 (Recall that
Corollary~\ref{fin-dim-base-koszul-duality-left-comodule-side},
in turn, is the opposite version of
Corollary~\ref{fin-dim-base-koszul-duality-comodule-side}.)

\begin{cor} \label{fin-dim-base-frobenius-comodule-side-cor}
 Assume that the left homological dimension of the ring $R$ is finite.
 Then the pair of adjoint triangulated functors
$\tA^\#\ot_R^{\rho'}{-}\,\:\Hot(B\modl)\rarrow\Hot(\tA^\#\modl)$
and $C^\#\ot_R^{\pi'}{-}\,\:\Hot(\tA^\#\modl)\rarrow
\Hot(B\modl)$ induces mutually inverse triangulated equivalences
\begin{equation} \label{fin-dim-base-frobenius-comodule-side-eqn}
 \sD(\tA^\#\modl)\simeq\sD^\co(B\modl)
\end{equation}
between the derived category of left $\tA^\#$\+modules and
the coderived category of left CDG\+modules over $(B,d,h)$.  \qed
\end{cor}

 The next theorem is the second main result of
Section~\ref{conversion-secn}.

\begin{thm} \label{frobenius-quadrality-theorem}
 Let $B=(B,d,h)$ be a relatively Frobenius Koszul CDG\+ring whose
degree-zero component $R=B^0$ is a ring of finite left homological
dimension.
 Let $\tA$ and $\tA^\#$ be two nonhomogeneous Koszul rings quadratic
dual to $B$ on the two sides, as above.
 Then there is a commutative square diagram of triangulated
equivalences
\begin{equation} \label{frobenius-quadrality-diagram}
\begin{tikzcd}
 \sD(\tA^\#\modl) \arrow[rrrr, bend right = 10, "{\Hom_R(T[-m],{-})}"']
 \arrow[dddd, "{C^\#\ot_R^{\pi'}{-}}"] &&&&
 \sD(\tA\modl) \arrow[llll, "{T[-m]\ot_R{-}}"']
 \arrow[dddd, "{\Hom_R^{\sigma'}(C,{-})}=B\ot_R^{\tau'}{-}"']
 \\ \\ \\ \\
 \sD^\co(B\modl) \arrow[rrrr, Leftrightarrow, no head, no tail]
 \arrow[uuuu, bend left = 20, "{\tA^\#\ot_R^{\rho'}{-}}"] &&&&
 \sD^\ctr(B\modl)
 \arrow[uuuu, bend right = 20, "{\Hom_R^{\tau'}(\tA,{-})}"']
\end{tikzcd}
\end{equation}
 Here the triangulated equivalence in the upper line is induced by
the equivalence of abelian categories from
Theorem~\ref{conversion-equivalence-theorem}(a) (up to
the cohomological shift by~$[-m]$).
 The triangulated equivalence in the leftmost column is the assertion
of Corollary~\ref{fin-dim-base-frobenius-comodule-side-cor}.
 The triangulated equivalence in the rightmost column is
a particular case of
Corollary~\ref{fin-dim-base-koszul-duality-contramodule-side}.
 The triangulated equivalence in the lower line is provided by
Theorem~\ref{frobenius-co-contra-thm}.
\end{thm}

\begin{proof}
 Recall that, by the definition, for any complex of left $\tA$\+modules
$L^\bu$ we have $\Hom_R^{\sigma'}(C,L^\bu)=\Ksp_{e'}(B,\tA)\ot_{\tA}
L^\bu=B\ot_R^{\tau'}L^\bu$, since $B^n=0$ for $n>m$
(see Sections~\ref{koszul-duality-functors-contramodule-side-subsecn}
and~\ref{revisited-subsecn}).
 Hence the equality in the label at the rightmost straight arrow in
the diagram.
 We also have $B\contra=B\modl=B\comodl$ because $B^n=0$ for $n>m$.

 It remains to explain why the diagram is commutative.
 It is convenient to check commutativity of the diagram of triangulated
functors between the homotopy categories formed by the three straight
arrows and the equality in the lower horizontal line.
 This follows from the natural closed isomorphism of left CDG\+modules
in Lemma~\ref{quadrality-commutativity-lemma}(b).
 Then one observes that both sides of the latter closed isomorphism
are CDG\+modules of the induced and coinduced type; so they are
adjusted to the derived functors in the construction of the triangulated
equivalence in Theorem~\ref{frobenius-co-contra-thm}.
\end{proof}

\Section{Examples} \label{examples-secn}

\subsection{Symmetric and exterior algebras}
\label{symmetric-exterior-subsecn}
 Let us start with the tensor ring.
 Let $R$ be an associative ring and $V$ be an $R$\+$R$\+bimodule.
 Consider the tensor ring $T_R(V)=
\bigoplus_{n=0}^\infty V^{\ot_R\,n}$, as defined in
Section~\ref{quadratic-duality-secn} (see the notation in
Section~\ref{relative-bar-subsecn}).
 By the definition, the graded ring $A=T_R(V)$ is quadratic.

 It is clear from, e.~g., Theorem~\ref{flat-koszul-theorem}(b)
that the graded ring $A$ is left flat Koszul whenever $V$ is a flat
left $R$\+module.
 By Theorem~\ref{projective-koszul-theorem}(a), the graded ring $A$
is left finitely projective Koszul whenever $V$ is a finitely
generated projective left $R$\+module.
 In the latter case, the quadratic dual right finitely projective
Koszul graded ring $B$ has the components $B_0=R$, \ $B_1=\Hom_R(V,R)$,
and $B_n=0$ for all $n\ge2$.

 Now let $R$ be a commutative ring and $V$ be an $R$\+module, viewed
as an $R$\+$R$\+bi\-module in which the left and right actions
of $R$ coincide.
 Then $T_R(V)$ is a graded $R$\+algebra (or in other words,
the degree-zero component $R=T_{R,0}(V)$ lies in the center of
$T_R(V)$).
 By the definition, the \emph{symmetric algebra} $\Sym_R(V)=
\bigoplus_{n=0}^\infty\Sym^n_R(V)$ is the largest commutative
quotient algebra of $T_R(V)$.
 Equivalently, $\Sym_R(V)=T_R(V)/(I)$ is the quadratic algebra
over $R$ with the submodule of quadratic relations $I\subset V\ot_RV$
spanned by all the tensors $v\ot w-w\ot v\in V\ot_RV$ with
$v$, $w\in V$.

 The \emph{exterior algebra} $\Lambda_R(V)=\bigoplus_{n=0}^\infty
\Lambda_R^n(V)$ is the largest \emph{strictly} graded commutative
quotient algebra of $T_R(V)$, that is, the largest quotient algebra
of $T_R(V)$ in which the identities $ab=(-1)^{|a||b|}ba$ and $c^2=0$
hold for all elements $a$~of degree~$|a|$, \ $b$~of degree~$|b|$,
and $c$~of odd degree.
 Equivalently, $\Lambda_R(V)=T_R(V)/(I)$ is the quadratic algebra
over $R$ with the submodule of quadratic relations $I\subset V\ot_RV$
spanned by all the tensors $v\ot v\in V\ot_RV$ with $v\in V$.
 When the element $2\in R$ is invertible, the submodule
$I\subset V\ot_RV$ is also spanned by the tensors $v\ot w+w\ot v$
with $v$, $w\in V$.

 When $V$ is a (finitely generated) free $R$\+module, the grading
components of the quadratic algebras $\Sym_R(V)$ and $\Lambda_R(V)$
are also (finitely generated) free $R$\+modules with explicit basises
which are easy to construct.
 It follows by passing to retracts that the $R$\+modules
$\Sym_R^n(V)$ and $\Lambda_R^n(V)$ are (finitely generated) projective
whenever the $R$\+module $V$ is (finitely generated) projective.
 Passing to the filtered direct limits, one can see that the grading
components of $\Sym_R(V)$ and $\Lambda_R(V)$ are flat $R$\+modules
whenever $V$ is a flat $R$\+module.

 It is worth noticing that the na\"\i ve definition of an exterior
algebra as the quotient algebra of the tensor algebra by the relations
$v\ot w+w\ot v$ with $v$, $w\in V$ (equivalently, the maximal quotient
algebra of $T_R(V)$ in which the identity $ab=(-1)^{|a||b|}ba$ holds)
does \emph{not} have these flatness/projectivity properties.
 For a counterexample, it suffices to take $R$ to be the ring of
integers $\boZ$ or the ring of $2$\+adic integers $\boZ_2$, and $V$
to be the free $R$\+module with one generator.

 Let $V$ be a finitely generated projective $R$\+module, and let
$V\spcheck=\Hom_R(V,R)$ denote the dual finitely generated projective
$R$\+module.
 Then the quadratic graded rings $A=\Sym_R(V)$ and
$B=\Lambda_R(V\spcheck)$ are quadratic dual to each other
(in the sense of Propositions~\ref{2-fin-proj-quadratic-duality}
and~\ref{3-fin-proj-quadratic-duality}).
 Moreover, both the quadratic algebras $\Sym_R(V)$ and
$\Lambda_R(V\spcheck)$ are (left and right) finitely projective Koszul,
as one can see, e.~g., from exactness of the Koszul complex associated
with any chosen basis in the module of generators $V$ of the polynomial
algebra $\Sym_R(V)$ (use Theorem~\ref{projective-koszul-theorem}(d)
or~(e)).

 As the class of (left or right) flat Koszul graded rings is closed
under filtered direct limits, it follows that the quadratic graded
rings $\Sym_R(V)$ and $\Lambda_R(V)$ are (left and right) flat Koszul
for any flat $R$\+module~$V$.

 Let $V$ be a finitely generated projective $R$\+module everywhere
of rank~$m$ (i.~e., for every prime ideal $\mathfrak p$ in $R$,
the localization $V_{\mathfrak p}$ is a free $R_{\mathfrak p}$\+module
of rank~$m$).
 Then one has $\Lambda_R^n(V)=0$ for $n>m$, and $\Lambda_R(V)$ is
a relatively Frobenius Koszul graded ring over $R$ in the sense of
Section~\ref{relatively-frobenius-koszul-graded-subsecn}.
 In other words, $\Lambda_R^m(V)$ is an invertible $R$\+module and
the multiplication map $\Lambda_R^n(V)\ot_R\Lambda_R^{m-n}(V)
\rarrow\Lambda_R^m(V)$ is a perfect pairing for every $0\le n\le m$
(in the sense of Section~\ref{rel-frobenius-subsecn}).

\subsection{Algebraic differential operators}
\label{algebraic-diffoperators-subsecn}
 Let $X$ be a smooth affine algebraic variety over a field~$k$ (which
we will eventually assume to have characteristic~$0$ in this section).
 Denote by $O(X)$ the finitely generated commutative $k$\+algebra
of regular functions on~$X$.
 \emph{Differential operators on~$X$} (in the sense of Grothendieck)
form a subring in the ring $\End_k(O(X))$ of $k$\+linear endomorphisms
of the $k$\+vector space $O(X)$.

 The differential operators of order~$0$ are the $O(X)$\+linear
endomorphisms of $O(X)$, that is the operators of multiplication
with regular functions $f\in O(X)$.
 A $k$\+linear map $D\:O(X)\rarrow O(X)$ is said to be
a \emph{differential operator of order}~$\le n$ if, for any
$f\in O(X)$, the operator $[D,f]=D\circ f-f\circ D\:O(X)\rarrow
O(X)$ is a differential operator of order~$\le n-1$.
 A \emph{differential operator} $D\:O(X)\rarrow O(X)$ is a map which
is a differential operator of some finite order~$n\ge0$.

 We denote the subspace of differential operators by $\Diff(X)\subset
\End_k(O(X))$ and the subspace of differential operators of
order~$\le n$ by $F_n\Diff(X)\subset\Diff(X)$ for every integer $n\ge0$.
 It is straightforward to check that $\Diff(X)$ is a subring in
$\End_k(O(X))$ (with respect to the composition multiplication of
linear operators) and $F$ is a multiplicative increasing filtration on
$\Diff(X)=\bigcup_n F_n\Diff(X)$.

 Let $T$ denote the tangent bundle to $X$ and $T^*$ denote
the cotangent bundle.
 The global sections of the tangent bundle are called the \emph{vector
fields} on $X$ and the global sections of the cotangent bundle are
called the \emph{differential $1$\+forms}.

 The $k$\+vector space of vector fields $T(X)$ can be constructed as
the subspace $T(X)\subset F_1\Diff(X)$ of all differential
operators~$v$ of order~$\le 1$ on $X$ which annihilate the constant
functions, that is $v(1)=0$.
 Equivalently, $T(X)$ is the space of all \emph{derivations} of
the $k$\+algebra $O(X)$, that is $k$\+linear maps $v\:O(X)\rarrow O(X)$
such that $v(fg)=v(f)g+fv(g)$ for all $f$, $g\in O(X)$.
 The subspace $T(X)\subset F_1\Diff(X)$ is preserved by the left
(but not right) multiplications of differential operators by
the functions; so the rule $(fv)(g)=f(v(g))$ defines
an $O(X)$\+module structure on $T(X)$.

 The $O(X)$\+module of differential $1$\+forms $T^*(X)$ is produced by
the construction of \emph{K\"ahler differentials}.
 Consider the free $O(X)$\+module spanned by the symbols $d(f)$ with
$f\in O(X)$ and take its quotient module by the submodule spanned by
all elements of the form $d(fg)-fd(g)-gd(f)$ with $f$, $g\in O(X)$
and $d(a)$ with $a\in k$; this quotient module is $T^*(X)$.

 Since we are assuming that $X$ is smooth, both $T(X)$ and $T^*(X)$
are finitely generated projective $O(X)$\+modules (of the rank equal
to the dimension of $X$ over~$k$).
 They are also naturally dual to each other: one has
$T(X)=\Hom_{O(X)}(T^*(X),O(X))$.

 We are interested in the symmetric powers of the $O(X)$\+module
$T(X)$ and the exterior powers of the $O(X)$\+module $T^*(X)$.
 Following the notation in Section~\ref{symmetric-exterior-subsecn},
the former are denoted by $\Sym^n_{O(X)}T(X)$ and the latter by
$\Lambda^n_{O(X)}(T^*(X))=\Omega^n(X)$, where $n\ge0$ (so,
in particular, $\Sym^0_{O(X)}T(X)=O(X)=\Omega^0(X)$, \
$\Sym^1_{O(X)}T(X)=T(X)$, and $\Omega^1(X)=T^*(X)$).
 These are the global sections of the symmetric/exterior powers of
the vector bundles $T$ and $T^*$ on~$X$.
 The elements of $\Omega^n(X)$ are called the \emph{differential
$n$\+forms} on~$X$.

 There exists a unique odd derivation~$d$ of degree~$1$ with $d^2=0$
on the graded algebra $\Omega(X)$ whose restriction to
$\Omega^0(X)=O(X)$ is the map $d\:O(X)\rarrow\Omega^1(X)$ appearing in
the above construction of the $O(X)$\+module $\Omega^1(X)=T^*(X)$
as the module of K\"ahler differentials.
 The differential $d\:\Omega^n(X)\rarrow\Omega^{n+1}(X)$, \,$n\ge0$,
is called the \emph{de~Rham differential}.
 So $(\Omega(X),d)$ is a DG\+ring.
 More precisely, it is a DG\+algebra over~$k$; it is called
the \emph{de~Rham DG\+algebra}.

 When the characteristic of~$k$ is equal to~$0$, the associated
graded ring $\gr^F\Diff(X)=\bigoplus_{n=0}^\infty
F_n\Diff(X)/F_{n-1}\Diff(X)$ is naturally isomorphic to the graded
$O(X)$\+algebra $\Sym_{O(X)}T(X)$.
 Since the latter is left (and right) finitely projective Koszul,
the filtered ring $(\Diff(X),F)$ is left finitely projective
nonhomogeneous Koszul in the sense of
Section~\ref{pbw-theorem-subsecn}.
 Furthermore, the ring $\Diff(X)$ is left augmented over its subring
$O(X)$ in the sense of Section~\ref{augmented-subsecn}.
 The natural left action of $\Diff(X)$ in $O(X)$ (by the differential
operators) provides the augmentation.

 One can check that the left augmented left finitely projective
nonhomogeneous Koszul ring $\Diff(X)$ (with the above filtration
and augmentation) corresponds to the right finitely projective
Koszul DG\+ring $\Omega(X)$ (with the de~Rham differential) under
the anti-equivalence of categories from
Corollary~\ref{left-augmented-koszul-duality-anti-equivalence}.

 Notice that the commutative ring $R=O(X)$ has finite homological
dimension (since $X$ is a smooth algebraic variety by assumption).
 Furthermore, the graded ring $\Omega(X)$ has only finitely many
grading components (indeed, $\Omega^n(X)=0$ for $n>\dim_kX$); so
there is no difference between graded \emph{comodules},
graded \emph{contramodules}, and the conventional graded
\emph{modules} over $\Omega(X)$.
 Hence the results of Sections~\ref{comodule-side-secn}
and~\ref{contramodule-side-secn} lead to the following theorem.

\begin{thm} \label{algebraic-differential-koszul-duality}
 For any smooth affine algebraic variety $X$ over a field~$k$ of
characteristic~$0$, the construction of
Corollary~\ref{fin-dim-base-koszul-duality-comodule-side} provides
a triangulated equivalence between the derived category of right
modules over the ring of algebraic differential operators $\Diff(X)$
and the coderived category of right DG\+modules over the de~Rham
DG\+algebra $(\Omega(X),d)$,
\begin{equation} \label{algebraic-differential-comodule-side}
 \sD(\modr{\Diff(X)})\simeq\sD^\co(\modr(\Omega(X),d)).
\end{equation}
 In the same context, the construction of
Corollary~\ref{fin-dim-base-koszul-duality-contramodule-side}
provides an equivalence between the derived category of left
$\Diff(X)$\+modules and the contraderived category of left
DG\+modules over $(\Omega(X),d)$,
\begin{equation} \label{algebraic-differential-contramodule-side}
 \sD(\Diff(X)\modl)\simeq\sD^\ctr((\Omega(X),d)\modl).
\end{equation} \qed
\end{thm}

 The abelian categories of left and right modules over the ring
$\Diff(X)$ are naturally equivalent (see the discussion in
Section~\ref{introd-conversion}).
 The results of Section~\ref{conversion-secn} provide a square
diagram of triangulated equivalences
connecting~\eqref{algebraic-differential-comodule-side}
with~\eqref{algebraic-differential-contramodule-side}.

\begin{thm} \label{algebraic-differential-quadrality}
 For any smooth affine algebraic variety $X$ over a field~$k$ of
characteristic~$0$, the constructions of
Theorem~\ref{frobenius-quadrality-theorem} (with
Remark~\ref{commutative-cdg-ring-remark} taken into account)
provide a commutative square diagram of triangulated equivalences
\begin{equation}
\begin{tikzcd}
 \sD(\modr{\Diff(X)}) \arrow[r, bend right = 8] \arrow[dd] &
 \sD(\Diff(X)\modl) \arrow[l] \arrow[dd]
 \\ \\
 \sD^\co((\Omega(X),d)\modl)
 \arrow[r, Leftrightarrow, no head, no tail]
 \arrow[uu, bend left = 20] &
 \sD^\ctr((\Omega(X),d)\modl) \arrow[uu, bend right = 20]
\end{tikzcd}
\end{equation}
where the equivalence of derived categories in the upper line is
induced by the conversion equivalence of abelian categories
$\modr{\Diff(X)}\simeq\Diff(X)\modl$ (up to a cohomological shift
by $[-\dim_kX]$), while the co-contra correspondence in the lower line
is the result of Theorem~\ref{frobenius-co-contra-thm}.
 The vertical equivalences
are~\eqref{algebraic-differential-comodule-side}
and~\eqref{algebraic-differential-contramodule-side}.  \qed
\end{thm}

 Modules over the ring of differential operators $\Diff(X)$ (or more
generally, sheaves of modules over the sheaf of rings of differential
operators over a nonaffine smooth algebraic variety~$X$) are known
colloquially as ``$D$\+modules''~\cite{Bern,Ginz,BB}.
 An approach to the theory of $D$\+modules based on DG\+modules over
the de~Rham DG\+algebra was developed in the paper~\cite{Ryb}.

\subsection{Crystalline differential operators}
\label{crystalline-subsecn}
 Let $X$ be a smooth affine algebraic variety over a field~$k$ of
arbitrary characteristic.
 The ring of \emph{crystalline differential operators} $\CDiff(X)$
is defined by generators and relations
as follows~\cite[Section~1.2]{BMR}.

 The generators are the elements of the ring of functions $O(X)$
and the $O(X)$\+mod\-ule of vector fields $T(X)$.
 The sum of any two elements of $O(X)$ in $\CDiff(X)$ equals their sum
in $O(X)$; and the sum of any two elements of $T(X)$ in $\CDiff(X)$
equals their sum in $T(X)$.
 Concerning the product, let us use the notation of
Section~\ref{self-consistency-subsecn} and denote by~$*$ the product
of any two elements in $\CDiff(X)$, to be distinguished from
their product in $O(X)$ or $T(X)$.
 Then the relations
\begin{equation} \label{crystalline-diff-relations}
\begin{gathered}
 f*g=fg \qquad\text{for all $f$, $g\in O(X)$,} \\
 f*v=fv \qquad\text{for all $f\in O(X)$ and $v\in T(X)$,} \\
 v*f=fv+v(f) \qquad\text{for all $f\in O(X)$ and $v\in T(X)$}
\end{gathered}
\end{equation}
are imposed.
 Here $fg\in O(X)$ denotes the product in $O(X)$ and $fv\in T(X)$
denotes the action of elements of $O(X)$ in the $O(X)$\+module $T(X)$,
while $v(f)\in O(X)$ denotes the action of vector fields by
differential operators (or more precisely derivations) on the functions;
so $v(f)$ is ``the derivative of~$f$ along~$v$''.

 Finally, we need to recall that derivations form a Lie algebra:
for any vector fields $v$ and $w\in T(X)$, there exists a unique
vector field $[v,w]\in T(X)$, called the \emph{commutator} of $v$
and~$w$, such that $[v,w](f)=v(w(f))-w(v(f))$ for all $f\in O(X)$.
 The relation
\begin{equation} \label{commutator-vector-fields-relation}
 v*w-w*v=[v,w] \qquad\text{for all $v$, $w\in T(X)$}
\end{equation}
is also imposed in $\CDiff(X)$.

 One considers the increasing filtration $F$ on the ring $\CDiff(X)$
generated by $F_1\CDiff(X)=O(X)\oplus T(X)$ over $F_0\CDiff(X)=O(X)$
(cf.\ the discussion of generated filtrations in
Section~\ref{nonhomogeneous-quadratic-subsecn}).
 Then, irrespectively of the characteristic of~$k$, the associated
graded ring $\gr^F\CDiff(X)=\bigoplus_{n=0}^\infty
F_n\CDiff(X)/F_{n-1}\CDiff(X)$ is naturally isomorphic to
the symmetric algebra $\Sym_{O(X)}T(X)$ (this is provable as
a particular case of Theorem~\ref{pbw-theorem-thm}).

 There is a natural homomorphism of filtered rings $\CDiff(X)\rarrow
\Diff(X)$ uniquely defined by the condition that it acts by
the identity maps on the subring $O(X)\subset\CDiff(X)$ (taking
it to the subring $O(X)\subset\Diff(X)$) and on the subspace
$T(X)\subset\CDiff(X)$ (taking it to the subspace $T(X)\subset
\Diff(X)$).
 Over a field~$k$ of characteristic~$0$, this is a ring isomorphism,
and in fact an isomorphism of filtered rings.
 But it is \emph{neither} surjective \emph{nor} injective in prime
characteristic.

 To give an example of noninjectivity, let $Y$ be the affine line
over a field~$k$ of prime characteristic~$p$; so $O(Y)=k[y]$ is
the ring of polynomials in one variable.
 Then $d/dy\:k[y]\rarrow k[y]$ is a vector field on~$Y$.
 One would expect $(d/dy)^p\:k[y]\rarrow k[y]$ to be a differential
operator of order~$p$, but in fact it is a zero map.
 So $(d/dy)^p\in\CDiff(Y)$ is an element of $F_p\CDiff(Y)$ not
belonging to $F_{p-1}\CDiff(Y)$, but belonging to the kernel of
the ring homomorphism $\CDiff(Y)\rarrow\Diff(Y)$.

 To give an example of nonsurjectivity, consider the ring of polynomials
with integer coefficients $R=\boZ[y]$.
 Then, for any element $f\in\boZ[y]$, the element $f^{(p)}(y)=
d^pf/dy^p\in\boZ[y]$ is divisible by~$p$; so $\frac{1}{p}(d/dy)^p$ is
a well-defined map $R\rarrow R$.
 Taking the tensor product $k\ot_{\boZ}{-}$, one obtains a differential
operator of order~$p$ on $k[y]$ which can be denoted by
``$\frac{1}{p}(d/dy)^p$\,''$\:k[y]\rarrow k[y]$.
 This is an element of $F_p\Diff(Y)$ not belonging to the sum of
$F_{p-1}\Diff(Y)$ with the image of the ring homomorphism
$\CDiff(Y)\rarrow\Diff(Y)$.

 To give another example, let $Z$ be the punctured affine line
over~$k$; so $O(Z)=k[z,z^{-1}]$.
 Then $z\frac{d}{dz}$ is a vector field on~$Z$.
 One would expect $\bigl(z\frac{d}{dz}\bigr)^p\:k[z,z^{-1}]\rarrow
k[z,z^{-1}]$ to be a differential operator of order~$p$, but in fact
it is the same map as $z\frac{d}{dz}$.
 So $\bigl(z\frac{d}{dz}\bigr)^p-z\frac{d}{dz}\in\CDiff(Z)$ is
an element of $F_p\CDiff(Z)$ not belonging to $F_{p-1}\CDiff(Z)$,
but belonging to the kernel of the ring homomorphism $\CDiff(Z)
\rarrow\Diff(Z)$.
 Similarly to the construction above, one can define a differential
operator ``$\frac{1}{p}\bigl(\bigl(z\frac{d}{dz}\bigr)^p-
z\frac{d}{dz}\bigr)$''$\:k[z,z^{-1}]\rarrow k[z,z^{-1}]$.
 This is an element of $F_p\Diff(Z)$ not belonging to the sum of
$F_{p-1}\Diff(Z)$ with the image of the ring homomorphism
$\CDiff(Z)\rarrow\Diff(Z)$.

 Returning to the general case, let point out that the map
$F_1\CDiff(X)\rarrow F_1\Diff(X)$ is still an isomorphism,
for any smooth affine variety $X$ over a field $k$ of any
characteristic.
 In fact, the map $F_{p-1}\CDiff(X)\rarrow F_{p-1}\Diff(X)$ is
an isomorphism when $k$~has characteristic~$p$.
 But the ring $\Diff(X)$ is not generated by $F_1\Diff(X)$ in
the latter case: the subring in $\Diff(X)$ generated by $F_1\Diff(X)$
does not contain $F_p\Diff(X)$ (whenever the dimension of $X$ is
more than zero); moreover, the ring $\Diff(X)$ is not finitely
generated~\cite[Section~3]{Smi}.
 The ring $\CDiff(X)$ has very different properties: over a field~$k$
of prime characteristic, $\CDiff(X)$ is finitely generated as
a module over its center~\cite[Sections~1.3 and~2]{BMR}.

 On the other hand, $(\CDiff(X),F)$ is a left finitely projective
nonhomogeneous Koszul ring in the sense of
Section~\ref{pbw-theorem-subsecn}, irrespectively of
the characteristic of~$k$.
 Composing the action of $\Diff(X)$ in $O(X)$ with the ring
homomorphism $\CDiff(X)\rarrow\Diff(X)$, one defines a left action
of $\CDiff(X)$ in $O(X)$ making $\CDiff(X)$ a left augmented ring
over its subring $O(X)$.
 The left augmented left finitely projective nonhomogeneous Koszul
ring of crystalline differential operators $\CDiff(X)$ (with
the above filtration and augmentation) corresponds to the right
finitely projective Koszul DG\+ring of differential forms
$\Omega(X)$ (with the above differential) under the anti-equivalence
of categories from
Corollary~\ref{left-augmented-koszul-duality-anti-equivalence}.

 Similarly to the characteristic~$0$ case of
Section~\ref{algebraic-diffoperators-subsecn}, the ring of functions
$R=O(X)$ on any smooth affine algebraic variety over a field~$k$ of
any characteristic has finite homological dimension, and the graded
ring $\Omega(X)$ has only finitely many grading components.
 Hence we have the following generalization of
Theorem~\ref{algebraic-differential-koszul-duality}.

\begin{thm} \label{crystalline-differential-koszul-duality}
 For any smooth affine algebraic variety $X$ over a field~$k$,
the construction of
Corollary~\ref{fin-dim-base-koszul-duality-comodule-side} provides
a triangulated equivalence between the derived category of right
modules over the ring of crystalline differential operators
$\CDiff(X)$ and the coderived category of right DG\+modules over
the de~Rham DG\+algebra $(\Omega(X),d)$,
\begin{equation} \label{crystalline-differential-comodule-side}
 \sD(\modr{\CDiff(X)})\simeq\sD^\co(\modr(\Omega(X),d)).
\end{equation}
 In the same context, the construction of
Corollary~\ref{fin-dim-base-koszul-duality-contramodule-side}
provides an equivalence between the derived category of left
$\CDiff(X)$\+modules and the contraderived category of left
DG\+modules over $(\Omega(X),d)$,
\begin{equation} \label{crystalline-differential-contramodule-side}
 \sD(\CDiff(X)\modl)\simeq\sD^\ctr((\Omega(X),d)\modl).
\end{equation} \qed
\end{thm}

 Similarly to Section~\ref{algebraic-diffoperators-subsecn},
there is a natural equivalence between the abelian categories of
left and right modules over $\CDiff(X)$, provided by the mutually
inverse functors of tensor product with $\Omega^m(X)$ and
$\Lambda^m_{O(X)}(T(X))$ over $O(X)$ (where $m=\dim_kX$).
 The results of Section~\ref{conversion-secn} provide a comparison
between~\eqref{crystalline-differential-comodule-side}
and~\eqref{crystalline-differential-contramodule-side}.

\begin{thm} \label{crystalline-differential-quadrality}
 For any smooth affine algebraic variety $X$ over a field~$k$,
the constructions of Theorem~\ref{frobenius-quadrality-theorem}
(with Remark~\ref{commutative-cdg-ring-remark} taken into account)
provide a commutative square diagram of triangulated equivalences
\begin{equation}
\begin{tikzcd}
 \sD(\modr{\CDiff(X)}) \arrow[r, bend right = 8] \arrow[dd] &
 \sD(\CDiff(X)\modl) \arrow[l] \arrow[dd]
 \\ \\
 \sD^\co((\Omega(X),d)\modl)
 \arrow[r, Leftrightarrow, no head, no tail]
 \arrow[uu, bend left = 20] &
 \sD^\ctr((\Omega(X),d)\modl) \arrow[uu, bend right = 20]
\end{tikzcd}
\end{equation}
where the equivalence of derived categories in the upper line is
induced by the conversion equivalence of abelian categories
$\modr{\CDiff(X)}\simeq\CDiff(X)\modl$ (up to a cohomological shift
by $[-\dim_kX]$), while the co-contra correspondence in the lower line
is the result of Theorem~\ref{frobenius-co-contra-thm}.
 The vertical equivalences
are~\eqref{crystalline-differential-comodule-side}
and~\eqref{crystalline-differential-contramodule-side}.
\end{thm}

\begin{proof}
 Notice that $\Omega(X)$ is a relatively Frobenius Koszul graded ring
(as per the discussion at the end of
Section~\ref{symmetric-exterior-subsecn}); so
Theorem~\ref{frobenius-quadrality-theorem} is indeed applicable.
\end{proof}

\subsection{Differential operators in a vector bundle}
\label{vector-bundle-subsecn}
 Let $E$ be a vector bundle over a smooth affine algebraic variety $X$
(over a field~$k$).
 Then the space of global sections $E(X)$ of the vector bundle $E$ is
a finitely generated projective $O(X)$\+module.

 Given two vector bundles $E'$ and $E''$ over $X$, one can consider
differential operators acting from $E'(X)$ to $E''(X)$.
 The space of such differential operators $\Diff(X,E',E'')$ is
a filtered $k$\+vector subspace in $\Hom_k(E'(X),E''(X))$.

 Specifically, the differential operators of order~$0$ are
the $O(X)$\+linear maps $E'(X)\allowbreak\rarrow E''(X)$.
 A $k$\+linear map $D\:E'(X)\rarrow E''(X)$ is said to be
a \emph{differential operator of order}~$\le n$ if, for every
regular function $f\in O(X)$, the map $[D,f]=D\circ f-f\circ\nobreak D
\:E'(X)\rarrow E''(X)$ is a differential operator of order~$\le n-1$.
 Here $f$~acts in $E'(X)$ and $E''(X)$ as in $O(X)$\+modules.
 A \emph{differential operator} $E'(X)\rarrow E''(X)$ is a map which
is a differential operator of some finite order $n\ge0$.

 We denote the subspace of differential operators of order~$\le n$
by $F_n\Diff(X,E',E'')\allowbreak\subset\Diff(X,E',E'')$.
 When the two vector bundles are the same, $E'=E=E''$, we write simply
$\Diff(X,E)$ instead of $\Diff(X,E',E'')$.
 The $k$\+vector space $\Diff(X,E)$ with its increasing filtration $F$
is a filtered ring and a subring in $\End_k(E(X))$.
 The subring of differential operators of order~$0$ in $\Diff(X,E)$ is
$F_0\Diff(X,E)=\Hom_{O(X)}(E(X),E(X))$; it can be described as
the ring of global sections $\End(E)(X)$ of the vector bundle $\End(E)$
of endomorphisms of the vector bundle $E$ over~$X$.

 Assuming that the characteristic of~$k$ is equal to~$0$,
the associated graded ring $\gr^F\Diff(X,E)$ is naturally isomorphic
to the tensor product of the $O(X)$\+algebra $\End(E)(X)$ and
the graded $O(X)$\+algebra $\Sym_{O(X)}(T(X))$,
$$
 \gr^F\Diff(X,E)\simeq\End(E)(X)\ot_{O(X)}\Sym_{O(X)}(T(X)).
$$
 Both $\End(E)(X)$ and the grading components of $\Sym_{O(X)}(T(X))$
are finitely generated projective $O(X)$\+modules, while
the filtration components of the ring $\Diff(X,E)$ have the left
and the right $O(X)$\+module structures, induced by the inclusion
$O(X)\rarrow\End(E)(X)=F_0\Diff(X,E)$.
 In particular, we have a short exact sequence of
$O(X)$\+$O(X)$\+bimodules
\begin{equation} \label{diff-operators-F0-F1-sequence}
 0\lrarrow\End(E)(X)\lrarrow F_1\Diff(X,E)\lrarrow
 \End(E)(X)\ot_{O(X)}T(X)\lrarrow0.
\end{equation}
 In fact, this is even a short exact sequence of bimodules over
the (noncommutative) $O(X)$\+algebra $\End(E)(X)$, which exists
irrespectively of the characteristic of~$k$ (cf.\ the discussion
in Section~\ref{crystalline-subsecn}).

 A \emph{connection} $\nabla$ in a vector bundle $E$ is a splitting
of~\eqref{diff-operators-F0-F1-sequence} as a short exact sequence
of \emph{left} $\End(E)(X)$\+modules.
 All the three terms of~\eqref{diff-operators-F0-F1-sequence} are
projective as left (as well as right) $\End(E)(X)$\+modules; so
a connection in a vector bundle $E$ over a smooth affine algebraic
variety $X$ always exists.
 In fact, the connections in $E$ form an affine space (in a different
language, a principal homogenenous space) over the vector space
$T^*(X)\ot_{O(X)}\End(E)(X)$.
 In other words, the difference $\nabla''-\nabla'$ of any two
connections in $E$ is an element of $T^*(X)\ot_{O(X)}\End(E)(X)$,
and conversely, to any connection $\nabla$ one can add any element of
$T^*(X)\ot_{O(X)}\End(E)(X)$ and obtain a new connection in~$E$.

 The inclusion of rings $O(X)\rarrow\End(E)(X)$ induces an inclusion
of $O(X)$\+modules $T(X)\rarrow\End(E)(X)\ot_{O(X)}T(X)$.
 Taking the pull-back of the short exact
sequence~\eqref{diff-operators-F0-F1-sequence} with respect to
the latter map, we obtain a short exact sequence of
$O(X)$\+$O(X)$\+bimodules
\begin{equation} \label{reduced-F0-F1-sequence}
 0\lrarrow\End(E)(X)\lrarrow \overline{F}_1\Diff(X,E)\lrarrow
 T(X)\lrarrow0.
\end{equation}
 Here $\overline{F}_1\Diff(X,E)$ is a certain
$O(X)$\+$O(X)$\+subbimodule in $F_1\Diff(X,E)$.
 A splitting of~\eqref{diff-operators-F0-F1-sequence} as a short
exact sequence of left $\End(E)(X)$\+modules is equivalent to
a left $O(X)$\+linear splitting of~\eqref{reduced-F0-F1-sequence}.

 We arrive to the definition of a connection $\nabla$ in $E$ as
a map assigning to every vector field $v\in T(X)$ a $k$\+linear
operator $\nabla_v\:E(X)\rarrow E(X)$ in such a way that
the following two equations are satisfied:
\begin{enumerate}
\renewcommand{\theenumi}{\roman{enumi}}
\item $\nabla_{fv}(e)=f\nabla_v(e)$ for all $f\in O(X)$, \,$v\in T(X)$,
and $e\in E(X)$;
\item $\nabla_v(fe)=v(f)e+f\nabla_v(e)$ for all
$f\in O(X)$, \,$v\in T(X)$, and $e\in E(X)$.
\end{enumerate}
 Here $v(f)\in O(X)$ is the derivative of~$f$ along~$v$.

 A connection $\nabla$ in $E$ can be also interpreted as a $k$\+linear
map $E(X)\rarrow\Omega^1(X)\ot_{O(X)}E(X)$ defined by the rule
$\lan v,\nabla(e)\ran =\nabla_v(e)$, where $\lan\ , \ \ran$ denotes
the $O(X)$\+linear map $T(V)\ot_{O(X)}T^*(X)\ot_{O(X)}E(X)\rarrow E(X)$
induced by the natural pairing $T(X)\ot_{O(X)}T^*(X)\rarrow O(X)$.
 Then the identity~(ii) takes the form $\nabla(fe)=d(f)\ot e
+f\nabla(e)$, where $d$~denotes the de~Rham differential
$d\:O(X)\rarrow\Omega^1(X)$.

 If vector bundles $E'$ and $E''$ over $X$ are endowed with
connections $\nabla'$ and $\nabla''$, then any vector bundle
produced naturally from $E'$ and $E''$, such as $E'\oplus E''$
and $E'\ot E''$, acquires an induced connection.
 In particular, the vector bundle $E'\ot E''$ is defined
by the rule $(E'\ot E'')(X)=E'(X)\ot_{O(X)}E''(X)$, and the induced
connection $\nabla$ on $E'\ot E''$ is given by the rule
$\nabla_v(e'\ot e'')=\nabla'_v(e')\ot e''+e'\ot\nabla''_v(e'')$
for all $e'\in E'(X)$, \ $e''\in E''(X)$, and $v\in T(X)$.
 Similarly, the vector bundle $\Hom(E',E'')$ is defined by
the rule $\Hom(E',E'')(X)=\Hom_{O(X)}(E'(X),E''(X))$, and
the induced connection $\nabla$ on $\Hom(E',E'')$ is given
by the rule $\nabla_v(g)(e')=\nabla''_v(g(e'))-g(\nabla'_v(e'))$
for all $g\in\Hom_{O(X)}(E'(X),E''(X))$, \ $e'\in E'(X)$
and $v\in T(X)$.

 Let $\nabla=\nabla_E$ be a connection in a vector bundle $E$ over~$X$.
 Consider the graded left $\Omega(X)$\+module $\Omega(X)\ot_{O(X)}E(X)$
of \emph{differential forms on $X$ with the coefficients in}~$E$.
 Then the connection $\nabla$ on $E$ induces an odd
derivation~$d_\nabla$ on the graded module $\Omega(X)\ot_{O(X)}E(X)$
compatible with the odd derivation~$d$ (the de~Rham differential)
on the graded ring $\Omega(X)$, in the sense of
Section~\ref{cdg-modules-subsecn}.
 The map $d_\nabla\:\Omega^n(X)\ot_{O(X)}\nobreak E(X)\allowbreak
\rarrow\Omega^{n+1}(X)\ot_{O(X)}\nobreak E(X)$ is given by the formula
\begin{equation} \label{connection-de-Rham}
 d_\nabla(\omega\ot e)=d(\omega)\ot e+(-1)^n\omega\wedge\nabla(e)
\end{equation}
for all $\omega\in\Omega^n(X)$ and $e\in E(X)$.
 Here the wedge~$\wedge$ denotes the map $\Omega^n(X)\ot_k\Omega^1(X)
\ot_{O(X)}E(X)\rarrow\Omega^{n+1}(X)\ot_{O(X)}E(X)$ induced by
the multiplication map $\Omega^1(X)\ot_k\Omega^n(X)\rarrow
\Omega^{n+1}(X)$.
 One needs to check that the map $d_\nabla=d_{\nabla_E}$ is
well-defined, that is $d_\nabla(f\omega\ot e)=d_\nabla(\omega\ot fe)$
for all $f\in O(X)$.

 The square~$d_\nabla^2$ of the differential~$d_\nabla$ on
$\Omega(X)\ot_{O(X)}E(X)$ is an $\Omega(X)$\+linear map.
 Hence there exists an element $h_\nabla\in
\Omega^2(X)\ot_{O(X)}\End(E)(X)$ such that $d_\nabla^2(\phi)=h(\phi)$
for all $\phi\in\Omega(X)\ot_{O(X)}E(X)$.
 Here the left action of $\Omega(X)\ot_{O(X)}\End(E)(X)$ in
$\Omega(X)\ot_{O(X)}E(X)$ is induced by the multiplication in $O(X)$
and the left action of $\End(E)(X)$ in $E(X)$.
 The element $h_\nabla=h_{\nabla_E}\in\Omega^2(X)\ot_{O(X)}\End(E)(X)$
is called the \emph{curvature} of the connection~$\nabla$ in
a vector bundle~$E$.

 The tensor product $\Omega(X)\ot_{O(X)}\End(E)(X)$ of the graded
$O(X)$\+algebra $\Omega(X)$ with the $O(X)$\+algebra $\End(E)(X)$
has a natural structure of graded $O(X)$\+algebra.
 The connection $\nabla=\nabla_E$ on the vector bundle $E$ induces
a connection $\nabla_{\End(E)}$ on the vector bundle $\End(E)$ on $X$,
as explained above.
 The induced differential $d_{\nabla_{\End(E)}}$ on the graded
$\Omega(X)$\+module $\Omega(X)\ot_{O(X)}\End(E)(X)$ is, in fact,
an odd derivation of the graded algebra $\Omega(X)\ot_{O(X)}\End(E)(X)$
(since the connection $\nabla_{\End(E)}$ is compatible with
the composition multiplication on $\End(E)$, in the appropriate sense).
 The triple
\begin{equation} \label{connection-de-Rham-cdg-ring}
 (\Omega(X)\ot_{O(X)}\End(E)(X),\>d_{\nabla_{\End(E)}},\>h_{\nabla_E})
\end{equation}
is a curved DG\+ring.
 The pair
\begin{equation} \label{connection-de-Rham-cdg-module}
 (\Omega(X)\ot_{O(X)} E(X),\>d_{\nabla_E})
\end{equation}
is a left CDG\+module over~\eqref{connection-de-Rham-cdg-ring}.

 The graded ring $A=\End(E)(X)\ot_{O(X)}\Sym_{O(X)}(T(X))$ is
left (and right) finitely projective Koszul (over its degree-zero
component $R=\End(E)(X)$).
 The quadratic dual right (and left) finitely projective Koszul
graded ring to $A$, as per the construction of
Propositions~\ref{2-fin-proj-quadratic-duality}
and~\ref{3-fin-proj-quadratic-duality},
is $B=\Omega(X)\ot_{O(X)}\End(E)(X)$.

 Assume that the characteristic of~$k$ is equal to~$0$.
 Then the filtered ring $(\Diff(X,E),F)$ is a left finitely projective
nonhomogeneous Koszul ring.
 The choice of a connection~$\nabla$ on $E$ means the choice of
a left $\End(E)(X)$\+linear splitting of the short exact
sequence~\eqref{diff-operators-F0-F1-sequence}; in the terminology
of Section~\ref{self-consistency-subsecn}, this is the choice
of a submodule of strict generators for $\Diff(X,E)$.
 The corresponding right finitely projective Koszul CDG\+ring produced
by the construction of Proposition~\ref{nonhomogeneous-dual-cdg-ring}
is the CDG\+ring
$\Omega(X)\ot_{O(X)}\End(E)(X)$~\,\eqref{connection-de-Rham-cdg-ring}.
 Replacing a connection~$\nabla_E$ in $E$ with another
connection~$\nabla'_E$ leads to a CDG\+ring
$(\Omega(X)\ot_{O(X)}\End(E)(X),\>d_{\nabla'_{\End(E)}},
\>h_{\nabla'_E})$ connected with~\eqref{connection-de-Rham-cdg-ring}
by a natural change-of-connection isomorphism of CDG\+rings, 
as per the discussion in Sections~\ref{curved-dg-rings-subsecn}
and~\ref{change-of-strict-gens-subsecn}.

 The ring of endomorphisms $R=\End(E)(X)$ of any vector bundle $E$ on
a smooth affine algebraic variety $X$ has finite left and right
homological dimensions; in fact, assuming that $E$ is nonzero on all
the connected components of $X$, both the abelian categories of left
and right $R$\+modules are equivalent to $O(X)\modl$.
 Furthermore, just as in Section~\ref{algebraic-diffoperators-subsecn},
the graded ring $\Omega(X)\ot_{O(X)}\End(E)(X)$ has only finitely
many grading components.
 Hence we obtain the following theorem.
\begin{thm} \label{vector-bundle-diffoperators-koszul-duality}
 For any smooth affine algebraic variety $X$ over a field~$k$ of
characteristic~$0$ and any vector bundle $E$ over $X$ with a chosen
connection~$\nabla_E$, the construction of
Corollary~\ref{fin-dim-base-koszul-duality-comodule-side} provides
a triangulated equivalence between the derived category of right
modules over the ring of differential operators $\Diff(X,E)$ acting
in the sections of $E$ and the coderived category of right CDG\+modules
over the CDG\+algebra~\eqref{connection-de-Rham-cdg-ring} of
differential forms on $X$ with the coefficients in $\End(E)$,
\begin{equation} 
 \sD(\modr{\Diff(X,E)})\,\simeq\,\sD^\co(\modr
 (\Omega(X)\ot_{O(X)}\End(E)(X),\,d_{\nabla_{\End(E)}},\,h_{\nabla_E})).
\end{equation}
 In the same context, the construction of
Corollary~\ref{fin-dim-base-koszul-duality-contramodule-side}
provides an equivalence between the derived category of left
$\Diff(X,E)$\+modules and the contraderived category of left
CDG\+modules over the CDG\+ring~\eqref{connection-de-Rham-cdg-ring}
\begin{equation} \label{vector-bundle-contramodule-side}
 \sD(\Diff(X,E)\modl)\,\simeq\,
 \sD^\ctr((\Omega(X)\ot_{O(X)}\End(E)(X),\,d_{\nabla_{\End(E)}},
 \,h_{\nabla_E})\modl).
\end{equation} \qed
\end{thm}

 For example, the equivalence of
categories~\eqref{vector-bundle-contramodule-side}
assigns the left CDG\+module $(\Omega(X)\ot_{O(X)}E(X),\,d_{\nabla_E})$
\,\eqref{connection-de-Rham-cdg-module} over
the CDG\+ring~\eqref{connection-de-Rham-cdg-ring} to the left
$\Diff(X,E)$\+mod\-ule~$E(X)$.
 The result of Theorem~\ref{vector-bundle-diffoperators-koszul-duality}
is the affine, characteristic~$0$ particular case of
the $\mathcal D$--$\Omega$ duality theorem
of~\cite[Theorem~B.2]{Pkoszul}.

\subsection{Twisted differential operators}
\label{twisted-differential-subsecn}
 Let us start with specializing the discussion in
Section~\ref{vector-bundle-subsecn} to the case of a \emph{line bundle}
$E=L$.
 In this case $L(X)$ in an \emph{invertible} finitely generated
projective $O(X)$\+module; specifically, one has
$L(X)\ot_{O(X)}L^*(X)\simeq O(X)$, where $L^*$ is the dual line
bundle to~$L$.
 Notice that for any vector bundle $E$ one has $End(X)=E\ot E^*$,
so $\End(E)(X)=E(X)\ot_{O(X)}E^*(X)$ (where $E^*$ is the dual
vector bundle to~$E$, so $E^*(X)=\Hom_{O(X)}(E(X),O(X))$).
 For a line bundle $L$, this means that $\End(L)$ is the trivial
line bundle, $\End(L)(X)=O(X)$.

 Choose a connection $\nabla=\nabla_L$ in~$L$.
 Then the induced connection $\nabla_{\End(L)}$ is the trivial
(canonical) connection in the trivial line bundle $\End(L)$.
 Therefore, the graded ring $\Omega(X)\ot_{O(X)}\End(L)(X)$ is simply
the ring of differential forms $\Omega(X)$, and
the differential~$d_{\nabla_{\End(L)}}$ in
the CDG\+ring~\eqref{connection-de-Rham-cdg-ring} is equal to
the standard de~Rham differential, $d_{\nabla_{\End(L)}}=d$.
 However, the curvature form $h_{\nabla_L}\in\Omega^2(X)$ of
the connection~$\nabla_L$ in $L$ can well be nontrivial.
 It is always a \emph{closed} differential form:
$d(h_{\nabla_L})=0$ in $\Omega^3(X)$ (cf.\ the equation~(iii) in
the definition of a CDG\+ring in Section~\ref{curved-dg-rings-subsecn}).

 Replacing $\nabla_L$ with another connection $\nabla'_L$ in $L$
replaces the differential $2$\+form $h_{\nabla_L}$ with
$h_{\nabla'_L}=h_{\nabla_L}+d(\alpha)$, where
$\alpha=\nabla'_L-\nabla_L\in\Omega^1(X)$ is
the change-of-connection $1$\+form (see the discussion
in Section~\ref{vector-bundle-subsecn};
cf.\ Section~\ref{change-of-strict-gens-subsecn}).
 So the cohomology class of the $2$\+form $h_{\nabla_L}$, viewed
as an element of the cohomology ring of the DG\+ring
$(\Omega(X),d)$, does not depend on a connection $\nabla_L$, but
only on the line bundle $L$ itself.
 It is called the \emph{first Chern class} of $L$ and denoted by
$c_1(L)\in H^2(\Omega(X),d)$.

 For any integer $n\in\boZ$, one can consider the line bundle
$L^{\ot n}$ over $X$, defined in the obvious way for $n\ge0$
and by the rule $L^{\ot n}=L^*{}^{\ot\,-n}$ for $n\le0$.
 The curvature form of the induced connection $\nabla_{L^{\ot n}}$
in $L^{\ot n}$ is given by the rule $h_{\nabla_{L^{\ot n}}}=
nh_{\nabla_L}\in\Omega^2(X)$.

 Now assume for a moment that $k$~is a field of characteristic~$0$,
and let $z\in k$ be any element.
 Then there is \emph{no} such thing as ``a line bundle $L^{\ot z}$
over~$X$''.
 However, one \emph{can} define a filtered ring $\Diff(X,L^{\ot z})$
of ``differential operators acting in the sections of $L^{\ot z}$\,''.
 For this purpose, one simply chooses a connection $\nabla$ in~$L$,
considers the right finitely projective Koszul CDG\+ring
$(\Omega(X),d,zh_\nabla)$, and constructs
$\Diff(X,L^{\ot z})$ as the left finitely projective nonhomogeneous
Koszul ring corresponding to $(\Omega(X),d,zh_\nabla)$ under
the anti-equivalence of categories from
Corollary~\ref{nonhomogeneous-koszul-duality-anti-equivalence}.

 Replacing the connection~$\nabla$ with another connection~$\nabla'$
in~$L$ corresponds to a natural change-of-connection isomorphism
between the CDG\+rings $(\Omega(X),d,h_\nabla)$ and
$(\Omega(X),d,h_{\nabla'})$, and consequently also
a change-of-connection isomorphism between the CDG\+rings
$(\Omega(X),d,zh_\nabla)$ and $(\Omega(X),d,zh_{\nabla'})$.
 This simply means that one has $h_{\nabla'_L}=h_{\nabla_L}+d(\alpha)$,
and consequently $zh_{\nabla'_L}=zh_{\nabla_L}+d(z\alpha)$, where
$\alpha=\nabla'_L-\nabla_L$ (as $\alpha^2=0$ in $\Omega^2(X)$).
 An isomorphism of right finitely projective Koszul CDG\+rings
induces an isomorphism of the corresponding left finitely projective
nonhomogeneous Koszul rings; so the filtered ring 
$(\Diff(X,L^{\ot z}),F)$ is defined uniquely up to a natural
isomorphism.
 Similarly, given two line bundles $L_1$ and $L_2$ over $X$ and two
scalars $z_1$ and $z_2\in k$, one can define a filtered ring
$\Diff(X,\>L_1^{\ot z_1}\ot L_2^{\ot z_2})$, etc.

 Quite generally, let $X$ be a smooth affine variety over a field~$k$
of arbitrary characteristic, and let $h\in\Omega^2(X)$ be a closed
$2$\+form (i.~e., $d(h)=0$).
 Then the triple $(\Omega(X),d,h)$, with the standard de Rham
differential~$d$ and the curvature element~$h$, is a curved DG\+ring
(because the graded ring $\Omega(X)$ is graded commutative, so in
particular $h$~is a central element in $\Omega(X)$ and the equation~(ii)
of Section~\ref{curved-dg-rings-subsecn} is satisfied).
 The curved DG\+ring $(\Omega(X),d,h)$ is right finitely projective
Koszul.
 The corresponding left finitely projective nonhomogeneous Koszul ring
under the anti-equivalence of categories from
Corollary~\ref{nonhomogeneous-koszul-duality-anti-equivalence} is
denoted by $\CDiff(X,h)$ and called the \emph{ring of twisted
crystalline differential operators on}~$X$ (twisted by~$h$).
 (See~\cite[Section~2]{BB} or~\cite[Chapter~II]{Ginz} for a much more
abstract discussion of twisted differential operators over nonaffine
varieties.)

 Explicitly, $\CDiff(X,h)$ is the filtered ring generated by
the functions on $X$ (placed in the filtration component~$F_0$)
and the vector fields on $X$ (placed in the filtration component~$F_1$),
subject to the same relations as in
Section~\ref{crystalline-subsecn}, except that
the relation~\eqref{commutator-vector-fields-relation} is replaced with
\begin{equation} \label{twisted-commutator-vector-fields}
 v*w-w*v=[v,w]+\lan v\wedge w,\,h\ran,
\end{equation}
where $\lan\ , \ \ran$ denotes the natural pairing
$\Lambda^2_{O(X)}(T(X))\ot_{O(X)}\Lambda^2_{O(X)}(T^*(X))\rarrow O(X)$
given by the rule $\lan v\wedge w,\>\alpha\wedge\beta\ran=
\lan v,\alpha\ran\lan w,\beta\ran-\lan v,\beta\ran\lan w,\alpha\ran$
for $v$, $w\in T(X)$ and $\alpha$, $\beta\in T^*(X)$
(cf.\ formula~\eqref{maps-p-and-h} in
Section~\ref{self-consistency-subsecn}).
 So $[v,w]\in T(X)$ is a vector field and $\lan v\wedge w,\>h\ran
\in O(X)$ is a function on~$X$.
 It is claimed, based on Theorem~\ref{pbw-theorem-thm}, that
the associated graded ring $\gr^F\CDiff(X,h)=
\bigoplus_{n=0}^\infty F_n\CDiff(X,h)/F_{n-1}\CDiff(X,h)$ is
naturally isomorphic to the symmetric algebra $\Sym_{O(X)}(T(X))$.
 The assumption that $d(h)=0$ is needed for this to be true.

 When two closed $2$\+forms $h'$ and $h''\in\Omega^2(X)$ represent
the same de~Rham cohomology class in $H^2(\Omega(X),d)$, the related
filtered rings of twisted (crystalline) differential operators
$\CDiff(X,h')$ and $\CDiff(X,h'')$ are isomorphic, but
\emph{not yet naturally isomorphic}.
 The choice of a $1$\+form $\alpha\in\Omega^1(X)$ such that $h''-h'=
d(\alpha)$ leads to a concrete isomorphism of filtered rings
$\CDiff(X,h')\simeq\CDiff(X,h'')$.

 The results of Sections~\ref{comodule-side-secn}
and~\ref{contramodule-side-secn} lead to the following theorem.

\begin{thm} \label{twisted-differential-koszul-duality}
 For any smooth affine algebraic variety $X$ over a field~$k$
and any closed differential $2$\+form $h\in\Omega^2(X)$, \,$d(h)=0$,
the construction of
Corollary~\ref{fin-dim-base-koszul-duality-comodule-side} provides
a triangulated equivalence between the derived category of right
modules over the ring of twisted (crystalline) differential operators
$\CDiff(X,h)$ and the coderived category of right CDG\+modules over
the de~Rham CDG\+algebra $(\Omega(X),d,h)$,
\begin{equation} \label{twisted-differential-comodule-side}
 \sD(\modr{\CDiff(X,h)})\simeq\sD^\co(\modr(\Omega(X),d,h)).
\end{equation}
 In the same context, the construction of
Corollary~\ref{fin-dim-base-koszul-duality-contramodule-side}
provides an equivalence between the derived category of left
$\CDiff(X,h)$\+modules and the contraderived category of left
CDG\+modules over $(\Omega(X),d,h)$,
\begin{equation} \label{twisted-differential-contramodule-side}
 \sD(\CDiff(X,h)\modl)\simeq\sD^\ctr((\Omega(X),d,h)\modl).
\end{equation} \qed
\end{thm}

 For any closed $2$\+form~$h$ on $X$, there is a natural equivalence
between the abelian categories of left modules over $\CDiff(X,h)$
and right modules over $\CDiff(X,-h)$, provided by the mutually
inverse functors of tensor product with the $O(X)$\+modules modules
of top differential forms and top polyvector fields on~$X$.
 This Morita equivalence can be obtained as a particular case of
Theorem~\ref{conversion-equivalence-theorem}
with Remark~\ref{commutative-cdg-ring-remark}.
 The results of Section~\ref{conversion-secn} provide a comparison
between~\eqref{twisted-differential-comodule-side}
and~\eqref{twisted-differential-contramodule-side}.

\begin{thm} \label{twisted-differential-quadrality}
 For any smooth affine algebraic variety $X$ over a field~$k$ and
any closed differential $2$\+form $h\in\Omega^2(X)$,
the constructions of Theorem~\ref{frobenius-quadrality-theorem}
(with Remark~\ref{commutative-cdg-ring-remark}) provide a commutative
square diagram of triangulated equivalences
\begin{equation}
\begin{tikzcd}
 \sD(\modr{\CDiff(X,-h)}) \arrow[r, bend right = 7] \arrow[dd] &
 \sD(\CDiff(X,h)\modl) \arrow[l] \arrow[dd]
 \\ \\
 \sD^\co((\Omega(X),d,h)\modl)
 \arrow[r, Leftrightarrow, no head, no tail]
 \arrow[uu, bend left = 20] &
 \sD^\ctr((\Omega(X),d,h)\modl) \arrow[uu, bend right = 20]
\end{tikzcd}
\end{equation}
where the equivalence of derived categories in the upper line is
induced by the conversion equivalence of abelian categories
$\modr{\CDiff(X,-h)}\simeq\CDiff(X,h)\modl$ (up to a shift by
$[-\dim_kX]$), while the co-contra correspondence in
the lower line is the result of Theorem~\ref{frobenius-co-contra-thm}.
 The vertical equivalences
are~\eqref{twisted-differential-comodule-side}
and~\eqref{twisted-differential-contramodule-side}.  \qed
\end{thm}

\subsection{Smooth differential operators}
\label{smooth-differential-subsecn}
 Let $X$ be a smooth compact real manifold.
 Denote by $O(X)$ the ring of smooth global functions
$X\rarrow\mathbb R$.

 We will consider smooth locally trivial vector bundles $E$ on~$X$.
 Then the $\mathbb R$\+vector space of smooth global sections $E(X)$
has a natural $O(X)$\+module structure.
 Moreover, the $O(X)$\+module $E(X)$ is finitely generated and
projective.
 The correspondence $E\longmapsto E(X)$ is an equivalence between
the category of (smooth locally trivial) vector bundles on $X$
and the category of finitely generated projective $O(X)$\+modules.
 In particular, any short exact sequence of vector bundles over
$X$ splits (as one can show using a partition of unity on~$X$).
 In this sense, smooth compact real manifolds are analogues of
\emph{affine} algebraic varieties.

 Given a vector bundle $E$ on $X$, the module of global sections of
the dual vector bundle $E^*$ can be obtained as the dual finitely
generated projective module, $E^*(X)=\Hom_{O(X)}(E(X),O(X))$.
 Given two vector bundles $E'$ and $E''$, the global sections of
the tensor product bundle $E'\ot E''$ are the tensor product of
the global sections, $(E'\ot\nobreak E'')(X)=E'(X)\ot_{O(X)}E''(X)$.
 Similarly, the global sections of the symmetric and exterior powers
of a vector bundle $E$ are computable as the symmetric and exterior
powers of the global sections, $(\Sym^nE)(X)=\Sym^n_{O(X)}(E(X))$ and
$(\Lambda^nE)(X)=\Lambda^n_{O(X)}(E(X))$.

 In particular, the tangent bundle $T$ and its dual (cotangent)
bundle $T^*$ are smooth locally trivial vector bundles on~$X$.
 In a local coordinate system $x_1$,~\dots, $x_m$ defined on an open
subset $U\subset X$, the sections of $T$ (called the \emph{vector
fields}) are represented by expressions like $\sum_{i=1}^m f_i\,
\d/\d x_i$ and the sections of $T^*$ (called the \emph{differential
$1$\+forms}) are represented by expressions like $\sum_{i=1}^m
f_i\,dx_i$, where $f_i\:U\rarrow\mathbb R$ are local functions.
 The natural map $d\:O(X)\rarrow T^*(X)$ (the differential) is
defined locally by the rule $d(f)=\sum_{i=1}^m \d f/\d x_i\, dx_i$,
and the action of vector fields in the functions (the derivative
$v(f)$ of a function~$f$ along a vector field~$v$) is given locally
by the formula $(\sum_{i=1}^m f_i\,\d/\d x_i)(f)=
\sum_{i=1}^m f_i\,\d f/\d x_i$.

 Let $E'$ and $E''$ be two vector bundles on $X$ of ranks
(the dimensions of the fibers) $r'$ and $r''\ge0$.
 Then a \emph{differential operator} $D\:E'(X)\rarrow E''(X)$ of
order~$\le n$ is an $\mathbb R$\+linear map such that, for any open
subset $U\subset X$ with a coordinate system $x_1$,~\dots, $x_m$
and any chosen trivializations of $E'$ and $E''$ over~$U$,
the operator $D$ can be expressed locally over $U$ as a linear
combination of compositions of at most~$n$ partial derivatives
$\d/\d x_i$ (acting in vector functions $U\rarrow\mathbb R^{r'}$
component-wise) with $r''\times r'$\+matrices of smooth functions
$U\rarrow\mathbb R$ as the coefficients.
 In particular, a differential operator $D\:O(X)\rarrow O(X)$ of
order~$\le n$ is represented locally over $U$ as a linear combination
of compositions of at most~$n$ partial derivatives $\d/\d x_i$, \
$1\le i\le m$, with smooth local functions $U\rarrow\mathbb R$ as
the coefficients.
 A differential operator $E'(X)\rarrow E''(X)$ of order~$0$ is
the same thing as a global section of the vector bundle
$\Hom(E',E'')$, or an $O(X)$\+linear map $E'(X)\rarrow E''(X)$.

 We denote the $\mathbb R$\+vector space of smooth differential
operators $E'(X)\rarrow E''(X)$ by $\Diff(X,E',E'')$, and the subspace
of differential operators of order~$\le n$ by $F_n\Diff(X,E',E'')
\subset\Diff(X,E',E'')$; so $\Diff(X,E',E'')=\bigcup_{n=0}^\infty
F_n\Diff(X,E',E'')$.
 When the two vector bundles $E'=E=E''$ are the same, we write
simply $\Diff(X,E)$.
 Then $\Diff(X,E)$ is a subalgebra in the $\mathbb R$\+algebra
$\End_{\mathbb R}(E(X))$ of all $\mathbb R$\+linear endomorphisms
of the vector space $E(X)$ (with respect to the composition).
 Furthermore, $\Diff(X,E)$ is a filtered ring with an increasing
filtration~$F$.
 When $E$ is the trivial line bundle on $X$ (so $E(X)=O(X)$),
we denote simply by $\Diff(X)$ the filtered ring of differential
operators $O(X)\rarrow O(X)$.

 Global sections of the exterior power $\Lambda^n(T^*)$ of
the cotangent bundle $T^*$ are called the \emph{differential
$n$\+forms} on~$X$.
 The vector space (in fact, $O(X)$\+module) of differential
$n$\+forms is denoted by $\Omega^n(X)=\Lambda^n_{O(X)}(T^*(X))$.
 There is a natural differential operator $d\:\Omega^n(X)\rarrow
\Omega^{n+1}(X)$ of order~$1$, defined for every $n\ge0$ and called
the \emph{de~Rham differential}.
 The operator~$d$ is uniquely defined by the property of being
and odd derivation of the graded algebra $\Omega(X)$ with $d^2=0$
whose restriction to the ring of functions, $d\:O(X)\rarrow
\Omega^1(X)=T^*(X)$, is the map mentioned above.

 The left action of $\Diff(X)$ by differential operators in $O(X)$
makes $\Diff(X)$ a left augmented ring over its subring
$O(X)=F_0\Diff(X)$, in the sense of Section~\ref{augmented-subsecn}.
 The associated graded ring $\gr^F\Diff(X)=\bigoplus_{n=0}^\infty
F_n\Diff(X)/F_{n-1}\Diff(X)$ is naturally isomorphic to
the graded ring $\Sym_{O(X)}(T(X))$, which is left and right
finitely projective Koszul (according to the discussion in
Section~\ref{symmetric-exterior-subsecn}).
 Similarly to Sections~\ref{algebraic-diffoperators-subsecn}\+-%
\ref{crystalline-subsecn}, the anti-equivalence of categories from
Corollary~\ref{left-augmented-koszul-duality-anti-equivalence}
assigns the right finitely projective Koszul DG\+ring $(\Omega(X),d)$
to the left augmented left finitely projective nonhomogeneous
Koszul ring $\Diff(X)$ (with the filtration~$F$).

 Hence the results of Sections~\ref{comodule-side-secn}
and~\ref{contramodule-side-secn} lead us to the following theorem.
 It is a rather rough algebraic version of derived relative
nonhomogeneous Koszul duality for the ring of smooth differential
operators $\Diff(X)$, in that the base ring of smooth functions
$O(X)$ is considered as an abstract ring \emph{without} any
additional structures (such as a Banach space metric or a topology).
 Notice that, similarly to the algebraic examples above, the graded
ring $\Omega(X)$ has only finitely many grading components (as
$\Omega^n(X)=0$ for $n>\dim_\boR X$).

\begin{thm} \label{smooth-differential-koszul-duality-theorem}
 For any smooth compact real manifold $X$, the construction of
Theorem~\ref{comodule-side-koszul-duality-theorem} provides
a triangulated equivalence between
the $\Diff(X)/O(X)$\+semicoderived category of right modules over
the ring of smooth differential operators $\Diff(X)$
and the coderived category of right DG\+modules over the de~Rham
DG\+algebra $(\Omega(X),d)$,
\begin{equation}
 \sD^\sico_{O(X)}(\modr{\Diff(X)})\simeq\sD^\co(\modr(\Omega(X),d)),
\end{equation}
while Theorem~\ref{reduced-koszul-duality-comodule-side-thm}
establishes a triangulated equivalence between the derived
category of right $\Diff(X)$\+modules and the reduced coderived
category of right DG\+modules over $(\Omega(X),d)$ relative to~$O(X)$,
\begin{equation}
 \sD(\modr{\Diff(X)})\simeq\sD^\co_{O(X)\red}(\modr(\Omega(X),d)).
\end{equation}
 In the same context, the construction of
Theorem~\ref{contramodule-side-koszul-duality-theorem} provides
a triangulated equivalence between
the $\Diff(X)/O(X)$\+semicontraderived category of left modules over
the ring $\Diff(X)$ and the contraderived category of left
DG\+modules over the DG\+algebra $(\Omega(X),d)$,
\begin{equation}
 \sD^\sictr_{O(X)}(\Diff(X)\modl)\simeq\sD^\ctr((\Omega(X),d)\modl),
\end{equation}
while Theorem~\ref{reduced-koszul-duality-contramodule-side-thm}
establishes a triangulated equivalence between the derived category
of left $\Diff(X)$\+modules and the reduced contraderived category
of left DG\+modules over $(\Omega(X),d)$ relative to~$O(X)$,
\begin{equation}
 \sD(\Diff(X)\modl)\simeq\sD^\ctr_{O(X)\red}((\Omega(X),d)\modl).
\end{equation} \qed
\end{thm}

 There is also a classical natural equivalence between
the abelian categories of left and right $\Diff(X)$\+modules,
$\Diff(X)\modl\simeq\modr{\Diff(X)}$, which was discussed in
Section~\ref{introd-conversion}.
 It can be obtained as a particular case of
Theorem~\ref{conversion-equivalence-theorem} with
Remark~\ref{commutative-cdg-ring-remark} (since in fact
$\Omega(X)$ is a relatively Frobenius Koszul graded ring).

\medskip

 More generally, the associated graded ring $\gr^F\Diff(X,E)=
\bigoplus_{n=0}^\infty F_n\Diff(X,E)/\allowbreak
F_{n-1}\Diff(X,E)$ is naturally isomorphic to the tensor product
$\End(E)(X)\ot_{O(X)}\Sym_{O(X)}(T(X))$ of the $O(X)$\+algebra
$\End(E)(X)=\Hom_{O(X)}(E(X),E(X))$ of sections of the vector bundle
$\End(E)$ and the graded $O(X)$\+algebra $\Sym_{O(X)}(T(X))$.
 This graded ring is left and right finitely projective Koszul over
its degree-zero component $\End(E)(X)$.
 So the filtered ring $(\Diff(X,E),F)$ is left finitely projective
nonhomogeneous Koszul.
 In order to describe the corresponding CDG\+ring, the notion of
a connection in a smooth vector bundle $E$ over $X$ is needed.

 A (\emph{smooth}) \emph{connection} $\nabla$ in a smooth vector bundle
$E$ is a map assigning to every smooth vector field $v\in T(X)$ and any
smooth section $e\in E(X)$ a smooth section $\nabla_v(e)\in E(X)$ such
that the equations~(i) and~(ii) from Section~\ref{vector-bundle-subsecn}
are satisfied.
 Using a partition of unity, one can show that a smooth connection
exists in any smooth vector bundle over a smooth manifold.

 Similarly to the algebraic case discussed in
Section~\ref{vector-bundle-subsecn}, a connection $\nabla$ in $E$
can be interpreted as an $\boR$\+linear map
$E(X)\rarrow\Omega^1(X)\ot_{O(X)}E(X)$ defined
by the rule $\lan v,\nabla(e)\ran = \nabla_v(e)$.
 In fact, it is a differential operator $E(X)\rarrow (T^*\ot E)(X)$ of
order~$1$.
 The trivial line bundle on $X$ has a canonical \emph{trivial
connection} defined by the rule $\nabla_v(f)=v(f)$ for all $v\in T(X)$
and $f\in O(X)$.
 If vector bundles $E'$ and $E''$ over $X$ are endowed with
connections $\nabla'$ and~$\nabla''$, then the vector bundles
$E'\oplus E''$, \ $E'\ot E''$, and $\Hom(E',E'')$ acquire
the induced connections.

 Similarly to Section~\ref{vector-bundle-subsecn}, a connection $\nabla$
in a smooth vector bundle $E$ induces an odd derivation~$d_\nabla$ of
degree~$1$ on the graded left $\Omega(X)$\+module
$\Omega(X)\ot_{O(X)}E(X)$ compatible with the de~Rham differential~$d$
on the graded ring $\Omega(X)$, in the sense of
Section~\ref{cdg-modules-subsecn}.
 The map $d_\nabla\:\Omega^n(X)\ot_{O(X)}E(X)\rarrow
\Omega^{n+1}(X)\ot_{O(X)}E(X)$, \ $n\ge1$ (in fact, a differential
operator of order~$1$ acting between the global sections of the vector
bundles $\Lambda^n(T^*)\ot E$ and $\Lambda^{n+1}(T^*)\ot E$ on~$X$) is
given by the formula~\eqref{connection-de-Rham}.

 The \emph{curvature $2$\+form} $h_\nabla$ of a smooth connection
$\nabla$ in $E$ is an element of the vector space (or $O(X)$\+module)
$\Omega^2(X)\ot_{O(X)}\End(E)(X)$ of \emph{differential $2$\+forms
on $X$ with the coefficients in} $\End(E)$, that is, global sections of
the vector bundle $\Lambda^2(T^*)\ot\End(E)$ on~$X$.
 The element $h_\nabla\in\Omega^2(X)\ot_{O(X)}\End(E)$ is defined by
the identity
$$
 \nabla_v(\nabla_w(e))-\nabla_w(\nabla_v(e))=
 \nabla_{[v,w]}(e)+\lan v\wedge w,\,h\ran (e)
$$
holding for all $v$, $w\in T(X)$ and $e\in E(X)$.
 Here $[v,w]\in T(X)$ is the commutator of the vector fields
$v$ and~$w$ on $X$ (defined as in Section~\ref{crystalline-subsecn}),
while $\lan\ , \ \ran$ denotes the $O(X)$\+linear map
$\Lambda^2_{O(X)}(T(X))\ot_{O(X)}\Omega^2(X)\ot_{O(X)}\End(E)(X)
\rarrow \End(E)(X)$ induced by the pairing
$\Lambda^2_{O(X)}(T(X))\ot_{O(X)}\Lambda^2_{O(X)}(T^*(X))\rarrow O(X)$
from Section~\ref{twisted-differential-subsecn}.

 The choice of a smooth connection $\nabla=\nabla_E$ in $E$ defines
a left $\End(E)(X)$\+submod\-ule of strict generators in
$F_1\Diff(X,E)$ (in the sense of Section~\ref{self-consistency-subsecn})
consisting of all the $\boR$\+linear combinations
of differential operators of the form $g\nabla_v\:E(X)\rarrow E(X)$,
where $g\in\End(E)(X)$ and $v\in T(X)$.
 Let $\nabla_{\End(E)}$ denote the connection in the vector bundle
$\End(E)$ induced by the connection $\nabla$ in~$E$.

 Similarly to the algebraic case considered in
Section~\ref{vector-bundle-subsecn}, the anti-equivalence of categories
from Corollary~\ref{nonhomogeneous-koszul-duality-anti-equivalence}
assigns the right finitely projective Koszul CDG\+ring
\begin{equation} \label{smooth-connection-de-Rham-cdg-ring}
 (\Omega(X)\ot_{O(X)}\End(E)(X),\>d_{\nabla_{\End(E)}},\>h_{\nabla_E})
\end{equation}
as in the formula~\eqref{connection-de-Rham-cdg-ring} to the left
finitely projective nonhomogeneous Koszul ring $\Diff(X,E)$ with
the filtration~$F$.
 Hence we obtain the following Koszul duality theorem.
 
\begin{thm} \label{smooth-vector-bundle-koszul-duality-theorem}
 For any smooth compact real manifold $X$ and vector bundle $E$
on $X$, the construction of
Theorem~\ref{comodule-side-koszul-duality-theorem} provides
a triangulated equivalence between
the $\Diff(X,E)/\End(E)(X)$\+semicoderived category of right modules
over the ring of smooth differential operators $\Diff(X,E)$ acting
in the sections of $E$ and the coderived category of right CDG\+modules
over the CDG\+algebra~\eqref{smooth-connection-de-Rham-cdg-ring}
of differential forms on $X$ with the coefficients in $\End(E)$,
\begin{multline}
 \sD^\sico_{\End(E)(X)}(\modr{\Diff(X,E)}) \\ \,\simeq\,
 \sD^\co(\modr(\Omega(X)\ot_{O(X)}\End(E)(X),
 \,d_{\nabla_{\End(E)}},\,h_{\nabla_E})),
\end{multline}
while Theorem~\ref{reduced-koszul-duality-comodule-side-thm}
establishes a triangulated equivalence between the derived
category of right $\Diff(X,E)$\+modules and the reduced coderived
category of right CDG\+modules over
the CDG\+ring~\eqref{smooth-connection-de-Rham-cdg-ring}
relative to the ring $\End(E)(X)$,
\begin{multline}
 \sD(\modr{\Diff(X,E)}) \\ \,\simeq\,\sD^\co_{\End(E)(X)\red}
 (\modr(\Omega(X)\ot_{O(X)}\End(E)(X),
 \,d_{\nabla_{\End(E)}},\,h_{\nabla_E})).
\end{multline}
 In the same context, the construction of
Theorem~\ref{contramodule-side-koszul-duality-theorem} provides
a triangulated equivalence between
the $\Diff(X,E)/\End(E)(X)$\+semicontraderived category of left
modules over the ring $\Diff(X,E)$ and the contraderived category
of left CDG\+modules over
the CDG\+algebra~\eqref{smooth-connection-de-Rham-cdg-ring},
\begin{multline} \label{smooth-vector-bundle-unreduced-contra-side}
 \sD^\sictr_{\End(E)(X)}(\Diff(X,E)\modl) \\ \,\simeq\,
 \sD^\ctr((\Omega(X)\ot_{O(X)}\End(E)(X),
 \,d_{\nabla_{\End(E)}},\,h_{\nabla_E})\modl),
\end{multline}
while Theorem~\ref{reduced-koszul-duality-contramodule-side-thm}
establishes a triangulated equivalence between the derived category
of left $\Diff(X,E)$\+modules and the reduced contraderived category
of left CDG\+modules over
the CDG\+ring~\eqref{smooth-connection-de-Rham-cdg-ring} relative to
the ring $\End(E)(X)$,
\begin{multline} \label{smooth-vector-bundle-reduced-contra-side}
 \sD(\Diff(X,E)\modl) \\ \,\simeq\,
 \sD^\ctr_{\End(E)(X)\red}((\Omega(X)\ot_{O(X)}\End(E)(X),
 \,d_{\nabla_{\End(E)}},\,h_{\nabla_E})\modl).
\end{multline} \qed
\end{thm}

 For example, the equivalences of
categories~(\ref{smooth-vector-bundle-unreduced-contra-side}\+-%
\ref{smooth-vector-bundle-reduced-contra-side})
assign the left CDG\+module $(\Omega(X)\ot_{O(X)}E(X),\,d_{\nabla_E})$
(as in~\eqref{connection-de-Rham-cdg-module}) over
the CDG\+ring~\eqref{smooth-connection-de-Rham-cdg-ring} to the left
$\Diff(X,E)$\+mod\-ule~$E(X)$.

 Furthermore, similarly to Section~\ref{twisted-differential-subsecn},
for any closed $2$\+form $h\in\Omega^2(X)$ (that is a smooth
differential $2$\+form such that $d(h)=0$), the triple
$(\Omega(X),d,h)$ is a right finitely projective Koszul CDG\+ring.
 The corresponding left finitely projective nonhomogeneous Koszul
filtered ring $\Diff(X,h)$, defined by
the relations~\eqref{crystalline-diff-relations}
and~\eqref{twisted-commutator-vector-fields}, is called
the \emph{ring of twisted smooth differential operators on}~$X$
(twisted by~$h$).
 It is claimed, on the basis of Theorem~\ref{pbw-theorem-thm}, that
the associated graded ring $\gr^F\CDiff(X,h)=
\bigoplus_{n=0}^\infty F_n\CDiff(X,h)/F_{n-1}\CDiff(X,h)$ is
naturally isomorphic to the symmetric algebra $\Sym_{O(X)}(T(X))$.

\begin{thm} \label{smooth-twisted-diffoperators-koszul-duality-theorem}
 For any smooth compact real manifold $X$ and any closed differential
$2$\+form $h\in\Omega^2(X)$, \,$d(h)=0$, the construction of
Theorem~\ref{comodule-side-koszul-duality-theorem} provides
a triangulated equivalence between
the $\Diff(X,h)/O(X)$\+semicoderived category of right modules over
the ring of twisted smooth differential operators $\Diff(X,h)$
and the coderived category of right CDG\+modules over the de~Rham
CDG\+algebra $(\Omega(X),d,h)$,
\begin{equation}
 \sD^\sico_{O(X)}(\modr{\Diff(X,h)})\simeq\sD^\co(\modr(\Omega(X),d,h)),
\end{equation}
while Theorem~\ref{reduced-koszul-duality-comodule-side-thm}
establishes a triangulated equivalence between the derived
category of right $\Diff(X,h)$\+modules and the reduced coderived
category of right CDG\+modules over $(\Omega(X),d,h)$ relative
to~$O(X)$,
\begin{equation}
 \sD(\modr{\Diff(X,h)})\simeq\sD^\co_{O(X)\red}(\modr(\Omega(X),d,h)).
\end{equation}
 In the same context, the construction of
Theorem~\ref{contramodule-side-koszul-duality-theorem} provides
a triangulated equivalence between
the $\Diff(X,h)/O(X)$\+semicontraderived category of left modules over
the ring $\Diff(X,h)$ and the contraderived category of left
CDG\+modules over the CDG\+algebra $(\Omega(X),d,h)$,
\begin{equation}
 \sD^\sictr_{O(X)}(\Diff(X,h)\modl)\simeq\sD^\ctr((\Omega(X),d,h)\modl),
\end{equation}
while Theorem~\ref{reduced-koszul-duality-contramodule-side-thm}
establishes a triangulated equivalence between the derived category
of left $\Diff(X,h)$\+modules and the reduced contraderived category
of left CDG\+modules over $(\Omega(X),d,h)$ relative to~$O(X)$,
\begin{equation}
 \sD(\Diff(X,h)\modl)\simeq\sD^\ctr_{O(X)\red}((\Omega(X),d,h)\modl).
\end{equation} \qed
\end{thm}

 Similarly to Section~\ref{twisted-differential-subsecn},
there is a natural equivalence between the abelian categories of
left modules over $\Diff(X,h)$ and right modules over $\Diff(X,-h)$,
provided by the mutually inverse functors of tensor product with
the $O(X)$\+modules $\Omega^m(X)$ and $\Lambda^m_{O(X)}(T(X))$,
where $m=\dim_\boR X$.
 This equivalence can be obtained as a particular case of
the construction of Theorem~\ref{conversion-equivalence-theorem}
with Remark~\ref{commutative-cdg-ring-remark}.

\subsection{Dolbeault differential operators}
\label{dolbeault-subsecn}
 For a reference on the basics of complex analytic geometry,
see, e.~g., \cite[Section~I.2]{Voi}.

 In this section, we denote by~$i$ a chosen imaginary unit
(a square root of~$-1$) in the field of complex numbers~$\boC$.
 The complex conjugation map is denoted by $z=x+iy\longmapsto
\bar z=x-iy$ (where $x$, $y\in\boR$).

 For any vector space $V$ over $\boC$, we denote by $V_\boR$
the underlying real vector space of~$V$.
 So $\dim_\boR V_\boR = 2\dim_\boC V$.
 Similarly, for any vector space $U$ over $\boR$, we denote by
$U_\boC$ the $\boC$\+vector space $U_\boC=\boC\ot_\boR U$.
 So $\dim_\boC U_\boC = \dim_\boR U$.
 We denote the elements of $U_\boC$ by $u+iv$, where $u$, $v\in U$
and $i$~is the imaginary unit.

 We also denote by $\overline V$ the vector space $V$ with the conjugate
complex structure.
 Denoting by $\bar v\in \overline V$ the element corresponding to
an element $v\in V$, the $\boC$\+vector space structure on $\overline V$
is defined by the rule $\bar z\bar v=\overline{zv}$ for all
$z\in\boC$ and $v\in V$, where $z\longmapsto \bar z\:\boC\rarrow\boC$
is the complex conjugation.
 The vector space $U_\boC$ is endowed with a natural isomorphism
of $\boC$\+vector spaces $\overline{U_\boC}\simeq U_\boC$;
the isomorphism is provided by the complex conjugation map
$w\longmapsto \overline w\:U_\boC\rarrow U_\boC$ defined by
the obvious rule $\overline{u+iv}=u-iw$ for all $u$, $v\in U$.

 For any $\boC$\+vector space $V$, the $\boC$\+vector space
$(V_\boR)_\boC$ decomposes canonically into a direct sum of two
$\boC$\+vector spaces, one of them isomorphic naturally to $V$
and the other one to~$\overline V$.
 Specifically, denote by $J\:V\rarrow V$ the action of the imaginary
unit $i\in\boC$ in~$V$.
 Let $V^{(1,0)}\subset (V_\boR)_\boC$ denote the set of all elements
of the form $v-iJv$, where $v\in V$, and let $V^{(0,1)}\subset
(V_\boR)_\boC$ denote the set of all elements of the form
$v+iJv$, where $v\in V$.
 Then both $V^{(1,0)}$ and $V^{(0,1)}$ are $\boC$\+vector subspaces
in $(V_\boR)_\boC$, and $(V_\boR)_\boC=V^{(1,0)}\oplus V^{(0,1)}$.
 Furthermore, there are natural isomorphisms of $\boC$\+vector spaces
$v\longmapsto v-iJv\:V\rarrow V^{(1,0)}$ and $\bar v\longmapsto
v+iJv\:\overline V\rarrow V^{(0,1)}$.
 The complex conjugation map $w\longmapsto\overline w\:(V_\boR)_\boC
\rarrow(V_\boR)_\boC$ switches the two sides of this direct
sum decomposition, taking $V^{(1,0)}$ into $V^{(0,1)}$ and
$V^{(0,1)}$ into $V^{(1,0)}$.

 For any $\boC$\+vector space $V$, one can consider the dual
$\boC$\+vector space $\Hom_\boC(V,\boC)$.
 One can also consider the dual $\boR$\+vector space
$\Hom_\boR(V,\boR)$.
 Let us endow the $\boR$\+vector space $\Hom_\boR(V,\boR)$ with
a $\boC$\+vector space structure by the rule $\lan v, Jf\ran
=\lan Jv,f\ran$ for all $v\in V$ and $f\in\Hom_\boR(V,\boR)$
(where $\lan\ , \ \ran\:V\ot_\boR\Hom_\boR(V,\boR)\rarrow\boR$
is the natural pairing).
 Then there is a natural isomorphism of $\boC$\+vector spaces
$\Hom_\boR(V,\boR)\simeq\Hom_\boC(V,\boC)$ given by the rule
$\Hom_\boR(V,\boR)\ni f\longmapsto \hat f\in\Hom_\boC(V,\boC)$ with
$\lan\lan v,\hat f\ran\ran=\lan v,f\ran-i\lan Jv,f\ran\in\boC$ for
all $v\in V$ (where $\lan\lan\ , \ \ran\ran\:V\ot_\boC\Hom_\boC(V,\boC)
\rarrow\boC$ is the notation for the $\boC$\+valued pairing).

 Let $X$ be a complex manifold (which we will eventually assume
to be compact).
 We denote the underlying smooth real manifold of $X$ by~$X_\boR$,
and refer to the previous Section~\ref{smooth-differential-subsecn}
for the discussion of differential operators and differential forms
on~$X_\boR$.
 For any vector bundle $V$ on $X_\boR$ and a point $s\in X$, we
denote by $V_s$ the fiber of $V$ at~$s$.
 So $V_s$ is a finite-dimensional (real or complex) vector space.

 Let $T$ and $T^*$ denote the tangent and cotangent bundles
to $X_\boR$; so $T_s$ is the tangent space to $X_\boR$ at~$s$
and $T^*_s=\Hom_\boR(T_s,\boR)$ is the cotangent space
for every $s\in X$.
 For every point $s\in X$, the $\boR$\+vector space $T_s$ has
a natural $\boC$\+vector space structure, which is constructed
as follows.

 Let $z_1$,~\dots, $z_m$ be a holomorphic local coordinate system
in an open subset $U\subset X$, \ $s\in U$.
 Denote by $x_k$ and~$y_k$ the real and imaginary parts of
the complex variable~$z_k$, \ $1\le k\le m$; so $z_k=x_k+iy_k$
and $x_k$, $y_k$ are $\boR$\+valued functions on~$U$.
 Let $w_1$,~\dots, $w_m$ be another holomorphic coordinate system
in $U$; put $w_l=u_l+iv_l$, where $u_l$, $v_l$~are $\boR$\+valued
local functions for all $1\le l\le m$.

 Then, for every fixed point $s\in U\subset X$,
the $2m\times 2m$\+matrix with real entries
$$
 \begin{pmatrix}
  \frac{\d u_l}{\d x_k}(s) & \frac{\d u_l}{\d y_k}(s) \\[5pt]
  \frac{\d v_l}{\d x_k}(s) & \frac{\d v_l}{\d y_k}(s)
 \end{pmatrix}
$$
can be obtained by applying the functor $V\longmapsto V_\boR$
to the $m\times m$\+matrix with complex entries
$\bigl(\frac{\d w_l}{\d z_k}(s)\bigr)$.
 Consequently, the $\boC$\+vector space structure on the vector
space $T_s$ defined by the rules $J\bigl(\frac{\d}{\d x_k}\bigr)=
\frac{\d}{\d y_k}$ and $J\bigl(\frac{\d}{\d y_k})=-\frac{\d}{\d x_k}$
does not depend on the choice of a holomorphic coordinate
system~$(z_k)_{k=1}^m$ in a neighborhood of~$s$ in~$X$.

 Following the above rule for the passage to the dual vector space,
there is the induced $\boC$\+vector space structure on
the $\boR$\+vector space $T^*_s$ for every $s\in X$,
defined in holomorphic local coordinates~$z_k$ by the formulas
$J(dx_k)=-dy_k$ and $J(dy_k)=dx_k$.
 Notice that $T_s^*$ is the vector space of all differentials $d_sf$
at the point $s\in X$ of \emph{real-valued} smooth functions
$f\:X\rarrow\boR$; still it has a natural $\boC$\+vector space
structure.
 If $g\:X\rarrow\boC$ is a complex-valued smooth function, then
the differential $d_sg$ belongs to the \emph{complexified} cotangent
space $T^*_{s,\boC}=\boC\ot_\boR T^*_s$.
 
 According to the above discussion, there is a natural direct sum
decomposition $T^*_{s,\boC}=T^*_s{}^{(1,0)}\oplus T^*_s{}^{(0,1)}$.
 Here the \emph{holomorphic cotangent space} $T^*_s{}^{(1,0)}$ to $X$
at~$s$ is a $\boC$\+vector space with a basis $dz_k=dx_k+i\,dy_k$,
\ $1\le k\le m$.
 It consists of the differentials $d_sh$ of holomorphic local functions
$h\:U\rarrow\boC$ at the point $s\in U$.
 The \emph{anti-holomorphic cotangent space} $T^*_s{}^{(0,1)}$ is
a $\boC$\+vector space with a basis $d\bar z_k=dx_k-i\,dy_k$.
 It consists of the differentials $d_s\bar h$ of anti-holomorphic local
functions $\bar h\:U\rarrow\boC$ (where $h$~is holomorphic).

 The complexified cotangent space $T^*_{s,\boC}$ is the dual
$\boC$\+vector space to the complexified tangent space $T_{s,\boC}=
\boC\ot_\boR T_s$.
 The natural pairing $\lan\ , \ \ran\:T_{s,\boC}
\ot_\boC T^*_{s,\boC}\rarrow\boC$ is constructed as the unique
$\boC$\+linear extension of the natural $\boR$\+linear pairing
$\lan\ , \ \ran\: T_s\ot_\boR T^*_s\rarrow\boR$.
 The dual basis to $(dz_k,\,d\bar z_k\in T^*_{s,\boC})_{k=1}^m$ with
respect to this pairing is denoted by $(\d/\d z_k,\,\d/\d\bar z_k\in
T_{s,\boC})_{k=1}^m$; so $\lan \d/\d z_k,\,dz_l\ran = \delta_{k,l} =
\lan \d/\d\bar z_k,\,d\bar z_l\ran$ and $\lan \d/\d z_k,\,d\bar z_l\ran
= 0 = \lan \d/\d\bar z_k,\,d z_l\ran$ for $1\le k$, $l\le m$.
 Then one has $dg=\sum_{k=1}^m\frac{\d g}{\d z_k}dz_k+
\sum_{k=1}^m\frac{\d g}{\d\bar z_k}d\bar z_k$
for any smooth function $g\:U\rarrow\boC$.
 Explicitly, $\d/\d z_k=\frac{1}{2}(\d/\d x_k-i\,\d/\d y_k)$ and
$\d/\d\bar z_k=\frac{1}{2}(\d/\d x_k+i\,\d/\d y_k)$.

 According to the same discussion above, there is also a natural direct
sum decomposition $T_{s,\boC}=T_s^{(1,0)}\oplus T_s^{(0,1)}$.
 Here the \emph{holomorphic tangent space} $T_s^{(1,0)}$ to $X$
at~$s$ is a $\boC$\+vector space with a basis $(\d/\d z_k)_{k=1}^m$
and the \emph{anti-holomorphic tangent space} $T_s^{(0,1)}$ is
a $\boC$\+vector space with a basis $(\d/\d\bar z_k)_{k=1}^m$.
 Smooth sections of the vector bundle $T^{(1,0)}$ on $X$ are called
\emph{$(1,0)$\+vector fields}, and smooth sections of the vector
bundle $T^{(0,1)}$ on $X$ are called \emph{$(0,1)$\+vector fields}.
 Holomorphic local functions $h\:U\rarrow\boC$ are distinguished
by the condition that $v(h)=0$ in $U$ for any $(0,1)$\+vector field
$v\in T^{(0,1)}(X)$.
 Similarly, anti-holomorphic local functions $\bar h\:U\rarrow\boC$
are characterized by the condition that $u(\bar h)=0$ in $U$ for
any $(1,0)$\+vector field $u\in T^{(1,0)}(X)$.

 For every $n\ge0$, the complexified space of exterior forms
$\Lambda^n_\boR(T^*_s)_\boC$ at a point $s$ in $X_\boR$ decomposes
naturally as
\begin{multline*}
 \boC\ot_\boR\Lambda^n_\boR(T^*_s)\simeq
 \Lambda^n_\boC(T^*_{s,\boC})\simeq
 \Lambda^n_\boC(T^*_s{}^{(1,0)}\oplus T^*_s{}^{(0,1)}) \\ \simeq
 \bigoplus\nolimits^{p,q\ge0}_{p+q=n}\Lambda^p_\boC(T^*_s{}^{(1,0)})
 \ot_\boC\Lambda^q_\boC(T^*_s{}^{(0,1)}).
\end{multline*}
 This is a direct sum decomposition of the smooth complex vector
bundle $\Lambda^n(T^*)_\boC$ on~$X_\boR$.
 The space of global sections decomposes accordingly,
\begin{equation} \label{p-q-forms-decomposition}
 \Omega^n(X_\boR)_\boC\simeq\Lambda^n(T^*)_\boC(X_\boR)\simeq
 \bigoplus\nolimits^{p,q\ge0}_{p+q=n}\Omega^{(p,q)}(X).
\end{equation}

 This is a decomposition of the $\boC$\+vector space
$\Omega^n(X_\boR)_\boC$ of smooth $\boC$\+valued differential
$n$\+forms on $X_\boR$ into the direct sum of the spaces of
\emph{differential $(p,q)$\+forms}.
 In holomorphic local coordinates~$z_k$, a differential $(p,q)$\+form
is an expression like
$$
 \sum\nolimits_{1\le k_1<\dotsb<k_p\le m}
 \sum\nolimits_{1\le l_1<\dotsb<l_q\le m}
 f_{k_1,\dotsc,k_p;l_1,\dotsc,l_q}
 dz_{k_1}\wedge\dotsb \wedge dz_{k_p}\wedge
 d\bar z_{l_1}\wedge\dotsb \wedge d\bar z_{l_q},
$$
with $p$~holomorphic differentials $d z_{k_1}$,~\dots, $dz_{k_p}$
and $q$~anti-holomorphic differentials $d \bar z_{l_1}$,~\dots,
$d\bar z_{l_q}$ in the exterior product.
 Here the coefficients $f_{k_1,\dotsc,k_p;l_1,\dotsc,l_q}$
are smooth complex-valued local functions $U\rarrow\boC$.
 In fact, this is even a direct sum decomposition of
$\Omega^n(X_\boR)_\boC$ as a module over the ring $O(X_\boR)_\boC$
of smooth complex-valued global functions on~$X_\boR$.
 This direct sum decomposition makes $\Omega(X_\boR)_\boC$
a \emph{bigraded algebra} over the ring $O(X_\boR)_\boC$.

 The de~Rham differential~$d$ on the graded algebra of differential
forms $\Omega(X_\boR)$ was discussed in
Section~\ref{smooth-differential-subsecn}.
 Taking the tensor product with $\boC$ over $\boR$, we obtain
an odd derivation of degree~$1$ on the graded algebra of
$\boC$\+valued differential forms $\Omega(X_\boR)_\boC$; we denote
this differential also by~$d$.
 With respect to the direct sum
decomposition~\eqref{p-q-forms-decomposition}, the differential
$d\:\Omega^n(X_\boR)_\boC\rarrow\Omega^{n+1}(X_\boR)_\boC$ has
two components:
$$
 \d\:\Omega^{p,q}(X)\lrarrow\Omega^{p+1,q}(X)
 \quad\text{and}\quad
 \bar\d\:\Omega^{p,q}(X)\lrarrow\Omega^{p,q+1}(X),
$$
making $\Omega^{\bu,\bu}(X)$ a \emph{bicomplex}.
 In particular, the differential $d\:O(X_\boR)_\boC\rarrow
\Omega^1(X_\boR)_\boC\simeq\Omega^{(1,0)}(X)\oplus\Omega^{(0,1)}(X)$
decomposes as $d=\d+\bar\d$, where in holomorphic local
coordinates~$z_k$ one has $\d(f)=\sum_{k=1}^m\frac{\d f}{\d z_k}d z_k$
and $\bar\d(f)=\sum_{k=1}^m\frac{\d f}{\d \bar z_k} d\bar z_k$
for any smooth complex-valued function~$f$.

 The $p=0$ part of the bicomplex $\Omega^{p,q}(X)$
\begin{equation} \label{dolbeault-complex}
 0\lrarrow\Omega^{0,0}(X)\overset{\bar\d}\lrarrow\Omega^{0,1}(X)
 \overset{\bar\d}\lrarrow
 \dotsb\overset{\bar\d}\lrarrow\Omega^{0,m}(X)\lrarrow0
\end{equation}
is called the \emph{Dolbeault complex} of a complex manifold~$X$.
 Here $m=\dim_\boC X$ is the dimension of $X$ as a complex manifold.
 The Dolbeault complex $\Omega^{0,\bu}(X)$ is a DG\+algebra over
the field~$\boC$.

 The degree-zero component $\Omega^{0,0}(X)$ is the $\boC$\+algebra
$O(X_\boR)_\boC$ of smooth $\boC$\+valued functions on~$X_\boR$.
 The kernel $H^0_{\bar\d}(\Omega^{0,\bu}(X))=O(X)$ of the differential
$\bar\d\:\Omega^{0,0}(X)\rarrow\Omega^{0,1}(X)$ in the Dolbeault
complex is the subalgebra $O(X)\subset O(X_\boR)_\boC$ of
\emph{holomorphic} global functions $X\rarrow\boC$.
 For comparison, the kernel $H^0_d(\Omega^\bu(X_\boR))$ of
the differential $d\:\Omega^0(X_\boR)\rarrow\Omega^1(X_\boR)$ in
the de~Rham complex is the subalgebra of \emph{locally constant}
global functions in the $\boR$\+algebra $O(X_\boR)$ of smooth
functions $X_\boR\rarrow\boR$.

 Consider the ring (or $\boR$\+algebra) of differential operators
$\Diff(X_\boR)$ on the smooth real manifold $X_\boR$, as discussed in
Section~\ref{smooth-differential-subsecn}.
 The tensor product $\boC\ot_\boR\Diff(X_\boR)=\Diff(X_\boR)_\boC$
is a $\boC$\+algebra acting naturally in the $\boC$\+vector space
of smooth complex-valued global functions $O(X_\boR)_\boC$.
 We are interested in the following subring (or $\boC$\+subalgebra)
$\Diff^{\bar\d}(X)\subset\Diff(X_\boR)_\boC$, which we call the ring
of \emph{$\bar\d$\+differential operators} (or \emph{Dolbeault
differential operators}) on a complex manifold~$X$.

 By the definition, $\Diff^{\bar\d}(X)$ is the subring generated by
the subring of smooth $\boC$\+valued functions $O(X_\boR)_\boC\subset
\Diff(X_\boR)_\boC$ and the left $O(X_\boR)_\boC$\+submodule of
smooth $(0,1)$\+vector fields $T^{(0,1)}(X)\subset T(X)_\boC\subset
\Diff(X_\boR)_\boC$ in the ring of smooth $\boC$\+valued
differential operators $\Diff(X_\boR)_\boC$ on~$X_\boR$.
 The increasing filtration $F$ by the order of differential operators
on $\Diff^{\bar\d}(X)$ is induced by the filtration $F$ on
$\Diff(X_\boR)_\boC$ (which, in turn, is induced by the filtration
$F$ on $\Diff(X_\boR)$).

 The main property of the ring $\Diff^{\bar\d}(X)$ is that its action
in the ring of smooth complex-valued functions $O(X_\boR)_\boC$
commutes with the operators of multiplication by holomorphic functions.
 This holds true because holomorphic functions have zero derivatives
along $(0,1)$\+vector fields.

 More generally, a smooth complex vector bundle $E$ over $X_\boR$ is
said to be \emph{holomorphic} (or ``have a holomorphic structure'') if
the notion of a holomorphic local section of $E$ is defined.
 For example, the $(1,0)$\+cotangent bundle $T^*{}^{(1,0)}$ on $X$
has a natural holomorphic structure in which $\boC$\+linear
combinations of the expressions $fd(g)$ with holomorphic local
functions $f$, $g\:U\rarrow\boC$ are the holomorphic local sections
(these are called the \emph{holomorphic differential $1$\+forms}
on~$X$).
 Similarly, the $(1,0)$\+tangent bundle $T^{(1,0)}$ on $X$ has
a natural holomorphic structure in which the expressions
$\sum_{k=1}^m f_k\frac{\d}{\d z_k}$ with holomorphic local functions
$f_k\:U\rarrow\boC$ are the holomorphic local sections (one has to
check that the class of \emph{holomorphic vector fields}, defined
in this way, does not depend on the choice of a local coordinate
system~$z_k$).
 We denote the $\boC$\+vector space of smooth sections of
a holomorphic vector bundle $E$ on $X$ by $E(X_\boR)$ and the subspace
of holomorphic sections by $E(X)\subset E(X_\boR)$.

 For any holomorphic vector bundle $E$ over $X$, the ring of
$\bar\d$\+differential operators $\Diff^{\bar\d}(X)$ acts naturally
in the $\boC$\+vector space of \emph{smooth} sections $E(X_\boR)$;
so $E(X_\boR)$ is a left $\Diff^{\bar\d}(X)$\+module.
 This action is constructed as follows.
 The ring of smooth $\boC$\+valued functions $O(X_\boR)_\boC$ acts
in $E(X_\boR)$ by the usual multiplications of smooth sections of
a vector bundle with smooth functions.
 To define the action of $(0,1)$\+vector fields in $E(X_\boR)$,
choose a local basis of holomorphic sections $e_l\in E(U)$,
\ $1\le l\le n$.
 This means that~$e_l$ are holomorphic sections of $E$ over $U$ and
for every point $s\in U$ the vectors $e_l(s)$, \ $1\le l\le n$,
form a basis of the $\boC$\+vector space~$E_s$.

 Let $e=\sum_{l=1}^n f_le_l\in E(U_\boR)$ be a smooth section of $E$
over $U$; so $f_l\in O(U_\boR)_\boC$ are smooth complex-valued
local functions.
 Furthermore, let $v=\sum_{k=1}^m g_k\,\d/\d\bar z_k$ be
a smooth $(0,1)$\+vector field in $U$; so $g_k\:U\rarrow\boC$ are
smooth complex-valued local functions as well.
 Then we put
$$
 v(e)=\sum\nolimits_{l=1}^n v(f_l)e_l=\sum\nolimits_{l=1}^n
 \sum\nolimits_{k=1}^m g_k\frac{\d f_l}{\d\bar z_k}e_l\,\in\,E(U_\boR).
$$
 Using the fact that $v(h)=0$ for any holomorphic local function
$h\:U\rarrow\boC$, one can check that the above definition of
$v(e)$ does not depend on the choice of a local basis of holomorphic
sections~$e_l$ of $E$ over~$U$; so the local constructions glue
together to a well-defined differential operator of order~$1$
providing the action of~$v$ in $E(X_\boR)$.

 Similarly one can show that, for any morphism $E'\rarrow E''$ of
holomorphic vector bundles over $X$ (that is, a morphism of smooth
compex vector bundes taking holomorphic local sections to holomorphic
local sections), the induced map of the spaces of smooth sections
$E'(X_\boR)\rarrow E''(X_\boR)$ is a morphism of left
$\Diff^{\bar\d}(X)$\+modules.
 This generalizes the above assertion about the commutativity of
the action of $\Diff^{\bar\d}(X)$ in $E(X_\boR)$ with
the multiplications by holomorphic functions.

 Now let us assume that $X$ is a \emph{compact} complex manifold.
 Then, according to the discussion in
Section~\ref{smooth-differential-subsecn}, the category of smooth
real vector bundles over $X_\boR$ is equivalent to the category of
finitely generated projective modules over $O(X_\boR)$.
 Similarly, the category of smooth complex vector bundles over $X_\boR$
is equivalent to the category of finitely generated projective
modules over the ring $O(X_\boR)_\boC$.

 The associated graded ring $\gr^F\Diff^{\bar\d}(X)=
\bigoplus_{n=0}^\infty F_n\Diff^{\bar\d}(X)/F_{n-1}\Diff^{\bar\d}(X)$
is isomorphic to the symmetric algebra
$A=\Sym_{O(X_\boR)_\boC}(T^{(0,1)}(X_\boR))$ of the finitely generated
projective module $T^{(0,1)}(X_\boR)$ of smooth $(0,1)$\+vector fields
on $X$ over the ring $O(X_\boR)_\boC$ of smooth complex-valued
functions.
 According to Section~\ref{symmetric-exterior-subsecn}, this graded ring
is left and right finitely projective Koszul.
 Hence the filtered ring $\Diff^{\bar\d}(X)$ is a left finitely
projective nonhomogeneous Koszul ring.
 Furthermore, the natural left action of $\Diff^{\bar\d}(X)$ in
$O(X_\boR)_\boC$ makes $\Diff^{\bar\d}(X)$ a left augmented ring over
its subring $F_0\Diff^{\bar\d}(X)=O(X_\boR)_\boC$.
 So $(\Diff^{\bar\d}(X),F)$ is a left augmented left finitely projective
nonhomogeneous Koszul ring.

 The graded ring $B=\bigoplus_{q=0}^m\Omega^{0,q}(X)=
\Lambda_{O(X_\boR)_\boC}(T^*{}^{(0,1)}(X_\boR))$ is left and right
finitely projective Koszul as well.
 It is also the right finitely projective Koszul ring quadratic dual
to the left finitely projective Koszul ring~$A$, in the sense of
Corollary~\ref{homogeneous-koszul-duality-cor}.
 The (right finitely projective Koszul)
Dolbeault DG\+ring $(\Omega^{0,\bu}(X),\bar\d)$
\,\eqref{dolbeault-complex} corresponds to the (left augmented left
finitely projective Koszul) filtered ring $(\Diff^{\bar\d}(X),F)$
under the anti-equivalence of categories from
Corollary~\ref{left-augmented-koszul-duality-anti-equivalence}.

 Similarly to the smooth real manifold case of
Theorem~\ref{smooth-differential-koszul-duality-theorem},
we obtain the following rough algebraic version of Koszul duality
for the ring of $\bar\d$\+differential operators $\Diff^{\bar\d}(X)$.

\begin{thm} \label{Dolbeault-koszul-duality-theorem}
 For any compact complex manifold $X$, the construction of
Theorem~\ref{comodule-side-koszul-duality-theorem} provides
a triangulated equivalence between
the $\Diff^{\bar\d}(X)/O(X_\boR)_\boC$\+semico\-derived category of
right modules over the ring of $\bar\d$\+differential operators
$\Diff^{\bar\d}(X)$ and the coderived category of right DG\+modules
over the Dolbeault DG\+algebra $(\Omega^{0,\bu}(X),\bar\d)$,
\begin{equation}
 \sD^\sico_{O(X_\boR)_\boC}(\modr{\Diff^{\bar\d}(X)})
 \simeq\sD^\co(\modr(\Omega^{0,\bu}(X),\bar\d)),
\end{equation}
while Theorem~\ref{reduced-koszul-duality-comodule-side-thm}
establishes a triangulated equivalence between the derived
category of right $\Diff^{\bar\d}(X)$\+modules and the reduced
coderived category of right DG\+modules over
$(\Omega^{0,\bu}(X),\bar\d)$ relative to~$O(X_\boR)_\boC$,
\begin{equation}
 \sD(\modr{\Diff^{\bar\d}(X)})\simeq
 \sD^\co_{O(X_\boR)_\boC\red}(\modr(\Omega^{0,\bu}(X),\bar\d)).
\end{equation}
 In the same context, the construction of
Theorem~\ref{contramodule-side-koszul-duality-theorem} provides
a triangulated equivalence between
the $\Diff^{\bar\d}(X)/O(X_\boR)_\boC$\+semicontraderived category
of left modules over the ring $\Diff^{\bar\d}(X)$ and
the contraderived category of left DG\+modules over the DG\+algebra
$(\Omega^{0,\bu}(X),\bar\d)$,
\begin{equation} \label{Dolbeault-unreduced-contra-side}
 \sD^\sictr_{O(X_\boR)_\boC}(\Diff^{\bar\d}(X)\modl)\simeq
 \sD^\ctr((\Omega^{0,\bu}(X),\bar\d)\modl),
\end{equation}
while Theorem~\ref{reduced-koszul-duality-contramodule-side-thm}
establishes a triangulated equivalence between the derived category
of left $\Diff^{\bar\d}(X)$\+modules and the reduced contraderived
category of left DG\+modules over $(\Omega^{0,\bu}(X),\bar\d)$
relative to~$O(X_\boR)_\boC$,
\begin{equation} \label{Dolbeault-reduced-contra-side}
 \sD(\Diff^{\bar\d}(X)\modl)\simeq
 \sD^\ctr_{O(X_\boR)_\boC\red}((\Omega^{0,\bu}(X),\bar\d)\modl).
\end{equation} \qed
\end{thm}

 In particular, for every holomorphic vector bundle $E$ over $X$,
we have the left $\Diff^{\bar\d}(X)$\+module $E(X_\boR)$, as
explained above.
 The corresponding left DG\+module over the Dolbeault
DG\+ring~\eqref{dolbeault-complex}, under
the equivalences of
categories~(\ref{Dolbeault-unreduced-contra-side}\+-%
\ref{Dolbeault-reduced-contra-side}),
is the \emph{Dolbeault complex}
$(\Omega^{0,\bu}(X)\ot_{O(X_\boR)_\boC}E(X_\boR),\,\bar\d)$
\emph{with the coefficients in a holomorphic vector bundle
$E$ on~$X$}.
 The the $\boC$\+vector spaces
$H^*_{\bar\d}(\Omega^{0,\bu}(X)\ot_{O(X_\boR)_\boC}E(X_\boR))$
of cohomology of the latter complex
are the sheaf cohomology spaces $H^*(X,\mathcal E)$ of the sheaf
of holomorphic local sections $\mathcal E$ of the vector bundle~$E$.

 Furthermore, there is a natural equivalence between the abelian
categories of left and right modules over the ring $\Diff^{\bar\d}(X)$,
provided by the mutually inverse functors of tensor product with
the invertible modules $\Omega^{0,m}(X)$ and
$\Lambda_{O(X_\boR)_\boC}^m(T^{(0,1)}(X))$ over the ring
$O(X_\boR)_\boC$ (where $m=\dim_\boC X$).
 This Morita equivalence can be obtained as a particular case
of the construction of Theorem~\ref{conversion-equivalence-theorem}
with Remark~\ref{commutative-cdg-ring-remark} (since
$\Omega^{0,\bu}(X)$ is a relatively Frobenius Koszul graded ring).

\subsection{Relative differential operators}
\label{relative-diffoperators-subsecn}
 Let $\iota\:S\rarrow R$ be a homomorphism of commutative rings.
 We will consider DG\+rings $(B,d)$ endowed with a ring homomorphism
$R\rarrow B^0$ such that the following conditions are satisfied:
\begin{enumerate}
\renewcommand{\theenumi}{\roman{enumi}}
\item the graded ring $B=\bigoplus_{n\in\boZ}B^n$ is \emph{strictly} 
graded commutative (in the sense of
Section~\ref{symmetric-exterior-subsecn});
\item $B$ is a DG\+algebra over $S$, that is, in other words,
the image of the composition $S\rarrow R\rarrow B^0$ is annihilated
by the differential $d_0\:B^0\rarrow B^1$.
\end{enumerate}
 The initial object in the category of all such DG\+rings $(B,d)$ with
ring homomorphisms $R\rarrow B^0$ is called the \emph{strictly graded
commutative DG\+algebra over $S$ freely generated by~$R$} and denoted
by $(\Omega_{R/S},d)$.
 The other name for $(\Omega_{R/S},d)$ is the \emph{de~Rham DG\+algebra
of $R$ over~$S$} \cite[Tag~0FKF]{SP}.

 Clearly, the graded ring $\Omega_{R/S}$ is generated by the images
of elements from $R$ and their differentials.
 It follows that one has $\Omega^n_{R/S}=0$ for $n<0$ and
$\Omega^0_{R/S}=R$.
 Furthermore, the graded ring $\Omega_{R/S}$ is generated by its
first-degree component $\Omega^1_{R/S}$ over $\Omega^0_{R/S}=R$.
 The $R$\+module $\Omega^1_{R/S}$ can be constructed as the module
of \emph{K\"ahler differentials} of $R$ over~$S$ \cite[Tag~00RM]{SP}
(cf.\ the discussion of the particular case when $S$ is a field in
Section~\ref{algebraic-diffoperators-subsecn} above).
 Specifically, $\Omega^1_{R/S}$ is the $R$\+module generated by
the symbols $d(r)$ with $r\in R$ subject to the relations
$d(fg)=fd(g)-gd(f)$ for all $f$, $g\in R$ and
$d(\iota(s))=0$ for all $s\in S$.

 The key observation is that the graded ring $\Omega_{R/S}$ is, in
fact, the exterior algebra of the $R$\+module $\Omega^1_{R/S}$,
that is $\Omega_{R/S}=\Lambda_R(\Omega^1_{R/S})$.
 In other words, this means that there exists a well-defined odd
derivation~$d$ of degree~$1$ on the graded ring
$\Lambda_R(\Omega^1_{R/S})$ whose restriction to
$R=\Lambda^0_R(\Omega^1_{R/S})$ is the natural map
$d\:R\rarrow\Omega^1_{R/S}$ and whose square vanishes.
 Such a derivation is clearly unique, because the exterior algebra
$\Lambda^0_R(\Omega^1_{R/S})$ is generated by $\Omega^1_{R/S}$ over
$R$ and the action of~$d$ on $\Omega^1_{R/S}$ is computable as
$d(fd(g))=d(f)\wedge d(g)+fd^2(g)=d(f)\wedge d(g)$.

 Concerning the existence, one needs to check that odd derivations
defined on the generators and satisfying obvious compatibilities
extend well to freely generated strictly graded commutative rings
(cf.\ the discussion in the last paragraph of the proof of
Proposition~\ref{nonhomogeneous-dual-cdg-ring}
and Lemma~\ref{odd-derivations-lemma}).
 In particular, one can compute that $d(c^2)=d(c)c-cd(c)=0$ since $c$
commutes with~$d(c)$ for any element~$c$ of odd degree in
$\Lambda_R(\Omega^1_{R/S})$ (as it should be).
 Furthermore, one has $d^2(ab)=d(d(a)b+(-1)^{|a|}ad(b))=d^2(a)b
+(-1)^{|a|+1}d(a)d(b)+(-1)^{|a|}d(a)d(b)+ad^2(b)=0$ provided
that $d^2(a)=0=d^2(b)$ for a given pair of elements
$a$, $b\in\Lambda_R(\Omega^1_{R/S})$; so the square of an odd
derivation vanishes whenever it vanishes on the generators.

 A simpler approach may be to construct $\Omega_{R/S}$ as
the graded commutative graded ring generated by the symbols~$f$ of
degree~$0$ and $d(f)$ of degree~$1$ for all $f\in R$, with
the relations that the addition and multiplication of the symbols~$f$
in $\Omega_{R/S}$ agrees with their addition and multiplication in $R$,
and also $d(f+g)=d(f)+d(g)$, \ $d(fg)=d(f)g+fd(g)$ for all $f$,
$g\in R$, and $d(\iota(s))=0$ for all $s\in S$.
 Then the isomorphism $\Omega_{R/S}=\Lambda_R(\Omega^1_{R/S})$
follows simply from the fact that all the relations imposed have
degrees~$0$ or~$1$.
 Having observed that, one needs to convince oneself that there exists
an (obviously unique) odd derivation~$d$ on $\Omega_{R/S}$ taking $f$
to $d(f)$ and $d(f)$ to~$0$ for all $f\in R$.
 Indeed, $d(d(fg)-d(f)g-fd(g))=0+d(f)d(g)-d(f)d(g)=0$, as it should be;
so the relation is preserved by the desired odd derivation.

 When $2$ is invertible in $R$ (so there is no difference between
graded commutativity and strict graded commutativity for~$B$),
the condition~(i) can be replaced by its weaker form
\begin{enumerate}
\renewcommand{\theenumi}{\roman{enumi}$'$}
\item the image of the map $R\rarrow B^0$ is contained in
the (graded) center of~$B$.
\end{enumerate}
 The universal graded commutative DG\+ring $(E,d)$ with a map
$R\rarrow E^0$ satisfying (i$'$) and~(ii) is the same as
the universal graded commutative DG\+ring with a similar map
satisfying (i) and~(ii).

 Indeed, the graded ring $E$, being universal, is clearly generated by
$E^0=R$ and $d(R)\subset E^1$.
 So it remains to check that $d(f)d(g)+d(g)d(f)=0$ in $B$ under~(i$'$)
for all $f$, $g\in R$.
 For this purpose, it suffices to compute that $d(f)d(g)+d(g)d(f)=
d(fd(g)-d(g)f)=d(0)=0$.
 More generally, one shows that the graded center $Z$ of any graded
ring $B$ is preserved by any odd derivation~$d$ on $B$, as one has
$d(z)b=d(zb)-(-1)^{|z|}zd(b)=(-1)^{|z||b|}d(bz)-(-1)^{|z||b|}d(b)z=
(-1)^{(|z|+1)|b|}bd(z)$ for all $b\in B$ and $z\in Z$.

 Elements of the $R$\+module $T_{R/S}=\Hom_R(\Omega^1_{R/S},R)$ are
interpreted as derivations of the $S$\+algebra $R$ (i.~e.,
$S$\+linear maps $v\:R\rarrow R$ such that $v(fg)=v(f)g+fv(g)$ for
all $f$, $g\in R$).
 Specifically, given an $R$\+linear map $v\:\Omega^1_{R/S}\rarrow R$,
the action of~$v$ on $R$ is defined by the rule $v(f)=v(df)\in R$.
 In fact, this rule defines a natural $R$\+module isomorphism between
$T_{R/S}$ and the $R$\+module of all $S$\+linear derivations of~$R$.
 Consequently, the underlying $S$\+module of $T_{R/S}$ acquires
the structure of a Lie algebra over $S$: the bracket $[v,w]$ of
two elements $v$, $w\in T_{R/S}$ is defined by the usual rule
$[v,w](f)=v(w(f))-w(v(f))$.

 The construction of the de~Rham DG\+algebra $\Omega_{R/S}$ is
well-behaved for \emph{smooth} morphisms of rings $S\rarrow R$.
 In order not to delve into the intricacies of various definitions of
a smooth morphism, let us define and impose the minimal condition
that we will actually need.
 Let us say that a morphism of commutative rings $S\rarrow R$ is
\emph{weakly smooth of relative dimension~$m$} if $\Omega^1_{R/S}$
is a finitely generated projective $R$\+module everywhere of
rank~$m$ (in the sense of Section~\ref{symmetric-exterior-subsecn}).
 In this case, $T_{R/S}$ is also a finitely generated projective
$R$\+module everywhere of rank~$m$.

 Notice that the relative dimension of a weakly smooth morphism can
exceed the (relative) Krull dimension.
 For example, given a field~$k$, denote by $k(x)$ the field of
rational functions in one variable~$x$ with the coefficients in~$k$.
 Then the natural inclusion $k\rarrow k(x)$ is a weakly smooth
morphism of relative dimension~$1$.
 Moreover, for a field~$k$ of prime characteristic~$p$, the inclusion
$k(x^p)\rarrow k(x)$ is a weakly smooth morphism of relative
dimension~$1$ in the sense of our definition (even though it is
a finite, algebraic field extension).

 The ring of \emph{relative crystalline differential operators}
$\CDiff_{R/S}$ is now defined similarly to
Section~\ref{crystalline-subsecn}.
 The ring $\CDiff_{R/S}$ is generated by elements of the ring $R$
and the $R$\+module $T_{R/S}$, subject to
the relations~(\ref{crystalline-diff-relations}\+-%
\ref{commutator-vector-fields-relation}) imposed for all
$f$, $g\in R$ and $v$, $w\in T_{R/S}$.
 There is a natural structure of left $\CDiff_{R/S}$\+module on
the ring $R$, with the elements $g\in R\subset\CDiff_{R/S}$ acting
in $R$ by the multiplication maps $f\longmapsto gf$ and the elements
$v\in T_{R/S}\subset\CDiff_{R/S}$ acting in $R$ by the derivations
$f\longmapsto v(f)$.
 The operators with which the ring $\CDiff_{R/S}$ acts in $R$ are
$S$\+linear, but not $R$\+linear.
 There is also a natural increasing filtration $F$ on the ring
$\CDiff_{R/S}$ generated by $F_1\CDiff_{R/S}=R\oplus T_{R/S}\subset
\CDiff_{R/S}$ over $F_0\CDiff_{R/S}=R$.

 Assume that the morphism of commutative rings $S\rarrow R$ is
weakly smooth.
 Then the filtered ring $\CDiff_{R/S}$ together with its action in
$R$ can be also defined as the left augmented left finitely projective
projective nonhomogeneous Koszul ring corresponding to
the right finitely projective Koszul DG\+ring $(\Omega_{R/S},d)$
under the anti-equivalence of categories from
Corollary~\ref{left-augmented-koszul-duality-anti-equivalence}.
 By the Poincar\'e--Birkhoff--Witt Theorem~\ref{pbw-theorem-thm},
the associated graded ring $\gr^F\CDiff_{R/S}=
\bigoplus_{n=0}^\infty F_n\CDiff_{R/S}/F_{n-1}\CDiff_{R/S}$
is naturally isomorphic to the symmetric algebra
$\Sym_R(T_{R/S})$ of the $R$\+module $T_{R/S}$.

 Hence the results of Sections~\ref{comodule-side-secn}
and~\ref{contramodule-side-secn} specialize to the following theorem.

\begin{thm} \label{relative-differential-koszul-duality-theorem}
 For any weakly smooth morphism of commutative rings $S\rarrow R$,
the construction of Theorem~\ref{comodule-side-koszul-duality-theorem}
provides a triangulated equivalence between
the $\CDiff_{R/S}/R$\+semicoderived category of right modules over
the ring of relative crystalline differential operators $\CDiff_{R/S}$
and the coderived category of right DG\+modules over the relative
de~Rham DG\+algebra $(\Omega_{R/S},d)$,
\begin{equation} \label{relativediff-semicoderived}
 \sD^\sico_R(\modr{\CDiff_{R/S}})\simeq\sD^\co(\modr(\Omega_{R/S},d)),
\end{equation}
while Theorem~\ref{reduced-koszul-duality-comodule-side-thm}
establishes a triangulated equivalence between the derived category of
right $\CDiff_{R/S}$\+modules and the reduced coderived category of
right DG\+modules over $(\Omega_{R/S},d)$,
\begin{equation} \label{relativediff-reduced-comodule-side}
 \sD(\modr{\CDiff_{R/S}})\simeq\sD^\co_{R\red}(\modr(\Omega_{R/S},d)).
\end{equation}
 In the same context, the construction of
Theorem~\ref{contramodule-side-koszul-duality-theorem} provides
a triangulated equivalence between
the $\CDiff_{R/S}/R$\+semicontraderived category of left modules over
the ring $\CDiff_{R/S}$ and the contraderived category of left
DG\+modules over the DG\+algebra $(\Omega_{R/S},d)$,
\begin{equation} \label{relativediff-semicontraderived}
 \sD^\sictr_R(\CDiff_{R/S}\modl)\simeq\sD^\ctr((\Omega_{R/S},d)\modl),
\end{equation}
while Theorem~\ref{reduced-koszul-duality-contramodule-side-thm}
establishes a triangulated equivalence between the derived category
of left $\CDiff_{R/S}$\+modules and the reduced contraderived category
of left DG\+modules over $(\Omega_{R/S},d)$,
\begin{equation} \label{relativediff-reduced-contramodule-side}
 \sD(\CDiff_{R/S}\modl)\simeq\sD^\ctr_{R\red}((\Omega_{R/S},d)\modl).
\end{equation} \qed
\end{thm}

 Furthermore, whenever the ring $R$ has finite homological dimension,
Corollaries~\ref{fin-dim-base-koszul-duality-comodule-side}
and~\ref{fin-dim-base-koszul-duality-contramodule-side} are applicable;
so there is no difference between~\eqref{relativediff-semicoderived}
and~\eqref{relativediff-reduced-comodule-side}, and similarly
there is no difference between~\eqref{relativediff-semicontraderived}
and~\eqref{relativediff-reduced-contramodule-side}.
 Moreover, the results of Section~\ref{conversion-secn} lead to
the following theorem in this case.
 Here the natural equivalence between the abelian categories of
left and right $\CDiff_{R/S}$\+modules is provided by the mutually
inverse functors of tensor product with the invertible
$R$\+modules $\Omega^m_{R/S}$ and $\Lambda^m_R(T_{R/S})$, as
described in Theorem~\ref{conversion-equivalence-theorem}.

\begin{thm} \label{relative-differential-quadrality}
 For any weakly smooth morphism of commutative rings $S\rarrow R$
of relative dimension~$m$ such that the ring $R$ has finite
homological dimension, the constructions of
Theorem~\ref{frobenius-quadrality-theorem}
(with Remark~\ref{commutative-cdg-ring-remark})
provide a commutative square diagram of triangulated equivalences
\begin{equation}
\begin{tikzcd}
 \sD(\modr{\CDiff_{R/S}}) \arrow[r, bend right = 8] \arrow[dd] &
 \sD(\CDiff_{R/S}\modl) \arrow[l] \arrow[dd]
 \\ \\
 \sD^\co((\Omega_{R/S},d)\modl)
 \arrow[r, Leftrightarrow, no head, no tail]
 \arrow[uu, bend left = 20] &
 \sD^\ctr((\Omega_{R/S},d)\modl) \arrow[uu, bend right = 20]
\end{tikzcd}
\end{equation}
where the equivalence of derived categories in the upper line is
induced by the conversion equivalence of abelian categories
$\modr{\CDiff_{R/S}}\simeq\CDiff_{R/S}\modl$ (up to a cohomological
shift by~$[-m]$), while the co-contra correspondence in the lower
line is the result of Theorem~\ref{frobenius-co-contra-thm}.
 The vertical equivalences
are~\textup{(\ref{relativediff-semicoderived}=%
\ref{relativediff-reduced-comodule-side})}
and~\textup{(\ref{relativediff-semicontraderived}=%
\ref{relativediff-reduced-contramodule-side})}.
\end{thm}

 Given a closed relative $2$\+form~$h$, i.~e., an element
$h\in\Omega^2_{R/S}$ such that $d(h)=0$, one can also consider
the right finitely projective Koszul CDG\+ring $(\Omega_{R/S},d,h)$
and construct its nonhomogeneous quadratic dual left finitely
projective nonhomogeneous Koszul filtered ring $\CDiff_{R/S,h}$ of
\emph{twisted relative crystalline differential operators},
using Theorem~\ref{pbw-theorem-thm}
or Corollary~\ref{nonhomogeneous-koszul-duality-anti-equivalence}.
 The results of
Theorems~\ref{relative-differential-koszul-duality-theorem}
and~\ref{relative-differential-quadrality} can be then extented
to the twisted case (similarly to
Section~\ref{twisted-differential-subsecn}).

\subsection{Lie algebroids}
 The setting in this section is a common generalizations of
Sections~\ref{algebraic-diffoperators-subsecn}\+-%
\ref{crystalline-subsecn} and \ref{smooth-differential-subsecn}\+-%
\ref{relative-diffoperators-subsecn}.

 A \emph{Lie algebroid} (known also as a \emph{Lie--Rinehart
algebra})~\cite{Rin,MM} is a pair of abelian groups $(R,\g)$ endowed
with the following structures:
\begin{itemize}
\item $R$ is a commutative ring (with unit);
\item $\g$ is a Lie algebra over $\boZ$, i.~e., it is endowed with
an additive map of \emph{Lie bracket} $[{-},{-}]\:\Lambda^2_\boZ\g
\rarrow\g$ satisfying the Jacobi identity;
\item $\g$ is an $R$\+module, so a commutative ring action map
$R\ot_{\boZ}\g\rarrow\g$ is given, denoted by
$a\ot x\longmapsto ax$ for all $a\in R$ and $x\in\g$;
\item $R$ is a $\g$\+module, so a Lie action map
$\g\ot_{\boZ}R\rarrow R$ is given, denoted by
$x\ot a\longmapsto x(a)$ for all $a\in R$ and $x\in\g$.
\end{itemize}

 In addition to the usual Jacobi identity on the bracket in~$\g$,
the identity involved in the notion of a $\g$\+module,
the associativity and commutativity equations on the multiplication
in $R$, and the associativity equation involved in the notion of
an $R$\+module, the listed structures must also satisfy
the following equations:
\begin{enumerate}
\renewcommand{\theenumi}{\roman{enumi}}
\item $\g$ acts in $R$ by derivations of the commutative
multiplication, that is
$$
 x(ab)=x(a)b+ax(b) \qquad\text{for all $x\in\g$ and $a$, $b\in R$};
$$
\item the identity
$$
 (ax)(b)=ax(b) \qquad\text{for all $x\in\g$ and $a$, $b\in R$}
$$
holds in $R$, where $ax\in\g$, \ $(ax)(b)\in R$, \ $x(b)\in R$, and $ax(b)\in R$ are the elements obtained using the action of $R$ in~$\g$,
the action of~$\g$ in $R$, the action of~$\g$ in $R$, and
the multiplication in $R$, respectively;
\item the identity
$$
 [x,ay]=x(a)y+a[x,y] \qquad\text{for all $x$, $y\in\g$ and $a\in R$}
$$
holds in~$\g$, where $x(a)\in R$ is the element obtained using
the action of~$\g$ in $R$, \ $ay$ and $x(a)y\in\g$ are the elements
obtained using the action of $R$ in~$\g$, \ $[x,y]$ and $[x,ay]\in\g$
are the Lie brackets in~$\g$, and $a[x,y]\in\g$ is the element obtained
using the action of $R$ in~$\g$.
\end{enumerate}

 Notice that $\g$~is \emph{not} a Lie algebra over $R$, as the Lie
bracket in~$\g$ is not $R$\+linear.
 The identity~(iii) describes the obstacle term to $R$\+linearity
of the bracket in~$\g$.

 For example, for any smooth affine algebraic variety $X$ over
a field~$k$, the commutative ring of functions $R=O(X)$ and
the Lie algebra of vector fields $\g=T(X)$ form a Lie algebroid
with $O(X)$ acting in $T(X)$ as in the module of sections of
a vector bundle (namely, the tangent bundle) on $X$ and $T(X)$ acting
in $O(X)$ by the derivations of functions along vector fields
(cf.\ Section~\ref{crystalline-subsecn}).
 More generally, for any homomorphism of commutative rings $S\rarrow R$,
the ring $R$ and the Lie algebra $\g=T_{R/S}$ of $S$\+linear
derivations of $R$ (as in Section~\ref{relative-diffoperators-subsecn})
form a Lie algebroid.

 Furthermore, for any smooth real manifold $X$, the ring $R=O(X)$
of smooth $\boR$\+valued functions on $X$ together with the Lie
algebra $\g=T(X)$ of smooth vector fields on $X$ (as in
Section~\ref{smooth-differential-subsecn}) form a Lie algebroid.
 Similarly, for any real manifold $X$, the ring $R=O(X)_\boC$ of
smooth $\boC$\+valued functions on $X$ together with the Lie algebra
$\g'=T(X)_\boC$ of smooth $\boC$\+valued vector fields on $X$
form a Lie algebroid.
 For a complex manifold $X$, the same ring of smooth $\boC$\+valued
functions $R=O(X_\boR)_\boC$ together with the Lie subalgebra
$\g=T^{(0,1)}(X)\subset\g'=T(X_\boR)_\boC$ of smooth $(0,1)$\+vector
fields on $X$ form a Lie algebroid as well (see
Section~\ref{dolbeault-subsecn}).

 The \emph{universal enveloping ring} $U(R,\g)$ of a Lie algebroid
$(R,\g)$ is an associative ring defined by generators and relations
as follows (see~\cite{Rin,MM} for a differently worded, but equivalent
version of this construction).
 The set of generators $R\sqcup\g$ is the disjoint union of $R$
and~$\g$.
 The sum of any two elements of $R$ in $U(R,\g)$ equals their sum
in $R$, and the sum of any two elements of~$\g$ in $U(R,\g)$ equals
their sum in~$\g$.
 Denoting the product in $U(R,\g)$ by~$*$, one imposes
the multiplicative relations similar to the ones
in Section~\ref{crystalline-subsecn}:
\begin{equation} \label{lie-algebroid-enveloping-relations}
\begin{gathered}
 a*b=ab \qquad\text{for all $a$, $b\in R$,} \\
 a*x=ax \qquad\text{for all $a\in R$ and $x\in\g$,} \\
 x*a=ax+x(a) \qquad\text{for all $a\in R$ and $x\in\g$}
\end{gathered}
\end{equation}
and
\begin{equation} \label{lie-algebroid-commutator-relation}
 x*y-y*x=[x,y] \qquad\text{for all $x$, $y\in\g$}.
\end{equation}

 The commutative ring $R$ has a natural structure of left module over
the associative ring $U(R,\g)$.
 To define this action, one lets the generators $b\in R$ of
the ring $U(R,\g)$ act in $R$ by the multiplication maps $a\longmapsto
ba$ and the generators $x\in\g$ of the ring $U(R,\g)$ act in $R$
as the Lie algebra~$\g$ acts in $R$, that is $a\longmapsto x(a)$.
 Then one needs to check that the assignment of such endomorphisms
of the abelian group $R$ to the generators of $U(R,\g)$ respects
the relations~(\ref{lie-algebroid-enveloping-relations}\+-%
\ref{lie-algebroid-commutator-relation}), so the resulting
action of $U(R,\g)$ in $R$ is well-defined.

 Let $(R,\g)$ be a Lie algebroid such that the $R$\+module~$\g$
is projective and finitely generated.
 Denote by $\g\spcheck=\Hom_R(\g,R)$ the dual finitely generated
projective $R$\+module.
 Consider the symmetric algebra $A=\Sym_R(\g)$ of the $R$\+module~$\g$
and the exterior algebra $B=\Lambda_R(\g\spcheck)$.
 According to Section~\ref{symmetric-exterior-subsecn}, $A$ and $B$ are
left and right finitely projective Koszul graded rings over~$R$.
 Furthermore, the right finitely projective Koszul graded ring $B$ is
quadratic dual to the left finitely projective Koszul graded ring~$A$,
in the sense of Section~\ref{quadratic-duality-secn}
and Proposition~\ref{finitely-projective-Koszul-duality}.

 Put $\tA=U(R,\g)$, and endow the ring $\tA$ with the increasing
filtration $F$ generated by $F_1\tA=\im(R\oplus\g\to\tA)$
over $F_0\tA=\im(R\to\tA)$.
 Consider the associated graded ring
$\gr^F\tA=\bigoplus_{n=0}^\infty F_n\tA/F_{n-1}\tA$.
 Then there is a unique homomorphism of graded rings
$A\rarrow\gr^F\tA$ forming commutative triangle diagrams with
the natural isomorphisms $R\simeq A_0$ and $\g\simeq A_1$
and the natural surjective maps $R\rarrow F_0\tA$ and
$\g\rarrow F_1\tA/F_0\tA$.
 This graded ring homomorphism is obviously surjective.

 We would like to show, as a particular case of
the Poincar\'e--Birkhoff--Witt Theorem~\ref{pbw-theorem-thm}, that
the graded ring homomorphism $A\rarrow\gr^F\tA$ is, in fact,
an isomorphism.
 In particular, it will follow that the natural surjective maps
$R\rarrow F_0\tA$ and $R\oplus\g\rarrow F_1\tA$ are, in fact, bijective.
 Let us emphasize that these assertions in the conclusion of
the Poincar\'e--Birkhoff--Witt theorem only have a chance to hold due
to~(i\+-iii) and other equations imposed on the structure maps of
a Lie algebroid (such as the Jacobi identity for the bracket in~$\g$).

 Using the notation of Section~\ref{self-consistency-subsecn},
put $V=\g=A_1$ and let $I\subset V\ot_RV$ denote the kernel of
the multiplication map $A_1\ot_R A_1\rarrow A_2$.
 So $I$ is the $R$\+submodule of skew-symmetric tensors in $V\ot_RV$.
 Put $q(x,a)=x(a)$ for every $x\in V=\g$ and $a\in R$.
 Then the relations~\eqref{lie-algebroid-enveloping-relations} take
the form~\eqref{map-q}.

 Following further the notation of
Section~\ref{self-consistency-subsecn}, denote by
$\hI\subset V\ot_{\boZ}V$ the full preimage of the submodule
$I\subset V\ot_RV$ under the natural surjective map
$V\ot_{\boZ}V\rarrow V\ot_RV$.
 So we get a surjective map $\hI\rarrow I$.

 The abelian group $\hI$ is spanned by the tensors of the form
$x\ot y-y\ot x$ and $ax\ot y-x\ot ay$, where $x$, $y\in\g$
and $a\in R$.
 The rules $p(x\ot y-y\ot x)=[x,y]\in\g$ and
$p(ax\ot y-x\ot ay)=-x(a)y\in\g$ define an abelian group
homomorphism $p\:\hI\rarrow\g$.
 Let $h\:\hI\rarrow R$ be the zero map, $h=0$.
 Then the relation~\eqref{lie-algebroid-commutator-relation}
takes the form~\eqref{maps-p-and-h} for all $\hi=x\ot y-y\ot x\in\hI$,
and it follows from
the relations~\eqref{lie-algebroid-enveloping-relations}
that~\eqref{maps-p-and-h} also holds for all $\hi=ax\ot y-x\ot ay
\in\hI$, i.~e., $ax*y-x*ay=-x(a)y$ in $U(R,\g)$.

 Using the equations imposed on the structure maps of a Lie algebroid,
one can check that the maps $q$, $p$, and~$h$ satisfy
the equations~(a\+-k) in
Proposition~\ref{self-consistency-equations-prop}.
 Following the proof of
Proposition~\ref{nonhomogeneous-dual-cdg-ring}, one can then
conclude that the formulas~(\ref{d-zero-defined}\+-\ref{d-one-defined})
define an odd derivation $d\:B\rarrow B$ of degree~$1$ with
zero square, $d^2=0$.
 So one obtains a DG\+ring $(B,d)$, called the \emph{cohomological
Chevalley--Eilenberg complex} of a Lie algebroid $(R,\g)$
(with trivial coefficients).
 Alternatively, it may be easier to check directly from the definition
of a Lie algebroid that
the formulas~(\ref{d-zero-defined}\+-\ref{d-one-defined}) define
a DG\+ring structure on the graded ring $B=\Lambda_R(\g\spcheck)$.

 Now Theorem~\ref{pbw-theorem-thm} is applicable to the right finitely
projective Koszul DG\+ring $(B,d)$, and the filtered ring $(\tA,F)$
together with the above action of $\tA=U(R,\g)$ in $R$ is
the left augmented left finitely projective nonhomogeneous Koszul ring
corresponding to $(B,d)$ under the anti-equivalence of categories
from Corollary~\ref{left-augmented-koszul-duality-anti-equivalence}.
 Hence the desired isomorphism $A\simeq\gr^F\tA$.
 A more general version of this result can be found in the classical
paper~\cite[Theorem~3.1]{Rin}.

 According to Example~\ref{contramodule-side-koszul-duality-examples}(2)
and Remark~\ref{contramodule-side-koszul-dual-complex-computes},
the Chevalley--Eilenberg complex $(\Lambda_R(\g\spcheck),d)$
computes the Ext groups/ring $\Ext^*_{U(R,\g)}(R,R)$
(cf.~\cite[Section~4]{Rin}).

 The results of Sections~\ref{comodule-side-secn}
and~\ref{contramodule-side-secn} specialize to the following derived
nonhomogeneous Koszul duality theorem.

\begin{thm} \label{lie-algebroid-koszul-duality-theorem}
 For any Lie algebroid $(R,\g)$ such that the $R$\+module\/ $\g$ is
projective and finitely generated,
the construction of Theorem~\ref{comodule-side-koszul-duality-theorem}
provides a triangulated equivalence between
the $U(R,\g)/R$\+semicoderived category of right modules over
the universal enveloping ring $U(R,\g)$
and the coderived category of right DG\+modules over
the Chevalley--Eilenberg DG\+ring $(\Lambda_R(\g\spcheck),d)$,
\begin{equation} \label{lie-algebroid-semicoderived}
 \sD^\sico_R(\modr U(R,\g))\simeq
 \sD^\co(\modr(\Lambda_R(\g\spcheck),d)),
\end{equation}
while Theorem~\ref{reduced-koszul-duality-comodule-side-thm}
establishes a triangulated equivalence between the derived category of
right $U(R,\g)$\+modules and the reduced coderived category of
right DG\+modules over $(\Lambda_R(\g\spcheck),d)$,
\begin{equation} \label{lie-algebroid-reduced-comodule-side}
 \sD(\modr U(R,\g))\simeq
 \sD^\co_{R\red}(\modr(\Lambda_R(\g\spcheck),d)).
\end{equation}
 In the same context, the construction of
Theorem~\ref{contramodule-side-koszul-duality-theorem} provides
a triangulated equivalence between
the $U(R,g)/R$\+semicontraderived category of left modules over
the ring $U(R,\g)$ and the contraderived category of left
DG\+modules over the DG\+algebra $(\Lambda_R(\g\spcheck),d)$,
\begin{equation} \label{lie-algebroid-semicontraderived}
 \sD^\sictr_R(U(R,\g)\modl)\simeq
 \sD^\ctr((\Lambda_R(\g\spcheck),d)\modl),
\end{equation}
while Theorem~\ref{reduced-koszul-duality-contramodule-side-thm}
establishes a triangulated equivalence between the derived category
of left $U(R,\g)$\+modules and the reduced contraderived category
of left DG\+modules over $(\Lambda_R(\g\spcheck),d)$,
\begin{equation} \label{lie-algebroid-reduced-contramodule-side}
 \sD(U(R,\g)\modl)\simeq
 \sD^\ctr_{R\red}((\Lambda_R(\g\spcheck),d)\modl).
\end{equation} \qed
\end{thm}

 Whenever the ring $R$ has finite homological dimension,
Corollaries~\ref{fin-dim-base-koszul-duality-comodule-side}
and~\ref{fin-dim-base-koszul-duality-contramodule-side} are applicable.
 So there is no difference between~\eqref{lie-algebroid-semicoderived}
and~\eqref{lie-algebroid-reduced-comodule-side}, and similarly
there is no difference between~\eqref{lie-algebroid-semicontraderived}
and~\eqref{lie-algebroid-reduced-contramodule-side} in this case.

 Assume that the finitely generated projective $R$\+module~$\g$ is
everywhere of the same rank $m\ge0$ (in the sense of the definition
at the end of Section~\ref{symmetric-exterior-subsecn}). 
 Then there is a natural equivalence between the abelian categories of
left and right $U(R,\g)$\+modules, provided by the mutually
inverse functors of tensor product with the invertible
$R$\+modules $T=\Lambda^m_R(\g\spcheck)$ and
$\Hom_R(T,R)=\Lambda^m_R(\g)$, as described in
Theorem~\ref{conversion-equivalence-theorem}.
 Moreover, the results of Section~\ref{conversion-secn} lead to
the following relative nonhomogeneous Koszul quadrality theorem.

\begin{thm} \label{lie-algebroid-differential-quadrality}
 For any Lie algebroid $(R,\g)$ such that\/ $\g$~is a finitely
generated projective $R$\+module everywhere of rank~$m$ and
the homological dimension of $R$ is finite, the constructions of
Theorem~\ref{frobenius-quadrality-theorem}
(with Remark~\ref{commutative-cdg-ring-remark})
provide a commutative square diagram of triangulated equivalences
\begin{equation}
\begin{tikzcd}
 \sD(\modr U(R,\g)) \arrow[r, bend right = 8] \arrow[dd] &
 \sD(U(R,\g)\modl) \arrow[l] \arrow[dd]
 \\ \\
 \sD^\co((\Lambda_R(\g\spcheck),d)\modl)
 \arrow[r, Leftrightarrow, no head, no tail]
 \arrow[uu, bend left = 20] &
 \sD^\ctr((\Lambda_R(\g\spcheck),d)\modl)
 \arrow[uu, bend right = 20]
\end{tikzcd}
\end{equation}
where the equivalence of derived categories in the upper line is
induced by the conversion equivalence of abelian categories
$\modr U(R,\g)\simeq U(R,\g)\modl$ (up to a cohomological
shift by~$[-m]$), while the co-contra correspondence in the lower
line is the result of Theorem~\ref{frobenius-co-contra-thm}.
 The vertical equivalences
are~\textup{(\ref{lie-algebroid-semicoderived}=%
\ref{lie-algebroid-reduced-comodule-side})}
and~\textup{(\ref{lie-algebroid-semicontraderived}=%
\ref{lie-algebroid-reduced-contramodule-side})}.
\end{thm}

\subsection{Noncommutative differential forms}
 Let $\iota\:S\rarrow R$ be a homomorphism of associative rings.
 We will consider DG\+rings $(B,d)$ endowed with a ring homomorphism
$R\rarrow B^0$ such that the image of the composition $S\rarrow R
\rarrow B^0$ is annihilated by the differential $d_0\:B^0\rarrow B^1$.

 The initial object in the category of all such DG\+rings $(B,d)$ with
ring homomorphisms $R\rarrow B^0$ is called the \emph{DG\+ring over
$S$ freely generated by~$R$} and denoted by $(NC_{R/S},d)$.
 The other name for $(NC_{R/S},d)$ is the \emph{DG\+ring of
noncommutative differential forms for $R$ over~$S$}
(cf.~\cite[proof of Proposition~II.1]{Con}).

 This definition is the noncommutative version of the one in
the beginning of Section~\ref{relative-diffoperators-subsecn}.
 The difference is that, even when the rings $S$ and $R$ happen to be
commutative, the DG\+rings $(B,d)$ considered in this section do not
need to be graded commutative.
 In other words, we are dropping the condition~(i) or~(i$'$) of
Section~\ref{relative-diffoperators-subsecn} and keeping only
the condition~(ii).

 As in Section~\ref{relative-diffoperators-subsecn}, it is clear that
the graded ring $NC_{R/S}$ is generated by the images of elements
from $R$ and their differentials.
 It follows that one has $NC_{R/S}^n=0$ for $n<0$ and $NC_{R/S}^0=R$.
 Furthermore, the graded ring $NC_{R/S}$ is generated by its
first-degree component $NC_{R/S}^1$ over $NC_{R/S}^0=R$.
 The next lemma provides a precise, explicit description.

\begin{lem}
 The maps
\begin{equation} \label{coefficient-on-left}
 R\ot_S R/\iota(S)\ot_S R/\iota(S)\ot_S\dotsb\ot_S R/\iota(S)
 \lrarrow NC^n_{R/S}
\end{equation}
and
\begin{equation} \label{coefficient-on-right}
 R/\iota(S)\ot_S R/\iota(S)\ot_S\dotsb\ot_S R/\iota(S)\ot_SR
 \lrarrow NC^n_{R/S}
\end{equation}
given by the formulas $f\ot\bar g_1\ot\dotsb\ot\bar g_n
\longmapsto fd(g_1)d(g_2)\dotsb d(g_n)$
and $\bar g_1\ot\dotsb\ot \bar g_n\ot f\longmapsto
d(g_1)d(g_2)\dotsb d(g_n)f$ are isomorphisms of
$S$\+$S$\+bimodules for all $n\ge0$.
 More precisely, the map~\eqref{coefficient-on-left} is an isomorphism
of $R$\+$S$\+bimodules, while the map~\eqref{coefficient-on-right}
is an isomorphism of $S$\+$R$\+bimodules.

 Here the left-hand side of~\eqref{coefficient-on-left} is
the tensor product of one factor $R$ and $n$~factors $R/\iota(S)$.
 The left-hand side of~\eqref{coefficient-on-right} is
the tensor product of $n$~factors $R/\iota(S)$ and one factor~$R$.
 All the tensor products signs
in~\textup{(\ref{coefficient-on-left}\+-\ref{coefficient-on-right})}
mean tensor products of $S$\+$S$\+bimodules.
 The notation is $f$, $g_i\in R$ for all\/ $1\le i\le n$
and $\bar g_i\in R/\iota(S)$ is the image of~$g_i$ under
the natural surjection $R\rarrow R/\iota(S)$.
\end{lem}

\begin{proof}
 The maps~(\ref{coefficient-on-left}\+-\ref{coefficient-on-right})
are well-defined due to the condition that $d(\iota(s))=0$
in $NC^1_{R/S}$ for all $s\in S$.
 Having observed that, one can split the assertions of the lemma
in two parts.
 Firstly, it is claimed that the maps $R\ot_SR/\iota(S)\rarrow
NC^1_{R/S}\larrow R/\iota(S)\ot_SR$ given by the formulas
$f\ot\bar g\longmapsto fd(g)$ and $\bar g\ot f\longmapsto d(g)f$
are isomorphisms (where $f$, $g\in R$ and $\bar g\in R/\iota(S)$
is the image of~$g$).

 Secondly, notice that
\begin{multline*}
 (R\ot_SR/\iota(S))\ot_R(R\ot_SR/\iota(S))\ot_R\dotsb\ot_R
 (R\ot_SR/\iota(S)) \\ \simeq
 R\ot_SR/\iota(S)\ot_SR/\iota(S)\ot_S\dotsb\ot_SR/\iota(S).
\end{multline*}
 So, in order to deduce the assertions of the lemma for an arbitrary
$n\ge1$ from such assertions for $n=1$, one needs to show that
the multiplication map
\begin{equation} \label{noncommutative-forms-tensor-ring-map}
 NC^1_{R/S}\ot_R\dotsb\ot_RNC^1_{R/S}\lrarrow NC^n_{R/S}
\end{equation}
($n$~factors in the left-hand side) is an isomorphism for $n\ge1$.
 In other words, this means that the graded ring $NC_{R/S}$ is
the tensor ring of the $R$\+$R$\+bimodule $NC^1_{R/S}$, that is
$NC_{R/S}\simeq T_R(NC^1_{R/S})$.

 Concerning the first part, one observes that the $R$\+$R$\+bimodule
$NC^1_{R/S}$ is spanned by the symbols $d(f)$, \,$f\in R$, with
the relations $d(fg)=d(f)g+fd(g)$ for $f$, $g\in R$ and $d(\iota(s))=0$
for $s\in S$.
 It follows immediately that the maps $R\ot_SR/\iota(S)\rarrow
NC^1_{R/S}\larrow R/\iota(S)\ot_SR$ are surjective.
 In order to show that these maps are isomorphisms, it suffices to
define $R$\+$R$\+bimodule structures on $R\ot_SR/\iota(S)$ and
$R/\iota(S)\ot_SR$ in such a way that the above relations are satisfied.

 Concerning the second part, in order to prove that
the maps~\eqref{noncommutative-forms-tensor-ring-map} are
isomorphisms, it suffices to show that the map $d\:R=NC^0_{R/S}
\rarrow NC^1_{R/S}$ extends to a well-defined odd derivation with
zero square on the tensor ring $T_R(NC^1_{R/S})$.
 Lemma~\ref{odd-derivations-lemma} is a suitable tool here.

 Similarly to Section~\ref{relative-diffoperators-subsecn},
a simpler alternative approach might be to define $NC_{R/S}$
as the graded associative ring generated by the symbols $f$ and $d(f)$
with $f\in R$, subject to the relations that the addition and
multiplication of elements $f$ in $R$ agrees with their addition
and multiplication in $NC_{R/S}$, and also $d(f+g)=d(f)+d(g)$, \
$d(fg)=d(f)g+fd(g)$ for $f$, $g\in R$, and $d(\iota(s))=0$ for $s\in S$.
 Then the isomorphism $NC_{R/S}\simeq T_R(NC^1_{R/S})$ follows simply
from the fact that all the relations imposed have degrees~$0$ or~$1$.
 Subsequently one needs to see that there exists an odd derivation~$d$
on $NC_{R/S}$ taking $f$ to $d(f)$ and $d(f)$ to~$0$ for all $f\in R$.
\end{proof}

 Now let us assume that the map $\iota\:S\rarrow R$ is injective
and the right $S$\+module $R/\iota(S)$ is finitely generated and
projective.
 Then it is clear from
the isomorphism~\eqref{coefficient-on-right} that the right
$R$\+module $NC^n_{R/S}$ is finitely generated and projective
for every $n\ge1$.
 Moreover, $NC^n_{R/S}\simeq T_R(NC^1_{R/S})$ is the tensor algebra
of the $R$\+$R$\+bimodule $NC^1_{R/S}$, as we have seen
in~\eqref{noncommutative-forms-tensor-ring-map}.
 Thus the DG\+ring of noncommutative differential forms $(NC_{R/S},d)$
is right finitely projective Koszul.

 Consider the following filtered ring $(\tA,F)$.
 Let $\tA=\Hom_{S^\rop}(R,R)$ be the ring of endomorpisms of
the right $S$\+module~$R$; so $R$ is a left $\tA$\+module.
 This is our ``ring of noncommutative differential operators acting
$S$\+linearly in~$R$''.
 The left action of $R$ in itself defines a natural ring inclusion
$R\rarrow\tA$.
 So $\tA$ is a left augmented ring over its subring $R$, in
the sense of Section~\ref{augmented-subsecn}.
 Let $\tA^+\subset\tA$ be the augmentation ideal, that is,
the subgroup of all elements in $\tA$ whose action annhilates
the unit element $1\in R$.
 So $\tA^+$ is a left ideal in $\tA$, and the left $R$\+module $\tA$
decomposes naturally as the direct sum $\tA=R\oplus\tA^+$.

 Put $F_0\tA=R$ and $F_n\tA=\tA$ for all $n\ge1$.
 Then the underlying left $R$\+module of the $R$\+$R$\+bimodule
$F_1\tA/F_0\tA$ is isomorphic to~$\tA^+$.
 One easily computes this left $R$\+module as $\tA^+\simeq
\Hom_{S^\rop}(R/\iota(S),R)\simeq\Hom_{R^\rop}(R/\iota(S)\ot_SR,\>R)$.
 In our assumptions, this is a finitely generated projective
left $R$\+module.
 It follows that the associated graded ring $A=\gr^F\tA=F_0\tA\oplus
F_1\tA/F_0\tA$ is left finitely projective Koszul and quadratic dual
to the right finitely projective Koszul graded ring $B=NC_{R/S}$
(in the sense of Section~\ref{quadratic-duality-secn}
and Proposition~\ref{finitely-projective-Koszul-duality}).

 The right finitely projective Koszul DG\+ring $(NC_{R/S},d)$
corresponds to the left augmented left finitely projective
nonhomogeneous Koszul filtered ring $(\tA,F)$ under the equivalence
of categories from
Corollary~\ref{left-augmented-koszul-duality-anti-equivalence}.
 
 The results of Sections~\ref{comodule-side-secn}
and~\ref{contramodule-side-secn} specialize to the following theorem.
 Notice that the graded ring of noncommutative differential forms
$NC_{R/S}$ usually has \emph{infinitely many} nonzero grading
components in our assumptions.
 So one has to distinguish $NC_{R/S}$\+comodules and
$NC_{R/S}$\+contramodules (in the sense of
Section~\ref{comodules-and-contramodules-secn}) from objects of their
ambient categories of arbitrary $NC_{R/S}$\+modules.

 We use the natural terminology ``DG\+comodules'' and
``DG\+contramodules'' for CDG\+comodules and CDG\+contramodules,
in the sense of Sections~\ref{coderived-cdg-comodules-subsecn}
and~\ref{contraderived-cdg-contramodules-subsecn},
over a DG\+ring (i.~e., a CDG\+ring with a vanishing
curvature element~$h$).

\begin{thm} \label{noncommutative-differential-koszul-duality-theorem}
 For any injective morphism of associative rings $\iota\:S\rarrow R$
such that $R/\iota(S)$ is a finitely generated projective right
$S$\+module,
the construction of Theorem~\ref{comodule-side-koszul-duality-theorem}
provides a triangulated equivalence between
the\/ $\Hom_{S^\rop}(R,R)/R$\+semicode\-rived category of right modules
over the endomorphism ring\/ $\Hom_{S^\rop}(R,R)$
and the coderived category of right DG\+comodules over
the nonnegatively graded DG\+ring of noncommutative differential
forms $(NC_{R/S},d)$,
\begin{equation} \label{noncommutativediff-semicoderived}
 \sD^\sico_R(\modr{\Hom_{S^\rop}(R,R)})
 \simeq\sD^\co(\comodr(NC_{R/S},d)),
\end{equation}
while Theorem~\ref{reduced-koszul-duality-comodule-side-thm}
establishes a triangulated equivalence between the derived category of
right\/ $\Hom_{S^\rop}(R,R)$\+modules and the reduced coderived category
of right DG\+comodules over $(NC_{R/S},d)$,
\begin{equation} \label{noncommutativediff-reduced-comodule-side}
 \sD(\modr{\Hom_{S^\rop}(R,R)})\simeq
 \sD^\co_{R\red}(\comodr(NC_{R/S},d)).
\end{equation}
 In the same context, the construction of
Theorem~\ref{contramodule-side-koszul-duality-theorem} provides
a triangulated equivalence between
the\/ $\Hom_{S^\rop}(R,R)/R$\+semicontraderived category of left modules
over the ring\/ $\Hom_{S^\rop}(R,R)$ and the contraderived category of
left DG\+contramodules over the nonnegatively graded DG\+ring
$(NC_{R/S},d)$,
\begin{equation} \label{noncommutativediff-semicontraderived}
 \sD^\sictr_R(\Hom_{S^\rop}(R,R)\modl)
 \simeq\sD^\ctr((NC_{R/S},d)\contra),
\end{equation}
while Theorem~\ref{reduced-koszul-duality-contramodule-side-thm}
establishes a triangulated equivalence between the derived category
of left $\Hom_{S^\rop}(R,R)$\+modules and the reduced contraderived
category of left DG\+contramodules over $(NC_{R/S},d)$,
\begin{equation} \label{noncommutativediff-reduced-contramodule-side}
 \sD(\Hom_{S^\rop}(R,R)\modl)\simeq
 \sD^\ctr_{R\red}((NC_{R/S},d)\contra).
\end{equation} \qed
\end{thm}

 Furthermore, whenever the ring $R$ has finite homological dimension,
Corollaries~\ref{fin-dim-base-koszul-duality-comodule-side}
and~\ref{fin-dim-base-koszul-duality-contramodule-side} are applicable;
so there is no difference
between~\eqref{noncommutativediff-semicoderived}
and~\eqref{noncommutativediff-reduced-comodule-side}, and similarly
there is no difference
between~\eqref{noncommutativediff-semicontraderived}
and~\eqref{noncommutativediff-reduced-contramodule-side}.

\end{document}